\documentclass[oneside]{amsbook}
\usepackage{mathtools}
\usepackage{amssymb}
\usepackage[bookmarks=true,
            bookmarksnumbered=true, breaklinks=true,
            pdfstartview=FitH, hyperfigures=false,
            plainpages=true, naturalnames=true,             
            colorlinks=true,
            pagebackref=false,
	    pdfpagelabels
]{hyperref}

\usepackage{comment}
\usepackage{dsfont}
\usepackage{amsfonts}
\usepackage{graphics}
\usepackage{color}
\usepackage{amssymb}
\usepackage{verbatim}
\usepackage{amsmath}
\usepackage{esint}
\usepackage{amssymb}
\usepackage{amsthm}
\usepackage{amscd}
\usepackage{euscript}
 \usepackage{tikz}
\usepackage{mathtools}
\usepackage[all,2cell,knot,cmtip]{xy}

\UseAllTwocells

\usepackage{extarrows}
\usepackage{hyperref}
\usepackage{enumitem}
\usepackage{mathrsfs}  
\usepackage{undertilde}

\newtheorem{Theorem}{Theorem}

\newtheorem{Corollary}[Theorem]{Corollary}
\newtheorem{Lemma}[Theorem]{Lemma}
\newtheorem{Definition}[Theorem]{Definition}

\newtheorem{Notation}[Theorem]{Notation}

\theoremstyle{remark}

\newtheorem{Example}[Theorem]{Example}
\newtheorem{Remark}[Theorem]{Remark}

\newtheorem{Proposition}[Theorem]{Proposition}
\newtheorem{Fundamental Theorem}{Fundamental Theorem}

\usepackage{stmaryrd}

\makeatletter
\newcommand{\xMapsto}[2][]{\ext@arrow 0599{\Mapstofill@}{#1}{#2}}
\def\Mapstofill@{\arrowfill@{\Mapstochar\Relbar}\Relbar\Rightarrow}
\makeatother

\newcommand{\Cc}{{\mathcal{C}}}
\newcommand{\W}
{
\mathbf{{W}}
}

\newcommand{\Ext}{\mathrm{Extend}}
\newcommand{\Exp}{\mathrm{Expand}}
\newcommand{\Res} {\mathrm{Restrict}}

\newcommand{\inc}{{\mathrm{inc}}}

\newcommand{\Orb}{{\mathrm{Orb}}}

\newcommand{\PC}{{\mathrm{PC}}}

\renewcommand{\t}{{\triangleright}}

\newcommand{\proj}{\mathrm{proj}}

\newcommand{\fin}{\mathrm{fin}}
\newcommand{\etab}{\overline{\eta}}

\newbox\pullbackbox
\setbox\pullbackbox=\hbox{\xy 0;<1mm,0mm>: \POS(3,0)\ar@{-} (0,0) \ar@{-} (3,3)
\endxy}
\def\pullback{\copy\pullbackbox}

\newbox\pushoutbox
\setbox\pushoutbox=\hbox{\xy 0;<1mm,0mm>: \POS(0,3)\ar@{-} (0,0) \ar@{-} (3,3)
\endxy}
\def\pushout{\copy\pushoutbox}
\makeatother
\newcommand{\mort}[5]
{
#1 \ar@/^0.7pc/[rr]^{#3}\ar@/_0.7pc/[rr]_{#4} & \Downarrow #5 & #2
}

\newcommand{\mortS}[5]
{
#1 \ar@/^1pc/[rr]^{#3}\ar@/_1pc/[rr]_{#4} & \Downarrow #5 & #2
}

\newcommand{\mortdash}[5]
{
#1 \ar@/^0.7pc/[rr]|{\not} \ar@/^0.7pc/[rr]^{#3}\ar@/_0.7pc/[rr]_{#4} \ar@/_0.7pc/[rr]|{\not} & \Downarrow #5 & #2
}

\newcommand{\bigt}[4]
{
#1 \ar@/^1pc/[rr]^{#3}\ar@/_1pc/[rr]_{#4} &  & #2
}
\makeatletter

\newcommand{\x}{\overline{x}}
\newcommand{\y}{\overline{y}}
\newcommand{\z}{\overline{z}}

\newcommand{\trl}{\triangleleft}
\newcommand{\tl}{\blacktriangleleft}
\newcommand{\tr}{\blacktriangleright}
\newcommand{\op}{\mathrm{op}}

\DeclareMathOperator{\colim}{colim}

\makeatother
 
\makeatletter

\def \smt {\mathbf{smt}}

\newcommand{\K} {{\mathcal{K}}}
\newcommand{\sk} {{\mathrm{sk}}}

\renewcommand{\trl} {\triangleleft}

\newcommand{\X} {{\mathcal{X}}}

\newcommand{\Z} {{\mathcal{Z}}}

\newcommand{\N} {{\mathcal{N}}}
\newcommand{\J} {{\mathcal{J}}}

\newcommand{\A} {\mathcal {A}}
\newcommand{\B} {\mathcal{B}}
\newcommand{\ra} {\xrightarrow}
\newcommand{\Ra} {\xRightarrow}
\newcommand{\E} {\mathcal{E}}
\newcommand{\id} {{\rm id}}
\newcommand{\Gc} {\mathcal{G}}

\newcommand{\I}{\mathcal{I}}
\newcommand{\C} {\mathbb{C}}
\newcommand{\Q} {\mathbb{Q}}

\newcommand{\Fc} {{\mathcal{F}}}

\newcommand{\g} {{\gamma}}
\newcommand{\Sing} {\mathrm{Sing}}
\newcommand{\Comb} {\mathrm{Comb}}

\newcommand{\To} {\Rightarrow}
\newcommand{\lTo} {\Longrightarrow}
\newcommand{\fr} {\mathrm{fr}}
\newcommand{\bd} {\mathrm{bd}}
\renewcommand{\d} {\partial}

\newcommand{\cob}[1] {\mathrm{\mathbf{Cob}}{}^{#1} }
\newcommand{\tcob}[1] {\mathrm{\mathbf{2Cob}}{}^{#1}}
\newcommand{\tdcob}[2] {\overline{\mathrm{\mathbf{2Cob}}} {}^{#1}_{\mathrm{dec}}}

\newcommand{\tcobp}[1] {\mathrm{\mathbf{2Cob'}}^{#1}}

\newcommand{\trcob}[1]{\overline{\mathrm{\mathbf{2Cob}}} {}^{#1}_{\mathrm{st}}}
\newcommand{\tFQtr}[1]{\overline{\boldsymbol{{2\mathcal{Q}}}} {}_{#1}^{\,\mathrm{st} }}
\newcommand{\tFQmortr}[1]{\overline{\boldsymbol{2\mathcal{Q}}} {}^{{\mathrm{Mor}}}_{#1, \mathrm{st}}}

\newcommand{\tFQmorch}[1]{{\boldsymbol{\widehat{2\mathcal{Q}}}{}^{{\mathrm{Mor}}}_{#1}}}

\newcommand{\Bc} {{\boldsymbol{\mathbf{B}}}}

\newcommand{\FQp}[2]{\boldsymbol{\mathcal{Q}}{}_{#1}^{{ #2 } }} 
\newcommand{\FQ}[1]{\boldsymbol{\mathcal{Q}}{}_{#1}} 

\newcommand{\FRp}[1]{\boldsymbol{\mathcal{R}}{\hskip-1pt}^{(#1)}}
\newcommand{\FRh}[1]{\overline{\boldsymbol{\mathcal{R}}}{\hskip-1pt}^{(#1)}}
\newcommand{\tFQmor}[1]{{\boldsymbol{\overline{2\mathcal{Q}}}{}^{{\mathrm{Mor}}}_{#1}}}

\newcommand{\tFQ}[1]{\boldsymbol{2\mathcal{Q}}_{#1}} 
\newcommand{\tFQd}[1]{{\overline{\boldsymbol{2\mathcal{Q}}} {}_{#1}^{\mathrm{dec} }} }
\newcommand{\tFQch}[1]{\boldsymbol{\widehat{2\mathcal{Q}}}_{#1}}

\newcommand{\tLin}{{\mathrm{Lin}_2}}
\newcommand{\Lin}{\mathrm{Lin}}

\newcommand{\im}{\mathrm{im}}

\makeatother
\newcommand{\cspnu}[5]
{\xymatrix@R=2pt{\\
 #1\ar[dr]_{#2} && #3 \ar[dl]^{#4} \\
& #5
}
}
\newcommand{\cspn}[5]
{\xymatrix@R=-5pt{
 #1\ar[dr]^{#2} && #3 \ar[dl]_{#4} \\
& #5&
}
}

\newcommand{\spnt}[5]
{
\xymatrix@R=2pt
{& #5\ar[dl]_{#2}\ar[dr]^{#4}&\\
 #1 && #3 
}
}

\newcommand{\spnd}[5]
{
\xymatrix@R=-5pt
{& #5\ar[dl]_{#2}\ar[dr]^{#4}&\\
 #1 && #3
}
}

\newcommand{\cspnc}[5]
{\xymatrix@C=35pt@R=-5pt{
& #1\ar[dr]_<<<<{#2} && #3 \ar[dl]^<<<<{#4} \\
&& #5
}
}

\newcommand{\cspnd}[5]
{\xymatrix@R=-5pt{
& #1\ar[dr]_<<<<{#2} && #3 \ar[dl]^<<<<{#4} \\
&& #5
}
}

\newcommand{\cspnb}[5]
{\xymatrix@R=-5pt{
 #1\ar[dr]_<<<<{#2} && #3 \ar[dl]^<<<<{#4} \\
& #5
}
}

\newcommand{\cspna}[5]
{\xymatrix@R=13pt{
& #1\ar[dr]_<<<<{#2} && #3 \ar[dl]^<<<<{#4} \\
&& #5
}
}
\newcommand{\ca}[2]
{
\Big(\,\, \raisebox{-8pt}
{{\hskip-1cm
\xymatrix{%
    & #1 \ar@(ul,ur)^{#2} 
}
}}\,\,\Big)
}

\newcommand{\uca}[2]
{
\Big(\,\, \raisebox{-4pt}
{{\hskip-1cm
\xymatrix{%
    & #1 \ar@{-}@(ul,ur)^{#2} 
}
}}\,\,\Big)
}

\newcommand{\spn}[5]
{\xymatrix@R=2pt{
& #1\ar@{<-}[dr]_{#2} && #3 \ar@{<-}[dl]^{#4} \\
&& #5
}
}
\makeatletter

\newcommand{\Fil} {\mathbf{Fil}}
\newcommand{\CW} {\mathbf{CW}}

\newcommand{\sCW} {\mathbf{sCW}}

\newcommand{\Simp} {\mathbf{Simp}}

\newcommand{\SIMP} {\mathrm{SIMP}}

\newcommand{\Sets} {{\mathbf{Set}}}
\newcommand{\Crs} {{ \mathbf{Crs}}}
\newcommand{\Ds}  {\underline{\Delta}}

\newcommand{\CRS} {\mathbf{\mathrm{CRS}}}
\newcommand{\hpiz} {\widehat{\pi}_0}
\newcommand{\hpi}{\widehat{\pi}}

\renewcommand{\L}{{\EuScript{L}}}

\makeatletter
\newcommand*{\doublerightarrow}[2]{\mathrel{
  \settowidth{\@tempdima}{$\scriptstyle#1$}
  \settowidth{\@tempdimb}{$\scriptstyle#2$}
  \ifdim\@tempdimb>\@tempdima \@tempdima=\@tempdimb\fi
  \mathop{\vcenter{
    \offinterlineskip\ialign{\hbox to\dimexpr\@tempdima+1em{##}\cr
    \rightarrowfill\cr\noalign{\kern.5ex}
    \rightarrowfill\cr}}}\limits^{\!#1}_{\!#2}}}
\newcommand*{\triplerightarrow}[1]{\mathrel{
  \settowidth{\@tempdima}{$\scriptstyle#1$}
  \mathop{\vcenter{
    \offinterlineskip\ialign{\hbox to\dimexpr\@tempdima+1em{##}\cr
    \rightarrowfill\cr\noalign{\kern.5ex}
    \rightarrowfill\cr\noalign{\kern.5ex}
    \rightarrowfill\cr}}}\limits^{\!#1}}}
\makeatother

\newcommand{\CGWH} {{\boldsymbol{\mathrm{CGWH}}}}
\newcommand{\CG} {{\boldsymbol{\mathrm{CG}}}}

\newcommand{\TOP}  {\boldsymbol{\mathrm{TOP}}}
\newcommand{\Mor}{\boldsymbol{\mathrm{Mor}}}

\newcommand{\Vect} {{\boldsymbol{\mathrm{Vect}}}}

\newcommand{\Hp}{\mathbf{H}}
\newcommand{\tH}{\mathbf{2H}}
\newcommand{\Hpb}{\mathbf{\overline{H}}}

\newcommand{\Eta}{\boldsymbol{\eta}}

\newcommand{\Id}{\boldsymbol{\mathrm{Id}}}

\newcommand{\tHp}[1]{{\mathbf{2H}^{#1}}}
\newcommand{\tHpc}[3]{{\mathbf{2H}^{#1}_{(#2,#3)}}}
\newcommand{\tHpb}[3]{{\mathbf{\overline{2H}}^{#1}_{(#2,#3)}}}

\newcommand{\bto}{\nrightarrow}

\renewcommand{\op} {{\mathrm{op}}}

\newcommand{\soml}[2]{{\langle  {#1}, {#2} \rangle}}

\newcommand{\const} {{\mathrm{const} }}
\newcommand{\HFb} {{{\mathrm{HF}^{\mathrm{span}}}}}

\newcommand{\HFbiso} {{{\mathrm{HF}}^{\mathrm{iso}}}}
\newcommand{\Diff}[1]{{\boldsymbol{\mathrm{Diff}}^{#1}}}
\newcommand{\Otimes} {{\boldsymbol{{\otimes}}}}
\newcommand{\Times}
{
\boldsymbol{{\times}}
}
\newcommand{\Sqcup} {{\boldsymbol{{\sqcup}}}}

\newcommand{\jfm}[1]{\textcolor{magenta}{#1}}

\newcommand{\Sdp} {{\boldsymbol{\Sigma'}}}

\newcommand{\Sd} {{\boldsymbol{\Sigma}}}
\newcommand{\fd} {\overline{f}}
\newcommand{\gd} {\overline{g}}

\newcommand{\Grp} {\boldsymbol{\mathrm{Grp}}}

\newcommand{\Prof} {\boldsymbol{\mathrm{Prof}}}
\newcommand{\vProf} {\boldsymbol{\mathrm{vProf}}}

\newcommand{\vProfGrp}{{\boldsymbol{\mathrm{vProf}}_{\boldsymbol{\Grp}}}}
\newcommand{\vProfGrphf}{{\boldsymbol{\mathrm{vProf}}_{\boldsymbol{\hf}}} }
 \newcommand{\vProfGrpfin}{{\boldsymbol{\mathrm{vProf}}_{\boldsymbol{\fin}}}}

\newcommand{\hf} {\mathrm{hf}}
\renewcommand{\fin} {\mathrm{fin}}

\usepackage[normalem]{ulem}

\renewcommand{\P}{{\mathcal{P}}}

\usepackage{pinlabel}
\newcommand{\Mu}{\mathsf{M}}
\newcommand{\Du}{\mathsf{D}}
\newcommand{\Au}{\mathsf{A}}
\newcommand{\Nu}{\mathsf{N}}

\newcommand{\T}{\mathcal{T}}

\newcommand{\oo}{\mathrm{or}}

\newcommand{\Mp}{\EuScript{M}}

\newcommand{\Np}{\EuScript{N}}

\newcommand{\Ap}{\EuScript{A}}

\newcommand{\Ip}{\EuScript{I}}

\newcommand{\Bp}{\EuScript{B}}

\newcommand{\Ep}{\EuScript{E}}

\usepackage{todonotes}

\begin{document}

\counterwithout*{section}{chapter}

\frontmatter

\title[A categorification of Quinn's TQFT]
{A categorification of Quinn's finite total homotopy  TQFT with application to TQFTs and once-extended TQFTs derived from  strict omega-groupoids}

\author[Faria Martins]{Jo\~ao Faria Martins}
\email[Jo\~ao Faria Martins]{j.fariamartins@leeds.ac.uk}

\author[Porter]{Timothy Porter}
\email[Timothy Porter ]{t.porter.maths@gmail.com
}

\begin{abstract}

 We first revisit the construction of Quinn's finite total homotopy TQFT, which depends on the choice of a homotopy finite space, $\Bc$. This constitutes a vast generalisation of the Dijkgraaf-Witten TQFT, with a trivial cocycle, and of Yetter's homotopy 2-type TQFT. We build our construction directly from homotopy theoretical techniques, and hence, as in Quinn's original notes from 1995, the construction works in all dimensions.

Our aim in this is to provide background for giving in detail the construction of a once-extended TQFT categorifying Quinn's finite total homotopy TQFT, in the form of a symmetric monoidal bifunctor from the bicategory of manifolds, cobordisms and extended cobordisms, first to the symmetric monoidal bicategory of profunctors (enriched over vector spaces), and then to the Morita bicategory of algebras, bimodules and bimodule maps. These once-extended versions of Quinn's finite total homotopy TQFT likewise are defined for all dimensions, and, as with the original version, depend on the choice of a homotopy finite space, $\Bc$.

To show the utility of this approach, we explicitly compute both Quinn's finite total homotopy TQFT, and its extended version, for the case when $\Bc$ is the classifying space of a homotopy finite omega-groupoid, in this paper taking the form of a crossed complex, following Brown and Higgins.

The constructions in this paper include, in particular, the description of once-extended TQFTs derived from the classifying space of a finite strict 2-group, of relevance for modelling     discrete higher gauge theory, although  the techniques involved are  considerably more general.

\end{abstract}

\date{\today}

\dedicatory{{\quad \\ This paper is dedicated to  Ronnie Brown, born 4 January 1935, died 6 December 2024.
The influence of his work in this paper is enormous.}}

\address[Faria Martins]{{School of Mathematics, University of
Leeds, Leeds, LS2 9JT, United Kingdom}}
\address[Porter]{Ynys M\^{o}n /  Anglesey, Cymru / Wales, ex-University of Bangor, United Kingdom.}



\maketitle

 \setcounter{tocdepth}{2}

\tableofcontents

\mainmatter

\section{Introduction}


In  Lecture 4 of his lecture notes,  \cite[1995]{Quinn},  on axiomatic topological quantum field theory, Quinn described, what he called the \emph{finite total homotopy TQFT},  following on from  a suggestion of Kontsevich, \cite{Kontsevich:rational:1988}.  This family of TQFTs, whose construction works in any spatial dimension, had, in special cases, been studied by Dijkgraaf and Witten, \cite{DW}, and also by Segal, \cite{segal:definition:1988}.  In that lecture, Quinn sketches the construction, starting from  a space, $\Bc$, which has `finite total homotopy' or, as we will say, `is homotopy finite'.
This finiteness  is to ensure that the resulting theory takes values in the category of finite dimensional vector spaces.

The basic construction used is quite simple in its main idea. Let ${n}$ be any non-negative integer. Let $\cob{{n}}$ be the symmetric monoidal category of closed smooth ${n}$-manifolds, and diffeomorphism classes of $(n+1)$-cobordisms between them. Given $\Bc$, and a (smooth and closed) ${n}$-manifold, ${{\Sigma}}$, then the TQFT,
which we will denote by $\FQ{\Bc} \colon \cob{{n}} \to \Vect_\mathbb{Q},$ assigns $\mathbb{Q}[{{\Sigma}},\Bc]$ to ${{\Sigma}}$.  The vector space corresponding to ${{\Sigma}}$ is thus based on the set, $[{{\Sigma}},\Bc]$, of homotopy classes of maps from ${{\Sigma}}$ to $\Bc$.  Important examples are when $\Bc$ is the classifying space of a finite group, or of a  finite (strict) 2-group, where one retrieves well known examples of TQFTs, such as Dijkgraff-Witten's TQFT, \cite{DW}, and the Yetter-Porter TQFT, \cite{Martins_Porter,Porter,Yetter}, but there are many other possibilities.  Given a cobordism, $M$, from ${{\Sigma}}$ to another manifold, ${{\Sigma'}}$, the construction gives a matrix, and hence a linear transformation from $\FQ{\Bc}({{\Sigma}})$ to $\FQ{\Bc}({{\Sigma'}})$, the matrix being with respect to the given bases of the two vector spaces.

Our main purpose in this paper is to \emph{categorify} this construction of Quinn to get what we call the \emph{once-extended Quinn TQFT}, which will be formulated in three different, closely related, ways.

What do we mean by `categorification'?  In very general terms, when categorifying a theory, one wants to try to replace sets by categories or groupoids, categories themselves by 2-categories, or better bicategories, functions between sets by functors between categories,  etc., and, when all that is done, to add another layer corresponding to natural transformations.  Here, for instance, we want to replace the category, $ \cob{{n}}$, by a bicategory / weak 2-category, $\tcob{{n}}$, incorporating a form of 2-cobordism, or cobordism between cobordisms, between manifolds. We want to replace $[{{\Sigma}},\Bc]$, which is the same as $\pi_0(\Bc^{{\Sigma}})$, the set of {path}-components of the mapping space\footnote{which may sometimes be more conveniently written as  $\TOP({{\Sigma}},\Bc)$.}, $\Bc^{{\Sigma}}$, by the fundamental groupoid, $\pi_1(\Bc^{{\Sigma}},\Bc^{{\Sigma}})$, and then do corresponding adaptations of Quinn's methodology to obtain $\Vect$-valued profunctors between these groupoids, associated to cobordisms. Finally, we then construct appropriate transformations between profunctors  to be associated to extended cobordisms connecting cobordisms. All this we want to linearise, forming the free linear categories on the groupoids, etc.

To make all this work, we need to start by taking apart Quinn's original method, and, noting that in the published version, \cite{Quinn}, a lot is merely sketched, we have included a more detailed rendition of his theory, and, in fact, will give a parametrised family of variants of his theory. Partially because of this, we work not only over $\mathbb{Q}$, but over a more general subfield, $\kappa$, of $\mathbb{C}$, as on occasion we will need the extra freedom that this gives us.

In developing this theory, we hoped that it would allow calculations that will generalise known ones, and, to this end, we develop methods of \emph{explicit calculation}, of both Quinn's finite total homotopy TQFT, and its categorified versions, in a particular family of cases, namely when $\Bc$ is the classifying space of a homotopy finite strict $\infty$-groupoid, using the crossed complex model of such algebraic objects, as developed by Brown--Higgins--Sivera, \cite{brown_higgins_sivera}, and Tonks, \cite{tonks:JPAA:2003}.  In particular, our framework includes the case when $\Bc$ is the classifying space of a strict 2-group, which is relevant for understanding TQFTs and extended TQFTs derived from discrete higher gauge theory, \cite{Baez_Huerta:2011,Companion}.

The framework for explicit calculations developed here can likely be extended in order to allow for combinatorial calculation of Quinn's finite total homotopy TQFT, and its once-extended versions, whenever the homotopy finite space, $\Bc$, is represented combinatorially, for instance when $\Bc$ is the classifying space of a finite simplicial group. (Note that finite simplicial groups are considerably more general than crossed complexes, and do model all homotopy finite spaces \cite{Ellis}.) Such a study will be deferred to a future paper.

We also expect that the categorification constructed here of Quinn's finite total homotopy TQFT, to a once-extended TQFT, can be further categorified to a doubly-extended, perhaps even fully-extended TQFT, \cite{Lurie}. This analysis will likewise be deferred to a subsequent paper.

\

\centerline{\sc{... in a bit more detail} }

\subsection{The `classical' Quinn finite total homotopy TQFT}
Let $\Vect_\kappa=\Vect$ be the category of $\kappa$-vector spaces and linear maps, which will usually  be considered together with its usual symmetric monoidal category structure. Throughout the paper, given a non-negative integer, ${n}$, the symmetric monoidal category of closed ${n}$-manifolds, and equivalence classes of cobordisms between them, will be denoted $\cob{{n}}$. By a \emph{${n}$-dimensional TQFT}, or an \emph{$(n,n+1)$-TQFT,} we will mean a symmetric monoidal functor from   $\cob{{n}}$  to $\Vect$, as in, for instance,   \cite{Lurie,Lec_TQFT}. Note that there is no assumption made, nor needed, that our manifolds or cobordisms be  oriented, or even orientable. However $\kappa$ must  be a subfield of $\C$ as we need to be able to invert positive integers.

In this paper, we will need to work over the category, $\CGWH$, of compactly generated, weak Hausdorff spaces, \cite{May,Strom}. Such a space, $B$, is called \emph{homotopy finite} if $B$ has only a finite number of path-components, each of which has only a finite number of non-trivial homotopy groups, all of which are finite, \cite{Ellis}. Also recall that given a homotopy finite space, its \emph{homotopy content}, $\chi^\pi(B)$, is defined by   the formula below (cf. also \cite[Lecture 4]{Quinn}, \cite{BaezDolan} and \cite[\S 3]{Galvezetal}),
$$\chi^\pi(B) = \sum_{[x] \in \pi_0(B)} 
\frac{ \big |\pi_2(B,x) \big |\, \big |\pi_4(B,x) \big |\, \big |\pi_6(B,x) \big |\ldots}{ \big |\pi_1(B,x) \big |\, \big |\pi_3(B,x) \big |\,  \big |\pi_5(B,x) \big |\ldots} \in \Q.$$

Given a  fixed homotopy finite space $\Bc$,  Quinn  defined what he called the \emph{finite total homotopy TQFT}, denoted here by $\FQ{\Bc} \colon \cob{{n}} \to \Vect_\mathbb{Q}$, defined for all $n\ge 0$. In Chapter \ref{Quinn-rend} of this paper
, we provide a thorough description of the construction of Quinn's finite total homotopy TQFT, giving full mathematical details, in particular defining, more generally, a parametrised version, $\FQp{\Bc}{s} \colon \cob{{n}} \to \Vect_\mathbb{\C}$, of it, where $s$ is a complex parameter. All TQFTs, $\FQp{\Bc}{s}$, for fixed $\Bc$, are related by  natural isomorphisms. The latter parameter, $s$, was not present in Quinn's original construction,  It is closely related to the parameter  appearing in the  similar `degroupoidification' functor in \cite{baezetal}.

Chapter \ref{Quinn-rend} of this paper is subdivided into two sections. In Section \ref{Classical Quinn}, we show the construction of Quinn's finite total homotopy TQFT. Prior to that, in Section \ref{Sec:prelim Quinn}, which is considerably more technical, we formulate the necessary homotopy-theoretical setting underpinning the construction of the TQFT. In particular, we introduce one of the main technical tools used  in this paper, which is the idea of a \emph{fibrant span}, $(p,M,p')\colon B \to B'$, of homotopy finite spaces, meaning that we have a diagram of homotopy finite spaces,  \begin{equation*} \vcenter{\xymatrix@R=-5pt{ &&M\ar[dl]_p\ar[dr]^{p'}\\  & B & &B',  }}
 \end{equation*}
 and, crucially, the   induced map, $\big \langle p,p'\big \rangle\colon M \to B\times B'$, is a  Hurewicz fibration.
 
 These fibrant spans of homotopy finite spaces can be composed by performing the obvious pullback.
Moreover, we have a category, $\HFb$, whose objects are homotopy finite spaces, and morphisms are fibred homotopy  classes of fibrant spans connecting them\footnote{A `dual' category of cofibrant cospans is treated in \cite{Torzewska,Torzewska-HomCobs}.}. The identity on a homotopy finite space, $B$, in $\HFb$, is given by the fibred homotopy class of the  fibrant span of homotopy finite spaces, $(s_B,B,t_B)\colon B \to B$, obtained from the path-space fibration, that is, \begin{equation*} \vcenter{\xymatrix@R=-5pt{ &&B^I\ar[dl]_{s_B}\ar[dr]^{_{t_B}}\\  & B & &B.  }}
 \end{equation*}
 Here $B^I$ is the space of all maps  from the unit interval, $I=[0,1]$, to $B$, and $s_B(\gamma)=\gamma(0)$ and $t_B(\gamma)=\gamma(1)$.

A crucial step for the construction of Quinn's finite total homotopy TQFT is a  family of functors, $\FRp{s}\colon \HFb\to \Vect$, (working now over $\C$), sending each homotopy finite space, $B$, to the free vector space, $\C(\pi_0(B))$, over the set $\pi_0(B)$,  and, given a fibrant span,  $(p,M,p')\colon B \to B'$, the matrix elements of the linear map, $\FRp{s} (p,M,p')\colon \FRp{s}(B) \to \FRp{s}(B')$,  are given by the equation
 \begin{multline*}
\big \langle \PC_b(B) \big | \FRp{s}(p,M,p')\big | \PC_{b'}(B') \big \rangle\\=\chi^\pi\big( \langle p,p'\rangle^{-1}(b,b') \big)\,\, \big({\chi^\pi(\PC_{b}{(B)})}\big)^{s}\,\, \big ({\chi^\pi(\PC_{b'}{(B')})}\big)^{1-s}.
\end{multline*}
 Here $b\in B$, $b'\in B'$, and we  denote the path component of $b$ in $B$ by $\PC_b(B)$, and the same for $b'$.
The existence of the  functor, $\FRp{s}\colon \HFb\to \Vect$, is implicit in the construction in \cite{Quinn}, and it is also  addressed by  G\'{a}lvez-Carrillo, Kock, and Tonks, \cite{Galvezetal}, albeit in the framework  of $\infty$-groupoids, and using homotopy pullbacks instead of the usual pullbacks (which we can use here, since we are working with fibrant spans only).
This also generalises the ``degroupoidification'' functor in \cite{baezetal}.

Let $\Bc$ be a fixed homotopy finite space. The final step of the construction of Quinn's finite total homotopy TQFT,  $\FQp{\Bc}{s}\colon \cob{{n}} \to \Vect$, made explicit in Section \ref{Classical Quinn}, relies on the existence of a functor,    $\Fc_\Bc\colon \cob{{n}} \to \HFb$, sending a ${n}$-manifold, $\Sigma$, to the space, $\Bc^\Sigma$, of continuous functions from $\Sigma$  to $\Bc$, and sending  the equivalence  class of an $(n+1)$-cobordism, $(i,S,j)\colon \Sigma \to \Sigma'$, from $\Sigma$ to $\Sigma'$, seen as a cospan,
 $$ \cspn{\Sigma}{i}{\Sigma'}{j}{S},$$
 to the equivalence class of the following fibrant span of homotopy finite spaces,
  \begin{equation*} \vcenter{\xymatrix@R=-5pt{ &&\Bc^S\ar[dl]_{i^*}\ar[dr]^{j^*}\\  & \Bc^\Sigma & &\Bc^{\Sigma'}.  }}
 \end{equation*}
 (Here given $f\colon S \to \Bc$, then $i^*(f)=f\circ i$ and $j^*(f)=f\circ j$.)
 This allows us to define Quinn's finite total homotopy TQFT, $\FQp{\Bc}{s} \colon \cob{{n}} \to \Vect$,  as the composite  functor,
 $$\cob{{n}} \ra{\Fc_\Bc}   \HFb \ra{ \FRp{s}} \Vect. $$
This functor,  $\FQp{\Bc}{s} \colon \cob{{n}} \to \Vect$, can be canonically given the structure of a symmetric monoidal functor, and hence of a TQFT.

In Section \ref{Classical Quinn}, we also show some properties of  $\FQp{\Bc}{s}$,  as we change $\Bc$.  For instance, we show that $\FQp{\Bc}{s}$ depends only on the homotopy type of $\Bc$, up to natural isomorphisms, which, given that $\FQp{\Bc}{s}$ is not functorial in $\Bc$, is  not immediate. Furthermore, given a homotopy finite space, $\Bc$, we have an action of the group, $\E(\Bc),$ of homotopy classes of homotopy equivalences of $\Bc,$ on the TQFTs $\FQp{\Bc}{s} \colon \cob{{n}} \to \Vect$, by natural isomorphisms. These are given in Theorems \ref{Th:changeB} and \ref{Thm:EB}.

We will show, much later on, in Subsection \ref{Quinn for finite A},  how to  calculate $\FQp{\Bc}{s}$, explicitly, in the case in which $\Bc$ is the classifying space of a strict $\infty$-groupoid, \cite{AraMetayer}, noting that such a structure is often, more classically,  called an $\omega$-groupoid as in \cite{brown_higgins_sivera}, and is often considered (as we will do here) in its form as a crossed complex, in the sense of \cite{brown_higgins_sivera}. The cases for  $\FQp{\Bc}{s}$ treated here, when $\Bc$ is the classifying space of a crossed complex,  include  that in which $\Bc$ is the classifying space of a finite 2-group, and are therefore relevant for understanding discrete higher gauge theory, \cite{Companion,Morton_Picken_2groupactions_and_moduli}. The explicit formulae that we will give for the TQFTs derived from crossed complexes, via the special case of 2-groups / crossed modules, complement and generalise those of our earlier paper, \cite{Martins_Porter}. The latter reference  considered only invariants of closed manifolds derived from crossed complexes, with cohomology classes, however providing a homotopy interpretation of the Yetter homotopy 2-type TQFT \cite{Yetter,Porter} for the case of closed manifolds.

\subsection{The once-extended versions of Quinn's finite total homotopy TQFT}

Chapter \ref{part:once_ext} of this paper, which again is subdivided into two sections, is devoted to the construction of once-extended versions of Quinn's finite total homotopy TQFT.

Let ${n}$ be a non-negative integer.
We let $\tcob{{n}}$ be the symmetric monoidal bicategory of ${n}$-dimensional closed (and smooth) manifolds, $(n+1)$-cobordisms between closed ${n}$-manifolds, and  diffeomorphism classes of $(n+2)$-extended cobordisms connecting $(n+1)$-cobordisms; see \cite{Schommer-Pries}, for instance, for precise definitions.

We let $\Prof_\kappa=\vProf$ be the symmetric monoidal bicategory with objects small linear categories (meaning categories enriched over $\Vect$). Given two small linear categories, $C$ and $C'$, 1-morphisms from $C$ to $C'$ are  $\Vect$-enriched profunctors, $\Hp\colon C \bto C'$, so, according to our conventions, they are enriched functors, $\Hp\colon C^{\op}\times C' \to \Vect$. The 2-morphisms  are natural transformations of such enriched functors from $C^\op\times C'$ to $\Vect$; see e.g. \cite{Bartlett_etal,Hansen-Shulman:constructing:2019} for complete definitions. 
 We review the definition of the bicategory, $\vProf$, in Subsection \ref{sec:prof_conventions}.

In this paper,  a \emph{once-extended TQFT}, {sometimes called here,} more briefly and more vaguely, an \emph{extended TQFT},  or an \emph{$(n,n+1,n+2)$-extended TQFT}, will be, by definition, a symmetric monoidal bifunctor, $\tcob{{n}}\to \vProf,$ thought of as a categorified version of a classical TQFT.
We will also consider once-extended TQFTs formulated as  symmetric monoidal bifunctors,
$\tcob{{n}}\to \Mor,$
where $\Mor$ is the bicategory of algebras, bimodules and bimodule maps. In this paper, the latter constructions however always originate from bifunctors with target $\vProf$ by `linearisation', in a similar way to \cite{morton:cohomological:2015}.

We will work with the full sub-bicategory,  $\vProfGrp$, of $\vProf$, whose objects are groupoids, $G=(s,t\colon G_1 \to G_0)$, each  made into a linear category  by applying the free vector space functor, $\Lin\colon \Sets \to \Vect$, to the $\hom$-sets of $G$.   Given groupoids, $G$ and $G'$, 1-morphisms in $\vProfGrp$, from $G$ to $G'$, are hence, by definition, \emph{$\Vect$-profunctors},  $\Hp\colon G\bto G'$,  in our conventions these being functors, $\Hp\colon G^\op\times G' \to \Vect$.  Given $\Vect$-profunctors, $\Hp,\Hp'\colon G^{\op}\times G' \to \Vect$, a 2-morphism, $\Eta\colon \Hp \To \Hp'$, in $\vProfGrp$, is a natural transformation of functors, $G^{\op}\times G' \to \Vect$.

There are two  sub-bicategories of  $\vProfGrp$ that we will consider here, namely  $\vProfGrphf$,  the full sub-bicategory of $\vProf$, whose objects are the homotopy finite groupoids, and  $\vProfGrpfin$, the full sub-bicategory of $\vProfGrp$ whose objects are the finite groupoids.
\subsection{The first version of the construction of the once-extended Quinn TQFT}\label{intro_2HF}
Our first construction of a once-extended version of Quinn's finite total homotopy TQFT is a (symmetric monoidal) bifunctor,
$$\tFQ{\Bc}\colon \tcob{{n}} \to \vProfGrphf,$$
constructed in Section \ref{sec:def-once-extended}, particularly Subsection \ref{sec:once extended}.

As in the construction of Quinn's finite total homotopy TQFT, we will factor $\tFQ{\Bc}$ by an intermediate homotopy-theoretical construction, which we describe in Section \ref{sec:homotopy_underpinning}, where we develop most of the  homotopy-theoretical underpinning  for the  once-extended  Quinn TQFT.  In particular, we consider a bicategory-like object (however not quite a bicategory), denoted $\mathbf{2span}(HF)$. The objects of $\mathbf{2span}(HF)$  are homotopy finite spaces, and the 1-cells of $\mathbf{2span}(HF)$,  from $X$ to $Y$, are fibrant spans, $(p,M,p')\colon X \to Y$, of homotopy finite spaces. The (not strictly associative) composition, $\bullet$, of 1-cells in $\mathbf{2span}(HF)$ is again  given by the obvious pullback. Each homotopy finite space, $X$, has its `horizontal identity' given by the path-space fibrant span, $(s_X,X^I,t_X)\colon X \to X$.

 Given two fibrant spans of homotopy finite spaces, $(p,M,p'),(q,N,q')\colon X \to Y$, 2-cells in $\mathbf{2span}(HF)$, connecting them, are given by
\emph{homotopy finite fibrant resolved 2-spans},  $\W\colon (p,M,p')\To (q,N,q')$, by definition consisting of diagrams in the category of topological spaces, as shown below,
 \begin{equation}\label{eq:W-intro}
 \W= \vcenter{ \xymatrix@C=40pt@R=20pt{
  X   &  M\ar[l]_p\ar[r]^{p'}  & Y  \\
  X^I\ar[u]^{s_X} \ar[d]_{t_X} & L \ar[l]_{l_X}\ar[r]^{r_Y}\ar[u]^P\ar[d]_Q& Y^I \ar[u]_{s_Y} \ar[d]^{t_Y}\\
  X  &  N\ar[l]^q\ar[r]_{q'}  & Y.}}
  \end{equation}
Here,  $X,Y, M,N, L$  are  homotopy finite spaces, and, crucially for the construction to work, the induced map below, called the \emph{filler of $\W$},  is a Hurewicz fibration:
 \begin{equation}\label{PL-intro}
L \stackrel{P_L}{\longrightarrow}  M  \underset{X \times Y}{\times} (X^I \times Y^I)  \underset{X \times Y}{\times} N.\end{equation}
(On the right-hand-side, we have the obvious pullback arising from the limit along the exterior faces of the diagram defining $\W$.)
Again, homotopy finite fibrant resolved $2$-spans  compose horizontally and vertically, though not associatively. 

In the crucial cases arising in the once-extended Quinn TQFT, some particular  1-cells, $(p,M,p')\colon X \to Y$, have `vertical units', $\id_{(p,M,p')}\colon (p,M,p') \To (p,M,p')$, and we, moreover, have `unitor 2-cells', all inside  $\mathbf{2span}(HF)$, denoted
\begin{align*}
 \boldsymbol{\rho}^{(p,M,p')}_{Y} \colon (X \ra{(p,M,p')} Y) \bullet (Y \ra{ (s_Y,Y^I,t_Y)} Y)   &\lTo (X \ra{(p,M,p')} Y),\\ \intertext{and}
\boldsymbol{\lambda}^{(p,M,p')}_{X} \colon (X \ra{(s_X,X^I,t_X)} X) \bullet (X \ra{(p,M,p')} Y)  &\lTo (X \ra{(p,M,p')} Y).
\end{align*}

Throughout Section \ref{sec:homotopy_underpinning}, we construct an `assignment',  $$\boldsymbol{\mathcal{H}}=\big(\pi_1(-,-),\Hp,\tH) \big)\colon \mathbf{2span}(HF) \to \vProfGrphf,$$
that gives the following.
\begin{enumerate}[leftmargin=1cm]
 \item Each homotopy finite space, $X$, is sent to its fundamental groupoid $\pi_1(X,X)$.
 \item 
 Given a, homotopy finite, fibrant span, $(p,M,p')\colon X \to Y$, we have a $\Vect$-profunctor,  $$\smash{\Hp\big (X \ra{(p,M,p')} Y \big) \colon \pi_1(X,X)^{\op} \times \pi_1(Y ,Y ) \to \Vect,}$$
 such that:
  \begin{enumerate}[leftmargin=0.5cm]
                \item  given  $x \in X$ and $y\in Y $, the (by construction, finite dimensional) vector space
                $\Hp\big ((p,M,p')\colon X \to Y \big)(x,y)$ is the free vector space on the path components of the fibre $\langle p,p' \rangle^{-1}(x,y)\subseteq M $, of $\langle p,p' \rangle\colon M \to X\times Y$.
                \item 
                Given morphisms  in $\pi_1(X,X)$ and $\pi_1(Y ,Y )$, i.e. equivalence classes of paths, in $X$ and $Y$, under path-homotopy,
$[\gamma^X] \colon x\to x'$ and $[\gamma^Y] \colon y\to y',$ the linear map,
$$ \Hp\big(X \ra{(p,M,p')} Y \big)\big([\gamma^X],[\gamma^Y]\big)\colon \Hp\big (X \ra{(p,M,p')} Y \big) (x',y) \to \Hp\big( X \ra{(p,M,p')} Y \big) (x,y'),$$ is induced by any of the homotopy equivalences, between fibres,
$$\langle p,p' \rangle^{-1}(x',y)\to \langle p,p' \rangle^{-1}(x,y'),$$
arising by applying the homotopy lifting property of $\langle p,p' \rangle\colon M \to X\times Y$ to $\overline{\gamma}^X$ and $\gamma^Y$, together. Here, $\overline{\gamma}^X$ is the reverse path of $\gamma^X$.
\end{enumerate}

\item Finally, given $\W\colon (p,M,p')\To (q,N,q')$, as in \eqref{eq:W-intro} above, we have a natural transformation of functors $\pi_1(X,X)^{\mathrm{op}}\times \pi_1(Y,Y) \to \Vect$, denoted $$\smash{\tHp{\W} \colon \Hp\big ( X\ra{(p,M,p')} Y\big)  \lTo \Hp\big ( X\ra{(q,N,q')} Y\big).}$$ Explicitly, given objects $x \in X$ and $y \in Y$, the linear map,
$$\smash{\tHpc{\W}{x}{y}\colon \Hp\big (X \ra{(p,M,p')} Y \big)(x,y) \to  \Hp\big(X \ra{(q,N,q')} Y \big)(x,y),}$$ has  the  following matrix elements, if
 $m\in \langle p,p' \rangle^{-1}(x,y) $, $n\in \langle q,q' \rangle^{-1}(x,y)
$, and where $P_L$ is the filler of $\W$, defined in Equation \eqref{PL-intro},
\begin{multline*}
 \big \langle \PC_{m}\big( \langle p,p' \rangle^{-1}(x,y) \big) \mid \tHpc{\W}{x}{y} \mid  \PC_{n}\big( \langle q,q' \rangle^{-1}(x,y) \big) \big\rangle
 \\= \chi^\pi \big( P_L^{-1}(m,\const_x,\const_y,n ) \big) \,\, \chi^\pi\big (\PC_n(  \langle q,q' \rangle^{-1}(x,y) ) \big). 
\end{multline*}
Here, $\const_x$ and $\const_y$ are the constant paths at $x$ and $y$.

 The proof that $\tHp{\W}$, defined this way, is indeed a natural transformation is far from being immediate and requires a wealth of careful verifications.
\end{enumerate}

Our main result in Section \ref{sec:homotopy_underpinning} is that  the assignment,  $$\boldsymbol{\mathcal{H}}\colon \mathbf{2span}(HF) \to \vProfGrphf,$$  preserves all various compositions, and the horizontal identities, up to applying appropriate natural isomorphisms. Moreover,  vertical units and unitors are preserved by $\boldsymbol{\mathcal{H}}$, whenever they exist. The hardest calculation is that indeed the natural transformations, $\tHp{\W}$, are well behaved with respect to the horizontal composition of fibrant resolved 2-spans of homotopy finite spaces. This is done in \S \ref{Hor-comp-resolved 2-cells}.

Having developed the homotopy-theoretical framework for the  once-extended Quinn TQFT, Section \ref{sec:def-once-extended} is devoted to its explicit construction, in three different forms.
Let $\Bc$ be a homotopy  finite space. Similarly to the case of Quinn's finite total homotopy TQFT, we have an assignment, $\Bc^{(-)}\colon\tcob{n} \to \mathbf{2span}(HF)$, sending each manifold, cobordism, or extended cobordism to its space of maps  to $\Bc$. This preserves all compositions, identities, and unitors, up to natural homeomorphisms.  Finally, the once-extended Quinn TQFT,  $$\tFQ{\Bc}\colon \tcob{n} \to \vProfGrphf,$$
is defined from the composite,
$$\smash{\tcob{n} \ra{\Bc^{(-)}} \mathbf{2span}(HF) \ra{\boldsymbol{\mathcal{H}}} \vProfGrphf.}$$
This is treated in Subsection \ref{sec:once extended}. We  check later, in Subsection \ref{Quinn_is_sym_mon}, that indeed this bifunctor can naturally  be given the structure of a  symmetric monoidal bifunctor.
\subsection{The finitary once-extended Quinn TQFT} If $\Sigma$ is a ${n}$-dimensional closed smooth manifold, then the groupoid that $\tFQ{\Bc}$ associates to $\Sigma$  is $\tFQ{\Bc}(\Sigma)=\pi_1(\Bc^\Sigma,\Bc^\Sigma)$. This is a homotopy finite groupoid, however its set of objects is, in general, uncountable. In Subsection \ref{sec:fin_ext}, we will explain how the  size of the image groupoids under  $\tFQ{\Bc}$ can be reduced by considering a closely related bifunctor, $$\tFQd{\Bc}\colon \tdcob{n}{\Bc} \to \vProfGrpfin,$$
that we call the \emph{finitary once-extended Quinn TQFT}.

Here, the objects of the bicategory, denoted $\tdcob{n}{\Bc} $, are now $\Bc$-decorated $n$-manifolds, $(\Sigma,\fd_\Sigma)$. These are, by definition, closed (and smooth) ${n}$-manifolds, $\Sigma$, equipped with a $\Bc$-decoration,  $\fd_\Sigma$, that is, a finite subset, $\fd_\Sigma$, of $\Bc^\Sigma$, containing at least one function in each path-component of $\Bc^\Sigma$. The rest of the bicategory structure of $\tdcob{n}{\Bc}$ is induced from that of $\tcob{n}$, in the obvious way. In particular 1-morphisms, $(\Sigma,\fd_\Sigma) \to (\Sigma',\fd_{\Sigma'})$,  are given  by $(n+1)$-cobordisms, $\Sigma\to \Sigma'$, with no further structure. Similarly  2-morphisms, are given by $(n+2)$-extended cobordisms $(\Sigma \to \Sigma')\To (\Sigma \to \Sigma')$.

On objects, the finitary once-extended Quinn TQFT gives
$$\tFQd{\Bc}(\Sigma,\fd_\Sigma)=\pi_1(\Bc^\Sigma,\fd_\Sigma),$$
and on 1-morphisms and 2-morphisms, we make use of the obvious restrictions of the profunctors and natural transformations given by $\tFQ{\Bc}$.

The groupoids that $\tFQd{\Bc}$ associates to a $\Bc$-decorated manifold, $(\Sigma,\fd_\Sigma)$, explicitly depend on the $\Bc$-decoration, $\fd_\Sigma$, of $\Sigma$. However this dependence is only up to a canonically defined invertible profunctor, which is functorial (up to natural isomorphism) with respect to further changes in the $\Bc$-decoration, and also natural with respect to the profunctors associated to cobordisms.

\subsection{The Morita-valued once-extended Quinn TQFT}
In Subsection \ref{sec:Mor_ext}, we  change the target bicategory of our categorification of Quinn's finite total homotopy TQFT from $\vProfGrpfin$ to $\Mor$, the symmetric monoidal bicategory of (finite dimensional) algebras, with $1$, bimodules and bimodule maps.

Our starting point will be the discussion of a naturally defined \emph{linearisation bifunctor}, 
 $\tLin\colon \vProfGrpfin \to \Mor,$
 essentially defined  in \cite{Mitchell72}, as part of a Morita equivalence between a linear category $\mathcal{C}$ and the algebra $[\mathcal{C}]$ that is associated to it.
 On objects, $\tLin$ sends a groupoid $G$ to its groupoid algebra, \cite{Willerton,loopy}, here denoted $\tLin(G)$. At the level of 1-morphisms, a $\Vect$-profunctor, $\Hp\colon G\bto G'$, of groupoids is then easily converted into a $\big(\tLin(G),\tLin(G')\big)$-bimodule, whose underlying vector space is $$\bigoplus_{x \in ob(G),\,\, y \in ob(G')} \Hp(x,y).$$ Likewise, natural transformations of profunctors naturally linearise to bimodule maps. This construction is closely related to that of Morton in \cite{morton:cohomological:2015}.
 
  These simple observations allow us to define yet one more version of the once-extended Quinn TQFT, the \emph{Morita valued once-extended Quinn TQFT}, $$\tFQmor{\Bc}\colon \tdcob{n}{\Bc} \to \Mor, $$
by considering the following composite of bifunctors:
$$ \smash{\tdcob{n}{\Bc} \ra{\tFQd{\Bc}} \vProfGrpfin \ra{\tLin}\Mor.} $$
 This is done in \S\ref{sec:Mor_valued-eTQFT}. 
 
The algebra, $\tFQd{\Bc} (\Sigma,\fd_\Sigma)$, that is associated to a $\Bc$-decorated ${n}$-manifold depends on the decoration, $\fd_\Sigma$, of $\Sigma$.  However,  this dependence is up to a canonically defined Morita equivalence, which is functorial with respect to further changes in the decoration, and natural with respect to the bimodules associated to cobordisms, a result rooted in \cite[Subsection 10.3]{Bullivant-thesis}.

 \subsection{Explicit calculations for classifying spaces of crossed complexes}
 Quinn's finite total homotopy TQFT, $\FQp{\Bc}{s}$, and its `finitary' once-extended versions, $\tFQd{\Bc}$ and  $\tFQmor{\Bc}$,
 can, in theory, be combinatorially calculated by passing to one of the existing combinatorial models for homotopy theory, for instance, simplicial sets, or simplicial groups,  \cite{Curtis,MaySimplicial}.
  However, notice that the calculation of homotopy contents of function spaces may in theory still present significant challenges since it requires computing all non-trivial homotopy groups. Explicit formulae in this general setting will be deferred to a future paper.

 In the last part of this paper, Chapter \ref{Quinn Calc},  we will work within  a  `truncation' of homotopy theory, obtained by passing to the category of strict infinity-groupoids, or $\omega$-groupoids in the nomenclature of \cite{brown_higgins_sivera}, which we will consider in their equivalent form as crossed complexes, following Brown and Higgins, \cite{brown_higgins_cubes}. (See also the more recent monograph, \cite{brown_higgins_sivera}, by Brown, Higgins and Sivera.)

We will give explicit formula for Quinn's finite total homotopy TQFT,  $\FQp{\Bc}{s}$,  and the two finitary versions of the once-extended Quinn TQFT, for cases in which $\Bc$ is the classifying space, $B_\A$, of a homotopy finite crossed complex, $\A$.  In this case, the computation of Quinn's finite homotopy TQFT, and its extended version, is quite simple since it does not require the explicit calculation of homotopy groups of function spaces. This is due to an `alternating product' formula for the homotopy content of finite crossed complexes, which is an analogue of the formula for the Euler characteristic of a finite complex in terms of an alternating sum of cardinalities of sets of simplices. This is treated in Subsection \ref{sec:homfinXcomp}.

We note, however, that the spaces of the form $B_\A$, where $\A$ is a homotopy finite crossed complex, do not include all possible homotopy classes of   homotopy finite spaces, but include, for instance, those that are 2-types (i.e. whose homotopy groups, $\pi_i$, vanish for $i \ge 3$).
 
 In order to work towards the explicit formulae for $\FQp{\Bc}{s}$, and its extended versions, in Section \ref{Sec:crossed_complexes}, we review the homotopy theory of crossed complexes, closely following work of Brown, Higgins, Sivera, \cite{brown_higgins_sivera,brown_higgins_classifying}, and Tonks, \cite{tonks:JPAA:2003}. Our main new results are in \S \ref{sec:fib-Xcompvs-FibMSpaces}, and,  given a subsimplicial set, $Y$, of a simplicial set, $X$, and a crossed complex, $\A$, they give a crossed complex model for the fibres of the induced fibration, $\TOP(|X|, B_\A) \to \TOP(|Y|,B_\A)$, obtained by restricting a function,  $f\colon |X|\to B_\A$, to $|Y|$. This has direct application to giving explicit formulae for Quinn's finite total homotopy TQFT, and its extended versions.

 In Section \ref{sec:TQFTS_xcomp}, we  finally give explicit formulae for $\FQp{B_\A}{s}$, $\tFQd{B_\A}$, and $\tFQmor{B_\A}$, where $\A$ is a homotopy finite crossed complex. The formulae are mainly given in terms of what we call \emph{simplicial stratifications}, $\zeta_\Sigma\colon |X_\Sigma| \to \Sigma$, of manifolds, $\Sigma$, and analogously  for cobordisms between manifolds, and extended cobordisms. Here $X_\Sigma$ is a simplicial set and $\zeta_\Sigma$ is a homeomorphism. Simplicial stratifications  are  more general than triangulations of manifolds, and typically allow us to decompose a manifold utilising a smaller number of simplices. We will also show computations for CW-decompositions of manifolds, in order to simplify formulae even further.

 It will not be necessary to prove that the formulae given do not depend on the chosen simplicial stratifications, since they are instead proved to coincide with quantities that are, by construction, topological invariants, except when it comes to  what the once-extended TQFTs assign to ${n}$-manifolds, where the dependence on a simplicial stratification, only up to naturally defined invertible profunctors, or bimodules, is naturally a feature of the construction.
 
 On that token,
 we will address, in \S  \ref{eTQFT-colA} and \S \ref{sec:mor-x-complexes}, yet two more versions of the once-extended Quinn TQFT, derived from a crossed complex $\A$, denoted,
\begin{align*}
 &\tFQtr{\A} \colon {\trcob{{n}}} \to \vProfGrpfin
  &&& \textrm{ and }  &&&&
 \tFQmortr{\A} \colon {\trcob{{n}}}  \to \Mor.\end{align*}
Here the bicategory, ${\trcob{{n}}}$, has objects pairs, $(\Sigma,\zeta_{\Sigma})$, with $\Sigma$  a closed ${n}$-manifold with a simplicial stratification, $\zeta_{\Sigma}$. The 1- and 2-morphisms, of ${\trcob{{n}}}$, are cobordisms, and  extended cobordisms, without any choice of simplicial stratification.

This will, in turn, give rise to the construction of (albeit non canonical) once-extended TQFTs,  obtained by picking a simplicial stratification of each path-connected closed ${n}$-manifold, denoted
 \begin{align*}
 &\tFQch{\A} \colon {\tcob{{n}}} \to \vProfGrpfin &&&\textrm{ and }
 &&&&\tFQmorch{\A} &\colon {\tcob{{n}}} \to \Mor.
\end{align*}
This latter construction uses the Axiom of Choice for classes. However the full force of the Axiom of Choice is not required when the domain bicategory of a once-extended TQFT is restricted to a `finitary' sub-bicategory of $\tcob{n}$,  for instance when considering finite presentations of the  symmetric monoidal bicategory $\tcob{0}$, as done in \cite{Schommer-Pries}, or of $\tcob{1}$, as done for example in \cite{Bartlett_etal,Bartlett_Goosen}.
 
One useful  general theorem proved  in this paper is the  following (see Theorem \ref{thm.abs}, in Section \ref{sec:TQFTS_xcomp}).

 \noindent \textbf{Theorem.} \emph{Let $\A$ be a finite crossed complex, and ${n}$  a non-negative integer. We have once-extended TQFTs,
 \begin{align*}
 &\tFQch{\A} \colon {\tcob{{n}}} \to \vProfGrpfin &&& \textrm{ and }
 &&&& \tFQmorch{\A} \colon {\tcob{{n}}} \to  \Mor.
\end{align*}
They can be `normalised' so that, if $\{\Sigma_k\}_{k \in \mathcal{K}}$ is any chosen set of path-connected closed smooth ${n}$-manifolds, and we have selected simplicial stratifications of each manifold $\Sigma_k$, namely $\zeta_{\Sigma_k}\colon |X_{\Sigma_k}| \to \Sigma_k$, then, for each $k$,
\begin{itemize}[leftmargin=1cm]
 \item
  $\tFQch{\A}(\Sigma_k)$ is the groupoid whose objects  are the crossed complex maps from the fundamental crossed complex, $\Pi(X_{\Sigma_k})$, of the simplicial set, $X_{\Sigma_k}$, to $\A$, and whose morphisms are crossed complex homotopies (considered up to 2-fold homotopy) between such  crossed complex maps $\Pi(X_{\Sigma_k})\to \A$,
  \item[]\hspace*{-12mm}and
  \item  $\tFQmorch{\A}(\Sigma_k)$ is the groupoid algebra of $\tFQch{\A}(\Sigma_k)$.
\end{itemize}}

 This is quite a general result, of which we will give some representative examples in Subsection    \ref{sec:exp_formulae_Quinn_groups}.  The category of crossed complexes includes that of groupoids and of strict 2-groups, as full subcategories. Taking $\A$ to be a finite group, $G$, or,  more generally, a finite groupoid, the theorem above gives a homotopy-theoretical interpretation, and a proof of existence, of the $(0,1,2)$-extended TQFTs derived, as in  \cite[\S 3.8]{Schommer-Pries}, from the fact that the groupoid  algebra of a finite groupoid is a ``separable symmetric Frobenius algebra'', see \cite[Example 5.1.]{LaudaPfeiffer}.  Passing to the $(1,2,3)$-extended TQFTs context, and considering a simplicial stratification of $S^1$, with single $0$- and $1$-simplices, and with  $\A$ again derived from a finite group,  $G$, $\tFQmorch{\A}$ associates the quantum double of the group algebra of $G$ to $S^1$. This gives a new proof of, and a homotopy-theoretical interpretation for, the fact that there exists a Morita-valued $(1,2,3)$-extended TQFT, sending $S^1$ to the quantum double of the group algebra of $G$, which is
due to Morton, \cite{morton:cohomological:2015}; see also \cite{Bartlett_etal,Qiu_Wang,MNS}.

We note that the overall construction is considerably more general, and it works in all dimensions, and for all finite crossed complexes. In particular, we also develop, at the end of the paper, the case when $\A=\Gc$, a crossed module of finite groups, which is of relevance for higher gauge theory, \cite{Baez_Huerta:2011,Baez-Schreiber:2004,Martins_Picken:2011}. Concretely, we write down, in Subsection \ref{sec:exp_formulae_Quinn_groups}, some explicit formulae for the $(1,2,3)$- and $(2,3,4)$-extended TQFTs derived from $\Gc$. Passing to the language of discrete higher gauge theory, as treated in \cite{Companion,Morton_Picken_2groupactions_and_moduli}, the  algebras that   $\tFQmorch{\Gc} $ associates to   $S^1$  and to the torus coincide with the `tube algebras' considered in \cite[Sections 10 and 13]{Bullivant-thesis}, \cite{Bullivant_Tube} and \cite[Section 3]{Bullivant_Excitations}, in the context of models for excitations of topological phases, derived from higher gauge theory. These algebras were one of the initial motivations for the  work in this paper.

A general result, directly following from the theorem above, is that if $\Sigma$ is an $n$-manifold, with a simplicial stratification, then there exists an $(n,n+1,n+2)$-extended TQFT that  sends $\Sigma$ to the groupoid of discrete $\Gc$-connections in $\Sigma$ and gauge transformation (considered up to 2-gauge transformation) between them \cite{Companion}.

We expect that, if $G$ is a finite simplicial group -- so that $G$ can represent any finite homotopy type by Ellis' theorem, \cite{Ellis}, -- then there will similarly exist once-extended TQFTs,
\begin{align*}
 &\tFQtr{G} \colon {\trcob{{n}}} \to \vProfGrpfin &&&\textrm{ and } &&&&\tFQmortr{G} \colon {\trcob{{n}}} \to \Mor,
 \end{align*}
sending $(\Sigma, \zeta_{\Sigma}\colon |X_\Sigma| \to \Sigma)$, to the groupoid of simplicial maps from $X_\Sigma$ to $\overline{W}(G)$, the simplicial classifying space of $G$, \cite[Definition (3.20)]{Curtis}, and homotopy classes of maps between them, up to 2-fold homotopy. In that case, $\tFQtr{G}(\Sigma, \zeta_\Sigma)$ will be the fundamental groupoid of the simplicial function space, $\overline{W}(G)^{X_{\Sigma}}$. We hope to  address this in a future publication. In particular, we expect that the recent construction, in \cite{TQFTs-2xmodules}, of topological invariants of 4-manifolds derived from 3-groups (2-crossed modules \cite{Gohla_Martins}), i.e. simplicial groups with Moore complex of length three \cite{Conduche}, is a particular case of the Quinn finite total homotopy TQFT, using a 3-type, $\Bc$, represented by a 2-crossed module of finite groups \cite{Martins-2CW,Baues_4D}, and therefore it can be categorified to a once-extended TQFT.

 In a future publication, we also expect to address the construction of homotopy quantum field theories, including extended ones, derived from crossed complexes. This should be closely related to the construction in \cite{Sozer-Virelizier:HQFT:ArXiv:2022}.
We also hope to address whether Quinn's finite total homotopy TQFT can be further categorified, and also explore its twisting by  cohomology classes. This should lead to a categorification of fully fledged Dijkgraaf-Witten TQFT, \cite{DW}, and its 2-group version, \cite{Martins_Porter}, taking as input more general homotopy finite spaces, enhanced with cohomology classes.

More immediate projects are the construction of the modular tensor categories giving the $(1,2,3)$-extended TQFTs \cite{Bartlett_etal}, derived from homotopy finite spaces $\Bc$, and the description of the representations of motion groups \cite{Qiu_Wang} derived from $\tFQ{\Bc}$.

\medskip
\noindent \textbf{Acknowledgement:}  JFM would like to express his gratitude to: Alex Bullivant, for discussions that framed many of the ideas leading to the construction of the once-extended version of Quinn's finite total homotopy TQFT; to Fiona Torzewska, for discussions on cobordism categories and on cofibrant cospans; to Ben Horton, for discussions on bicategories of extended cobordisms; to Jack Rom\"{o} for discussions on modular categories derived from homotopy finite spaces; and also to Paul Martin, Vincent Koppen, and Eric Rowell for  discussions. We both thank Ronnie Brown for discussions on classifying spaces of crossed complexes, Nicola Gambino for pointing out some important references, and Catherine Meusburger for useful suggestions.

\smallskip

\noindent \textbf{Funding statement:}
The  initial research for this paper was financed by the Leverhulme Trust research project grant RPG-2018-029:
\emph{Emergent Physics from Lattice Models of Higher Gauge Theory}, and then by the project PTDC/MAT-PUR/31089/2017: \emph{Higher Structures and Applications}, of FCT (Portugal). The (years long) final stages of the writing of this paper were funded by EPSRC, via the Programme Grant EP/W007509/1: \emph{Combinatorial Representation Theory: Discovering the Interfaces of Algebra with Geometry and Topology.}

\smallskip

\noindent \textbf{Data Access Statement:} No data was created while producing this publication.

\smallskip

\noindent \textbf{Authorship:} The two authors, Jo\~{a}o Faria Martins and Timothy Porter, contributed equally to this publication. The order of authors is alphabetical.

\smallskip

\noindent \textbf{Open access statement:} For the purpose of open access, the authors have applied a Creative Commons Attribution (CC BY) licence to any Author Accepted Manuscript version arising from this submission.


\chapter{Review of Quinn's Finite Total Homotopy TQFT}\label{Quinn-rend}

\section{The homotopy  underpinning of Quinn's finite total homotopy TQFT}\label{Sec:prelim Quinn}
\sectionmark{The homotopy underpinning of Quinn's  TQFT}

\subsection{Preliminaries: some general conventions and notation}
We will review some of the background theory in the various areas that will feed into this paper. Many readers will not need to read these short sections and need only refer to them when the ideas  and results, mentioned here, are needed in later sections.

\subsubsection{General notation and conventions}
\begin{itemize}[leftmargin=0.6cm]
\item  \label{matrix elements} Let $V$ and $W$ be vector spaces 
with given bases, $X$ and $Y$, respectively. We will  specify a linear map, $f\colon V \to W$, by giving its \textit{matrix elements}, denoted $\langle x \mid f \mid y \rangle \in \kappa$, for $x \in X$ and $y\in Y$. Hence $f(x)=\sum_{y \in Y}  \langle x \mid f \mid y \rangle y$, if $x \in X$.

\item If $X$ is a finite set, then its cardinality is denoted $|X|$.

\item The category of sets is here denoted $\Sets$. The category of $\kappa$-vector spaces is denoted $\Vect$, or $\Vect_\kappa$. The category of topological spaces is  denoted  $\mathbf{Top}$.

\item If $\Cc$ is a monoidal category, the tensor product functor is denoted by $\Otimes_\Cc\colon\Cc\times \Cc \to \Cc$. In the cases where the tensor product arises from a coproduct or a product in $\Cc$, we will also use the notations (respectively),  $$ \Times_\Cc\colon \Cc\times \Cc\to \Cc \textrm{ \quad and \quad }  \Sqcup_\Cc\colon \Cc\times \Cc\to \Cc.$$
\item Given a category $\Cc$ with products, and morphisms $f\colon X \to Y$, $g\colon X \to Z$, in $\Cc$, the map derived from the universal property of the product is denoted $\langle f,g \rangle\colon X \to Y\times Z$. Similarly, if $\Cc$ has coproducts, given maps $f\colon Y \to X$ and $g\colon Z \to X$, the map derived from the universal property of the coproduct is denoted $\langle  f,g \rangle\colon Y\sqcup Z \to X$.
\item Throughout the paper we will implicitly choose a  particular realisation of all limits, colimits, coends, etc, appearing; so, for instance, if $X$ and $Y$ are sets, $X\sqcup Y=X\times\{0\} \cup Y \times \{1\}$.
\end{itemize}
\subsubsection{Conventions and notation for groupoids}\label{sec:convention_groupoids}
The category of groupoids and functors between them is denoted $\Grp$. In this paper, groupoids $G$  will be denoted $G=(s,t\colon G^1 \to G^0)$, where $G^1$ and $G^0$ are, respectively, the set of morphisms and the set of objects of $G$. Morphisms of $G$ are frequently denoted as $\big(s(g) \ra{g} t(g)\big)$, or as $(g\colon s(g) \to t(g))$,  so $s(g)$ is the source of $g$ and $t(g)$ its target. The identity of $x \in G^0$ is denoted $1_x=(1_x\colon x \to x)$, or $1_x^G$. Our convention for notation for  composition in this context is that the composition of $(g\colon x \to y)$ and $(h\colon y \to z)$ is $\big((gh)\colon x\ \to z\big)$. The set of arrows from $x$ to $y$ is denoted $G(x,y)$, or $\hom_G(x,y)$. The vertex group at $x \in G$ is $G(x):=\hom_G(x,x)$.\label{G(x,y)} We always identify a  groupoid having just a single object with its corresponding group of morphisms.

A \emph{totally disconnected groupoid} is a groupoid for which the source and target maps coincide. Totally disconnected groupoids are frequently denoted in the style $A=(\beta\colon A^1  \to A^0)$, where $\beta:=s=t$. If $X$ is a set and for each $x\in X$ we have a group, $E_x$, we can build a totally disconnected groupoid, $(\beta\colon \bigsqcup_{x \in X} E_x \to X)$, where $\beta$ is the obvious map, identifying each component of the disjoint union, and where given $x\in X$, the composition in $\hom(x,x)\cong E_x$ is the product in $E_x$.

A groupoid, $G$,  is said to be \emph{discrete} if it has no non-identity arrows. In this case, it is more or less indistinguishable from its set, $G^0$, of objects. 
We often identify a set, $X$, with the corresponding discrete groupoid having $X$ as its set of objects, and, of course, just the identity arrows as the arrows. This gives an inclusion of the category of sets into that of groupoids.  This inclusion functor  has a left adjoint, sending a  groupoid, $G$, to the set of connected components, $\pi_0(G)$.  For basic information on the theory of groupoids, see \cite{Brown}.

We will often think of groupoids as modelling very simple homotopy types (1-types). We will also recall the notion of homotopy finite space; see Subsection \ref{HF space}.  Combining the two notions, we will have a notion of \emph{homotopy finite}
groupoid.
This is just one of several related finiteness conditions on groupoids used here, namely:
\begin{itemize}[leftmargin=0.6cm]
\item A groupoid, $G$, will be said to be \emph{finite} if both $G^0$ and $G^1$ are finite sets,
\item $G$ is called \emph{locally finite} if each `hom-set' $G(x,y)$ is a finite set,
\item $G$ will be called \emph{homotopy finite} (or briefly to be a \emph{HF-groupoid}) if it has finitely many connected components and each vertex group, $G(x)$, is  finite.
\end{itemize}

\subsubsection{Conventions for compactly generated topological spaces} \label{sec:conventions_top}  We will require a certain background of concepts and notation when handling topological spaces, not all of which is considered in many sources on topology.
\begin{enumerate}[leftmargin=0.6cm]

\item Recall that a space $X$ is called \emph{weak Hausdorff}, see   \cite{May,Strickland} and \cite[\S 7.9]{Dieck}, if given any continuous map, $f\colon K \to X$, where $K$ is compact Hausdorff, then $f(K)$ is closed in $X$.

\item A space $X$ is called \emph{compactly generated} if a set $F\subseteq X$ is closed if, and only if, $u^{-1}(F)$ is closed, for any continuous map $u\colon K \to X$, where $K$ ranges over the class of all compact Hausdorff spaces; see  \cite{Strickland}, \cite[\S 2.4]{Hovey}, \cite[\S 7.9]{Dieck}  or \cite[page 242]{FP}. (Note that compactly generated spaces are called $k$-spaces in \cite{May,FP,Hovey}.) We will denote by $\CG$ the full subcategory of the category $\mathbf{Top}$ with objects the compactly generated  spaces.
\item We have a $k$-ification {functor}, denoted $k\colon \mathbf{Top} \to \CG$. Definitions are in \cite[page 242]{FP} and \cite[\S 7.9]{Dieck}. It is a right adjoint to the inclusion functor $\CG \to \mathbf{Top}$; see \cite[page 243]{FP}. If $X$ is a space,
then the map $k(X) \to X$ given by the identity function,  which we will sometimes denote $\epsilon_X\colon k(X) \to X$, is continuous. This gives the counit of the adjunction. If $f\colon X \to Y$ is a continuous map between topological spaces, then $k(f)\colon k(X) \to k(Y)$ is $f$ itself, at the level of maps between sets.

\item If $K$ is a compact Hausdorff space, then a set map, $f\colon K \to k(X)$, is continuous if, and only if, the same map $f\colon K \to X$ is continuous. (The same holds if $K$ is compactly generated.) In particular, since all disks, $D^n$, are compact Hausdorff, the map $\epsilon_X \colon k(X) \to X$ is a weak homotopy equivalence.
\item By the discussion above, if $X$ is weak Hausdorff, then so is $k(X)$.
\item As in \cite{May,FP,Strom}, we will work in the category $\CGWH$, the full subcategory of $\mathbf{Top}$ with objects the compactly generated and weak Hausdorff  topological spaces, (which we will refer to as \emph{CGWH  spaces}). These include all compact Hausdorff spaces and all metric spaces. Recall that $\CGWH$ has all small limits and colimits,  \cite{Lewis,Strickland}. Note further, see  \cite[Proposition 2.4.22]{Hovey}, the limits in $\CGWH$ are computed by computing the limits in $\mathbf{Top}$, and then applying the $k$-ification functor, so, for example, given a pair of CGWH spaces, their product is $X\times Y=k(X\times_0 Y),$ where $\times_0$ is the product in $\mathbf{Top}$. The colimits in $\CGWH$ are computed as those in $\mathbf{Top}$.

\item Most importantly, we recall moreover, that $\CGWH$ is a cartesian closed  category, \cite{Lewis,Strickland}. Given CGWH spaces, $X$ and $Y$, the  space of maps from $X$ to $Y$, will {be  denoted both by $Y^X$ and by $\TOP(X,Y)$. If $X$ is compact Hausdorff, the topology on $Y^X$ is the $k$-ification of the compact-open topology on the set of maps from $X$ to $Y$.} \label{mapping space}

\item As in \cite{Strickland,Lewis}, a subset, $A$, of a CGWH space, $X$, will be always be given  the $k$-ification of the topology induced by $X$, called the \emph{CGWH induced topology}.  Note that
\begin{itemize}
\item if $F$ is closed in $X$, then $F$ with the induced topology from $X$ is already CGWH, so $k$ will not modify the topology, hence the CGWH induced topology on $F$ is the usual induced topology  on $F$ as a subspace of $X$;
\item\label{coherence_k-induced}  if $A\subseteq B \subseteq X$, then the $k$-ification of the topology that $X$ induces on $A$ coincides with the $k$-ification of the topology that $B$, with the $k$-ification of the induced topology from $X$, induces on $A$;
\item the inclusion, $A \to X$, is continuous,
\item[]\hspace*{-8mm}and
\item if $A\subseteq X$, and $f\colon Y \to A$, with both $X$ and $Y$ being CGWH spaces, then $f$ is continuous (where $A$ has the $k$-ification of the induced topology) if, and only if, $f$ is continuous when considered as a map from $Y $ to $X$.
\end{itemize}
\item\label{itemref:X^I} Given a CGWH space, $X$, then $X^I$ denotes the space of maps from $I=[0,1]$ to $X$, with the $k$-ification of the compact-open topology. We have continuous maps, which we will often  denote, $s:=s_X, t:=t_X\colon X^I \to X$ with $s(\gamma)=\gamma(0)$ and $t(\gamma)=\gamma(1)$. The notation emphasises that these pick out the \emph{source} and \emph{target} of a path.

\item Given a CGWH space $X$, and $j \in \{0,1\}$, also define the inclusions,  $\iota_j^X:=\iota_j\colon \ X \to X \times I$, by $\iota_j(x)=(x,j)$. These are continuous. (We may occasionally use simplified notation for these end inclusions.)
 \item  Given a CGWH space $X$, we can define $\pi_0(X)$  as the underlying set of the coequaliser, in $\CGWH$, of the maps $s,t\colon X^I\to  X $.
\item  \label{Def_PC} Given a CGWH space, $X$, and an element, $x \in X$,  the path-component that $x$ belongs to will be denoted by $\PC_x(X)$. Each path component of $X$ is given the $k$-ification of the topology induced by $X$.

\item The category with objects the CGHW spaces and morphisms from $X$ to $Y$, the maps, $X\to Y$, considered up to homotopy, will be denoted $\CGWH/\simeq$.

\item\label{Def_hatpi}
We will consider a functor, $\hpi_0\colon \CGWH/\simeq  \,\to \Sets$. This sends a CGWH space, $X$, to the set, $\hpi_0(X)$,  of ($k$-ified) path components in $X$, in other words to the set of  $\PC_x(X)$, for $x\in X$.
Let $X$ and $Y$ be CGWH spaces. Given a homotopy class, $[f]\colon X \to Y$, of maps from $X$ to $Y$, we put $$\hpi_0(f)\big(\PC_x(X)\big):=\PC_{f(x)}(Y).$$
While there is an obvious one-to-one correspondence, $\pi_0(X) \leftrightarrow\hpi_0(X)$,
throughout the paper, it will be useful to distinguish between $\pi_0(X)$ and $\hpi_0(X)$. We will write ${\hpi_0}(X)=\{\PC_x(X)\mid x\in X\}$, but note that different $x$ in $X$ may correspond to the same element $\PC_x(X) \in {\hpi_0}(X)$.

\end{enumerate}

\subsubsection{Review of fibrations}\label{HureFib}

Let us recall, for instance from  \cite[Chapter 7]{May} or \cite[Chapter 5]{Strom} the definition.
\begin{Definition}[Fibration]Let $E$ and $B$ be CGWH spaces. We say that $p\colon E \to B$ is a \emph{Hurewicz fibration} (abbr. \emph{fibration}) if the following  \emph{homotopy lifting property holds:}  given any CGWH space $X$, any homotopy, $H\colon X \times I \to B$, and any map,
$f\colon X \to E$, making the diagram with solid arrows, below, commutative, then there exists a map, $H'\colon X \times I \to E$, making the full diagram commutative.
$$\xymatrix@R=15pt{
  X \ar[d]_{\iota^X_0} \ar[rr]^{{f}}&&E\ar[d]^p \\
X \times I \ar@{-->}[rru]^{H'}\ar[rr]_{H} &&B
.}
$$
This  $H'\colon X \times I \to E$ is called a \emph{lifting of $H$ starting at ${f}$}.
\end{Definition}
 Differently from the conventions in \cite[Chapter 7]{May}, we do not impose that fibrations are surjective\footnote{This surjectivity condition was dropped in the subsequent \cite{May-Ponto}; see footnote on page 25.},  hence, given a space $B$,  $\emptyset \to B$ is  a fibration.

 The following is well known, and has an easy proof; see e.g. (\cite[Chapter 6]{May}).
\begin{Lemma}The composite of fibrations is a fibration. If $X$ and $Y$ are CGWH-spaces, then the projections $p_X\colon X \times Y\to X$ and  $p_Y\colon X \times Y\to Y$  are both fibrations.
Finally, \emph{pullbacks of fibrations are fibrations.} This means that if $p\colon E \to B$ is a fibration, and $f\colon X \to B$ is any map (of CGWH spaces), then the map, $q\colon X \times_B E \to X$, appearing in the pullback diagram below is a fibration,
 $$\xymatrix@R=15pt{&X\times_B E\ar@{}[dr]|<<<{\pullback} \ar[d]_q\ar[r] &E\ar[d]^p\\& X \ar[r]_f & B.}$$
\end{Lemma}

If $p\colon E \to X$ is a fibration, then, given $x \in X$, the \emph{fibre} of $p$ at $x$ is
$E_x:=p^{-1}(x).$
Following our conventions in Subsection \ref{sec:conventions_top},  the fibre $E_x$ is given the induced CGWH-topology from $E$.
Since $E_x$ is closed in $E$, it is compactly generated already, with the induced topology, so the $k$-ification step does not modify the topology in $E_x$.
We also note that we have the following pullback diagram in $\CGWH$, $$\xymatrix@R=15pt{& E_x \ar@{}[dr]|<<<{\pullback}\ar[d] \ar[r]^{\inc} & E \ar[d]^p\\ & \{x\} \ar[r]_{\inc} & X ,}$$
where the $\inc$ denote the obvious inclusion maps.
More generally, let $A$ be a subset of $X$, and let $E_A:=p^{-1}(A)$,  and consider the induced map, $p_A\colon E_A \to A$. Give both $A$ and $E_A$ the induced CGWH topology. We have a pullback diagram in $\CGWH$, where $p_A\colon E_A \to A$ is a fibration, namely,
$$\xymatrix@R=15pt{& E_A \ar@{}[dr]|<<<{\pullback}
\ar[d]_{p_A} \ar[r]^{\inc} & E \ar[d]^p\\ & A \ar[r]_{\inc} & X .}$$

As a closely related case, suppose $E \neq \emptyset$ and let  $e \in E$. Put $x=p(e)$. It is easy to see that $p(\PC_e(E))=\PC_x(X)$, and that  the induced map $p_e\colon \PC_e(E) \to \PC_x(X)$ is a surjective fibration; see \cite[Lemma 2.3.1]{Spanier}.

We will make extensive use of the fact that if $p\colon E \to B$ is a fibration, and $x,y \in B$ are in the same path-component, then the fibres, $p^{-1}(x)$ and $p^{-1}(y)$, are homotopy equivalent; see e.g. \cite[Chapter 7]{May}. We will also need that if $E$ is path-connected and $x\in X$, it follows that all  path-components of  $E_x$ are homotopy equivalent, \cite[Proposition 3]{Fadell}.
We will review some of these results in more detail, later, starting with Lemma  \ref{main_techfib}, page \pageref{main_techfib}.

\subsubsection{Cofibrations}
Looking at the dual setting, recall that a map, $f\colon A \to X$, of CGWH spaces is a \emph{cofibration}, \cite[Chapter 6]{May}, or  \cite[\S 5.1]{Strom}, if it satisfies the homotopy extension property:

 \begin{itemize}[leftmargin=0.4cm]
  \item[]
\textit{For any CGWH space, $B$, any map, $g\colon X \to B$, and any homotopy, $H\colon A \times I \to B$, as in the solid arrows of the  diagram,
$$\xymatrix@R=13pt@C=40pt
{A\ar[dd]_{\iota_0^A}  \ar[r]^f & X\ar[d]_{\iota_0^X} \ar[r]^g & B\\
               & X\times I \ar@{-->}[ur]^{H'}\\ A \times I  \ar[ru]^{f \times \id_I} \ar@/_1.7pc/[rruu]_>>>>>>>{H} } $$
            there is a homotopy, $H'\colon X\times I\to B$, making the diagram commute.}
\end{itemize}
The following two well-known results will be used without further comment.
\begin{itemize}[leftmargin=0.6cm]
 \item Let $f\colon A \to X$ be a cofibration and let $B$ be a CGWH space, then the induced map on mapping spaces, $f^*\colon B^X \to B^A$, {sending $\phi\colon X \to B$ to $\phi\circ f \colon A \to B$,}  is a fibration. (For example, see  \cite[Section 7.2]{May}.)
\item If $(X,Y)$ is a CW-complex pair, meaning that $Y$ is a subcomplex of the CW-complex, $X$, then the inclusion map $i\colon Y \to X$ is a cofibration; see, for instance, \cite[Corollary 1.4.7]{FP}.
 \end{itemize}


\subsubsection{Fibre homotopy}\label{sec:fhe}
Let $f,g\colon X \to Y$ be  maps of CHWH spaces. A homotopy, $H\colon X \times I \to Y$,  connecting $f$ to $g$, will frequently be denoted by $f\Ra{H} g$.

Given two fibrations, $p\colon X \to B$ and $q\colon Y \to B$, over the same space, a \emph{fibre map}, or \emph{fibred map}, $f\colon X \to Y$, is a map such that the diagram below commutes,
\begin{equation}\label{diagfibs}
\vcenter{\xymatrix@R=13pt{X\ar[dr]_p \ar[rr]^f &&Y.\ar[dl]^q\\ &B&}}
\end{equation}
Two fibre maps, $f,g\colon X \to Y$, are \emph{fibre homotopic} if there exists a homotopy, $H\colon X \times I\to Y$, called a \emph{fibre} or \emph{fibred} homotopy, connecting $f$ and $g$, and  such that for each $(x,t) \in X \times I$, $p(x)=q(H(x,t))$.  We write $f\Ra[B]{H} g$.

We say that a pair of fibre maps, $f\colon X \to Y$ and $f'\colon Y \to X$, realises a \emph{fibre homotopy equivalence} if we have fibre homotopies, $$f\circ f'\Ra[B]{H} \id_Y \qquad \textrm{ and } \qquad f'\circ f\Ra[B]{H'} \id_X .$$
This means that $H\colon X\times I \to X$ satisfies $p(H(x,t))=p(x)$, for each $(x,t) \in X \times I$, and similarly for $H'$.

The following non-immediate, but well-known, result will be extensively used later.  It is the dual version of Dold's Theorem; see the discussion and references  in \cite[Chapter I, section 6, p. 33]{KHK-TP:1997}.
For a proof, see \cite[Chapter 7.5]{May}, \cite[Theorem 3.4]{Brown_Heath}. Very thorough discussions in the dual case of cofibrations appear in \cite[7.4.2: Addendum]{Brown} and in \cite{KHK-TP:1997}. 
\begin{Lemma}[Dual Dold Theorem]\label{miracle}
Suppose that  in \eqref{diagfibs}, $f\colon X \to Y$ is a homotopy equivalence, then there exists a homotopy inverse, $f'\colon Y \to X$, of $f$, which is a fibre map,  and such that $f$ and $f'$ realise a fibre homotopy equivalence.
\end{Lemma}
Almost by definition, we have:
\begin{Lemma}\label{fibre_hom}
 Suppose that $f\colon X \to Y$, as in \eqref{diagfibs}, is a fibre homotopy equivalence. Given  $b \in B$, the map, $f$, restricts to a homotopy equivalence, $p^{-1}(b) \to q^{-1}(b)$.
\end{Lemma}
We will also need a `relative' version of the result from Lemma \ref{miracle}.
\begin{Definition}[Map of fibrations]\label{map of fibs}
Given two fibrations, $p\colon D\to A$ and $q\colon E\to B$, a \emph{map from $p$ to $q$} is a pair, $(g,f)$, of maps as in the square,
$$\xymatrix@R=18pt{D\ar[r]^g\ar[d]_p&E\ar[d]^q\\
A\ar[r]_f&B,
}$$making that square commute. We write $(g,f)\colon p\to q$.
\end{Definition}
\begin{Definition}\label{he of fibs}Given two fibrations, $p\colon D\to A$ and $q\colon E\to B$, a map, $(g,f)\colon p\to q$, as above, is a \emph{homotopy equivalence  of fibrations} if there are homotopy inverses, $f'$ of $f$ and $g'$ of $g$, such that $p\circ g'=f'\circ q$, and, in addition, there are homotopies, $H\colon g'\circ g\simeq \id_D$, and $K\colon g\circ g'\simeq \id_E$, that cover homotopies, $h\colon f'\circ f\simeq \id_A$ and $k\colon f\circ f'\simeq \id_B$.
\end{Definition}
The relative version of  Lemma \ref{miracle}, to be found in \cite[p. 53]{May}, is then:
\begin{Proposition}\label{relative miracle}If $(g,f)\colon p\to q$ is a map of fibrations, where 
both $g$ and $f$ are homotopy equivalences, then $(g,f)\colon p\to q$ is a homotopy equivalence of fibrations. %
\end{Proposition}
\begin{Corollary}\label{h.e. on fibres}
Let $p\colon D\to A$ and $q\colon E\to B$ be fibrations. Let $(g,f)\colon p\to q$ be a homotopy equivalence  of fibrations, then, for any $a\in A$, the induced map on fibres, $f\colon p^{-1}(a)\to q^{-1}(f(a))$ is a homotopy equivalence.
\end{Corollary}
\subsection{Homotopically finite (HF) spaces}\label{HF space}
This is a key notion for this paper.
\begin{Definition}[Homotopically finite (HF) space]
 A  space, $B$, is called \emph{homotopy finite} (abbr. HF) if $B$ is CGWH  (see \S\ref{sec:conventions_top}) and, moreover, $B$ has only a finite set of path-components, each of which has only a finite set of non-trivial homotopy groups, all of which are finite.
\end{Definition}
Clearly, finite disjoint unions and finite products of HF-spaces are HF.
Each path-component of a HF-space is also HF (after possibly applying the $k$-ification functor in order to make it a CGWH space).

\begin{Remark}In the literature, one finds some alternative terminology used for homotopy finite spaces.
Lurie,  \cite[Appendix E]{lurie:spectral:2018}, calls them `$\pi$-finite spaces' whilst Anel, \cite{anel:elementary:2021}, uses the term `truncated coherent space'.
\end{Remark}
The following is essentially in \cite[Lemma 3.4]{Galvezetal}, albeit stated in the context of $\infty$-groupoids.
\begin{Lemma}\label{2of3}
Let $p\colon E \to B$ be a
fibration. Given $b \in B$, we let $E_b:=p^{-1}(b)$.
\begin{enumerate}[leftmargin=1cm]
 \item  Suppose that $B$ is path-connected, and that $p$ is surjective.  If any two of $B$, $E$ and $E_b$ are HF, then so is the third.
\item Let $B$ be any space. If $B$ and $E$ are HF, then so is $E_b$, whenever $b \in B$.
\item If  $B$ and
each $E_b$ are HF (for each $b \in B$), then so is $E$.
\end{enumerate}

In particular, if $p\colon E \to B$ is a fibration, and $E$ and $B$ each are HF, then each fibre of $p$ is HF.
\end{Lemma}
 \begin{proof} Follows from the homotopy long exact sequence of $p\colon E \to B$; see Equation \eqref{hles} below.
 \end{proof}
The second point of the following result will be crucial for what follows. This is stated in \cite[Lemma 3.13]{Galvezetal} for $\infty$-groupoids.
\begin{Lemma}\label{pull_HF}
 Consider a pullback diagram of spaces, where $p\colon E \to B$ is a fibration, $$\xymatrix@R=18pt{&X\times_B E\ar@{}[dr]|<<<<{\pullback} \ar[d]_q\ar[r] &E\ar[d]^p\\& X \ar[r]_f & B,}
 $$
 then  $q$ is a fibration. If $X$, $E$ and $B$ are HF, then $X\times_B E$ is HF.
\end{Lemma}
\begin{proof}
 That $q$ is a fibration follows from the standard fact  that pullbacks of fibrations are fibrations; see e.g. \cite[\S 6.1 Lemma]{May}.

 Let us prove that $X\times_B E$ is HF. We  use the previous lemma. By assumption, $X$ is $HF$. We only need to prove that the fibres of $q$ are HF. Given $x \in X$, the fibre of $q$ at $x$ is  homeomorphic to $p^{-1}(f(x))$, which is HF since $B$ and $E$ are.\end{proof}
 \subsubsection{{The homotopy content  of a HF-space}} This is a key notion for this paper.
\begin{Definition}[{Homotopy content}]\label{def:hcont}
Let $B$ be a path connected HF-space. The homotopy content  of $B$ is defined as ,
\begin{align*}
\chi^\pi(B)
=
\frac{ \big |\pi_2(B,x) \big |\, \big |\pi_4(B,x) \big |\, \big |\pi_6(B,x) \big |\ldots}{ \big |\pi_1(B,x) \big |\, \big |\pi_3(B,x) \big |\,  \big |\pi_5(B,x) \big |\ldots} \in \Q,
\end{align*}
where $x\in B$ is any point.
In general, if $X$ is a HF space, recalling the notation in item \eqref{Def_hatpi} on page \pageref{Def_hatpi}, define
$$\chi^\pi(X)=\sum_{B \in \hpi_0(X)} \chi^\pi(B).$$
We also put $\chi^\pi(\emptyset)=0$. (Note that for all other HF spaces $F$, we have $\chi^\pi(F)>0$.)
\end{Definition}

We note that what we have called `homotopy content' is called  `homotopy order' in  \cite[Lecture 4]{Quinn},  and `homotopy cardinality' in \cite{BaezDolan} and in \cite[\S 3]{Galvezetal}.

The homotopy  content of a space also appeared in \cite{MartinsCW}, without being given a name,  and, there, was also considered for crossed complexes. We will consider that form separately a bit later on here. The case of $\infty$-groupoids is treated in \cite{Galvezetal}, which proves similar results to the two lemmas below, in that context.

 Note that homotopic HF spaces, and, more generally, weakly homotopic HF spaces, have the same homotopy  content.

\begin{Example} The customary examples are (i) when $X$ is a finite set, thought of as a discrete space, then $\chi^\pi(X)$ is the usual cardinality of $X$, and (ii) when $X$ is the classifying space of a finite groupoid, $G$, then $\chi^\pi(X)= \sum_{[x]\in \pi_0(G)} \frac{1}{|G(x)|}$, which is the \emph{groupoid cardinality of $G$}, in the sense of Baez and Dolan, \cite{BaezDolan}, and \cite{baezetal}.  \end{Example}
\begin{Lemma}\label{lem:sums.prods}
If $B$ and $B'$ are HF-spaces, then
$$
 \chi^\pi(B \sqcup B')=\chi^\pi(B)+\chi^\pi(B')\quad \textrm{ and }\quad \chi^\pi(B \times B')=\chi^\pi(B) \chi^\pi(B').
$$
\end{Lemma}
\begin{proof}
The first equation is straightforward.  For the second note that
$$\hpi_0(B\times B')\cong \{ A \times A' | A \in \hpi_0(B), A'\in \hpi_0(B')\},$$
and that, if $x \in B$ and $x'\in B'$, then $\pi_n\big (B\times B',(x,x')\big)\cong \pi_n(B,x)\times \pi_n (B',x').$
\end{proof}

More generally,
\begin{Lemma}[Quinn, \cite{Quinn}, Baez--Dolan, \cite{BaezDolan}, and Galv\'{e}z-Carillo--Kock--Tonks, \cite{Galvezetal}]\label{main1}
Suppose that $p\colon E \to B$ is a   fibration of HF-spaces and that $B$ is path-connected. Let $b \in B$ be arbitrary, then, recalling  $E_b=p^{-1}(b)$ is the fibre at $b$,
$$\chi^\pi(E)=\chi^\pi(B) \, \chi^\pi(E_b).$$
\end{Lemma}
The proof we give below is as hinted at in the above references, with some crucial technical details added.

\begin{proof}
If $E$ is empty, then so is $E_b$, and in this case there is nothing to prove. If $B$ is empty, then so are $E$ and $E_b$, and there is nothing to prove either.

 We are left with the case  that $E,B \neq \emptyset$. In this case, it follows that $p\colon E \to B$ is surjective, as $B$ is path-connected. More generally, if $E'$ is a path-component of $E$, the restriction $p'\colon E' \to B$ of $p$ is also surjective.

 Suppose, firstly, that $E$ is path-connected. Let $x \in E$ and $b=p(x)$, then, cf. \cite[p. 52]{May} or \cite[p. 376]{Hatcher}, the homotopy long exact sequence of $p\colon E \to B$, at $b$ and $x$  reads
\begin{multline}\label{hles}
\to \pi_i(E_b,x) \stackrel{\iota}{\to} \pi_i(E,x) \stackrel{\d}{\to} \pi_i(B,b) \stackrel{\delta}{\to} \pi_{i-1}(E_b,x) \to \dots\\  \stackrel{\iota}{\to} \pi_1(E,x) \stackrel{\partial}{\to} \pi_1(B,b) \stackrel{\delta_x}{\to} \pi_0(E_b) \ra{\iota} \pi_0(E)=\{0\}.  
\end{multline}
Here, for the last stages of the sequence, the exactness means the following:
\begin{itemize}[leftmargin=0.6cm]
\item we have a left-action, $\triangleright$, of $\pi_1(B,b)$ on $\pi_0(E_b)$ (reviewed in Lemma  \ref{main_techfib_loops}), whose stabiliser subgroup at the path-component, $\PC_x(E_b)$, of $x \in E_b$, is $\d(\pi_1(E,x))$;
\item the map $\delta_x\colon \pi_1(B,b) \to \pi_0(E_x)$, which is defined as $\delta_x(g)= g\triangleright \PC_x(E_b)$, is surjective,
\item $\iota\colon \pi_0(E_b) \to \pi_0(E)$ is surjective, and descends to a  bijection $\pi_0(E_b)/\pi_1(B,b)\cong \pi_0(E)$.
\end{itemize}
Also note that, by the orbit-stabiliser theorem, $|\pi_0(E_b)|=|\pi_1(B,b)|/|\d( \pi_1(E,x))|$.

The  exactness of the sequence \eqref{hles} yields that:
\begin{align*}
 |\pi_i(E,x)|&=|\d(\pi_i(E,x))|\, |\iota(\pi_i(E_b,x))|, \textrm{ if } i \ge 1,\\
| \pi_i(B,b)|&=|\d(\pi_i(E,x))|\, |\delta(\pi_i(B,b))|, \textrm{ if } i \ge 2,\\ 
 |\pi_1(B,b)|&=|\d(\pi_1(E,x))|\, |\pi_0(E_b)|,\\
| \pi_i(E_b,x)|&=|\delta(\pi_{i+1}(B,b))| |\iota(\pi_i(E_b,x))|, \textrm{ if } i \ge 1. 
 \end{align*}
Therefore, noting that $B$ and $E$ are by assumption path-connected,
\begin{align*}
\chi^\pi(B)&= \frac{1}{|\d(\pi_1(E,x))|\, |\pi_0(E_b)|} \prod_{k=2}^{+\infty} \big( \big |\d(\pi_k(E,x)) \big |\, \big |\delta(\pi_k(B,b)) \big | \big)^{( (-1)^{-k})},\\
\chi^\pi(E)&= \prod_{k=1}^{+\infty} \big( \big |\partial(\pi_k(E,x)) \big |\, \big |\iota(\pi_k(E_b,x)) \big | \big)^{ ( (-1)^{-k})},\intertext{and also,}
\chi^\pi(E_b)&= \big |\pi_0(E_b) \big | \prod_{k=1}^{+\infty} \big( \big |i(\pi_k(E_b,x)) \big |\, \big |\delta(\pi_{k+1}(B,b)) \big | \big)^{( (-1)^{-k})}.
\end{align*}
Crucially, in the last equation, we also used the fact that given that $p\colon E \to B$ is a  fibration, and $E$ is  path-connected,    all  path-components of  $E_b=p^{-1}(b)$ are homotopy equivalent,  \cite[Proposition 3]{Fadell}. This is reviewed in Lemma  \ref{main_techfib}.

We thus have
\begin{align*}
  \chi^\pi(E_b) \chi^\pi(B)
&=  \prod_{k=1}^{+\infty} \big( \big |i(\pi_k(E_b,x)) \big | \,  \big |\d(\pi_k(E,x)) \big |\big)^{ ( (-1)^k)}=\chi^\pi(E).
\end{align*}

Suppose now that $E$ may have more that one path-component (but recall that we still take  $B$ to be path-connected). Let $E^1,\dots, E^n$ be the path-components of $E$. Let $p_k\colon E^k \to B$ be the restriction of $p$ to $E^k$, for each $k=1,\dots, n$. Each $p_k$ is itself a fibration, and  is surjective. Let $F_k= p_k^{-1}(b)=E_b\cap E^k$. Note that we have an obvious continuous bijection $\sqcup_{k=1}^n F_k \to E_b$, which is always a weak homotopy equivalence. We therefore have:
\begin{align*}
\chi^\pi(E)&=\sum_{k =1}^n \chi^\pi(E^k)=\sum_{k=1}^n \chi^\pi(F_k )\, \chi^\pi (B)
\\&=  \chi^\pi\big( {\sqcup}_{k=1}^n F_k \big)\, \chi^\pi (B)=\chi^\pi(E_b)\,\chi^\pi(B).
\end{align*}
(Note that we do not necessarily have a homeomorphism ${\sqcup}_{k=1}^n F_k\cong E_b$.)
\end{proof}
We have the following, which is very useful later on.
\begin{Theorem}\label{main2}
Let $p\colon E \to B$ be a  fibration, where $B$ and $E$ are HF. If $b \in B$, and $E_b=p^{-1}(b)$, then
$$\chi^\pi(E)=\sum_{[b] \in \pi_0(B)} \chi^\pi(E_b) \,\chi^{\pi}({\PC_b(B)}).$$\end{Theorem}
\noindent (Here we have chosen a representative of each path-component of $B$, noting that if $b$ and $b'$ are in the same path-component then $E_b$ is homotopic to $E_{b'}$.)

\begin{proof} If $B$ is empty, then so is $E$, so the result follows trivially, so we suppose that $B\neq \emptyset$.
That if $b \in B$, then $E_b$ is HF follows from Lemma \ref{2of3}.
Given $[b] \in \pi_0(B)$,  put $E_{[b]}=p^{-1}({\PC_b(B)})$. The restriction, $p_b \colon E_{[b]} \to {\PC_b(B)}$, of $p\colon E \to B$, is a  fibration. We  have weak homotopy equivalences,
$$\bigsqcup_{[b] \in \pi_0(B)} \PC_b(B) \to B \quad \textrm{ and }\quad \bigsqcup_{[b] \in \pi_0(B)} E_{[b]} \to E,$$
therefore
\begin{align*}
 \chi^\pi(E)&=\chi^\pi\Big (\bigsqcup_{[b] \in \pi_0(B)} E_{[b]} \Big) =\sum_{[b] \in \pi_0(B)}\chi^\pi\big ( E_{[b]} \big)=\sum_{[b] \in \pi_0(B)} \chi^\pi(E_b) \,\chi^{\pi}(\PC_b(B)).
\end{align*}

\end{proof}

\subsection{Fibrant spans of HF spaces and their composition: {the category $\HFb$}}\label{sec-fib-span}

Before we introduce fibrant spans in detail, we should briefly motivate why we are going to use them.  The objects considered in basic TQFTs are manifolds of some type, and the cobordisms between them.  Such a set-up gives a cospan {of CGWH spaces},
$$ \vcenter{\cspn{\Sigma}{i}{\Sigma'}{j}{S}}$$ and we have that the induced map, $\Sigma\sqcup \Sigma'\to S$, is an inclusion, and furthermore a cofibration; see later in Subsection \ref{cobordism categories}, starting on page \pageref{cobordism categories}, for a more detailed discussion. Such \emph{cofibrant cospans} of spaces are studied in detail in \cite{Torzewska,Torzewska-HomCobs}.

To study the state spaces associated to the manifolds, we form the space of maps from such  manifolds to a `classifying space' $\Bc$, which we will take in Section \ref{Classical Quinn} to be homotopy finite, but, in so doing, we convert a cospan, as above, to a span,$$\xymatrix@R=-5pt{ &\Bc^S\ar[dl]_{{i^*}}\ar[dr]^{{j^*}}\\  \Bc^\Sigma & &\Bc^{\Sigma'},  }$$ where $i^*$ and $j^*$ denote the obvious restriction maps,
and we note that the induced map from $\Bc^S$ to $\Bc^\Sigma \times \Bc^{\Sigma'}$ is a fibration. To study that type of situation, we need to understand \emph{fibrant spans} and we will examine them in some generality, not just in this particular function space set-up.

\subsubsection{Fibrant spans  and HF fibrant spans}
Let $B$ and $B'$ be CGWH spaces.
\begin{Definition}[Fibrant span]\label{fib_span}
A  \emph{fibrant span}, \smash{$B \ra{(p,M,p')} B'$},  from $B$ to  $B'$, also denoted  $(p,M,p')\colon B \to B'$ is a diagram in $\CGWH$ of form,
\begin{equation}\label{hff} \vcenter{\xymatrix@R=-5pt{ &M\ar[dl]_p\ar[dr]^{p'}\\   B & &B' , }}
 \end{equation}
where the induced map $\big \langle p,p'\big \rangle\colon M \to B\times B'$ is a   fibration. If the spaces, $B,B'$ and $M$  are all HF, we will say this is a HF fibrant span or a fibrant span of HF spaces.
\end{Definition}
\begin{Remark}\label{Spans and Lambda}
Consider the  Hurewicz / Str\o m model structure on $\CGWH$; see \cite{Str}. Let $\Lambda$ be the category $\{-1 \leftarrow 0 \rightarrow 1\}$. This is an inverse category in the sense used in, for instance, \cite[\S 5.1]{Hovey}. If we give  the injective model structure to $\CGWH^{\Lambda}$, then weak equivalences and cofibrations are given objectwise, see \cite[Theorem 5.1.3]{Hovey}, whilst the fibrant objects are precisely our fibrant spans.
\end{Remark}
\begin{Example}\label{ids_ex} Let $X$ be a  space. The trivial span on $X$ is then $$\xymatrix@R=-5pt{ &X\ar[dl]_{\id_X}\ar[dr]^{\id_X}\\   X & &X.  }$$This is clearly not a fibrant span if $X$ is non-empty, so we want to replace it by a fibrant one. Consider, where $X^I$ is the space of functions from $I=[0,1]$ to $X$,
\begin{equation}\label{fibspan_unit} 
\vcenter{\xymatrix@R=-5pt{ &X^I\ar[dl]_{s_X}\ar[dr]^{t_X}\\   X & &X . }}
 \end{equation}
Here if $\gamma\colon I \to X$, then $s_X(\gamma)=\gamma(0)$ and $t_X(\gamma)=\gamma(1)$.  We can see directly that we have a fibration, $$\langle s_X,t_X \rangle\colon X^I \to  X \times X.$$ This  follows, for instance, from the fact that the inclusion, $\iota\colon \{0,1\} \to I$, is a cofibration, and hence the induced map, $\iota^*\colon X^I \to X^{\{0,1\}}\cong X\times X$, is a fibration.
If $X$ is HF, then so is $X^I$, as it is homotopy equivalent to $X$.\end{Example}
\begin{Lemma}\label{All-fibrations}
 Consider a fibrant span,  $(p,M,p')\colon B \to B'$. Both maps, $p\colon M \to B$ and $p'\colon M \to B'$, are   fibrations, and moreover, given  $b \in B$ and $b' \in B'$,  both of the induced maps, $p^{-1}(b) \to B'$ and ${p'}^{-1}(b') \to B$, are   fibrations.
\end{Lemma}
\begin{proof}
 For the first point, given the fact that both projections, $B\times B' \to B$ and $B\times B' \to B'$, are fibrations, and also that the composite of fibrations is a fibration, it follows that both $p$ and $p'$ are fibrations.
The second point follows from the fact that pullbacks of fibrations are fibrations.
\end{proof}
 \begin{Lemma}\label{all-HF}
 Suppose that the fibrant  span, $(p,M,p')\colon B \to B'$,
 is HF. Let $b\in B$ and $b' \in B'$. The spaces $p^{-1}(b)$ and ${p'}^{-1}(b')$, and also the fibre of $\langle p,p^\prime \rangle\colon M \to B\times B'$, over $(b,b^\prime)$, i.e., $\big \langle p,p^\prime\big \rangle^{-1}(b,b')$, are all HF.
 \end{Lemma}
 \begin{proof}The second statement follows from Lemma \ref{2of3}, since $M$ and $B\times B'$ are both HF. (Recall that the product of two HF spaces is HF.) For the first, by Lemma \ref{All-fibrations}  we have a   fibration $p'\colon p^{-1}(b) \to B'$. The fibres have the form $\big \langle p,p^\prime\big \rangle^{-1}(b,b')$, and hence they are  also fibres for the fibration, $\langle p,p'\rangle\colon M \to B\times B'$, so they must be HF. Since $B$ is HF it hence follows that $p^{-1}(b) \to B'$ is HF,  by Lemma \ref{2of3}.
\end{proof}

 \begin{Lemma}\label{main3}
Let $B$, $B'$ and $B''$ be HF-spaces. Consider HF fibrant spans,   $${(p,M,p')}\colon B \to  B'  \qquad \textrm{ and  } \qquad  (p'',M',p''')\colon B' \to B''.$$ We form the obvious pullback, as in the diamond in the  diagram below,
\begin{equation}\label{compHFfibrations}\vcenter{\xymatrix@R=0pt{ &&& 
M\times_{B'} M'
\ar[dd]^P \ar[dl]_{q}\ar[dr]^{q'}\\ &&M\ar[dl]_p\ar[dr]^{p'}&&M'\ar[dl]_{p''}\ar[dr]^{p'''}\\  & B & &B'  & &B'',  }}
\end{equation}
where $P:=p'\circ q=p''\circ q'$, then the span,
\begin{equation}\label{eq:comp-fib-spans}
\big(B \ra{(p,M,p')  \bullet (p'',M',p''')} B''\big):=\big( B \ra{ (p\circ q,  M\times_{B'} M'
,p'''\circ q' )} B''\big),
\end{equation}
is a fibrant span of HF spaces. 

We also have that $\big \langle p\circ q, P, p'''\circ q'\big \rangle\colon M\times_{B'} M' \to B \times B' \times B''$ is  a   fibration.
\end{Lemma}
\begin{proof}
That $\big \langle p\circ q, P, p'''\circ q'\big \rangle$ is  a fibration is clear from the fact that $\big \langle p,p' \big \rangle$ and $\big \langle p'',p''' \big \rangle$ are fibrations, and from the universal property of pullbacks. It follows that $\big \langle p\circ q,p'''\circ q' \big \rangle$ is also a fibration, for the projection, $B\times B'\times B''\to B \times B''$, is a fibration.
To prove that  $M\times_{B'} M'$ is HF, it suffices (by Lemma \ref{2of3}) to observe that $B \times B' \times B''$ is HF and that the fibres of the fibration, $\big \langle p\circ q, P, p'''\circ q'\big \rangle$, have the form $(\big \langle p,p'\big \rangle^{-1}(b,b')) \times (\big \langle p'',p'''\big \rangle^{-1}(b',b''))$. Each fibre is thus  HF, since both components of the product are.
\end{proof}
 \begin{Definition}[{Composition of HF fibrant spans}]\label{main3def}
 The HF fibrant span in \eqref{eq:comp-fib-spans}
is  called the \emph{composite of $(p,M,p')\colon B \to B'$ and $(p'',M',p''')\colon B' \to B''$.}
\end{Definition}
 
\subsubsection{The category $\HFb$}
We note that this will be  a non-locally small category, as we have a class of maps between objects. This will, however, not cause any difficulties.

The class of objects of $\HFb$ is the class of all HF spaces. Given HF-spaces, $B$ and $B'$, the class of morphisms from $B$ to $ B'$ is given by equivalence classes of HF fibrant spans, $(p,M,p')\colon B\to B'$, 
as we now explain. We will make use of the materials on fibre homotopy equivalence recalled in \S\ref{sec:fhe}.

\begin{Definition}[Equivalent and isomorphic HF fibrant spans]\label{equiv-of-HF-spans} Let $B$ and $B'$ be HF spaces. Two HF fibrant spans,
$ (p,M,p')\colon  B \to B'$ and $(q,N,q')\colon  B\to B'$,
  are said to be equivalent if there exist fibred maps, $\Psi \colon M \to N$ and $\Psi'\colon N \to M$, i.e. maps making the diagrams below commute,
 \begin{equation}\label{equiv}\vcenter{\xymatrix@R=5pt{ &&M\ar[dd]^\Psi\ar[dl]_p\ar[dr]^{p'}\\  & B & &B' \\ &&N\ar[ul]^q\ar[ur]_{q'} }}
 \qquad \textrm{ and }  \vcenter{\xymatrix@R=5pt{ &&M\ar[dl]_p\ar[dr]^{p'}\\  & B & &B', \\ &&N\ar[uu]_{\Psi'}\ar[ul]^q\ar[ur]_{q'} }}
 \end{equation}
realising a fibre homotopy equivalence, with respect to the fibrations, $\langle p,p'\rangle \colon M \to B\times B'$ and $\langle q,q'\rangle \colon N \to B \times B'$. This  means that homotopies, $H\colon M \times I \to M$ and $H'\colon N \times I \to N$, exist such that:
\begin{enumerate}[leftmargin=1cm]
 \item\label{1} $H(m,1)=\Psi'(\Psi(m))$ and $H(m,0)=m$ for each $m \in M$;
 \item $p(H(m,t))=p(m)$ and $p'(H(m,t))=p'(m)$, for each $m \in M$ and $t \in I$;
 \item $H'(n,1)=\Psi(\Psi'(n))$ and $H'(n,0)=n$, for each $n \in N$;
 \item\label{4} $q(H(n,t))=q(n)$ and $p'(H(n,t))=q'(n)$, for each $n \in N$ and $t \in I$.
\end{enumerate}
If $\Psi$ and $\Psi'$ are inverses of each other, then the HF fibrant spans are said to be \emph{isomorphic.}
\end{Definition}
\noindent Standard arguments prove that indeed this defines an equivalence relation on the class of all HF fibrant spans, from $B$ to $B'$.
An equivalence class of HF fibrant spans, from $B$ to $B'$, will usually be  denoted $[(p,M,p')]\colon B \to B'$.

Using the context and notation of Definition \ref{equiv-of-HF-spans}, we recall  that pullbacks along fibrations are homotopy limits.
 Given  that $\langle p,p' \rangle\colon M \to B \times B'$ and  $\langle q,q' \rangle\colon N \to B \times B'$ are fibrations, several conditions in the definition of equivalence between HF fibrant spans are, in fact, redundant. By Lemma \ref{miracle}, it follows that:
 \begin{Lemma}\label{miracle_2} Two HF fibrant spans,  
$(p,M,p')\colon B\to B'$ and $(q,N,q')\colon B\to B'$, are equivalent if there exists a map, $\Psi\colon M \to N$, making the left-most diagram of \eqref{equiv} commute and such that $\Psi\colon M \to N$ is a homotopy equivalence of spaces.
\end{Lemma}

\begin{Definition}\label{idsHFb}
Given a HF space, $B$,  we define \begin{equation}\label{hff_unit} \id_B^{\HFb}:=  \vcenter{ \xymatrix@R=-5pt{ &B^I\ar[dl]_{s_B}\ar[dr]^{t_B}\\  B & &B.}}
 \end{equation}
 This is a fibrant span, as discussed in Example \ref{ids_ex}.
\end{Definition}

\begin{Lemma}[The category $\HFb$]\label{main_desc}The composition of HF fibrant spans  in Definition \ref{main3} descends to the quotient under the equivalence relation in Definition \ref{equiv-of-HF-spans}, and, with this, the identities satisfy the evident rules.

We thus have a  category, $\HFb$, whose objects are the HF-spaces, and where morphisms from $B$ to $B'$ are  equivalence classes of HF fibrant spans, connecting $B$ and $B'$.  Given a HF space, $B$, the identity in $B$ is given by  $\id_B^{\HFb}\colon B \to B$.
\end{Lemma}
\begin{proof}
That the composition descends to the quotient follows from the universal property of pullbacks\footnote{The arguments are essentially  identical to those proving that cospans of spaces and maps between them can be arranged into a bicategory; see  \cite{Dawson_et_al}.}. More precisely, suppose that  $\Psi_1\colon M_1 \to N_1$ and  $\Psi_1'\colon N_1 \to M_1$ realise a fibred homotopy equivalence between the HF fibrant spans, $(p_1,M_1,q_1)\colon B_1 \to B$ and  $(p_1',N_1,q_1')\colon B_1 \to B$.  Suppose that  $\Psi_2\colon M_2 \to N_2$ and  $\Psi_2'\colon N_2 \to M_2$ realise a fibred homotopy equivalence between the HF fibrant spans $(p_2,M_2,q_2)\colon B \to B_2$ and  $(p_2',N_2,q_2')\colon B \to B_2$. This is as in the diagram,
\begin{equation}\vcenter{\xymatrix{
 && M_1\ar[dl]_{p_1} \ar[dr]^{q_1 }\ar@/_1pc/[dd]_{\Psi_1} \ar@{<-}@/^1pc/[dd]^{\Psi'_1}  &  M_1{\times_B} M_2\ar[l]\ar[r] & M_2
\ar@/_1pc/[dd]_{\Psi_2} \ar@{<-}@/^1pc/[dd]^{\Psi'_2}
\ar[dl]_{p_2} 
\ar[dr]^{q_2 }\\
&B_1 && B && B_2 \,.\\
 && N_1\ar[ul]^{p_1'} \ar[ur]_{q_1' } &N_1{\times_B} N_2\ar[l]\ar[r]& N_2\ar[ul]^{p_2'} \ar[ur]_{q_2' }}}
\end{equation}
 The universal property of pullbacks gives maps, $(\Psi_1{\times_B} \Psi_2) \colon M_1{\times_B} M_2 \to   N_1{\times_B} N_2$, co-gluing $\Psi_1$ and $\Psi_2$, and
 $(\Psi_1'{\times_B} \Psi_2') \colon N_1{\times_B} N_2 \to    M_1{\times_B} M_2$ doing the same for the other pair.
 
Choose fibred homotopies (using the notation in \S\ref{sec:fhe}):
\begin{align*}
\Psi_1'\circ \Psi_1 \Ra[B_1 \times B]{H_1} \id_{M_1}, \quad \quad\quad \Psi_1\circ \Psi'_1 \Ra[B_1 \times B]{H_1'} \id_{N_1}, \\\Psi_2'\circ \Psi_2 \Ra[B \times B_2]{H_2} \id_{M_2}, \quad \quad\quad  \Psi_2\circ \Psi_2' \Ra[B \times B_2]{H_2'} \id_{N_2}.
 \end{align*}
 Conditions \ref{1} to \ref{4} of Definition \ref{equiv-of-HF-spans}  imply that these homotopies can be (co)glued to homotopies, ${J}\colon (M_1{\times_B} M_2) \times I \to M_1{\times_B} M_2$ and ${J'\colon} (N_1{\times_B} N_2) \times I \to N_1{\times_B} N_2$. By construction, they are such that
\begin{align*}
 (\Psi_1'{\times_B} \Psi_2')\circ (\Psi_1{\times_B} \Psi_2)& \Ra[{B_1 \times B_2}]{{J}} \id_{M_1{\times_B} M_2}, \\ (\Psi_1{\times_B} \Psi_2)\circ (\Psi_1'{\times_B} \Psi_2')&\Ra[{B_1 \times B_2}]{{J'}} \id_{N_1{\times_B} N_2}.
 \end{align*}

To handle the point about identities, let $B$ and $B'$ be HF spaces. Consider a HF fibrant span,  $(p,M,q)\colon B \to B'$. Let us prove that  we have maps $\Psi$ and $\Psi'$, as below, realising an equivalence of HF fibrant spans,
 $$\vcenter{\xymatrix@R=5pt{ &&M\ar[dd]^\Psi\ar[dl]_p\ar[dr]^{q}\\  & B & &B', \\ && B^I {\times}_B M \ar[ul]^{p'}\ar[ur]_{q'} }} \quad \textrm{ and } \quad \vcenter{\xymatrix@R=5pt{ &M\ar@{<-}[dd]^{\Psi'}\ar[dl]_p\ar[dr]^{q}\\   B & &B'. \\ & B^I {\times}_B M \ar[ul]^{p'}\ar[ur]_{q'} }}$$
Here we consider the obvious pull-back, appearing as the diamond in the diagram:
$$\xymatrix@R=3pt{ &&& 
B^I {\times}_B M
\ar[dl]_{}\ar[dr]^{} \ar@/^1.4pc/[ddrr]^>>>>>>>>>{q'} \ar@/_1.4pc/[ddll]_>>>>>>>>>{p'}\\
&&B^I\ar[dl]^{s_B}
\ar[dr]_{t_B}&&M \ar[dl]^{p}\ar[dr]_{q}\\  & B & &B  & &B'.  }
$$

We put $\Psi(m)=(\widehat{p(m)},m)$, where $\widehat{p(m)}$ is the constant path at $p(m)\in B$.  Clearly $\Psi$ is fibred. By Lemma \ref{miracle}, in the context of Lemma \ref{miracle_2}, we only need to prove that $\Psi\colon M \to B^I {\times}_B M$ is a homotopy equivalence of spaces (as opposed to a homotopy equivalence of fibred spaces).  A homotopy inverse of  $\Psi\colon M \to B^I {\times}_B M$ is given by the map $\Phi\colon B^I {\times}_B M \to M$ such that   $\Phi(\gamma,m)=m$. (Note that this is not a fibred map.)  We have that
$\Phi\circ \Psi =\id_M$. On the other hand $\Psi(\Phi(\gamma,m))=(\widehat{\gamma(1)},m)$, for each $(\gamma,m) \in B^I {\times}_B M$. (Here $\widehat{\gamma({1})})$ is the constant path at $\gamma(1) \in B.$)  The following homotopy connects $\Psi\circ \Phi$ and $\id_{B^I {\times}_B M}$, where $s \in I$,
$$ (B^I {\times}_B M) \times I \ni (\gamma,m,t)  \mapsto \Big( s \mapsto \gamma\big(t+ (1-t)s\big), m \Big)\in  B^I {\times}_B M.$$

That the resolved identity spans given by the mapping spaces  are also identities on the right is dealt with similarly.
\end{proof}

\begin{Remark} A `dual' category to  $\HFb$, whose objects are spaces, and morphisms are cofibred homotopy equivalence classes of cofibrant cospans was constructed in \cite{Torzewska,Torzewska-HomCobs}. Our methods of proofs here are very similar, but, of course, needed switching from  cofibred to fibred homotopy equivalences.
 \end{Remark}

\begin{Definition}\label{def:HFiso} We let $\HFbiso$ be the subcategory of $\CGWH$ with objects the HF spaces, and  homeomorphisms  of HF spaces as morphisms.
\end{Definition}

\begin{Lemma}\label{lem:functions_to_spans} We have a functor, $\I\colon \HFbiso \to \HFb$, given by, if $X$ is a HF space, then $\I(X)=X$, and if $f\colon X \to Y$ is a homeomorphism of HF spaces, then $\I(f)$ is the equivalence class of the span,
$$
\smash{\xymatrix@R=-5pt{ &&X^I\ar[dl]_{s_X}\ar[dr]^{f\circ t_X}\\  & X & &Y.}}
 $$
 \end{Lemma}
\noindent (It is likely that a similar functor will map the category with objects the HF spaces, and morphisms the homotopy classes of homotopy equivalences of HF spaces, to $\HFb$, but we will not consider this, nor do we  need it.)
 \begin{proof}If $f\colon X\to Y$ is a homeomorphism of HF spaces, it follows that $\I(f)=\big(s_X,X^I,f\circ t_Y\big)\colon X \to Y$ is a HF fibrant span, since $\big(s_X,X^I, t_X\big)\colon X \to X$ is a HF fibrant span. Clearly, $\I$ sends the identities in $\HFbiso$ to the identities in $\HFb$.
 
Let $f\colon X \to Y$ and $g\colon Y \to Z$ be \label{I as pseudo-functor}homeomorphisms of HF spaces. We check that $\I(g\circ f)=\I(f)\bullet \I(g)$. To see this, look at the diagram below, where the top diamond is a pullback, defining the composite $\I(f)\bullet \I(g)$:
$$ \xymatrix@R=7pt{
&& X^I \times_Y Y^I\ar[dl]_>>>>>>{\proj_1} \ar[dr]^>>>>>>{\proj_2}\\
& X^I \ar[dl]_{s_X}\ar[dr]_{f\circ t_X} && Y^I \ar[dl]^{s_Y}\ar[dr]^{g\circ t_Y} \\
 X&& Y && Z.\\
&&X^I \ar[ull]^{s_X} \ar[urr]_{g\circ f\circ t_X}
}
$$
The map, $\Psi\colon X^I\times_Y Y^I \to X^I$, such that,
$$\Psi(\gamma,\gamma')(t)=\begin{cases}
                           \gamma(2t),& \textrm{ if } t \in [0,1/2],\\
                           f^{-1}(\gamma'(2t-1)), &\textrm{ if } t \in [1/2,1],\\
                          \end{cases}
$$
is a homeomorphism that makes the obvious diagram commute. This shows that $\I(g\circ f)=\I(f)\bullet \I(g)$.
\end{proof}

\begin{Remark}
There is a more general version of the notion of equivalence of (HF) spans, as given in Definition \ref{equiv-of-HF-spans}, that will be useful slightly later on. Recall that  spans form a category, $\CGWH^\Lambda$, as noted in Remark \ref{Spans and Lambda}, in which a morphism is simply a natural transformation, of functors $\Lambda \to \CGWH$,
$$\xymatrix{(p,M,p')\ar@{=>}[d]_{(f_{-1},f_0,f_1)}\\ (q,N,q')}
\hspace*{5mm}\xymatrix{B\ar[d]_{f_{-1}}&M\ar[l]_{p}\ar[d]^{f_0}\ar[r]^{p'}&B'\ar[d]^{f_1}\\
C&N\ar[l]^{q}\ar[r]_{q'}&C'.
}
$$
As before, we take  the Hurewicz / Str\o m model structure on $\CGWH$, and the injective model structure on $\CGWH^\Lambda$. A morphism, such as $(f_{-1},f_0,f_1)$, is thus a cofibration, in that model structure, if each of $f_{-1}$, $f_0$, and $f_1$ is a cofibration in $\CGWH$, and is a weak equivalence if each of these maps is a weak equivalence. As, in the Hurewicz / Str\o m model  category structure, the weak equivalences are, in fact, `strong' homotopy equivalences, we make the following definition:
\begin{Definition} Two fibrant spans,  
$$(p,M,p')\colon B\to B'\quad\textrm{  and  }\quad(q,N,q')\colon C\to C',$$ 
are said to be homotopy equivalent if there is a morphism,
$$(f_{-1},f_0,f_1):(p,M,p')\Rightarrow (q,N,q'),$$
in which each $f_i$ is a homotopy equivalence.
\end{Definition}\end{Remark}
Of course, in this case, $(f_{-1},f_0,f_1)$ is a homotopy  equivalence\footnote{more precisely a weak equivalence in the injective model structure on the category of spans.}, and, in the setting in which $f_{-1}$ and $f_1$ are the respective identities, we retrieve the notion of equivalence given in Definition \ref{equiv-of-HF-spans}. We also note that fibrant spans are cofibrant-fibrant objects  in the injective model structure in  $\CGWH^\Lambda$, so a homotopy equivalence of fibrant spans will actually be a strong homotopy equivalence.

Given a map, $(f_{-1},f_0,f_1):(p,M,p')\Rightarrow (q,N,q'),$ of fibrant spans, we get an induced map,
$$(f_0,\langle f_{-1},f_1\rangle): \langle p,p'\rangle\to \langle q,q'\rangle,$$
of fibrations (in the sense of Definition \ref{map of fibs}), so, if $(f_{-1},f_0,f_1)$ is a homotopy equivalence, then $(f_0,\langle f_{-1},f_1\rangle)$ will be a homotopy equivalence of fibrations, (Definition \ref{he of fibs}), and by Proposition \ref{relative miracle}, there will be a homotopy inverse.

Of course, if $b\in B$ and $b'\in B'$, then by Corollary \ref{h.e. on fibres}, or by using the fact that $(f_{-1},f_0,f_1)$ is a strong homotopy equivalence, we have the following:
\begin{Proposition}\label{cor:h.e. of spans}
If $(f_{-1},f_0,f_1):(p,M,p')\Rightarrow (q,N,q')$  is a homotopy equivalence of fibrant spans, then, for any $(b,b')\in B\times B'$, the induced map on fibres,
$$\langle p,p'\rangle^{-1}(b,b')\to \langle q,q'\rangle^{-1}(f_{-1}(b),f_1(b')),$$
is a homotopy equivalence.
\end{Proposition}

\subsection{A family of functors, $\FRp{s}\colon \HFb\to \Vect$, derived from the homotopy  content}
The results in this subsection are closely related to some  given in \cite{Galvezetal}, where they are stated in the language of $\infty$-groupoids.  They were, in fact, essentially implicit in \cite[Section 4]{Quinn}. The indexation of the \emph{family} of functors, with $s \in \C$, is, however, a generalisation of the \emph{$\alpha$-degroupoidification} set-up introduced by Baez, Hoffnung and Walker in \cite{baezetal}, Proposition 3.3. 

The setting here is particularly suited  to constructing Quinn's finite total homotopy TQFT, and explicitly to compute it in a number of cases, as well as moving towards extended versions of Quinn's finite total homotopy TQFT.  The point about the parameter, $s$, is then that, for $s=0$, one has  Quinn's theory, as we will see shortly, but, for other values of $s$,  one also gets functors linked to other TQFTs, in the normalisations in which they were initially constructed.
\subsubsection{Spatial slices}
 We  now introduce some nomenclature and notation that we will use later on to simplify some formulae.

 Let $ (p,M,p')\colon B \to  B'$ be a HF fibrant span. Let $b \in B$, $b' \in B'$.

\begin{Notation}\label{space-els} We  denote the fibre $\langle p,q \rangle^{-1}(b,b')$ of $\langle p,q \rangle\colon M \to B\times B'$ as:
$$
\{b \big |(p,M,p') \big |b'\}\stackrel{abbr.}{=} \{b \big |M \big |b'\}.$$ To $\{b|M|b'\}$ we call the \emph{spatial slice of} $(p,M,p')\colon B \to  B'$, \emph{at $b\in B$ and $b'\in B'$}. \end{Notation}
We note that  the abbreviated notation, $\{b|M|b'\}$, does not show the dependence on $p\colon M \to B$ and $p'\colon M \to B'$, but we will use it more often than the complete $\{b \big |(p,M,p') \big |b'\}$ or $\langle p,q \rangle^{-1}(b,b')$, so as not to overload the various formulae.

\emph{We also define the following spaces, also called \emph{spatial slices}, but \emph{over the various subsets of $B$, or $B'$,} as indicated:
\begin{align*}
\{b \big |M \big |\PC_{b'}{(B')}\}&:=\{m \in M\colon p(m)= b \textrm{ and } p'(m)\in \PC_{b'}{(B')}\},\\
\{\PC_{b}{(B)}\big |M \big |b'\}&:=\{m \in M\colon p(m)\in \PC_{b}{(B)} \textrm{ and } p'(m)=b'\},\\
\intertext{and}
\{\PC_{b}{(B)} \big |M \big |\PC_{b'}{(B')}\}&:=\{m \in M\colon p(m)\in \PC_{b}{(B)}\textrm{ and } p'(m)\in \PC_{b'}{(B')}\}.
\end{align*}}
\begin{Remark}\label{space_els_props} We collect some useful facts about spatial slices.
\begin{itemize}[leftmargin=0.6cm]
\item The fibration $\langle p,p' \rangle \colon M \to B \times B'$ restricts to a fibration,
  $$\{\PC_{b}{(B)} \big |M \big |\PC_{b'}{(B')}\} \to B \times B',$$ and, by  Lemma \ref{All-fibrations}, $p\colon M \to B$ and $p'\colon M \to B'$ restrict to fibrations
 $$\{\PC_{b}{(B)} \big |M \big |\PC_{b'}{(B')}\} \to B\textrm{ and  }\{\PC_{b}{(B)} \big |M \big |\PC_{b'}{(B')}\} \to B',$$ the fibres, or more generally, the inverse images are, respectively, the spaces,
 $\{b \big |M \big |b'\}$, $\{b \big |M \big |\PC_{b'}{(B')}\}$ and $\{\PC_{b}{(B)}\big |M \big |b'\}.$ 
\item In particular, the homotopy type of spaces,   $\{b \big |M \big |b'\}$, $\{b \big |M \big |\PC_{b'}{(B')}\}$ and $\{\PC_{b}{(B)}\big |M \big |b'\}$,  depends only on the path-components of $b\in B$ and $b'\in B'$.

\item Lemma \ref{All-fibrations} also gives that $p$ and $p'$ restrict to fibrations,
$$\{b | M | \PC_{b'}(B')\} \to \PC_{b'}(B') \textrm{ and }\{\PC_b(B) | M | b'\} \to \PC_b(B).$$ The fibres, again,  have the form $\{b |M | b'\}$. 

\item All the spaces  below are HF,  which follows from Lemma \ref{all-HF},  $$\{b \big |M \big |b'\}, \quad  \{b \big |M \big |\PC_{b'}{(B')}\}, \quad  \{\PC_{b}{(B)}\big |M \big |b'\} \quad \textrm{ and } \quad \{\PC_{b}{(B)}\big |M \big | \PC_{b'}{(B')}\}.$$ We can therefore take their homotopy  content.
 \end{itemize}
 
 \end{Remark}

\begin{Lemma}\label{rels_matrix_els}
 Let $(p,M,p')$ be a HF fibrant span. We have (in $\Q$):
\begin{align*}
 \chi^\pi(\{\PC_{b}{(B)} \big |M \big |\PC_{b'}{(B')}\})&= \chi^\pi(\{b \big |M \big |\PC_{b'}(B)\}) \,\,\chi^\pi(\PC_{b}{(B)})\\
 &= \chi^\pi(\{\PC_{b}(B) \big |M \big |b'\}) \,\,\chi^\pi(\PC_{b'}{(B')})\\
  &= \chi^\pi(\{b \big |M \big |b'\}) \,\,\chi^\pi(\PC_{b}{(B)}) \,\,\chi^\pi(\PC_{b'}{(B')}).
\end{align*}
\end{Lemma}
\begin{proof}Follows from Lemma \ref{main1} applied to the  fibrations  in Remark \ref{space_els_props}.
\end{proof}
\subsubsection{Matrix elements}
 
 First recall Definition \ref{space-els} and  Lemma \ref{rels_matrix_els}, and, as we want to define some linear maps, recall also the general notation and terminology relating to matrix elements, mentioned at the start of the paper, on page \pageref{matrix elements}.

Let $(p,M,p')\colon B \to  B'$ be a HF fibrant span.  We  introduce a matrix over $\C$, parametrised by a complex valued index, $s$.

\begin{Definition}\label{matrix_els} 
Let $s \in \C$. Given (non-empty) path-components, $\PC_b(B)$, $\PC_{b'}(B')$,  of $B$ and $B'$, define the following complex-valued `matrix-elements',
\begin{align*}
\big \langle &\PC_b(B) \big | \FRh{s}(p,M,p') \big | \PC_{b'}(B') \big \rangle\\&:=\chi^\pi\big(\{b \big |M \big |b'\}\big)\,\, \big({\chi^\pi(\PC_{b}{(B)})}\big)^{s}\,\, \big ({\chi^\pi(\PC_{b'}{(B')})}\big)^{1-s}\\
&=\chi^\pi(\{\PC_b(B) \big |M \big |\PC_{b'}{(B')}\})\,\, \big({\chi^\pi(\PC_{b}{(B)})}\big)^{s-1}\,\, \big ({\chi^\pi(\PC_{b'}{(B')})}\big)^{-s}\in \C.
\end{align*}
\end{Definition}
We note, as well, that the homotopy type of $\{b|(p,M,p')|b'\}$ depends only on the path-components, in $B$ and $B'$, that $b$ and $b'$,  respectively, belong to, so $\chi^\pi(\{b \big |M \big |b'\})$ is indeed a function of $\PC_b(B)$ and  $\PC_{b'}(B')$, only.

The following result is essentially in \cite{Quinn}.  A version of this result for groupoids and spans appears in \cite[Theorem 41]{baezetal}, whilst a  version for $\infty$-groupoids is in \cite[Proposition 8.2]{Galvezetal}. We abbreviate $\PC_b(B)$ as $[b]_B$ to save space.

\begin{Lemma}\label{mult} Let $s\in \C$. The matrix elements corresponding to  ${\FRh{s}}$ are multiplicative with respect to composition of HF fibrant spans. Explicitly,
consider  HF fibrant spans, $(p,M,p')\colon B \to B'$, $(p'',M',p''')\colon B' \to B''$, and their composition, connecting $B$ to $B''$,  $(P_L,M\times_{B'} M',P_R)=(p,M,p')\bullet (p'',M',p''')$, defined from the diagram below, where the middle diamond is a pullback,
\begin{equation}\vcenter{
\xymatrix@R=3pt{ &&&
M\times_{B'} M'\ar@/_1.4pc/[ddll]_>>>>>>>{P_L} \ar@/^1.4pc/[ddrr]^>>>>>>>{P_R}
\ar[dd]^P \ar[dl]_{q}\ar[dr]^{q'}\\ &&M\ar[dl]^p\ar[dr]_{p'}&&M'\ar[dl]^{p''}\ar[dr]_{p'''}\\  & B & &B'  & &B'',  }}
\end{equation}
 then, given $b \in B$ and $b'' \in B''$, we have
\begin{multline*}
 \big \langle
 [b]_B  \big | \FRh{s}(P_L ,M\times_{B'} M',P_R) \big | [b'']_{B''} \big \rangle=\\ \sum_{ [b'] \in \pi_0(B') } \big \langle [b]_{B} \big |\FRh{s}(p,M,p') \big | [{b'}]_{B'} \big \rangle   \,\big \langle [{b'}]_{B'} \big |\FRh{s}(p',M',p'') \big |[b'']_{B''}\big \rangle.
 \end{multline*}
\end{Lemma}
\begin{proof}
We  apply Theorem \ref{main2} to the fibration,
 $$P_{b,b''}\colon \{b  \big | M \times_{B'} M'  \big | b''\}\to B',$$ obtained by restricting $P \colon M\times_{B'} M' \to B'$ to $\{b  \big | M \times_{B'} M'  \big | b''\}=\langle P_L,P_R\rangle^{-1}(b,b')$. This map, $P_{b,b''}$, is a fibration, since, by Lemma \ref{main3},  the map $$\big \langle P_L, P,  P_R\big \rangle\colon M \times_{B'} M'\to B \times B' \times B'' $$ is a fibration; cf. the proof of the second part of Lemma \ref{All-fibrations}.

We then have
\begin{align*}
\big \langle
 \PC_b(B) & \big | \FRh{s}(P_L,M\times_{B'} M',P_R) \big | \PC_{b''}(B'') \big \rangle
\\&=\chi^\pi(\PC_{b}(B))^s\,\,\chi^\pi(\PC_{b''}(B''))^{1-s}\,\,\chi^\pi (\{
 b  \big | M\times_{B'} M' \big | b''\})\\
 &= \chi^\pi(\PC_{b}(B))^s\,\,\chi^\pi(\PC_{b''}(B''))^{1-s} \sum_{ [b'] \in \pi_0(B') } \chi^\pi(P_{b,b''}^{-1}(b'))\,\, \chi^\pi(\PC_{b'}(B')).
 \end{align*}
Now note that, as spaces, we have, $$P_{b,b''}^{-1}(b')=\{b \big |(p,M,p') \big |b'\} \times \{b' \big |(p',M',p'') \big |b''\}, $$ so
$$\chi^\pi(P_{b,b''}^{-1}(b') ) = \chi^\pi(\{b \big |(p,M,p') \big |b'\})\,\chi^\pi( \{b' \big |(p',M',p'') \big |b''\}).$$
This yields the main formula in the statement of the lemma.
\end{proof}

\begin{Lemma}\label{OmegaiHF} Let $B$ be a path-connected HF space. Choose a base-point $*$.
Let $\Omega_*(B)$ denote  the pointed loop space of $(B,*)$, that is, the space of all loops in $B$ starting and ending in $*$, then $\Omega_*(B)$ is HF, and   $\chi^\pi(\Omega_*(B))=1/\chi^\pi(B)$.
\end{Lemma}
This is the analogue of \cite[Lemma 3.10]{Galvezetal}, which is for $\infty$-groupoids.
\begin{proof}
 Let $\mathcal{F}_*(B)$ denote
 the space of all paths, $\alpha\colon [0,1]\to B$, starting at $*$. There is, of course, a fibration, $ \mathcal{F}_*(B)\to B$, sending a path to its second end-point,  and we note that  $\mathcal{F}_*(B)$ is contractible
. Since both  $ \mathcal{F}_*(B)$ and $ B$ are HF, it follows that the fibre at $*$, which is $\Omega_*(B)$, is HF; see Lemma \ref{2of3}. Since  $\mathcal{F}_*(B)$ is contractible, we have that $\chi^\pi(\mathcal{F}_*(B))=1$, and  Lemma \ref{main1} gives that $\chi^\pi(\Omega_*(B))\, \chi^\pi(B)=\chi^\pi(\mathcal{F}_*(B))=1$.
\end{proof}

\begin{Lemma}\label{ids}
Suppose $B$  is  a HF space and let $b,b'\in B$.  Let $s \in \C$, then 
\begin{align*}
 \big \langle \PC_{b}(B) \big |\FRh{s}(s_X, B ^{[0,1]},t_Y) \big |  \PC_{b'}(B) \big\rangle &=\delta\big(\PC_b(B),  \PC_{b'}(B)\big).
\end{align*}
More generally, let $f\colon B \to B'$ be homeomorphism of HF spaces, then,
\begin{align*}
 \big \langle \PC_{b}(B) \big |\FRh{s}(s_X, B ^{[0,1]},f\circ t_Y) \big |  \PC_{b'}(B') \big\rangle &=\delta\big(\PC_{b}(B),  \PC_{f^{-1}(b')}(B)\big).
 \end{align*}
\end{Lemma}
\noindent Here, if $X$ is a set, which we will need  to be the set of path-components of $B$, we take $\delta\colon X \times X \to \{0,1\}$, to be such that $\delta(x,y)$ is $0$, if $x\neq y$, and $\delta(x,x)=1$.
\begin{proof} We prove the most general case. First of all note that if $\PC_b(B) \neq \PC_{f^{-1}(b')}(B)$, then
\[
\big \langle \PC_{b}(B) \big |\FRh{s}(s_X, B^{[0,1]},f\circ t_Y) \big |  \PC_{b'}(B') \big\rangle=0,
\]
as in this case  $\{ b \big | (s_B,B ^{[0,1]},f\circ t_B) \big | b'\}$ is empty. On the other hand,
\begin{align*}\big \langle \PC_{b}(B) &\big |\FRh{s}(s_X, B ^{[0,1]},f\circ t_Y) \big |  \PC_{f(b)}(B') \big\rangle\\&=
\chi^\pi\big(\{ b \big | (s_B,B ^{[0,1]},f\circ t_B) \big | f(b)\}\big)\,\,\chi^\pi(\PC_{b}(B))^s\,\, \,\,\chi^\pi(\PC_{f(b)}(B'))^{1-s}\\&=
\chi^\pi\big(\{ b \big | (s_B,B ^{[0,1]},  t_B) \big | b\}\big)\,\,\chi^\pi(\PC_{b}(B))^s\,\, \,\,\chi^\pi(\PC_{b}(B))^{1-s}\\
&=\chi^\pi\big(\Omega_b(\PC_{b}(B))\big)\chi^\pi(\PC_{b}(B))=1,
\end{align*}
where we have used Lemma \ref{OmegaiHF}.
\end{proof}

We thus have fixing $s \in \C$ (the notation $\hpi_0(B)$ is explained at the end of \S\ref{sec:conventions_top}):

\begin{Theorem}\label{thm:FR}
 There is  a  functor, $\FRp{s}\colon \HFb \to \Vect_\C$, such that:
 \begin{itemize}[leftmargin=1cm]
 \item if $B$ is a HF space, then $\FRp{{s}}(B) =\C(\hpi_0(B))$,  the free vector space over the set of all path-components of $B$;
 
 \item for  $[(p,M,p')]\colon B \to B'$,  a  morphism in $\HFb$,  the matrix elements for the linear map, $$\FRp{s}\big([(p,M,p')] \big) \colon \C(\hpi_0(B))\to \C(\hpi_0(B')),$$ are given by
$$ \big \langle \PC_b(B) \big | \FRp{s}\big( [ (p,M,p')] \big) \big | \PC_{b'}(B') \big \rangle:= \big \langle \PC_b(B) \big | \FRh{s}(p,M,p') \big | \PC_{b'}(B') \big \rangle,$$
 for given path-components $\PC_b(B)\in \hpi_0(B)$ and $\PC_{b'}(B')\in \hpi_0(B')$.
\end{itemize} 
\end{Theorem}
\begin{proof}
Compatibility of $\FRp{s}$ with the composition and  identities in $\HFb$  follows from Lemmas \ref{mult} and \ref{ids}, respectively. What we have not yet shown is that the matrix elements $\FRh{s}$ are invariant under equivalence of spans. This follows by using lemma  \ref{fibre_hom}, combined with Definition \ref{equiv-of-HF-spans} and  Proposition \ref{cor:h.e. of spans}.
\end{proof}

\begin{Lemma}\label{lem:mon0} We have a functor, $ \Times \colon \HFb\times \HFb \to \HFb$, sending the equivalence class of,
$$\left(\vcenter{\xymatrix@R=-5pt{&M\ar[dl]_p\ar[dr]^{q}&\\ X && Y  } },\vcenter{\xymatrix@R=-5pt{&M'\ar[dl]_{p'}\ar[dr]^{q'}\\X' & &Y'  }}\right),$$
 to the equivalence class of,
 $$\xymatrix@R=-5pt{ &M\times M'\ar[dl]_{p\times p'}\ar[dr]^{q\times q'}\\   X\times X' & &Y\times Y'.}$$
 (Note that the latter span is fibrant, since the product of  two  fibrations is a fibration.)
\end{Lemma}
\begin{proof}
This follows from straightforward calculations using Lemma \ref{lem:sums.prods}.
\end{proof}
It can furthermore be proved that $\HFb$ is a monoidal category, with this tensor product, but we will not use this here. For details in the dual case of cofibrant cospans, see \cite{Torzewska,Torzewska-HomCobs}.

\begin{Lemma}\label{lem:mon1}
 Let $\Otimes\colon  \Vect \times \Vect \to \Vect$ denote the tensor product functor in $\Vect$. We have a natural isomorphism of functors from $\HFb\times\HFb$ to $\Vect$,
 $$\eta\colon \Otimes \circ \big (\FRp{s}\times \FRp{s}\big) \To \FRp{s}\circ \Times. $$
It is such that,  given HF spaces $X$ and $X'$, and $x\in X, x'\in X'$, then
 $$\eta_{X,X'} \big (\PC_{x}(X) {\otimes} \PC_{x'}(X') \big)=\PC_{(x,x')}(X\times X').$$
\end{Lemma}
\begin{proof}
 If $X$ and $X'$ are spaces, then we have a natural bijection from $\hpi_0(X) \times \hpi_0(X')$ to $\hpi_0(X \times X')$, sending $\big (\PC_{x}(X),  \PC_{x'}(X') \big)$ to $\PC_{(x,x')}(X\times X')$. The naturality of $\eta$ then follows from the calculation below.

Consider two HF fibrant spans,
$$\vcenter{\xymatrix@R=-5pt{&M\ar[dl]_p\ar[dr]^{q}&\\ X && Y  } }  \quad  \textrm{ and } \quad \vcenter{\xymatrix@R=-5pt{&M'\ar[dl]_{p'}\ar[dr]^{q'}\\X' & &Y'.  }}$$
If $x \in X, x'\in X', y \in Y$ and $y'\in Y'$, we have:
\begin{multline*}
 \Big \langle  \PC_{(x,x')} (X \times X') \,\, | \FRh{s}\big ( (p,p'), M\times M', (q,q') \big) \,|  \PC_{(y,y')} (Y\times Y')    \Big \rangle \\
 =\chi^\pi\left (  \big \{ (x,x') \, |\,\big ( (p,p'), M\times M', (q,q') \big)  \,|\,  (y,y') \big\} \right ) \\
 \chi^\pi(\PC_{(x,x')} (X \times X') \big)^s\,\,  \chi^\pi\big( \PC_{(y,y')} (Y \times Y') \big)^{1-s}.
 \end{multline*}
Note that we have homeomorphisms,
\begin{align*}
  \big \{ (x,x' | \big ( (p,p'), M\times M', (q,q') \big) |  (y,y') \big\} &\cong
  \{ x |  (p,M,q) | y\} \times \{x' | (p',M',q') | y'\},\\
 \PC_{(x,x')}(X\times X') &\cong \PC_x(X) \times \PC_{x'}(X'),\\ \PC_{(y,y')}(Y\times Y') &\cong \PC_y(Y) \times \PC_y(Y').\end{align*}
It therefore follows that,
\begin{multline*}
 \big \langle  \PC_{(x,x')} (X \times X') \, |\, \FRh{s}\big ( (p,p'), M\times M', (q,q') \big) \,| \,  \PC_{(y,y')} (Y\times Y')    \big \rangle =\\
 \big \langle  \PC_{x} (X ) \, |\, \FRh{s}\big ( (p, M, q ) \big) \,|  \PC_{y} (Y)   \big \rangle \,
 \big \langle  \PC_{x'} (X') \, | \, \FRh{s}\big (( p', M',  q') \big) \,|  \PC_{y'} (Y')    \big \rangle.
\end{multline*}
The last equation follows from Lemma \ref{lem:sums.prods}.
\end{proof}

Finally, note the following result, that follows from  straightforward calculations.
\begin{Proposition}\label{nat:trans_st}
Let $s,t$ be complex numbers. We have a natural isomorphism, $\eta^{s,t}\colon \FRp{s}  \To \FRp{t}$, of functors, from $\HFb$ to $\Vect_\C$, which is such that, if $X$ is a space, and $x \in X$, then, $$\eta^{s,t}_X(\PC_x(X))=\chi^\pi\left( \PC_x(X)\right)^{s-t} \,\,\PC_x(X) .$$
\end{Proposition}

 \section{The construction of Quinn's finite total homotopy TQFT}\label{Classical Quinn}


After the homotopical constructions in the previous section, we now explain the construction of Quinn's finite total homotopy TQFT, $\FQ{\Bc} \colon \tcob{n} \to \Vect_\C$, where $n$ is a non-negative integer, and $\Bc$ is a homotopy finite space.  As in \cite[3.3 Proposition]{baezetal}, this is one of a family of such constructions, $\FQp{\Bc}{s} \colon \tcob{n} \to \Vect_\C$, one for each $s\in \C$, albeit all related by monoidal natural isomorphisms.

First let us recall some basic definitions.  Most of the time, we will be able to do the necessary constructions without this level of detail, just as Quinn does in the primary source, \cite[Lecture 4]{Quinn}, where he discusses a version of the theory using just CW complexes, but very occasionally, these results, or their consequences,  are needed.  They also help to tie Quinn's theory into the general theory of TQFTs.

\subsection{Cobordism categories}\label{cobordism categories}
A topological manifold of dimension $n$ is a Hausdorff and second countable topological space, $S$, such that each point of $S$ has a neighbourhood homeomorphic to an open subset of the upper half-plane of $\mathbb{R}^n$. A smooth manifold, $(S,\smt_S)$, is a pair, consisting of a topological manifold, $S$,  and a smooth structure, $\smt_S$, on $S$; see, for instance, {\cite{Hirsch}  or \cite[\S 1]{Milnor}}.  We call $S$ the underlying topological manifold of $(S,\smt_S)$, and  will usually abbreviate $(S,\smt_S)$ to $S$, when the context  makes this unambiguous.
We note that a topological manifold being smooth is a structure, not a property, and some topological manifolds do not have a smooth structure at all.

If $M$ is a compact smooth manifold, then it can be given a finite triangulation, and, in particular, it can be given the structure of a finite CW-complex, see \cite{MunkresDiff}.
We also  note that if $M$ is a smooth manifold with border, then we can find, again see \cite{MunkresDiff},  a triangulation of the \emph{pair}, $(M,\d M )$, making $\d M $ a subcomplex of  $M$.  In particular, the inclusion, $\iota\colon \d M \to M$, is a cofibration.

For a positive integer $n$, we let $\cob{n}$ denote the monoidal category of compact smooth manifolds and cobordisms between them. Details are  discussed in many places in the literature,  e.g. \cite{Lec_TQFT,Kock:TQFT:2003,Milnor}. We note  that we make no assumption that orientations on manifolds and cobordisms are given, or even that they exist.

The class of objects of $\cob{n}$ is given by all compact smooth-manifolds of dimension $n$. Given compact smooth $n$-manifolds, $\Sigma$ and $\Sigma'$, morphisms in  $\cob{n} $ from $\Sigma$ to $\Sigma'$, are equivalence classes of cobordisms, $(i,S,j)\colon \Sigma \to \Sigma'$. Here a cobordism is a cospan of compact smooth manifolds and smooth maps, as below,
\begin{equation}\label{eq:cob-1st}
\cspn{\Sigma}{i}{\Sigma'\, ,}{j}{S}
\end{equation}
where $i$ and $j$ are smooth maps inducing a diffeomorphism, $\langle i,j\rangle\colon \Sigma \sqcup \Sigma' \to \partial S$.

Two cobordisms, $(i,S,j),(i',S',j')\colon \Sigma \to \Sigma'$,
are considered equivalent if a smooth diffeomorphism, $f\colon S \to S'$, exists, making the diagram below commute,
$$\xymatrixrowsep{.2pc}\xymatrix{ & S\ar[dd]^f &\\
              \Sigma\ar[ur]^i\ar[dr]_{i'} &&  \Sigma'\, .\ar[ul]_j\ar[dl]^{j'}\\ & S'&}
$$
We will use a hopefully evident notation for the equivalence classes of cobordisms.

The composition of morphisms, $[(i,S,j)]\colon \Sigma \to \Sigma'$ and $[(i',S',j')]\colon \Sigma' \to \Sigma''$, is done as follows.  We first consider the pushout, $S\sqcup_{\Sigma'} S'$, in $\CGWH$, as in the diagram below (the nodes in the first and second rows contain the underlying topological manifolds of the corresponding smooth manifolds):
$$\xymatrix@R=5pt{&\Sigma\ar[dr]_{i} && \Sigma'\ar[dl]_{j} \ar[dr]^{i'} && \Sigma'' \ar[dl]^{j'}\\
&& S\ar[dr]_{k} \ar@{}[rr]_{\widehat{\qquad}}&& S'\ar[dl]^{k'}\\
&&& S\sqcup_{\Sigma'} S'
}.$$

The topological space,  $S\sqcup_{\Sigma'} S'$, is a topological manifold, see  {\cite[{\S 1}]{Milnor}} or  \cite[\S 8.2]{Hirsch},  and $\langle  k\circ i, k'\circ j' \rangle\colon \Sigma \sqcup \Sigma''\to\partial( S\sqcup_{\Sigma'} S')$ is a homeomorphism. This yields the cospan below in $\CGWH$, (the nodes, as yet, only denote topological manifolds), $$ \cspnc{\Sigma}{k\circ i}{\Sigma''}{k'\circ j'}{S\sqcup_{\Sigma'} S'}.$$

As is well known, $S \sqcup_{\Sigma'} S'$ can be given a smooth structure, which `restricts' to the smooth structures in $S$ and in $S'$. This smooth structure, despite not being unique, as it depends on the choice of a collar of $\Sigma'$ in $S$ and in $S'$, is unique up to a diffeomorphism, which is the identity on $\partial (S\sqcup_{\Sigma'} S') $;  {for discussion see \cite[{\S 3}]{Milnor}.

The composition, $\bullet$, of morphisms, $[(i,S,j)]\colon \Sigma \to \Sigma'$ and $[(i',S',j')]\colon \Sigma' \to \Sigma''$, in $\cob{n}$ is then given as
$$ \big([(i,S,j)]\colon \Sigma \to \Sigma'\big)\bullet \big([(i',S',j')]\colon \Sigma' \to \Sigma''\big):=\big ( [(k\circ i, S\sqcup_{\Sigma'} S', k'\circ j')]\colon \Sigma \to \Sigma''\big).$$

Given a closed smooth   $n$-manifold, $\Sigma$, the identity, $\id_\Sigma\colon\Sigma \to \Sigma$, in $\cob{n}$, is the equivalence class of the cobordism below, where $\iota_i^\Sigma(x)=(x,i)$, for all $i \in \{0,1\}$,
 $$\cspnc{\Sigma}{\iota_0^\Sigma}{\Sigma}{\iota_1^\Sigma}{\Sigma \times I}.$$

To finalise this discussion, we note that if we forget the smooth structures in a cobordism, \eqref{eq:cob-1st},   each  cobordism gives a cospan in $\CGWH$. Crucially for what follows, note that  the induced map $\langle i,j\rangle\colon \Sigma \sqcup \Sigma' \to S$ is
 a cofibration, as it is given by the inclusion of the boundary of a smooth manifold into the manifold. We say that the cospan is \emph{cofibrant}. Such cofibrant cospans are treated in \cite{Torzewska,Torzewska-HomCobs}.

\subsubsection{The symmetric monoidal structures in $\Diff{n}$ and $\cob{n}$ }\label{moncob}

Let $\Diff{n}$ denote the category of closed $n$-manifolds and diffeomorphism between them. We
have a functor, $\I'\colon\Diff{n} \to \cob{n}$, {which is the identity on objects,} and such that $\I'(f\colon \Sigma \to \Sigma')$ is the equivalence class of the cobordism below,\label{I-construction}
$$\cspnc{\Sigma}{\iota_0^\Sigma}{\Sigma'\,.}{\iota_1^{\Sigma'}\circ f^{-1}}{\Sigma \times I} $$
The proof of this fact is dual to that of Lemma \ref{lem:functions_to_spans}. This, of course, implies that each $\I'(f)$ will be an invertible morphism in $\cob{n}$, i.e., that the cobordism is invertible up to equivalence.
It is well known that this functor, $\I'\colon \Diff{n} \to \cob{n}$, descends to the category with objects the closed $n$-manifolds and morphisms isotopy classes of diffeomorphisms. This will, however, not be used in the following.

\begin{Remark} \label{left coleg}There is another way  to obtain a cobordism from a diffeomorphism. Instead of using $f^{-1}$ in the right (co)leg, we  can use $f$ in the left one,
$$\cspnc{\Sigma}{{\iota_0^{\Sigma}\circ f}}{\Sigma'\,.}{\iota_1^{\Sigma'}}{\Sigma' \times I} $$
This gives an equivalent cobordism, as the one used when constructing $\I'(f)$.
We will later on `categorify' this second $\I'$-construction, when we are considering the symmetric monoidal bicategory structure on $\tcob{n}$; see \S \ref{sec:mon_tcob}.
\end{Remark}

 Recall that both $\cob{n}$ and $\Diff{n}$ are symmetric monoidal categories, where the tensor product on objects is given by the disjoint union, $\Sigma \sqcup \Sigma'$, of closed $n$-manifolds, $\Sigma$ and $\Sigma'$. For both symmetric monoidal bicategories, the  unit object is the empty manifold, $\emptyset$. In  $\Diff{n}$,  the tensor product of morphisms is achieved as in $(\CGWH,\sqcup)$, i.e. by performing the disjoint union of diffeomorphisms, namely
 $$(f_1 \colon \Sigma_1 \to \Sigma_1') \sqcup (f_2 \colon \Sigma_2 \to \Sigma_2')=(f_1 \sqcup f_2)\colon \Sigma_1\sqcup \Sigma_2 \to \Sigma'_1\sqcup \Sigma'_2.  $$
 The associativity constraints, braiding, etc., in $\Diff{n}$ are   as those in  $(\CGWH,\sqcup)  $. So, for instance, given closed smooth manifolds $\Sigma, \Sigma',\Sigma''$, the associativity contraint is given by the obvious diffeomorphism $\alpha_{\Sigma,\Sigma',\Sigma''}\colon (\Sigma \sqcup \Sigma')\sqcup \Sigma'' \to \Sigma \sqcup ( \Sigma'\sqcup \Sigma'')$.

In $\cob{n}$, the monoidal structure is based on the functor, $$\Sqcup\colon \cob{n}\times\cob{n} \to \cob{n},$$   so is obtained from the disjoint union of cobordisms, which descends to their equivalence classes. (This is dual to the construction in Lemma \ref{lem:mon0}.)
 Crucially for what follows, the associativity and unit
constraints, and the braiding in $\cob{n}$, are obtained from those of $\Diff{n}$ by applying $\I'\colon \Diff{n} \to \cob{n}$.

\subsection{Quinn's results on HF function spaces}\label{sec:QuinnFunSpace}In his original reference \cite{Quinn}, Quinn uses various results on mapping spaces whose  proof is not immediate.

Let $\Bc$ be a HF space. 
\begin{Lemma}[Quinn]\label{lem:mappingfromCW}
 Let $X$ be a finite CW-complex, then \smash{$\Bc^X$} is a HF-space.
\end{Lemma}
\begin{proof}The proof follows from an induction on the number of cells of $X$, by making use of the following lemma in each induction step; cf.  \cite[Chapter 4]{Quinn}.
\end{proof}

\begin{Lemma} Let $i$ be a non-negative integer. Let a space, $Y$, be obtained from the CW-complex $X$ by attaching an $i$-cell. Suppose $\Bc^X$ is HF, then so is $\Bc^Y$.
\end{Lemma}
\begin{proof}
Let $Y$ be obtained from  $X$ by attaching an $i$-cell along $f\colon S^{i-1} \to X$. (Here, by convention $S^{-1}=\emptyset$.)
We therefore have a pushout diagram,
$$\xymatrix@R=14pt{
& S^{i-1} \ar[d] \ar[rr]^f  & \ar@{}[dr]|>>>{\pushout}& X\ar[d]\\
&D^i \ar[rr] &&Y\, ,
} 
$$
where both vertical arrows are induced by inclusion of subcomplexes, hence they are cofibrations. Passing to function spaces, and using the fact that $\CGWH$ is  monoidal closed,
we have a pullback diagram, where the vertical arrows are moreover fibrations, given that they are `dual' to cofibrations, 
$$\xymatrix@R=14pt{
&\Bc^Y\ar[rr] \ar@{}[dr]|<<<{\pullback} \ar[d]&& \Bc^{D^{i}}\ar[d]  \\
&\Bc^X \ar[rr]^{{f^*}} &&\Bc^{S^{i-1}}
\,. }
$$
Concretely, each vertical arrow is obtained by restricting a function  defined on a CW-complex to a subcomplex.

We now apply Lemma \ref{pull_HF}. Since $\Bc^X$ and $\Bc^{D^i}$ are HF, the first by assumption, the second since $\Bc^{D^i}$ is contractible,  the proof is reduced to proving that  $\Bc^{S^{i-1}}$ is HF. This is proved in the following lemma.
\end{proof}
\begin{Lemma}
 Given any non-negative integer ${n}$, the space, $\Bc^{S^{{n}}}$, is HF.
\end{Lemma}
\begin{proof}
 Again the proof is by induction in ${n}$. The base case follows from the fact that  $\Bc^{S^{0}}=\Bc^{\{0,1\}}\cong \Bc \times \Bc$. The induction step follows by observing that we have the following pullback diagram, where the vertical arrows are fibrations:
 $$ \xymatrix
@R=14pt {&\Bc^{S^n}\ar@{}[dr]|<<<{\pullback}\ar[rr]\ar[d] && \Bc^{D^n}\ar[d]\\
 &\Bc^{D^n} \ar[rr]&& \Bc^{S^{n-1}}
 }
 $$
 and noting, once again,  that $\Bc^{D^n}$ is $HF$, since it is contractible.
\end{proof}

\begin{Lemma}\label{1dual} Let $n$ be a non-negative integer, and $\Bc$ be a HF-space.
\begin{enumerate}[leftmargin=1cm]
\item\label{it1} There is a  functor, $\Fc_\Bc^0 \colon (\Diff{n})^{\mathrm{op}}\to \HFbiso$ (see Definition \ref{def:HFiso}), that sends a closed and smooth $n$-manifold, $\Sigma$, to $\Bc^\Sigma$ and a diffeomorphism,  $f\colon \Sigma \to \Sigma'$, to the induced map, $f^*\colon \Bc^{\Sigma'} \to \Bc^\Sigma$, sending $\phi\colon \Sigma' \to \Bc$ to  $\phi\circ f \colon \Sigma \to \Bc$.
\item\label{it2} There is a functor, $\Fc_\Bc\colon \cob{n} \to \HFb$, that sends the equivalence class of a cobordism,
$(i,S,j)\colon \Sigma \to \Sigma'$
 to the equivalence class of the HF fibrant span,
 \begin{equation}\label{eq:fib_span_local}
 \vcenter{\spnd{\Bc^\Sigma}{i^*}{\Bc^{\Sigma'}.}{j^*}{\Bc^S}}
 \end{equation}
\end{enumerate}
\end{Lemma}
 
\begin{proof}
If $\Sigma$ is a compact smooth manifold, then $\Sigma$ has a finite triangulation, and, in particular, it can be given the structure of a finite CW-complex, \cite{MunkresDiff}, so, by Lemma \ref{lem:mappingfromCW}, $\Bc^\Sigma$ is HF.  The rest of Item \eqref{it1} follows from the fact that $\CGWH$ is a cartesian closed category.

For the second point, also note that if ${S}$ is a smooth manifold with boundary, then we can find,  again by, for instance, \cite{MunkresDiff},  a triangulation of the pair $({S},\d {S} )$ making $\d {S}$ a subcomplex of  ${S}$, so the inclusion $\iota\colon \d{S} \to {S}$ is a cofibration\footnote{An alternative proof of this fact is in \cite{Torzewska,Torzewska-HomCobs}.}. As a consequence, the induced map, $\iota^*\colon \Bc^{{S}} \to \Bc^{\d {S}}$, is a fibration. To see that the HF span in \eqref{eq:fib_span_local} is fibrant,  note that $\Bc^{\d {S}}\cong \Bc^{\Sigma \sqcup \Sigma'}\cong \Bc^\Sigma \times \Bc^{\Sigma'}$, where we used the fact that $\CGWH$ is cartesian closed in the last step.
Hence $\Fc_\Bc$ sends (equivalence classes of) cobordisms of manifolds to (equivalence classes of) HF fibrant spans.

That $\Fc_\Bc$ preserves the compositions of $\cob{n}$ and $\HFb$ follows again from the fact that $\CGWH$ is cartesian closed, and, in particular, that the contravariant functor, $\Bc^{(-)}\colon \CGWH \to \CGWH$, sends colimits to limits.
Finally units are preserved by definition.
\end{proof}

\subsection{Quinn's finite total homotopy TQFT}\label{sec:QuinnTQFT}
We will work over the complex number field $\C$. In the following, we let $n$ be a non-negative integer,
 $\Bc$ be an arbitrary, but fixed, HF space, and we fix some $s \in \C$.

\begin{Definition}[(The $s$-indexed form of) Quinn's finite total homotopy TQFT]\label{def:quinnTQFT}  Quinn's finite total homotopy TQFT,  abbr. \emph{Quinn's TQFT},
$\FQp{\Bc}{s}\colon \cob{n} \to {\Vect}$, is defined to be the composite of the functors,
$$\Fc_\Bc\colon \cob{n} \to \HFb \quad \textrm{ and } \quad \FRp{s}\colon \HFb \to {\Vect}.   $$
We will write $\FQ{\Bc}$ for $\FQp{\Bc}{0}$. {(This was Quinn's original normalisation.)}

\end{Definition}

For the second of these functors, see Theorem \ref{thm:FR}, and use  Lemma \ref{1dual} to allow its application here.

Taking this apart,  and in more detail, using the notation in \S\ref{sec:conventions_top},
\begin{itemize}[leftmargin=1cm]\item given a closed smooth $n$-manifold, $\Sigma$, then
$\FQp{\Bc}{s}(\Sigma)=\C (\hpi_0(\Bc^\Sigma)); $
\item given an 
$(n+1)$-cobordism, ${(i_1,S,i_2)}\colon \Sigma_1 \to \Sigma_2$,
between the closed smooth $n$-manifolds, $\Sigma_1$ and $\Sigma_2$, we have that the matrix elements of the resulting linear operator are given by the equation below, for continuous functions
 $f_1\colon \Sigma_1 \to \Bc$ and $f_2\colon \Sigma_2 \to \Bc$,
\begin{multline*}
\big \langle \PC_{f_1}(\Bc^{\Sigma_1}) \mid \FQp{\Bc}{s}( [{(i_1,S,i_2  })] )\mid \PC_{f_2}(\Bc^{\Sigma_2})\big \rangle
\\= \chi^\pi\big(\{f_1|(i_1^*,\Bc^S,i_2^*)|f_2\}\big)\, \big(\chi^\pi(\PC_{f_1}(\Bc^{\Sigma_1}))\big) ^s
\big(\chi^\pi(\PC_{f_2}(\Bc^{\Sigma_2}))\big)^{1-s}.
\end{multline*}
We are using the notation of Notation \ref{space-els}, so we have the HF space:
$$ \{f_1|(i_1^*,\Bc^S,i_2^*)|f_2\}
=
\Bigg\{
H\colon S \to \Bc \Big{ |}
\vcenter{
\xymatrix@R=2pt{ &\Bc\\
 \Sigma_1\ar[ur]^{f_1} \ar[dr]_{i_1} && \Sigma_2\ar[ul]_{f_2} \ar[dl]^{i_2} \\
& S\ar[uu]^H }} \textrm{ commutes}\Bigg \}.
$$
\end{itemize}

 
\subsubsection{The monoidality of Quinn's finite total homotopy TQFT}
The functor, $\FQp{\Bc}{s}\colon \cob{n} \to \Vect,$ can be upgraded to be a symmetric monoidal functor, and hence defines a fully-fledged TQFT.
The following elementary lemma will be used in the proof of this fact. We use the notation in \S\ref{moncob}.

\begin{Lemma}\label{definitionT}
Let $\Bc$ be a homotopy finite space. Let $n$ be a non-negative integer.
\begin{enumerate}[leftmargin=1cm]
\item There is a symmetric monoidal functor, $T_\Bc\colon \Diff{n} \to \Vect$, such that
\begin{itemize}[leftmargin=.5cm]
 \item $T_\Bc(\Sigma)=\FQp{\Bc}{s}(\Sigma)=\C\big(\hpi_0(\Bc^\Sigma)\big)$,
 \item[\quad and]
 \item  given $f\colon \Sigma \to \Bc$ and $f'\colon \Sigma' \to \Bc$, and a diffeomorphism, $\phi\colon \Sigma \to \Sigma'$, then the matrix elements satisfy $$
 \big\langle \PC_f(\Bc^\Sigma)  \,|\, T_\Bc(\phi\colon \Sigma \to \Sigma') \, |\, \PC_{f'}(\Bc^{\Sigma'}) \big \rangle =
 \begin{cases}
1, \textrm{ if } \PC_f(\Bc^\Sigma)= \PC_{f'\circ \phi}(\Bc^{\Sigma}),\\
0, \textrm{otherwise}.
 \end{cases}$$
 \end{itemize}
\item If $\phi\colon \Sigma \to \Sigma'$ is a diffeomorphism, then
$\FQp{\Bc}{s}(\I'(\phi))=T_\Bc(\phi)$.
\end{enumerate}
\end{Lemma}

\begin{proof}
The existence of $T_\Bc$, and that it can be upgraded to be a symmetric monoidal functor, follows from standard results from algebraic topology. The crucial point is that the functor $\pi_0\colon \CGWH \to \Sets$ is a symmetric monoidal functor, since, given CGWH spaces $X$ and $Y$, we have a natural bijection, $\eta''_{X,Y}\colon\pi_0(X)\times \pi_0(Y) \to \pi_0(X\times Y)$, such that, for $x\in X$ and $y \in Y$,
 $(\PC_x(X),\PC_y(Y)) \mapsto \PC_{(x,y)}(X \times Y)$.  Furthermore, the natural isomorphism, $\eta''\colon \Times \circ (\pi_0\times \pi_0) \To \pi_0\circ \Times$, is associative,  i.e.  the diagrams below commutes,  given spaces $X,Y$ and $Z$,
 $$
\xymatrix@R=18pt
{
\big(\pi_0(X) \times \pi_0(Y) \big)\times \pi_0(Z) \ar[d]_{\eta''_{X,Y}\otimes \pi_0(Z) }\ar[rrr]^{\alpha^\Sets_{\pi_0(X),\pi_0(Y),\pi_0(Z)}} &&& \pi_0(X) \times \big(\pi_0(Y)\times \pi_0(Z)\big) \ar[d]^{\pi_0 (X) \otimes \eta''_{Y,Z}}\\
 \Big(\pi_0 (X\times  Y)\Big) \otimes \pi_0 (Z) \ar[d]_{\eta''_{X\times Y,Z}}  &&& \pi_0 (X)\otimes \Big( \pi_0(Y\times Z) \Big) \ar[d]^{ \eta''_{X, Y\times Z}}\\
 \pi_0 \Big( (X\times  Y)\times Z) \Big)
\ar[rrr]_{\pi_0 \big(\alpha^{\CGWH}_{X,Y,Z} \big)}
&&& \pi_0 \Big (X\times (Y\times Z) \Big).
}
$$
If we consider the bijection, $\epsilon\colon \{*\} \to \pi_0(\{*\})$,  we can see that, given a CGHW space $X$, the two diagrams pertaining to the unitality of $\eta''$,  commute, for instance,
$$\xymatrix@R=18pt{
\pi_0(X)\times \{*\}\ar[d]_{\pi_0(X)\times \epsilon} \ar[rr]^{\rho_{\pi_0(X)}^\Sets}\ar[d] &&\pi_0(X)\\  \pi_0(X)\times \pi_0(\{*\}) \ar[rr]_{\eta''_{X,\{*\}}}  &&\pi_0(X\times \{*\})
\ar[u]_{\pi_0(\rho_{X}^\CGWH)},}
$$
 so the triple $(T_\Bc,\eta'',\epsilon)$ is a (strong) monoidal functor.

 Given that the diagram below clearly commutes, $$\xymatrix@R=18pt{\pi_0(X)\times \pi_0(Y) \ar[d]_{\eta''_{X,Y}}\ar[rr]^{\tau^\Sets_{\pi_0(X),\pi_0(Y)}} && \pi_0(Y)\times \pi_0(X)\ar[d]^{\eta''_{Y,X}}\\
\pi_0(X\times Y) \ar[rr]_{\pi_0(\tau^\CGWH_{X,Y})} && \pi_0(Y\times X),
}
$$
(where $\tau$ denotes the natural transformations obtained by exchange of coordinates),  $(T_\Bc,\eta'',\epsilon)$ is also a symmetric monoidal functor.

 The free vector space functor, $\Lin\colon \Sets \to \Vect,$ is  symmetric monoidal. By using  that $\CGWH$ is a monoidal closed category, the functor  $\Fc_\Bc^0 \colon (\Diff{n})^{\mathrm{op}}\to \HFbiso$ (for notation see Lemma \ref{definitionT}) is also symmetric monoidal, in a natural way.
 Now note that $T_\Bc$ is given by the following composition of functors:
$$\Diff{n}\ra{(-)^{-1}} (\Diff{n})^{\mathrm{op}} \ra{\Fc_\Bc^0} \HFbiso \ra{\mathrm{inc}} \CGWH \ra{\pi_0} \Sets \ra{\Lin}\Vect. $$

The second point of the lemma follows   from the second part of Lemma \ref{ids}.
\end{proof}

\begin{Remark}
Note that a part of the monoidal structure of $T_\Bc\colon \Diff{n} \to \Vect$, is a natural isomorphism,
$\Otimes \circ (T_\Bc \times T_\Bc) \To T_\Bc\circ \Sqcup$, which, by tracking the sequence of compositions above,  explicitly is such that, given closed smooth $n$-manifolds, $\Sigma$ and $\Sigma'$, and maps, $f\colon \Sigma \to \Bc$ and $f'\colon \Sigma' \to \Bc$, we have,
$$\PC_f(\Bc^\Sigma) \otimes \PC_{f'}(\Bc^{\Sigma'})\mapsto \PC_{\langle f,f'\rangle}\big (\Bc^{\Sigma \sqcup \Sigma'}\big). $$
\end{Remark}

\begin{Theorem}\label{Thm:Quinn-sym-mon} The functor $\FQp{\Bc}{s}\colon \cob{n} \to \Vect,$ can be upgraded to be a symmetric monoidal functor, which we also denote  $\FQp{\Bc}{s}\colon \cob{n} \to \Vect$.
\end{Theorem}
The symmetric monoidal structure of $\FQp{\Bc}{s}$ is clarified in the proof.
\begin{proof}
We first consider the natural isomorphism, $$\eta'\colon \Otimes \circ \big (\FQp{\Bc}{s}\times \FQp{\Bc}{s}\big) \To \FQp{\Bc}{s}\circ \Sqcup, $$ of functors from $\cob{n}\times \cob{n}$ to $\Vect$ defined as the composite,
\begin{multline*}
  \Otimes \circ \big (\FQp{\Bc}{s}\times \FQp{\Bc}{s}\big) =  \Otimes \circ \big (\FRp{s}\circ \Fc_\Bc\times \FRp{s}\circ \Fc_\Bc)\\\Ra{{\eta \circ (\Fc_\Bc \times \Fc_{\Bc}) }}\FRp{s}\circ \Times \circ (\Fc_\Bc \times \Fc_{\Bc})\stackrel{\cong } {\implies} \FQp{\Bc}{s}\circ\Fc_\Bc\circ \Sqcup.
\end{multline*} 
Here $\eta$ is  defined in Lemma \ref{lem:mon1}, and the last natural isomorphism follows from the fact that $\CGWH$ is cartesian closed.
Explicitly, given compact smooth $n$-manifolds $\Sigma$ and $\Sigma'$, the respective component of $\eta'$, as below
$$\eta'_{\Sigma, \Sigma'}\colon \FQp{\Bc}{s} (\Sigma)\otimes \FQp{\Bc}{s}(\Sigma') \to \FQp{\Bc}{s} (\Sigma\sqcup \Sigma')$$ 
is such that, given $f\colon \Sigma \to \Bc$ and $f'\colon \Sigma' \to \Bc$,
$$ \PC_f(\Bc^\Sigma) \otimes \PC_{f'}(\Bc^{\Sigma'})\mapsto \PC_{\langle f,f'\rangle}\big (\Bc^{\Sigma \sqcup \Sigma'}\big). $$

This natural isomorphism can easily be proved to be `associative', meaning that, given closed smooth $n$-manifolds $\Sigma, \Sigma',\Sigma''$, the following diagram commutes (where we omitted the labels in the associativity constraints in $\Vect$),
$$
\xymatrix@R=18pt@C=28pt
{
\Big(\FQp{\Bc}{s} (\Sigma)\otimes \FQp{\Bc}{s}(\Sigma')\Big) \otimes \FQp{\Bc}{s} (\Sigma'')\ar[d]_{\eta'_{\Sigma, \Sigma'}\otimes \FQp{\Bc}{s} (\Sigma'') }\ar[r]^{\alpha^{\Vect}} & \FQp{\Bc}{s} (\Sigma)\otimes \Big( \FQp{\Bc}{s}(\Sigma')  \otimes \FQp{\Bc}{s} (\Sigma'') \Big) \ar[d]^{\FQp{\Bc}{s} (\Sigma) \otimes \eta'_{\Sigma',\Sigma''}}\\
 \Big(\FQp{\Bc}{s} (\Sigma\sqcup  \Sigma')\Big) \otimes \FQp{\Bc}{s} (\Sigma'') \ar[d]_{\eta'_{\Sigma\sqcup \Sigma',\Sigma''}}  & \FQp{\Bc}{s} (\Sigma)\otimes \Big( \FQp{\Bc}{s}(\Sigma'\sqcup \Sigma'') \Big) \ar[d]^{ \eta'_{\Sigma, \Sigma'\sqcup \Sigma''}}\\
 \FQp{\Bc}{s} \Big( (\Sigma\sqcup  \Sigma')\sqcup \Sigma'') \Big)
\ar[r]_{\FQp{\Bc}{s} \big(\alpha^{\cob{n}}_{\Sigma,\Sigma',\Sigma''} \big)}
& \FQp{\Bc}{s} \Big (\Sigma\sqcup (\Sigma'\sqcup \Sigma'') \Big).
}
$$

That this diagram commutes, follows from the fact that, using the definition of the monoidal structure of $\cob{n}$ sketched in \S\ref{moncob}, we have
\begin{align*}
 \FQp{\Bc}{s} \big(\alpha^{\cob{n}}_{\Sigma,\Sigma',\Sigma''} \big)&=
 \FQp{\Bc}{s} \Big( \I'\big(\alpha^{\Diff{n}}_{\Sigma,\Sigma',\Sigma''}\big) \Big)\\
 &=T_\Bc \big(\alpha^{\Diff{n}}_{\Sigma,\Sigma',\Sigma''}\big).
\end{align*}
Here we have used the second point of Lemma \ref{definitionT} in the last step. We have that the diagram above commutes, since the functor $T_\Bc$ is monoidal, by the first point of Lemma  \ref{definitionT}. Note that $T_\Bc$ and $\FQp{\Bc}{s}$ coincide on objects.

The remaining bits of the proof of the fact that $\FQp{\Bc}{s}\colon \cob{n} \to \Vect$ can be turned into a symmetric monoidal functor follow exactly the same pattern.  
\end{proof}
\begin{Remark}
 We have monoidal natural isomorphisms connecting all normalisations of Quinn's finite total homotopy TQFT, obtained by applying Proposition \ref{nat:trans_st}. The details are left to the reader.
\end{Remark}

\begin{Remark}
The ground field for Quinn's finite total homotopy TQFT can be taken to be $\Q$, for $s=1$ or $s=0$, or the Galois closure of $\mathbb{Q}$, for $s=1/2$. The case $s=1/2$ coincides with the conventions in \cite{Yetter} where Yetter develops a (2+1)-dimensional TQFT, derived from a finite crossed module. Passing to its classifying space gives a particular case of Quinn TQFT; some details are in  \cite{Martins_Porter}.
\end{Remark}

\begin{Remark}[Cohomology twisting]\label{coho-twist-Quinn}\label{rem:coho_twist}
 We note that if we restrict to oriented $n$-manifolds and oriented cobordisms Quinn's finite total homotopy TQFT, $\FQp{\Bc}{s}$, can be `twisted' by a cohomology class in $H^n(\Bc,U(1))$. The details are left to the reader. If $\Bc$ is the classifying space of a finite group, this gives exactly Dijkgraaf-Witten TQFT \cite{DW}. Some details of the construction, in the closed case, when $\Bc$ is the classifying space of a  finite crossed module, can be found in \cite{Martins_Porter}.
\end{Remark}

\subsection{{Some  elementary properties of Quinn's finite total homotopy TQFT}}
Methods  for concrete calculations of Quinn's TQFT will be addressed in Chapter \ref{Quinn Calc}, when $\Bc$ is the classifying space of a homotopy finite  crossed complex. There we will also discuss  the  methods of calculation in the extended version.

For the moment, we will restrict our attention to some simple examples and observations. Fix a non-negative integer $n$ and $s \in \C$.

We recall the idea of direct sum of TQFTs, as defined initially in \cite{Durhuus-Jonsson:classification:1994} in the case of (2+1)-dimensional TQFTs, but later applied, in \cite{sawin:direct:1995}, to general TQFTs.

Suppose that $\mathcal{Z}_1$ and $\mathcal{Z}_2$ are two TQFTS.
\begin{Definition}
The \emph{direct sum}, $\mathcal{Z}_1\oplus \mathcal{Z}_2$, of $\mathcal{Z}_1$ and $\mathcal{Z}_2$ is the theory which:
\begin{itemize}[leftmargin=1cm]
\item associates, to each connected $\Sigma$, the vector space $\mathcal{Z}_1(\Sigma)\oplus \mathcal{Z}_2(\Sigma)$;
\item associates, to each disconnected $\Sigma$, the tensor product of the vector spaces associated to its components;
\item associates, to each connected cobordism, $(i,S,j)\colon \Sigma_1 \to \Sigma_2$, the linear map, $\mathcal{Z}_1(i,S,j)\oplus \mathcal{Z}_2(i,S,j)$, interpreted as an operator on the appropriate vector spaces,
\item[]\hspace{-1.2cm} and
\item associates to each disconnected cobordism, the tensor product of the values on the components.
\end{itemize}
\end{Definition}
\begin{Definition}
The \emph{tensor product}, $\mathcal{Z}_1\otimes \mathcal{Z}_2$, of TQFTs $\mathcal{Z}_1,\mathcal{Z}_2\colon \cob{n} \to \Vect$, is the theory obtained as the composite:
$$\cob{n}\ra{\langle \mathcal{Z}_1 , \mathcal{Z}_2 \rangle } \Vect \times \Vect \ra{\Otimes_{\Vect }} \Vect. $$
The \emph{trivial TQFT} assigns $\C$ to all $n$-dimensional manifolds and the identity map to all cobordisms.
\end{Definition}
The following follows from standard arguments, that are left to the reader.
\begin{Theorem} The following holds for Quinn's TQFT:
\begin{enumerate}[leftmargin=1cm]
 \item  If $\Bc=\{\ast\}$ then $\FQp{\Bc}{s}\colon \cob{n} \to \Vect$ is isomorphic to the trival TQFT;
 \item and given  homotopy finite spaces $\Bc$ and $\Bc'$ we have:
   $$\FQp{\Bc\sqcup \Bc'}{s}\cong \FQp{\Bc}{s} \oplus\FQp{\Bc'}{s} \qquad \textrm{ and } \qquad
   \FQp{\Bc\times \Bc'}{s}\cong \FQp{\Bc}{s} \otimes\FQp{\Bc'}{s}.$$
\end{enumerate}
\end{Theorem}

\subsubsection{Changing $\Bc$}\label{Changing B}
Given  the previous theorem,  one might think that there was some possible functoriality of  $\FQp{\Bc}{s}$ with $\Bc$ itself, but recall, for instance from \cite[\S 2.5 and Appendix A2]{Lec_TQFT}, that if $\varphi\colon \mathcal{Z}_1\To \mathcal{Z}_2$ is a (monoidal) natural transformation between TQFTs, then it is a natural isomorphism.  If $f\colon \Bc_1\to \Bc_2$ is a general continuous map, we therefore should not expect that  there would be some sort of induced `morphism' between $\FQp{\Bc_1}{s}$ and $\FQp{\Bc_2}{s}$.

Under certain circumstances, however, a map, $f\colon \Bc_1\to \Bc_2$, between HF spaces, does induce a natural transformation between $\FQp{\Bc_1}{s}$ and $\FQp{\Bc_2}{s}$, which as we noted must, then, be a natural isomorphism.

\begin{Theorem}\label{Th:changeB}
If  $f\colon\Bc_1\to \Bc_2$ is a homotopy equivalence, then $f$ induces a monoidal natural isomorphism,  $f_*\colon\FQp{\Bc_1}{s}\To \FQp{\Bc_2}{s}$, between  $\FQp{\Bc_1}{s}$ and $\FQp{\Bc_2}{s}$.
This natural isomorphism is defined in the following way:

If $\Sigma$ is a closed $n$-manifold, then the linear map, ${(f_*)}_\Sigma\colon \FQp{\Bc_1}{s}(\Sigma)\to \FQp{\Bc_2}{s}(\Sigma)$, is such that, given $g\colon \Sigma \to \Bc_1$, then
$${(f_*)}_\Sigma\big ( \PC_g\big( \Bc_1^\Sigma \big)\big)= \PC_{f \circ g}( \Bc_2^\Sigma ).$$
\end{Theorem}
\begin{proof}
We  always  have this family of mappings, $(f_*)_{\Sigma}\colon \FQp{\Bc_1}{s}(\Sigma)\to \FQp{\Bc_2}{s}(\Sigma)$, induced by post-composition with $f$, but, {in general,} this need not define a natural transformation, due to possible incompatibility with the cobordisms.

Suppose $(i,S,j)\colon \Sigma_1\to \Sigma_2$ is a cobordism, thus giving us  fibrant spans,  $$\xymatrix@R=-5pt{ &&\Bc_i^S\ar[dl]_{i^*}\ar[dr]^{j^*}\\  & \Bc_i^{\Sigma_1} & &\Bc_i^{\Sigma_2} , }$$ for $i=1,2$. The function $f\colon \Bc_1 \to \Bc_2$ gives us  a commutative diagram, $$\xymatrix@R=18pt{\Bc_1^{\Sigma_1}\ar[d]_{f^{\Sigma_1}}&\Bc_1^{S}\ar[l]_{i^*}\ar[d]^{f^S}\ar[r]^{j^*}&\Bc_1^{\Sigma_2}\ar[d]^{f^{\Sigma_2}}\\
\Bc_2^{\Sigma_1}&\Bc_2^{S}\ar[l]^{i^*}\ar[r]_{j^*}&\Bc_2^{\Sigma_2},
}
$$
in which the vertical maps are all homotopy equivalences. (We note that this is \emph{not} a morphism in the category  $\HFb$, but does relate to a higher category structure  on the class of HF spans.)
The commutative diagram induces a map of fibrations,
\begin{equation}\label{Eq_diag_fibrations}
\vcenter{\xymatrix@R=18pt{
&\Bc_1^{S} \ar[rr]^>>>>>>>>>>{\langle i^*, j^*\rangle}\ar[d]_{f^{S}} && \Bc_1^{\Sigma_1}\times  \Bc_1^{\Sigma_2}\ar[d]^{f^{\Sigma_1}\times f^{\Sigma_2}} \\
&\Bc_2^{S} \ar[rr]_>>>>>>>>>>{\langle i^*, j^*\rangle} && \Bc_2^{\Sigma_1}\times  \Bc_2^{\Sigma_2},
}}
\end{equation}
where the vertical arrows are homotopy equivalences.

We also have a diagram of vector spaces and linear maps,
\begin{equation}\label{diag:natQf}
\vcenter{\xymatrix@R=18pt@C=30pt{\FQp{\Bc_1}{s}(\Sigma_1)\ar[r]^{ (f_*)_{\Sigma_1}}\ar[d]_{\FQp{\Bc_1}{s}(S)}&\FQp{\Bc_2}{s}(\Sigma_1)\ar[d]^{\FQp{\Bc_2}{s}(S)}\\
\FQp{\Bc_1}{s}(\Sigma_2)\ar[r]_{ (f_*)_{\Sigma_2}}&\FQp{\Bc_2}{s}(\Sigma_2).
}}
\end{equation}
To check that this diagram commutes, in this context, we pick basis elements in $\FQp{\Bc_1}{s}(\Sigma_1)$ and $\FQp{\Bc_1}{s}(\Sigma_2)$, and compare  the matrices corresponding to the left-hand side, with those on the right-hand side, with respect to the image basis.

First we note that $\hpi_0(\Bc_1^\Sigma)$ and $\hpi_0(\Bc_2^\Sigma)$ are related by the bijection, $(f_*)_\Sigma$,  induced from $f$. Consider arbitrary maps $g\colon \Sigma_1 \to \Bc_1$ and $g'\colon \Sigma_2 \to \Bc_1$. To prove that the diagram in \eqref{diag:natQf} commutes, it suffices to prove that
\begin{multline*}
\big\langle \PC_g(\Bc_1^{\Sigma_1}) \mid \FQp{\Bc_1}{s}( [(i,S,j)] )\mid \PC_{g'}(\Bc_1^{\Sigma_2})\big \rangle\\
= \big\langle \PC_{f\circ g}(\Bc_2^{\Sigma_1}) \mid \FQp{\Bc_2}{s}( [(i,S,j)] )\mid \PC_{f\circ g'}(\Bc_2^{\Sigma_2})\big \rangle.
\end{multline*}
 
Unpacking the notation, this amounts to comparing the corresponding fibres of  the horizontal fibrations in diagram in \eqref{Eq_diag_fibrations}. These fibres are homotopy equivalent by Corollary \ref{cor:h.e. of spans}, applied to the map of fibrations {given} in \eqref{Eq_diag_fibrations}.

We note that these isomorphisms respect the monoidal structure and also the composition, which completes the proof.
\end{proof}

\begin{Corollary} If $\Bc$ is contractible then the TQFT $\FQp{\Bc}{s}$ is trivial.
\end{Corollary}

We note that, from the proof of Theorem \ref{Th:changeB}, how the natural isomorphism, $f_*\colon \FQp{\Bc_1}{s} \to \FQp{\Bc_2}{s} $, depends only on the homotopy class of the homotopy equivalence, $f\colon \Bc_1 \to \Bc_2$.   This gives:
\begin{Theorem}\label{Thm:EB}
Given any HF-space, $\Bc$, there is an action of the group, $\mathcal{E}(\Bc)$, of homotopy classes of self homotopy equivalences of $\Bc$ on the TQFT $\FQp{\Bc}{s}$, by natural isomorphisms.
\end{Theorem}
\noindent N.B: The groups $\mathcal{E}(\Bc)$ are in general non trivial. E.g., if $B_G$ is the classifying space of a group $G$, then $\mathcal{E}(B_G)$ is isomorphic to the group of outer automorphism of $G$.

\chapter{Once-extended versions of Quinn's finite total homotopy TQFT}\label{part:once_ext}

Chapter \ref{part:once_ext} of this paper consists of two sections. The first, Section \ref{sec:homotopy_underpinning}, looks at the homotopy-theoretical and bicategorical underpinning of the  once-extended  Quinn TQFT. The second, Section \ref{sec:def-once-extended}, gives the detailed construction of that extended TQFT.
Throughout Chapter \ref{part:once_ext}, we will  work with an arbitrary subfield, $\kappa$, of $\C$.

Let $\Bc$ be an arbitrary, but fixed, homotopy finite space, and $n$ be a non-negative integer. In this chapter, we will see how the $s=0$ case  of Quinn's finite total homotopy TQFT  (abbr. `Quinn TQFT'), formulated  in Definition \ref{def:quinnTQFT} and Theorem \ref{Thm:Quinn-sym-mon} as a symmetric monoidal functor,  $\FQ{\Bc}\colon \cob{n} \to \Vect$,
 can be `categorified' to what we called the \emph{once-extended Quinn TQFT}. We will formulate this categorification  as a (symmetric monoidal) bifunctor, appearing in  Definition \ref{def:once_ext_Quinn}, and denoted $$\tFQ{\Bc}\colon \tcob{n} \to \vProfGrphf .$$

Here $\tcob{n}$ is the bicategory with objects the closed (and, by convention, smooth) $n$-manifolds, the 1-morphisms being the $(n+1)$-cobordisms between closed $n$-manifolds,  and the 2-morphisms the equivalence classes of extended $(n+2)$-cobordisms between $(n+1)$-cobordisms; see \cite{Mortoncospans,morton:cohomological:2015,Schommer-Pries}. Our convention are explained in Subsection \ref{sec:sym_mon}. On the other hand, the objects of the bicategory $\vProfGrphf$ are the homotopy finite groupoids, as defined in \S\ref{sec:convention_groupoids},   and given groupoids $G$ and $H$ the 1-morphisms, $G\bto H$, are $\Vect$-valued profunctors between groupoids, in other words functors $G^\op \times H\to \Vect$. The 2-morphisms of $\vProfGrphf$ are natural transformations of functors. These constructions are reviewed in Subsection \ref{sec:prof_conventions}.

Let $\Sigma$ be a closed and smooth $n$-manifold. Typically, the groupoid $\tFQ{\Bc}(\Sigma)$, despite being homotopy finite, is uncountable. In order to reduced the size  of the target groupoids, we will consider a bicategory, $\tdcob{n}{\Bc}$, with objects $\Bc$-decorated  $n$-manifolds, which are closed smooth $n$-manifolds $\Sigma$ equipped with a finite subset of the function space $\TOP(\Sigma,\Bc)$, containing at least one function for each homotopy class of maps $\Sigma \to \Bc$.  The rest of the bicategory structure of $\tdcob{n}{\Bc}$ is induced by that of $\tcob{n}$, in the obvious way, as discussed in Subsection \ref{sec:fin_ext}.

We will then consider another once-extended TQFT, called the \emph{finitary once-extended Quinn TQFT},
$$\tFQd{\Bc}\colon \tdcob{n}{\Bc} \to \vProfGrpfin;$$
 see Definition \ref{def:finitary_ext_TQFT}.
The bicategory, $\vProfGrpfin$,  is the full sub-bicategory of $\vProfGrphf$, with objects the finite groupoids; see   \S\ref{sum:conv_prof}.
This, in turn, gives rise to another once-extended TQFT, called the \emph{Morita-valued once-extended Quinn TQFT}, in Definition \ref{Sec:ext:Mor}, denoted,
$$\tFQmor{\Bc}\colon \tdcob{n}{\Bc} \to \Mor,$$
 where $\Mor$ is the bicategory of $\kappa$-algebras, bimodules and bimodule maps.

The algebraic construction showing how to go from  the bicategory $\vProfGrpfin$
to the bicategory  $\Mor$, starting from \emph{groupoid algebras}, may be of independent interest. This is laid out in
Subsection \ref{sec:Mor_ext}.

Depending on which setting is chosen, the groupoids, or algebras that $\tFQd{\Bc}$  and  $\tFQmor{\Bc}$ assign to a closed manifold, $\Sigma$, with a $\Bc$-decoration,  explicitly depend on the $\Bc$-decoration of $\Sigma$. However this dependence is up to a canonically defined, and invertible, profunctor or bimodule, which is functorial with respect to further changes in the decoration (up to natural isomorphism), and natural with respect to the profunctors, or bimodules, assigned to cobordisms. This is discussed in Subsection \ref{sec:fin_ext} and \S\ref{sec:Mor_valued-eTQFT}

We will show  explicit examples of calculations of these once-extended TQFTs later on in Chapter \ref{Quinn Calc}, Section \ref{sec:TQFTS_xcomp}, for the case in which $\Bc$ is the classifying space of a finite crossed complex. This includes the case of classifying spaces of finite 2-groups, as appearing in higher gauge theory, see e.g. \cite{Baez-Schreiber:2004,Baez_Huerta:2011,Companion,Martins_Picken:2011}.

\textbf{Remark:} It is possible to define a one-parameter categorification of
Quinn TQFT, 
 but we will not deal with that here. This would considerably increase the complexity of our formulae, without adding much  generality to our construction.

\section{The homotopy  underpinning of the  once-extended  Quinn TQFT}\label{sec:homotopy_underpinning}
\sectionmark{Homotopy underpinning of the once-extended  Quinn TQFT}

Similarly to  our exposition of Quinn's finite total homotopy TQFT, we will factor its once-extended version through a homotopy-theoretical bicategorical object, later denoted $\mathbf{2span}(HF)$, whose objects are HF spaces, 1-morphisms are fibrant HF spans, and  2-morphisms consist of {HF fibrant resolved 2-spans}.  This construction is the main topic of this section. A summary of the construction, with all terms explained, can be found in Subsection \ref{sec:comm_summary}, at the end of this section.

\subsection{Conventions on bicategories}\label{sec:bicategories}
We will frequently need to use the terminology of the theory of bicategories, also called weak 2-categories.
\subsubsection{The basics of bicategories}
A basic introduction to this theory can be found in Leinster's \cite{leinster:bicat:1998}, with a more thorough and complete description given in Borceux's \cite{borceux1}. We may also use the summary to be found in \cite{Schommer-Pries}, and the relevant parts of the draft book, \cite{johnson-yau:bimonoidal:2021}, by Johnson and Yau.

\begin{Definition}
A \emph{bicategory}, $\mathcal{B}$, is specified by the following:
\begin{itemize}[leftmargin=1cm]\item a collection of \emph{objects}, denoted $Ob(\mathcal{B})$, or sometimes $\mathcal{B}_0$;

\item for each pair of objects, $a,b$ in $\mathcal{B}$, a locally small category, $\mathcal{B}(a,b)$, whose objects are \emph{1-morphisms} from $a$ to $b$, whose morphisms are called \emph{2-morphisms} and whose composition is sometimes referred to as \emph{vertical composition};
\item for objects, $a,b,c$ in $\mathcal{B}$, there are composition \emph{functors},
$$\mathbf{c}_{a,b,c}:\mathcal{B}(a,b)\times \mathcal{B}(b,c)\to \mathcal{B}(a,c),$$
and, for each object $a$ in $\mathcal{B}$, a functor, $I_a:[0]\to \mathcal{B}(a,a),$ (where $[0]$ is the `singleton' category).
The functors, $\mathbf{c}$, are called \emph{horizontal compositions};
\item[]\hspace*{-1cm} and
\item natural isomorphisms,

 $\alpha: \mathbf{c}_{a,b,d}\circ (\mathbf{c}_{b,c,d}\times id)\Rightarrow \mathbf{c}_{a,c,d}\circ(id\times \mathbf{c}_{a,b,c}),$

 $\lambda:\mathbf{c}_{a,b,b}\circ(I_b\times id)\Rightarrow id$,\\
 and

 $\rho:  \mathbf{c}_{a,a,b} \circ (id\times I_a)\Rightarrow id,$\\ called, respectively, the \emph{associator} and the \emph{left} and \emph{right unitors}.
 \end{itemize}
 These are required to satisfy the pentagon and triangle identities, which we omit here, referring the reader to Borceux, \cite{borceux1}, and the many other existing sources.
\end{Definition}

\begin{Notation}
When discussing specific bicategories,  we will more often than not use generic composition symbols such as $\bullet$ or $\circ$, but  when it is clear whether we intend horizontal or vertical composition, it can be useful to have available some specific notation that distinguishes them. In such cases, we may use $\#_0$ for horizontal composition (including whiskering of 2-morphisms by 1-morphisms) and $\#_1$ for the vertical composition of 2-morphisms.
\end{Notation}

\begin{Definition}\label{def:bifunctor}
Let $\mathcal{A}$ and $\mathcal{B}$ be bicategories.  A \emph{bifunctor},  also called a \emph{homomorphism}, or a \emph{pseudo-functor},
$F=(F, \varphi)\colon \mathcal{A}\to \mathcal{B}$, consists of
\begin{enumerate}
\item a function, $F\colon \mathcal{A}_0\to \mathcal{B}_0$, mapping objects to objects;
\item for each pair of objects, $a, a'$ in $\mathcal{A}$, a functor,
$$F_{a, a'}\colon \mathcal{A}(a, a')\to \mathcal{B}(F(a), F(a'));$$
\item\label{natural-iso} natural isomorphisms, $\varphi_{a_0,a_1,a_2}$, for each triple, $ a_0, a_1,a_2$, of objects in $\mathcal{A}$,  as shown in the diagrams,
$$\xymatrix{ \mathcal{A}( a_0, a_1)\times\mathcal{A}( a_1, a_2)\ar[r]^{\mathbf{c}^\mathcal{A}}\ar[d]_{F_{a_0, a_1}\times F_{a_1, a_2}}\drtwocell<\omit>{^\hspace*{-1cm} \varphi_{a_0,a_1,a_2}}
&\mathcal{A}( a_0, a_2)\ar[d]^{F_{a_0, a_2}}\\
\mathcal{B}( F(a_0), F(a_1))\times\mathcal{B}( F(a_1), F(a_2))\ar[r]_{\hspace{15mm}\mathbf{c}^\mathcal{B}}&\mathcal{B}(F( a_0), F(a_2))
}$$
and, for each object, $a$ in $\mathcal{A}$,
$$\xymatrix{[0]\ar[r]^{I^\mathcal{A}_a}\ar[dr]_{I^\mathcal{B}_{F(a)}}\druppertwocell<\omit>{^<-2>{{\varphi_a\,\,\,}}}&\mathcal{A}(a,a)\ar[d]^{F_{a,a}}\\
&\mathcal{B}(F(a),F(a)),
}$$
such that certain diagrams, expressing compatibility with the corresponding associators and unitors, commute, and, again, we will not give them here as they can easily be found in the literature and we will give them in a simplified case slightly later.

If  $\varphi_{a_0,a_1,a_2}$ and $\varphi_a$ are all identities,  $F$ is said to be a \emph{strict homomorphism}.
\end{enumerate}\end{Definition}

\textbf{Notation:} Given a 1-morphism, $b_1\xrightarrow{h}b_2$, in a bicategory $\mathcal{B}$, we denote by $h_*$ and $h^*$, the natural transformations / induced morphisms,
$$h_*\colon\mathcal{B}(b,b_1)\to \mathcal{B}(b,b_2),$$
and $$h^*\colon\mathcal{B}(b_2,b)\to \mathcal{B}(b_1,b).$$

\begin{Definition}
Let $(F,\varphi), (G,\psi)\colon\mathcal{A}\to \mathcal{B}$ be two homomorphisms between bicategories. A \emph{transformation}, also called a \emph{pseudo-natural} transformation, $\sigma\colon F\Rightarrow G$, is given by
\begin{itemize}[leftmargin=1cm]\item 1-morphisms, $\sigma_a\colon F(a)\to G(a)$, for each object $a$ in $\mathcal{A}$;
\item  given objects $a$ and $a'$ of $\mathcal{A}$, natural isomorphism,  $\sigma_{a,a'}$, as in the diagram,
$$\xymatrix@C=30pt@R=18pt{\mathcal{A}(a,a')\ar[r]^<<<<<<<<<{F_{a,a'}}\ar[d]_{G_{a,a'}}\drtwocell<\omit>{^\hspace*{-0.3cm} \sigma_{a,a'}}& \mathbf{B}\big(F(a),F(a')\big)\ar[d]^{(\sigma_{a'})_*}&\\
 \mathbf{B}\big(G(a),G(a')\big)\ar[r]_{(\sigma_{a})^*}&\mathbf{B}\big(F(a),G(a')\big).
}$$
\end{itemize}
As before we omit the conditions for compatibility with the other structure, referring to the literature.
\end{Definition}

\begin{Definition}
Given $F,G\colon\mathcal{A}\to \mathcal{B}$, as before, and  $\sigma, \theta\colon F\To G$, two transformations, a \emph{modification}, $\Gamma\colon\sigma\Rrightarrow \theta$, consists of a 2-morphism, $\Gamma_a\colon\sigma_a\Rightarrow \theta_a$, for every object $a$ in $\mathcal{A}$.  These are required to make the following square commute,
$$\xymatrix@R=18pt@C=50pt{G_{a,a'}(f)\#_0\sigma_a\ar@{=>}[r]^{G_{a,a'}(f)\#_0 \Gamma_a}\ar@{=>}[d]_{(\sigma_{a,a'})_f}&G_{a,a'}(f)\#_0\theta_a\ar@{=>}[d]^{(\theta_{a,a'})_f}\\
\sigma_{a'}\#_0 F_{a,a'}(f)\ar@{=>}[r]_{ \Gamma_{a'} \#_0 F_{a,a'}(f)}&\theta_{a'}\#_0F_{a,a'}(f),}$$
for every 1-morphism $f\colon a\to a'$ in $\mathcal{A}$.
\end{Definition}
\begin{Remark}
We refer the reader to \cite{Schommer-Pries} for how to compose transformations, etc., so that one gets a bicategory $Bicat(\mathcal{A},\mathcal{B})$ with the resulting structure, provided the bicategories are small. \end{Remark}

We suppose that $\mathcal{A}$ is a bicategory.
\begin{Definition}
(i)  Given a pair of 1-morphisms, $f\colon A\to B$ and $u\colon B\to A$, in $\mathcal{A}$, we say $f$ is \emph{left adjoint} to $u$ {(and $u$ is \emph{right adjoint} to $f$),} and written $f\dashv u$, if there are two 2-morphisms,
$$\eta \colon1_A\Rightarrow u f, \qquad \textrm{ and } \qquad \varepsilon\colon fu\Rightarrow 1_B,$$
such that the following equations hold:
$$(u\xleftrightarrow{\cong}1_Au\xrightarrow{\eta\cdot u}uf\cdot u\xleftrightarrow{\cong}u\cdot fu\xrightarrow{u\cdot \varepsilon}u\cdot 1_B\xleftrightarrow{\cong}u)= id_u,$$
and
$$(f\xleftrightarrow{\cong}f1_A\xrightarrow{f\eta} f\cdot uf \xleftrightarrow{\cong}fu\cdot f \xrightarrow{\varepsilon f}1_B\cdot f\xleftrightarrow{\cong}f)= id_f.$$We have written $\cong$ to label the evident unitors and associators, or their inverses.

(ii)  A 1-morphism, $f\colon A\to B$, in $\mathcal{A}$ is an \emph{equivalence} if there is a 1-morphism, $g\colon B\to A$, and two 2-isomorphisms, $gf\xRightarrow{\cong}1_A$ and $fg\xRightarrow{\cong}1_B$. We have an \emph{adjoint equivalence} if $f\dashv g$ and both $\eta$ and $\varepsilon$ are  isomorphisms.
\end{Definition}

\subsubsection{{Pseudo-functors from categories to bicategories}}\label{sec:pseud-cat-bicat}
 In Section \ref{sec:def-once-extended}, we will consider several examples of bifunctors / pseudo-functors, $F\colon\mathcal{A}\to \mathcal{B}$, in which
 \begin{itemize}[leftmargin=0.6cm]
 \item the domain, $\mathcal{A}$, is `locally discrete', meaning that each $\mathcal{A}(x,y)$ is a discrete category, i.e., a set.  (We will often just say that $\mathcal{A}$ is a category.)
  \item[] \hspace*{-1.2cm}But
 \item  the codomain, $\mathcal{B}$, is a  bicategory, usually that of $\tcob{n}$, any of the span or cospan bicategories, $\vProf$, one of its variants, or $\Mor$.
 \end{itemize}

We repeat the specification of a pseudo-functor,  $F\colon\mathcal{A}\to \mathcal{B}$,   but  in the simplified form that this context allows.  We have\label{def:special pseudofunctor}
 \begin{itemize}[leftmargin=0.6cm]\item for each object, $a$ in $\mathcal{A}$, an object, $F(a)$, in $\mathcal{B}$;
 \item for each pair, $a_0,a_1$, of objects in $\mathcal{A}$, a functor,
 $$F_{a_0,a_1}\colon\mathcal{A}(a_0,a_1)\to \mathcal{B}(F(a_0),F(a_1)), $$
 and, because $\mathcal{A}(a_0,a_1)$ is discrete, this just means a family of 1-morphisms,
 $F(f)\colon F(a_0)\to F(a_1)$, where $f\colon a_0\to a_1 \in \mathcal{A}$;
 \item for each composable pair of morphisms, $a_0\xrightarrow{f}a_1$, $a_1\xrightarrow{g}a_2$, in $\mathcal{A}$, an invertible 2-morphism, $\varphi_{g,f}\colon F(g)F(f)\Rightarrow F(gf)$ in $\mathcal{B}$;
 \item for each object $a$ of $\mathcal{A}$, an invertible 2-morphism, $\varphi_a\colon id_{F(a)}\Rightarrow F(id_a)$, in $\mathcal{B}$.
 \end{itemize}
\noindent These must satisfy the following conditions:
\begin{enumerate}[leftmargin=1cm]\item compatibility with the associator in $\mathcal{B}$,  so given, in addition, $h\colon a_2\to a_3$, the following diagram commutes:
$$\xymatrix@R=4pt{
(F(h)F(g))F(f)\ar[r]\ar[dd]_{a_\mathcal{B}}&F(hg)F(f)\ar[dr]&\\
&&F((hg)f)=F(h(gf)),\\
F(h)(F(g)F(f))\ar[r]&F(h)F(gf)\ar[ur]&}$$
where the unlabelled arrow are derived from the various $\varphi$ 2-cells,
\item compatibility with the right and left unitors, (which are `equalities' in $\mathcal{A}$), so, for each $f\colon a_0\to a_1$ in $\mathcal{A}$, the diagrams below commute:
$$\vcenter{\xymatrix@R=13pt{F(f)\cdot id_{F(a_0)}\ar[r]\ar[dr]_<<<<<<<{\rho^\mathcal{B}_{F(a_0)}}&F(f)\cdot F(id_{a_0})\ar[d]\\
&F(f\cdot id_{a_0})=F(f),
}} \,\,
\vcenter{\xymatrix@R=13pt{id_{F(a_1)}\cdot F(f) \ar[r]\ar[dr]_<<<<<<<{\lambda^\mathcal{B}_{F(a_1)}}&F(id_{a_0})\cdot F(f)\ar[d]\\
&F( id_{a_1}\cdot f)=F(f),
}}$$where the unlabelled arrows are the evident ones. \end{enumerate}

When formalising or analysing  a structure in a category or bicategory, the structure is often expressed in terms of the commutativity of certain diagrams.  Suppose we have a commutative diagram in a category $\mathcal{A}$. We can think of this as a functor $D\colon \mathcal{I}\to \mathcal{A}$, where $\mathcal{I}$ is some `template' for the commutative diagram. We can then think of $D$ as a (trivially structured) pseudo-functor and compose it with our given $F\colon \mathcal{A}\to \mathcal{B}$.  The result will be a pseudo-functor from $\mathcal{I}$ to $\mathcal{B}$, so a `pseudo-commutative' diagram in $\mathcal{B}$. We note all the 2-cells in this diagram will be invertible. \label{page:pseudofunctors-to-diagrams}

As a more-or-less trivial example, we can take $\mathcal{I}=[2]$, the small category corresponding to the ordered set, $0\to 1\to 2$, and the $D$ will correspond to a commutative diagram of form $$\xymatrix@R=18pt{&a_1\ar[dr]^{a_{12}}&\\a_0\ar[ur]^{a_{01}}\ar[rr]_{a_{02}}&&a_2.}$$
The corresponding 2-diagram, $FD$, in $\mathcal{B}$ will be the pseudo-commutative one having an invertible 2-arrow from $F(a_{12})F(a_{01})$ to $F(a_{02}),$ together with three invertible 2-arrows, $id_{F(a_i)}\Rightarrow F(id_{a_1}),$ for $i=0,1,2$.

\subsubsection{Monoidal bicategories}
We will also need the bicategorical analogues of monoidal and symmetric monoidal categories.

One of the motivating examples for the notion of a symmetric monoidal bicategory is the following, where $R$ is a commutative ring.
\begin{Example}Let $Alg(R)$ be the bicategory such that \label{Alg(R)}
\begin{itemize}[leftmargin=1cm]
\item the objects are $R$-algebras, denoted $\mathcal{A}$, $\mathcal{B}$, etc.;
\item the morphisms from $\mathcal{A}$ to $\mathcal{B}$ are the left-right $(\mathcal{A},\mathcal{B})$-bimodules;
\item the 2-morphisms are the bimodule homomorphisms.
\end{itemize}
(This is all formally treated in \S\ref{sec:mor_defined}.)
The monoidal product is the tensor product over $R$, so the unit is $R$ itself, considered as an $R$-algebra. This bicategory is also often denoted $Alg_2(R)$ or, as later in this paper, by  $\Mor_R$, or simply by $\Mor$. It will be one of the main codomain bicategories for the once-extended Quinn theory.
\end{Example}

A monoidal bicategory is a bicategory that also has a monoidal structure, up to the equivalence inherent in the bicategorical context.
They can be defined in various ways, for instance, as a tricategory having just one object, \cite{GPS,Gurski_book}. Other definitions mention Gray categories, for which see, for example, \cite{Carqueville-Meusberger-Schaumann:2016}.  Each of these is fairly complex to give, and needs a few more definitions.  The following is one of the simpler ones in as much as it seems fairly clearly motivated by the definition of monoidal category suitably weakened with equality replaced by equivalence.  It does use some bicategorical language that we have not given earlier, but is, perhaps, fairly self explanatory\footnote{The structure and laws are well illustrated in A. S. Corner's thesis, \cite{corner:thesis:day:2016} in \S 1.6,  in the draft book by Johnson and Yau, \cite{johnson-yau:bimonoidal:2021} and in Mike Stay's article, \cite{stay:compact-closed-bicat:2016}.}.

\begin{Definition}\label{def-mon-bicat}
A \emph{monoidal bicategory}, $\mathcal{A}$, consists of \begin{itemize}[leftmargin=1cm]
\item a bicategory, $\mathcal{A}$;
\item a pseudofunctor/homomorphism, $\otimes: \mathcal{A}\times \mathcal{A}\to \mathcal{A}$;
\item a pseudofunctor/homomorphism, $I:\mathbf{1}\to \mathcal{A}$, where $\mathbf{1}$ is the unit bicategory;
\item an adjoint equivalence, sometimes called the \emph{monoidal associator}, diagrammatically denoted
$$\xymatrix@R=15pt{\mathcal{A}^3\ar[r]^{\otimes\times\mathcal{A}}\ar[d]_{\mathcal{A}\times\otimes} \drtwocell<\omit>{\alpha}&\mathcal{A}^2\ar[d]^\otimes\\
\mathcal{A}^2\ar[r]_\otimes&\mathcal{A},}$$
in $Bicat(\mathcal{A}^3,\mathcal{A})$, corresponding to associativity of $\otimes$ in a monoidal category. This adjoint equivalence consists of $\alpha$, its adjoint, $\alpha^*$, with unit, $\eta^\alpha$, and counit, $\varepsilon^\alpha$. To see what these do, we take a triple, $A,B,C$, of objects in $\mathcal{A}$, so $(A,B,C)$ is in $\mathcal{A}^3$, and then we have  that
\begin{align*}
\alpha_{CBA}&\colon (C\otimes B)\otimes A\to C\otimes (B\otimes A), \intertext{whilst} \alpha^*_{CBA}& \colon  C\otimes (B\otimes A)\to (C\otimes B)\otimes A,
\end{align*}with unit and counit,
\begin{align*}
\eta^\alpha_{CBA}&\colon Id\Rightarrow \alpha^*_{CBA}\circ\alpha_{CBA}\intertext{ and }
\varepsilon^\alpha_{CBA}&\colon \alpha_{CBA}\circ\alpha^*_{CBA}\Rightarrow Id,
\end{align*} being isomorphisms. Furthermore, given 1-morphisms, $f\colon C \to C', g\colon B \to B'$ and $h\colon A \to A'$, we  have  natural 2-morphisms, for instance,
$$\xymatrixcolsep{2.5cm}\xymatrix@R=13pt{ (C\otimes B)\otimes A \ar[d]_{\alpha_{CBA}} \drtwocell<\omit>{\hspace*{8mm} \alpha_{(f,g,h)}}\ar[r]^{(f\otimes g)\otimes h} & (C' \otimes B')\otimes A' \ar[d]^{\alpha_{C'B'A'}}\\
 C\otimes (B\otimes A) \ar[r]_{f\otimes ( g\otimes h)} & C' \otimes (B'\otimes A');\\
}
$$
\item adjoint equivalences, sometimes called the \emph{monoidal unitors},
$$\vcenter{\xymatrix@R=13pt{&\mathcal{A}^2\ar[dr]^\otimes&\\
\mathcal{A}\ar[ur]^{I\times \mathcal{A}}\ar[rr]_{\mathcal{A}}\rruppertwocell<\omit>{<-2>\ell}&&\mathcal{A}
}}
\qquad \textrm{ and } \qquad \vcenter{\xymatrix@R=13pt{&\mathcal{A}^2\ar[dr]^\otimes&\\
\mathcal{A}\ar[ur]^{\mathcal{A}\times I}\ar[rr]_{\mathcal{A}}\rruppertwocell<\omit>{<-2>r}&&\mathcal{A},
}}$$
in $Bicat(\mathcal{A}^2,\mathcal{A})$, corresponding to left and right unitors,
\item[]\hspace*{-10mm} and
\item an invertible modification giving the analogue of the pentagon axiom for monoidal product in a monoidal category. This is called the \emph{pentagonator}.  Its component 2-morphisms, for objects $A,B,C, D$ in $\mathcal{A}$ looks like,\label{pentagonator}
$$\scriptsize
\xymatrix@R=10pt@C=-5pt{&(D\otimes (C\otimes B))\otimes A\ar[rr]&&D\otimes((C\otimes B)\otimes A)\ar[dr]&\\
((D\otimes C)\otimes B)\otimes A  \xtwocell[rrrr]{}<>{\qquad\,\,\,\,\pi_{DCBA}}  \ar[ur]\ar[drr] &&&&D\otimes(C\otimes(B\otimes A)),\\
&&(D\otimes C)\otimes (B\otimes A)\ar[urr]&}$$
where each of the unlabelled arrows corresponds to a use of an associator, possibly combined with an identity on an object,  as in the usual pentagon rule for monoidal categories;
\item[]\hspace{-1cm} and
\item invertible modifications, $\utilde{\mu}$, $\utilde{\lambda}$ and $\utilde{\rho}$, called the \emph{middle}, \emph{left} and \emph{right 2-unitors}, respectively, with component 2-morphisms, for objects $A$,  $B$ in $\mathcal{A}$,
$$\xymatrix@R=15pt{(B\otimes I)\otimes A\ar[r]\drtwocell<\omit>{\quad\,\,\mu_{B,A}}&B\otimes(I\otimes A)\ar[d]^{B\otimes \ell_A}\\
B\otimes A\ar[u]^{r_B^*\otimes A}\ar[r]_=&B\otimes A,}$$
$$\xymatrixcolsep{2.5mm}\xymatrix@R=15pt{(I\otimes B)\otimes A\ar[rr]^{\ell_B\otimes A}\ar[dr]_\alpha&\ar@{}[d]|{\Downarrow\utilde{\lambda}_BA}&B\otimes A,\\
&I\otimes (B\otimes A)\ar[ur]_{\ell _{B\otimes A}}&
} \qquad
\xymatrixcolsep{2.5mm}\xymatrix@R=15pt{B\otimes A\ar[rr]^{B\otimes r^*_A}\ar[dr]_{r^*_{B\otimes A}}&\ar@{}[d]|{\Downarrow\utilde{\rho}_BA}&B\otimes (A\otimes I).\\
& (B\otimes A)\otimes I\ar[ur]_\alpha&
}$$
\end{itemize}
This data is required to satisfy three pasting diagrams, which we omit, but which are well presented in Johnson and Yau, \cite{johnson-yau:bimonoidal:2021}, and in \cite{GPS,Gurski_book}, from the point of view of the more general tricategories. In a string diagram form, they are also to be found  in Corner, \cite{corner:thesis:day:2016}. These are easier to draw, but still quite complex to read.\end{Definition}

To ease our way towards a sketch of the definition of symmetric monoidal bicategory, we will briefly recall the corresponding definition of symmetric monoidal category. Although originally introduced directly by specifying that there was  a natural isomorphism, $X\otimes Y\cong Y\otimes X$, satisfying certain axioms, for our purposes it is slightly better to go via the definition of a braided monoidal category, so we briefly recall that first.
(We have adapted the definition given in Etingof, Gelaki, Nikshych and Ostrik,
\cite{etingof-et-al:tensor:2015}.)
\begin{Definition}\label{def:braided mon cat}
A \emph{braided monoidal category} is a monoidal category, $(\mathcal{C},\otimes,I)$, equipped with a natural isomorphism, $$R_{X,Y}:X\otimes Y\cong Y\otimes X,$$ called the \emph{braiding}, such that the diagrams below commute,
$$\vcenter{\xymatrix@R=19pt@C=15pt{X\otimes(Y\otimes Z)\ar[r] & (Y\otimes Z)\otimes X\ar[d]\\
(X\otimes Y)\otimes Z\ar[u]\ar[d]&Y\otimes (Z\otimes X)\\
(Y\otimes X)\otimes Z\ar[r]&Y\otimes(X\otimes Z)\ar[u]
}}\qquad \textrm{ and } \qquad \vcenter{\xymatrix@R=19pt@C=15pt{(X\otimes Y)\otimes Z\ar[r]&Z\otimes (X\otimes Y)\ar[d]\\
X\otimes (Y\otimes Z)\ar[u]\ar[d]&(Z\otimes X)\otimes Y\\
X\otimes (Z\otimes Y)\ar[r]&(X\otimes Z)\otimes Y,\ar[u]
}}$$
for all choices of objects, $X,Y, Z$, in $\mathcal{C}$, and where each arrow is an evident application of the associator, its inverse  or of the braiding.
\end{Definition}
 \begin{Definition}
 A braided monoidal category, $\mathcal{C}$, is said to be symmetric if, for all $X$, $Y$ in $\mathcal{C}$,
 $R_{Y,X}\circ R_{X,Y}= id_{X\otimes Y}.$
 \end{Definition}

 A symmetric monoidal \emph{bicategory} categorifies the above, so replacing equalities by structural morphisms. A complete description of symmetric monoidal bicategories can be found in \cite{Gurski_coherence,Gurski_Osorno}. Loosely  speaking, we have a bicategory, $\mathcal{B}$, which is monoidal as above (Definition \ref{def-mon-bicat}).  The monoidal structure is assumed to be \emph{braided}, so, we have a pseudo-natural transformation, $R\colon \otimes \to \otimes \circ \tau$, of bifunctors from  $\B \times \B$  to $\B$, where $\tau$ arises from swapping coordinates, \cite[page 4234]{Gurski_coherence}, so, in particular,
 for every $X,Y$ in $\mathcal{B}$, there is an \emph{equivalence} (within $\mathcal{B}$),
 $$R_{X,Y}:X\otimes Y\xrightarrow{\simeq} Y\otimes X,$$
 and also invertible 2-cells between the two obvious  composites from $(X\otimes Y)\otimes Z$ to $Y\otimes (Z\otimes X)$, and similarly from $X\otimes (Y\otimes Z)$ to $(Z\otimes X)\otimes Y$, so replacing equality, in the diagrams of Definition \ref{def:braided mon cat}, by invertible 2-cells.  (We will denote these 2-cells by $R_{X|YZ}$ and $R_{XY|Z}$ as seems to be the fairly standard notation currently in use; see \cite{johnson-yau:bimonoidal:2021} and  \cite{stay:compact-closed-bicat:2016}.) A full definition of braided monoidal bicategories can be found in \cite[Subsection 2.4]{Gurski_coherence}.

 There is another intermediate step before getting to the final form for `symmetric', rather than merely `braided'. Following \cite[1.1. Definitions]{Gurski_Osorno}, a \emph{sylleptic monoidal bicategory} is a braided monoidal bicategory with a \emph{syllepsis}. Such a structure is a natural isomorphism,  given by, for each pair $X,Y$ of objects, an isomorphism,
 $$\nu_{XY}:R_{YX}R_{XY}\xrightarrow{\cong}Id_{X\otimes Y}.$$
 The structure of a symmetric monoidal bicategory is, then, to satisfy  one additional axiom which says that the two ways of rewriting $R_{XY}R_{YX}R_{XY}$ to $R_{XY}$, one using $\nu_{XY}$, the other using $\nu_{YX}$, agree.

 \begin{Example} \label{Examples:symm.mon.bicat} An excellent list of examples of symmetric monoidal bicategories can be found on page 2 of Mike Stay's paper, \cite{stay:compact-closed-bicat:2016}.  We will select a few of most relevance to this work, adapting some to fit the context here.  We will also add a few others. Yet others will be included later on, once the necessary terminology has been introduced, and, for those here, we will simply mention them briefly, with a  reference to where they are discussed later.
 \begin{itemize}[leftmargin=1cm]

 \item $\mathcal{V}\!-\!\mathbf{Cat}$: If $\mathcal{V}$ is a symmetric monoidal category, then the 2-category, $\mathcal{V}\!-\!\mathbf{Cat}$, of $\mathcal{V}$-categories, enriched functors and enriched natural transformation, forms a symmetric monoidal 2-category, and thus a symmetric monoidal bicategory; see Kelly, \cite{Kelly:enriched:TAC}, page 12. In particular, this applies if $\mathcal{V}$ is the symmetric monoidal category of small categories, or if $\mathcal{V}=\Sets$.
 \item \textbf{Bicategories with finite products}: One has that any bicategory, $\mathcal{A}$,  with binary product, $-\times -$, and terminal object, $1$, underlies
a symmetric monoidal bicategory with $-\times -$  as its tensor product and $1$ as its unit object; see Theorem 2.15 of Carboni, Kelly, Walters and Wood, \cite{Carboni-Kelly-Walters-Wood:Cartesian-bicatII:2008}.  Essentially the same arguments work for bicategories with finite coproducts.\label{bicats with coproducts}

 \item $\mathbf{Span(\mathcal{C})}$:  As is well known, a span from $A$ to $B$ in a category, $\mathcal{C}$, is a diagram
$$ \xymatrix@R=-5pt{ &&C\ar[dl]_f\ar[dr]^g\\  & A & &B.  }$$
For each pair $A$, $B$, there is a category, $Span(\mathcal{C})(A,B)$,  with objects spans $(f,C,g)$, from $A$ to $B$, and morphisms, from $(f,C,g)$ to $(f',C',g')$, consisting of morphisms $C\to C'$ making the obvious diagram commute.  If $\mathcal{C}$ has pullbacks, then we can horizontally compose spans, as show e.g. in Lemma \ref{main3}.  This gives a bicategory, $\mathbf{Span(\mathcal{C})}$; see Borceux, \cite{borceux1}, Examples 7.7.3.

If $\mathcal{C}$ is a category with finite products, then the bicategory,  $\mathbf{Span(\mathcal{C})}$, is a symmetric monoidal bicategory\footnote{In fact, $Span(\mathcal{C})$ is a compact closed bicategory in the sense of Stay's paper.}, in which the tensor product on both objects and spans is given by the product.\label{Span(C)-fin prod}

This needs taking apart a little as there are subtleties that are important later on. Of course, the objects of  $\mathbf{Span(\mathcal{C})}$ are just the objects of $\mathcal{C}$.  Given any two objects, $A_1$ and $A_2$,  in $\mathbf{Span(\mathcal{C})}$, and thus in $\mathcal{C}$, their tensor product will be $A_1\otimes A_2:=A_1\times A_2$, whilst the tensor product of two spans is
$$(A_1\leftarrow C_1\rightarrow B_1)\otimes (A_2\leftarrow C_2\rightarrow B_2) :=(A_1\times A_2\leftarrow C_1\times C_2\rightarrow B_1\times B_2),$$
where the maps are as one would expect.

To define the associator on objects, we suppose that we have three objects, $A_1$, $A_2$ and $A_3$, and we need a `morphism' (in $\mathbf{Span(\mathcal{C})}$),
$$\alpha_{A_1A_2A_3}:(A_1\otimes A_2)\otimes A_3\to A_1\otimes (A_2\otimes A_3).$$
Within the base category, $\mathcal{C}$, we have an associator  (iso)morphism, (the usual one coming from the universal property of products),
$$a_{A_1A_2A_3}:(A_1\times A_2)\times A_3\to A_1\times (A_2\times A_3),$$which satisfies the requirements that the pentagon diagrams commute. (Remember $\mathcal{C}$ is a monoidal category with the product as tensor, so it is in a simpler setting than $\mathbf{Span(\mathcal{C})}$.)

In $\mathbf{Span(\mathcal{C})}$, as we said, the associator transformation is to be made of spans,  and the one that works is
$$(A_1\times A_2)\times A_3\xleftarrow{id} (A_1\times A_2)\times A_3\xrightarrow{a_{A_1A_2A_3}}  A_1\times (A_2\times A_3),$$
in other words, using the way that $\mathcal{C}$ can be thought if as being `embedded' in $\mathbf{Span(\mathcal{C})}$, using the second legs of the spans.

The associator is not just $\alpha_{A_1A_2A_3}$,  but has to be part of an adjoint equivalence, so we need a $\alpha^*_{A_1A_2A_3}$ going the other way, which is given by the reverse span, i.e. using the `first leg', and we also need $\eta$ and $\varepsilon$ as in Definition \ref{def-mon-bicat}. That `second leg' then quickly shows how to specify $\eta$ and $\varepsilon$ for $\mathbf{Span(\mathcal{C})}$ in terms of the corresponding ones in $\mathcal{C}$.

The unitors and the braidings are similarly handled giving them first in $\mathcal{C}$ before transferring them to $\mathbf{Span(\mathcal{C})}$ using the second leg process. More general results are described in much more detail in \cite{stay:compact-closed-bicat:2016}.

\item $\mathbf{Cospan(\mathcal{C})}$: If we replace $\mathcal{C}$ by the opposite category then we have that, if $\mathcal{C}$ has finite colimits,  $\mathbf{Cospan(\mathcal{C})}$  will be a symmetric monoidal bicategory, having coproduct, $\sqcup$, as its tensor product.

\item $\Prof$: The bicategory of profunctors, or distributeurs. This will be revisited after we have recalled the basic theory in the next section.
\item $\vProf$: The bicategory of $\Vect$-enriched categories, enriched profunctors and enriched natural transformations.
\item  $\tcob{d}$: Let $d$ be a non-negative integer.  It is well known that $\cob{d}$, the category of closed smooth $d$-manifolds and diffeomorphism classes of cobordisms between them, \cite{Milnor,Lec_TQFT}, forms a symmetric monoidal category with coproduct /  disjoint union, $\sqcup$, as the tensor product.

As proved in \cite{Schommer-Pries}, $\tcob{d}$, the bicategory of closed smooth $d$-manifolds, their cobordisms, and diffeomorphism classes of extended cobordisms between cobordisms,
 is a symmetric monoidal bicategory, again having $\sqcup$ as its tensor product. We will sketch the construction of the symmetric monoidal structure of $\tcob{d}$ in \S \ref{sec:mon_tcob}.
 \item  ${Alg(R)}$, also denoted $\Mor$, the bicategory of $R$-algebras, bimodules and bimodule maps. We mentioned this important example earlier in Example \ref{Alg(R)}. Here $R$ is a commutative ring.

  \end{itemize}
  \end{Example}

We now sketch the definition of a symmetric monoidal bifunctor, following  Definition 2.5 of \cite{Schommer-Pries}, and \cite{Gurski_coherence,Gurski_book,Gurski_Osorno}.

 \begin{Definition}\label{Sym. Mon. Bifunctor}
 A \emph{symmetric monoidal bifunctor}, $F\colon \mathcal{A}\to \mathcal{B}$, consists of: 
 \begin{itemize}[leftmargin=1cm]
 \item a homomorphism, (i.e. a bifunctor), $F\colon \mathcal{A}\to \mathcal{B}$, between the underlying bicategories;
 \item a transformation, (i.e. a pseudo-natural transformation), $$\chi\colon \otimes_\B\circ (F\times F) \Rightarrow F\circ \otimes_\A,$$ of bifunctors, from $\A\times \A$ to $\B$, so we have, given objects, $A_0$ and $A_1$, of $\A$, a 1-morphism, in $\B$,
 $$\chi_{A_0,A_1}:F(A_0)\otimes_\mathcal{B}F(A_1)\to F(A_0\otimes_\mathcal{A} A_1),$$
 and given 1-morphisms, $f_0\colon A_0 \to A_0'$ and $f_1\colon A_1\to \A_1'$, in $\A$, we have a natural 2-cell in $\B$,
 $$\xymatrixcolsep{2cm}\xymatrix@R=15pt{
&F(A_0)\otimes_\B F(A_1)   \ar[r]^{F(f_0) \otimes_\B F(f_1)}\ar[d]_{{\chi}_{(A_0,A_1)} } \drtwocell<\omit>{\hspace*{0.8cm}{\chi}_{(f_0,f_1)}} & F(A'_0)\otimes_\B F(A'_1) \ar[d]^{{\chi}_{(A_0',A_1')} } \\
& F(A_0\otimes_\A A_1)
 \ar[r]_{F ( f_0\otimes_A f_1)}  &  F(A'_0 \otimes_\A A'_1),}
$$
which is compatible with horizontal compositions and horizontal identities,
 and we also have a transformation,
 $$i\colon I_\mathcal{B}\Rightarrow F \circ I_\mathcal{A},$$
 together with  corresponding adjoint equivalence transformations, $\chi^*$ and $i^*$, (and with the relevant adjunction data);
 \item invertible modifications, $\omega$, $\gamma$ and $\delta$, measuring compatibility with the relevant associators and unitors (see Fig 2.5. of \cite{Schommer-Pries});
 \item[]\hspace*{-1cm}and
 \item an invertible modification, $u$, giving compatibility with the braiding (see Fig 2.6 of \cite{Schommer-Pries}, or \cite[page 4239]{Gurski_coherence}, for the relevant diagram, so that
 $$u:F(R_\mathcal{A})\circ \chi \Rightarrow \chi\circ R_\mathcal{B}:F(B)\otimes_\mathcal{B}F(A)\to F(A\otimes_\mathcal{A} B).$$
 \end{itemize}
 This data is to satisfy certain axioms, which we omit, referring to the discussion in \cite[Definition 2.5]{Schommer-Pries} for further details, including further references.  In particular the compatibility conditions in \cite[\S 4.3]{Gurski_book} / \cite[page 17]{GPS} hold, dealing with the preservation of the monoidal structure by $F$, and those of \cite[page 4239]{Gurski_coherence} hold, similarly describing  the symmetric monoidal structure of the bifunctor $F$.
 \end{Definition}

 We will sketch the construction of the entire structure of a particular symmetric monoidal bifunctor when we prove, in \S \ref{Quinn_is_sym_mon}, that the once-extended Quinn TQFT,  $\tFQ{\Bc}\colon \tcob{n} \to \vProfGrphf$, is symmetric monoidal.
\subsection{Conventions on profunctors}\label{sec:prof_conventions}

The detailed theory of profunctors can be found in many texts and on-line sources, for instance, \cite[Chapter 7]{borceux1}, \cite[Section 5]{loregian_2021}  and the nLab \cite{nLab}. The enriched version of profunctor bicategories can be found in e.g. \cite[1.3. Distributors]{GambinoJoyal}.  When searching for such theory, it is important to note that the terms `distributor' and `bimodule' are often-used alternative names for profunctors. We will give a very minimal sketch here.
\subsubsection{Some background and basic definitions}

In a general context, given two small categories, $\mathcal{A}$ and $\mathcal{B}$, we have:

\begin{Definition} A ($\Sets$-valued) \emph{profunctor} from $\mathcal{A}$ to $\mathcal{B}$ is a functor,
$$F\colon \mathcal{A}^{op}\times \mathcal{B}\to \Sets.$$
\end{Definition}
\noindent This will be sometimes written $F\colon \mathcal{A}\nrightarrow \mathcal{B}$, or   $F\colon \mathcal{A}\to \mathcal{B}$ in commutative diagrams.

\begin{Example}\label{Ex:maps_to_profunctors}
Given a functor, $F:\mathcal{A}\to \mathcal{B}$, we can define two profunctors,
$$\varphi^F:\mathcal{A}\bto \mathcal{B}, \qquad  \textrm{ and }  \qquad  \varphi_F:\mathcal{B}\bto \mathcal{A},$$
by $$\varphi^F(A,B)= \mathcal{B}(F(A),B),  \qquad  \textrm{
whilst }  \qquad  \varphi_F(B,A)=\mathcal{B}(B,F(A)).$$
\end{Example}

\begin{Remark} There are several different conventions used in the literature as to the `direction' of the profunctor. One of the most current, but not the one that we will use, is to say a profunctor, $F:\mathcal{A}\nrightarrow \mathcal{B}$, is a functor, $\mathcal{B}^{op}\times \mathcal{A}\to \Sets$. \end{Remark}

With a suitable notion of composition, these profunctors, {together with their natural transformations},  form a bicategory $\Prof$, see for instance \cite[Chapter 7]{borceux1}, or \cite[Section 5]{loregian_2021}, which is symmetric monoidal, for which see \cite[1.3. Distributors]{GambinoJoyal} or \cite{Bartlett_etal,Hansen-Shulman:constructing:2019}. We will briefly recall the construction of this bicategory here, however in the slightly enriched context of $\Vect$-profunctors between groupoids.

\begin{Definition}  Given $\Vect$-enriched categories, $\mathcal{A}$ and  $\mathcal{B}$, a \emph{$\Vect$-enriched profunctor, from $\mathcal{A}$ to $\mathcal{B}$,} is an enriched functor,
$F\colon \mathcal{A}^{op}\times \mathcal{B}\to \Vect.$

Given groupoids,  $G$ and $G'$, a \emph{$\Vect$-profunctor from $G$ to $G'$}, denoted $\Hp\colon G \bto G'$, is a functor $\Hp\colon G^\op\times G' \to \Vect$.
\end{Definition}
Let $G$ and $G'$ be groupoids. Any functor, $F\colon {G}^{op}\times G'\to \Sets$, determines  a functor,  $\mathbf{F}\colon {G}^{op}\times G'\to \Vect ,$ by composing $F$ with the free vector space functor, $\Lin\colon \Sets \to \Vect$.
Note  also that a groupoid $G$ can be `converted' to a $\Vect$-enriched category $Lin(G)$,  by applying the free vector space functor to the $\hom$-sets of $G$. Any $\Vect$-profunctor $G \bto G'$ gives rise to a $\Vect$-enriched profunctor, from  $Lin(G)$ to  $Lin(G')$, canonically.

 \begin{Example}\label{identity profunctor} Let $\mathcal{A}$ be a (small) category, then we have the bivariant hom-functor, $\mathcal{A}(-,-):\mathcal{A}^{op}\times \mathcal{A}\to \Sets$, which is a $\Sets$-valued profunctor from $\mathcal{A}$ to itself.  If $\mathcal{A}$ is a $\Vect$-enriched category, the natural analogue of the above is $\mathcal{A}(-,-):\mathcal{A}^{op}\times \mathcal{A}\to \Vect$, and so is a $\Vect$-enriched profunctor. We denote the latter by  $\Id_\mathcal{A}$, and call it the \emph{identity profunctor} on the (linear) category, $\mathcal{A}$. We may shorten this to $\Id_G$, if $\mathcal{A}$ is the linearisation, $Lin(G)$, of a groupoid, $G$.
\end{Example}

\begin{Definition} Given $\Vect$-profunctors, $\Hp,\Hp'\colon G^{\op}\times G' \to \Vect$,  a 2-morphism, or 2-cell, $\Eta\colon \Hp \To \Hp'$,  between them, is a natural transformation of functors, $G^{\op}\times G' \to \Vect$. (Hence, we have, given $x\in G_0$ and $y\in G'_0$, a linear map, $\Eta_{x,y}\colon \Hp(x,y) \to \Hp'(x,y)$, which is natural in both $x$ and $y$.)
\end{Definition}

Let $G,H,K$ be groupoids, with sets of objects $G_0,H_0$ and $K_0$. Given $\Vect$-profunctors, $\Hp\colon G \bto H $, and $\Hp'\colon H \bto K$, their composite, $\Hp\bullet \Hp'\colon G \bto K$, will be  the $\Vect$-profunctor such that, if $x \in G_0$ and $z \in K_0$, then\footnote{Note that this coend is, \textit{a priori}, defined only up to isomorphism. In this paper, we always implicitly choose a natural realisation for all limits, colimits and coends appearing.}
\begin{equation}\label{Profunctor composition}(\Hp\bullet \Hp')(x,z):=\int^{y \in H_0} \Hp(x,y)\otimes\Hp'(y,z)=\Big(\bigoplus_{y \in H_0}  \Hp(x,y)\otimes\Hp'(y,z)\Big )/\simeq.\end{equation}
Here, fixing $x \in G$ and $z \in K$, the equivalence relation, $\simeq $, is generated (as a linear equivalence relation\footnote{i.e. as an equivalence relation whose quotient is a vector space.}) by
 \\ for $y, y' \in Y$, $v_{x,y} \in \Hp(x,y)$ and  $v'_{y',z} \in \Hp(y',z)$ and an arrow, $y\ra{h} y'$, in $H$,
$$  v_{x,y} \otimes \Hp'( y\ra{ h} y', z \ra{1_z} z)( v'_{y',z}) \simeq  \Hp(x \ra{1_x} x, y\ra{ h} y') (v_{x,y} ) \otimes  v'_{y',z},  $$
or, more informally,
\begin{equation}\label{profunctor comp as bimodule tensor}v_{x,y} \otimes h\cdot  v'_{y',z} \simeq v_{x,y}\cdot h\otimes  v'_{y',z}.\end{equation}

We note the convention on the order of composition that we are using. This convention is used because it reflects the geometric intuition, being a concatenation order of composition. It also reflects a useful convention for the bicategory, $\Mor$, of algebras, bimodules and bimodule maps, to which $\Prof$ is closely related.

If we just have $\Sets$-valued profunctors,  this formula for composition still makes sense by interpreting $\otimes$ as $\times$, and we note that $\Lin\colon \Sets \to \Vect$ preserves that composition in the evident way.

Of course, there is a projection, which we will need later on,
\begin{equation}\label{proj:prof}proj \colon \bigoplus_{y \in H_0}  \Hp(x,y)\otimes\Hp'(y,z)\to (\Hp\bullet \Hp')(x,z).\end{equation}  Given any element in $ (\Hp\bullet \Hp')(x,z),$ we can represent it by an element in $\Hp(x,y)\otimes\Hp'(y,z)$, for some $y$, but, working with that, just as in the very similar setting of tensor product of bimodules, any resulting calculation has to be shown to be  invariant under the action of the arrows of $H$.

\begin{Example} \label{varphi is pseudo}Suppose we have functors, $\mathcal{A}\xrightarrow{F}\mathcal{B}\xrightarrow{G}\mathcal{C}$, then we have corresponding profunctors, $\varphi^F$ and $\varphi^G$,  as in Example \ref{Ex:maps_to_profunctors}, so we can form their composite, $\varphi^F\bullet \varphi^G$.  It is not hard to check that   we have a natural isomorphism, $\varphi^F\bullet \varphi^G\To \varphi^{GF}$. This is part of the data that says that $\varphi^{(-)}$ is a pseudo-functor, from the bicategory of small categories, functors and natural transformations, to the bicategory, $\Prof$, of small categories, profunctors and natural transformations between profunctors. This is explained e.g. in \cite[1.3. Distributors]{GambinoJoyal}.
\end{Example}

The composition of general profunctors, including that of $\Vect$-profunctors, has left and right (lax) identities. Suppose that we have $\Hp\colon G \bto H $, then we can compose it with the  hom-functor, $G(-,-)\colon G^{op}\times G\to \Sets$, or to $\Vect$ after applying $\Lin$. This acts like a left identity on $G$. We have natural isomorphisms, called the `left unitor' and  the `right unitor',
\begin{equation}\label{left unitor}\lambda_G^\Hp\colon G(-,-)\bullet \Hp\To \Hp \quad \textrm{ and } \quad
\rho_H^\Hp\colon \Hp\bullet H(-,-)\To \Hp,
\end{equation}  These are easy to write down, for example,  using Equation \eqref{profunctor comp as bimodule tensor}.

Given  natural transformations, $\Eta\colon \Hp_1 \To \Hp_2$, between $\Vect$-profunctors, $G \bto H$, and $\Eta'\colon \Hp_1' \To \Hp_2'$, between $\Vect$-profunctors, $H \bto K$, we have a natural transformation, $(\Eta\bullet \Eta')\colon \Hp_1\bullet \Hp_1' \to \Hp_2\bullet \Hp_2'$.
Explicitly, given $x \in G_0$ and $z \in K_0$, then  $(\Eta\bullet\Eta')_{(x,z)}$ sends the equivalence class of $v_{x,y}\otimes v'_{y,z}$ to the equivalence class of $\Eta_{(x,y)}(v_{x,y})\otimes \Eta'_{(y,z)}(v'_{y,z})$. Here $y\in H_0$, $v_{x,y}\in \Hp_1(x,y)$ and $v'_{y,z}\in \Hp_2(y,z)$.
In other words, given $x \in G_0$ and $z \in K_0$, the linear map,   $(\Eta\bullet\Eta')_{(x,z)}$, is the unique map that makes the diagram below commute:
\begin{equation}\label{hor_comp:etas} 
\vcenter{\xymatrix@R=13pt
{ & \save[]+<0cm,0cm>*{ {\displaystyle \bigoplus_{y \in H_0}  \Hp_1(x,y)\otimes\Hp_1'(y,z)} \,\,} \ar@<4pt>[rr]^{\proj} \restore
\ar[dd]_{\displaystyle \bigoplus_{y \in H_0} \Eta_{(x,y)} \otimes \Eta'_{(y,z)}}  && \displaystyle  \int^{y \in H_0} \Hp_1(x,y)\otimes\Hp_1'(y,z) \ar[dd]^{ \displaystyle (\Eta\bullet\Eta')_{(x,z)}}\\
\\
& \displaystyle \bigoplus_{y \in H_0}  \Hp_2(x,y)\otimes\Hp_2'(y,z)  \ar[rr]_{\proj} && \displaystyle \int^{y \in H_0} \Hp_2(x,y)\otimes\Hp_2'(y,z)\, .
}}
\end{equation}

\subsubsection{The symmetric monoidal bicategory $\vProfGrp$}\label{sum:conv_prof}
This bicategory will be a primary target for our once-extended Quinn TQFT.
Our starting point to construct $\vProfGrp$ is that we have a bicategory,  $\vProf$, the bicategory of $\Vect$-enriched categories, $\Vect$-enriched profunctors and $\Vect$-enriched natural transformations, which is symmetric monoidal; see \cite[Corollary 6.6.]{Hansen-Shulman:constructing:2019}, \cite[Remark 1.3.4]{GambinoJoyal}, and also \cite{cattani-winskel:profunctors:2005}. We  define $\vProfGrp$, as the full sub-bicategory of $\vProf$,  whose objects are groupoids, $G$, each made into a $\Vect$-enriched category by applying the free vector space functor to the sets of morphisms, $\hom_{G}(x,y)$, of $G$.

Here is an explicit description of $\vProfGrp$.

\begin{Definition}[$\vProfGrp$]Let $\kappa$ be a field. The bicategory $\vProfGrp$ is defined as the bicategory such that:
\begin{enumerate}[leftmargin=1cm]
 \item the objects are groupoids;
 \item given groupoids, $G$ and $G'$, 1-morphisms are $\Vect$-profunctors, $\Hp\colon G \bto G'$;
 \item the identity profunctor on a groupoid $G$ is the $\Vect$-profunctor, $\Id_G\colon G \bto G$, defined in Example \ref{identity profunctor};
 \item composition of 1-morphisms is defined in Equation \eqref{Profunctor composition};
 \item 2-morphisms between $\Vect$-profunctors, $\Hp,\Hp'\colon G \to G'$, are natural transformations, $\Eta\colon \Hp \To
 \Hp'$, of functors  $G^{op}\times G' \to \Vect$, \item the vertical composition of 2-morphisms is the usual composition of natural transformations, of functors $G^{op}\times G' \to \Vect$;
 \item  the left and right unitors are as in equation \eqref{left unitor};
 \item the horizontal composition of 2-morphisms is explained in Equation \eqref{hor_comp:etas}.
\end{enumerate}
\end{Definition}
\noindent The associators in $\vProfGrp$ can be written down with no difficulty.

The bicategory $\vProfGrp$ inherits a symmetric monoidal structure from $\vProf$, where:
\begin{itemize}[leftmargin=1cm]
 \item the tensor product of groupoids $G$ and $G'$ is their usual cartesian product,
 \item the tensor product of two $\Vect$-profunctors $\Hp_1\colon G_1  \times G_1' \to \Vect$ and $\Hp_2\colon G_2   \times G_2'$, written $(\Hp_1\otimes \Hp_1)\colon G_1 \times G_2 \bto  G_1' \times G_2'$, is defined as the composite:
 $$(G_1\times G_2)^{op} \times (G_1'\times G_2') \xrightarrow{\cong } (G_1^{\op}\times G_1') \times (G_2^{op}\times G_2') \xrightarrow{\Hp_1\times \Hp_2} \Vect \times \Vect \xrightarrow{\Otimes} \Vect. $$
 \end{itemize}
From the previous two items, we can then build a bifunctor $\Otimes \colon \vProfGrp\times \vProfGrp \to \vProfGrp$.
\begin{itemize}[leftmargin=1cm]
\item The  associativity and braiding profunctors are inherited from  the associativity and braiding morphisms of the symmetric monoidal 2-category of small categories, functors and natural transformation,  by applying the $\varphi^{(-)}$-construction of Examples \ref{Ex:maps_to_profunctors} and \ref{varphi is pseudo}, and then linearising. In particular, the braiding in $\vProfGrp$ is given explicitly as follows.

Suppose $G$ and $H$  are groupoids, then there is an isomorphism   $$R: G\times H\xrightarrow{\cong} H\times G.$$
We therefore have a (set-valued) profunctor, $\varphi^R:G\times H\bto H\times G$. Applying $\Lin\colon \Sets \to \Vect$, we then get a $\Vect$-profunctor $\boldsymbol{\varphi^R}\colon G\times H\bto H\times G$.
\end{itemize}

Finally define, using the nomenclature in \S\ref{sec:convention_groupoids}:
\begin{itemize}[leftmargin=0.6cm]
\item $\vProfGrphf$: the full  sub-bicategory of $\vProfGrp$, whose objects are the homotopy finite groupoids.
\item $\vProfGrpfin$: the full  sub-bicategory of $\vProfGrp$ with objects the finite groupoids.
\end{itemize}
These are likewise symmetric monoidal bicategories, whose structure is inherited from that of $\vProfGrp$, in the obvious way.
\subsubsection{Matrix elements for natural transformations}
Let $G$ and $H$ be groupoids. In most examples appearing in this paper, $\Vect$-profunctors, $G^{\op}\times H \to \Vect$,
 are linearisations, $\mathbf{F}=\Lin\circ F$, of set-valued profunctors, $F\colon G^{\op}\times H \to \Sets$.  However, the   natural transformations we construct between such $\Vect$-profunctors will not, in general, arise from linearising natural transformations between $\Sets$-valued profunctors. Nevertheless, they have a simple description, as we now explain.

Given a profunctor, $F{\colon} G\bto H$, its linearisation, $\mathbf{F}\colon G\bto H$, comes with given bases on each of its constituent vector spaces. We  make the assumption that all the $F(x,y)$ are finite set, as this holds in all situations met in this paper.

An object of $G^{op}\times H$ is, of course, a pair, $(x,y)\in G_0\times H_0$. For $g\colon x\to x'$ in $G_1$, we then have, for each $y\in H_0$, a linear map, $\mathbf{F}(g,y)\colon \mathbf{F}(x',y)\to \mathbf{F}(x,y)$. Similarly, if $h{\colon}y\to y'\in H_1$, we have $\mathbf{F}(x,h){\colon} \mathbf{F}(x,y)\to \mathbf{F}(x,y')$.  This looks like a `many object' bimodule, and we will recall the relationship with the theory of bimodules over algebras more fully in \S \ref{Groupoid algebra}.

For the moment, we need to examine the way of describing the natural transformations between such profunctors.  Suppose $\mathbf{F},\mathbf{F}'{\colon}G^{\op}\times H \to \Vect$ are two $\Vect$-valued profunctors, linearised from some $\Sets$-valued ones, $F$ and $F'$, as above.  Further suppose $\varphi{\colon}\mathbf{F}\Rightarrow\mathbf{F}'$ is a natural transformation from $\mathbf{F}$ to $\mathbf{F}'$.  For each $(x,y)\in G_0\times H_0$, we then have a linear mapping, $\varphi(x,y){\colon}\mathbf{F}(x,y)\to \mathbf{F}'(x,y)$, and hence, for each $f\in F(x,y)$ and $f'\in F'(x,y)$,  a matrix element, $\langle f\mid \varphi(x,y)\mid f'\rangle$, so that we have a `state sum' description, $$\varphi(x,y)(f) =\sum_{f'\in F'(x,y)} \langle f\mid \varphi(x,y)\mid f'\rangle f'.$$  The fact that $\varphi$ is a natural transformation means that it must be compatible with changes along any $g{\colon}x'\to x$ and $h{\colon}y\to y'$, and so must satisfy equations involving the various $\mathbf{F}(g,h){\colon} \mathbf{F}(x,y)\to \mathbf{F}(x',y')$\label{nat-trans as matrices}. This will be extensively used in  \S\ref{sec_resolvedspans_to_nat_tranfs}.
\subsubsection{{Some lemmas on coends of functors from groupoids}}\label{some:lemmas:coends}
For the convenience of the reader, we collect a few elementary lemmas, whose explicit formulation can be difficult to find in the literature. They will be useful when giving some explicit details in the  proof that the once-extended Quinn TQFT is indeed  a bifunctor,   when it comes to preservation of horizontal compositions, see \S \ref{sec:path-components-pull} and \S \ref{Hor-comp-resolved 2-cells}.

Let $G=(s,t\colon G_1 \to G_0)$ be a groupoid. Consider a functor $F\colon G^{\op} \times G \to \Sets$. The coend of $F$ is a universal wedge, making the diagram,
$$\xymatrixcolsep{10mm}\xymatrix{ F(y,x)\ar[r]^{F(g, \id_x)}\ar[d]_{F(\id_y , g)} &F(x,x)\ar[d]^<<<<<{p_x} \\
             F(y,y) \ar[r]_<<<<<{p_y} &  \displaystyle\int^{z \in G_0}F(z,z),}
$$
commute, for all choices of morphisms $g\colon x \to y$ in $G$. Therefore, as written above,
$$\int^{z \in G_0}F(z,z)= \Big(\bigsqcup_{z \in G_0} F(z,z)\Big)\big/\sim,$$
where $\sim$ is the smallest equivalence relation that makes the diagram above commute for all choices of $g\colon x \to y$.
As $G$ is a groupoid, given any $g\colon x \to y$ and $g'\colon x'\to y'$, the map,
$$F(g, g')\colon F(y,x') \to F(x,y'),$$
is a bijection. This gives:
\begin{Lemma} We have
$$\int^{z \in G_0}F(z,z)= \Big (\bigsqcup_{z \in G_0} F(z,z)\Big )\big /\simeq,$$
where $\simeq$ is the equivalence relation in which $$u_x \in F(x,x) \simeq u_y \in F(y,y)$$ if there exists $g\colon x \to y$ such that:
$$F\left( g^{-1}\colon y \to x,g\colon x \to y\right)(u_x)=u_y.   $$
\end{Lemma}
\noindent In fact, this shows that the equivalence relations $\sim$ and $\simeq$ are really  the same.

Note that any groupoid $G$ comes with a \emph{contravariant} functor $(-)^{-1}\colon G \to G$, that is the identity on objects and sends $g\colon x \to y$ to $g^{-1}\colon y \to x$.
There is also the diagonal functor, $\Delta\colon G \to G\times G$, sending $x \in G_0$ to $(x,x)\in G_0\times G_0$, and with $$\Delta\big (x\ra{g} y\big)=\big ((x,x) \ra{(g,g)} (y,y)\big). $$  Hence we have a functor,
$$ F\circ \big ( (-)^{-1} \times \id\big)\circ \Delta \colon G \to \Sets.$$

The previous lemma gives:
\begin{Lemma}\label{lem:colimstoprof}
There is a canonical bijection
$$ \int^{z \in G_0}F(z,z)\cong \mathrm{colim} \Big (F\circ \big ( (-)^{-1} \times \id\big)\circ \Delta\Big).$$
\end{Lemma}

Let $A$ and $B$ be sets. Consider a linear map, $f\colon \kappa(A) \to \kappa(B)$, between free vector spaces and  equivalence relations, $\sim_A$ and $\sim_B$, on $A$ and $B$, such that $f$ descends to a map $f'\colon \kappa(A/\sim_A) \to \kappa(B/\sim_B),$ then, given $a \in A$ and $b \in B$, the matrix elements of $f'$ satisfy
$$\big\langle [a] | f'| [b]\big\rangle=\sum_{b' \in [b]} \big\langle a | f| b'\big \rangle.$$
Combined with the previous discussion, this gives the following.
\begin{Lemma}\label{lem:profmaps}
 Let $F,F'\colon G^\op\times G \to \Sets$ be functors. Consider, for each $(x,y) \in G_0 \times G_0$, a linear map, $\eta_{(x,y)}\colon \kappa\big ( F(x,y)\big) \to \kappa \big( F'(x,y)\big)$, such that putting all of the $\eta_{(x,y)}$ together gives a  natural transformation, $\eta\colon \Lin \circ F {\To} \Lin\circ F'$, where $\Lin\colon\Sets \to \Vect$ is the free vector space functor. Further, consider the induced map (see \cite[Notation 1.1.15]{loregian_2021}) between coends, as below
 $$\int^{x \in G_0} \eta_{(x,x)} \colon  \int^{x \in G_0} \kappa ( F(x,x) ) \to \int^{x \in G_0} \kappa (  F'(x,x) ),$$
that is (since the free vector space functor preserves colimits, and with a minor abuse of notation),
 $$\int^{x \in G_0} \eta_{(x,x)} \colon \kappa \Big( \int^{x \in G_0} F(x,x) \Big) \to \kappa\Big( \int^{x \in G_0} F'(x,x) \Big).$$
Its matrix elements satisfy, for each $z \in G_0$, $u_z \in F(z,z)$ and $w_z \in F'(z,z)$,
 $$
  \Big \langle [u_z] \Big | \int^{x \in G_0} \eta_{(x,x)} \Big| [w_z] \Big \rangle= \sum_{w_z' \in \mathrm{Orb}(w_z) } \left \langle u_z  |  \eta_{(z,z)} \ | w_z'   \right \rangle.
  $$
Here $\mathrm{Orb}(w_z)$ is the orbit of $w_z$, under the action of the group $\hom_G(z,z)$, on $F'(z,z)$, defined as ,
$$v_z\trl (g^{-1}\colon z \to z):= F'( g^{-1}\colon z \to z,g\colon z \to z)(v_z),  $$where $v_z\in F'(z,z)$.

On the other hand, if $z'\in G_0$ belongs to a different connected component from $z$ in $G$, and if $t_{z'} \in F'(z',z')$, then
 $$
  \Big \langle [u_z] \Big | \int^{x \in G_0} \eta_{(x,x)} \Big | [t_{z'}] \Big \rangle =0.
  $$
\end{Lemma}
\begin{proof}
This follows from the previous discussion, and the fact that the diagram below commutes,
 $$\xymatrix@C=40pt{ &\displaystyle \bigoplus_{x \in G_0} \kappa(F(x,x)) \ar[d]\ar[rr]^{ \bigoplus_{x \in G_0} \eta_{(x,x)}} && \displaystyle \bigoplus_{x \in G_0} \kappa(F'(x,x))\ar[d]\\
  &\kappa \Big(\displaystyle \int^{x \in G_0} F(x,x) \Big) \ar[rr]_{   \int^{x \in G_0} \eta_{(x,x)}} &&\kappa \Big( \displaystyle\int^{x \in G_0} F'(x,x) \Big),
 }
 $$
  where the vertical arrows are the canonical projections.
\end{proof}

\subsection{Notation and some more basic results about fibrations} We need some additional results and notation about fibrations, {as defined in \S\ref{HureFib}.}

\subsubsection{Holonomy maps and the functor, $\Fc^M\colon \pi_1(B,B) \to \CGWH/\simeq$, associated to a fibration $p\colon M \to B$}\label{sec:more_fibrations}
 
Recall from Subsection \ref{sec:conventions_top} that given a CGWH space, $B$, $B^I$ is the space of paths $\gamma\colon [0,1] \to B$, and that we defined the maps, $$s_B=s,t_B=t\colon B^I \to B,$$ such that $s(\gamma)=\gamma(0)$ and $t(\gamma)=\gamma(1)$. We have a fibration, $\langle s,t \rangle \colon B^I \to B\times B$, from which we constructed the identities in $\HFb$; see Lemma  \ref{idsHFb} and Definition \ref{main_desc}.   Given $\gamma \in B^I$, the reverse path to $\gamma$ will be  denoted $\overline{\gamma}$, so $\overline{\gamma}(u)=\gamma(1-u)$, for each $u \in [0,1]$.
We also consider the map, $\const\colon B \to B^I$, sending $x \in B$ to the constant path, $\const_x$, at $x$.

Let $p\colon M \to B$ be a fibration. Recall that the fibre of $x\in B$ is denoted $$M_x:=p^{-1}(x).$$ 
Let $ M\times_B B^I$ denote the pullback of the  maps, $p\colon M \to B$ and $s\colon B^I \to B.$ Consider also the canonical projections, $\proj_1\colon M\times_B B^I \to M\textrm{ and } \proj_2\colon M\times_B B^I \to B^I,$
so we have a pullback diagram,
$$\xymatrix@R=13pt{ && M\times_B B^I\ar[dl]_{\proj_1} \ar[dr]^{\proj_2}  \\ &M \ar[dr]_p  && B^I . \ar[dl]^{s_B}\\ && B } $$

We have a continuous function, $\lambda^M\colon I\times ( M\times_B B^I)\to M $, arising from the diagram below and the homotopy lifting property of $p\colon M \to B$,
$$ \xymatrix{ &M\times_B B^I \ar[rrrr]^{\proj_1} \ar[d]_{\{0\} \times (\_)} &&&& M\ar[d]^p\\
              & I \times (M\times_B B^I)\ar@{-->}[rrrru]^{\lambda^M} \ar[rr]_{\id_I\times \proj_2}  &&I \times B^I \ar[rr]_{(u,\gamma) \mapsto \gamma(u)} && B . }$$
(Here $\big(\{0\} \times (\_)\big)(m,\gamma)=(0,m,\gamma)$.)
\begin{Definition}[
Holonomy Map]
A  function, $\lambda^M$, making the diagram above commute will be called a \emph{holonomy map on the fibration $p\colon M \to B$}.
\end{Definition}
\noindent The nomenclature ``holonomy map'' is borrowed from differential geometry. We will frequently write ``holonomy'' rather than ``holonomy map''. The holonomy maps considered here are equivalent to the ``path-lifting functions'' in  \cite[Chapter 7]{May}, and the ``lifting functions'' in  \cite{Fadell}.

The following string of classical results are  to be found, essentially, in  \cite{Fadell} or \cite[Chapter 7]{May}. They follow from simple application of the appropriate homotopy lifting property, and are `well known', but we give a reference for each one.
\begin{Lemma}\label{main_techfib} Let $p\colon M \to B$ be a fibration.
 Consider a fixed holonomy map, $\lambda^M$, on $p\colon M \to B$. Let $\gamma$ be a path, in $B$, from $x\in B$ to $y\in B$. Consider the map $\Gamma_{\gamma}^M:M_x \to M_y$, defined  by $$\Gamma_\gamma^M(m):=\lambda^M(1,m,\gamma), \textrm{ for all } m \in M_x.$$
 Up to homotopy of  maps from $M_x$  to $M_y$, the map, $\Gamma^M_\gamma\colon M_x \to M_y$, then, depends only on the homotopy class of $\gamma$ (and, in particular, not on the chosen holonomy map, $\lambda^M$). In fact,
 $\Gamma^M_{\gamma}:M_x \to M_y$ 
 is a homotopy equivalence between {the fibres} $M_x$ and $M_y$, with a homotopy inverse to  $\Gamma^M_\gamma\colon M_x \to M_y$ given by  $\Gamma_{\overline{\gamma}}^M$. \end{Lemma}
 \begin{proof} See \cite[\S 7.6 (Change of fiber)]{May}.\end{proof}

\begin{Lemma}\label{pphi-cmpt} Let $\phi$ be a path in $M$,  and $\gamma=p\circ \phi$, its image path in $B$, then $\Gamma_\gamma^M(\phi(0))$ and $\phi(1)$ are in the same path-component of $M_{\gamma(1)}=p^{-1}(\gamma(1))$.

 In particular, if $M$ is path-connected and $x,y \in B$, then if $A \in\hpi_0( M_x)$ and $A' \in\hpi_0( M_y)$ are path components of the chosen fibres, it follows that $A$ and $A'$ are homotopy equivalent. Concretely, choose $m \in A$ and $m' \in A'$, and a path, $\phi$, in $M$ connecting $m$ to $m'$, then $\Gamma_{\gamma}^M\colon M_x \to M_y$ restricts to a map, $A\to A'$, giving the desired homotopy equivalence. (Here $\gamma=p\circ \phi$.)
 \end{Lemma}
 \begin{proof} See \cite[page 3]{Fadell}.\end{proof}

 Suppose that we have paths, $\gamma\colon x \to y $ and $\gamma'\colon y \to z$, then $\Gamma^M_{\gamma \gamma'}$ is homotopic to $\Gamma^M_{\gamma'} \circ \Gamma^M_\gamma$, as maps from $M_x$ to $M_z$; see \cite[\S 7.6]{May}. Moreover, given $x \in B$, the map, $\Gamma_{\const_x}^M \colon M_x \to M_x$, is homotopic to the identity.

Recall that  $\CGWH/\simeq$ denotes  the category with objects the CGHW spaces, with morphisms being homotopy classes of maps; see \S \ref{sec:conventions_top}.

\begin{Lemma} \label{main_techfibfunctor}There is a functor, $$\Fc^M\colon \pi_1(B,B) \to \CGWH/\simeq,$$
(where $\pi_1(B,B)$ is the fundamental groupoid of $B$). Given $x \in B$, $\Fc^M(x):= M_x$, and  given a path, $x \ra{\gamma} y$, in $B$, then,  $$\Fc^M(x \ra{[\gamma]} y): = [\Gamma^M_{\gamma}]\colon M_x \to M_y.$$
Here $[\Gamma^M_{\gamma}]$ is the homotopy class of $\Gamma^M_{\gamma}\colon M_x \to M_y$. 
This functor, $\Fc^M$, depends only on the fibration, $p\colon M \to B$, and not on the chosen holonomy map, $\lambda^M$.
\end{Lemma}
\begin{proof} See \cite[\S 7.6]{May}.\end{proof}

We, thus,  have a functor  $\hpiz\circ \Fc^M \colon \pi_1(B,B) \to \Sets.$ It sends $x \in B$ to $\hpi_0(M_x)$.
Given a path, $x\ra{\gamma} y$, in $B$, the functor is such that, if $m \in M_x$,\label{hpiz}
$$\big((\hpi_0\circ \Fc^M)(x\ra{[\gamma]} y)\big)(\PC_m(M_x))=\PC_{\lambda^M(1,m,\gamma)}(M_y)= \PC_{\Gamma^M_\gamma(m)}(M_y),$$
where we recall that $\PC_m(M_x)$ denotes the path-component of $m$ in $M_x$; see  Subsection \ref{sec:conventions_top}.
This functor depends only on the fibration,  $p\colon M \to B$, and not on the chosen holonomy map, $\lambda^M$.

Finally, in this string of lemmas, we have:
 \begin{Lemma}\label{main_techfib_loops}
Let $x \in B$. There are left and right actions of $\pi_1(B,x)$ on $\hpi_0(M_x)$. These are  such that,  if  $m \in M_x$, $\gamma \in \Omega_x(B)$, the loop space of $B$ based at $x$, and  $[\gamma]$ is the associated element of $\pi_1(B,x)$,  then:
\begin{align*}
[\gamma] \t \PC_m(M_x) &=\PC_{\Gamma_{\overline{\gamma}}^M(m)}(M_x),\\
\PC_m(M_x) \trl [\gamma]&=\PC_{\Gamma_\gamma^M(m)}(M_x).
\end{align*}
\end{Lemma}

The following result will be needed when addressing why the once-extended Quinn TQFT can be given the structure of a symmetric monoidal bifunctor.
\begin{Lemma}\label{Lem:FExE'}
Let $p\colon E \to X$ and $p'\colon E' \to X'$ be fibrations. We thus have a fibration, $(p\times p') \colon E \times E' \to X \times X'$. The functor, $$\Fc^{E \times E'}\colon \pi_1(X\times X',X\times X') \to \CGWH/\simeq,$$ provided by  $(p\times p')\colon E \times E' \to X \times X'$, is given by the composition of the functors below,
 \begin{multline*}
 \pi_1(X \times X',X \times X') \cong \pi_1(X,X) \times \pi_1(X',X') \\ \ra{\Fc^E \times \Fc^{E'}}( \CGWH/\simeq) \times (\CGWH/\simeq) \ra{\Times_\CGWH} \CGWH/\simeq.
 \end{multline*}
\end{Lemma}
\noindent Here $\Times_\CGWH$ is  the product monoidal  structure on $\CGWH$; see \S \ref{sec:path-components-pull} below.
\begin{proof}
On objects, this follows from the fact that, if $x \in X$ and $x'\in X'$, then
 $(p\times p')^{-1}(x,x')= p^{-1}(x) \times p'^{-1}(x')$. On morphisms, this follows from the fact that a holonomy map for the fibration, $(p\times p')\colon E\times E' \to X \times X'$, can be obtained by doing the ``product'' of those of $p\colon E \to X$ and $p\colon E' \to X'$, namely,
 $$ I \times (E \times_X X^I) \times (E' \times_{X'} {X'}^I) \ni (t,e,\gamma,e',\gamma') \mapsto \left (\lambda^E(t,e,\gamma), \lambda^{E'}(t,e',\gamma')\right) \in E \times E'.$$
\end{proof}

\begin{Definition}
Let $p\colon M \to B$ be a fibration. Choose a subset, $\x_B$, of $B$. The functor, $$\Fc^M_{\x_B}\colon \pi_1(B,\x_B) \to \CGWH/\simeq,$$   is defined by restricting the functor,  $\Fc^M\colon \pi_1(B,B) \to \CGWH/\simeq$, to $\pi_1(B,\x_B)$,  the full subgroupoid of $\pi_1(B,B)$, with set of objects $\x_B$.

We have, thus, also defined a functor, $\hpi_0\circ \Fc^M_{\x_B}\colon \pi_1(B,\x_B) \to \Sets$.
\end{Definition}

This latter functor, $\hpi_0\circ \Fc^M_{\x_B}\colon  {\pi_1(B,\x_B)}   \to \Sets$, is a key ingredient for the construction of the once-extended Quinn TQFT. The set,  $\x_B$, is a choice of a set of base-points, typically, but not exclusively, for the components of $B$, or, at the other extreme, we could take $\x_B$ to be the set of all the elements of $B$. We can choose. If we pick $\x_B$ finite, this can be used to reduce the groupoids, $\pi_1(B,B)$, and their action (interpreted as a functor) to the more classical setting  of finite group algebras together with categories of bimodules over them.

 We will need such many pointed extensions of quite a few otherwise classical results, which are not that easy to find given in an explicit form in the literature, and so will give them in a bit of detail.

Recall that,
\begin{Definition}A pair, $(X,\x_X)$, of topological spaces is said to be $0$-connected if the set, $\x_X$, has at least one point in each path-component of $X$.
\end{Definition}

\begin{Lemma}\label{lem:colimfib}Let $p\colon M \to B$ be a fibration. Choose a subset, $\x_B$, of $B$ such that $(B,\x_B)$ is 0-connected. We have a natural bijection,
$$F\colon \colim  (\hpi_0\circ \Fc^M_{\x_B}) =\Big(\,\bigsqcup_{x \in \x_B} \hpi_0(M_x) \Big)\bigg/\sim \, \to \hpi_0(M).$$
Given $x \in \x_B$ and $m\in M_x$, this sends the equivalence class of $\PC_m(M_x)$ to $\PC_m(M)$.
\end{Lemma}
\begin{proof}
 By construction, the map, $$F'\colon \bigsqcup_{x \in \x_B} \hpi_0(M_x) \to \hpi_0(M),$$ such that, if  $x \in \x_B$ and $m\in M_x$,  then $\PC_m(M_x)\mapsto \PC_m(M)$, descends to $\colim  (\hpi_0\circ \Fc^M_{\x_B})$. Let us explain this a bit more. Given ${[\gamma]}\colon x\to y$ in $\pi_1(B,\x_B)$, if $m \in M_x$ and $n\in N_y$, and we have that
 $$\big(\hpi_0\circ \Fc^M_{\x_B}([\gamma]\colon x \to y)\big)\big(\PC_m(M_x)\big)=\PC_n(M_y), $$
 then this means that $\lambda^M(1,m,\gamma)$ is in $\PC_n(M_y).$ From this, it follows that $m$ and $n$ are in the same path-component in $M$.  (Note that the $0$-connectedness of $(B,\x_B)$ was not used.)

Now to prove that $F$ is injective. Suppose that, given $x,x'\in \x_B$, and $m\in M_x$ and $m'\in M_{x'}$, we have $\PC_m{(M)}=\PC_{m'}(M)$. Choose any path, $\phi$, in $M$, starting in $m$ and ending in $m'$. Let $\gamma=p\circ \phi$, then, by using  Lemma \ref{pphi-cmpt}, it follows that $\PC_{\Gamma^M_\gamma(m)}{(M)}=\PC_{m'}(M)$, so $$\big((\hpi_0\circ \Fc^M_{\x_B} )\big({[\gamma]}\colon  x\to  y)\big)(\PC_m(M_x) )=\PC_{m'}(M_y).$$
  (Again note that the $0$-connectedness of $(B,\x_B)$ was not used.)
  
 That $F$ is surjective follows  analogously, but here we use the fact that $(B,\x_B)$ is 0-connected. If we are given $m \in M$,  there is a path, $\gamma$, in $B$ connecting  $p(m) \in B$ to some $x \in \x_B$.
 Put $m'=\Gamma^M_\gamma(m)$. We then have $F([\PC_{m'}(M_{x})])=\PC_m(M)$.
\end{proof}


\subsubsection{The path components of pullbacks along fibrations}\label{sec:path-components-pull}

Let $p\colon M \to B$ and $q\colon N \to B$ be  fibrations and consider the  pullback diagram in $\CGWH$ given by the diamond in the diagram below, where we put $P=q\circ {\proj_2}=p\circ {\proj_1}$,
\begin{equation}\label{eq:pull_fibs}\vcenter{
\xymatrix@R=15pt{&&M{\times}_B N\ar[dr]^{\proj_2 } \ar[dl]_{\proj_1} \ar[dd]^P\\&M
\ar[dr]_p && N.\ar[dl]^q\\
&&B&}}
\end{equation}
It is clear, e.g.,  by the universal property of pullbacks, that $P$ is a fibration.

Let $\lambda^M$ and $\lambda^N$ be holonomy maps for the fibrations, $p\colon M \to B$ and $q\colon N \to B$, then a holonomy map, $\lambda^{M\times_B N}$, for the fibration,   $P\colon M\times_B N \to B$, can  be given  such that, if we have a path $x\ra{\gamma} y$ in $B$, and  a point, $(m,n)\in M_x\times N_x  \subseteq  M\times_B N$,
\begin{equation}\label{hol_pull}
 \lambda^{M\times_B N}(t,m,n,\gamma)=\big(\lambda^M(t,m,\gamma),\lambda^N(t,n,\gamma)\big).
\end{equation}

Given $b\in B$, the fibre $P^{-1}(b)$, of $P$ at $b$, is homeomorphic to $M_b \times N_b$.  We have a bijection, from $\hpi_0(M_b) \times \hpi_0(N_{b})$ to $\hpi_0(P^{-1}(b))$. This bijection sends
$\big(\PC_m(M_b),\PC_n(N_b)\big)$ to $\PC_{(m,n)}(P^{-1}(b))$, where $m\in M_b$ and $n \in N_b$.

 We have functors, $\Fc^M,\Fc^N\colon \pi_1(B,B) \to \CGWH/\simeq $, given by the fibrations $p\colon M \to B$ and $q\colon N \to B$.
By construction, and using \eqref{hol_pull}, the functor, $$\Fc^{M\times_B N}\colon \pi_1(B,B) \to \CGWH/\simeq,$$ given by the fibration, $P\colon M\times_B N \to B$,  is naturally isomorphic to the composite of functors below, where $\left \langle\Fc^M,\Fc^N\right \rangle$ is given by the universal property of a product,
$$\pi_1(B,B) \ra{\left \langle\Fc^M,\Fc^N\right \rangle} (\CGWH/\simeq) \times (\CGWH/\simeq) \ra{\Times_\CGWH}\CGWH/\simeq,$$
and where $\Times_\CGWH\colon \CGWH\times \CGWH \to \CGWH$ denotes the product functor in $\CGWH$, which descends to a functor, also denoted $$\Times_\CGWH\colon (\CGWH/\simeq) \times (\CGWH/\simeq) \to \CGWH /\simeq.$$ (Explicitly,  $\Times_\CGWH$ sends a pair, $(X,Y)$, of CGWH spaces to  their product, $X\times Y$, and analogously for maps between spaces.)

The functor,  $\hpi_0\circ \Fc^{M\times_B N}\colon \pi_1(B,B) \to \Sets, $ is, thus, naturally isomorphic to
$$\Times_\Sets \circ \left\langle \hpi_0\circ \Fc^M,\hpi_0\circ \Fc^N \right\rangle\colon \pi_1(B,B) \to \Sets, $$ where $\Times_\Sets\colon \Sets\times \Sets \to \Sets$ is the product functor in $\Sets$.

Given $\x_B \subseteq  B$, it also follows that  the functor, $$\hpi_0\circ \Fc^{M\times_B N}_{\x_B}\colon \pi_1(B,\x_B) \to \Sets ,$$ is naturally isomorphic to
$$\Times_\Sets \circ \left \langle \hpi_0\circ \Fc^M_{\x_B},\hpi_0\circ \Fc^N_{\x_B} \right\rangle\colon \pi_1(B,\x_B) \to \Sets.$$
\begin{Lemma}
Suppose that the pair $(B,\x_B)$ is $0$-connected. There is a bijection,
$$
\colim \big ( \Times_\Sets \circ  \left \langle \hpi_0\circ \Fc^M_{\x_B},\hpi_0\circ \Fc^N_{\x_B} \right \rangle  \big) \rightarrow \hpi_0(M\times_B N).$$
This bijection is such that, if $x \in \x_B$, $m \in M_x$ and $n \in N_x$, then,
                                                                $$\big[\big(\PC_m(M_x),\PC_n(N_x)\big)\big] \mapsto \PC_{(m,n)}(M\times_B N)  . $$
\end{Lemma}
\begin{proof}
Follows from Lemma \ref{lem:colimfib} combined with the previous discussion.
\end{proof}

Any groupoid, $G$, comes with a \emph{contravariant} functor $(-)^{-1}\colon G\to G$, that is the identity on objects and sends each morphism to its inverse. In particular, we have a functor, $$\Times_{\Sets}\circ \big (\hpi_0\circ \Fc^M_{\x_B}\circ (-)^{-1} \times \hpi_0\circ \Fc^N_{\x_B}\big)\colon \pi_1(B,\x_B)^{\op} \times \pi_1(B,\x_B) \to \Sets .$$
This gives us the following. (We note that a generalisation of this lemma, written in the context of $\infty$-groupoids, is in \cite[Lemma 3.8]{Galvezetal}.)
\begin{Lemma}\label{techfib_2}
Let $p\colon M \to B$ and $q\colon N \to B$ be fibrations. Choose $\x_B\subseteq  B$ such that the pair $(B,\x_B)$ is 0-connected. There is a bijection,
$$\int^{x \in \x_B} \big( (\hpi_0\circ \Fc^M)\circ (-)^{-1}(x)\big) \times \big ((\hpi_0\circ \Fc^N)  (x)\big) \to \hpi_0(M\times_B N).$$
Noting that,
\begin{multline*}
\int^{x \in \x_B} \big( (\hpi_0\circ \Fc^M)\circ (-)^{-1}(x)  \big) \times \big ((\hpi_0\circ \Fc^N)  (x)\big)\\=\Big(\bigsqcup_{x \in \x_B}  (\hpi_0\circ \Fc^M)(x) \times  (\hpi_0\circ \Fc^N)(x)\Big) \bigg/ \sim,
\end{multline*}
given $x \in \x_B$, the bijection sends the equivalence class of $\big(\PC_m(M_x),\PC_n(N_x)\big)$ to $\PC_{(m,n)}(M\times_B N)$, where $m\in M_x$ and $n \in N_x$.
\end{Lemma}
\begin{proof}
 This follows from the previous lemma, since any arrow in $\pi_1(B,\x_B)$ is invertible. (Lemma \ref{lem:colimstoprof}, is useful to translate between coends and colimits.)
\end{proof}

\subsubsection{The homotopy content of path-components of pullbacks along fibrations}\label{sec:defT}
The results here will be  crucial to formulate  the once-extended Quinn TQFT.
\begin{Lemma}\label{lem:coher_path_comp}
 Let $f\colon E \to X$ be a fibration, with $E \neq \emptyset$. Let $e\in E$ and  $x =f(e)$. 
 \begin{enumerate}[leftmargin=1cm]
 \item The induced map, $f_e\colon \PC_e(E) \to \PC_x(X)$, is a fibration.
 \item If $k\in f_e^{-1}(x)$, then
 $\PC_k(f_e^{-1}(x))=\PC_k(f^{-1}(x))$.
 \end{enumerate}

 \end{Lemma}
\begin{proof}
The first point is immediate from the homotopy lifting property. For the second point, note that clearly $\PC_k(f_e^{-1}(x))\subseteq  \PC_k(f^{-1}(x))$ as sets. The reverse inclusion also holds. This is because  a path in $f^{-1}(x)$, starting in $k \in   f_e^{-1}(x)$, cannot leave  $\PC_e(E)$, so it is a path in  $f_e^{-1}(x)$.
\end{proof}

Consider two fibrations, $p\colon M \to B$, and $q\colon N \to B$, and the resulting fibration, $P\colon M\times_B N\to B$, of diagram \eqref{eq:pull_fibs}. Suppose that $M,N$ and $B$ are homotopy finite.
\begin{Lemma}
 The space, $M\times_B N $, is homotopy finite.
\end{Lemma}
\begin{proof}This lemma is a particular case of Lemma \ref{main3}.
\end{proof}

Let $b \in B$. By Lemma  \ref{main_techfib_loops}, we have a right action of $\pi_1(B,b)$ on $\hpi_0(M_b) \times \hpi_0(N_b)\cong \hpi_0(P^{-1}(b))$.
Given the form of the holonomy map for the fibration $P\colon M\times_B N\to B$ in  \eqref{hol_pull},
this action is such that, if $\gamma \colon I \to B$ connects $b$ to $b$, then, given $m' \in M_b$ and $n' \in N_b$, so $(m',n')\in P^{-1}(b)$,  we have,
\begin{equation}\label{eq:prop-hol}
\begin{split}
\big(\PC_{m'}(M_b), \PC_{n'}(N_b)\big)\trl [\gamma]&=\big (\PC_{\Gamma^M_\gamma(m')}(M_b), \PC_{\Gamma^N_\gamma(n')}(N_b)\big)\\
&=\big(\PC_{m'}(M_b) \trl [\gamma] , \PC_{n'}(N_b) \trl [\gamma]\big).
\end{split}
\end{equation}
In the last equation, the actions are derived from the fibrations $M \to B$ and $N \to B$.

Now fix $b \in B$, and elements, $m\in M_b$, $n\in N_b$, in the fibres of the two fibrations, $p$ and $q$.
The fibration, $P\colon M\times_B N \to B$, restricts to a map, $$P_{(m,n)}\colon \PC_{(m,n)}(M\times_B N) \to \PC_b(B),$$ which  is a fibration, by Lemma \ref{lem:coher_path_comp}.
 Assuming that $M,N$ and $B$ are homotopy finite, then $P_{(m,n)}\colon \PC_{(m,n)}(M\times_B N) \to \PC_b(B)$ is  a fibration of homotopy finite spaces. In particular, the  fibre of  $P_{(m,n)}\colon \PC_{(m,n)}(M\times_B N) \to \PC_b(B)$ has only a finite number of path-components.
 \begin{Notation}[$T_{(m,n)}^{M\times_B N} $]
We write $T_{(m,n)}^{M\times_B N} $
for the number of path-components of the fibre, at $b\in B$, of the fibration, $P_{(m,n)}\colon \PC_{(m,n)}(M\times_B N) \to \PC_b(B)$.
\end{Notation}
\noindent By definition, or from the long exact sequence \eqref{hles} of $P\colon (M\times_B N) \to B$, we have:

\begin{Lemma}\label{Lem:T_is_orbit} 
The value of $T_{(m,n)}^{M\times_B N} $ is equal to the cardinality of the orbit of the path-component $\PC_{(m,n)}(P^{-1}(b))$ under the right-action of $\pi_1(B,b)$ on $\hpi_0(P^{-1}(b))$.
\end{Lemma}
\begin{proof}
We prove that   $\hpi_0(P_{(m,n)}^{-1}(b)) \subseteq  \hpi_0(P^{-1}(b))$ coincides with the  $\pi_1(B,b)$-orbit, of  $\PC_{(m,n)}(P^{-1}(b))$, inside $\hpi_0(P^{-1}(b))$.

Let $m'\in M_b$ and $n'\in N_b$. If $\PC_{(m',n')}(P^{-1}(b))\in  \hpi_0(P_{(m,n)}^{-1}(b))$,
then, in particular, $(m',n')\in M_b\times N_b \subseteq  M\times_B N$ is in the same path-component as $(m,n)$ in $M\times_B N$.
Choose a path, $\phi$, in $M\times_B N$ connecting $(m,n)$ and $(m',n')$. Applying  Lemma \ref{pphi-cmpt},
it follows that $\PC_{(m,n)}(P^{-1}(b))\trl [p(\phi))]=\PC_{(m',n')}(P^{-1}(b))$. The rest follows by construction.
\end{proof}

\begin{Lemma}\label{Lem:chi_pull}
 Let $b\in B$, $m\in M_b$, and $n \in N_b$, then
$$\chi^\pi\big( \PC_{(m,n)}(M\times_B N)\big)=T_{(m,n)}^{M\times_B N} \,\, \chi^\pi(\PC_b(B))\, \chi^\pi(\PC_m(M_b))\, \chi^\pi(\PC_n(N_b)).$$
\end{Lemma}
\begin{proof}
The decomposition of $P^{-1}(b)= M_b\times N_b$, into path-components, gives a weak homotopy equivalence,
$$\bigsqcup_{(A,A')\in \hpi_0(M_b) \times\hpi_0(N_b)} A\times A' \to M_b\times N_b .$$
Each path component of the fibre $P_{(m,n)}^{-1}(b)$, of $P_{(m,n)}\colon \PC_{(m,n)}(M\times_B N) \to \PC_b(B)$, at $b$, is also a  path-component of  $P^{-1}(b)= M_b\times N_b$.  (We are using Lemma \ref{lem:coher_path_comp}.). \emph{A priori}, however, there may be fewer path-components in {$P^{-1}_{(m,n)}(b)$} than there are in $M_b\times N_b$.

 If we use Lemma \ref{main1}, applied to the HF fibration,   $P_{(m,n)}\colon \PC_{(m,n)}(M\times_B N) \to \PC_b(B)$,
 we obtain
 \begin{align*}
  \chi^\pi\big( \PC_{(m,n)}(M\times_B N)\big)&=\chi^\pi(\PC_b(B))\,\,\chi^\pi\big(P_{(m,n)}^{-1}(b)\big).
 \end{align*}
Now, by Lemma \ref{pphi-cmpt}, all path-components of $P_{(m,n)}^{-1}(b)$ are homotopy equivalent. The result follows from the fact that, by definition, we have $T_{(m,n)}^{M\times_B N} $ such path components, each of which is homotopic to
\begin{align*}
\PC_{(m,n)}\big(P_{(m,n)}^{-1}(b)\big)&=\PC_{(m,n)}\big(P^{-1}(b)\big) = \PC_{(m,n)}(M_b\times N_b)\\ & \cong \PC_m(M_b) \times \PC_n(N_b).
\end{align*}
We now apply Lemma \ref{lem:sums.prods}.
\end{proof}
\subsection{The profunctor construction}\label{proffibspan}
Using the results of the previous sections, we will show that each fibrant span gives us a profunctor with $\Vect$-values, and that the composition of fibrant spans in Definition \ref{main3def} translates under this to composition of profunctors, as described in Subsection \ref{sec:prof_conventions}.

We note that our results are related to those of \cite[8. Cardinality as a functor]{Galvezetal}, which were written in the language of $\infty$-groupoids.
\subsubsection{The profunctor associated to a fibrant span}\label{sec:profunctor_fib_span}
Consider a fibrant span from $X$ to $Y$, then we have a diagram in $\CGWH$ of form,
\begin{equation*}\label{hff2} \xymatrix@R=-5pt{ &M\ar[dl]_p\ar[dr]^{p'}\\   X & &Y,  }
 \end{equation*}
where the induced map, $\big \langle p,p'\big \rangle\colon M \to X\times Y$, is a  fibration. Given $x\in X$ and $y \in Y$, recall, from Notation \ref{space-els}, that we defined the spatial slice at $x$ and $y$ as:
$$\{x| (p,M,p')|y\}\stackrel{abbr.}{=}\{x| M|y\}=\langle p,p'\rangle^{-1}(x,y).$$

We have a holonomy map, $\lambda^M$, for the fibration,  $\big \langle p,p'\big \rangle\colon M \to X\times Y$, of form,
  $$\lambda^M \colon I \times (M \times_{X \times Y} (X \times Y)^I) \to M.$$
  (We are using the notation of \S \ref{sec:more_fibrations}.)
  Let $x,x'\in X$ and $y,y'\in Y$.  Given paths, ${\gamma^X}\colon x\to  x'$ in $ X$, and  $\gamma^Y\colon y\to  y'$
  in $Y$, the holonomy map, $\lambda^M$,  induces a homotopy equivalence,
 \begin{align*} 
  \Gamma_{{\langle \gamma^X,\gamma^Y\rangle}}^M\colon \{x |M |y\} &\to \{x' |M |y'\},\\
                m &\mapsto \lambda^M(1,m,\gamma^X,\gamma^Y).
  \end{align*}
  Here $\langle \gamma^X,\gamma^Y\rangle$ is the path in $X \times Y$ such that $I \ni u \mapsto \big(\gamma^X(u),\gamma^Y(u)\big) \in X \times Y$.

  The homotopy class of $ \Gamma_{{\langle \gamma^X,\gamma^Y\rangle}}^M\colon \{x |M |y\} \to \{x' |M |y'\}$ depends only on the fibration, $\langle p,p'\rangle \colon M\to X\times Y$, and not on the chosen holonomy map. This yields, by Lemma \ref{main_techfibfunctor}, a functor,  \begin{align*}
  \Fc^M\colon \pi_1\big( X\times Y, X\times Y\big) \cong \pi_1(X,X) \times \pi_1(Y,Y) &\to \CGWH/\simeq,
  \end{align*}
where $\Fc^M(x,y):=\{x|M|y\}$ and $\Fc^M\big (x \ra{[\gamma^X]} x', y \ra{[\gamma^Y]} y'\big):= [\Gamma_{(\gamma^X,\gamma^Y)}^M]$.

 For convenience, we will repeat our conventions, from Subsection \ref{sec:prof_conventions},  for profunctors between groupoids.
  Fix a subfield, $\kappa$, of the complex field, $\C$. Recall that $\Vect_k=\Vect$ denotes the category of $\kappa$-vector spaces and linear maps. 

\begin{Definition}[$\Sets$-profunctor and $\Vect$-profunctor] Consider groupoids $G$ and $G'$.
A \emph{$\Sets$-profunctor}, $H\colon G \bto{} G'$, is a functor,$$H\colon G^{\op}\times G' \to \Sets.$$ A \emph{$\Vect$-profunctor}, or \emph{$\Vect$-valued profunctor}, $\Hp\colon G \bto{} G'$, is a functor,  $$\Hp\colon G^{\op}\times G' \to \Vect.$$
\end{Definition}
The free vector space functor, from $\Sets$ to $\Vect$, is  denoted $$\Lin=\Lin_\kappa\colon \Sets \to \Vect.$$ Each $\Sets$-profunctor, $H\colon G^{\op} \times G'\to \Sets$, gives rise to a $\Vect$-valued profunctor, $\Hp:=\Lin\circ H\colon G^\op\times G'\to \Vect$.

\begin{Definition}[The $\Vect$-profunctor associated to a fibrant span]\label{def:fibrantSpanstoprof} Consider a fibrant span, of GCWH spaces, connecting the HF spaces $X$ and $Y$, as below
$$X \ra{(p,M,p')} Y= \left(\vcenter{\xymatrix@R=-5pt{ &M\ar[dl]_p\ar[dr]^{p'}\\   X & &Y}}\right).  $$ Its associated $\Vect$-profunctor, denoted,  $$\Hp(X \ra{(p,M,p')} Y ) \colon \pi_1(X,X) \bto{} \pi_1(Y ,Y ),$$
which we will frequently abbreviate to $$\Hp^M\colon \pi_1(X,X) \bto{} \pi_1(Y ,Y ),$$ 
    is, by definition, the following composite of functors, recalling Definition \ref{main_techfibfunctor},
    \begin{multline*}
    \pi_1(X,X)^\op \times \pi_1(Y,Y)\ra{(-)^{-1} \times \id} \pi_1(X,X) \times \pi_1(Y,Y)\ra{\cong}\pi_1(X\times Y,X\times Y)\\\ra{\Fc^M} \CGWH/\simeq \ra{\hpiz} \Sets \ra{\Lin} \Vect.
    \end{multline*}
    \end{Definition}
    
     Taking this apart,  given  $x \in X$ and $y\in Y $, $\Hp^M(x,y)$ is the free vector space over $\hpi_0(\{x|M|y\})$, the set of path-components of the fibre of $\langle p,p' \rangle\colon M \to X\times Y$, at $(x,y)$.  Given morphisms,  in $\pi_1(X,X)$ and $\pi_1(Y ,Y )$, say
$$[\gamma^X]\colon x\to x' \quad \textrm{ and } \quad [\gamma^Y]\colon y\to y',$$ the linear map, $$\Hp^M\big([\gamma^X],[\gamma^Y]\big)\colon \Hp^M(x',y) \to \Hp^M(x,y'),$$ is induced by the homotopy equivalence, between fibres,
$$\Gamma_{{\langle \overline{\gamma^X},\gamma^Y\rangle}}^M\colon \{x'|M|y\} \to \{x|M|y'\},$$
by applying $\hpi_0\colon \CGWH \to \Sets$, and then $\Lin\colon \Sets \to \Vect$.
Here $\langle \overline{\gamma^X},\gamma^Y\rangle$ is the path in $X \times Y$, such that $ I \ni u \mapsto \big(\gamma^X(1-u),\gamma^Y(u)\big) \in X \times Y$.

The following result will implicitly be used a number of times.
\begin{Lemma}\label{lem:finiteness} Suppose that the fibrant span,  $(p,M,p')\colon X  \to Y$, is homotopy finite, so  that $X,Y$ and $M$ are homotopy finite spaces, then the profunctor, $$\Hp^M\colon \pi_1(X,X) \bto{} \pi_1(Y ,Y ),$$ is a 1-morphism in the bicategory, $\vProfGrphf$, of homotopy finite groupoids and $\Vect$-profunctors between them; see \S \ref{sum:conv_prof}.

Given $x \in X$ and $y \in Y$, the vector space, $\Hp^M(x,y)$, is finite dimensional.
     \end{Lemma}
\begin{proof} Since $X$ and $Y$ are homotopy finite, it follows that the groupoids, $\pi_1(X,X)$ and $\pi_1(Y,Y)$, each are homotopy finite. Given $x \in X$ and $y\in Y$, then by Lemma \ref{2of3}, the fibre $\{x|M|y\}=\langle p,p' \rangle^{-1}(x,y)$ is homotopy finite, and  thus it only has a finite number of path-components. 
\end{proof}

\begin{Notation}\label{not_barH}
More generally, choose subsets, $\x_X$ and $\y_{Y }$, of $X$ and $Y $, respectively. The restriction, of $\Hp^M\colon \pi_1(X,X)^{\op}\times  \pi_1(Y ,Y )\to \Vect$,  to $\pi_1(X,\x_X)\times \pi_1(Y ,\y_{Y })$,
is a $\Vect$-profunctor,  $$\pi_1(X,\x_X)\bto{} \pi_1(Y,\y_{Y }).$$

We will use three different notations for it:
$$\Hpb_{(\x_X,\y_Y)}(p,M,p')\colon \pi_1(X,\x_X)^{\op}\times \pi_1(Y ,\y_{Y })\to \Vect,
$$ 
$$\Hpb_{(\x_X,\y_Y)}\Big( X  \ra{ (p,M,p')} Y \Big) \colon \pi_1(X,\x_X)^{\op}\times \pi_1(Y ,\y_{Y })\to \Vect,
$$ 
and finally,
$$\Hpb^M_{(\x_X,\y_{Y })}\colon \pi_1(X,\x_X)^{\op}\times \pi_1(Y ,\y_{Y })\to \Vect,
$$ depending on the context and the amount of detail needed. 
\end{Notation}

We will be mainly interested in the case when \smash{$X \ra{(p,M,p')} Y$} is homotopy finite, and furthermore both sets, $\x_X$ and $\y_Y$, are finite. In this case, we therefore have that $\Hpb^M_{(\x_X,\y_{Y })}\colon \pi_1(X,\x_X)\bto  \pi_1(Y ,\y_{Y })$ is a 1-morphism in the bicategory $\vProfGrpfin$, of finite groupoids and $\Vect$-profunctors between them.

\subsubsection{The symmetric monoidal-like structure of $\Hpb^M_{(-,-)}$}\label{sec:sim_mon_0}

The following result will be used when proving that  the constructions of once-extended Quinn TQFTs, given here,   do indeed give bifunctors, which are symmetric monoidal.
\begin{Lemma}\label{Sym_mon:1}
 Consider two fibrant spans of homotopy finite spaces, $$(p,M,q)\colon X \to Y \textrm{ and } (p',M',q')\colon X' \to Y'.$$ Let $\x_X {\subseteq}  X$, $\overline{x}'_{X'} {\subseteq}  X'$,   $\y_Y {\subseteq}  Y$ and  $\overline{y}'_{Y'}{\subseteq}  Y' $. Form the product HF fibrant span, $$(p\times p',M\times M',q \times q')\colon X\times X' \to Y \times Y'.$$
There is a natural isomorphism from the profunctor,
 \begin{multline*}
 \Hpb_{(\x_X\times \overline{x}'_{X'} ,\y_Y\times \overline{y}'_{Y'})}
 \big ( p\times p',M\times M',q \times q'\big)
 \colon\\  \pi_1(X\times X',\x_X\times \overline{x}'_{X'})^{\op}\times \pi_1(Y \times Y',\y_{Y } \times \overline{y}'_{Y' })\to \Vect,
 \end{multline*}
 to the profunctor obtained from the following composition of functors,
\begin{multline*}
 \pi_1(X\times X',\x_X\times \overline{x}'_{X'})^{\op}\times \pi_1(Y \times Y',\y_{Y } \times \overline{y}'_{Y' })\\ \ra{\cong} \big( \pi_1(X,\x_X)^{\op} \times \pi_1(Y,\y_Y)\big) \times \big( \pi_1(X',\overline{x}'_{X'})^{\op}\times \pi_1(Y,\overline{y}'_{Y'})\big) \\ \ra{  \Hpb_{(\x_X,\y_Y)}(p,M,q) \times  \Hpb_{(\overline{x}'_{X'},\overline{y}'_{Y'})}(p',M',q') }\Vect\times \Vect \ra{\Otimes_\Vect} \Vect.
\end{multline*}
Let $x \in \x_X$, $x' \in \overline{x}'_{X'}$,   $y\in \y_Y$  and $y'\in \overline{y}'_{Y'}$. This natural isomorphism is such that, if $m \in  \{x|M|y\}$ and $m' \in \{x'|M'|y'\}$, then
$$\PC_m\big (\{x|M|y\}\big )\otimes \PC_{m'}\big (\{x'|M'|y'\}\big ) {\longleftrightarrow} \PC_{(m,m')} \big ( \{(x,x')| M\times M'| (y,y')\}\big).$$
\end{Lemma}
\begin{proof}
 This follows from Lemma \ref{Lem:FExE'}.
\end{proof}

In order to prove that the once-extended Quinn TQFT, in its various forms, gives a symmetric monoidal bifunctor, it is convenient to change slightly the language of the previous result, approximating that of  Definition \ref{Sym. Mon. Bifunctor}. (This point will be made concrete later, in Subsection \ref{sec:sym_mon}.)
Consider the canonical natural isomorphisms, of groupoids,
\begin{align*}
m_{(X,X')}\colon \pi_1(X,\x_X)\times \pi_1(X',\overline{x}'_{X'}) &\to \pi_1(X\times X',\x_X\times \overline{x}'_{X'}),\\
m_{(Y,Y')}\colon \pi_1(Y,\y_Y)\times \pi_1(Y',\overline{y}'_{Y'})& \to \pi_1(Y\times Y',\y_Y\times \overline{y}'_{Y'}).
\end{align*}
By using the notation of Example \ref{Ex:maps_to_profunctors}, they yield profunctors,
\begin{align*}
\varphi^{m_{(X,X')}}\colon \pi_1(X,\x_X)\times \pi_1(X',\overline{x}'_{X'}) &\bto \pi_1(X\times X',\x_X\times \overline{x}'_{X'}),\\
\varphi^{m_{(Y,Y')}}\colon \pi_1(Y,\y_Y)\times \pi_1(Y',\overline{y}'_{Y'})& \bto \pi_1(Y\times Y',\y_Y\times \overline{y}'_{Y'}).
\end{align*}

Continuing the notation of Lemma \ref{Sym_mon:1}, we have
\begin{Lemma}\label{def_chi}
There is  a diagram of $\Vect$-profunctors, and a natural isomorphism, $\chi_{(M,M')}$, or in full $\chi_{( (p,M,q),(p',M',q'))}$, of $\Vect$-profunctors, as shown below,
$$
\hskip-2.6mm\xymatrixcolsep{2.5cm}\xymatrix{\pi_1(X,\x_X)\times \pi_1(X',\overline{x}'_{X'})\ar[rr]^{\Hpb_{(\x_X,\y_Y)}(p,M,q) \otimes  \Hpb_{(\overline{x}'_{X'},\overline{y}'_{Y'})}(p',M',q')  }\ar[d]_{ \varphi^{m_{(X,X')}}} \drrtwocell<\omit>{\hspace*{1cm}\chi_{(M,M')}}&& \pi_1(Y,\y_Y)\times \pi_1(Y',\overline{y}'_{Y'})\ar[d]^{\varphi^{m_{(Y,Y')}}}\\
 \pi_1(X \times X',\x_X\times \overline{x}'_{X'})\ar[rr]_{\Hpb_{(\x_X\times \overline{x}'_{X'} ,\y_Y\times \overline{y}'_{Y'})}\left  (p\times p',M\times M',q \times q'\right)} &&  \pi_1(Y \times Y',\y_Y\times \overline{y}'_{Y'}).}
$$
Let $x \in \x_X$, $x' \in \overline{x}'_{X'}$,   $(y_1\ra{\gamma} y_2)\in \pi_1(Y,\y_Y)$  and $(y'_1\ra{\gamma'} y'_2)\in \pi_1(Y,\y_Y)$. This natural isomorphism is such that, referring to the notation in Equation \eqref{Profunctor composition}, if $m \in  \{x|M|y_1\}$ and $m' \in \{x'|M'|y'_1\}$, then the equivalence class of
$$
\big(\PC_m\big (\{x|M|y_1\}\big )\otimes \PC_{m'}\big (\{x'|M'|y'_1\}\big ) \big) \otimes m_{(Y,Y')}(\gamma , \gamma'),
$$
is sent to the equivalence class of
\begin{multline*}
\id_{(x,x')} \otimes\\ \Hpb_{(\x_X\times \overline{x}'_{X'} ,\y_Y\times \overline{y}'_{Y'})}^{M\times M'}\big(\id_{(x,x')},m_{(Y,Y')}(\gamma,\gamma')\big)\big(
\PC_{(m,m')} \big ( \{(x,x')| M\times M'| (y_2,y'_2)\}\big)\big).
\end{multline*}
\end{Lemma}
\begin{proof} This follows from the previous lemma, and elementary properties of profunctors.
\end{proof}

\subsubsection{The $\Hp$  construction preserves the composition of fibrant spans}\label{sec:IsHfunctor}

Now we know how to construct profunctors and $\Vect$-profunctors from fibrant spans, we should ask how that construction behaves with respect to composition of fibrant spans, and also what does it do to identity spans.  We will examine preservation of composition here, whilst preservation of identities\label{preservation of identities} will be  discussed later, being a consequence of  Lemma \ref{preserves unitors}, as it is more convenient to package it with similar results later on.

Consider HF fibrant spans, $$(p_1,M_1,p_1')\colon X \to Y \textrm{ and } (p_2,M_2,p_2')\colon Y \to Z.$$ Recall the definition of the composition ({see Lemma \ref{main3} and Definition \ref{main3def}),
 $$(p_1,M_1,p_1')\bullet (p_2,M_2,p_2')=(\overline{p_1},M_1\times_Y M_2,\overline{p_2'})\colon X\to Z,$$
 which is,  itself, a HF fibrant span. 
 To recall and extend notation from earlier, we repeat  the relevant commutative diagram, in Equation \eqref{compHFfibrations},
in which the middle diamond is a pullback,
$$\xymatrix@R=10pt{ &&& 
M_1{\times}_Y M_2\ar[dd]^P  \ar[dl]_{}\ar[dr]^{} \ar@/_1.5pc/[ddll]_{\overline{p_1}} \ar@/^1.5pc/[ddrr]^{\overline{p_2'}}\\ &&M_1\ar[dl]^{p_1}\ar[dr]_{p_1'}
&&M_2\ar[dl]^{p_2}\ar[dr]_{p_2'}\\  & X & &Y  & &Z, } 
$$

We recall, from Lemma \ref{main3}, that $\langle \overline{p_1},P,\overline{p'_2}\rangle\colon M_1\times_Y M_2 \to X\times Y \times Z$ is a fibration. We also note that, given holonomies, for $\langle p_1,p_1' \rangle \colon M_1 \to X\times Y$ and for $\langle p_2,p_2' \rangle\colon  M_2 \to Y\times Z$, denoted\footnote{Note that $(X\times Y)^I\cong X^I\times Y^I$ and  $(Y\times Z)^I\cong Y^I\times Z^I$, canonically.},
\begin{align*}
\lambda^{M_1}&\colon I\times \big ( M_1\times_{X\times Y} (X^I\times Y^I) \big) \to M_1,\\ \intertext{and}
\lambda^{M_2}&\colon I\times \big (M_2\times_{Y\times Z} (Y^I\times Z^I)\big) \to M_2,
\end{align*} 
(respectively), then a holonomy, $$\lambda^{M_1\times_Y M_2}\colon I \times \big ( (M_1\times_YM_2)\times_{X \times Y \times Z} ( X^I\times Y^I \times Z^I)  \big) \to M_1\times_YM_2,$$ for  the fibration $\langle \overline{p_1},P,\overline{p'_2}\rangle\colon M_1\times_Y M_2 \to X\times Y \times Z$, is obtained from the holonomies, $\lambda^{M_1}$ and $\lambda^{M_2}$, in the obvious way, namely,
$$\big(t,(m_1,m_2),(\gamma^X,\gamma^Y,\gamma^Z)\big)\mapsto \Big ( \lambda^{M_1}\big(t,m_1, (\gamma^X,\gamma^Y)\big), \lambda^{M_2}\big(t,m_2, (\gamma^Y,\gamma^Z)\big)\Big).$$

The following result shows that the profunctor construction is compatible with composition of fibrant spans.

 \begin{Proposition}\label{techfib3}
Choose subsets, $\x_X$, $\y_Y$ and $\z_Z$, of $X,Y$ and $Z$, respectively. Suppose that $(Y,\y_Y)$ is $0$-connected.
We have a canonical isomorphism of $\Vect$-profunctors from $\pi_1(X,\x_X)$ to $\pi_1(Z,\z_Z)$,
\begin{equation*}
\eta^{M_1,M_2}_{(\x_X,\y_Y,\z_Z)}\colon \Hpb_{(\x_X,\y_Y)}^{M_1}\bullet \Hpb^{M_2}_{(\y_Y,\z_Z)}
\To \Hpb_{(\x_X,\z_Z)}^{M_1\times_Y M_2},
\end{equation*}
 such that, if $x \in \x_X$ and $z \in \z_Z$, then, given any $y\in \y_Y$, $m_1 \in \{x|M_1|y\}$ and $m_2 \in \{y|M_2|z\}$, we have that the linear map,
$$(\eta^{M_1,M_2}_{(\x_X,\y_Y,\z_Z)})_{(x,z)}\colon \int^{y \in \y_Y} \Hpb_{(\x_X,\y_Y)}^{M_1}(x,y) \otimes \Hpb_{(\y_Y,\z_Z)}^{M_2}(y,z){\to} \Hpb_{(\x_X,\z_Z)}^{M_1\times_Y M_2}(x,z),
$$
 (using the notation in Equation \eqref{Profunctor composition}),
sends the equivalence class of
$$\PC_{m_1}(\{x|M_1|y\}) \otimes \PC_{m_2}(\{y|M_2|z\}) \in  \bigoplus_{y \in \y_Y} \Hpb_{(\x_X,\y_Y)}^{M_1}(x,y) \otimes \Hpb_{(\y_Y,\z_Z)}^{M_2}(y,z) $$
to 
$\PC_{(m_1,m_2)}(\{x| M_1\times_Y M_2|z\}).$
\end{Proposition}
\begin{proof}
This follows  by combining the previous discussion with Lemma \ref{techfib_2}.
\end{proof}
Note that, in a situation in which we have two pairs of composable fibrant spans, say, $$(p_1,M_1,p_1')\colon X \to Y \textrm{ and } (p_2,M_2,p_2')\colon Y \to Z,$$ and also, $$(\hat{p_1},\hat{M_1},\hat{p_1'})\colon \hat{X} \to \hat{Y} \textrm{ and } (\hat{p_2},\hat{M_2},\hat{p_2'})\colon \hat{Y} \to \hat{Z},$$
then the natural isomorphisms in Lemma \ref{techfib3} are compatible with those of Lemma  \ref{def_chi}.
We leave it to the reader to write down the corresponding commutative diagram of natural transformations between profunctors.
\subsection{HF fibrant resolved 2-spans connecting fibrant spans}\label{sec:fibrant_resolved_2spans}
So far we have, in the main,  been handling only the 1-categorical structure related to fibrant spans. We now introduce `spans between spans', that is `2-spans', or, as we will call them `windows'. We will later be investigating how this second level of structure on the spatial side is reflected, via the profunctor construction, in the `linear algebra', and, of course, this is the beginning of the extension of Quinn's theory.
\subsubsection{HF fibrant windows}

\begin{Definition}[Window]\label{def:window}
 By a  window, $\W$, we will mean a diagram of, as usual, CGWH spaces of the form below, so $\W$ is a `span of spans',
 \begin{equation}\label{window}
\W=\vcenter{\xymatrix@R=18pt{
  X  &  M\ar[l]_p\ar[r]^{p'}  & Y\\
  Z\ar[u]^{P_l} \ar[d]_{Q_l} & L \ar[l]_{l}\ar[r]^{r}\ar[u]^P\ar[d]_Q& W \ar[u]_{P_r} \ar[d]^{Q_r}\\
  X'  &  N\ar[l]^q\ar[r]_{q'}  & Y'.}}
  \end{equation}

The \emph{boundary}, $\bd{(\W)}$, of the above window,  is the following diagram,
 \begin{equation}\label{boundary}
\bd(\W)=\vcenter{\xymatrix@R=18pt{
  X  &  M\ar[l]_p\ar[r]^{p'}  & Y\\
  Z\ar[u]^{P_l} \ar[d]_{Q_l} & & W \ar[u]_{P_r} \ar[d]^{Q_r}\\
  X'  &  N\ar[l]^q\ar[r]_{q'}  & Y'.}}
  \end{equation}
  By the \emph{frame}, $\fr(\W)$, of the window, $\W$, above,  we will mean  the following limit,
\begin{equation}\label{frame}
 \fr(\W)=  M  \underset{X \times Y}{\times} (Z \times W)  \underset{X' \times Y'}{\times} N\cong \lim(\bd(\W)),
\end{equation}
and the \emph{filler}, $P_L$, of the window, $\W$, is given by the naturally defined map,
\begin{equation}\label{filler}
P_L\colon L \to \fr(\W).
\end{equation}

The restrictions of the diagram, \eqref{window}, to each of its  four boundary spans, are called the \emph{top}, \emph{bottom}, \emph{left} and \emph{right boundary spans} of $\W$. These are:
\begin{equation}\label{diag:bound-spans}
\vcenter{\xymatrix@R=5pt@C=12pt{ &M\ar[dl]_{p} \ar[dr]^{p'} \\
X && Y,}} \quad
\vcenter{\xymatrix@R=5pt@C=12pt
{ &N\ar[dl]_{q} \ar[dr]^{q'} \\
X' && Y',}}\quad
\vcenter{\xymatrix@R=5pt@C=12pt
{ &Z\ar[dl]_{P_l} \ar[dr]^{Q_l} \\
X && X', }}\quad
\vcenter{\xymatrix@R=5pt@C=12pt
{ &W\ar[dl]_{P_r} \ar[dr]^{Q_r} \\
Y && Y'.}}
\end{equation}
The \emph{middle horizontal} and \emph{middle vertical} spans of the window $\W$, are defined as:
\begin{equation}\label{diag:middle-spans}\vcenter{\xymatrix@R=5pt@C=12pt
{ &L\ar[dl]_{P} \ar[dr]^{Q} \\
M && N}} \qquad \textrm{ and } \qquad  \vcenter{\xymatrix@R=5pt@C=12pt
{ &L\ar[dl]_{l} \ar[dr]^{r} \\
Z && W.}}
\end{equation}
\end{Definition}
In this paper, we will see windows as being `oriented' from top to bottom and from left to right.

\begin{Definition}[HF fibrant window]\label{def:HF-fibrant-window}
A \emph{fibrant window} is a window, $\W$, as in \eqref{window}, such that:
\begin{enumerate}[leftmargin=1cm]
 \item the filler, $P_L\colon L \to \fr(\W)$, is a fibration,
 \item[]\hspace{-12.5mm} and, 
\item the four boundary spans (top, bottom, left and right) of $\W$ are all fibrant.
\end{enumerate}

If, in addition, all the spaces appearing in diagram \eqref{window} are HF, then the window,  $\W$, will be  called a \emph{HF fibrant window}.
\end{Definition}
The following are some immediate consequences of the definition of HF fibrant windows.
\begin{enumerate}[leftmargin=1cm]
 \item Suppose that a fibrant window, $\W$, as in \eqref{window}, is HF, then its frame   $$\fr(\W)=M  \underset{X \times Y}{\times} (Z \times W)  \underset{X' \times Y'}{\times} N$$ is a HF space. This follows by  applying Lemma \ref{pull_HF} to the pullback diagram,
 \begin{equation}\label{crucial:pullback}\vcenter{
 \xymatrix@R=18pt{ &
 M  \underset{X \times Y}{\times} (Z \times W)  \underset{X' \times Y'}{\times} N\ar[rr]\ar[d]  && M\times N \ar[d]^{\langle p,p'\rangle\times \langle q,q'\rangle} \\
 &Z\times W\ar[rr]_<<<<<<<<<<<<<<<<<{\tau\circ (\langle P_l,Q_l\rangle\times \langle P_r,Q_r\rangle ) }&& (X\times Y) \times (X'\times Y'),
 }}
 \end{equation}
 where $\tau$ is the obvious transposition. This uses that the top and bottom boundary spans, of $\W$,  are fibrant.
\item There are  two naturally defined maps, $\fr(\W) \to Z\times W$ and $\fr(\W)\to M \times N$.
Both are fibrations. For the first map,  $\fr(\W) \to Z \times W,$ this follows from the pullback diagram \eqref{crucial:pullback} above together with Lemma \ref{pull_HF}, and similarly for the  map $\fr(\W)\to M \times N$, by symmetry.
\item Composing with the filler, $P_L\colon L \to \fr(\W)$, of $\W$, which by definition is  a fibration, this gives:
 \end{enumerate}

\begin{Lemma}\label{lem:midhor_is_fib} If a window, $\W$,  is fibrant, then so are the middle horizontal and middle vertical spans, in Equation \eqref{diag:middle-spans}.\end{Lemma}

 \subsubsection{Isomorphic  windows and equivalent  fibrant windows}
 \begin{Definition}[Isomorphic windows and equivalent fibrant windows]\label{Def:iso_windows}
  Two windows, $\W_1$ and $\W_2$, as below, so with the same boundary,
    \begin{equation}
 \W_1=\vcenter{\xymatrix@R=18pt{
  X  &  M\ar[l]_p\ar[r]^{p'}  & Y\\
  Z\ar[u]^{P_l} \ar[d]_{Q_l} & L_1 \ar[l]_{l_1}\ar[r]^{r_1}\ar[u]^{P_1}\ar[d]_{Q_1}& W \ar[u]_{P_r} \ar[d]^{Q_r}\\
  X'  &  N\ar[l]^q\ar[r]_{q'}  & Y',}}
 \textrm{ and }\quad
\W_2=\vcenter{\xymatrix@R=18pt{
  X  &  M\ar[l]_p\ar[r]^{p'}  & Y\\
  Z\ar[u]^{P_l} \ar[d]_{Q_l} & L_2 \ar[l]_{l_2}\ar[r]^{r_2}\ar[u]^{P_2}\ar[d]_{Q_2}& W \ar[u]_{P_r} \ar[d]^{Q_r}\\
  X'  &  N\ar[l]^q\ar[r]_{q'}  & Y',}}
  \end{equation}
and thus, in particular, $\fr(\W_1)=\fr(\W_2)$, are said to be \emph{isomorphic} if there exists a homeomorphism, $F\colon L_1 \to L_2$, making the obvious three dimensional diagram commute.
  This is equivalent to saying that there exist maps, $F\colon L_1 \to L_2$ and $F'\colon L_2 \to L_1$, making the diagrams below commute,
  $$\vcenter{\xymatrix@R=10pt{&L_1\ar[rr]^{F}\ar[dr]_{P_{L_1}} & &L_2,\ar[dl]^{P_{L_2}} \\ && \fr(\W_2)}} \xymatrix{ \quad \textrm{and } \quad }\hskip-1cm
  \vcenter{\xymatrix@R=10pt{&L_1\ar@{<-}[rr]^{F'}\ar[dr]_{P_{L_1}} & &L_2,\ar[dl]^{P_{L_2}} \\ && \fr(\W_2)}}
  $$
  such that $F\circ F'=\id_{L_2}$ and $F'\circ F=\id_{L_1}$.

  More generally, if $\W_1$ and $\W_2$ are fibrant, then $\W_1$ and $\W_2$ are called \emph{equivalent} if there exist  $F\colon L_1 \to L_2$ and $F'\colon L_2 \to L_1$ making the diagrams above  commute, together with fibred homotopies (see Subsection \ref{sec:fhe}),
  $$F'\circ F\Ra[\fr(\W_2)]{H_1}\id_{L_1} \quad \textrm{ and } \quad F\circ F'\Ra[\fr(\W_2)]
  {H_2}\id_{L_2}.$$
 \end{Definition}
\noindent It is easy to see that equivalence between  fibrant windows is an equivalence relation.

 \subsubsection{HF fibrant resolved 2-spans}\label{sec:HF_f_r_2_s} We recall the notation in item \eqref{itemref:X^I} in page \pageref{itemref:X^I}, and also Definition \ref{idsHFb}.
\begin{Definition}[HF fibrant resolved 2-span]
Given HF fibrant spans, from $X$ to $Y$, namely $(p,M,p'), (q,N,q')\colon X \to Y$,
a \emph{HF fibrant resolved 2-span}, $$\W=(l_X,P,L,Q,r_Y)\colon (p,M,p')\Longrightarrow (q,N,q'),$$
from $(p,M,p')$ to  $(q,N,q')$,
also written
 $$\xymatrix{\mort{X}{Y,}{(p,M,p')}{(q,N,q')}{\W}}
 $$
 is a HF fibrant window of  form, as below
 \begin{equation}\label{2Morphism}
  \W=\vcenter{\xymatrix@R=18pt{
  X  &  M\ar[l]_p\ar[r]^{p'}  & Y\\
  X^I\ar[u]^{s_X} \ar[d]_{t_X} & L \ar[l]_{l_X}\ar[r]^{r_Y}\ar[u]^P\ar[d]_Q& Y^I \ar[u]_{s_Y} \ar[d]^{t_Y}\\
  X  &  N\ar[l]^q\ar[r]_{q'}  & Y.}}
  \end{equation}
   This means that $X,Y, M,N$, and $L$  are  HF spaces, (and, thus, so are $X^I$ and $Y^I$), and the filler of $\W$ (see Definition \ref{def:window}),  below, is a fibration,
 \begin{equation}\label{fib_bound}
 L \stackrel{P_L}{\longrightarrow} \fr(\W)= M  \underset{X \times Y}{\times} (X^I \times Y^I)  \underset{X \times Y}{\times} N.\end{equation}
 We will  frequently identity a  HF-fibrant resolved 2-span with its filler,  $L \stackrel{P_L}{\longrightarrow} \fr(\W)$.
\end{Definition}
\begin{Remark}Note that the top and bottom boundary spans in \eqref{2Morphism} are already fibrant, by assumption. The left and right boundary spans are also fibrant, by construction, see Example \ref{ids_ex}.
Crucially for our constructions later on, the middle horizontal and middle vertical spans in \eqref{2Morphism} are fibrant, by Lemma \ref{lem:midhor_is_fib}.
\end{Remark}

 The `resolved' terminology arises from the fact that we allow the left and right boundaries of a 2-span to take values in the respective path spaces.  For 2-spans in a usual setting, the left and right vertical edges would be identity spans, but to ensure the result is fibrant, we  must `resolve' those vertical edges, replacing them with their `fibrant replacements', as in  Example \ref{ids_ex}.  This is in line with the definition of extended cobordisms between cobordisms of manifolds. Our definition was also designed  so that the horizontal composition of resolved 2-spans, defined in Subsection \ref{sec:hor_comp_2spans}, is a resolved 2-span.

\begin{Definition}[Equivalent and isomorphic HF fibrant resolved 2-spans]\label{Def_iso-2-spans} Let $X$ and $Y$ be HF spaces. Two HF resolved 2-spans, $$\W_1,\W_2\colon \big ( (p,M,p')\colon X \to Y) \implies \big ((q,N,q')\colon X \to Y),$$ will be said to be \emph{equivalent} if they are equivalent as HF fibrant windows, and, similarly, $\W_1$ and $\W_2$ are \emph{isomorphic} if they are isomorphic as   HF fibrant  windows.  
\end{Definition}


\subsection{HF resolved fibrant 2-spans and natural transformations of profunctors}
The next step is to study these HF resolved fibrant 2-spans before turning to how their properties are reflected by the profunctor construction.
\subsubsection{The spatial 2-slices of a HF resolved fibrant 2-span }\label{sec:2-els}
Let $X$ and $Y$ be HF-spaces. Consider HF fibrant spans, $(p,M,p'),(q,N,q')\colon X \to Y,$
and a  HF fibrant resolved 2-span,
connecting them, $$\W=(l_X,P,L,Q,r_Y)\colon (p,M,p')\Longrightarrow (q,N,q').$$  Its underlying HF fibrant window, $\W$, is the commutative diagram of solid arrows in Equation \eqref{diag:W}, just below.  (The dashed arrows showing the inclusion, $x\mapsto \const_x$, of a space, $X$, into the corresponding path space, $X^I$, via constant paths, do not necessarily commute with the rest of the diagram. They will, however, have an important role later.)
 \begin{equation}\label{diag:W}
 \W= \,\, \vcenter{ \xymatrix{
  X \ar@{-->}@/_2pc/[d]_{\const}  &  M\ar[l]_p\ar[r]^{p'}  & Y \ar@{-->}@/^2pc/[d]^{\const} \\
  X^I\ar[u]^{s_X} \ar[d]_{t_X} & L \ar[l]_{l_X}\ar[r]^{r_Y}\ar[u]^P\ar[d]_Q& Y^I \ar[u]_{s_Y} \ar[d]^{t_Y}\\
  X \ar@{-->}@/^2pc/[u]^{\const}  &  N\ar[l]^q\ar[r]_{q'}  & Y \ar@{-->}@/_2pc/[u]_{\const}}}
  \end{equation}

  Let $x\in X$ and $y\in Y$.  Recall, from Notation \ref{space-els}, that the spatial slices
 of the fibrant spans, $(p,M,p')$ and  $(q,N,q')$,  are defined as
  \begin{align*}
  \{x| M| y\}&=\langle p,p' \rangle^{-1}(x,y)=\{m \in M: p(m)=x \textrm{ and } p'(m)=y\},\\
  \{x| N| y\}&=\langle q,q' \rangle^{-1}(x,y)=\{n \in N: q(n)=x \textrm{ and } q'(n)=y\},
  \end{align*}
 and these spatial slices will be homotopy finite; see Lemma \ref{all-HF}.
  
  \begin{Definition}[Spatial 2-slices]\label{space-2els} Let $\W=(l_X,P,L,Q,r_Y)\colon (p,M,p')\Longrightarrow (q,N,q')$ be a HF resolved 2-span,  as in Equation \eqref{diag:W}. Let $x \in X$, $y \in Y$. We define\footnote{Note the square brackets rather than braces here. That distinction will be needed shortly.} the following space, which we call the \emph{spatial 2-slice},  of $\W$, at $(x,y)$,
  \begin{equation}\label{eq:s2slxy}
   \begin{split}
  [ x | L | y]&:=\langle  l_X , r_Y\rangle^{-1}(\const_x,\const_y)
  \\&= \big\{l \in L:  l_X(l)=\const_x, r_Y(l)=\const_y\big\}.
  \end{split}
  \end{equation}
Given $m \in  \{x| M| y\}$ and $n\in  \{x| N| y\}$,  also consider the following space, which we call the \emph{spatial 2-slice}, of $\W$, at $(m,x,y,n)$,
\begin{equation}\label{eq:s2slxymn}
\left[\hskip-1pt\hskip-1pt\begin{array}{c|c|c}
&m&\\
x&L & y\\
&n&\end{array}\hskip-1pt\right]:=P_{L}^{-1}(m,\const_x,\const_y,n)
\end{equation}

More generally, given paths,
$\gamma^X\colon x \to x'$ and $\gamma^Y\colon y \to y'$,
in $X$ and $Y$, and elements, $m \in \{x| M |y\}$ and $n\in \{x'| N |y'\}$, we define the following space, also called a spatial 2-slice of $\W$, but at $(m,\gamma^X,\gamma^Y,n)$,
\begin{equation}\label{eq:s2slgamma}
\left[\hskip-1pt\begin{array}{c|c|c}
&m&\\
\gamma^X&L & \gamma^Y\\
&n&\end{array}\hskip-2pt\right]  :=P_{L}^{-1}(m,\gamma^X,\gamma^Y,n).
\end{equation}
\end{Definition}

 In the context of this definition, note that, given $x \in X$ and $y\in Y$, the maps, $P\colon L \to M$ and $Q\colon L \to N$,  canonically  restrict to maps, as below,
 \begin{align*}
 P_{x,y} \colon  \left [ x \left | L  \right |y\right ] \to  \{x| M| y\} &&\textrm{ and }\quad\qquad
 Q_{x,y} \colon  \left [ x \left |  L \right |y\right ] \to  \{x| N| y\}.
 \end{align*}
Moreover, given $m \in \{x|M|y\}$ and $n\in \{x|N|y\}$, we  have the following,
\begin{equation}\label{eq:Mels-vertspan}\left[\hskip-1pt\begin{array}{c|c|c}
&m&\\
x&L & y\\
&n&\end{array}\hskip-1pt\right]  = \langle P_{x,y},Q_{x,y}\rangle^{-1}(m,n).
\end{equation}

\begin{Lemma}\label{lem:matrix2els}
Let $x\in X$, and $y \in Y$. \begin{enumerate}[leftmargin=1cm]
 \item  The induced map, below, is a fibration, $$\langle P_{x,y},Q_{x,y}\rangle \colon   \left [ x \left |L \right |y\right ]  \to \{x |M|y\} \times \{x|N|y\}.$$
 \item\label{item:independence} If  $m \in \{x|M|y\}$ and $n\in \{x|N|y\}$, the homotopy type of the spatial 2-slice of  $\W$, at $(m,x,y,n)$ in \eqref{eq:s2slxymn},
depends only on the path-components, in $ \{x|M|y\}$, resp. $ \{x|N|y\}$, containing $m$, resp. $n$.
\end{enumerate}
\end{Lemma}
\begin{proof} The first item follows by direct application of the homotopy lifting property  of the fibration, $P_L\colon L \to \fr(\W)$. The second follows from the fact that all fibres of a fibration over a path-connected space are homotopy equivalent.
\end{proof}
Recalling that we are assuming that $X,Y,L,M,N$ are HF, we have:
\begin{Lemma}\label{lem;2elsHF}  All spaces appearing in Definition \ref{space-2els} are HF spaces.\end{Lemma}
\begin{proof}
For the spatial 2-slice, of $\W$, at $(m,x,y,n)$ in \eqref{eq:s2slxymn},
this follows from the fact that both $L$ and $\fr(\W)$ are HF, (for the latter fact, see the discussion just after Definition \ref{def:HF-fibrant-window}), and Lemma \ref{2of3}, 
applied to the fibration $P_L\colon L \to \fr(\W)$. 
The same argument works for the general spatial 2-slices, in \eqref{eq:s2slgamma}.

 We have a fibration, $\langle P_{x,y},Q_{x,y}\rangle \colon   \left [ x \left |L \right |y\right ]  \to \{x |M|y\} \times \{x|N|y\},$
 in which the spaces, $\{x|M|y\}$ and  $\{x|N|y\}$, are both HF (see Lemma \ref{all-HF}), and all of whose fibres, i.e., all of the     $P_{L}^{-1}(m,\const_x,\const_y,n)$, with $m\in \{x|M|y\}$ and $n \in \{x|N|y\}$,
 are HF. We thus have that $ \left [ x \left |L \right |y\right ] $ is also HF, by the last item of Lemma \ref{2of3}.
\end{proof}

We define:
\begin{Definition}[Vertical span of slices]\label{def:vspan_el} Let  $\W\colon (p,M,p')\Longrightarrow (q,N,q')$ be a HF resolved 2-span, as in \eqref{diag:W}. Let $x \in X$ and $y \in Y$. The fibrant span, below, will be called the \emph{vertical span of slices, of $\W$, at $x$ and $y$},
  $$[x|\W|y]:=\left( \vcenter{\xymatrix@R=0pt@C=20pt
{ &[x|L| y]\ar[dl]_{P_{x,y}} \ar[dr]^{Q_{x,y}} \\
 \{x|M|y\} && \{x|N|y\}}}\right).$$
\end{Definition}
\noindent  By the discussion just given, $[x|\W|y]$ is a HF fibrant span.

\subsubsection{Homotopy invariance of spatial 2-slices}
We continue to fix a HF resolved 2-span, $\W$, as in \eqref{diag:W} and prove that, in several important cases, the spatial 2-slices of $\W$  are homotopy equivalent.

Firstly consider  holonomy maps for the fibrations, $\langle p,p' \rangle\colon M \to X \times Y$ and  $\langle q,q' \rangle\colon N \to X \times Y$, which will be denoted, respectively, by:
  $$\lambda^M \colon I \times (M \times_{X \times Y} (X \times Y)^I) \to M\textrm{ and } \lambda^N \colon I \times (N  \times_{X \times Y} (X \times Y)^I) \to N.$$

  Given paths, ${\gamma^X}\colon x\to  x'$ in $X$ and $\gamma^Y\colon y\to  y'$ in $Y$,
  these holonomy maps induce homotopy equivalences (where we are using the notation of Lemma \ref{main_techfib}),
   $$\Gamma_{{\langle \gamma^X,\gamma^Y\rangle}}^M\colon \{x |M |y\} \to \{x' |M |y'\}\quad \textrm{ and }\quad \Gamma_{{\langle \gamma^X,\gamma^Y\rangle}}^N\colon \{x |N |y\} \to \{x' |N| y'\}.$$
Here $\langle \gamma^X,\gamma^Y\rangle$ is the path, in $X \times Y$, such that, for $u\in I$, $u \mapsto \big(\gamma^X(u),\gamma^Y(u)\big)$.

\begin{Lemma}\label{lem:act2ls} Consider a HF resolved 2-span  as in \eqref{diag:W}. Let $x,x' \in X$ and $y,y'\in Y$.
Let $\gamma^X\colon I \to X$ be a path in $X$, from $x$ to $x'$, and  $\gamma^Y\colon I \to Y$ be one in $Y$, from $y$ to $y'$. Let $m \in \{x|M|y\}$, $n\in \{x | N | y\}$ and   $n' \in\{x' |N| y'\}$. We have:
 \begin{enumerate}[leftmargin=1cm]
 \item\label{lem:act2ls!} The two spaces, below, are homotopically equivalent,
 $$\left[\hskip-1pt\begin{array}{c|c|c}
&m&\\
x&L & y\\
&n&\end{array}\hskip-1pt\right] \quad \textrm{  and  } \quad  \left[\hskip-1pt\begin{array}{c|c|c}
&\Gamma_{{\langle \gamma^X,\gamma^Y\rangle}}^M(m)&\\
x'&L & y'\\
& \Gamma_{{\langle \gamma^X,\gamma^Y\rangle}}^N(n)&\end{array}\hskip-1pt\right].$$
 \item\label{lem:act2ls!!} The three spaces, below, are homotopically equivalent,
 $$\left[\hskip-1pt\begin{array}{c|c|c}
&m&\\
\gamma^X &L &\gamma^Y\\
&n'&\end{array}\hskip-1pt\right],\qquad
\left[\hskip-1pt\begin{array}{c|c|c}
&m&\\
x&L & y\\
&\Gamma_{ {\langle\overline{\gamma^X},\overline{\gamma^Y}\rangle}}^N(n')&\end{array}\hskip-1pt\right], \quad \textrm{  and   } \quad
\left[\hskip-1pt\begin{array}{c|c|c}
&\Gamma_{{\langle {\gamma^X},{\gamma^Y}\rangle}}^M(m)&\\
x'&L & y'\\
&n'&\end{array}\hskip-1pt\right].$$
 \end{enumerate}
(We recall that given a path $\gamma$, then $\overline{\gamma}$ denotes its reverse.)
\end{Lemma}
\begin{proof}
 The first point follows from the fact that the two spaces are fibres of the fibration, $P_L\colon L \to \fr(\W)$, over  points in the same path component of $\fr(\W).$ Indeed, the two spaces are
$$P_{L}^{-1}(m,\const_x,\const_y,n),\textrm{  and  }P_L^{-1}\big(  \Gamma_{{\langle \gamma^X,\gamma^Y\rangle}}^M(m) ,\const_{x'},\const_{y'}, \Gamma_{{\langle \gamma^X,\gamma^Y\rangle}}^N(n)\big),$$ respectively, for  points in $\fr(\W)$.
A path, in $\fr(\W)$,  connecting
$$(m,\const_x,\const_y,n) \qquad \textrm { and  } \qquad \big(  \Gamma_{{\langle \gamma^X,\gamma^Y\rangle}}^M(m) ,\const_{x'},\const_{y'}, \Gamma_{{\langle \gamma^X,\gamma^Y\rangle}}^N(n)\big),$$  is  given by
$$ 
t\mapsto \Big(\lambda^M\big(t,m,{\langle \gamma^X,\gamma^Y\rangle}\big),\const_{\gamma^X(t)},\const_{\gamma^Y(t)},  \lambda^N\big(t,n,
{\langle \gamma^X,\gamma^Y\rangle}\big)\Big). $$

The second point follows again by identifying the spaces with fibres of the fibration $P_L\colon L \to \fr(\W)$, and examining different points of the same path component. For instance, a path, in $\fr(\W)$, connecting the points,
$$(m,\gamma^X,\gamma^Y,n')\qquad  \textrm{  and
} \qquad \big(m,\const_x,\const_y,\Gamma_{
{\langle \overline{\gamma^X},\overline{\gamma^Y}\rangle}}^N(n')\big),$$
is
$$
s \mapsto  \Big (m,\gamma^X_s,
\gamma^Y_s,
\lambda^N\big(s,n',{\langle\overline{\gamma^X},\overline{\gamma^Y}\rangle }\big )\Big).$$
Here, given a path $\gamma\colon I \to B$, where $B$ is a space, and $s \in [0,1]$, we have written  $\gamma_s\colon I \to B$ for the path $t\mapsto \gamma\big((1-s)t\big )$.
\end{proof}

\subsubsection{The natural transformation of profunctors associated to a HF resolved fibrant 2-span}\label{sec_resolvedspans_to_nat_tranfs}
As usual, we fix $\kappa$,  a subfield of the complex field $\C$.
Let  
$$\W=(l_X,P,L,Q,r_Y)\colon \big ( (p,M,p')\colon X \to Y \big) \Longrightarrow \big ( (q,N,q')\colon X \to Y\big)$$ be a HF fibrant resolved 2-span, as in \eqref{diag:W}, connecting the HF fibrant spans, $(p,M,p')$ and $(q,N,q')$. 
We assume given, or chosen, subsets $\x_X {\subseteq}  X$ and $\y_Y {\subseteq}  Y$.

We have $\Vect$-profunctors, as  defined in \S\ref{sec:profunctor_fib_span}, in particular  in Notation \ref{not_barH},
\begin{align*}\Hpb^M_{(\x_X,\y_Y)}\colon \pi_1(X,\x_X)^\op \times \pi_1(Y,\y_Y) \to \Vect, \\ \intertext{and}
 \Hpb^N_{(\x_X,\y_Y)}\colon \pi_1(X,\x_X)^\op \times \pi_1(Y,\y_Y) \to \Vect. 
\end{align*}
These are the restrictions, to $\pi_1(X,\x_X)^\op\times \pi_1(Y,\y_Y)$, of the profunctors,
\begin{align*}\Hp^M\colon \pi_1(X,X)^\op \times \pi_1(Y,Y) \to \Vect, \\ \intertext{and}
 \Hp^N\colon \pi_1(X,X)^\op \times \pi_1(Y,Y) \to \Vect,
\end{align*}
in Definition \ref{def:fibrantSpanstoprof}.

 We want to  define a natural transformation of profunctors, denoted $$\tHpb{\W}{\x_X}{\y_Y}\colon \Hpb^M_{(\x_X,\y_Y)} \To \Hpb^N_{(\x_X,\y_Y)}.$$
This natural transformation will, itself,  be the restriction of a natural transformation,
which we will define first, and
  denote $$\tHp{\W} \colon \Hp^M\  \To \Hp^N.$$
 
 Of course, given $x\in X$ and $y\in Y$, $\Hp^M(x,y)$ and $\Hp^N(x,y)$ are both $\kappa$-vector spaces, so to specify a natural transformation, as required, we have to specify a linear transformation, $\tHpc{\W}{x}{y}$, from  $\Hp^M(x,y)$ to  $\Hp^N(x,y)$, depending on $x$ and $y$ in a `natural' way. Given that, by Lemma \ref{lem:finiteness}, both vector spaces are finite dimensional, we may specify this linear map by giving its matrix elements with respect to the evident bases, consisting of path-components, $\PC_m(\{x|M|y\})$ and  $\PC_n(\{x|N|y\})$, of $m \in  \{x|M|y\}$ and of $n \in  \{x|N|y\}$.

\begin{Definition} Consider a HF resolved 2-span, $\W\colon (p,M,p')\Longrightarrow (q,N,q'),$ as in \eqref{diag:W}.
Given $x\in X$ and $y\in Y$, we define the linear map,
$$\tHpc{\W}{x}{y}\colon \Hp^M(x,y)\to  \Hp^N(x,y),$$
where
$$\Hp^M(x,y)=\Lin \big(\hpi_0(\{x|M|y\})\big) \qquad \textrm{ and } \qquad \Hp^N(x,y)=\Lin\big(\hpi_0(\{x|N|y\})\big).$$
This is to  have the following matrix elements, with respect to the usual bases.

Given $m\in \{x|M|y\}$ and $n\in \{x|N|y\}$,
\begin{multline}\label{eq:2matrixels}
\ \big \langle \PC_{m}(\{x|M|y\}) \mid \tHpc{\W}{x}{y} \mid  \PC_{n}(\{x|N|y\}) \big\rangle \\:= \chi^\pi \left(\left[\hskip-1pt\begin{array}{c|c|c}
&m&\\
x&L & y\\
&n&\end{array}\hskip-1pt\right]\right) \chi^\pi\big (\PC_n(\{x|N|y\} ) \big).
\end{multline}\end{Definition}
\noindent We refer to Definition \ref{space-2els} for notation.
Note that, by Lemma \ref{lem;2elsHF}, all the spatial 2-slices met here are HF spaces, so we can consider their homotopy content.

\begin{Lemma}\label{def_2HW} Consider a HF resolved 2-span, as in \eqref{diag:W},
$$\W\colon \big ( (p,M,p')\colon X \to Y \big) \Longrightarrow \big ( (q,N,q')\colon X \to Y\big).$$
\begin{itemize}[leftmargin=1cm]
 \item 
Let $x \in X$ and $y \in Y$. Given $m\in \{x|M|y\}$ and $n\in \{x|N|y\}$, the value of the right-hand-side of
Equation \eqref{eq:2matrixels} depends only on the path-components, in  $\{x|M|y\}$ and in $\{x|N|y\}$, respectively, to which $m$ and $n$ belong.
  \item The family of  linear maps, $\tHpc{\W}{x}{y}\colon \Hp^M(x,y)\to  \Hp^N(x,y)$, for all $x \in X$ and $y \in Y$, together defines a natural transformation of $\Vect$-profunctors, $$\tHp{\W}\colon \big (\Hp^M\colon  \pi_1(X,X) \bto \pi_1(Y,Y)\big ) \Longrightarrow \big( \Hp^N\colon  \pi_1(X,X) \bto \pi_1(Y,Y)\big ) ,$$
  and, therefore (by Lemma \ref{lem:finiteness}) a 2-morphism in the bicategory $\vProfGrphf$.
\end{itemize}
  
\end{Lemma}
\begin{proof}
 The first statement follows directly from Item \eqref{item:independence} of Lemma \ref{lem:matrix2els}, since  the homotopy content of a homotopy finite space is a homotopy invariant.

 The second statement follows   from point \eqref{lem:act2ls!} of Lemma \ref{lem:act2ls}, given the explicit forms of  $\Hp^M$ and  $\Hp^N$  in Definition \ref{def:fibrantSpanstoprof}. We also use the fact that, given
   paths  ${\gamma^X}\colon x \to x'$ in $ X$, and  ${\gamma^Y}\colon y \to y'$ in $Y$, the holonomy map, $\lambda^N$, for $\langle q,q' \rangle\colon N \to X\times Y$, gives rise to a homotopy equivalence, between fibres,
   $$\Gamma_{{\langle \gamma^X,\gamma^Y\rangle}}^N\colon \{x |N |y\} \to \{x' |N| y'\},$$
   and, in particular, induces a bijection between the sets of path components. Moreover, given $n \in  \{x |N |y\}$, $\Gamma_{{\langle \gamma^X,\gamma^Y\rangle}}^N$  therefore restricts to a homotopy equivalence,
   $$\PC_n\big ( \{x|N|y\}\big)\cong  \PC_{\Gamma_{{\langle \gamma^X,\gamma^Y\rangle}}^N(n)}\big ( \{x'|N|y'\}\big).$$
The remaining details are left to the reader.
\end{proof}

\begin{Definition}\label{def:2matrixels} Choose $\x_X {\subseteq}  X$ and $\y_Y {\subseteq}  Y$. The natural transformation,
 $$\tHpb{\W}{\x_X}{\y_Y}\colon \Hpb^M_{(\x_X,\y_Y)} \To \Hpb^N_{(\x_X,\y_Y)},$$  is defined by restricting $\tHp{\W}\colon\Hp^M\Rightarrow \Hp^N$ to $\pi_1(X,\x_X)^{\op} \times \pi_1(Y,\y_Y)$.
\end{Definition}
Explicitly, given $x\in \x_X$ and $y\in \y_Y$, the linear map,
$$(\tHpb{\W}{\x_X}{\y_Y})_{(x,y)}\colon \Hpb^M_{(\x_X,\y_Y)}(x,y) \to  \Hpb^N_{(\x_X,\y_Y)}(x,y), $$
is, therefore,  $\tHpc{\W}{x}{y}\colon \Hp^M(x,y)\to  \Hp^N(x,y)$.

\begin{Remark}When $\kappa=\C$, there is a 1-parameter version of $\tHp{\W}\colon\Hp^M \Rightarrow \Hp^N,$  denoted  $\tHp{(\W,s)}\colon\Hp^M\Rightarrow \Hp^N,$ where $s \in \C$. This has as matrix elements
 \begin{multline*}
 \left \langle \PC_{m}(\{x|M|y\}) \mid \tHpc{(\W,s)}{x}{y} \mid  \PC_{n}(\{x|N|y\}) \right\rangle \\ := \chi^\pi \left(\left[\hskip-1pt\begin{array}{c|c|c}
&m&\\
x&L & y\\
&n&\end{array}\hskip-1pt\right]\right) \chi^\pi\big (\PC_n(\{x|N|y\} ) \big)   \chi^\pi(\PC_x(X))^{1-s}  \chi^\pi(\PC_y(Y))^{s}.
\end{multline*}
All results go through with this extra generality. This is a special case of a more general 2-parameter version, also involving the vertical direction of the spatial 2-slices. This is left to the reader to explore.
\end{Remark}

\begin{Remark}\label{rem:simpler-matrix-2-els}Going back to Definition \ref{space-2els},  let $x \in X$, $y \in Y$,
 $m \in  \{x| M| y\}$ and $n\in  \{x| N| y\}$. Define the following space,
$$
 \left[\hskip-1pt\hskip-1pt\begin{array}{c|c|c}
&m&\\
x&L & y\\
&[n]&\end{array}\hskip-1pt\right]:=P_{L}^{-1}\big (\{m\}\times \{\const_x\} \times \{\const_y\}\times \PC_n(\{x|N|y\})\big).$$
We then have, by using the same argument as in Lemma \ref{rels_matrix_els}, the following,
\begin{equation}\label{eq:simpler-matrix-2-els}
\big \langle \PC_{m}(\{x|M|y\}) \mid \tHpc{\W}{x}{y} \mid  \PC_{n}(\{x|N|y\}) \big\rangle = \chi^\pi \left(\left[\hskip-1pt\begin{array}{c|c|c}
&m&\\
x&L & y\\
&[n]&\end{array}\hskip-1pt\right]\right)
.
\end{equation}
\end{Remark}

\subsubsection{The symmetric monoidal like structure of $\tHpb{\W}{-}{-}$}\label{sec:H-is-sym-monoidal}
The following string of results will  be used later on (Subsection \ref{sec:sym_mon})  to prove that the constructions of the once-extended Quinn TQFTs, defined in this paper, give bifunctors     which, furthermore, can be given symmetric monoidal structures. Our starting point is Lemmas \ref{Sym_mon:1} and \ref{def_chi} of \S \ref{sec:sim_mon_0}.

In order to simplify the notation we will temporarily denote $\PC_x(B)=[x]_B$.
\begin{Lemma}\label{Sym_mon:2}
 Consider HF fibrant resolved 2-spans, denoted
 \begin{align*}
 \W\colon \big ( (p,M,q)\colon X \to Y \big) &\To  \big (( f,N,g)\colon X \to Y \big),\\
\W'\colon \big ( (p',M',q')\colon X' \to Y' \big) & \To  \big ( (f',N',g')\colon X' \to Y' \big).
 \end{align*}
The diagrams, for $\W$ and $\W'$, will be \begin{equation}
 \vcenter{\xymatrix@R=17pt{
  &X  &  M\ar[l]_{p}\ar[r]^{q}  & Y\\
  &X^I\ar[u]^{s_X} \ar[d]_{t_X} & L \ar[l]_{l}\ar[r]^{r}\ar[u]^{P}\ar[d]_{Q}& Y^I \ar[u]_{s_Y} \ar[d]^{t_Y}\\
  &X  &  N\ar[l]^{f}\ar[r]_{g}  & Y}} \textrm{ \quad and }  \vcenter{
  \xymatrix@R=17pt{
  &X'  &  M'\ar[l]_{p'}\ar[r]^{q'}  & Y'\\
  &{X'}^I\ar[u]^{s_{X'}} \ar[d]_{t_{X'}} & L' \ar[l]_{l'}\ar[r]^{r'}\ar[u]_{P'}\ar[d]_{Q'}& {Y'}^I \ar[u]_{s_{Y'}} \ar[d]^{t_{Y'}}\\
  &X'  &  N'\ar[l]^{f'}\ar[r]_{g'}  & Y'.}}
  \end{equation}
The following hold.
\begin{enumerate}[leftmargin=1cm]
 \item The window, $\W\times\W'$, below, is a HF fibrant resolved 2-span,
 \begin{multline*}
 \big((p\times p',M\times M',q\times q')\colon X\times X' \to Y\times Y' \big) \\ \To  \big ((f\times f',N\times N',g\times g')\colon X\times X' \to Y\times Y' \big),
 \end{multline*}
 $$\W\times\W':=\vcenter{\xymatrix@C=50pt{
  X \times X' &  M\times M'\ar[l]_{p\times p'}\ar[r]^{q\times q'}  & Y\times Y'\\
  (X\times X')^I \ar[u]^{s_X\times s_{X'}} \ar[d]_{t_X\times t_{X'}} & L\times L' \ar[l]_<<<<<<<<<{l\times l'}\ar[r]^>>>>>>>>>{r\times r'}\ar[u]^{P\times P'}\ar[d]_{Q\times Q'}&  (Y\times Y')^I \ar[u]_{s_Y\times s_{Y'}} \ar[d]^{t_Y\times t_{Y'}}\\
  X\times X'  &  N\times N'\ar[l]^{f\times f'}\ar[r]_{g\times g'}  & Y\times Y'.}}$$
We used that
  $ (X\times X')^I \cong X^I\times {X'}^I$ and $(Y\times Y')^I \cong Y^I\times {Y'}^I$.
\item Let $x \in X$, $x'\in X'$, $y \in Y$ and $y'\in Y'$. Given $m\in \{x|M|y\}, n\in \{x|N|y\}$, $m'\in \{x'|M'|y'\}$ and $n'\in \{x'|N'|y'\}$, we have that,
\begin{multline*}
\left \langle [{(m,m')}]_{\{(x,x')|M\times M'|(y,y')\}} \mid \tHpc{\W\times \W'}{(x,x')}{(y,y')} \mid  [{(n,n')}]_{\{(x,x')|N\times N'|(y,y')\}} \right\rangle
\\=\left \langle [{m}]_{\{x|M|y\}} \mid \tHpc{\W}{x}{y} \mid  [{n}]_{\{x|N|y\}} \right\rangle \, \left \langle [{m'}]_{\{x'|M'|y'\} } \mid \tHpc{\W'}{x'}{y'} \mid  [{n'}]_{\{x'|N'|y'\}} \right\rangle.
\end{multline*}
\end{enumerate}
\end{Lemma}

\begin{proof}
The first point follows from the fact that the product of fibration is a fibration. The second follows from the fact that, clearly,
$$\left[\hskip-1pt\begin{array}{c|c|c}
&(m, m') &\\
(x, x') \, &L \times L' & \, (y, y')\\
&(n, n') &\end{array}\hskip-1pt\right]\cong \left[\hskip-1pt\begin{array}{c|c|c}
&m&\\
x\, &L  & \, y\\
&n&\end{array}\hskip-1pt\right]\times \left[\hskip-1pt\begin{array}{c|c|c}
&m'&\\
x' \, & L' & \, y'\\
&n'&\end{array}\hskip-1pt\right],$$
and
$$\PC_{(n,n')}\big(\{(x,x')|N\times N'|(y,y')\}  \big)\cong \PC_{n}\big ( \{x|N|y\}\big)\times  \PC_{n'}\big ( \{x'|N'|y'\}\big) .$$
The formula therefore follows from the fact that the homotopy content of  HF spaces is multiplicative with respect to their product.
\end{proof}

In order to prove that the once-extended Quinn TQFT is a symmetric monoidal bifunctor, it is  convenient to change the language of the previous result, to a language closer to that of  Definition \ref{Sym. Mon. Bifunctor}. In particular, combining Lemma \ref{Sym_mon:2} with Lemma \ref{def_chi}, whose notation we follow, gives the following.
\begin{Lemma}\label{lem:chi-is-natural} Let $\x_X {\subseteq}  X$, $\overline{x}'_{X'} {\subseteq}  X'$,   $\y_Y {\subseteq}  Y$ and  $\overline{y}'_{Y'} {\subseteq}  Y'.$
The two natural transformations of $\Vect$-profunctors, i.e. the two 2-morphisms in $\vProfGrphf$, obtained by pasting the diagrams in $\vProfGrphf$, below, coincide,
$$\xymatrix@C=42.5pt@R=30pt{\pi_1(X,\x_X)\times \pi_1(X',\overline{x}'_{X'})\ar[rrr]_<<<<<<<<<<<<<<<<<<<<<<<<<{\Hpb_{(\x_X,\y_Y)}(f,N,g)
\otimes  \Hpb_{(\overline{x}'_{X'},\overline{y}'_{Y'})}(f',N',g')  }
\ar@/^2pc/[rrr]^{\Hpb_{(\x_X,\y_Y)}(p,M,q)
\otimes  \Hpb_{(\overline{x}'_{X'},\overline{y}'_{Y'})}(p',M',q')  }
\xtwocell[rr]{}<>{<-2>\qquad \qquad \qquad \qquad\quad\big(\tHpb{\W}{\x_X}{\y_Y} \otimes \tHpb{\W'}{\overline{x}'_{X'}}{\overline{y}'_{Y'}}\big)}
\ar[d]_{ \varphi^{m_{(X,X')}}} \xtwocell[drrr]{}<>{<1>\qquad\,\,\chi_{(N,N')}}&&& \pi_1(Y,\y_Y)\times \pi_1(Y',\overline{y}'_{Y'})\ar[d]^{\varphi^{m_{(Y,Y')}}}\\
 \pi_1(X \times X',\x_X\times \overline{x}'_{X'})\ar[rrr]_{\Hpb_{(\x_X\times \overline{x}'_{X'} ,\y_Y\times \overline{y}'_{Y'})}\left  (f\times f',N\times N',g \times g'\right)} &&& \pi_1(Y \times Y',\y_Y\times \overline{y}'_{Y'}),}$$
and
$$\xymatrix@C=42.5pt@R=30pt{\pi_1(X,\x_X)\times \pi_1(X',\overline{x}'_{X'})\ar[rrr]^{\Hpb_{(\x_X,\y_Y)}(p,M,q)
\otimes  \Hpb_{(\overline{x}'_{X'},\overline{y}'_{Y'})}(p',M',q')  }
\ar[d]_{ \varphi^{m_{(X,X')}}} \xtwocell[rrrd]{}<>{<-0.8>\qquad\quad\chi_{(M,M')}}&&& \pi_1(Y,\y_Y)\times \pi_1(Y',\overline{y}'_{Y'})\ar[d]^{\varphi^{m_{(Y,Y')}}}\\
 \pi_1(X \times X',\x_X\times \overline{x}'_{X'})\ar[rrr]^{\Hpb_{(\x_X\times \overline{x}'_{X'} ,\y_Y\times \overline{y}'_{Y'})}\left  (p\times q,M\times M',q \times q'\right)}
 \ar@/_2pc/[rrr]_{\Hpb_{(\x_X\times \overline{x}'_{X'} ,\y_Y\times \overline{y}'_{Y'})}\left  (f\times f',N\times N',g \times g'\right)}
\xtwocell[rr]{}<>{<2>\qquad \qquad \qquad \qquad \big(\tHpb{\W\times \W'}{\x_X \times \overline{x}'_{X'}}{\y_Y \times\overline{y}'_{Y'}} \big)}
&&&  \pi_1(Y \times Y',\y_Y\times \overline{y}'_{Y'}).}$$
\end{Lemma}

We continue to follow the notation in Definition \ref{Sym. Mon. Bifunctor}.
\begin{Notation}\label{def_omega}
Given pairs of spaces, $(X,\x_X)$, $(X',\overline{x}'_{X'})$ and  $(X'',\overline{x}''_{X''})$, with $X,X'$ and $X''$ homotopy finite, we have a canonical invertible 2-morphism in $\vProfGrphf$, as shown in the diagram below. (We have condensed the notation, so $\pi(X)$ means $\pi_1(X,\x_X)$, $X\times X'$ means $(X\times X',\x_X\times \overline{x}'_{X'})$, and so on.)
$$
\hskip-1cm
\xymatrix@R=15pt@C=17pt{&&\big(\pi(X)\times \pi(X')\big)\times \pi(X'')
 \ar[dl]_>>>>>>>>{\varphi^{\alpha^{\Grp}_{(\pi(X),\pi(X'),\pi(X''))}\qquad}}\ar[dr]^>>>>>>>>{\qquad \varphi^{m_{(X,X')}}\times \id_{\pi(X'')}}\\
&\pi(X) \times \big (\pi(X')\times\pi( X'') \big )\ar[d]^{\id_{\pi(X)} \times \varphi^{m_{(X',X'')}}}
& \ar@{}[dd]\dtwocell<\omit>{\omega_{(X,X',X'')}} &\pi(X \times X' )\times \pi(X'')\ar[d]_{\varphi^{m_{(X\times X',X'')}}} \\&\pi(X) \times \big (\pi(X' \times  X'') \big)\ar[dr]_{\varphi^{m_{(X,X'\times X'')}}\quad}
&&\pi\big ((X \times X') \times X''\big) \ar[dl]^<<<<<<{\varphi^{\pi\left (\alpha^{\CGWH}_{(X,X',X'')}\right )}}\\
&&\pi\big (X \times (X' \times X'') \big   ).
}
$$
\end{Notation}
This 2-morphism of profunctors arises from the fact that, if we switch back the profunctors in the arrows of the diagram above, to the functors that gave rise to them, then, applying Example \ref{Ex:maps_to_profunctors}, gives rise to a commutative diagram of functors. Indeed, note that, in general, if $F\colon C \to C'$ and $G\colon C'\to C''$ are functors, then we have a canonical natural isomorphism from the profunctor, $\varphi^{G\circ F}\colon C \bto C''$, to the composition of the profunctors,
$C \stackrel{\varphi^F}{\bto} C' \stackrel{\varphi^{G}}{\bto} C'';$ for more details see Example \ref{varphi is pseudo}.

The 2-morphisms, $\omega_{(X,X',X'')}$ in $\vProfGrphf$ satisfy an obvious cocycle condition, given pairs of spaces,  $(X,\x_X)$, $(X',\overline{x}'_{X'})$, $(X'',\overline{x}''_{X''})$, and $(X''',\overline{x}'''_{X'''})$. The equations satisfied are those in \cite[\S 4.3]{Gurski_book} / \cite[page 17]{GPS}. This follows from an explicit calculation.

These cocycles, $\omega_{(X,X',X'')}$, are furthermore compatible with the natural transformations, $\tHpb{\W}{-}{-}$ and the $\chi_{(-,-)}$, as we now explain.

Suppose that we are given three HF spans,
$$
(p,M,q)\colon X \to Y, \qquad (p',M',q')\colon X' \to Y', \quad \text{ and \, } (p'',M'',q'')\colon X'' \to Y''.$$
The 2-morphisms, $\chi_{(-,-)},$ in  $\vProfGrphf$, arising from Lemma \ref{def_chi}, can be pasted in two different ways, as written below (where we have again condensed the notation, in the obvious way),
$$\xymatrixcolsep{2.7cm}
\hskip-2.5cm
\xymatrix{&\big(\pi(X)\times \pi(X')\big)\times \pi(X'')\ar[d]_{\varphi^{m_{(X,X')}}\times \id_{\pi(X'')}} \drtwocell<\omit>{\qquad\qquad \quad\chi_{(M,M')} \otimes  \Hpb^{M''}}
\ar[r]^{(\Hpb^{M} \otimes \Hpb^{M'}) \otimes \Hpb^{M''}} & \big (\pi(Y)\times \pi(Y')\big)\times \pi(Y'')\ar[d]^{\varphi^{m_{(Y,Y')}}\times \id_{\pi(X'')}}  \\
&\pi(X \times X' )\times \pi(X'') \ar[r]|{\Hpb^{M\times M'} \otimes \Hpb^{M''}}
\ar[d]_{\varphi^{m_{(X\times X',X'')}}} \drtwocell<\omit>{\qquad\qquad \,\chi_{(M\times M',M'')}}
& \pi(Y \times Y')\times \pi(Y'')\ar[d]^{\varphi^{m_{(Y\times Y',Y'')}}}
\\
&\pi\big ((X \times X') \times X''\big)
\ar[d]_{\varphi^{\pi\left(\alpha^{\CGWH}_{(X,X',X'')}\right)}} \drtwocell<\omit>{\cong}
\ar[r]|{\Hpb^{(M\times M')\times M''}} & \pi\big((Y \times Y')\times Y''\big) \ar[d]^{\varphi^{\pi\left(\alpha^{\CGWH}_{(Y,Y',Y'')}\right)}}\\
&\pi\big(X \times ( X' \times X'')\big) \ar[r]_{\Hpb^{M\times (M' \times M'')}} & \pi\big (Y \times ( Y' \times Y'')\big),
}
$$
and
$$\xymatrixcolsep{2.7cm}
\hskip-2cm
\xymatrix{&\big (\pi(X)\times \pi(X')\big)\times \pi(X'')\ar[d]_{\varphi^{\alpha^{\Grp}_{ (\pi(X),\pi(X'),\pi(X''))}}} \drtwocell<\omit>{\cong}
\ar[r]^{(\Hpb^{M} \otimes \Hpb^{M'}) \otimes \Hpb^{M''}} & \big (\pi(Y)\times \pi(Y')\big)\times \pi(Y'')\ar[d]^{\varphi^{\alpha^{\Grp}_{(\pi(Y),\pi(Y'),\pi(Y''))}}} \\
&\pi(X) \times \big ( \pi(X' )\times \pi(X'') \big) \ar[r]^{\Hpb^{M} \otimes (\Hpb^{M'} \otimes \Hpb^{M''})}
\ar[d]_{\pi(X)\times\varphi^{m_{(X',X'')}}} \drtwocell<\omit>{\qquad\qquad \quad\Hpb^{M}\otimes\chi_{(M',M'')}}
& \pi(Y) \times \big ( \pi(Y' )\times \pi(Y'')\big) \ar[d]^{\pi(Y)\times\varphi^{m_{(Y',Y'')}}}
\\
&\pi (X) \times\pi( X' \times X'')
\ar[d]_{  \varphi^{m_{(X,X'\times X'')}}} \drtwocell<\omit>{\qquad \qquad \,\chi_{(M,M'\times M'')}}
\ar[r]|{\Hpb^{M} \otimes \Hpb^{ M' \times M''}} & \pi(Y) \times \pi( Y' \times Y'') \ar[d]^{{  \varphi^{m_{(Y,Y'\times Y'')}}}}\\
&\pi\big (X \times ( X' \times X'')\big) \ar[r]_{\Hpb^{M\times (M' \times M'')}} & \pi\big (Y \times ( Y' \times Y'')\big).
}
$$
We note that the first of these diagrams fits together with the diagram for $\omega_{(X,X',X'')}$ on the left, and the various  2-morphisms compose well. The second lower diagram likewise composes, this time with $\omega_{(Y,Y',Y'')}$ and on the right, again giving a second composite 2-morphism, which has the same source and target composite 1-morphisms as the first.
\begin{Lemma}\label{Omega_is_natural}
The two composite 2-morphisms obtained as above by pasting,  respectively, $\omega_{(X,X',X'')}$  or   $\omega_{(Y,Y',Y'')}$
 to the two diagrams, are equal.
\end{Lemma}
\begin{proof}This follows from an explicit calculation.
\end{proof}

Given HF spans $(p,M,q)\colon X \to Y$ and  $(p',M',q')\colon X' \to Y',$
the profunctor,$$
 \Hpb^{M\times M'}_{(\x_X\times \overline{x}'_{X'} ,\y_Y\times \overline{y}'_{Y'})}\colon  \pi_1(X\times X',\x_X\times \overline{x}'_{X'})\bto \pi_1(Y \times Y',\y_{Y } \times \overline{y}'_{Y' }),$$ is similarly well behaved with respect to swapping the order of coordinates, and with products with trivial spaces, if $X=\{*\}$ or $Y=\{*\}$. We leave it to the reader to unpack what this means in terms of diagrams similar to those  just presented.

\subsection{The horizontal composition of HF resolved 2-spans}\label{sec:hor_comp_2spans}
HF resolved 2-spans can be composed  both horizontally and vertically. Here we first look at the horizontal composition before considering the effect of that composition on the corresponding natural transformations.

\subsubsection{The horizontal composition of HF resolved 2-spans in detail}\label{sec:horcomp}
Consider HF fibrant spans, $$(p_1,M_1,p_1')\colon X \to Y \textrm{ and } (p_2,M_2,p_2')\colon Y \to Z.$$  Their composition from  Lemma \ref{main3} and Definition \ref{main3def},  written as  $$(p_1,M_1,p_1')\bullet (p_2,M_2,p_2')=(\overline{p_1},M_1\times_Y M_2,\overline{p_2'})\colon X\to Z$$ is itself a HF fibrant span. It is defined by the pullback diagram, appearing as the diamond in the commutative diagram below,
$$\xymatrix@R=0pt{ && 
M_1{\times}_Y M_2  \ar[dl]_{}\ar[dr]^{} \ar@/_1.2pc/[ddll]_>>>>>>>{\overline{p_1}} \ar@/^1.2pc/[ddrr]^>>>>>>>{\overline{p_2'}}\\ &M_1\ar[dl]^{p_1}\ar[dr]_{p_1'}
&&M_2\ar[dl]^{p_2}\ar[dr]_{p_2'}\\   X & &Y  & &Z. }
$$

Consider a diagram of spaces, HF spans, and HF resolved 2-spans, as below,
$$\xymatrix{&X \ar@/^1pc/[rr]^{(p_1,M_1,p_1')}\ar@/_1pc/[rr]_{ {(q_1,N_1,q_1')}} & \Downarrow \W_1 & Y \ar@/^1pc/[rr]^{(p_2,M_2,p_2')}\ar@/_1pc/[rr]_{ {(q_2,N_2,q_2')}} & \Downarrow \W_2 & Z
}
 ,$$
 where the diagrams for $\W_1$ and $\W_2$ are the windows shown below,
 \begin{equation}\label{ash0}
 \vcenter{\xymatrix@R=15pt{
  &X  &  M_1\ar[l]_{p_1}\ar[r]^{p'_1}  & Y\\
  &X^I\ar[u]^{s_X} \ar[d]_{t_X} & L_1 \ar[l]_{l_1}\ar[r]^{r_1}\ar[u]^{P_1}\ar[d]_{Q_1}& Y^I \ar[u]_{s_Y} \ar[d]^{t_Y}\\
  &X  &  N_1\ar[l]^{q_1}\ar[r]_{q'_1}  & Y}} \textrm{ \quad and }  \vcenter{\xymatrix@R=15pt{
  &Y  &  M_2\ar[l]_{p_2}\ar[r]^{p'_2}  & Z\\
  &Y^I\ar[u]^{s_Y} \ar[d]_{t_Y} & L_2 \ar[l]_{l_2}\ar[r]^{r_2}\ar[u]^{P_2}\ar[d]_{Q_2}& Z^I \ar[u]_{s_Z} \ar[d]^{t_Z}\\
  &Y  &  N_2\ar[l]^{q_2}\ar[r]_{q'_2}  & Z.}}
  \end{equation}
  We will also need to consider  the fillers, Definition \ref{def:window}, of $\W_1$ and $\W_2$, as usual denoted by
  $$P_{L_1}\colon L_1 \to \fr(\W_1) \textrm{ \quad   and   \quad }  P_{L_2}\colon L_2 \to \fr(\W_2). $$

   We will define the horizontal composite, $\W_1\#_0 \W_2$, of $\W_1$ and $\W_2$, in such a way  that $\W_1\#_0 \W_2 $ is a HF resolved 2-span, which fits inside the diagram below,
   \begin{equation}\label{diag:W1W2}{
  \vcenter{ \xymatrix{\mort{X}{Z}{(\overline{p_1}, M_1\times_Y M_2,\overline{p_2'})}{(\overline{q_1}, N_1\times_Y N_2,\overline{q_2'})}{\W_1\#_0\W_2}} }  =
  \vcenter{ \xymatrix{\mort{X}{Z}{(p_1,M_1,p_1')\bullet (p_2,M_2,p_2')}{(q_1,N_1,q_1')\bullet (q_2,N_2,q_2')}
   {\W_1\#_0\W_2}.}}}
   \end{equation}
This is  done by considering the  obvious pullback along the common vertical HF fibrant span in \eqref{ash0}. Explicitly the horizontal composition  of $\W_1$ and $\W_2$ will be given by the window
 \begin{equation}\label{def:W1W2}
  \W_1 \#_0 \W_2:= \vcenter{\xymatrix@C=40pt@R=18pt{
  X  &  M_1\times_Y M_2\ar[l]_{\overline{p_1}}\ar[r]^{\overline{p_2'}}  & Z\\
  X^I\ar[u]^{s_X} \ar[d]_{t_X} & L_1\times_{Y^I} L_2 \ar[l]_<<<<<<<{\overline{l_X}}\ar[r]^<<<<<<<{\overline{r_Z}}\ar[u]^{\mathcal{P}}\ar[d]_{\mathcal{Q}}& Z^I \ar[u]_{s_Z} \ar[d]^{t_Z}\\
X  &  N_1\times_Y N_2\ar[l]^{\overline{q_1}}\ar[r]_{\overline{q_2'}}  & Z.
  }}
\end{equation}
Here we will need to consider the pullback diagram included (as the middle diamond) in the commutative diagram below,  in which, \\(i)  $\Psi_{L_1,L_2}\colon L_1\times_{Y^I} L_2\to Y^I$ is ${\Psi}_{L_1,L_2}=l_2\circ \proj_2=r_1\circ\proj_1$, \\(ii)  $\overline{l_X}=l_1\circ \proj_1$,\\(iii) $\overline{r_Z}=r_2\circ \proj_2$,
\begin{equation}\label{Equation:psiL1L2}\vcenter{
\xymatrix@R=8pt{&&& L_1 \times_{Y^I} L_2\ar@/_2pc/[ddll]_>>>>>>>>>{\overline{l_X}} \ar@/^2pc/[ddrr]^>>>>>>>>>{\overline{r_Z}} \ar[dd]|{\Psi_{L_1,L_2}}\ar[dl]_>>>>{\proj_1} \ar[dr]^>>>>{\proj_2}  \\& & L_1\ar[dr]_{r_1}\ar[dl]^{l_1} && L_2\ar[dl]^{l_2} \ar[dr]_{r_2}\\& X^I && Y^I &&Z^I\,.}}
\end{equation}

We have obvious maps, $\P\colon L_1\times_{Y^I} L_2 \to  M_1 \times_Y M_2$, induced by $P_1\colon L_1 \to M_1$ and $P_2\colon L_2 \to M_2$, and  $\mathcal{Q}\colon L_1\times_{Y^I} L_2 \to  N_1 \times_Y N_2$, induced by $Q_1\colon L_1 \to M_1$ and $Q_2\colon L_2 \to N_2$.

Let us explain why  $\W_1\#_0 \W_2$, in Equation \eqref{def:W1W2}, is a HF fibrant resolved 2-span, fitting inside diagram \eqref{diag:W1W2}. This follows by a sequence of observations.
\begin{enumerate}[leftmargin=0.6cm]
 \item Lemma \ref{main3}, together with the fact that the top, bottom and middle horizontal (Lemma \ref{lem:midhor_is_fib}) spans of the diagrams in \eqref{ash0} are HF fibrant, implies that all spaces appearing in diagram \eqref{def:W1W2}   are HF.
 \item Secondly, we  note that the naturally defined map, $\mathcal{P}_{L_1,L_2}$ below, which, by definition, is the filler of $\W_1\#_0\W_2$, is a fibration,
$$
\mathcal{P}_{L_1,L_2} \colon L_1\times_{Y^I} L_2 \to \lim \left( \vcenter{ \xymatrix@C=30pt@R=13pt{
  X  &  M_1\times_Y M_2\ar[l]^>>>>>>>>{\overline{p_1}}\ar[r]_<<<<<<<{\overline{p_2'}}  & Z\\
  X^I\ar[u]^{s_X} \ar[d]_{t_X} &  & Z^I \ar[u]_{s_Z} \ar[d]^{t_Z}\\
X  &  N_1\times_Y N_2\ar[l]_>>>>>>>{\overline{q_1}}\ar[r]^<<<<<<{\overline{q_2'}}  & Z
  }}\right).
$$
The argument to prove this  is very similar to that  in the proof of Lemma \ref{main3}:
\begin{enumerate}[label=(\alph*),leftmargin=1cm]
 \item First note that we can put the fillers, $P_{L_1}\colon L_1 \to \fr(\W_1)$ and $P_{L_2}\colon L_2\to \fr(\W_2)$, together, to obtain a map,
 $P_{L_1,L_2}$, as below,
 \begin{equation}\label{eq:defPL1L2} 
  L_1\times_{Y^I} L_2 \xrightarrow{P_{L_1,L_2}} \lim \left(\vcenter{
 \xymatrix@C=25pt@R=13pt{
  X  &  M_1\ar[l]^{p_1}\ar[r]_{p'_1}  & Y    &  M_2\ar[l]^{p_2}\ar[r]_{p'_2}  & Z
  \\
  X^I\ar[u]^{s_X} \ar[d]_{t_X} &
  & Y^I \ar[u]_ {s_Y} \ar[d]^{t_Y}
  & 
  & Z^I \ar[u]_{s_Z} \ar[d]^{t_Z}
  \\
  X  &  N_1\ar[l]_{q_1}\ar[r]^{q'_1}  & Y    &  N_2\ar[l]_{q_2}\ar[r]^{q'_2}  & Z }}\right).
\end{equation}
 The fact that both the fillers, $P_{L_1}$ and $P_{L_2}$, are fibrations, together with the universal property of pullbacks, then gives that $P_{L_1,L_2}$ is a fibration.

 \item We also have a naturally defined map, $P_{out}$, from the limit above to:
$$\lim \left( \vcenter{ \xymatrix@C=35pt@R=13pt{
  X  &  M_1\times_Y M_2\ar[l]^{\overline{p_1}}\ar[r]_{\overline{p_2'}}  & Z\\
  X^I\ar[u]^{s_X} \ar[d]_{t_X} &  & Z^I \ar[u]_{s_Z} \ar[d]^{t_Z}\\
X  &  N_1\times_Y N_2\ar[l]_{\overline{q_1}}\ar[r]^{\overline{q_2'}}  & Z
  }}\right).$$
  
From the fact that $\langle s_Y,t_Y\rangle\colon Y^I \to Y \times Y$ is a fibration, it follows that $P_{out}$ is a fibration.  Hence $\mathcal{P}_{{L_1,L_2}}=P_{out}\circ P_{L_1,L_2}$ is a fibration.
\end{enumerate}
\end{enumerate}

We now  analyse the composite window, $\W_1\#_0\W_2,$ via its spatial 2-slices,
\begin{equation}\label{eq:sp-2-slice}
\left[\hskip-1pt\begin{array}{c|c|c}
&(m_1,m_2)&\\
x\, &L _1\times_{Y^I}L_2& \, z\\
&(n_1,n_2)&\end{array}\hskip-1pt\right],\end{equation}
with, of course, $x\in X$, $z\in Z$, $(m_1,m_2)\in M_1\times_{Y}M_2$ and $(n_1,n_2)\in N_1\times_YN_2$. Taking this apart a  bit, just listing  elementary properties, we get, with the conventions as in the diagram in Equation \eqref{ash0}:\begin{itemize}[leftmargin=0.6cm]\item[]
$x=p_1(m_1)$, $z=p_2'(m_2)$, and there is some $y= p_1'(m_1)=p_2(m_2) \in Y$;
\item[] $x=q_1(n_1)$, $z=q_2'(n_2)$, and there is some $y'= q_1'(n_1)=q_2(n_2) \in Y$.\end{itemize}
The two elements, $y$ and $y'$, in $Y$, need not be the same here. However note that
\begin{itemize}[leftmargin=0.6cm]
\item[] 
$m_1 \in \{x|M_1|y\}$, \quad   $m_2 \in \{y|M_2|z\}$,\quad $n_1 \in \{x|N_1|y'\}$, and $n_2 \in \{y'|N_2|z\}$.
\end{itemize} Moreover, if $\ell=(\ell_1,\ell_2)$ is in the spatial 2-slice in \eqref{eq:sp-2-slice} then it must satisfy:\begin{itemize}[leftmargin=1cm]\item[--]
$\mathcal{P}(\ell)= (m_1,m_2)$;
\item[--] $\mathcal{Q}(\ell)= (n_1,n_2)$;
\item[--] $\ell= (\ell_1,\ell_2) \in L _1\times_{Y^I}L_2$, so $r_1(\ell_1)=l_2(\ell_2)$, which is in $Y^I$, of course, and so it is some path, $\gamma:[0,1]\to Y$; see diagram \eqref{ash0}.
\item[--]  Referring again to Equation \eqref{ash0}, the path $\gamma = r_1(\ell_1)$ starts at $s_Y\big(r_1(\ell_1)\big)= p'_1(P_1(\ell_1))=p'_1(m_1)= p_2(m_2)= y$, and, similarly, it ends at  $q_2(n_2)=y'$.
\end{itemize}
The following simple result then holds:
\begin{Lemma}\label{2-slice non-empty}
If we have that the element, $\ell= (\ell_1,\ell_2)\in L_1 \times L_2$, belongs to the spatial 2-slice in \eqref{eq:sp-2-slice}, then $r_1(\ell_1)$ is a path in $Y$, from $y:=p_2(m_2)$ to $y':=q_2(n_2)$, and so $y$ and $y'$ will be in the same path component of $Y$.

Furthermore, if the spatial 2-slice in \eqref{eq:sp-2-slice} is  non-empty,  then it has the same homotopy type as a spatial 2-slice of the form,
$$\left[\hskip-1pt\begin{array}{c|c|c}
&(m_1,m_2)&\\
x \, &L _1\times_{Y^I}L_2& \, z\\
&(n'_1,n'_2)&\end{array}\hskip-1pt\right],$$  with  $n'_1 \in \{x|N_1|y\}$ and $n'_2 \in \{y|N_2|z\}$, where $y=p_1'(m_1)=p_2(m_2).$
\end{Lemma}
\begin{proof} The first statement of the lemma follows from the discussion before. For the second, if $y$ and $y'$ are in the same path-component of $Y$,  we can always find a representative, $(n_1',n_2')$, of the path-component, $ \PC_{(n_1,n_2)}(\{x|N_1\times_Y N_2|z\})$, with $n_1'\in\{x|N_1|y\}$ and $n_2'\in\{y|N_2|z\}$. (This uses the last point of Lemma \ref{main3}, and applying the homotopy lifting property of the fibration,
$\{x|N_1\times_Y N_2|z\} \to Y,$
 to a path connecting $y'$ to $y$.) Since
$$\PC_{(n_1,n_2)}(\{x|N_1\times_Y N_2|z\}) = \PC_{(n_1',n_2')}(\{x|N_1\times_Y N_2|z\}),$$ and using the second point of Lemma \ref{lem:matrix2els}, we then have a homotopy equivalence:
$$\left[\hskip-1pt\begin{array}{c|c|c}
&(m_1,m_2)&\\
x\, &L _1\times_{Y^I}L_2& \, z\\
&(n_1,n_2)&\end{array}\hskip-1pt\right] \cong \left[\hskip-1pt\begin{array}{c|c|c}
&(m_1,m_2)&\\
x\, &L _1\times_{Y^I}L_2& \, z\\
&(n'_1,n'_2)&\end{array}\hskip-1pt\right].$$
\end{proof}
Finally we have the following, which will be used shortly. Recall that we denote $\Omega_y(Y)=\{ \gamma \in Y^{[0,1]}: \gamma(0),\gamma(1)=y\}$, the loop space of $Y$  at $y$.
\begin{Lemma}\label{lem:fibtoOmega}
 Let $x \in X$, $y\in Y$ and $z \in Z$. Also fix $m_1 \in \{x|M_1|y\}$,  $m_2 \in \{y|M_2|z\}$, $n_1 \in \{x|N_1|y\}$ and $n_2 \in \{y|N_2|z\}$.
The map, $\Psi_{L_1,L_2}\colon L_1\times_{Y^I} L_2 \to Y^I$, in diagram \eqref{Equation:psiL1L2}, induces (by restricting domain and codomain) a fibration,
 $$\overline{\Psi}_{L_1,L_2}\colon \left[\hskip-1pt\begin{array}{c|c|c}
&(m_1,m_2)&\\
x\, &L _1\times_{Y^I}L_2& \, z\\
&(n_1,n_2)&\end{array}\hskip-1pt\right] \to \Omega_y(Y).
 $$
\end{Lemma}
\begin{proof}
Follows from the fact that the map, $P_{L_1,L_2}$, in \eqref{eq:defPL1L2} is a fibration.
\end{proof}
\subsubsection{The natural transformation associated to the horizontal composition of HF resolved 2-spans}\label{Hor-comp-resolved 2-cells}
We assume given  a diagram of spaces, HF spans, and HF resolved 2-spans, and their horizontal composite, as shown\footnote{The  underlying windows, $\W_1$ and $\W_2$, are  in \eqref{ash0}, and their horizontal composite  in \eqref{def:W1W2}.},
$$\xymatrix@C=15pt{X \ar@/^1pc/[rr]^{(p_1,M_1,p_1')}\ar@/_1pc/[rr]_{ {(q_1,N_1,q_1')}} & \Downarrow \W_1 & Y \ar@/^1pc/[rr]^{(p_2,M_2,p_2')}\ar@/_1pc/[rr]_{ {(q_2,N_2,q_2')}} & \Downarrow \W_2 & Z,
} \quad \textrm{ and } \quad \xymatrix{\mortS{X}{Z}{(\overline{p_1}, M_1\times_Y M_2,\overline{p_2'})}{(\overline{q_1}, N_1\times_Y N_2,\overline{q_2'})}{\W_1\#_0\W_2}}.$$

 Let $\x_X{\subseteq}  X$, $\y_Y {\subseteq}   Y$ and $\z_Z{\subseteq}  Z$ be subsets.
Definition \ref{def:2matrixels} gives natural transformations of profunctors, therefore 2-morphisms in $\vProfGrphf$, \begin{align*}
 \tHpb{\W_1}{\x_X}{\y_Y}&\colon   \Hpb^{M_1}_{(\x_X,\y_Y)} \To \Hpb^{N_1}_{(\x_X,\y_Y)} ,\\
  \tHpb{\W_2}{{\y_Y}}{\z_Y}&\colon  \Hpb^{M_2}_{({\y_Y},\z_Y)}  \To  \Hpb^{N_2}_{(\y_Y,\z_Z)},\\
 \intertext{and}
\tHpb{\W_1\#_0\W_2}{\x_X}{\z_Z}&\colon   \Hpb^{M_1\times_Y M_2}_{(\x_X,\z_Z)}  \To \Hpb^{N_1\times_Y N_2}_{(\x_X,\z_Z)}, \end{align*}
where, in more notational detail, $ \Hpb^{M_1}_{(\x_X,\y_Y)}\colon \pi_1(X,\x_X) \bto \pi_1(Y,\y_Y)  ,$ etc.

We will now prove the important fact that, with the assumption that $(Y,\y_Y)$ is $0$-connected,  and noting Proposition \ref{techfib3},  the natural transformation arising from ${\W_1\#_0\W_2}$ is obtained by horizontally composing the natural transformations given by $\W_1$ and $\W_2$.

A crucial fact that we use is  Lemma \ref{lem:fibtoOmega}. This will be used together with Theorem \ref{main2}. The notation and results in \S\ref{sec:defT} will also play a key role.

 \begin{Proposition}\label{lemhorcomp} Let $\x_X {\subseteq}  X$, $\y_Y {\subseteq}  Y$ and $\z_Z {\subseteq}  Z$ be subsets with $(Y,\y_Y)$ 0-connected, then
 $$\tHpb{\W_1\#_0 \W_2}{\x_X}{\z_Z}=\tHpb{\W_1}{\x_X}{\y_Y} \bullet \tHpb{\W_2}{\y_Y}{\z_Z},$$
as natural transformations of $\Vect$-profunctors,
$$\Hpb^{M_1}_{(\x_X,\y_Y)}\bullet \Hpb^{M_2}_{({\y_Y},\z_Y)} \To \Hpb^{N_1}_{(\x_X,\y_Y)}\bullet \Hpb^{N_2}_{(\y_X,\z_Y)}.  $$
 \end{Proposition}
Here we have abused notation, and noted that, by Proposition \ref{techfib3}, we have canonical natural isomorphisms,
 of $\Vect$-profunctors,
 \begin{align*}
 \Hpb_{(\x_X,\y_Y)}^{M_1}\bullet \Hpb^{M_2}_{(\y_Y,\z_Z)}  &\xRightarrow{\eta_{(\x_X,\y_Y,\z_Z)}^{M_1,M_2}} \Hpb_{(\x_X,\z_Z)}^{M_1\times_Y M_2},\intertext{and} \Hpb_{(\x_X,\y_Y)}^{N_1}\bullet \Hpb^{N_2}_{(\y_Y,\z_Z)} &\xRightarrow{\eta_{(\x_X,\y_Y,\z_Z)}^{N_1,N_2}} \Hpb_{(\x_X,\z_Z)}^{N_1\times_Y N_2}.
 \end{align*}
 \begin{proof}
We will prove that the following diagram of natural transformations,
$$\xymatrix@R=15pt{
&\displaystyle\int^{y \in {\y_Y}} \Hpb_{(\x_X,\y_Y)}^{M_1}(-,y) \otimes \Hpb_{(\y_Y,\z_Z)}^{M_2}(y,-)\ar@{=>}[dd]|{\tHpb{\W_1}{\x_X}{\y_Y} \bullet \tHpb{\W_2}{\y_Y}{\z_Z}}
\ar@{=>}[rrr]^>>>>>>>>>>>>>>>>>{\eta_{(\x_X,\y_Y,\z_Z)}^{M_1,M_2}}&&& \Hpb_{(\x_X,\z_Z)}^{M_1\times_Y M_2}
\ar@{=>}[dd]|{\tHpb{\W_1\#_0 \W_2}{\x_X}{\z_Z}}
\\\\
&\displaystyle\int^{y \in {\y_Y}} \Hpb_{(\x_X,\y_Y)}^{N_1}(-, y) \otimes \Hpb_{(\y_Y,\z_Z)}^{N_2}(y,-)\ar@{=>}[rrr]_>>>>>>>>>>>>>>>>>{\eta_{(\x_X,\y_Y,\z_Z)}^{N_1,N_2}}&&& \Hpb_{(\x_X,\z_Z)}^{N_1\times_Y N_2}\,, }
$$
commutes. To this end, we prove that, given $x \in \x_X$ and $z \in \z_Z$, the following diagram of linear maps commutes,
\begin{equation}\label{Ldiag}\xymatrixrowsep{0.3cm}
\begin{gathered}\xymatrix{{} \save[]+<0cm,0cm>*{
{\displaystyle {\bigoplus_{y \in {\y_Y}}} \Hp^{M_1}(x_,y) \otimes \Hp^{M_2}(y,z)}} \restore 
 \\\\
\displaystyle \int^{y \in {\y_Y}} \Hpb_{(\x_X,\y_Y)}^{M_1}(x,y) \otimes \Hpb_{(\y_Y,\z_Z)}^{M_2}(y,z)\ar@{->}[ddd]|{\big(\tHpb{\W_1}{\x_X}{\y_Y} \bullet \tHpb{\W_2}{\y_Y}{\z_Z}\big)_{(x,z)}}
\ar@{->}[rrr]^<<<<<<<<<<<<<<<{\big(\eta_{(\x_X,\y_Y,\z_Z)}^{M_1,M_2}\big)_{(x,z)}}
\ar@{<-}[uu]^>>>>>{proj}
&&& \Hpb_{(\x_X,\z_Z)}^{M_1\times_Y M_2}(x,z)
\ar@{->}[ddd]|{\big(\tHpb{\W_1\#_0 \W_2}{\x_X}{\z_Z}\big)_{(x,z)}}
\\ \\\\
\displaystyle \int^{y \in {\y_Y}} \Hpb_{(\x_X,\y_Y)}^{N_1}(x,y) \otimes \Hpb_{(\y_Y,\z_Z)}^{N_2}(y,z)\ar@{->}[rrr]_<<<<<<<<<<<<<<<<{\big(\eta_{(\x_X,\y_Y,\z_Z)}^{N_1,N_2}\big)_{(x,z)}}&&& \Hpb_{(\x_X,\z_Z)}^{N_1\times_Y N_2}(x,z),
}
\end{gathered}\end{equation}where $proj$ is the projection mentioned earlier, in Equation \eqref{proj:prof} on page \pageref{proj:prof}.

We prove that the diagram in \eqref{Ldiag} commutes by explicitly computing matrix elements.
So let $y,y'\in \y_Y$.  Let also $m_1\in \{x|M_1|y\}$ and  $m_2\in \{y|M_2|z\}$. Finally let $n_1\in \{x|N_1|y'\}$  and $n_2\in
 \{y'|N_2|z\}$. 
 
Consider the path in   \eqref{Ldiag} that passes through the top right corner,  yielding $$F^{tr}\colon \bigoplus_{y \in {\y_Y}} \Hp^{M_1}(x_,y) \otimes \Hp^{M_2}(y,z)  \to \Hpb_{(\x_X,\z_Z)}^{N_1\times_Y N_2}(x,z),  $$ a linear map. Combining  Equation \eqref{eq:2matrixels} with  Proposition \ref{techfib3}, the corresponding matrix elements of $F^{tr}$ are given by
\begin{multline} \label{eq:form:tr}
 \big\langle \PC_{m_1}(\{x|M_1|y\})\otimes \PC_{m_2}(\{y|M_2|z\})\, \big| F^{tr} \big |\, \PC_{(n_1,n_1)}(\{x|N_1\times_Y N_2 | z\}) \big \rangle\\=\chi^\pi \left (\left[\hskip-1pt\begin{array}{c|c|c}
&(m_1,m_2)&\\
x\, &L _1\times_{Y^I}L_2& \, z\\
&(n_1,n_2)&\end{array}\hskip-1pt\right]\right) \chi^\pi\big (\PC_{(n_1,n_2)}(\{x|N_1\times_Y N_2|z\} ) \big).
\end{multline}
By Lemma \ref{2-slice non-empty}, the spatial 2-slices on the right-hand-side of \eqref{eq:form:tr} are
empty if  $y$ and $y'$ are not in the same path component of $Y$. Hence,  if  $y$ and $y'$ are  in different components of $Y$, the  matrix elements of $F^{tr}$ in \eqref{eq:form:tr}  have value 0.

On the other hand, again by Lemma \ref{2-slice non-empty}, if $y$ and $y'$ are in the same path-component of $Y$,  we can  find a representative, $(n_1',n_2')$, of the path-component, $ \PC_{(n_1,n_2)}(\{x|N_1\times_Y N_2|z\})$, with $n_1'\in\{x|N_1|y\}$ and $n_2'\in\{y|N_2|z\}$.  Consequently, when the matrix elements in \eqref{eq:form:tr} are non zero,  we can suppose that $y'=y$.

We therefore  consider $y\in \y_Y$, $m_1\in \{x|M_1|y\}$,  $m_2\in \{y|M_2|z\}$,  $n_1\in \{x|N_1|y\}$  and $n_2\in
 \{y|N_2|z\}$, and compute the value of \eqref{eq:form:tr} in this case.
To this end,  we  apply Theorem \ref{main2} to the fibration, $\overline{\Psi}_{L_1,L_2}$, in Lemma \ref{lem:fibtoOmega}.

Let $\gamma_1,\dots,\gamma_r$, where $r=|\pi_1(Y,y))|=|\pi_0(\Omega_y(Y))|$, be representatives of the different path-components of $\Omega_y(Y)$, which we recall is homotopy finite (see Lemma \ref{OmegaiHF}). We will suppose, with no loss of generality, that $\gamma_1=\const_y$, the constant path at $y$. Using the notation in Definition \ref{space-2els}, we have
\begin{align}\allowdisplaybreaks\nonumber
  & \chi^{\pi}\left(\left[\hskip-1pt\begin{array}{c|c|c}
&(m_1,m_2)&\\
x\,&L _1\times_{Y^I}L_2& \,z\\
&(n_1,n_2)&\end{array}\hskip-1pt\right]\right) =\sum_{i=1}^r \chi^\pi\big(\overline{\Psi}_{L_1,L_2}^{-1}(\gamma_i)\big)\,\chi^\pi\big(\PC_{\gamma_i}(\Omega_y(Y))\big)\\
  \nonumber&=\sum_{i=1}^r \chi^\pi\left( \left[\hskip-1pt\begin{array}{c|c|c}&m_1&\\ \const_x&L_1&\gamma_i\\&n_1&\end{array}\hskip-1pt\right]\times   \left[\hskip-1pt\begin{array}{c|c|c}&m_2&\\ \gamma_i&L_2&\const_z\\&n_2&\end{array}\hskip-1pt\right] \right)\chi^\pi\big(\PC_{\gamma_i}(\Omega_y(Y))\big)\\\nonumber
   &=\sum_{i=1}^r \chi^\pi\left( \left[\hskip-1pt\begin{array}{c|c|c}&m_1&\\ \const_x&L_1&\gamma_i\\&n_1&\end{array}\hskip-1pt\right]\right)\chi^\pi\left(  \left[\hskip-1pt\begin{array}{c|c|c}&m_2&\\ \gamma_i&L_2&\const_z\\&n_2&\end{array}\hskip-1pt\right] \right)\chi^\pi(\PC_{\gamma_i}(\Omega_y(Y)))\\
 &=\sum_{i=1}^r  \chi^\pi\left( \left[\hskip-1pt\begin{array}{c|c|c}&m_1&\\ x&L_1&y\\&\Gamma^{N_1}_{{\langle c_x,\overline{\gamma_i}\rangle}}(n_1)&\end{array}\hskip-1pt\right]\right)\chi^\pi\left(  \left[\hskip-1pt\begin{array}{c|c|c}&m_2&\\ y&L_2&z\\&\Gamma^{N_2}_{{\langle\overline{\gamma_i},c_z\rangle }}(n_2)&\end{array}\hskip-1pt\right] \right)\chi^\pi(\PC_{\gamma_i}(\Omega_y(Y))).
   \label{calc1:end}\end{align}
In  step
   \eqref{calc1:end} we put $\const_x=c_x$ and $\const_z=c_z$, and  used Item \eqref{lem:act2ls!!} of Lemma   \ref{lem:act2ls}.
 
Next we apply Lemma \ref{Lem:chi_pull} to calculate the next term in the right-hand-side of \eqref{eq:form:tr}, that is, $\chi^\pi\big( \PC_{(n_1,n_2)}(\{x|N_1\times_Y N_2|z\})\big)$.  To this end, we consider  the commutative diagram, below, where the middle diamond is a pullback,
\begin{equation}\label{eq:-def P}
\vcenter{\xymatrix@R=4pt{ &&&
N_1{\times}_Y N_2\ar[dd]^P  \ar[dl]_{}\ar[dr]^{} \ar@/_1.5pc/[ddll]_>>>>>{\overline{q_1}} \ar@/^1.5pc/[ddrr]^>>>>>{\overline{q_2'}}\\ &&N_1\ar[dl]^{q_1}\ar[dr]_{q_1'}
&&N_2\ar[dl]^{q_2}\ar[dr]_{q_2'}\\  & X & &Y  & &Z  \,. }}
\end{equation}
Put $\{x|N_1\}:=q_1^{-1}(x)$ and $\{N_2|z\}:={q_2'}^{-1}(z)$,
 then $q_1'\colon N_1 \to Y$ and $q_2\colon N_2\to Y$ restrict to fibrations, $q_r\colon \{N_2|z\} \to Y$ and $q_l\colon \{x| N_1\} \to Y$; see Lemma \ref{All-fibrations}. Moreover, we have a pullback diagram, where $P_{x,z}\colon \{x|N_1\times_Y N_2|z\} \to Y$ is the unique map making the diagram commute,
$$\xymatrix@R=7pt{
&& \{x|N_1\times_Y N_2|z\}\ar[dd]^{P_{x,z}}\ar[dl]_>>>>{\proj_1} \ar[dr]^>>>>{\proj_2} \\ & \{x| N_1\} \ar[dr]_{q_l}&& \{N_2|z\}  \ar[dl]^{q_r} \\ && Y.
}$$

Note that $P_{x,z}$ is a fibration (see the proof of Lemma \ref{mult}). Its fibre at $y\in Y$ is
$$P_{x,z}^{-1}(y)= \{x|N_1|y\} \times \{y|N_2|z\}. $$

By using Lemma \ref{Lem:chi_pull}, given $n_1\in \{x|N_1|y\}$ and $n_2\in \{y|N_2|z\}$, then
\begin{multline}\label{calc2:end}
\chi^\pi\big( \PC_{(n_1,n_2)}(\{x|N_1\times_Y N_2|z\})
\big)\\=T_{(n_1,n_2)}^{ \{x|N_1\times_Y N_2|z\}} \, \chi^\pi\big(\PC_y(Y)\big)\, \chi^\pi\big(\PC_{n_1}(\{x|N_1|y\})\big)\, \chi^\pi\big(\PC_{n_2}(\{y|N_2|z\})\big),
\end{multline}
where  $T_{(n_1,n_2)}^{ \{x|N_1\times_Y N_2|z\}}$ is
the cardinality of the orbit of the path-component, $$\PC_{(n_1,n_2)}\big( P_{x,z}^{-1}(y) \big) \in \hpi_0\big(P_{x,z}^{-1}(y)\big),$$ under the
 action, $\trl$, of $\pi_1(Y,y)$ on  $\hpi_0\big(P_{x,z}^{-1}(y)\big)$, as in Lemma \ref{main_techfib_loops}.
Note that $$\hpiz\big ( P_{x,z}^{-1}(y)\big) = \hpiz\big (  \{x|N_1|y\} \times \{y|N_2|z\}\big) \cong \hpiz\big(\{x|N_1|y\}\big) \times \hpiz\big( \{y|N_2|z\}\big),$$
and that the action of $\pi_1(Y,y)$ on $\hpi_0\big(P_{x,z}^{-1}(y)\big)$, is derived, in the obvious way, from the product of the actions, of $\pi_1(Y,y)$ on $\hpi_0(\{x|N_1|y\})$ and on $\hpi_0(\{y|N_2|z\})$, arising from the fibrations  $q_l\colon \{x| N_1\} \to Y$ and $q_r\colon \{N_2|z\} \to Y$; see Equation \eqref{eq:prop-hol}.

Going back to \eqref{eq:form:tr}, we now put \eqref{calc1:end} and \eqref{calc2:end} together. Two more observations are needed. Recall $r=|\pi_0(\Omega_y(Y))|=|\pi_1(Y,y)|$.

\begin{enumerate}[leftmargin=0.6cm]

\item We have homotopy equivalences, for  each $i \in \{1,\dots, r\}$,
\begin{align*}
\PC_{n_1}(\{x|N_1|y\}) &\cong  \PC_{{\Gamma^{N_1}_{{\langle \const_x,\overline{\gamma_i}\rangle }}(n_1)}}(\{x|N_1|y\}), 
\intertext{and}
\PC_{n_2}(\{y|N_2|z\}) &\cong  \PC_{{\Gamma^{N_2}_{{\langle \overline{\gamma_i},\const_z\rangle}}(n_2)}}(\{y|N_2|z\}). 
\end{align*}
These are induced by the homotopy equivalences, 
$$\Gamma^{N_1}_{{\langle \const_x,\overline{\gamma_i}\rangle}}\colon \{x|N_1|y\}\to  \{x|N_1|y\} \quad \textrm{ and } \quad \Gamma^{N_2}_{{\langle \overline{\gamma_i},\const_z\rangle}}
\colon \{y|N_2|z\}\to  \{y|N_2|z\} .$$
{(We are using Lemma \ref{pphi-cmpt} here.)}

\item By Lemma \ref{path_omega} below, all path-components, $\PC_{\gamma_i}(\Omega_y(Y))$,
of the loop space, of $Y$ at $y$, are homotopic. We therefore have
$$\chi^\pi\big(\PC_{\gamma_i}(\Omega_y(Y))\big)=\chi^\pi(\Omega_y(Y))/|\pi_1(Y,y)|, \textrm{ for each } i\in \{1, \dots, r\}. $$
\end{enumerate}

Putting everything together, we have:
\begin{align*}
\big\langle& \PC_{m_1}(\{x|M_1|y\})\otimes \PC_{m_2}(\{y|M_2|z\}) \, |F^{tr} | \, \PC_{(n_1,n_1)}(\{x | N_1\times_Y N_2 | z\}) \big \rangle\\&=
T_{(n_1,n_2)}^{ \{x|N_1\times_Y N_2|z\}} \chi^\pi(\PC_y(Y)) \frac{\chi^\pi(\Omega_y(Y))}{|\pi_1(Y,y)|}
\\ &\hspace*{12mm}\sum_{i=1}^{|\pi_1(Y,y)|} \big \langle \PC_{m_1}(\{x|M_1|y\})|  (\tHpb{\W_1}{\x_X}{\y_Y})_{(x,y)} |   \PC_{{\Gamma^{N_1}_{(\const_x,\overline{\gamma_i})}(n_1)}}(\{x|N_1|y\}) \big\rangle\\
&\,\,\hspace*{25mm}\big\langle \PC_{m_2}(\{y|M_2|z\})|  \big(\tHpb{\W_2}{\y_Y}{\z_Z}\big)_{(y,z)} |   \PC_{{\Gamma^{N_2}_{(\overline{\gamma_i},\const_z)}(n_2)}}(\{y|N_2|z\}) \big\rangle.
\end{align*}

Using Equation \eqref{eq:prop-hol}, the action, $\trl$, of $\pi_1(Y,y)$ on the set of path-components of the fibre, $P_{x,z}^{-1}(y)$, at $y$, of the fibration,
$P_{x,z}\colon  \{x | N_1 \times_Y N_2 | z\} \to Y$, is the product of the  actions, $\trl$, of $\pi_1(Y,y)$ on $\hpi_0(\{x|N_1|y\})$ and on $\hpi_0(\{y|N_2|z\}$, derived from the fibrations  $q_l\colon \{x| N_1\} \to Y$ and $q_r\colon \{N_2|z\} \to Y$. We hence have
\begin{multline*}
\big\langle \PC_{m_1}(\{x|M_1|y\})\otimes \PC_{m_2}(\{y|M_2|z\})\, | F^{tr} |\,  \PC_{(n_1,n_1)}(\{x | N_1\times_Y N_2 | z\}) \big \rangle\\ =\frac{T_{(n_1,n_2)}^{ \{x|N_1\times_Y N_2|z\}}}{|\pi_1(Y,y)|} \sum_{g\in \pi_1(Y,y)} \phantom{--------------------}\\ \big \langle \PC_{m_1}(\{x|M_1|y\})|    (\tHpb{\W_1}{\x_X}{\y_Y})_{(x,y)}|   \PC_{n_1}(\{x|N_1|y\})\trl g \big \rangle\\
\qquad \ \qquad  \ \quad \big \langle \PC_{m_2}(\{y|M_2|z\})| \big(\tHpb{\W_2}{\y_Y}{\z_Z}\big)_{(y,z)} |   \PC_{n_2}(\{y|N_2|z\})\trl g \big\rangle.
\end{multline*}
Note that we also  have used  that, by Lemma \ref{OmegaiHF},  $\chi^\pi(\PC_y(Y)) {\chi^\pi(\Omega_y(Y))}=1$.

We now let $\mathcal{T}_{(n_1,n_2)}^{ \{x|N_1\times_Y N_2|z\}}$  denote the $\pi_1(Y,y)$-orbit of the element,  $$\big(\PC_{n_1}(\{x|N_1|y\}), \PC_{n_2}(\{y|N_2|z\})\big) \in \hpi_0(\{x|N_1|y\})\times\hpi_0(\{y|N_2|z\}), $$
hence,
$|\mathcal{T}_{(n_1,n_2)}^{ \{x|N_1\times_Y N_2|z\}}|={T}_{(n_1,n_2)}^{ \{x|N_1\times_Y N_2|z\}}$.
Using the elementary fact that, if a finite group, $G$,   acts on a set,  then given any pairs of elements, $k$ and $l$, in the same orbit, the cardinality of $\{g \in G: k \trl g=l\}$ is that of the stabiliser subgroup of $k$, then, on applying this to the case of $G$ being  $\pi_1(Y,y)$, and, invoking  the orbit-stabiliser theorem, we have,
\begin{multline}\label{eq:latter-horcomp}
\big\langle \PC_{m_1}(\{x|M_1|y\})\otimes \PC_{m_2}(\{y|M_2|z\}) \, | F^{tr} |\, \PC_{(n_1,n_1)}(\{x| N_1\times_Y N_2 |y\}) \big \rangle
\\ =\sum \big  \langle \PC_{m_1}(\{x|M_1|y\})|  \big(\tHpb{\W_1}{\x_X}{\y_Y}\big)_{(x,y)} |   \PC_{n_1'}(\{x|N_1|y\}) \big \rangle\\
\hspace*{3cm} \big\langle \PC_{m_2}(\{y|M_2|z\})|  \big(\tHpb{\W_2}{\y_Z}{\z_Z}\big)_{(y,z)} |   \PC_{n_2'}(\{y|N_2|z\})\big \rangle,
\end{multline}
where  the sum is indexed by the set of elements of the $\pi(Y,y)$-orbit, namely those $$\big (\PC_{n_1'}(\{x|N_1|y\}), \PC_{n_2'}(\{y|N_2|z\})\big)\in  \mathcal{T}_{(n_1,n_2)}^{ \{x|N_1\times_Y N_2|z\}}.$$

Recall that $x\in \x_X$, $y\in \y_Y$ and  $z \in \z_Z$, and  also  $m_1\in \{x|M_1|y\}$,  $m_2\in \{y|M_2|z\}$, and $n_1\in \{x|N_1|y\}$,  $n_2\in
 \{y|N_2|z\}$. The formula in Equation \eqref{eq:latter-horcomp} gives exactly the corresponding matrix element,
$$ \big\langle \PC_{m_1}(\{x|M_1|y\})\otimes \PC_{m_2}(\{y|M_2|z\})\, \big| F^{bl} \big |\, \PC_{(n_1,n_2)}(\{x| N_1\times_Y N_2 | z\}) \big \rangle,$$
of the linear map associated to  the  path in \eqref{Ldiag} passing through
 the bottom left corner. This follows from \eqref{hor_comp:etas}, because if  $x \in \x_X$ and $z \in \z_Z$,  the linear bijection,
 $$\hskip -1cm\xymatrix{
 &\displaystyle \int^{y \in {\y_Y}} \Hpb_{(\x_X,\y_Y)}^{N_1}(x,y) \otimes \Hpb_{(\y_Y,\z_Z)}^{N_2}(y,z)\ar@{->}[rrr]^>>>>>>>>>>>>>>>>>{\big(\eta_{(\x_X,\y_Y,\z_Z)}^{N_1,N_2}\big)_{(x,z)}}&&& \Hpb_{(\x_X,\z_Z)}^{N_1\times_Y N_2}(x,z),
}$$
  is such that given  $y\in {\y_Y}$, $n_1 \in \{x|N_1|y\}$ and $n_2 \in \{y|N_2|z\}$, it sends the equivalence class of
 $$\PC_{n_1}(\{x|N_1|y\}) \otimes \PC_{n_2}(\{y|N_2|z\})\in  \bigoplus_{y \in {\y_Y}} \Hpb_{(\x_X,\y_Y)}^{N_1}(x,y) \otimes \Hpb_{(\y_Y,\z_Z)}^{N_2}(y,z) $$
 to 
 $\PC_{(n_1,n_2)}(\{x|N_1\times_Y N_2|z\}).$
Lemma \ref{lem:profmaps} in \S\ref{some:lemmas:coends} is  useful to translate between the categorical and the combinatorial languages.

Now suppose that $y,y'\in \y_Y$ are not in the same path-component in $Y$. If $m_1\in \{x|M_1|y\}$,  $m_2\in \{y|M_2|z\}$, $n_1\in \{x|N_1|y'\}$,  and $n_2\in
 \{y'|N_2|z\}$, we already saw that
  $$ \big\langle \PC_{m_1}(\{x|M_1|y\})\otimes \PC_{m_2}(\{y|M_2|z\})\, \big| F^{tr} \big |\, \PC_{(n_1,n_2)}(\{x | N_1\times_Y N_2 | z\}) \big \rangle=0.$$
 Applying the second point of Lemma \ref{lem:profmaps}, it also follows that,
 $$ \big\langle \PC_{m_1}(\{x|M_1|y\})\otimes \PC_{m_2}(\{y|M_2|z\})\, \big| F^{bl} \big |\, \PC_{(n_1,n_2)}(\{x | N_1\times_Y N_2 | z\}) \big \rangle=0.$$

Therefore, diagram \eqref{Ldiag} commutes as required. 
 \end{proof}

In the previous proof, we used the following lemma.
\begin{Lemma}\label{path_omega} Let $Y$ be any CGWH space {(which need not be  HF)}.
 All path components of $\Omega_y(Y)=\{\gamma \in Y^I\colon s_Y(\gamma)=t_Y(\gamma)=y\}$ are homotopic.
\end{Lemma}
 
\begin{proof}  Let $\mathcal{P}_y(Y)=\{\gamma \in Y^I: s_Y(\gamma)=y\}$, a path-connected space. We have a fibration,
 $t_y\colon \mathcal{P}_y(Y) \to Y$, induced by $t_Y\colon Y^I \to Y$. Clearly $\Omega_y(Y)$ is the fibre of  $t_y$ at $y$. By Lemma \ref{pphi-cmpt}, all path-components of $\Omega_y(Y)$ are homotopy equivalent.
 \end{proof}

 \begin{Remark}Recall Equation \eqref{eq:simpler-matrix-2-els} for the matrix elements of $\tHpc{\W}{x}{y}$. That the diagram in  \eqref{Ldiag} commutes can also be proven in the following way.
We use the notation in \eqref{ash0}.
Let
 $y\in \y_Y$,  $m_1\in \{x|M_1|y\}$ and  $m_2\in \{y|M_2|z\}$, and $n_1\in \{x|N_1|y\}$, $n_2\in \{y|N_2|z\}$. Let $\mathcal{P}_{L_1,L_2}$ be the filler of $\W_1\#_0\W_2$. Since the fillers, $P_{L_1}\colon L_1 \to \fr(\W_1)$  and   $P_{L_2}\colon L_2 \to \fr(\W_2), $   are fibrations, it follows that we have a fibration, $Q_!$, where the constant paths in $x$ and $z$ are denoted $c_x$ and $c_z$,
 \begin{multline*}
 Q_!\colon \left[\hskip-1pt\begin{array}{c|c|c}
&(m_1,m_2)&\\
x\, &L _1\times_{Y^I}L_2& \, z\\
&[(n_1,n_2)]&\end{array}\hskip-1pt\right]=\mathcal{P}_{L_1,L_2}^{-1}\big ( (m_1,m_2), c_x, c_z, \PC_{(n_1,n_2)}(\{x | N_1\times_Y N_2 | z\}) \big)\\\to  \mathcal{P}_y(Y),
\end{multline*}
sending $(\ell_1,\ell_2)$ to  $r_1(\ell_1)= l_2(\ell_2)$. The base space, $\mathcal{P}_y(Y)$, of the fibration contracts onto the constant path $\const_y$, so $\chi^\pi\left( \mathcal{P}_y(Y)\right) =1$. We  have another fibration,
$Q_{-}\colon Q_{!}^{-1}(\const_y)\to \{x|N_1|y\}\times \{y|N_2|z\}$, such that $(\ell_1,\ell_2)\mapsto \big ( Q_1(\ell_1),Q_2(\ell_2)\big)$.

Recall that images of path-components under fibrations are path-components themselves, of the base. Choose  representatives, $(a_1^i,a_2^i) \in \{x|N_1|y\}\times \{y|N_2|z\} $, of the path components of $\{x|N_1|y\}\times \{y|N_2|z\}$ that are images, under $Q_-$, of path-components of  $Q_{!}^{-1}(\const_y)$.  There are  only finite of these,
 and they are all in the  $\pi(Y,y)$-orbit, $\mathcal{T}_{(n_1,n_2)}^{ \{x|N_1\times_Y N_2|z\}}$, of $\PC_{n_1}\big (\{x|N_1|y\}\big)\times \PC_{n_2}\big(\{y|N_2|z\}\big)\in \hpiz(\{x|N_1|y\})\times \hpiz(\{y|N_2|z\})$.
 We then have:
$$Q_-\big (  Q_{!}^{-1}(\const_y)\big)=\bigsqcup_i \PC_{a_1^{i}}\big (\{x|N_1|y\}\big)\times \PC_{a_2^{i}}\big(\{y|N_2|z\}\big).  $$
\begin{Lemma}For each $i$, we have, using the notation in Remark \ref{rem:simpler-matrix-2-els},
$$Q_-^{-1} \Big( \PC_{a_1^{i}}\big (\{x|N_1|y\}\big)\times \PC_{a_2^{i}}\big(\{y|N_2|z\}\big)\Big ) =\left [ \begin{array}{c|c|c}
&m_1&\\
x\, &L_2  & \, y\\
&[a_1^i]&\end{array}\hskip-1pt\right]\times \left[\hskip-1pt\begin{array}{c|c|c}
&m_2&\\
y \, & L_2 & \, z\\
&[a_2^i]&\end{array}\hskip-1pt\right].$$\end{Lemma}
\begin{proof}The $\supseteq$ incusion  follows since
 \begin{align*}
 \PC_{a_1^i}(\{x|N_1|y\})\times \PC_{a_2^i}(\{y|N_2|z\})&\subseteq  \PC_{(a_1^i,a_2^i)}(\{x | N_1\times_Y N_2 | z\})\\ &=  \PC_{(n_1,n_2)}(\{x | N_1\times_Y N_2 | z\}).
 \end{align*}
  The $\subseteq$ inclusion is immediate.
\end{proof}

As images of path-components under fibrations are path-components, we have:
\begin{Lemma}
 An element $\PC_{n_1'}(\{x|N_1|y\})\times \PC_{n_2'}(\{y|N_2|z\})\in \mathcal{T}_{(n_1,n_2)}^{ \{x|N_1\times_Y N_2|z\}}$   is a path-component of  $Q_-\big(Q_{!}^{-1}(\const_y)\big)$ if, and only if:
 $$\left [ \begin{array}{c|c|c}
&m_1&\\
x\, &L_2  & \, y\\
&[n'_1]&\end{array}\hskip-1pt\right] \neq \emptyset \textrm{ and }\ \left[\hskip-1pt\begin{array}{c|c|c}
&m_2&\\
y \, & L_2 & \, z\\
&[n'_2]&\end{array}\hskip-1pt\right]\neq \emptyset.$$
\end{Lemma}

\noindent Using Lemma \ref{main1} to $Q_-$, and applying the proof of Lemma \ref{Lem:chi_pull}, we then get:
 \begin{align*}
 \big\langle & \PC_{m_1}(\{x|M_1|y\})\otimes \PC_{m_2}(\{y|M_2|z\})\, \big| F^{tr} \big |\, \PC_{(n_1,n_2)}(\{x | N_1\times_Y N_2 | z\}) \big \rangle\\&=\chi^\pi \left (\left[\hskip-1pt\begin{array}{c|c|c}
&(m_1,m_2)&\\
x\, &L _1\times_{Y^I}L_2& \, z\\
&[(n_1,n_2)]&\end{array}\hskip-1pt\right]\right)=\chi^\pi\big ( Q_{!}^{-1}(\const_y )\big)\\
&=\chi^\pi\Big ( Q_-^{-1} \big (\bigsqcup_i \PC_{a_1^{i}}(\{x|N_1|y\})\times \PC_{a_2^{i}}(\{y|N_2|z\})\big) \Big)\\
&=\sum_{i} \chi^{\pi}\Big (Q_{-}^{-1} \big( \PC_{a_1^i}(\{x|N_1|y\})\times  \PC_{a_2^i}(\{y|N_2|z\}) \big )\Big)\\
&=\sum_i \chi^\pi \left (
 \left[\hskip-1pt\begin{array}{c|c|c}
&m_1&\\
x\, &L_2  & \, y\\
&[a_1^i]&\end{array}\hskip-1pt\right]\times \left[\hskip-1pt\begin{array}{c|c|c}
&m_2&\\
y \, & L_2 & \, z\\
&[a_2^i]&\end{array}\hskip-1pt\right]\right ) =\sum_i \chi^\pi \left (
 \left[\hskip-1pt\begin{array}{c|c|c}
&m_1&\\
x\, &L_2  & \, y\\
&[a_1^i]&\end{array}\hskip-1pt\right]\right) \chi^\pi \left ( \left[\hskip-1pt\begin{array}{c|c|c}
&m_2&\\
y \, & L_2 & \, z\\
&[a_2^i]&\end{array}\hskip-1pt\right]\right )\\
&=
\sum \chi^\pi \left (
 \left[\hskip-1pt\begin{array}{c|c|c}
&m_1&\\
x\, &L_2  & \, y\\
&[n'_1]&\end{array}\hskip-1pt\right]\right) \chi^\pi \left ( \left[\hskip-1pt\begin{array}{c|c|c}
&m_2&\\
y \, & L_2 & \, z\\
&[n'_2]&\end{array}\hskip-1pt\right]\right )
\\ &=\sum \Big  \langle \PC_{m_1}(\{x|M_1|y\})|  \big(\tHpb{\W_1}{\x_X}{\y_Y}\big)_{(x,y)} |   \PC_{n_1'}(\{x|N_1|y\}) \Big \rangle\\
& \hspace*{3cm} \Big\langle \PC_{m_2}(\{y|M_2|z\})|  \big(\tHpb{\W_2}{\y_Z}{\z_Z}\big)_{(y,z)} |   \PC_{n_2'}(\{y|N_2|z\})\Big \rangle,
\end{align*}
where  the last two sums are indexed by the elements of $\mathcal{T}_{(n_1,n_2)}^{ \{x|N_1\times_Y N_2|z\}}$.
 \end{Remark}

\subsection{The vertical composition of HF resolved 2-spans}
 
 Checking that horizontal composition translates via the profunctor construction to horizontal composition of  the corresponding natural transformations required some counting arguments, the corresponding checks for the vertical composition require other methods.
\subsubsection{Preliminaries for the vertical composition of HF resolved 2-spans}\label{sec:prelims_vert}
Let $X$ be a CGWH space.  By Example \ref{ids_ex}, there is a fibrant span, ${(s_X,X^I,t_X)} \colon X  \to  X$, from which we constructed the identity of $X$ in the category $\HFb$.

The composite, (see Definition \ref{main3def}),
$(s_X,X^I,t_X)\bullet (s_X,X^I,t_X)\colon X \to X$ is  the fibrant span,
$(\overline{s_X}, X^I \times_X X^I, \overline{t_Y})\colon X\to X.$ As before,
 $$X^I \times_X X^I=\{(\gamma,\gamma') \in X^I \times X^I\mid \gamma(1)=\gamma'(0)\},$$
and we recall that  $\overline{s_X}(\gamma,\gamma')=\gamma(0)$ and $\overline{t_X}(\gamma,\gamma')=\gamma'(1)$.

We consider the homeomorphism, $\mathcal{F}_X\colon X^I\times_X X^I \to X^I$, defined as
\begin{equation}\label{eq:concatenation}
 \mathcal{F}_X(\g,\g')(t)=\begin{cases}
                         \gamma(2t), &t \in [0,1/2],\\
                          \gamma'(2t-1),& t \in [1/2,1].
                        \end{cases}
 \end{equation}
Clearly  $\mathcal{F}_X$ makes the diagram below commute,
$$\xymatrix@R=10pt{&&X^I{\times}_{X} X^I\ar[dd]^{\mathcal{F}_X}\ar[dl]_{\overline{s_X}} \ar[dr]^{\overline{t_X}}  \\
&X & & X.\\ && X^I \ar[ul]^{s_X} \ar[ur]_{t_X}
}
$$
We thus have that $\mathcal{F}_X$ is an isomorphism (of fibrations) over $X\times X$.

\subsubsection{The vertical composition of HF resolved 2-spans}\label{sev_vert_compresolved}
Let $X$ and $Y$ be HF spaces. Consider HF fibrant resolved 2-spans of form
\begin{align*}
 \W_2\colon \big ((p_2,M_2,p_2')\colon X \to Y\big) \To\big ((p_1,M_1,p_1')\colon X \to Y\big),
 \intertext{and}
 \W_1\colon \big ((p_1,M_1,p_1')\colon X \to Y\big) \To\big ((q_1,N_1,q_1')\colon X \to Y\big),
\end{align*}
where the windows, $\W_1$ and $\W_2$, are:
\begin{equation}\label{WindowsW1W2}
\W_1=\vcenter{\xymatrix@R=15pt{
  X  &  M_1\ar[l]_{p_1}\ar[r]^{p'_1}  & Y\\
  X^I\ar[u]^{s_X} \ar[d]_{t_X} & L_1 \ar[l]_{l_1}\ar[r]^{r_1}\ar[u]^{P_1}\ar[d]_{Q_1}& Y^I \ar[u]_{s_Y} \ar[d]^{t_Y}\\
  X  &  N_1\ar[l]^{q_1}\ar[r]_{q'_1}  & Y}} \quad \textrm{and} \quad  \W_2= \vcenter{\xymatrix@R=15pt{
  X  &  M_2\ar[l]_{p_2}\ar[r]^{p'_2}  & Y\\
  X^I\ar[u]^{s_X} \ar[d]_{t_X} & L_2 \ar[l]_{l_2}\ar[r]^{r_2}\ar[u]^{P_2}\ar[d]_{Q_2}& Y^I \ar[u]_{s_Y} \ar[d]^{t_Y}\\
  X  &  M_1\ar[l]^{p_1}\ar[r]_{p'_1}  & Y \,.} }
\end{equation}

 We want to  construct a vertical composite, $\W_2 \#_1 \W_1$, which should be  a HF resolved 2-span, such that $\W_2 \#_1 \W_1\colon (p_2,M_2,p_2') \To  (q_1,N_1,q_1')$.  We do this in two steps.
Firstly, as when we constructed the horizontal composition of HF fibrant resolved 2-spans, we perform the obvious pullback along the common
horizontal spans of $\W_1$ and $\W_2$.
This yields the HF fibrant window,
 $$
  \W_2 \#_1' \W_1 =\vcenter{\xymatrix@C=50pt@R=15pt{
  X  &  M_2\ar[l]_{p_2}\ar[r]^{p'_2}  & Y\\
  X^I\times_XX^I\ar[u]^{\overline{s_X}} \ar[d]_{\overline{t_X}} & {L_2\times_{M_1} L_1} \ar[l]_{\overline{l_1}}\ar[r]^{\overline{r_2}}\ar[u]^{\overline{P_2}}\ar[d]_{\overline{Q_1}}& Y^I\times_Y Y^I \ar[u]_{\overline{s_Y}} \ar[d]^{\overline{t_Y}}\\
  X  &  N_1\ar[l]^{q_1}\ar[r]_{q'_1}  & Y.}}
$$
Here, given $(\ell_2,\ell_1) \in L_2\times_{M_1} L_2$, we have written
\begin{align*}\overline{l_1}\big( (\ell_2,\ell_1)\big)&=\big(l_2(\ell_2),l_1(\ell_1)\big),  &\overline{r_2}\big( (\ell_2,\ell_1)\big)&=\big(r_2(\ell_2),r_1(\ell_1)\big),\\
\overline{P_2}\big( (\ell_2,\ell_1)\big)&=P_2(\ell_2), &\overline{Q_1}\big( (\ell_2,\ell_1)\big)&=Q_1(\ell_1).
\end{align*}

To prove that $\W_2 \#_1'\W_1$ is a HF fibrant window, we can use  the same argument that we used for  the horizontal composition of HF resolved 2-spans; see \S \ref{sec:horcomp}.

We now need to `adjust' the left and right vertical spans of $\W_2 \#_1'\W_1$. To this end, we use the homeomorphisms, below, of fibrant spans; as in   \S\ref{sec:prelims_vert},
$$\vcenter{\xymatrix@R=7pt{&X^I{\times}_X X^I\ar[dd]^{\mathcal{F}_X}\ar[dl]_{\overline{s_X}} \ar[dr]^{\overline{t_X}}  \\
X & & X\\ & X^I \ar[ul]^{s_X} \ar[ur]_{t_X}
}}
 \textrm{\quad and \quad}
\vcenter{\xymatrix@R=7pt{&Y^I{\times}_Y Y^I\ar[dd]^{\mathcal{F}_Y}\ar[dl]_{\overline{s_Y}} \ar[dr]^{\overline{t_Y}}  \\
Y & & Y,\\ & Y^I \ar[ul]^{s_Y} \ar[ur]_{t_Y}
}}
$$
and the commutative diagram,
$$
  \xymatrix@C=30pt@R=15pt{&
  &X  &  M_2\ar[l]_{p_2}\ar[r]^{p'_2}  & Y&\\
& X^I \ar[ur]^{s_X} \ar[dr]_{t_Y}   &X^I\times_X X^I\ar[u]^{\overline{s_X}} \ar[d]_{\overline{t_X}}\ar[l]_{\mathcal{F}_X} & {L_2\times_{M_1} L_1} \ar[l]_{\overline{l_1}}\ar[r]^{\overline{r_2}}\ar[u]^{\overline{P_2}}\ar[d]_{\overline{Q_1}}& Y^I\times_Y Y^I \ar[u]_{\overline{s_Y}} \ar[d]^{\overline{t_Y}}
\ar[r]^{\mathcal{F}_Y}& Y^I. \ar[ul]_{s_Y} \ar[dl]^{t_Y}
\\
 & &X  &  N_1\ar[l]^{q_1}\ar[r]_{q'_1}  & Y}
$$
This yields what will be called the \emph{vertical composite of the fibrant resolved 2-spans,} $\W_2$ and $\W_1$, as displayed below,
 $$ \W_2\#_1 \W_1:=  \vcenter{\xymatrix@C=60pt@R=18pt{
  X  &  M_2\ar[l]_{p_2}\ar[r]^{p'_2}  & Y\\
  X^I\ar[u]^{{s_X}} \ar[d]_{{t_X}} & {L_2\times_{M_1} L_1} \ar[l]_{\mathcal{F}_X\circ\overline{l_1}}\ar[r]^{\mathcal{F}_Y\circ \overline{r_2}}\ar[u]^{\overline{P_2}}\ar[d]_{\overline{Q_1}}& Y^I \ar[u]_{{s_Y}} \ar[d]^{{t_Y}}\\
  X  &  N_1\ar[l]^{q_1}\ar[r]_{q'_1}  & Y.}}$$
By construction, the window, $\W_2\#_1 \W_1$,   is a HF fibrant resolved 2-span such that $$\W_2\#_1\W_1\colon \big ((p_2,M_2,p_2')\colon X \to Y\big) \Longrightarrow \big ((q_1,N_1,q_1')\colon X \to Y\big).$$

Now recall Definitions \ref{equiv-of-HF-spans}  and  \ref{def:vspan_el}. Let $x \in X$ and $y\in Y$.
\begin{Lemma}\label{vert_span_comp}  We have an isomorphism   of HF spans, from $\{x|M_2|y\}$ to $\{x|N_1|y\}$, $$[x| \W_2\#_1\W_1|y] \cong  [x|\W_2|y] \bullet  [x|\W_1|y]. $$
 \end{Lemma}
\begin{proof} Follows from the fact that the concatenation of two paths, $\gamma\colon a \to  a$ and $\gamma'\colon a \to a$, see \eqref{eq:concatenation}, is a constant path if, and only if, $\gamma,\gamma'=\const_a$.
\end{proof}
Let now $\x_X$ and $\y_Y$ be subsets of $X$ and $Y$. Consider the  $\Vect$-profunctors,
\begin{align*}
\Hpb^{M_2}_{(\x_X,\y_Y)}\colon \pi_1(X,\x_X)^{\op} \times \pi_1(Y,\y_Y) &\to \Vect,\\
\Hpb^{M_1}_{(\x_X,\y_Y)}\colon \pi_1(X,\x_X)^{\op} \times \pi_1(Y,\y_Y) &\to \Vect,\\ 
\Hpb^{N_1}_{(\x_X,\y_Y)}\colon \pi_1(X,\x_X)^{\op} \times \pi_1(Y,\y_Y) &\to \Vect.
\end{align*} The
 HF fibrant resolved 2-spans $\W_2$ and $\W_1$
give rise to natural transformations of $\Vect$-profunctors, as in Definition \ref{def:2matrixels}.
\begin{Lemma}\label{Lem:_vert_comp}
The diagram of natural transformations of profunctors commutes,
$$\xymatrix@C=45pt{\Hpb^{M_2}_{(\x_X,\y_Y)} \ar@{=>}[rd]_{\tHpb{\W_2\#_1\W_1}{\x_X}{\y_Y}\,\,\,\,\,} \ar@{=>}[r]^>>>>>>>>>>{\tHpb{\W_2}{\x_X}{\y_Y}} &\Hpb^{M_1}_{(\x_X,\y_Y)}\ar@{=>}[d]^{\tHpb{\W_1}{\x_X}{\y_Y}} \\ & \Hpb^{N_1}_{(\x_X,\y_Y)}. } $$
\end{Lemma}
\begin{proof}
 Let $(x,y) \in \x_X \times \y_Y$. We claim that the diagram in $\Vect$ commutes:
$$\xymatrix@C=70pt{\Hpb^{M_2}_{(\x_X,\y_Y)}(x,y) \ar@{->}[dr]_{\big(\tHpb{\W_2\#_1\W_1}{\x_X}{\y_Y}\big)_{(x,y)}\,\,\,\,} \ar@{->}[r]^<<<<<<<<<<<<{\big(\tHpb{\W_2}{\x_X}{\y_Y}\big)_{(x,y)}} &\Hpb^{M_1}_{(\x_X,\y_Y)}(x,y)  \ar@{->}[d]^{\big(\tHpb{\W_1}{\x_X}{\y_Y}\big)_{(x,y)}} \\ &\Hpb^{N_1}_{(\x_X,\y_Y)}(x,y). } $$
This follows by combining Lemma \ref{vert_span_comp} with the $s=0$ case of Lemma \ref{mult}, and using Definition \ref{def:2matrixels} and Equation \eqref{eq:2matrixels}.
\end{proof}

\subsection{Towards horizontal and vertical identities}
We still have to examine if the suggested compositions, both horizontal and vertical, have identities.

\subsubsection{The vertical identity}\label{vert_idsspan}
Let $(p,M,q)\colon X \to Y$ be an HF fibrant span. We define the following window,
$$ \id^2_{(p,M,q)}:=\vcenter{\xymatrix@R=13pt@C=30pt{
  X  &  M\ar[l]_p\ar[r]^{q}  & Y\\
  X^I\ar[u]^{s_X} \ar[d]_{t_X} & M^I  \ar[l]_{l_X}\ar[r]^{r_Y}\ar[u]^{s_M}\ar[d]_{t_M}& Y^I \ar[u]_{s_Y} \ar[d]^{t_Y}\\
  X  &  M\ar[l]^p\ar[r]_{q}  & Y.}}
$$
  Here, given $\gamma\colon I \to M$, we put $l_X(\gamma)=p\circ \gamma$ and $r_X(\gamma)=q\circ \gamma$.

  This definition is motivated by the construction of the bicategory $\tcob{n}$, below in Subsection \ref{tcobn}.  In particular, the  diagram above is a function space counterpart of the  vertical identity of a  cobordism, as given  in Item \eqref{vert_id} on page \pageref{vert_id}.

 \begin{Remark}\label{rem:vertical_id_naive} We do not know whether $\id^2_{(p,M,q)}$ is, in general, a fibrant window or not. Whenever $\id^2_{(p,M,q)}$ is a fibrant window (which holds in all cases required in the construction of the once-extended Quinn TQFT, in Section \ref{sec:def-once-extended}),  we note that it will be a HF fibrant resolved 2-span, connecting $(p,M,q)$ to itself. This is because $M^I$, $X^I$ and $Y^I$ are all HF, as they are homotopic to $M$, $X$ and $Y$, respectively.
 \end{Remark}

   \begin{Lemma}\label{lem:vert_ids} Let $(p,M,q)\colon X \to Y$ be a HF fibrant span. Suppose that  $\id^2_{(p,M,q)}$  is a HF fibrant resolved 2-span, therefore connecting $(p,M,q)$ to itself. In this situation, given any subsets $\x_X \subseteq  X$ and $\y_Y \subseteq Y$, the natural transformation,
   $$\tHpb{\id^2_{(p,M,q)}}{\x_X}{\y_Y} \colon \Hpb^M_{(\x_X,\y_Y)}\implies \Hpb^M_{(\x_X,\y_Y)} , $$
is the identity natural transformation.    
   \end{Lemma}
\begin{proof}
 Let $x \in \x_X$ and $y \in \y_Y$. If $m,n \in \{x|M|y\}$, then
 \begin{multline*}\big\langle  \PC_{m}(\{x|M|y\}) \,|\, (\tHpb{\id^2_{(p,M,q)}}{\x_X}{\y_Y} )_{(x,y)} \, | \,\PC_{n}(\{x|M|y\}) \big\rangle
 \\=
 \chi^\pi \left( \left[\hskip-1pt\begin{array}{c|c|c}
&m&\\
x&M^I& y\\
&n&\end{array}\hskip-1pt\right] \right)\chi^\pi\big (\PC_n(\{x|M|y\} ) \big).
 \end{multline*}
Now note that, by Example \ref{ids_ex} and Remark \ref{space_els_props}, we have an HF fibrant span,
$$\xymatrix@R=1pt{   & \{x|M|y\}^I \ar[dl]_{s_{\{x|M|y\}}} \ar[dr]^{t_{\{x|M|y\}}} \\ 
\{x|M|y\} && \{x|M|y\}, }$$
and that
$$ \left[\hskip-1pt\begin{array}{c|c|c}
&m&\\
x&M^I& y\\
&n&\end{array}\hskip-1pt\right]=\big\langle s_{\{x|M|y\}}, t_{\{x|M|y\}}\big\rangle^{-1}(m,n),$$ so we only need to apply the first part of Lemma \ref{ids}, for the case $s=0$.
\end{proof}

\subsubsection{Horizontal identities and unitors}\label{horidsspan}

Let $X$ be a HF space and  consider the HF fibrant span, $(s_X,X^I,t_X)\colon X \to X$. Let $\x_X\subseteq  X$.
The $\Vect$-profunctor,
$$\Hpb_{(\x_X,\x_X)}^{X^I}\colon \pi_1(X,\x_X)^\op \times \pi_1(X,\x_X)\to \Vect, $$is  such that, given $x,y \in \x_X$,  where $\Lin\colon \Sets \to \Vect$ denotes free vector space,
$$\Hpb_{(\x_X,\x_X)}^{X^I}(x,y)=\Lin \big(\hpi_0(\{x|X^I|y\})\big)\cong \Lin\big(\hom_{\pi_1(X,\x_X)}(x,y)\big).$$

A holonomy map, \smash{$\lambda^{X^I}$},  for the fibration,  \smash{$\langle s_X,t_X\rangle\colon  X^I \to X \times X$}, can  be constructed so that, given paths in $X$,  $\g\colon x \to y$, $\g_l\colon x' \to x$ and $\g_r\colon y \to y'$, then
$$\Gamma_{(\overline{\gamma_l}, \gamma_r)}^{X^I}(\gamma)={\gamma_l}*\gamma *\gamma_r, $$
the concatenation of ${\gamma_l}, \gamma$ and $\gamma_r$, each fitting into a third of $[0,1]$.
In particular,
$$\Hp^{X^I}_{(\x_X,\x_X)}([\gamma_l],[\gamma_r])([\gamma])=[\gamma_l][\gamma][\gamma_r].$$  The profunctor associated with the identity span,  $(s_X,X^I,t_X)\colon X \to X,$ is,   therefore, canonically isomorphic to the horizontal identity, $\Id_{\pi_1(X,\x_X)}$, of $\pi_1(X,\x_X)$, in the bicategory $\vProfGrp$;  see Example \ref{identity profunctor}.

Continuing this approach, we now discuss a type of `would be' unitor for HF fibrant spans, given by  certain HF fibrant resolved 2-spans, and also how the \textit{bona fide} unitors in the bicategory $\vProfGrp$ can be obtained from the former by computing the associated natural transformations of  profunctors. This will be crucial for constructing the once-extended Quinn TQFT in Section \ref{sec:def-once-extended}. We will only discuss left unitors, as the case of right unitors is  analogous.

 Let $X$ and $Z$ be HF spaces. Consider a HF fibrant span,  $(p,M,q)\colon X \to Z$. We suppose, and this will be the case in all settings that we need for constructing the once-extended Quinn TQFT in Section \ref{sec:def-once-extended}, that the following conditions, \eqref{first-condition} -- \eqref{final-condition}, are satisfied.
   (These conditions may seem a bit mysterious at this stage. However, as we will see later, they arise naturally from the construction of the unitors of a cobordism in the bicategory $\tcob{n}$, when looking at their function space counterparts; see Subsection \ref{tcobn}, especially item  \eqref{l_unitor}, starting on page \pageref{l_unitor}.)
\begin{enumerate}[leftmargin=0.6cm]\item\label{first-condition} We have a homeomorphism, $\Phi\colon X^I\times_X M \to M$, and a commuting diagram,
 $$\xymatrix@R=5pt{ &&M \ar[dl]_p\ar[dr]^{q}\\  & X & &Z. \\ && X^I {\times}_X M \ar[uu]^\Phi \ar[ul]^>>>>>{p'}\ar[ur]_>>>>>{q'} }
$$
(This homeomorphism is an analogue of a collar of the boundary of a manifold, when considering spaces of functions on manifolds.) Here we used the pullback diamond inside the commutative diagram,
$$\xymatrix@R=0pt{ &&& 
X^I {\times}_X M
\ar[dl]_{}\ar[dr]^{} \ar@/^1.3pc/[ddrr]^>>>>>{q'} \ar@/_1.3pc/[ddll]_>>>>>{p'}\\ &&X^I\ar[dl]^{s_X}\ar[dr]_{t_X}&&M \ar[dl]^{p}\ar[dr]_{q}\\  & X & &X  & &Z. }
$$
\item Let $x,x' \in X$ and $z \in Z$. Given a path, $x \ra{\gamma} x'$, in $X$ and $m\in \{x'|M|z\}$, then $\Phi(\gamma,m) \in \{x|M|z\}$ is in the same path-component as $\Gamma^{M}_{\langle \overline{\gamma},\const_z\rangle}(m)\in \{x|M|z\}$.
Here, recalling  the notation in Lemma \ref{main_techfib}, $$\Gamma^{M}_{\langle \overline{\gamma}, \const_z\rangle } \colon \{x'|M|z\} \to \{x|M|z\}$$ is defined from the fibration $\langle p,q \rangle \colon M \to X\times Z$. \item\label{final-condition} The  following window is fibrant, where
$l_X(\gamma)=p\circ \gamma$ and $r_Z(\gamma)=q \circ \gamma$,
$$\boldsymbol{\lambda}^{(p,M,q)}_{X}:=\vcenter{\xymatrix@R=15pt@C=30pt{
  X  &  X^I\times_X M\ar[l]_>>>>>{p'}\ar[r]^>>>>>{q'}  & Z\\
  X^I\ar[u]^{s_X} \ar[d]_{t_X} & M^I  \ar[l]^{l_X}\ar[r]_{r_Z}\ar[u]^{\Phi^{-1}\circ s_M}\ar[d]_{t_M}& Z^I \ar[u]_{s_Z} \ar[d]^{t_Z}\\
  X  &  M\ar[l]^p\ar[r]_{q}  & Z. }}
   $$
\end{enumerate}

\begin{Lemma}\label{preserves unitors}  Let $X$ and $Z$ be HF spaces. Let a HF fibrant span,  $(p,M,q)\colon X \to Z$, satisfy  the conditions, \eqref{first-condition} -- \eqref{final-condition}, just outlined. Suppose $\x_X\subseteq  X$ and $\z_Z \subseteq  Z$.

If $x,x' \in \x_X$, $\gamma \in \{x|X^I|x'\}$, $z \in \z_Z$, and $m \in \{x'|M|z\}$, $m'\in \{x|M|z\}$, then
\begin{multline*}
\big \langle \PC_{(\gamma,m)} \Big(\big\{x|X^I \times_X M|z\big\} \Big)\,\, \Big{|}\,\,
\Big(\tHpb{\boldsymbol{\lambda}^{(p,M,q)}_{X}}
{\x_X}{\z_Z} \Big)_{(x,z)}  \,\,\Big{|} \,\,\PC_{m'}\big(\{x|M|z \}\big)\big\rangle
\\=\begin{cases} 1,\textrm{ if }  \PC_{\Gamma^M_{\langle \overline{\gamma}, \const_z\rangle}(m)}\big(\{x|M|z\}\big)= \PC_{{m'}}\big(\{x|M|z\}\big),\\
 0, \textrm{ otherwise.}                                                                                                                                                                                                                                                             \end{cases}
 \end{multline*}
 In particular, when $(X,\x_X)$ is 0-connected, the natural transformation, \smash{$\tHpb{\boldsymbol{\lambda}^{(p,M,q)}_{X}}
{\x_X}{\z_Z} $}, of profunctors $\pi_1(X,\x_X) \bto \pi_1(Z,\z_Z)$,  gives the appropriate left-unitor
 in the bicategory $\vProfGrp$.
More precisely, the following diagram, of profunctors $\pi_1(X,\x_X) \bto \pi_1(Z,\z_Z)$ and natural transformations, commutes,
$$\xymatrix@C=30pt{ \Id_{\pi_1(X,\x_X)}\bullet {\Hpb^{M}_{(\x_X,\z_Z)}} \ar@{=>}[drrr]_{{ \lambda_{\pi_1(X,\x_X)}^{\Hpb^{M}_{(\x_X,\z_Z)}} }}\ar@{=>}[r]^{\cong}  & \Hpb_{(\x_X,\x_X)}^{X^I} \bullet {\Hpb^{M}_{(\x_X,\z_Z)}}  \ar@{=>}[rr]^>>>>>>>>>>>>>{\eta^{X^I,M}_{(\x_X,\x_X,\z_Z)}}_\cong && \Hpb_{(\x_X,\z_Z)}^{X^I \times_X M}\ar@{=>}[d]^{\tHpb{\boldsymbol{\lambda}^{(p,M,q)}_{X}}
{\x_X}{\z_Z}}\\
&&&
 \Hpb_{(\x_X,\z_Z)}^{M}.
} $$
\end{Lemma}
We note that here the first equivalence is discussed earlier in this section, and the second is in Lemma \ref{techfib3}.
\begin{proof}
The proof of the first statement is exactly as in the proof of Lemma \ref{lem:vert_ids}. The second statement follows by passing to the language of profunctors.
\end{proof}

\subsection{Comment and Summary}\label{sec:comment-summary}\label{sec:comm_summary}
There have been a lot of fairly technical results in this section and it is easy to lose track of what they say \emph{in toto}, so we will step back to look at why they were necessary in the form we gave.

In this section, we have constructed a bicategorical type of object, though not quite a bicategory, that we will from now on  denote by $\mathbf{2span}(HF)$, following the description starting in page \pageref{intro_2HF} of the Introduction. Similar constructions are in \cite{grandis:collaredII:2007,Mortoncospans,Torzewska,Torzewska-HomCobs}.

The objects of
 $\mathbf{2span}(HF)$ are homotopy finite spaces. Given homotopy finite spaces, $X$ and $Y$, the 1-cells, from $X$ to $Y$, are homotopy finite fibrant spans, $(p,M,q)\colon X \to Y$. We have a non-associative composition, $\bullet$, of 1-cells, obtained via the obvious pullback, see  Definition \ref{main3def}.  Each homotopy finite space $X$ has a `horizontal identity', given by the path-space fibrant span, $(s_X,X^I,t_X)\colon X \to X$.

  Given 1-cells $(p,M,p'),(q,N,q')\colon X \to Y$, the 2-cells, in $\mathbf{2span}(HF)$, connecting them, consist of homotopy finite resolved 2-spans (see \S \ref{sec:HF_f_r_2_s}),  $\W\colon  (p,M,p')  \To (q,N,q').$
Again, those 2-cells can be composed horizontally and vertically, as described in detail in  \S \ref{sec:horcomp} and \S \ref{sev_vert_compresolved}. None of these compositions is associative.

As discussed in \S \ref{vert_idsspan} and \S \ref{horidsspan}, if we apply certain restrictions on the 1-cells $(p,M,p')\colon X \to Y$ that we allow (which will be automatically satisfied in the cases arising in the construction of the once-extended Quinn TQFT, in Section \ref{sec:def-once-extended}), we then have `vertical identities', $\id_{(p,M,p')}\colon (p,M,p') \To (p,M,p')$, as well as `unitor 2-cells', as shown below, whenever a 1-cell comes equipped with the function space analogue of a collar neighbourbood of the boundary of a manifold,
\begin{align*}
 \boldsymbol{\rho}^{(p,M,q)}_{X} \colon (X \ra{(p,M,p')} Y) \bullet (Y \ra{ (s_Y,Y^I,t_Y)} Y)   &\To (X \ra{(p,M,p')} Y),\\
\boldsymbol{\lambda}^{(p,M,q)}_{X} \colon (X \ra{(s_X,X^I,t_X)} X) \bullet (X \ra{(p,M,p')} Y)  &\To (X \ra{(p,M,p')} Y).
\end{align*}
Throughout Section \ref{sec:homotopy_underpinning}, we constructed an `assignment', from now on denoted  $$\boldsymbol{\mathcal{H}}=\big(\pi_1(-,-),\Hp,\tH\big)\colon \mathbf{2span}(HF) \to \vProfGrphf,$$ more precisely a map of 2-truncated globular sets, that gives the following.
\begin{enumerate}[leftmargin=0.6cm]
 \item Each homotopy finite space, $X$, is sent to its fundamental groupoid, $\pi_1(X,X)$.
 \item  Given a homotopy finite fibrant span, $(p,M,p')\colon X \to Y$, we have a $\Vect$-profunctor, as defined in \S\ref{sec:profunctor_fib_span},
 $$\smash{\Hp\big (X \ra{(p,M,p')} Y \big) \colon \pi_1(X,X)^{\op} \times \pi_1(Y ,Y ) \to \Vect.}$$
\item Given a homotopy finite fibrant resolved 2-span, $\W\colon (p,M,p')\To (q,N,q')$, as above, we have a natural transformation, of functors $\pi_1(X,X)^{\mathrm{op}}\times \pi_1(Y,Y) \to \Vect$, as discussed  in \S \ref{sec_resolvedspans_to_nat_tranfs}, $$\tHp{\W} \colon \Hp\big ( X\ra{(p,M,p')} Y\big)  \To \Hp\big ( X\ra{(q,N,q')} Y\big).$$
\end{enumerate}

We proved in Section \ref{sec:homotopy_underpinning} that the assignment, $\boldsymbol{\mathcal{H}}\colon \mathbf{2span}(HF) \to \vProfGrphf,$  preserves all various compositions, and the horizontal identities in $\mathbf{2span}(HF)$ and in  $\vProfGrphf$, up to applying appropriate natural isomorphisms, and that  vertical identities, and unitors  likewise are preserved by $\boldsymbol{\mathcal{H}}$, whenever they exist.

There is a \textit{relative} variant of $\mathbf{2span}(HF)$, from now on denoted $\overline{\mathbf{2span}(HF)}$, where homotopy finite spaces, $X$, come equipped with subsets, $\x_X{\subseteq}  X$, such that $(X,\x_X)$ is 0-connected, and the rest of the  `bicategorical' structure of $\overline{\mathbf{2span}(HF)}$ is induced by that of $\mathbf{2span}(HF)$. We also saw in this section that  $\boldsymbol{\mathcal{H}}$ can be modified to a `assignment',  $\boldsymbol{\overline{\mathcal{H}}}\colon \overline{\mathbf{2span}(HF)} \to \vProfGrphf$, that gives the following.
\begin{enumerate}[leftmargin=0.6cm]
 \item Each pair, $(X,\x_X), $ is sent to the  fundamental groupoid, $\pi_1(X,\x_X)$.
 (We will, in the following section, furthermore suppose that $\x_X$ is finite, so given that $X$ is homotopy finite, it follows that $\pi_1(X,\x_X)$ is then a finite groupoid.)
 \item  Given a homotopy finite fibrant span, $(p,M,p')\colon X \to Y$, $\x_X {\subseteq}  X$, and $\y_Y {\subseteq}  Y$, we  have a 1-cell, $(p,M,p')\colon (X,\x_X) \to (Y,\y_Y)$, in $\overline{\mathbf{2span}(HF)}$, and a $\Vect$-profunctor, as defined in  {Notation} \ref {not_barH},
$$\Hpb_{(\x_X,\y_Y)}\big( X  \ra{ (p,M,p')} Y \big) \colon \pi_1(X,\x_X)^{\op}\times \pi_1(Y ,\y_{Y })\to \Vect,
$$
obtained by restricting \smash{$\Hp\big(X \ra{(p,M,p')} Y \big)$} to $\pi_1(X,\x_X)^{\op}\times \pi_1(Y ,\y_{Y })$.
\item Given $\x_X {\subseteq}  X$ and $\y_Y {\subseteq}  Y$, and a 2-cell in ${\mathbf{2span}(HF)}$,  $\W\colon (p,M,p')\To  (q,N,q')$, then Definition  \ref{def:2matrixels}, gives  a  natural transformation of profunctors,
 $$\tHpb{\W}{\x_X}{\y_Y}\colon \Hpb_{(\x_X,\y_Y)}\big ( X  \ra{ (p,M,p')} Y \big)  \To \Hpb_{(\x_X,\y_Y)}\big ( X  \ra{ (q,N,q')} Y \big).$$
\end{enumerate}
In this section, we also proved that, just as for $\boldsymbol{\mathcal{H}}\colon \mathbf{2span}(HF) \to \vProfGrphf$,  the relative version $\boldsymbol{\overline{\mathcal{H}}}\colon \overline{\mathbf{2span}(HF)} \to \vProfGrphf,$ preserves all the various compositions, plus the unitors and identities (when they exist) up to natural isomorphisms.

We also saw in \S \ref{sec:H-is-sym-monoidal}, buiding from lemmas \ref{Sym_mon:1} and \ref{def_chi}, that $\boldsymbol{{\mathcal{H}}}$, and similarly $\boldsymbol{\overline{\mathcal{H}}}$, takes the product / cartesian monoidal structure in $\mathbf{2span}(HF)$ to the tensor product in $\vProfGrphf$.  This will be discussed further in Subsection \ref{Quinn_is_sym_mon}, and written down in the language  of symmetric monoidal bifunctors. 

As already mentioned,  $\mathbf{2span}(HF)$ is not a symmetric monoidal bicategory. However, when we define our  once-extended versions of Quinn's TQFT, we will only use $\mathbf{2span}(HF)$ as a `half-way house' between the bicategory, $\tcob{n}$, of 2-cobordisms, that we will introduce in detail in the next section, and $\vProfGrphf$, the second part of this process being given by $\boldsymbol{\mathcal{H}}$. Even if $\mathbf{2span}(HF)$ is not a symmetric monoidal bicategory, the composite assignment, $\tcob{n}\to \vProfGrphf$ will still be  a symmetric monoidal bifunctor. More will revealed in the next section.

\section{Once-extended versions of Quinn's TQFT}\label{sec:def-once-extended}


Let $X$ be a space. In this section, we will frequently abbreviate $\iota^X_k\colon X \to X \times [0,1]$ to $\iota_k$, hence $\iota_k(x)=\iota_k^X(x)=(x,k)$, usually for $k=0$ or $k=1$.

\subsection{Conventions for the bicategory $\tcob{n}$}\label{tcobn} Let $n$ be a non-negative integer. The bicategory, $\tcob{n}$, is that of \emph{closed smooth $n$-manifolds, $(n+1)$-cobordisms between manifolds, and equivalence classes of $(n+2)$-extended cobordisms connecting $(n+1)$-cobordisms.} The details of the construction are in \cite[\S 3.1.2]{Schommer-Pries} and \cite{Mortoncospans}.
We will give an overview in what follows, so as to set out the conventions we  use. As in previous sections, we  make no assumption that orientations exist on the manifolds, cobordisms, nor now on the extended cobordisms.

The bicategory, $\tcob{n}$, is thus defined as follows.
\begin{enumerate}[leftmargin=0.41cm]
 \item The class of objects of $\tcob{n}$ is the class of all smooth closed (i.e. compact and with empty boundary) $n$-dimensional manifolds.
 \item Given objects, $\Sigma$ and $\Sigma'$, a 1-morphism, $(i,S,j)\colon \Sigma \to \Sigma'$, called a  \emph{$(n+1)$-cobordism}, is a cospan,  of smooth manifolds and smooth maps, as below,
 $$\xymatrix@R=-8pt{ \Sigma \ar[dr]^{i} &&  \Sigma' \ar[dl]_{j}\\ & S&}.
$$
This should be  such that $S$ is a compact smooth $(n+1)$-manifold, possibly with a non-empty boundary, and the universally defined map, $\langle i,j \rangle\colon \Sigma \sqcup \Sigma'\to S,$ gives a diffeomorphism, $\Sigma\sqcup \Sigma'\cong \d S$.
\item \label{compcob} The composition of the 1-morphisms $(i,S,j)\colon \Sigma \to \Sigma'$ and $(i',S',j')\colon \Sigma' \to \Sigma''$, denoted $(\overline{i},S\bullet S',\overline{j'}):=(\overline{i},S \sqcup_\Sigma S',\overline{j'})$, is given by considering the pushout, in $\CGWH$, included as the diamond in the commutative diagram, below,
$$\xymatrix@R=8pt{&\Sigma\ar@/_1.3pc/[ddrr]_<<<<<{\overline{i}}\ar[dr]^{i} && \Sigma'\ar[dl]_{j} \ar[dr]^{i'}&& \Sigma''. \ar@/^1.3pc/[ddll]^<<<<<{\overline{j'}} \ar[dl]_{j'}\\
&& S\ar[dr]^{k} && S'\ar[dl]_{k'}\\
&&& S\sqcup_{\Sigma'} S'
}$$
(As already mentioned, in this paper, we implicitly choose a natural realisation for all limits and colimits. In this case, we took the obvious choice, $S\sqcup_{\Sigma'} S'= \big((S\times \{0\})\cup (S' \times \{1\})\big)/j(s)\sim i'(s)$, for all $s \in \Sigma'$, with the quotient topology.)

We note that the  pushout is formed in $\CGWH$, so initially we forget the smooth structure on the given manifolds, and consider just their underlying topological spaces. The smooth structure on $S\sqcup_{\Sigma'} S'$ is then inserted afterwards.

As recalled in Subsection \ref{cobordism categories}, in order to give a smooth structure to  $S\sqcup_\Sigma S'$, we could, for instance, consider collars of $\Sigma'$ in $S$ and $S'$. However the collars are not part of the structure given here to cobordisms. This issue can be resolved as in \cite[\S 3.1.2 and \S3.2]{Schommer-Pries}, either by  considering ``halations'', where collars essentially become part of the cobordism information, or applying the  axiom of choice for classes, to endow each cobordism with appropriate collars. We will not say more on this issue (and essentially will ignore it when we come to compose  extended cobordisms, below).  We can safely do this as our constructions depend only on the underlying topological manifolds of the smooth manifolds.

\item Given closed smooth $n$-manifolds, $\Sigma_1$ and $\Sigma_2$, and cobordisms, $(i_1,S,i_2)\colon \Sigma_1 \to \Sigma_2$ and  $(i_1',S',i_2')\colon \Sigma_1 \to \Sigma_2$,
the 2-morphisms between them,
$$[\K]\colon \big((i_1,S,i_2)\colon \Sigma_1 \to \Sigma_2\big)\To \big((i_1',S',i_2')\colon \Sigma_1 \to \Sigma_2\big), $$  are given by equivalence classes of diagrams, in the category of manifolds and smooth maps, of the form \eqref{ex_extended} below, called \emph{$(n+2)$-extended cobordisms},
$$\K\colon \big((i_1,S,i_2)\colon \Sigma_1 \to \Sigma_2\big)\To \big((i_1',S',i_2')\colon \Sigma_1 \to \Sigma_2\big),$$
\begin{equation}\label{ex_extended}
\K= \vcenter{\xymatrix@C=2pc@R=2pc{
  \Sigma_1  \ar[r]^{i_1}\ar[d]_{\iota_0^{\Sigma_1}}  &  S\ar[d]^{i_S}   & \Sigma_2\ar[l]_{i_2} \ar[d]^{\iota_0^{\Sigma_2}}\\
  \Sigma_1\times I\ar[r]_{i_E}   &  K & \Sigma_2\times I\ar[l]^{i_W} \\
  \Sigma_1 \ar[u]^{\iota_1^{\Sigma_1}} \ar[r]_{i_1'}  &  S'\ar[u]_{i_{S'}}    & \Sigma_2.\ar[l]^{i_2'} \ar[u]_{\iota_1^{\Sigma_2}}
 }}
  \end{equation}

 Here $K$ is a compact smooth $(n+2)$-manifold with corners, called the \emph{support of $\K$}. (The $E$ on the middle right pointing map is there to indicate that the arrow  is `pointing' east in the diagram, and the $W$, similarly, is pointing west.)
  
  Dually to the ideas of windows and fibrant resolved 2-spans, see Definition \ref{def:window} and \S \ref{sec:HF_f_r_2_s}, the \emph{frame}, $\fr(\K)$, of an extended cobordism, $\K$, as in \eqref{ex_extended},
  is defined to be
  \begin{equation}\label{eq:frame} \
  \fr(\K):=
  \colim \left(
  \vcenter  
  {
  \xymatrix@R=2pc@C=1.7pc{
  \Sigma_1  \ar[r]_{i_1}\ar[d]_{\iota_0^{\Sigma_1}}  &  S   & \Sigma_2\ar[l]^{i_2} \ar[d]^{\iota_0^{\Sigma_2}}\\
  \Sigma_1\times I    &  & \Sigma_2\times I \\
  \Sigma_1 \ar[u]^{\iota_1^{\Sigma_1}} \ar[r]^{i_1'}  &  S'   & \Sigma_2\ar[l]_{i_2'} \ar[u]_{\iota_1^{\Sigma_2}}}
  }
  \right).
\end{equation}
  We have a canonically defined map, $f_\K\colon \fr(\K) \to K$, as  for HF resolved 2-spans, hence also called the \emph{filler of $\K$}. It is required that $f_\K$ provides a diffeomorphism $\fr(\K) \cong \d(K)$, the boundary of the manifold with corners, $K$.
  \item Two extended cobordisms, $\K,\K'\colon \big((i_1,S,i_2)\colon \Sigma_1 \to \Sigma_2\big)\To \big((i_1',S',i_2')\colon \Sigma_1 \to \Sigma_2\big) $, so with the same frame, are called \emph{equivalent} if there exists a diffeomorphism, $f\colon K \to K'$, between the supports of $\K$ and $\K'$, making the diagram  commute,
  $$\xymatrix@R=1.7pc{&\fr{(\K)}\ar[d]_{\id}\ar[r]^{f_\K} & K\ar[d]^f\\
                &\fr{(\K')}\ar[r]_{f_{\K'}} & K' . } $$
\item The horizontal and vertical compositions of extended cobordisms are done via the obvious horizontal and vertical pushouts, dually to the  case of HF resolved 2-spans, as treated  in \S \ref{sec:horcomp} and \S \ref{sev_vert_compresolved}; see also \cite[\S 3.1.2]{Schommer-Pries} and \cite{Mortoncospans,morton:cohomological:2015}.  

These compositions of extended cobordisms descend to  equivalence classes,  defining  horizontal and vertical compositions of 2-morphisms in $\tcob{n}$.

\item\label{vert_id} Given an $(n+1)$-cobordism, $(i_1,S,i_2)\colon \Sigma_1 \to \Sigma_2$,  its \emph{vertical identity} is the equivalence class of the extended cobordism,
$$\id^2_{(i_1,S,i_2)}:= \vcenter{\xymatrix@C=3pc@R=1.7pc{
  \Sigma_1  \ar[r]_{i_1}\ar[d]_{\iota_0^{\Sigma_1}}  &  S \ar[d]^{\iota_0^S}   & \Sigma_2\ar[l]^{i_2} \ar[d]^{\iota_0^{\Sigma_2}}\\
  \Sigma_1\times I\ar[r]_{i_E}   &  S\times I & \Sigma_2\times I\ar[l]^{i_W} \\
  \Sigma_1 \ar[u]^{\iota_1^{\Sigma_1}} \ar[r]^{i_1}  &  S\ar[u]_{\iota_1^S}    & \Sigma_2.\ar[l]_{i_2} \ar[u]_{\iota_1^{\Sigma_2}}
 }}$$
Here $i_E(s,t)=(i_1(s),t)$, for $s \in \Sigma_1$ and $t \in I$, and, similarly,  $i_W(s,t)=(i_2(s),t)$.

\item Given a smooth compact $n$-manifold, $\Sigma$, the \emph{horizontal identity} of $\Sigma$ is $$\id^1_\Sigma:=(\iota_0^\Sigma,\Sigma \times I,\iota_1^\Sigma)\colon \Sigma \to \Sigma.$$

\item \label{l_unitor} Given an $(n+1)$-cobordism, $(i_1,S,i_2)\colon \Sigma_1 \to \Sigma_2$, we have \emph{left} and \emph{right unitor} $(n+2)$-extended cobordisms,
$$\boldsymbol{\lambda'}_{\Sigma_1}^{(i_1,S,i_2)}\colon (\iota_0^{\Sigma_1},\Sigma_1 \times I,\iota_1^{\Sigma_1}) \bullet  (i_1,S,i_2) \To (i_1,S,i_2),$$\and
$$\boldsymbol{\rho'}_{\Sigma_2}^{(i_1,S,i_2)}\colon (i_1,S,i_2)\bullet (\iota_0^{\Sigma_2},\Sigma_2 \times I,\iota_1^{\Sigma_2}) \To (i_1,S,i_2).$$
The support of both is $S\times I$. We will explain the  construction of the left unitor extended cobordism,
$\boldsymbol{\lambda'}_{\Sigma_1}^{(i_1,S,i_2)}$. The construction of the right unitor extended cobordism is  similar.

Consider the $(n+1)$-cobordism,
$$\big (i_1',   (\Sigma_1 \times I) \sqcup_{\Sigma_1} S, i_2'\big)= (\iota_0^{\Sigma_1},\Sigma_1 \times I,\iota_1^{\Sigma_1}) \bullet  (i_1,S,i_2),$$ and also an explicit isomorphism of cospans,
$$\xymatrix@R=9pt{ & S\ar[dd]^{\overline{\Phi}}& \\
\Sigma_1\ar[ru]^{i_1}\ar[rd]_<<<<<<{i_1'} && \Sigma_2\ar[lu]_{i_2}\ar[ld]^<<<<<<{i_2'}\\
&( \Sigma_1 \times I) \sqcup_{\Sigma_1} S.&
}
$$
(Note that to construct such a homeomorphism, $\overline{\Phi}\colon S \to ( \Sigma_1 \times I) \sqcup_{\Sigma_1} S $, we need a collar for the inclusion of $\Sigma_1$ in  $S$.)
The left unitor extended cobordism is defined by tweaking the vertical identity of $ (i_1,S,i_2)\colon \Sigma_1 \to \Sigma_2$, as shown in the diagram below (cf. \S\ref{horidsspan} where a dual construction is discussed),
$$ \boldsymbol{\lambda'}_{\Sigma_1}^{(i_1,S,i_2)}  :=\vcenter{\xymatrix@C=3pc@R=1.7pc{
  \Sigma_1  \ar[r]_<<<<<<<{i_1'}\ar[d]_{\iota_0^{\Sigma_1}}  &   (\Sigma_1 \times I)\sqcup_{\Sigma_1} S \ar[d]^{\iota_0^S \circ \overline{\Phi}^{-1}}   & \Sigma_2\ar[l]^<<<<<<{i_2'} \ar[d]^{\iota_0^{\Sigma_2}}\\
  \Sigma_1\times I\ar[r]|{i_E}   &  S\times I & \Sigma_2\times I\ar[l]|{i_W} \\
  \Sigma_1 \ar[u]^{\iota_1^{\Sigma_1}} \ar[r]^{i_1}  &  S\ar[u]_{\iota_1^S}    & \Sigma_2. \ar[l]_{i_2} \ar[u]_{\iota_1^{\Sigma_2}} }}
 $$

The equivalence classes of the left and right unitor extended cobordisms give  the left and right unitors in $\tcob{n}$, as denoted below,
$$\boldsymbol{\lambda}_{\Sigma_1}^{(i_1,S,i_2)}=[\boldsymbol{\lambda'}_{\Sigma_1}^{(i_1,S,i_2)}]\colon (\iota_0^{\Sigma_1},\Sigma_1 \times I,\iota_1^{\Sigma_1}) \bullet  (i_1,S,i_2) \To (i_1,S,i_2),$$ and
$$\boldsymbol{\rho}_{\Sigma_2}^{(i_1,S,i_2)}=[\boldsymbol{\rho'}_{\Sigma_2}^{(i_1,S,i_2)}]\colon (i_1,S,i_2)\bullet (\iota_0^{\Sigma_2},\Sigma_2 \times I,\iota_1^{\Sigma_2}) \To (i_1,S,i_2).$$
 
\end{enumerate}


In addition to the above basic structure,  we note that, in the classical setting, the category, $\cob{n}$, has the structure of  a symmetric monoidal category with the coproduct /  disjoint union, $\sqcup$, as the tensor product, as recalled in \S \ref{moncob}, and that the extended form, $\tcob{n}$, similarly, has a symmetric monoidal bicategory structure, again having $\sqcup$ as its tensor product.   An explicit proof is in \cite[\S 3.1.4]{Schommer-Pries}. We will revisit this structure in Subsection \ref{sec:sym_mon}, particularly \S\ref{sec:mon_tcob}.

\subsection{A once-extended version of Quinn's TQFT}\label{sec:once extended}

As before, $n$ will be  a non-negative integer, and  $\Bc$   a homotopy finite space. These will be the standard assumptions throughout this section.

Consider an $(n+1)$-cobordism, between closed smooth $n$-manifolds, as in Subsections \ref{cobordism categories} and \ref{tcobn}, viewed  as a cospan in the category $\CGWH$,
 $$(i,S,j):=\Big(\vcenter{\xymatrix@R=-5pt{\Sigma \ar[dr]_{i} &&  \Sigma' \ar[dl]^{j}\\ & S} }\Big),
$$so the nodes only encode the data of the underlying topological manifolds.
Applying the contravariant mapping space  functor, $\Bc^{(-)}\colon \CGWH \to \CGWH$,  sends this cospan to a span in $\CGWH$, whose nodes contain the corresponding spaces of maps from the topological manifolds into $\Bc$,
$$(i^*,\Bc^S,j^*):=\Big(\vcenter{\xymatrix@R=-5pt{ & \Bc^{S}\ar[dl]_{i^*} \ar[dr]^{j^*}&\\
 \Bc^\Sigma &&  \Bc^{\Sigma'}}}\Big).
$$
\begin{Lemma}\label{span of function spaces} This span, $(i^*,\Bc^S,j^*)$, of function spaces  is a fibrant span in which all the  spaces appearing are homotopy finite.\end{Lemma}
 
Now consider an extended  $(n+2)$-cobordism with 2-cospan diagram as follows (as before this is a diagram in $\CGWH$, and similarly for all subsequent diagrams),
\begin{equation}\label{diag-cob-forR}
\K=  \vcenter{\xymatrix@C=3pc@R=2pc{
  \Sigma_1  \ar[r]_{i_1}\ar[d]_{\iota_0}  &  S\ar[d]^{i_N}   & \Sigma_2\ar[l]^{i_2} \ar[d]^{\iota_0}\\
  \Sigma_1\times I\ar[r]^{i_E}   &  K & \Sigma_2\times I\ar[l]_{i_W} \\
  \Sigma_1 \ar[u]^{\iota_1} \ar[r]^{i'_1}  &  S'\ar[u]_{i_S}    & \Sigma_2.\ar[l]_{i'_2} \ar[u]_{\iota_1}
 }}
\end{equation}
 Applying the same contravariant functor, $\Bc^{(-)}$, to $\K$ gives
a dual `window',
\begin{equation}\label{diag-dcob-forR}
\Bc^{\K}:=\vcenter{\xymatrix@C=5pc@R=2pc{
  \Bc^{\Sigma_1}  \ar@{<-}[r]_{i_1^*}\ar@{<-}[d]_{\iota_0^*}  &  \Bc^{S}\ar@{<-}[d]^{i_N^*}   & \Bc^{\Sigma_2}\ar@{<-}[l]^{i_2^*} \ar@{<-}[d]^{\iota_0^*}\\
  \Bc^{\Sigma_1\times I}\ar@{<-}[r]^{i_E^*}   &  \Bc^K & \Bc^{\Sigma_2\times I}\ar@{<-}[l]_{i_W^*} \\
  \Bc^{\Sigma_1} \ar@{<-}[u]^{\iota_1^*} \ar@{<-}[r]^{{i'^*_1}}  &  \Bc^{S'}\ar@{<-}[u]_{i_S^*}    & \Bc^{\Sigma_2}.\ar@{<-}[l]_{{i'^*_2}} \ar@{<-}[u]_{\iota_1^*}
 }}
\end{equation}

Recalling the definition of fibrant resolved 2-spans in \S \ref{sec:HF_f_r_2_s}, we obtain:
\begin{Lemma}\label{2-span of function spaces} The window, $\Bc^{\K}$, of mapping spaces is a fibrant resolved 2-span,  in which all the spaces appearing are homotopy finite.
{Furthermore, the application of  $\Bc^{(-)}$ preserves the compositions of all cobordisms and extended cobordisms, sending them to the corresponding compositions of fibrant spans, as in Definition \ref{main3def}, and of fibrant resolved 2-spans, as in \S \ref{sec:horcomp} and \S \ref{sev_vert_compresolved}.}
The vertical units in $\tcob{n}$, as well as the horizontal identities, and unitors, are also sent to those  of \S\ref{vert_idsspan} and \S\ref{horidsspan}, when passing to the mapping spaces.
\end{Lemma}
We note that $\CGWH$ is cartesian closed, so $\Bc^{\Sigma_1\times I}\cong (\Bc^{\Sigma_1})^I$, canonically.

\begin{proof} We will prove  Lemmas \ref{span of function spaces} and  \ref{2-span of function spaces} together as the proofs are related. It may be useful to compare with the proof of our earlier Lemma \ref{1dual}.  We will continue to refer to the mapping space picture as being `dual' to the other one.

The fibrancy of the dual span  follows from the fact that  the inclusion of $\Sigma \sqcup \Sigma'\cong \d  S $ into $S$ is a cofibration and, similarly, for the dual window,  the inclusion of $f_\K\colon \fr(\K)\cong \d K$ into $K$ is a cofibration, so the dual map, $f_\K^*\colon \Bc^K \to \Bc^{\fr(\K)}\cong \fr(\Bc^\K)$, is a fibration. (In that last step, we again used the fact that $\CGWH$ is cartesian closed, so the mapping space contravariant functor, $\Bc^{(\_)}\colon \CGWH \to \CGWH$, send colimits to limits.) For the latter reason, all compositions are preserved (up to isomorphism) when going from cobordisms and extended cobordisms to fibrant spans and {fibrant resolved 2-spans.}
 
 In order to prove that all spaces in $\Bc^\K$ are homotopy finite, we use the fact that all compact smooth manifolds can be given the structure of a  finite CW-complex, and use, once again, Lemma \ref{lem:mappingfromCW}.
 
 Vertical units, horizontal units, and unitors, are preserved by construction.
\end{proof}

We also note the following, that once again follows from the fact that $\CGWH$ is cartesian closed.
\begin{Lemma}\label{B is monoidal}Given two manifolds, $\Sigma_1$ and $\Sigma_2$, in $\CGWH$, we have a natural isomorphism,
$$\Bc^{\Sigma_1\sqcup \Sigma_2}\cong \Bc^{\Sigma_1}\times \Bc^{\Sigma_2},$$and this is true also for  1- and 2-cobordisms. \end{Lemma}

 Recall, now, the construction of the bicategory,  $\vProfGrphf$, defined in Subsection \ref{sec:prof_conventions}, particularly \S\ref{sum:conv_prof}.

\begin{Definition}\label{def:once_ext_Quinn}
The \emph{once-extended Quinn TQFT}, denoted $$\tFQ{\Bc}\colon \tcob{n} \longrightarrow \vProfGrphf,$$ is defined to be the bifunctor given by:
\begin{itemize}[leftmargin=7mm]
 \item  if $\Sigma$ is a closed smooth $n$-manifold, then $\tFQ{\Bc}^0(\Sigma):= \pi_1(\Bc^\Sigma,\Bc^\Sigma)$; 
 \item if $(i,S,j)\colon \Sigma \to \Sigma'$ is a cobordism, then  $\tFQ{\Bc}^1\big( {(i,S,j)}\colon \Sigma \to \Sigma'\big)\colon \tFQ{\Bc}^0(\Sigma) \bto \tFQ{\Bc}^0(\Sigma')$ is the profunctor,
$$\Hp{\big (\Bc^\Sigma\ra{(i^*,\Bc^S,j^*)} \Bc^{\Sigma'} \big )} \colon  \pi_1(\Bc^\Sigma,\Bc^\Sigma)\bto \pi_1(\Bc^{\Sigma'},\Bc^{\Sigma'}), $$  we are using the notation, $\Hp$, from  Definition \ref{def:fibrantSpanstoprof};
\item[] \hspace{-5.5mm}and
 \item the equivalence class of an extended cobordism,  as in Equation \eqref{diag-cob-forR},
 $$\K\colon \big((i_1,S,i_2)\colon \Sigma_1 \to \Sigma_2\big) \To \big((i_1',S',i_2')\colon \Sigma_1 \to \Sigma_2\big), $$ is sent to the natural transformation of profunctors,
 $$\tFQ{\Bc}^2([\K])\colon  \tFQ{\Bc}^1\big( {(i_1,S,i_2)}\colon \Sigma \to \Sigma'\big) \To\tFQ{\Bc}^1\big( {(i_1',S',i_2')}\colon \Sigma \to \Sigma'\big),$$
derived from $\Bc^{\K}$ in \eqref{diag-dcob-forR}, using  Lemma \ref{def_2HW} in \S \ref{sec_resolvedspans_to_nat_tranfs}, namely,
 $$\tHp{\Bc^\K}\colon  \Hp\big (\Bc^{\Sigma_1} \ra{(i_1^*,\Bc^S,i_2^*)} \Bc^{\Sigma_2} \big )\To \Hp\big( \Bc^{\Sigma_1} \ra{ ({i'^*_1},\Bc^{S'},{i'^*_2)}} \Bc^{\Sigma_2} \big).  $$
\end{itemize}
\end{Definition}
It should perhaps be noted that the name we have used here needs justifying. We have not as yet shown that the structure outline above does give a once-extended TQFT as  that will require a proof that the bifunctor \emph{is} symmetric monoidal.  That will be shown later (see Theorem \ref{eQuinn is an eTQFT} in Subsection \ref{Quinn_is_sym_mon}).

From the constructions in Section \ref{sec:homotopy_underpinning}, combined with the previous lemmas,  it follows that we do indeed have a bifunctor,
$\tFQ{\Bc}\colon \tcob{n} \to \vProfGrphf$.
The fact that  $\tFQ{\Bc}$ preserves the composition of cobordisms is in Proposition \ref{techfib3},  that $\tFQ{\Bc}$ preserves the horizontal composition of extended cobordisms follows from Proposition \ref{lemhorcomp},
and that  $\tFQ{\Bc}$ preserves the vertical composition of extended cobordisms is dealt with by Lemma \ref{Lem:_vert_comp}.
Preservation of vertical identities follows from Lemma \ref{lem:vert_ids}. Finally, preservation of horizontal identities and unitors follows from the discussion in \S\ref{horidsspan}, particularly Lemma
\ref{preserves unitors}.

Also note that,  if $\Sigma$ is a smooth closed manifold, then the groupoid,  $\tFQ{\Bc}^0(\Sigma)= \pi_1({\Bc}^\Sigma,{\Bc}^\Sigma)$,  is homotopy finite. This follows since the function space, $\Bc^\Sigma$, is homotopy finite (Lemma \ref{lem:mappingfromCW}) and so, given a pair of objects, $f, f'\colon \Sigma \to \Bc$, of $\pi_1(\Bc^\Sigma,\Bc^\Sigma)$, the set of arrows from $f$ to $f'$, in  $\pi_1(\Bc^\Sigma,\Bc^\Sigma)$, is finite.

\begin{Remark}\label{rem:tcobp} Let $\tcobp{n}$ be obtained from $\tcob{n}$, by considering the 2-cells to be extended cobordisms, therefore not considering the latter to be up to equivalence. Because the vertical composition of extended cobordisms is not associative,  $\tcobp{n}$ is then not a bicategory. However, we still have compositions, units and unitors, so $\tcobp{n}$ is, similarly to $\mathbf{2span}(HF)$, a 2-truncated cubical set with compositions, see Subsection \ref{sec:comm_summary}.

Moreover, the mapping space construction $\Bc^{-}$, in Lemma \ref{2-span of function spaces}, gives rise to a map of 2-truncated globular sets,
$\Bc^{(-)}\colon \tcobp{n} \to  \mathbf{2span}(HF)$,
which preserves all compositions, units and unitors, up to isomorphism.

The once-extended Quinn TQFT arises from the composite of assignments, below, using the notation of Subsection \ref{sec:comm_summary},
$$ \tcobp{n} \ra{\Bc^{(-)}}\mathbf{2span}(HF) \ra{\boldsymbol{\mathcal{H}}} \vProfGrphf  .$$
\end{Remark}

\subsubsection{Reduced extended cobordisms}\label{sec: reduced cobordism}
We  need some further notation.
Consider an extended cobordism $\K$ as in \eqref{diag-cob-forR}. The pushout diagram on the left-hand-side of \eqref{eq:reduceddef}, below, then induces the commutative diagram on the right, obtained from \eqref{diag-cob-forR} by squashing the vertical cylinders, $\Sigma_1\times I$ and $\Sigma_2\times I$, to $\Sigma_1$ and $\Sigma_2$,
\begin{equation}\label{eq:reduceddef}
\vcenter{ \xymatrix{(\Sigma_1\times I) \sqcup (\Sigma_2\times I)\ar[d]_{\langle i_E,i_W\rangle}\ar[r]^>>>>>{proj}  & \Sigma_1 \sqcup \Sigma_2\ar[d]^{\langle \hat{i_E},\hat{i_W}\rangle }\\
            K\ar@{{}{}{}}[ru]|{\text{push out }} \ar[r] &\hat{K}, }} \qquad \hat{\K}= \vcenter{\xymatrix@C=2.8pc@R=2pc{
   &  S\ar[d]^{\hat{i_N}}   & \\
  \Sigma_1\ar[r]_{\hat{i_E}}\ar[ur]^{i_1}  \ar[rd]_{{i_1'}}  &  \hat{K} & \Sigma_2.\ar[l]^{\hat{i_W}} \ar[ul]_{{i_2}}\ar[ld]^{{i_2'}}\\
    &  S'\ar[u]_{\hat{i_{S}}}    &
 }}
  \end{equation}
To $\hat{\K}$ we call the \emph{reduced cobordism} of $\K$, and the colimit $S\sqcup_{\Sigma_1\sqcup \Sigma_2} S'$ is called the \emph{reduced frame}, $\hat{\fr}(\K)$, of $\K$. The \emph{reduced filler}, of $\K$, is the universally defined map $\hat{f}_{\K}\colon \hat{\fr}({\K})\to \hat{K}$.

\subsubsection{A more explicit description of $\tFQ{\Bc}\colon \tcob{n} \to \vProfGrphf$}

In the discussion below, we use,  from items \eqref{Def_PC} and \eqref{Def_hatpi}  on  page \pageref{Def_hatpi}, that, if $X$ is a CGHW space, and $x \in X$, then the path-component in $X$ to which $x$ belongs, topologised with the induced CGWH topology, is denoted $\PC_x(X)$, and that we denote the set of all such path-components by $\hpi_0(X).$
We also note the definition of the homotopy content, $\chi^\pi(B)$, of a HF  space $B$ in Definition \ref{def:hcont}.

Here is an explicit explanation of $\tFQ{\Bc}\colon \tcob{n} \to \vProfGrphf$.

On objects, $\tFQ{\Bc}$ sends a closed $n$-manifold, $\Sigma$, to the fundamental groupoid, $\pi_1(\Bc^\Sigma,\Bc^\Sigma)$, of the space $B^\Sigma$ of functions $\Sigma \to B$.

Given a cobordism, $(i_1,S,i_2)\colon \Sigma_1 \to \Sigma_2$,  if  $f_1\colon \Sigma_1 \to \Bc$ and $f_2\colon \Sigma_2 \to \Bc$ are continuous functions,  hence objects of $\pi_1(\Bc^{\Sigma_1},\Bc^{\Sigma_1})$ and of  $\pi_1(\Bc^{\Sigma_2},\Bc^{\Sigma_2})$, we have:
 $$\tFQ{\Bc}^1(i_1,S,i_2)(f_1,f_2)= \Lin\big(\hpi_0(\{f_1|\Bc^{S}|f_2\})\big).$$
Here $\Lin\colon \Sets\to \Vect$ denotes the free vector space functor, and
$ \{f_1|\Bc^{S}|f_2\}$, in full $ \{f_1|\Bc^{(i_1^*,S,i_2^*)}|f_2\}$, is the space of maps,
$H\colon S \to \Bc$, making the diagram commute,  $$\xymatrix@R=2pt{ &\Bc\\
 \Sigma_1\ar[ur]^{f_1} \ar[dr]_{i_1} && \Sigma_2.\ar[ul]_{f_2} \ar[dl]^{i_2} \\
& S\ar[uu]^H } $$

Given paths, $\gamma_1 \colon f_1' \to f_1$ in $\Bc^{\Sigma_1}$ and $\gamma_2\colon f_2 \to f_2'$ in $\Bc^{\Sigma_2}$,  the linear map,
 \begin{multline*}
 \tFQ{\Bc}^1(i_1,S,i_2)(f_1'\ra{[\gamma_1]}f_1,f_2\ra{[\gamma_2]}f_2')\colon\\ \Lin\big(\hpi_0(\{f_1|\Bc^{S}|f_2\})\big) \to \Lin\big(\hpi_0(\{f_1'|\Bc^{S}|f_2'\})\big),\end{multline*}
 is defined from the functor, 
 $${\Fc^{(\Bc^{S})}\colon \pi_1(\Bc^{\Sigma_1} \times \Bc^{\Sigma_2},\Bc^{\Sigma_1} \times \Bc^{\Sigma_2} ) \to \CGWH/\simeq,}$$ obtained from the path-space fibration,  $$\langle i_1,i_2 \rangle^*\colon \Bc^{S} \to \Bc^{\Sigma_1} \times \Bc^{\Sigma_2}\cong \Bc^{\Sigma_1 \sqcup \Sigma_2},$$ in the usual way (see \cite[Chapter 7]{May}, as reviewed in Lemma \ref{main_techfibfunctor}), and then by applying $\hpi_0\colon \CGWH/\simeq \to \Sets$;  finally linearising by applying $\Lin\colon \Sets \to \Vect$.
 (Note that the path $\gamma_1 \colon f_1' \to f_1$ must be inverted before applying $\Fc^{(\Bc^{S})}$.) An explicit description can be obtained from the comments just after Definition \ref{def:fibrantSpanstoprof}.

For an extended cobordism,  as in \eqref{diag-cob-forR}, $$\K\colon \big((i_1,S,i_2)\colon \Sigma_1 \to \Sigma_2\big)\To \big((i_1',S',i_2')\colon \Sigma_1 \to \Sigma_2\big),$$
then the natural transformation, of profunctors, $$\tFQ{\Bc}^2([\K])\colon \tFQ{\Bc}^1( (i_1,S,i_2)\colon \Sigma_1 \to \Sigma_2)\implies \tFQ{\Bc}^1( (i'_1,S',i'_2)\colon \Sigma_1 \to \Sigma_2), $$
is such that, if   $f_1\colon \Sigma_1 \to \Bc$ and $f_2\colon \Sigma_2 \to \Bc$, given $H\in
\{f_1|\Bc^{S}|f_2\}$ and $H'\in  \{f_1|\Bc^{S'}|f_2\}$, then we have the following formula for the corresponding matrix elements, of $\big(\tFQ{\Bc}^2([\K])\big)_{(f_1,f_2)}\colon\tFQ{\Bc}^1(i_1,S,i_2)(f_1,f_2)\to \tFQ{\Bc}^1(i'_1,S',i'_2)(f_1,f_2)$,
\begin{multline}\label{morphism from ext cobord}
 \big \langle \PC_{H}\big({\{f_1|\Bc^{S}|f_2\}}\big) \mid \big(\tFQ{\Bc}^2([\K])\big)_{(f_1,f_2)}\mid  \PC_{H'}\big ({\{f_1|\Bc^{S'}|f_2\}}\big)\big \rangle\\[0.5ex]
 =\chi^\pi\left(\left \{ T\colon K \to \Bc\left |
 \begin{smallmatrix}  T\circ i_N=H, \quad T\circ i_S=H',\smallskip
 \\  \forall s \in \Sigma_1, \,\,\forall t \in [0,1]:\,\, T(i_E(s,t))=f_1(s), \smallskip\\ \forall  s'\in \Sigma_2, \,\, \forall t \in [0,1]:\,\,  T(i_W(s',t))=f_2(s'). \end{smallmatrix}\right.\right\}\right)
 \\ \  \qquad  \qquad\,\,\,\, \chi^\pi \big( \PC_{H'}\big({\{f_1|\Bc^{S'}|f_2\}}\big)\big) .
\end{multline}
 (It follows from the construction in Section \ref{sec:homotopy_underpinning}, in particular Lemma \ref{lem;2elsHF}, that we are indeed considering homotopy contents only of homotopy finite spaces.)

This simplifies further. We use the notation in \S \ref{sec: reduced cobordism}.
\begin{Theorem}\label{thm:red-form-2matrixeles}Given $f_1\colon \Sigma_1 \to \Bc$ and $f_2\colon \Sigma_2 \to \Bc$, the  matrix elements, of $\big(\tFQ{\Bc}^2([\K])\big)_{(f_1,f_2)}\colon\tFQ{\Bc}^1(i_1,S,i_2)(f_1,f_2)\to \tFQ{\Bc}^1(i'_1,S',i_2')(f_1,f_2)$, are equal to
\begin{multline}\label{morphism from ext cobord-simp}
 \big \langle \PC_{H}\big({\{f_1|\Bc^{S}|f_2\}}\big) \mid \big(\tFQ{\Bc}^2([\K])\big)_{(f_1,f_2)}\mid  \PC_{H'}\big ({\{f_1|\Bc^{\smash{S'}}|f_2\}}\big)\big \rangle\\[0.5ex]
 =\chi^\pi\left(\left \{ T\colon \hat{K} \to \Bc\left |
 \begin{smallmatrix}  T\circ \hat{i_N}=H, \quad T\circ \hat{i_S}=H', \smallskip
 \\ T\circ\hat{i_E}=f_1, \quad  T\circ\hat{i_W}=f_2.\end{smallmatrix}\right.\right\}\right)
  \chi^\pi \big( \PC_{H'}\big({\{f_1|\Bc^{S'}|f_2\}}\big)\big) .
\end{multline}
\end{Theorem}
\begin{proof}
 This follows from Equation  \eqref{morphism from ext cobord} and Lemma \ref{Lem:squash-functions}, where $A=(\Sigma_1\times I) \sqcup (\Sigma_2\times I)$, $A'=\Sigma_1\sqcup \Sigma_1$ and $h\colon A \to A'$ is the map that squashes $I$ to a point.
\end{proof}
\begin{Lemma}\label{Lem:squash-functions}
 Let $X,A'$ and $B$ be CGWH spaces. Let $A$ be a subspace of $X$, and let $i_A\colon A \to X$ be the inclusion. Let $P_A\colon \TOP(X,B)\to \TOP(A,B)$ be the restriction function. Let $h\colon A \to A'$ be a continuous function, and form the pushout, on the left-hand-side of the equation below. Let $P_{A'}\colon  \TOP(X',B)\to \TOP(A',B)$  denote the restriction function.
Finally,  let $f\colon A' \to B$ be any function.

There exists a homeomorphism between fibres: $P_{A'}^{-1}(f)\cong P_{A}^{-1}(f\circ h)$, so in between the two pullbacks on the right-hand side of the equation below:
 $$\vcenter{ \xymatrix{ A \ar[d]_{i_A}\ar[r]^h  & A'\ar[d]^{{i'_A} }\\
            X\ar@{{}{}{}}[ru]|{\text{pushout }} \ar[r]_{\hat{h}}  &X',}}\quad  \vcenter{\xymatrix{P_A^{-1}(f\circ h)\ar[r]^{inc}\ar[d] & \TOP(X,B)\ar[d]^{P_A}\\
                                                                         \{f\circ h\}\ar[r]_{inc} \ar@{{}{}{}}[ru]|{\text{pullback }} &\TOP(A,B),        }} \quad \vcenter{\xymatrix{P_{A'}^{-1}(f)\ar[r]^<<<<<<{inc}\ar[d] & \TOP(X',B)\ar[d]^{P_{A'}}\\
                                                                         \{f\}\ar@{{}{}{}}[ru]|{\text{pullback}}\ar[r]_<<<<<<<<{inc} &\TOP(A',B).        }} $$
\end{Lemma}
\begin{proof}Consider diagram below, and apply the pullback pasting lemma:
$$ \vcenter{\xymatrix{P_{A'}^{-1}(f)\ar[r]^<<<<<{inc}\ar[d] & \TOP(X',B)\ar[d]^{P_{A'}} \ar[r]^{\hat{h}^\ast}  & \TOP(X,B)\ar[d]^{P_A} \\
                                                                         \{f\}\ar@{{}{}{}}[ru]|{\text{pullback}}\ar[r]_<<<<<<<<{inc} &\TOP(A',B) \ar@{{}{}{}}[ru]|{\text{pullback} } \ar[r]_{h^\ast}  &   \TOP(A,B).  }} $$
(Note $h^*(f)=  f\circ h$.) The square on the right is a pullback since $\CGWH$ is cartesian closed.
\end{proof}

Note that, unless  $\Bc$ is a finite set with the discrete topology, then for $\Sigma$,  a closed smooth $n$-manifold, the groupoid, $\tFQ{\Bc}^0(\Sigma)=\pi_1(\Bc^\Sigma,\Bc^\Sigma)$, is uncountable, since it has an uncountable set of objects. However  $\tFQ{\Bc}^0(\Sigma)$    \emph{is} homotopy finite, so we can extract an equivalent finite subgroupoid from it. Namely, if we choose a finite subset, $\fd_\Sigma\subset \Bc^\Sigma$, containing at least one element for each path-component, then $\pi_1(\Bc^\Sigma,\fd_\Sigma)$ will be a finite groupoid, equivalent to $\tFQ{\Bc}^0(\Sigma)=\pi_1(\Bc^\Sigma,\Bc^\Sigma)$.

This latter fact will be used in the next subsection to construct a finitary version of the once-extended Quinn TQFT.

\subsection{A finitary version of the once-extended Quinn TQFT}\label{sec:fin_ext}

For technical and historical reasons in the applications of the above theory, it is often useful to replace groupoids that have possibly infinitely many objects, by  more finitary, but equivalent, ones. There are several useful ways of doing this, for instance, using triangulations of the manifolds, as we will do in Section \ref{sec:TQFTS_xcomp}. We will, in this section, reduce the size of the groupoids by a different means as follows.

As always, let $\Bc$ be a homotopy finite space, and  $n$ be a non-negative integer.

\begin{Definition}[$\Bc$-decorated manifold]\label{def:DecoratedMan}
 A $\Bc$-decorated $n$-manifold, denoted  $\Sd=(\Sigma,\fd_{\Sigma})$, is given by a closed smooth $n$-manifold, $\Sigma$, called the \emph{underlying manifold} of $\Sd$, together with a finite subset, $\fd_\Sigma$, of $\Bc^\Sigma$, containing at least one function, $f\colon \Sigma \to \Bc$, for each path component of the space, $\Bc^\Sigma$.
\end{Definition}
Let $\Sigma$ be a closed smooth manifold. We recall, \cite{MunkresDiff}, that  $\Sigma$ has a finite CW-decomposition. Since $\Bc$ is homotopy finite,  $\Bc^\Sigma$ is homotopy finite (Lemma \ref{lem:mappingfromCW}), and hence it has a finite number of path-components. In particular, we can see that  all closed (smooth) manifolds possess a $\Bc$-decoration.

\begin{Definition}
\label{FinQ} We define a  bicategory, \smash{$\tdcob{n}{\Bc}$}, as follows.
\begin{itemize}[leftmargin=10mm]
\item The objects are $\Bc$-decorated $n$-manifolds.

\item Given $\Bc$-decorated $n$-manifolds, $\Sd=(\Sigma,\fd_\Sigma)$ and $ \Sdp = (\Sigma',\fd_{\Sigma'})$, the 1-morphisms, $(i,S,j)\colon \Sd \to \Sdp, $ are given by $(n+1)$-cobordisms, $(i,S,j)\colon \Sigma \to \Sigma'$, with no additional structure on $S$.
(The  $(n+1)$-cobordism, $(i,S,j)\colon \Sigma \to \Sigma'$, associated to a 1-morphism $(i,S,j)\colon \Sd\to \Sdp$, will be  called the \emph{underlying $(n+1)$-cobordism} of that 1-morphism.) 

\item Given 1-morphisms, $(i,S,j),(i',S',j')\colon \Sd \to \Sdp $, the 2-morphisms, 
$${[\K]\colon \big( (i,S,j)\colon \Sd \to \Sdp \big) \implies\big( (i',S',j') \colon \Sd \to \Sdp\big)} ,$$ are given by equivalence classes of extended cobordisms,
$${\K\colon \big( (i,S,j)\colon \Sigma \to \Sigma' \big) \implies\big( (i',S',j') \colon \Sigma \to \Sigma'\big),} $$ with the equivalence relation as in $\tcob{n}$.
(As for 1-morphisms, we define the underlying 2-morphism of $\tcob{n}$, associated to such a 2-morphism of $\tdcob{n}{\Bc}$,  by forgetting the $\Bc$-decoration on the $n$-manifolds $\Sigma$ and $\Sigma'$.)
\end{itemize}
The rest of the bicategory structure in $\tdcob{n}{\Bc}$ is induced, in the obvious way, by that of the undecorated case,
 i.e. by composing the underlying cobordisms or underlying 2-morphisms of the 1- and 2-morphisms in $\tdcob{n}{\Bc}$. For instance, the composition below,
$$(\Sigma, \fd_\Sigma) \ra{(i,S,j) }( \Sigma', \fd_{\Sigma'}) \ra{(i',S',j') }  ( \Sigma'', \fd_{\Sigma''})$$
 simply gives
$$(\Sigma, \fd_\Sigma) \ra{(\overline{i},S\bullet S',\overline{j'})}( \Sigma'', \fd_{\Sigma''}),$$
where $(\overline{i},S\bullet S',\overline{j'})$ gives the composite of cobordisms as in item \eqref{compcob} on page \pageref{compcob}.
\end{Definition}

For convenience, we recall that the bicategory, $\vProfGrpfin$, defined in \S\ref{sum:conv_prof},
is the sub-bicategory of $\vProfGrp$, whose objects are the finite groupoids.

Given  a $\Bc$-decorated manifold, $ \Sd=(\Sigma,\fd_{\Sigma})$, the pair $(\Bc^\Sigma,\fd_{\Sigma})$ is, by definition,  0-connected. From Lemmas \ref{techfib3} and \ref{lemhorcomp}, it follows that the bifunctor, $$\tFQ{\Bc}\colon \tcob{n} \to \vProfGrphf,$$  induces a bifunctor, $$\tFQd{\Bc}\colon \tdcob{n}{\Bc} \to \vProfGrphf,$$
by restricting from $\pi_1(\Bc^\Sigma,\Bc^\Sigma)$ to $\pi_1(\Bc^\Sigma,\fd_\Sigma)$, and leaving the rest of the structure unaltered.
Since we assume that $\Sigma$ is a closed smooth manifold, as above, it follows that $\Bc^\Sigma$ is homotopy finite, and, thus, that the groupoid $\pi_1(\Bc^\Sigma,\fd_\Sigma)$ is finite.

This leads to the following definition in which we are again using the notation of Definition \ref{def:2matrixels} and from from notational comment  \ref{not_barH} in
 Subsection \ref{proffibspan}.

\begin{Definition}\label{def:finitary_ext_TQFT}
The \emph{finitary once-extended Quinn TQFT}, $$\tFQd{\Bc}\colon \tdcob{n}{\Bc} \to \vProfGrpfin,$$ is the bifunctor defined as follows.
 \begin{itemize}[leftmargin=10mm]
 \item If $\Sd=(\Sigma,\fd_\Sigma)$ is a $\Bc$-decorated $n$-manifold, then $\tFQd{\Bc}(\Sd):= \pi_1(\Bc^\Sigma,\fd_\Sigma).$   \item  If $(i,S,j)\colon \Sd=(\Sigma,\fd_\Sigma) \to \Sdp=(\Sigma',\fd_{\Sigma'})$ is a 1-morphism, then
 \begin{align*}
 \tFQd{\Bc}&\big( (\Sigma,\fd_{\Sigma}) \ra{(i,S,j)} (\Sigma',\fd_{\Sigma'}) \big)\\ &:=
 \Hpb_{(\fd_{\Sigma},\fd_{\Sigma'})}\big( \Bc^\Sigma\ra{(i^*,\Bc^S,j^*)}\Bc^{\Sigma'} \big)\colon  \pi_1(\Bc^\Sigma,\fd_{\Sigma})\bto \pi_1(\Bc^{\Sigma'},\fd_{{\Sigma'}}).
 \end{align*}
 Concretely the functor, $$\tFQd{\Bc}\big( (\Sigma,\fd_{\Sigma}) \ra{(i,S,j)} (\Sigma',\fd_{\Sigma'}) \big)\colon \pi_1(\Bc^\Sigma,\fd_{{\Sigma}})^{\op}\times \pi_1(\Bc^{\Sigma'},\fd_{{\Sigma'}})\to \Vect,$$ is  the restriction to $\pi_1(\Bc^\Sigma,\fd_{{\Sigma}})^{\op}\times \pi_1(\Bc^{\Sigma'},\fd_{{\Sigma'}})$, of the functor,
 $$\tFQ{\Bc}^1{\big ( \Sigma \ra{(i,S,j)} \Sigma'\big )}\colon \pi_1(\Bc^\Sigma,\Bc^\Sigma)^{\op}\times \pi_1(\Bc^{\Sigma'},\Bc^{\Sigma'})\to \Vect.$$
 \item A 2-morphism, of  $\tdcob{n}{\Bc}$,
 $$[\K]\colon\big( (i_1,S,i_2)\colon  (\Sigma_1,\fd_{{\Sigma_1}}) \to (\Sigma_2,\fd_{{\Sigma_2}})\big) \implies \big ((i_1',S',i_2')\colon (\Sigma_1,\fd_{{\Sigma_1}}) \to (\Sigma_2,\fd_{{\Sigma_2}})\big),$$
 is sent to the natural transformation,
 $$
 \tHpc{\Bc^\K}{\fd_{{\Sigma_1}}}{\fd_{{\Sigma_2}}}\colon  \Hpb_{(\fd_{{\Sigma_1}},\fd_{\Sigma_2})}(i_1^*,\Bc^S,i_2^*)\implies \Hpb_{(\fd_{{\Sigma_1}},\fd_{\Sigma_2})}({i'^*_1},\Bc^{S'},{i'^*_2}),$$
of functors, $\pi_1(\Bc^{\Sigma_1},\fd_{\Sigma_1})^{\op}\times  \pi_1(\Bc^{\Sigma_2},\fd_{{\Sigma_2}})\to \Vect$, as in Definition \ref{def:2matrixels}.
 \end{itemize}
\end{Definition}
The proof that we have indeed defined a bifunctor  $\tFQd{\Bc}$, follows as for the earlier case of $\tFQ{\Bc}$. As before, the crucial fact that $\tFQd{\Bc}$ preserves the composition of the 1-morphisms and the horizontal compositions of the 2-morphisms of $\tdcob{n}{\Bc}$, follows from Proposition \ref{techfib3} and Proposition \ref{lemhorcomp}. It is for these cases that we need to impose that each pair, $(\Bc^\Sigma,\fd_\Sigma)$, is $0$-connected.
The proof that  $\tFQd{\Bc}$ is symmetric monoidal will be done in Theorem \ref{Decorated case is eTQFT}.

\subsubsection{The dependence of $ \tFQd{\Bc}\colon \tdcob{n}{\Bc} \to \vProfGrpfin$  on decorations}\label{Sec-Decorations:ind}
Let,  as before, $\Sigma$ be a closed smooth $n$-manifold.
 The finitary once-extended Quinn TQFT, $\tFQd{\Bc}$ does not give a value to $\Sigma$ itself, except when $\Sigma$ is given a $\Bc$-decoration, $\fd_{{\Sigma}}$. However, given two $\Bc$-decorations, $\fd_{{\Sigma}}$ and $\fd'_{{\Sigma}}$, we have a  profunctor, $\Psi({\fd_{{\Sigma}},\fd'_{{\Sigma}}})\colon \tFQd{\Bc}(\Sigma, \fd_{{\Sigma}})\bto \tFQd{\Bc}(\Sigma, \fd'_{\Sigma})$, defined by
 $$\tFQd{\Bc}\big ( (\Sigma,\fd_\Sigma) \ra{ (\iota_0^\Sigma,\Sigma\times I,\iota_1^\Sigma)} (\Sigma,\fd'_\Sigma) \big)\colon \tFQd{\Bc}(\Sigma, \fd_\Sigma)\bto \tFQd{\Bc}(\Sigma, \fd'_{\Sigma}).$$

 By construction, the profunctors $\Psi({\fd_{{\Sigma}},\fd'_{{\Sigma}}})$ are compatible with the bifunctor,  $\tFQd{\Bc}\colon \tdcob{n}{\Bc} \to \vProfGrpfin$,  in that the following diagram commutes, given a cobordism $(i,S,j)\colon \Sigma_1 \to \Sigma_2$, and pairs of decorations in $\Sigma_1$ and $\Sigma_2$,   up to a natural isomorphism of profunctors,
 $$ \xymatrix@C=140pt{\tFQd{\Bc}(\Sigma_1, \fd_{\Sigma_1})\ar[r]^{\tFQd{\Bc}\big( (i,S,j)\colon (\Sigma_1, \fd_{\Sigma_1}) \to (\Sigma_2, \fd_{\Sigma_2}) \big) }\ar[d]_{\Psi({\fd_{{\Sigma_1}},\fd'_{{\Sigma_1}}})}
 & \tFQd{\Bc}(\Sigma_2, \fd_{\Sigma_2})\ar[d]^{\Psi({\fd_{{\Sigma_2}},\fd'_{{\Sigma_2}}})}\\
 \tFQd{\Bc}(\Sigma_1, \fd'_{\Sigma_1})
 \ar[r]_{\tFQd{\Bc}\big( (i,S,j)\colon (\Sigma_1, \fd'_{\Sigma_1}) \to (\Sigma_2, \fd'_{\Sigma_2}) \big) }&\tFQd{\Bc}(\Sigma_2, \fd'_{\Sigma_2}).
 }
 $$
 This is because the two paths in the diagram above yield diffeomorphic cobordisms.

The profunctors  $\Psi({\fd_{{\Sigma}},\fd'_{{\Sigma}}})$ associated to pairs of $\Bc$-decorations of $\Sigma$ are furthermore `functorial' in the following sense.
Consider decorations 
$\fd_\Sigma$, $\fd'_\Sigma$ and $\fd''_\Sigma$ of $\Sigma$. We have 1-morphisms in $\tdcob{n}{\Bc}$,
$$  (\Sigma,\fd_\Sigma) \ra{ (\iota_0^\Sigma,\Sigma\times I,\iota_1^\Sigma)}  (\Sigma,\fd'_\Sigma), \qquad (\Sigma,\fd'_\Sigma) \ra{ (\iota_0^\Sigma,\Sigma\times I,\iota_1^\Sigma)}  (\Sigma,\fd''_\Sigma) , $$ 
and the following, which is homeomorphic to their composite,
$$  (\Sigma,\fd_\Sigma) \ra{ (\iota_0^\Sigma,\Sigma\times I,\iota_1^\Sigma)}  (\Sigma,\fd''_\Sigma) . $$
Applying Lemma \ref{techfib3},
we have a natural isomorphism of profunctors,
$$\Psi({\fd_{{\Sigma}},\fd'_{{\Sigma}}}) \bullet \Psi({\fd'_{{\Sigma}},\fd''_{{\Sigma}}})\Longrightarrow \Psi({\fd_{{\Sigma}},\fd''_{{\Sigma}}}).$$
It can be proved that these natural isomorphisms satisfy appropriate relations when we consider four different $\Bc$-decorations of $\Sigma.$  In particular,  the profunctor associated to a change of $\Bc$-decoration is always invertible up to 2-isomorphism, i.e. is an adjoint equivalence.
\subsection{The Morita valued extended Quinn TQFT}\label{sec:Mor_ext}

We continue working with our chosen subfield,  $\kappa$, of $\C$.
The finitary theory, as  given in the previous section, takes values in a bicategory of $\Vect$-valued profunctors between finite $\kappa$-linear categories. To be more easily able to use more usual representation theoretic methods and ideas, it can be convenient to replace this bicategory by one that is better know within the representation theoretic setting, namely that of finite dimensional algebras (with 1), bimodules and morphisms between them.

Here we will first  review the construction of an algebra from a linear category, as given by Mitchell, \cite{Mitchell72}, \S 7, and then look at it in detail for $Lin(\Gamma)$, the linear category associated to a (finite) groupoid, $\Gamma$, obtained by applying the free vector space functor to the morphism sets of $\Gamma$.  In this case, the resulting algebra is the well known \emph{groupoid algebra}, \cite{Willerton}.  We look, in some detail, at the relationship between bimodules over a category algebra and profunctors. Some  of this is \textit{folklore}, and is quite difficult to find explicitly in the literature, yet it seems important for the understanding of the relationship between the $\Prof$-valued and the $\Mor$-valued extended TQFTs.

We will  define the Morita bicategory, $\Mor$, or, more exactly, $\Mor_\kappa$, (also sometimes denoted $\mathbf{Alg_2}$), of algebras, bimodules and the bimodule morphisms, these latter being often known as \emph{intertwiners} in a representation theoretic context, and will examine its relation to   $\vProf$.  We will then see how to define a \emph{Morita valued extended Quinn TQFT.
}

\subsubsection{The algebra of a small linear category}
 \label{Groupoid algebra}
 In Section 7, (page 33), of the classic paper, \cite{Mitchell72}, by Mitchell, it was shown how to associate a ring to a small additive category.  That construction easily extends to a $\kappa$-linear version for any $\kappa$-linear category.
 
 Let $\mathcal{C}$ be a (small) $\kappa$-linear category, having $\mathcal{C}_0$ as its set of objects. We set $[\mathcal{C}]$ to be the set of $\mathcal{C}_0\times \mathcal{C}_0$ matrices, $c$, where,  for $p, q\in\mathcal{C}_0$,  the $(p,q)$-entry, denoted $c_{p,q}$, is an element of the vector space, $\mathcal{C}(p,q)$, of arrows from $p$ to $q$,  and each row and column has only a finite number of non-zero entries.  Using the addition in each $\mathcal{C}(p,q)$, together with the composition from $\mathcal{C}$, we can give $[\mathcal{C}]$ the structure of a $\kappa$-algebra, which will not usually be commutative, nor, in general, unital.
 
 We thus have that, as a vector space, 
 $[\mathcal{C}]=\bigoplus_{p,q\in \mathcal{C}_0}\mathcal{C}(p,q)$,
 and the multiplication is given by 
 $$g\cdot f=\begin{cases}g\circ f& \textrm{ if } domain(g)=codomain(f)\\0&\textrm{ otherwise.}
 \end{cases}$$
 
 Although, in general, $[\mathcal{C}]$ will not have a multiplicative identity, each object $p\in \mathcal{C}_0$ gives an idempotent matrix, $\mathbf{1}_p$, namely the matrix having the identity morphism on $p$ in the $(p,p)$-position and zeroes elsewhere. 
 
 If $\mathcal{C}_0$ is finite, then $\sum_{p\in \mathcal{C}_0}\mathbf{1}_p$ is, however, a multiplicative identity for  $[\mathcal{C}]$.  
 The algebra, $[\mathcal{C}]$, is an example of a generalised matrix algebra. This algebra is called the \emph{category algebra} of $\mathcal{C}$.

Any element, $\mathbf{c}=(c_{p,q})$, in $[\mathcal{C}]$ can be written as a sum of matrices of form $\mathbf{c}_{p,q}$, where the matrix $\mathbf{c}_{p,q}$  is to be zero in all positions except the $(p,q)$ position, where it is,   no surprise, $c_{p,q}$. This sum is finite.  The element, $c_{p,q}$, clearly has domain equal to $p$ and codomain  equal to $q$, so, in a completely classical way,  the product 
$\mathbf{c}\cdot \mathbf{1}_p=\sum_{r,s}\mathbf{c}_{r,s}\cdot \mathbf{1}_p= \sum_s \mathbf{c}_{p,s}$, whilst $\mathbf{1}_q\cdot \mathbf{c}= \sum_r  \mathbf{c}_{r,q}$, so $\mathbf{1}_q\cdot \mathbf{c}\cdot \mathbf{1}_p= \mathbf{c}_{p,q}$, and, in particular, we have the useful equation: $\mathbf{1}_q\cdot \mathbf{c}_{p,q}\cdot \mathbf{1}_p= \mathbf{c}_{p,q}$.

  \begin{Example}\label{Matrix algebra}If $\mathbf{P}$ is a (finite) pre-ordered set, and $\mathcal{C} $ is the linearisation of the corresponding small category, then $[\mathcal{C}]$ is the incidence algebra of the poset. As specific examples, if $\mathbf{P}=\{1<\ldots <n\}$, then $[\mathcal{C}]$ is the algebra of $n\times n$ upper triangular matrices over $\kappa$.  If we replace the given preorder by the discrete preorder, so $p\leq q$ here means $p=q$, then the corresponding $[\mathcal{C}]$ is the algebra of diagonal matrices. If, on the other hand, we replace the preorder by the codiscrete preorder (in which $p\leq q$ for  every pair of elements $(p,q)$), then $[\mathcal{C}]= M_n(\kappa)$, the full algebra of $n\times n$ matrices over  $\kappa$.
 \end{Example}

The main case of category algebras $[\mathcal{C}]$ that we consider here arise in the case when the $\kappa$-linear category $\mathcal{C}$ is the  $\kappa$-linearisation  of a (usually finite) groupoid, $\Gamma$. This has the same objects as $\Gamma$, and the vector space of morphisms $x \to y$ is the free vector over $\hom_\Gamma(x,y)$. If $\Gamma$ is finite, more generally if $\Gamma$ has a  finite set of objects, then the resulting `\emph{groupoid algebra}' will be unital.  It has a well known description in terms of the arrows of $\Gamma$, which we include in case it is found easier to understand, as it is written in a slightly less abstract way; see also \cite{Bohm:Handbook:2009,bohm-et-al:weak-bialgebras:2014,nikshych-vainermann:finiteqntmgpds:2002,Willerton} and \cite{loopy} for various versions and the development of further theory.
 
 \begin{Definition}[Groupoid algebra]Let $\Gamma$ be a finite groupoid. The \emph{groupoid algebra}, $\tLin(\Gamma)$,\label{grpd-alg} has as its underlying vector space,  $\Lin(\Gamma_1)$, the free vector space over the set of morphisms of $\Gamma$.   The product in $\tLin(\Gamma)$ is given on generators by
 $$(x \ra{g} y) (x'\ra{g'} y):=\delta(y,x')(x \ra{g g'} y),$$ where  $\delta(y,x')$ is 1, if the two objects, $y$ and $x'$, are equal, and is zero otherwise.
 The multiplicative unit of $\tLin(\Gamma)$ is given by $\sum_{x \in \Gamma_0}( \id_x\colon x\to x)$.
 \end{Definition}
 \begin{Example}(The Quantum Double as a groupoid algebra) \label{ex.Qdoublefingroup}
  Let $G$ be a finite group. Consider the action of $G$ on itself by conjugation,  so we can form the action groupoid $G\sslash G$ of this action\footnote{We will give a more  general form of action groupoid later; see page \pageref{def:act_groupoid}.}.  This has the elements of $G$ as its objects and the arrows have form $(g,a)\colon g\to aga^{-1}$, where $g,a \in G$.  The groupoid algebra of $G\sslash G$ is given in detail in, for instance, \cite{Willerton} and \cite[ \S1.10]{loopy}.  The product on the basis elements, as above, is given by the formula below, where $g,g',a,a'\in G$,
 $$(g,a)(g',a')=\delta(aga^{-1},g')(g,aa').$$

 As noted in \cite{Willerton}, if we define a comultiplication,
 $$\Delta(x,g)=\sum_{yz=x}(y,g)\otimes (z,g),$$
and a counit, $\epsilon(x,g) = \delta(x,1_G),$ then, for suitable  definitions of an antipode and an $R$-matrix, the resulting object is a quasi-triangular Hopf algebra.  It is clear, from standard descriptions of the `double construction', that this  is $D(G)$ the Drinfel'd double or quantum double of the Hopf algebra, $\kappa[G]$.
   \end{Example}

\subsubsection{Non-functoriality of $[-]$ and Morita equivalences}\label{non-functoriality}

This subsection examines more properties of this situation, but they will not be immediately needed for the main theme of this paper, so can be left aside on first reading, moving on to section \ref{sec:mor_defined}, where the bicategory, $\Mor$ is discussed. We would recommend that they be at least skimmed in a later reading as they  provide further insights on the algebraic mechanisms involved later on.

A functor, $F\colon \mathcal{C}\to \mathcal{D}$, does induce a \emph{linear} map, $[F]\colon[\mathcal{C}]\to [\mathcal{D}]$, but,
in general, this map will not preserve multiplication, so $[-]$  is not a functor, from the category of $\kappa$-linear categories to the category of $\kappa$-algebras. The more-or-less minimal example for
this  is to take $\mathcal{C}$ to be the ($\kappa$-linearisation of the) discrete category on the set having just two elements, say $a$ and $b$, and $\mathcal{D}$ to be the corresponding construction on a singleton set, $\{c\}$.

We will examine this slightly odd situation in a bit more detail shortly, as it is perfectly manageable given the approach that we are using.  In any case,  the following  result of Mitchell, \cite{Mitchell72}, Theorem 7.1, makes one realise that there is a lot of power in the category algebra construction.  For the statement, we think of the category algebra as a linear category having just a single object. 
\begin{Theorem}[Mitchell]\label{Mitchell-Morita}
Suppose $\mathcal{C}$ is a linear category having only finitely many objects, then $\mathcal{C}$ and $[\mathcal{C}]$ are Morita equivalent categories. Explicitly, let $\mathcal{C}\!-\!Mod$ be $Func_\kappa(\mathcal{C},\Vect)$, the category of  $\kappa$-linear functors from $\mathcal{C}$ to $\Vect$. Then $\mathcal{C}\!-\!Mod$ and $[\mathcal{C}]\!-\!Mod$, the category of representations of $[\mathcal{C}]$,  are equivalent categories.
\end{Theorem}

We will sketch out  a proof of this as it contains some ideas that help one to understand what is happening here, and hence why the `linearisation' versus `categorification' process is so useful. (A full proof is given in \cite{Mitchell72} on page 34.  A discussion of the ideas can be found online in the \emph{$n$-Category Caf\'{e}}, (May 14, 2014), in a post, \emph{Categories vs. Algebras}, by Tom Leinster; see \cite{ncafe}.)

\begin{proof}(Sketch.)
First we construct a functor, $[-]$, from $\mathcal{C}\!-\!Mod$ to $[\mathcal{C}]\!-\!Mod$, so suppose that $M\colon \mathcal{C}\to \Vect$ is a $\kappa$-linear functor.  We set $$[M]=\oplus_{q\in \mathcal{C}_0}M(q),$$ the direct sum of all the image vector spaces of the functor, $M$. This is finite dimensional if each $M(q)$ is, as $\mathcal{C}_0$ is finite. It has a $[\mathcal{C}]$-module structure in a natural and fairly obvious way. If $c= (c_{p,q})$ is an element of $[\mathcal{C}]$, then $c_{p,q}:p\to q$ is in $\mathcal{C}$, so $M(c_{p,q}):M(p)\to M(q)$ is a linear map. Now, if $m=(m_p)$ is an element of $[M]$, we define $c\cdot m=\sum_{p,q}M(c_{p,q})(m_p)$.  The linearity and functoriality of $M$ ensures that this does give a $[\mathcal{C}]$-module structure to $[M]$. This is easily seen to define a functor, $[-]$, as claimed. This forms part of the claimed equivalence.

The other direction, starting from a $[\mathcal{C}]$-module  and ending up with a $\kappa$-linear functor from $\mathcal{C}$ to $\Vect$, is not so obvious, although it is, in fact,  a generalisation of a well known process from elementary linear algebra. 

 Suppose $N$ is a $[\mathcal{C}]$-module and that $p$ is an object of $\mathcal{C}$.  The element, $\mathbf{1}_p \in [\mathcal{C}]$, is idempotent, so multiplication by it gives an idempotent linear map from $N$ to itself. We can thus split $N$ as $\mathbf{1}_p N\oplus (1-\mathbf{1}_p)N$, and we can repeat this with each object. We get $N\cong\oplus_{p\in \mathcal{C}_0}\mathbf{1}_pN$. We set $\widetilde{N}(p)$ to be the summand, $\mathbf{1}_p N$, and show that this is the object part of the required functor. If $c_{p,q}:p\to q$ in $\mathcal{C}$, then, as before,  let $\mathbf{c}_{p,q}\in [\mathcal{C}]$ be  the matrix having $c_{p,q}$ in position $(p,q)$.  If $n\in \widetilde{N}(p)$, then $\mathbf{1}_p\cdot n= n$, and $\mathbf{c}_{p,q}\cdot n= \mathbf{1}_q\cdot \mathbf{c}_{p,q}\cdot\mathbf{1}_p\cdot n$, which is in $\mathbf{1}_q N=\widetilde{N}(q)$.  We define $\widetilde{N}(c_{p,q}):\widetilde{N}(p)\to \widetilde{N}(q)$ by  $\widetilde{N}(c_{p,q})(n)=\mathbf{c}_{p,q}\cdot n$.  Linearity of this map is automatic.  The proof that  $\widetilde{N}:\mathcal{C}\to \Vect$ is a functor is fairly routine, as is that of the functoriality of the construction,  $\widetilde{(-)}:[\mathcal{C}]\!-\!Mod\to \mathcal{C}\!-\!Mod$. Finally it should be fairly clear  that this is the required quasi-inverse for $[-]$.
 \end{proof}

We have left `to the reader' the detailed verification that the above constructions do yield an equivalence between $\mathcal{C}\!-\!Mod$ and $[\mathcal{C}]\!-\!Mod$, as it is fairly routine to give a direct proof. We will, in fact, investigate that equivalence by a separate route. For this, we recall that $[\mathcal{C}]$, as it is a $\kappa$-algebra, can be considered as a $\kappa$-linear category in its own right, namely one having a single object, $\ast$, and with $[\mathcal{C}](\ast,\ast)$ being the set of elements of $[\mathcal{C}]$ itself, with composition being the multiplication in $[\mathcal{C}]$.
   We will not make any notational distinction between the $\kappa$-algebra, $[\mathcal{C}]$, and the linear category, $[\mathcal{C}]$, at least where no confusion is likely to arise by so doing.

\subsubsection{The bimodules underpinning Mitchell's Morita equivalence}

It is well known that two $\kappa$-algebras, $R$ and $S$, are Morita equivalent if there are bimodules, $_RA_S$ and $_SB_R$, such that the functors, $A\otimes_S-$ and $B\otimes_R-$, form an adjoint equivalence. What is also clear it that this should generalise to $\kappa$-linear categories and it does.  It does, however, seem a bit difficult to find a simple published proof of this, as it is a special case of some very wide ranging generalisations, whose generality we do not need, or, in fact, want here as our aim is to justify and interpret some calculations in a specific case of that general theory.

It does, however, suggest that we try to find `bimodules', $_\mathcal{C}A_{[\mathcal{C}]}$ and $_{[\mathcal{C}]}B_\mathcal{C}$, with similar properties. What are such `bimodules' to be? They are just another name for $\Vect$-valued profunctors, which, in the case of interest, would give $\mathbf{A}\colon \mathcal{C}\bto [\mathcal{C}]$ and $\mathbf{B}\colon [\mathcal{C}]\bto \mathcal{C}$. This observation, and quite a bit of what follows, is adapted from the $n$-Category Café discussion, (May 14, 2014), \cite{ncafe}, as mentioned before. (The ideas, there and here, were largely given by Karol Szumi\l{}o, but with a few additional features and verifications added here. We should add that  any errors should be attributed to us, and not to him.)

Before that, however, we will give the  pair of profunctors as suggested above.
Earlier, on and around page \pageref{identity profunctor}, we saw that, in $\Prof$ or $\vProf$, the identity profunctor on a category, $\mathcal{C}$, was the double Yoneda embedding,
$\mathcal{C}(-,-)\colon \mathcal{C}^{op}\times \mathcal{C}\to \Sets$,
or, of course, with codomain $\Vect$ if $\mathcal{C}$ is a $\kappa$-linear category.  We need, here, a functor, $\mathbf{A}\colon\mathcal{C}^{op}\times [\mathcal{C}]\to \Vect$, and an obvious candidate can be derived from that Yoneda based functor, by applying the $[-]$-construction to one side of it.  We, therefore, define a functor, $\mathbf{A}$, as required, by
$$\mathbf{A}(p,\ast)=\oplus_{q\in \mathcal{C}_0} \mathcal{C}(p,q).$$
We recall that, here, $\ast$ is the unique object of the $\kappa$-linear category, $[\mathcal{C}]$.   The formula is clearly (contravariantly) functorial in $p$, so it remains to see how some $\ast\xrightarrow{\mathbf{c}} \ast$, acts on $\mathbf{A}(p,\ast)$, where $\mathbf{c}$ is a matrix, $(c_{r,s})$, and each $c_{r,s}\in \mathcal{C}(r,s)$.

Let $p$ be an object of $\mathcal{C}$. Let $\mathbf{x}_{p}$ be an element of $\mathbf{A}(p,\ast)$. As we wrote before, given another object $q$, its $q$ component is some $x_{p,q}\in \mathcal{C}(p,q)$, and then
$$(\mathbf{x}_{p}\cdot \mathbf{c})_{p,s}= \sum_{q\in \mathcal{C}}x_{p,q}\cdot c_{q,s}.$$

We  want to calculate the composite, \smash{$\mathcal{C}\stackrel{\mathbf{A}}{\bto} [\mathcal{C}]\xrightarrow{N}\Vect$}, for a $[\mathcal{C}]$-module, $N$. We will write $N$ for both the module, and the functor, $N\colon[\mathcal{C}]\to \Vect$, although we should remember that $N$ is also $N(\ast)$, i.e. the functor evaluated on the single object of the algebra (considered as a linear category).

What should this mean?  The composite of a profunctor and a functor?  We can interpret this as being $\mathbf{A}\bullet \varphi^N\colon\mathcal{C}\bto \Vect$, so giving us a profunctor. (The notation $\varphi^N$ is  explained in Example \ref{Ex:maps_to_profunctors}). That would be a first step, thus we want to examine  the corresponding functor, $\mathbf{A}\bullet \varphi^N\colon\mathcal{C}^{op}\times \Vect \to \Vect$. We evaluate it on a pair of objects, $(p,V)$, with $p\in \mathcal{C}_0$ and $V$ being a vector space over $\kappa$.  The formula for the composition in $\vProf$ gives
$$(\mathbf{A}\bullet \varphi^N)(p,V) = \int^\ast \mathbf{A}(p,\ast)\otimes \Vect(N,V),$$
but we note that the coend is `integrating' over the one object category corresponding to $[\mathcal{C}]$, so is just the term being `integrated' divided by the diagonal action of the algebra, i.e.,  it is
$$(\oplus_{q\in \mathcal{C}_0}\mathcal{C}(p,q)\otimes \Vect(N,V))/\simeq,$$
where the action of any $\mathbf{c}_{r,s}$, which is homogeneous with value ${c}_{r,s}$, on the right hand side, $\Vect(N,V)$, is by the action on $N$, so if $v\colon N\to V$, then $(\mathbf{c}_{r,s}\cdot v)(n)=v(\mathbf{c}_{r,s}\cdot n)$, whilst on the left hand side, it is  by post-composition by $c_{r,s}$. Any element, $(c_{p,q},v)$ in this direct sum is $\simeq$-equivalent to one of the form $(1_p,w)$, by factoring the $c_{p,q}$ as $1_p\cdot c_{p,q}$, and then shifting the $c_{p,q}$ across to the other side. If we do this to $1_p$ itself, we find that $(1_p,w)\simeq (1_p,\mathbf{1}_p\cdot w)$, so is determined by the restriction of the linear map, $w$, to the direct summand $\mathbf{1}_pN$, which we have denoted above by $\widetilde{N}(p)$, as $\mathbf{1}_p\cdot w$ is the composition of $w$ with the projection onto that direct summand.  In other words,
$$(\mathbf{A}\bullet \varphi^N)(p,V)\cong \Vect(\widetilde{N}(p),V).$$
It is easy to see that this isomorphism is natural in both $N$ and $V$, so $(\mathbf{A}\bullet \varphi^N)$ is a representable profunctor, represented by $\widetilde{N}$. To summarise, the composite profunctor, $\mathcal{C}\stackrel{\mathbf{A}}{\bto} [\mathcal{C}]\xrightarrow{N}\Vect$,  is `really' the functor $\widetilde{N}$, as we hoped.

We now turn to the profunctor, $\mathbf{B}\colon[\mathcal{C}]\bto \mathcal{C}$, so $\mathbf{B}\colon[\mathcal{C}]^{op}\times \mathcal{C}\to \Vect$.  Given the success of the formula for $\mathbf{A}$ above, the `obvious' formula for $\mathbf{B}$ is
$$\mathbf{B}(\ast,q)= \oplus_{p\in \mathcal{C}_0}\mathcal{C}(p,q).$$
This certainly gives a functor, $[\mathcal{C}]^{op}\times \mathcal{C}\to \Vect$, (and thus a profunctor) as hoped for, and, to get the analogue of our earlier calculation, we will think of a functor, $M\colon[\mathcal{C}]\to \Vect$, as a profunctor, $\varphi^M\colon[\mathcal{C}]\bto \Vect$.  We have
$$(\mathbf{B}\bullet \varphi^M)(\ast,V)\cong \int^q\oplus_{p\in \mathcal{C}_0}\mathcal{C}(p,q)\otimes \Vect(M(q),V),$$
for $q\in \mathcal{C}_0$ and a vector space, $V$.  This is isomorphic to $\oplus_{p\in C_0} \Vect(M(p),V)$ by one of the forms of the co-Yoneda lemma, and this, in turn, is $\Vect([M],V)$, up to isomorphism.  These isomorphisms are natural, so $(\mathbf{B}\bullet \varphi^M)\cong \varphi^{[M]}$.  The composite profunctor, $(\mathbf{B}\bullet \varphi^M)$, is thus representable, and is `really' $[M]$.

This sets up the two functors, $[-]$ and $\widetilde{(-)}$, on the categories of `modules', as being given by the profunctors $\mathbf{A}$ and $\mathbf{B}$, respectively.  The final steps to explore, in this investigation of Theorem \ref{Mitchell-Morita}, are to calculate the composites, $\mathbf{A\bullet B}\colon\mathcal{C}\bto \mathcal{C}$ and $\mathbf{B\bullet A}\colon[\mathcal{C}]\bto [\mathcal{C}]$.  These are (slightly careful) manipulations involving the coend formulation of profunctor composition.
\begin{Proposition}\label{A.B and B.A}
(i) $\mathbf{A\bullet B}\cong \mathcal{C}(-,-)$, the unit profunctor on $\mathcal{C}$.

(ii) $\mathbf{B\bullet A}\cong [\mathcal{C}](\ast,\ast)\cong [\mathcal{C}]$, the unit profunctor / bimodule on $[\mathcal{C}]$.
\end{Proposition} 
\begin{proof}
(i) We take $p,q\in \mathcal{C}_0$, then

\begin{align*}
\mathbf{A\bullet B}(p,q)&=\int^\ast \mathbf{A}(p,\ast)\otimes \mathbf{B}(\ast,q)\\
&=\int^\ast \oplus_{r\in \mathcal{C}_0}\mathcal{C}(p,r)\otimes \oplus_{s\in \mathcal{C}_0}\mathcal{C}(s,q).
\end{align*}
As this coend is over (the single object category) $[\mathcal{C}]$, it can be calculated as $$(\oplus_{r\in \mathcal{C}_0}\mathcal{C}(p,r))\otimes_{[\mathcal{C}]} (\oplus_{s\in \mathcal{C}_0}\mathcal{C}(s,q)),$$so as a  tensor product over (the algebra),  $[\mathcal{C}]$, then, given  the form of the multiplication in  $[\mathcal{C}]$, it is clear that this  tensor product  is isomorphic to the vector space, $\mathcal{C}(p,q)$, and as all the isomorphisms are natural in $p$ and $q$, we  thus have that $\mathbf{A\bullet B}\cong \mathcal{C}(-,-)$, the unit profunctor on $\mathcal{C}$ as required.

(ii) This part is easier:
\begin{align*}
\mathbf{B\bullet A}(\ast,\ast)&=\int^q \mathbf{B}(\ast,q)\otimes \mathbf{A}(q,\ast)=\int^q \oplus_{p\in \mathcal{C}_0}\mathcal{C}(p,q)\otimes \oplus_{r\in \mathcal{C}_0}\mathcal{C}(q,r)\\
&\cong \oplus_{p,r}\mathcal{C}(p,r)\cong [\mathcal{C}],
\end{align*}
which is, of course, the same as $ [\mathcal{C}](\ast,\ast)$, as required.
\end{proof}
\begin{Remark}
This resolves, at least in part, the problem that we noted earlier, namely that $[-]$ is not a functor as such, at least in the most obvious sense.  Suppose, however, that $F\colon\mathcal{C}\to \mathcal{D}$ is a $\kappa$-linear functor, then there is an `induced' way to get from $[\mathcal{C}]$ to $[\mathcal{D}]$. It can be given by the composite profunctor,
$$\mathbf{B}_\mathcal{C}\bullet \varphi^F\bullet \mathbf{A}_\mathcal{D}\colon[\mathcal{C}]\bto \mathcal{C}\bto \mathcal{D}\bto [\mathcal{D}],$$
where we have indicated the `versions' of the profunctors, $\mathbf{A}$ and $\mathbf{B}$, by adding suitable suffices, e.g., $\mathbf{A}_\mathcal{D}$ being the $\mathbf{A}$ profunctor for $\mathcal{D}$, and so on. As this is a profunctor between two single object linear categories, it is `just' a left $[\mathcal{C}]$-, right $[\mathcal{D}]$-bimodule (determined up to isomorphism).

This, in fact, shows clearly that the bicategory of algebras, bimodules, and bimodule morphisms has some better properties than the category of algebras and algebra homomorphisms. The original functor induces a bimodule, but, in general, not a homomorphism, between the two category algebras. We could have got to this induced bimodule without going via the Morita context, but the route we have taken has some advantages for what we will be needing.
\end{Remark}

\begin{Remark}\label{Mitchel:map}
We note that the above construction for a functor extends easily to handling a profunctor, $\mathbf{H}\colon\mathcal{C}\bto \mathcal{D}$. Any such profunctor can be whiskered by suitable units and counits, $\mathbf{A}$ and $\mathbf{B}$, to give
$$\mathbf{B}_\mathcal{C}\bullet \mathbf{H}\bullet \mathbf{A}_\mathcal{D}\colon[\mathcal{C}]\bto \mathcal{C}\bto \mathcal{D}\bto [\mathcal{D}].$$

The discussion that we just gave can be used to prove that we have a bifunctor from the bicategory of linear categories, $\Vect$-enriched profunctors between them, and enriched natural transformations, to the bicategory of algebras, bimodules and bimodule maps. It sends $\mathcal{C}$ to $[\mathcal{C}]$, and $\mathbf{H}\colon\mathcal{C}\bto \mathcal{D}$ to the composite profunctor above. This construction will be clarified in \S\ref{sec:Mor_monoidal}. We will deal with a particular case of this latter construction in the following section.
\end{Remark}

We next turn to the bicategory of algebras, and describe it in a bit more detail, as it is the target for our next version of the once-extended TQFT.

    
\subsubsection{The bicategory,  $\Mor$}\label{sec:mor_defined}
Let us give a detailed definition of the bicategory of algebras, bimodules, and bimodule morphisms / intertwiners, that we have been using in a fairly sketchy form for some time. In so doing, we will shift our notation to put the actions of the algebras on the bimodules into a more central role.
 
 This bicategory is sometimes denoted $\mathbf{Alg}$ or $\mathbf{Alg_2}$ in the literature, but we will denote it by $\Mor$, and refer to it as the \emph{Morita bicategory}, as it is the natural and classical setting for Morita equivalence,  an adjoint equivalence in $\Mor$, in the sense of bicategory theory, being precisely a classical Morita equivalence.

We follow \cite{Bartlett_etal,Hansen-Shulman:constructing:2019,Schommer-Pries}, as well as more classical sources on bicategories.

\begin{Definition}\label{Def:Mor}
The \emph{Morita bicategory}, $\Mor=\Mor_\kappa$, is  such that:
\begin{itemize}[leftmargin=9mm]
 \item the objects of $\Mor$ are unital $\kappa$-algebras;
 \item given algebras $\A$ and $\B$, 1-morphisms $\A \bto \B$ are $(\A,\B)$-bimodules, $M$.   So $M$ is  a $\kappa$-vector space equipped with a left $\A$-representation / action, $\t$, and a right $\B$-representation,  $\trl$, that are compatible, meaning that given $a\in \A$, $b\in \B$ and $m\in M$, we have that $(a\t m)\trl b=a\t (m \trl b)$;
 \item the 2-morphisms, $F\colon (M\colon \A \bto \B) \To  (N\colon \A \bto \B)$, are given by $(\A,\B)$-bimodule maps, $F\colon M \to N$;
 \item the composition of $M\colon \A \bto \B$ with $N\colon \B \bto \Cc$ is $(M\otimes_{\B} N)\colon \A \bto \Cc$, where $M\otimes_{\B} N$ is the usual tensor product over $\B$, with its $(\A,\Cc)$-bimodule structure;
 \item and, finally, the horizontal composite of a compatible pair of 2-morphisms is
 $$\xymatrix{\A \ar@/^0.7pc/[rr]|{\not} \ar@/^0.7pc/[rr]^{M_1} \ar@/_0.7pc/[rr]|{\not} \ar@/_0.7pc/[rr]_{ {N_1}} & \Downarrow F_1 & \B \ar@/^0.7pc/[rr]|{\not} \ar@/^0.7pc/[rr]^{M_2} \ar@/_0.7pc/[rr]|{\not} \ar@/_0.7pc/[rr]_{ {N_2}} & \Downarrow F_2 & \Cc
}=   \xymatrix{\mortdash{\A}{\Cc}{ M_1\otimes_\B M_2}{ N_1\otimes_\B N_2}{F_1\otimes_\B F_2}}.$$
\end{itemize}
\end{Definition}
There are also 'well known' horizontal units and unitors, completing the construction of the bicategory $\Mor$, whose explicit description is left to the reader.

 Recall that $\vProfGrpfin$ is the full sub-bicategory of $\vProfGrphf$, with objects the finite groupoids. The constructions in \S \ref{non-functoriality} give a bifunctor, $\tLin\colon  \vProfGrpfin \to \Mor$. This bifunctor sends:
\begin{itemize}[leftmargin=6mm]
 \item each finite groupoid, $\Gamma$ to its \emph{groupoid algebra}, $\tLin(\Gamma)$;
 \item[]\hspace{-1cm}and
 \item each $\Vect$-produnctor $\Hp\colon \Gamma^\op \times \Gamma' \to \Vect$, to the bimodule, $\tLin(\Hp)$, with
 $$\tLin(\Hp):=\bigoplus_{x \in \Gamma_0,  \hspace{1mm} y\in \Gamma'_0} \Hp(x,y) .$$

 To describe the bimodule structure on $\tLin(\Hp)$, we let $a \in \Gamma_0 $ and $b\in\Gamma_0'$. Below, we will not distinguish between an element, $v_{(a,b)} \in  \Hp(a,b)$, and its image under the obvious inclusion of $\Hp(a,b)$ into  $\tLin(\Hp)$.
 The left and right actions of the algebras, $\tLin(\Gamma)$ and $\tLin(\Gamma')$, on $\tLin(\Hp)$ are such that, given
 $v_{(a,b)} \in \Hp(a,b)$, $(g\colon x \to y)\in \Gamma_1$ and  $(g'\colon x' \to y') \in \Gamma'_1$, we have:
 $$(x \ra{g} y) \t v_{(a,b)}=\begin{cases}  \Hp\big(x \ra{g} y, b \ra{\id_{b}} b\big)(v_{(a,b)}), &\textrm{ if } y=a,\\
                        0, &\textrm{ if } y \neq a,
                       \end{cases}
$$
and
$$ v_{(a,b)} \trl (x' \ra{g'} y') =\begin{cases}  \Hp\big(a \ra{\id_{a}} a, x' \ra{g'} y' \big)(v_{(a,b)}), &\textrm{ if } x'=b,\\
                        0, &\textrm{ if } x' \neq b.
                       \end{cases}
$$ \end{itemize}
 \begin{Remark}\label{Lin^(2) as reflection}
 The bimodule, $\tLin(\Hp)\colon \tLin(\Gamma) \bto \tLin(\Gamma')$, is an instance of the general construction mentioned at the end of \S\ref{non-functoriality}, namely $\tLin(\Hp)$ is isomorphic to the composite $$\mathbf{B}_\mathcal{C}\bullet \mathbf{H}\bullet \mathbf{A}_\mathcal{D}\colon[\mathcal{C}]\bto \mathcal{C}\bto \mathcal{D}\bto [\mathcal{D}],$$ where, here, $\mathcal{C}$ is $\tLin(\Gamma)$ and $\mathcal{D}$ is $\tLin(\Gamma')$.  This can help when checking, for instance, preservation, up to invertible 2-morphisms,  of horizontal composition for the candidate bifunctor, $\tLin\colon \vProfGrpfin \to \Mor$, see below.
 \end{Remark}
\begin{Remark}\label{actby1}
 Note that if $a \in \Gamma_0$ and $b\in \Gamma'_0$, then, for $v_{(a,b)}\in \Hp(a,b)$, we have that
 $$(a \ra{\id_a} a) \t v_{(a,b)}=  \Hp\big(a \ra{\id_a} a, b \ra{\id_{b}} b\big)(v_{(a,b)})=v_{(a,b)}$$ and
  $$  v_{(a,b)} \trl (b \ra{\id_b} b)=  \Hp\big(a \ra{\id_a} a, b \ra{\id_{b}} n\big)(v_{(a,b)})=v_{(a,b)}.$$
\end{Remark}

The remaining details of the verification that the above construction does give a  bifunctor, $\tLin\colon \vProfGrpfin \to \Mor_\kappa$,  will mostly be left to the reader. The key property that $\tLin$ preserves horizontal compositions of 1-morphisms, up to a canonical natural equivalence, is given by the following lemma.
\begin{Lemma} Consider finite groupoids, $\Gamma,\Gamma'$, and $\Gamma''$, and $\Vect$-profunctors, $ \Hp\colon \Gamma \bto \Gamma'$ and $\Hp'\colon \Gamma'\bto \Gamma''$. We have an isomorphism of $\big(\tLin(\Gamma), \tLin(\Gamma'')\big)$-bimodules,
$$\I\colon \tLin(\Hp\bullet \Hp') \To \tLin(\Hp) \otimes_{{\tLin(\Gamma')}}\tLin(\Hp).  $$
\end{Lemma}
\begin{proof}
As mentioned above, in Remark \ref{Lin^(2) as reflection}, this follows from the calculations in the previous section, and in particular on the properties of the composite profunctors, $\mathbf{A\bullet B}$ and $\mathbf{B\bullet A}$, as given in Lemma \ref{A.B and B.A}.  This is sketched in \S\ref{sec:Mor_monoidal}. We now give a direct proof, so as to accustom the reader to the links between profunctor and bimodule composition arguments.

 We first see what happens at the level of underlying vector spaces. Let $x\in \Gamma_0$ and $z \in \Gamma''_0$, then
 $$(\Hp\bullet \Hp')(x,z)=\int^{y \in \Gamma'_0} \Hp(x,y)\otimes\Hp'(y,z)=\big(\bigoplus_{y \in \Gamma'_0}  \Hp(x,y)\otimes\Hp'(y,z)\big )/\simeq.$$
Here, fixing $x \in \Gamma_0$ and $z \in \Gamma_0''$, the linear equivalence relation\footnote{i.e., an equivalence relation whose quotient is a vector space.}, $\simeq$, is generated by, for $y, y' \in \Gamma_0'$, $v_{(x,y)} \in \Hp(x,y)$ and  $v'_{(y',z)} \in \Hp(y',z)$, and an arrow, $y\ra{g} y'$, in $\Gamma'_1$,
$$  v_{(x,y)} \otimes \Hp'( y\ra{ g} y', z \ra{1_z} z)( v'_{(y',z)}) \simeq  \Hp(x \ra{1_x} x, y\ra{ g} y') (v_{(x,y)} ) \otimes  v'_{(y',z)}.  $$
The latter relation  means exactly that, given $y, y' \in \Gamma_0'$, $v_{(x,y)} \in \Hp(x,y)$ and  $v'_{(y',z)} \in \Hp(y',z)$, and an arrow, $y\ra{g} y'$ in $\Gamma'_1$, we have
$$  v_{(x,y)} \otimes (( y\ra{ g} y')\t  v'_{(y',z)}) \simeq  (v_{(x,y)} \trl ( y\ra{ g} y')) \otimes  v'_{(y',z)}.  $$
We also note that
 $$\tLin(\Hp\bullet \Hp')=\bigoplus_{x \in \Gamma_0, z\in \Gamma''_0} (\Hp\bullet \Hp')(x,z)  .$$

On the other hand, we have 
$$\tLin(\Hp) \otimes_{{\tLin(\Gamma')}}\tLin(\Hp') = \Big(\bigoplus_{x \in \Gamma_0, z\in \Gamma''_0}\quad  \bigoplus_{y,y' \in \Gamma'_0}  \Hp(x,y)\otimes\Hp'(y',z)\Big )/\sim.$$
Here the linear equivalence relation, $\sim$,  is such that, given $x\in \Gamma_0$, $z \in \Gamma''_0$,  $y, y' \in \Gamma_0'$, $v_{(x,y)} \in \Hp(x,y)$ and  $v'_{(y',z)} \in \Hp(y',z)$,  we have
$$  v_{(x,y)} \otimes (( w\ra{ g} w')\t  v'_{(y',z)}) \sim  (v_{(x,y)} \trl ( w\ra{ g} w')) \otimes  v'_{(y',z)},  $$
for arbitrary $(w\ra{ g} w')\in\Gamma'_1$. 

Clearly we have a bimodule map, $$\I\colon \tLin(\Hp\bullet \Hp') \To \tLin(\Hp) \otimes_{{\tLin(\Gamma')}}\tLin(\Hp)  ,$$ sending the equivalence class of 
$$v_{(x,y)} \otimes v'_{(y,z)}\in \bigoplus_{x \in \Gamma_0, z\in \Gamma''_0}\quad \bigoplus_{y \in \Gamma'_0}  \Hp(x,y)\otimes\Hp'(y,z) , $$ under $\simeq$, to the  equivalence class of
$$v_{(x,y)} \otimes v'_{(y,z)}\in \bigoplus_{x \in \Gamma_0, z\in \Gamma''_0} \quad\bigoplus_{y,y' \in \Gamma'_0}  \Hp(x,y)\otimes\Hp'(y',z),$$ under $\sim$, and we claim that $\I$ is a bijection.

If $y,y'\in \Gamma'_0$ are not equal, and $v_{(x,y)} \in \Hp(x,y)$ and  $v'_{(y',z)} \in \Hp(y',z)$, then $v_{(x,y)} \otimes  v'_{(y',z)}\sim 0$. This is because, (on using Remark \ref{actby1}),
\begin{align*}
v_{(x,y)} \otimes  v'_{(y',z)} &=\big (v_{(x,y)}\trl (y\ra{1_y} y) \big) \otimes  v'_{(y',z)} \\
&\sim v_{(x,y)}\otimes  \big( (y\ra{1_y} y) \big) \t   v_{(y',z)}\big) =0.
\end{align*}
In particular, $\I$ is surjective. 

We now define a bimodule  map, $$\I'\colon  \tLin(\Hp) \otimes_{{\tLin(\Gamma')}}\tLin(\Hp')  \To \tLin(\Hp\bullet \Hp').$$ There is a  bilinear map,  $\I''\colon  \tLin(\Hp) \times \tLin(\Hp')  \to \tLin(\Hp\bullet \Hp')$, given by, if $x \in \Gamma_0$, $y,y'\in \Gamma_0$ and $z \in \Gamma''_0$,  and also $v_{(x,y)} \in \Hp(x,y)$ and  $v'_{(y',z)} \in \Hp(y',z)$, then,
$$\I''( v_{(x,y)} ,v'_{(y',z)} )=\begin{cases}
                            [v_{(x,y)} \otimes v'_{(y',z)}]_{\simeq}, &\textrm{ if } y=y',    \\
                            0,  &\textrm{ if } y\neq y' .
                             \end{cases}
 $$
This is clearly balanced, considering the right and left actions of $\tLin(\Gamma')$, so $\I''$ descends to a 
linear map, $\I'\colon  \tLin(\Hp) \otimes_{{\tLin(\Gamma')}}\tLin(\Hp')  \to \tLin(\Hp\bullet \Hp').$ By construction, $\I'\circ \I=\id$, and so, in particular, $\I$ is injective as well.

The rest of the details are left to the reader.
\end{proof}

\subsubsection{The Morita-valued once-extended Quinn TQFT}\label{sec:Mor_valued-eTQFT}
Using the results of the previous sections, we can take our finitary version of the once-extended Quinn TQFT, defined in Subsection \ref{sec:fin_ext}, and reflect it into the bicategory $\Mor$, as follows.

As always,  we let $n$ be a non-negative integer and  $\Bc$ be a homotopy finite space.

\begin{Definition}\label{Sec:ext:Mor}
 The \emph{Morita valued once-extended Quinn TQFT},
$$\tFQmor{\Bc}\colon \tdcob{n}{\Bc} \to \Mor, $$
is defined as the following composite of bifunctors,
$$ \tdcob{n}{\Bc} \ra{\tFQd{\Bc}} \vProfGrpfin \ra{\tLin}\Mor. $$
\end{Definition}
\begin{Remark}\label{Decorations:indMod} (We  follow here an approach found  in \cite[Subsection 10.3]{Bullivant-thesis}.) Let $\Sigma$ be a closed smooth $n$-manifold. Given two $\Bc$-decorations,  $\fd_\Sigma$ and $\fd_\Sigma'$, of $\Sigma$, the same discussion as in \S\ref{Sec-Decorations:ind} gives a canonically defined invertible bimodule $\overline{\Psi}\big( \fd_{\Sigma}, \fd_{\Sigma'}\big)\colon \tFQmor{\Bc}(\Sigma,\fd_\Sigma) \bto \tFQmor{\Bc}(\Sigma,\fd_\Sigma')$, defined as:
$$\overline{\Psi}\big( \fd_{\Sigma}, \fd_{\Sigma'}\big):=\tFQmor{\Bc}\big ( (\Sigma,\fd_\Sigma) \ra{ (\iota_0^\Sigma,\Sigma\times I,\iota_1^\Sigma)} (\Sigma,\fd'_\Sigma) \big).$$
The discussion in \S \ref{Sec-Decorations:ind} passes over to \smash{$\tFQmor{\Bc}$} with the obvious modifications.
\end{Remark}
From the previous remark, all algebras,
 \smash{$\tFQmor{\Bc}(\Sigma,\fd_\Sigma)$}, where $\fd_\Sigma$ is a $\Bc$-decoration of $\Sigma$, are Morita equivalent. Crucially,  appropriate Morita equivalences can be canonically and functorially chosen, for any pair of decorations of $\Sigma$.

\subsection{The symmetric monoidal structure in $\tcob{n}$}\label{sec:sym_mon} We fix a non-negative integer $n$ throughout this subsection and the following as well. The central result of this paper is that one can categorify the finite total homotopy TQFT of Quinn, \cite{Quinn}, in a sensible way to get a once-extended TQFT,
$\tFQ{\Bc}\colon \tcob{n} \longrightarrow \vProfGrphf$.  From there we have shown that the resulting theory can be cut down in size to be more finitary by various means such as the introduction of decorations, and can be linked up with better known `algebraic' bicategories such as $\Mor$, which are frequently met in representation theoretic contexts.

Following Schommer-Pries,  \cite{Schommer-Pries},  Lurie, \cite{Lurie}, and others, we have taken a once-extended TQFT to be a symmetric monoidal bifunctor, as above, but we remark that the existing definitions do not agree on the target / codomain bicategory.  We have defined $\tFQ{\Bc}$, and have shown it to be a bifunctor.  There is, however, one further step to complete the proof that these constructions give  once-extended TQFTs, and that is to prove $\tFQ{\Bc}$, and its cousins, are symmetric monoidal bifunctors.  For this, we have to specify the symmetric monoidal structures on the cobordism bicategory, $\tcob{n}$, and will also recall that of $\vProfGrphf$, which was formally proved to exist in \cite{Hansen-Shulman:constructing:2019}.

We note that being a \emph{symmetric monoidal} bicategory or a \emph{symmetric monoidal}  bifunctor is a \emph{structure}, not a property, and refer the reader to the sketch in  Definition \ref{Sym. Mon. Bifunctor} and to  \cite[Definition 2.5]{Schommer-Pries} for a more detailed description. Sometimes the extra structure, i.e., that beyond being a bifunctor,  is `evident', but in our case that extra categorical structure encodes some of the `geometric' structure, for instance cobordisms, and 2-cobordisms, and we do need to have the transition between the various contexts made explicit to allow the naturality of the constructions to be made clear.

\subsubsection{A preliminary result towards the construction of the symmetric monoidal structure in $\tcob{n}$.}
The details of the construction of the symmetric monoidal structure in the bicategory $\tcob{n}$, using the language of symmetric monoidal pseudo-double categories \cite{Hansen-Shulman:constructing:2019}, can be found in \cite[\S 3.1.4]{Schommer-Pries}.
 In \cite[Remark 1.2.7.]{Lurie}, it is stated that, in the case of $\tcob{n}$ (or, more exactly, Lurie's analogue of this), the monoidal structure is straightforward, as \textit{``the tensor product operation is simply given by disjoint union of manifolds''}, just as in the more classical case of $\cob{n}$. Although correct, this statement hides some important details. The disjoint union of manifolds, cobordisms and extended cobordism indeed gives rise to a bifunctor, by abuse of language denoted\footnote{It is important to note that this bifunctor is not strict, in the sense that the natural isomorphism, $\varphi$, in item \eqref{natural-iso} of Definition \ref{def:bifunctor} is non-trivial.}
  $$\sqcup \colon  \tcob{n}\times \tcob{n} \to \tcob{n}. $$
 A monoidal bicategory is however not just the tensor product and unit, but also the associator, unitors and with additional pentagonators, etc., as sketched in Definition \ref{def-mon-bicat}.  Moreover, we need the tensor product to be symmetric, so need to specify a braiding, etc.  This may seem excessive detail to give, but is needed as there is a slight trap that has to be avoided, as we will now see.

As we saw in Definition \ref{def-mon-bicat}, in  a monoidal bicategory, $(\mathcal{A},\otimes, I, \ldots)$, we have a bifunctor $\otimes\colon \mathcal{A}\times \mathcal{A}\to\mathcal{A}$ and an adjoint equivalence, $\alpha\colon \otimes\circ (\otimes \times \id_{\mathcal{A}}) \To  \otimes\circ ( \id_{\mathcal{A}} \times \otimes)$. In the situation when $\mathcal{A}=\tcob{n}$, in particular we need to specify cobordisms, for all triples of $n$-manifolds, $C,B,A$,} \label{assoc in 2cob}
 $$\alpha_{CBA}=\vcenter{\xymatrix@R=-5pt{              (C\sqcup B)\sqcup A \ar[dr]_{i} && C\sqcup (B\sqcup A),\ar[dl]^{j}\\ & M&}}
$$
This is not difficult, but does involve some technicalities.  The two ends of the required cobordisms are not \emph{equal}, although they are naturally homeomorphic (and diffeomorphic if we include consideration of the smooth structure).  The required $(n+1)$-manifold, $M$, will be, topologically, a cylinder, but some care is needed with the embeddings $i$ and $j$, which must be specified.

We cannot directly use the diffeomorphism between the two sides, $(C\sqcup B)\sqcup A$ and  $C\sqcup (B\sqcup A)$, as the associator, as that would be a morphism in  $\Diff{n}$, but not in $\tcob{n}$.  We have to convert that isomorphism to a cobordism before checking that it works. For this, and for similar later situations, we need   a result which is, in some sense, a dual of Lemma  \ref{lem:functions_to_spans}, the context for which we will set up next.

In Section \ref{Classical Quinn}, we saw, in \jfm{\S\ref{moncob}}, that denoting by $\Diff{n}$, the category of closed $n$-manifolds and diffeomorphisms between them, we have a functor, $\I'\colon\Diff{n} \to \cob{n}$, which is the identity on objects.
We need a  categorified version, $\I$, of $\I'$.

We suppose that $f$ is a diffeomorphism from $X$ to $Y$. Write $\I(f)$ for the cobordism written down below, where $\iota^Y_k(y)=(y,k)$
\begin{equation}\label{def:I(f)}
\I(f):=\vcenter{\xymatrix@R=-5pt{X\ar[dr]_<<<<{\iota^Y_0f}&&Y\ar[dl]^<<<<{\iota_1^Y}.\\
&Y\times I&}}
\end{equation}
Note that this is a cobordism and not just the equivalence class determined by it.

This $\I$ does not give a `functor' from $\Diff{n}$ to $\tcob{n}$.
The reason is, essentially, that the horizontal composition in $\tcob{n}$ is that of a bicategory, not a category, since we are now not taking cobordisms up to diffeomorphism.  However, we instead have a pseudo-functor $\I\colon\Diff{n} \to \tcob{n}$, in the sense that we now describe.

Suppose that we have diffeomorphisms, $X\xrightarrow{f}Y\xrightarrow{g}Z,$
and thus two cobordisms  $\I(f)\colon X\to Y$ and $\I(g)\colon Y\to Z$, as well as $\I(gf)\colon X\to X$.  We can form $\I(f)\bullet \I(g)$ by the usual pushout and can put all this into a diagram as follows:
\begin{equation}\label{I-construction diffeo}\vcenter{
\xymatrix@R=13pt{X\ar[dr]^{\iota^Y_0f}\ar@/_2pc/[dddrr]_{\iota^Z_0 gf}&&Y\ar[dr]^{\iota^Z_0g}\ar[dl]_{\iota^Y_1}&&Z\ar[dl]_{\iota^Z_1}\ar@/^2pc/[dddll]^{\iota^Z_1}
\\
&Y\times I\ar@/_1pc/[ddr]|{\Psi^Y_{g,f}}\ar[r]_\ell&   PO(g,f) \ar@{-->}[dd]^{\Psi_{g,f}} &Z\times I\ar@/^1pc/[ddl]|{\Psi^Z_{g,f}}\ar[l]^r&\\
&&&&\\
&&Z\times I&&
}}\end{equation}
The pushout, $PO(g,f)$, is given  by $(Y\times I)\sqcup(Z\times I)/\sim$, where, for all $y\in Y$,
$(y,1)\sim (g(y),0).$  We have a homeomorphism,\label{PO homeo Z times I}
$\Psi_{g,f}\colon PO(g,f)\to Z\times I,$
given by $$\Psi^Y_{g,f}(y,t)=(g(y),t/2), \textrm{ and } \Psi^Z_{g,f}(z,t)= (z,(t+1)/2).$$

We note that the pushout in \eqref{I-construction diffeo} is a pushout in $\CGWH$. In order for this construction to be usable in the context of $\tcob{n}$, we must put a smooth structure on $PO(g,f)$, and also possibly modify $\Psi_{g,f}:PO(g,f)\to Z\times I$, slightly, in order that it is smooth at the junction, where the cylinders $Y\times I$ and $Z\times I$ join. These however can be easily handled using the usual mechanisms of collars, etc., so we will not concern ourselves more with this aspect.

Up to now, of course, this construction is very similar to what we used in our earlier section, \S \ref{I-construction}, to show that the uncategorified version of the construction gave a functor from  $\Diff{n}$ to $\cob{n}$, except that, as we already mentioned, we are now not taking the quotient of cobordisms by diffeomorphism (relative to the boundary). However it is not yet quite in the right form to be used for extended cobordisms, so as to give an extended cobordism / 2-cobordism between $\I(f)\bullet\I(g)$ and $\I(gf)$.  For that we use an analogue of the $\I$-construction one dimension up.

In general, suppose we have an isomorphism of cobordisms,
i.e. a diffeomorphism, $f$, making the diagram below commute,
\begin{equation}\label{eq: iso-cobs}\vcenter{
\xymatrix@R=1.3pt{ & M\ar[dd]^{f}_\cong& \\
X\ar[ru]^{i}\ar[rd]_{i'} && Y.\ar[lu]_{j}\ar[ld]^{j'}\\
&N&}}
\end{equation}
We can expand this out as a map of cospans,
\begin{equation}\label{eq: iso-cobs-cospan}
\vcenter{\xymatrixcolsep{1.5cm}\xymatrix{X\ar[r]^i\ar[d]_{id_X}&M\ar[d]^{ f}& Y\ar[l]_j\ar[d]^{id_Y}\\
X\ar[r]_{i'}&N&Y   \ar[l]^{j'},
}}
\end{equation}
 to which we apply the same idea as in the $\I$-construction to each vertical diffeomorphism to get the following extended cobordism, of dimension $n+2$,
\begin{equation}\label{eq: iso-cobs-Jconst}\J(f):=\vcenter{\xymatrixcolsep{1.5cm}\xymatrix@R=15pt{X\ar[r]_i\ar[d]_{\iota^X_0}&M\ar[d]^{\iota^N_0\circ f}& Y\ar[l]^j\ar[d]^{\iota_0^Y}\\
X\times I\ar[r]_{{i'}\times I}&N\times I&Y\times I\ar[l]^{j'\times I}\\
X\ar[r]^{i'}\ar[u]^{\iota^X_1}&N\ar[u]_{\iota^N_1}&Y.  \ar[u]_{\iota^Y_1}\ar[l]_{j'}
}}
\end{equation}
Passing to equivalence classes, we get a 2-morphism,\label{def:J}
$$[\J(f)]:(i,M,j)\To (i',N,j'),$$in $ \tcob{n}$. This is a vertically invertible 2-morphism.

Now suppose that we have diffeomorphisms of cospans, $f\colon (i,M,j)\to (i',N,j')$ and $g\colon (i',N,j')\to (i'',P,j''),$
as below,
$$\vcenter{
\xymatrix@R=0.8pc{ & M\ar[dd]^{f}_\cong& \\
X\ar[ru]^{i}\ar[rd]_{i'} && Y\ar[lu]_{j}\ar[ld]^{j'}\\
&N&}} \qquad \textrm{ and } \qquad  \vcenter{
\xymatrix@R=0.8pc{ & N\ar[dd]^{g}_\cong& \\
X\ar[ru]^{i'}\ar[rd]_{i''} && Y.\ar[lu]_{j'}\ar[ld]^{j''}\\
&P&}}$$
We can then compose them to get $gf\colon (i,M,j)\to (i'',P,j'')$. The 2-cospans, $\J(f)\colon (i,M,j)\To (i',N,j')$ and $\J(g)\colon (i',N,j')\To (i'',P,j''),$ equally well compose, using the vertical composition  given by the obvious pushout diagram, which fits into a diagram analogous to
the diagram, (\ref{I-construction diffeo}), above, but, of course,  replacing $X$, $Y$, and $Z$, with $M$, $N$ and $P$, respectively. (We leave the enterprising reader to extend this diagram to include what happens to the vertical cospans, $X\rightarrow X\times  I\leftarrow X$, etc.) The composite 2-cospan will be of form, 
$$\xymatrixcolsep{1.5cm}\xymatrix@R=13pt{X\ar[r]\ar[d]&M\ar[d]& Y\ar[l]\ar[d]\\
(X\times I)\sqcup_X(X\times I)\ar[r]&L&(Y\times I)\sqcup_Y(Y\times I)\ar[l]\\
X\ar[r]\ar[u]&P\ar[u]&Y,    \ar[u]\ar[l]
}$$
in which $L$ is given by the pushout,
$$\xymatrix@R=13pt{N\ar[r]^{\iota_0^P\circ g} \ar@{}[dr]|>>>\pushout\ar[d]_{\iota_1^N}&P\times I\ar[d]\\
N\times I\ar[r]&L.
}$$ There are diffeomorphisms, $(X\times I)\sqcup_X(X\times I)\xrightarrow{\cong}X\times I$, extending the obvious one from $I\sqcup_{\{*\} }I\to I$, and the discussion given after \eqref{I-construction diffeo} carries over to the setting here, giving an equivalence between $\J(f)\#_1 \J(g)$ and $\J(gf)$, so, in $\tcob{n}$,
$$[\J(f)]\#_1 [\J(g)]=[\J(gf)].$$

It should, now, be more-or-less clear  that we have a pseudo-functor,
 $$\I\colon\Diff{n} \to \tcob{n},$$
so we refer back to page \pageref{def:special pseudofunctor} for a checklist of structure and properties needed.  (We note that $\I$ is contravariant due to our notational convention for composition of cobordisms.)  In this setting,
\begin{itemize}[leftmargin=10mm]\item for each manifold, $X$, considered as an object of $\Diff{n}$, we have that $\I(X)$ is that same object considered as an object of $ \tcob{n}$, but note that we will write $X$ instead of $\I(X)$ most of the time in this context;
\item for each diffeomorphism, $f\colon X\to Y$, we have a 1-morphism  $\I(f)\colon X\to Y$;
\item for each composable pair, 
$X\xrightarrow{f}Y\xrightarrow{g}Z,$
we have an invertible  2-morphism, $$[\J(\Psi_{g,f})]\colon \I(f)\bullet \I(g)\To \I(gf);$$
\item[]\hspace{-12mm} and
\item for each object, $X$, $\I(id_X)$ is the chosen identity cobordism on $X$.
\end{itemize}

This leaves us to check compatibility of $\I$ with associators in $\tcob{n}$, namely that, given a triple of composable diffeomorphisms, $$X\xrightarrow{f}Y\xrightarrow{g}Z\xrightarrow{h}W,$$ the diagram
\begin{equation}\label{compatibility cocycle}\vcenter{\xymatrixcolsep{15ex}\xymatrix@R=10pt{
(\I(f)\#_0 \I(g))\#_0 \I(h)\ar@{=>}[r]^<<<<<<<<<<<{[\J(g,f)]\#_0 \I(h)}\ar@{=>}[dd]_{a}&\I(
gf)\#_0 \I(h)\ar@{=>}[dr]^{\quad[\J(h,gf)]}&\\
&&\I(hgf),\\
\I(f)\#_0 ( \I(g)\#_0 \I(h))\ar@{=>}[r]_<<<<<<<<<<<{\I(f)\#_0[\J(h,g)]}&\I(f)\#_0 \I(hg)\ar@{=>}[ur]_{\quad[\J(hg,f)]}&}}\end{equation}
commutes. Here we have written $hgf$ for the value of $(hg)f$ and $h(gf)$, which, of course, are equal, and have abbreviated $\J(\Psi_{g,f})$ to $\J(g,f)$ for ease of labelling the diagram. Furthermore, to emphasise that, here, it is the horizontal composition that is being used, we have replaced the  convenient, but `generic', symbol for composition,  $\bullet$, by the more specific one, $\#_0$.

We will formalise this in a proposition for ease of reference,  whose remaining details are left to the reader.
\begin{Proposition}\label{lemma:higherI}
There is a pseudo-functor, $\I\colon\Diff{n} \to \tcob{n},$ given as the identity on objects, and where the rest of the structure is as sketched above.
\end{Proposition}

\begin{Remark} We note that, if $g\colon Y \to X$ is the inverse diffeomorphism of $f\colon X \to Y$, then,
in $\tcob{n}$, $$\mathcal{I}(f)\bullet\mathcal{I}(g)\cong id_X, \qquad \textrm{  and } \qquad \mathcal{I}(g)\bullet\mathcal{I}(f)\cong id_Y.$$ These can be used to prove that, in the bicategory  $\tcob{n}$, $\mathcal{I}(f)$ forms part of an adjoint equivalence.
\end{Remark}

Finally, if $f\colon A\to B$ and $g\colon C\to D$ are diffeomorphisms, then we can form $f\sqcup g\colon A\sqcup C\to B\sqcup D$, and it is easy to see that $\mathcal{I}(f\sqcup g)\cong \mathcal{I}(f)\sqcup \mathcal{I}(g)$, again by the diffeomorphism coming from  $(B\sqcup D)\times I\cong (B\times I)\sqcup (D\times I)$. This implies that the pseudo-functor, $$\I\colon \Diff{n} \to \tcob{n},$$ is compatible with the coproduct monoidal structure. In particular, we will use this  when one of the two diffeomorphisms is the identity on the corresponding object.
\subsubsection{A sketch of the construction of the symmetric monoidal structure in  the bicategory, $\tcob{n}$}\label{sec:mon_tcob}

After this technical diversion, we can return to the problem of the monoidal associators in the monoidal bicategory, $\tcob{n}$, that we started discussing on page \pageref{assoc in 2cob}.  We need to formalise things a little more. For this, it may be helpful to give a reference for a fairly standard form of the axioms for a monoidal bicategory.  We will use Johnson and Yau, \cite[ \S1.2]{johnson-yau:bimonoidal:2021}, as a basic reference and will, in general, use their terminology.

For each triple of objects, $A$, $B$, $C$ in $\tcob{n}$, we seek  a cobordism,
$$ \alpha_{CBA}=\vcenter{\xymatrix@R=-5pt{              (C\sqcup B)\sqcup A \ar[dr]_{i} && C\sqcup (B\sqcup A).\ar[dl]^{j}\\ & M&}}
$$
We have, for $n$-manifolds, $A,B,C$, a diffeomorphism,
$a_{CBA}:(C\sqcup B)\sqcup A\to C\sqcup (B\sqcup A),$
and note that, as $(\Diff{n},\sqcup, \emptyset)$, forms a monoidal category, these satisfy the pentagon axiom, so for $A$, $B$, $C$ and $D$, the diagram,
\begin{equation}\label{pentagon:s5}
\vcenter{\scriptsize\xymatrix@C=-8pt@R=13pt{&(D\sqcup (C\sqcup B))\sqcup A\ar[rr]&&D\sqcup((C\sqcup B)\sqcup A)\ar[dr]&\\
((D\sqcup C)\sqcup B)\sqcup A\ar[ur]\ar[drr]&&&&D\sqcup(C\sqcup(B\sqcup A)),\\
&&(D\sqcup C)\sqcup (B\sqcup A)\ar[urr]&}}
\end{equation}
of manifolds and diffeomorphism, commutes.  We now write $$\alpha_{CBA}:=\mathcal{I}(a_{CBA}),$$ using the notation defined in \eqref{def:I(f)}. Given that we have a pseudo-functor, $$\I\colon\Diff{n} \to \tcob{n},$$ as shown in Proposition \ref{lemma:higherI}, whenever we have  two composable diffeomorphisms, $f$ and $g$, we have a 2-morphism,
 $$[\J(\Psi_{g,f})]\colon \I(f)\bullet \I(g)\To \I(gf),$$
 which satisfies the cocycle identity in \eqref{compatibility cocycle}. Applying this to the arrows in \eqref{pentagon:s5}, we can then derive an expression for the required pentagonator:
$$\scriptsize\xymatrix@C=-8pt@R=13pt{&(D\sqcup (C\sqcup B))\sqcup A\ar[rr]&&D\sqcup((C\sqcup B)\sqcup A)\ar[dr]&\\
((D\sqcup C)\sqcup B)\sqcup A\ar[ur]\ar[drr]\ar@{}[rrrr]|{\Downarrow \pi_{DCBA}}&&&&D\sqcup(C\sqcup(B\sqcup A)).\\
&&(D\sqcup C)\sqcup (B\sqcup A)\ar[urr]&}$$
By construction, as $\I\colon\Diff{n} \to \tcob{n}$ is a pseudo-functor, this pentagonator then  satisfies a higher order cocycle identity, as in \cite[Page 61]{Gurski_book} and \cite[Page 10]{GPS}, when we have five (closed and smooth) $n$-manifolds.

Also, given $(n+1)$-cobordisms, $(i_A,K,j_{A'})\colon A \to A'$,  $(i_B,M,j_{B'})\colon B \to B'$, and  $(i_C,N,j_{C'})\colon C \to C'$, we have a natural 2-morphism in $\tcob{n}$, fitting inside the diagram below,
 \begin{equation}\label{2natur}{\hskip-4.5cm\vcenter{
 \xymatrix@C=13pc{& (C\sqcup B)\sqcup A \ar[d]_{\alpha_{CBA}} \drtwocell<\omit>{\quad \qquad\alpha^2_{NMK}}\ar[r]^{\big( ( i_C,N,j_{C'}) \sqcup (i_B,M,j_{B'}) \big )\sqcup (i_A,K,j_{A'})
 } & (C' \sqcup B')\sqcup A' \ar[d]^{\alpha_{C'B'A'}}\\
 & C\sqcup (B\sqcup A) \ar[r]_{( i_C,N,j_{C'}) \sqcup \big ( (i_B,M,j_{B'}) \sqcup (i_A,K,j_{A'})\big )}
  & C' \sqcup (B'\sqcup A').\\
 }}}
 \end{equation}
 We note that, in the diagram above, we have abbreviated the notation, putting $$\alpha^2_{NMK}=\alpha_{\big ( ( i_C,N,j_{C'}) , (i_B,M,j_{B'}) , (i_A,K,j_{A'}) \big)}.$$

This 2-morphism,  $\alpha^2_{N,M,K}$,  arises from the obvious diffeomorphism between the $(n+1)$-cobordisms obtained from the two paths, from  $(C\sqcup B)\sqcup A $ to   $C'\sqcup (B' \sqcup A')$, in the diagram above, together with the construction in Equation \eqref{eq: iso-cobs-Jconst}. That diffeomorphism  underpins the naturality of the associativity constraints in $\cob{n}$, where the diagram consisting of the 1-dimensional arrows in \eqref{2natur} would commute. In the monoidal bicategory $\tcob{n}$, this diffeomorphism is unsurprisingly promoted to being a part of the symmetric monoidal
bicategory structure.

Together with the associator 1-morphisms, $\alpha_{CBA}$, the class of all 2-morphisms, $\alpha^2_{NMK}$, defines a pseudo-natural transformation of bifunctors,  $\alpha\colon( \tcob{n})^3 \to  \tcob{n}$, called the \emph{associator pseudo-natural transformation}, as shown below,
$$\xymatrixcolsep{3.3cm}
\xymatrix{
( \tcob{n})^3\ar[r]^{\sqcup\times( \tcob{n})}\ar[d]_{( \tcob{n})\times\sqcup} \drtwocell<\omit>{\alpha}&( \tcob{n})^2\ar[d]^\sqcup\\
( \tcob{n})^2\ar[r]_\sqcup &\tcob{n}
.}
$$

We have two different  bifunctors, from  $( \tcob{n})^4$ to  $\tcob{n}$, defined as $\sqcup \circ  (  \sqcup  \times \id  ) \circ(  \sqcup \times \id \times \id)$ and as $ \sqcup \circ  ( \id \times \sqcup) \circ   (\id \times \id \times \sqcup).$ Two different  pseudo-natural transformations between these bifunctors can be constructed using the  associator pseudo-natural transformation, $\alpha$, above, by considering the two different paths in  diagram  \eqref{pentagon:s5}.  The class of all pentagonators, $\pi_{DCBA}$, then defines a modification between the corresponding pseudo-natural transformations. This ``pentagonator modification''  satisfies its own cocycle identity, where we have five copies of $\tcob{n}$. The equation satisfied is in \cite[p. 61]{Gurski_book} and \cite[p. 10]{GPS}.

We can similarly use that the unit object in $(\CGWH,\sqcup, \emptyset)$, comes with natural isomorphisms,
$$\emptyset \sqcup A\xrightarrow{\ell_A}A \textrm{  and  }A\sqcup \emptyset \xrightarrow{r_A}A,$$
to obtain cospans, $\lambda_A:=\mathcal{I}(\ell_A)$ and $\rho_A:=\mathcal{I}(r_A)$.  These are just the obvious ones, but linking them with the construction of the pseudo-functor explicitly means that certain diagrams will immediately do what we need, without further checking.

The Middle Unity Axiom gives that 
$$\xymatrix{(A\sqcup \emptyset)\sqcup B\ar[r]^{a_{A,\emptyset, B}}\ar[d]_{{r_A}\sqcup B}&A\sqcup (\emptyset\sqcup B)\ar[d]^{A\sqcup \ell_B}\\
A\sqcup B\ar[r]_=&A\sqcup B,}$$
commutes, so, on applying $\mathcal{I}$, we get a specific modification,
$$\mu_{A,B}:(id_A\sqcup \lambda_B)\circ \alpha_{A,\emptyset, B}\to \rho_A\sqcup id_B,$$
in which $\circ$ stands for the composition of cospans.

The evident commutative diagram,
$$\xymatrix@R=13pt@C=13pt{(\emptyset\sqcup A)\sqcup B\ar[rr]^{\ell_A\sqcup B}\ar[dr]_{a_{\emptyset,A,B}}&&A\sqcup B\\
&\emptyset\sqcup (A\sqcup B)\ar[ur]_{\ell_{A\sqcup B}}&
}$$
after application of $\mathcal{I}$ gives a left 2-unitor, and the reverse / adjoint of $r$, denoted $r^*_A:A\to A\sqcup \emptyset$, likewise gives the right 2-unitor.

The fact that the pasting diagrams for these modifications work as required follows from the (trivially commutative) diagrams in  $(\CGWH,\emptyset,\sqcup)$ itself, on application of $\mathcal{I}$.  All this works in $\CGWH$, but we note that if the objects are smooth manifolds, the structure gives corresponding cobordisms as required.

Turning to the braiding, $\mathbf{R}$, on  $\tcob{n}$, the structural 1-morphisms are obtained as the image under $\mathcal{I}$ of the braiding, $\tau_{A,B}:A\sqcup B\cong B\sqcup A,$
in $\CGWH$, given by the universal property of the coproduct, so  given closed smooth $n$-manifolds, $A$ and $B$, we put $R_{A,B}=\mathcal{I}(\tau_{A,B})$. As for the case of the associator pseudo-natural transformation, given cobordisms, $(i_A,M,j_{A'})\colon A \to A'$ and  $(i_B,K,j_{B'})\colon B \to B'$, we have an extended cobordism, $$R^2_{M,N}=R_{((i_A,M,j_{A'}), (i_B,K,j_{B'}))},$$ fitting into the commutative diagram,
$$\xymatrixcolsep{5cm}\hskip-5.5cm\xymatrix@R=15pt{& A\sqcup B \ar[d]_{R_{A,B}} \drtwocell<\omit>{\qquad  R^2_{M,N}}\ar[r]^{   (i_A,M,j_{A'}) \sqcup (i_B,K,j_{B'})
 } & A'\sqcup B' \ar[d]^{R_{A',B'}}\\
 & B\sqcup A  \ar[r]_{ (i_B,K,j_{B'}) \sqcup (i_A,M,j_{A'})}
  & B'\sqcup A'. }
 $$

 Again, this extended cobordism arises from the obvious diffeomorphism between the two composite cobordisms from $A\sqcup B$ to $B'\sqcup A'$, obtained from the two paths from  $A\sqcup B$ to $B'\sqcup A'$ in the diagram above,  on applying  the construction in diagram \eqref{eq: iso-cobs-Jconst}. (Similarly to the associator natural transformation, this diffeomorphism underpins the naturality of the braiding in $\cob{n}$, but is now promoted to a crucial bit of structure in the symmetric monoidal bicategory $\tcob{n}$.) Together with the $R_{A,B}$, the class of all $R^2_{M,N}$ defines a pseudo-natural transformation of bifunctors, $\mathbf{R}$, fitting into the diagram below,
 $$\xymatrixcolsep{0.2cm}\hskip-0.3cm\xymatrix@R=17pt{ \tcob{n} \times \tcob{n} \ar[drr]_{\sqcup} \ar[rr]^{  \tau} &\drtwocell<\omit>{ \, \,\mathbf{R}}& \tcob{n} \times \tcob{n} \ar[d]^{\sqcup} \\
 &&\tcob{n}.}
 $$
 (Here the bifunctor $\tau$ is obtained simply by swapping coordinates.)
 Moreover, this is part of an adjoint equivalence, as in \cite[page 4234]{Gurski_coherence}.

As in the case of the associator natural transformation, to finish constructing a braiding in the monoidal bicategory $\tcob{n}$, we still need to specify modifications as in \cite[page 4235]{Gurski_coherence}, which we will not need explicitly here, and also check the remaining axioms for a braided monoidal bicategory, see \textit{loc cit}. Finally, this braiding satisfies the axioms for a braided monoidal bicategory to be a symmetric monoidal bicategory, which can be found in \cite[1.1. Definitions]{Gurski_Osorno}.

This finishes the sketch of the construction of the symmetric monoidal structure on $\tcob{n}$.

\subsection{The symmetric monoidal structure of the bifunctor $\tFQ{\Bc}$}\label{Quinn_is_sym_mon}
As usual, let $\Bc$ be a homotopy finite space, and recall that we  fix a non-negative integer, $n$, throughout this section.
\subsubsection{The basic case}\label{sec:basicTQFT}
We now sketch the proof of the fact that the bifunctor, $$\tFQ{\Bc}\colon \tcob{n} \to \vProfGrphf,$$ can be given the structure of a symmetric monoidal bifunctor, with respect to the symmetric monoidal structure, $\sqcup$, in $\tcob{n}$, whose construction we just sketched
 and the symmetric monoidal structure  in $\vProfGrphf$, mentioned in \S \ref{sum:conv_prof}. The latter symmetric monoidal structure is a particular case of that of the bicategory of $\Vect$-enriched productors, which is discussed in  \cite[Corollary 6.6]{Hansen-Shulman:constructing:2019}.

 The monoidal structure in $\vProfGrphf$ is essentially given as follows:
\begin{itemize}[leftmargin=6mm]
\item  on objects  it is given by the usual cartesian product of groupoids;
\item if $F\colon \mathcal{A}_0\bto \mathcal{B}_0$ and $G\colon \mathcal{A}_1\bto \mathcal{B}_1$ are 1-morphisms in $\vProfGrphf$, i.e. $\Vect$-profunctors, then  $F\otimes G\colon \mathcal{A}_0\times \mathcal{A}_1\bto \mathcal{B}_0\times \mathcal{B}_1 $ is given by the composite functor,
$$(\mathcal{A}_0\times \mathcal{A}_1)^{op}\times (\mathcal{B}_0\times \mathcal{B}_1)\xrightarrow{\cong}(\mathcal{A}_0^{op}\times\mathcal{B}_0)\times (\mathcal{A}_1^{op}\times \mathcal{B}_1)\xrightarrow{F\times G}\Vect\times \Vect \xrightarrow{\Otimes_\Vect}\Vect;$$
\item on 2-morphisms, the rule is $$(\alpha\otimes \beta)_{(A_0,A_1), (B_0,B_1)}=\alpha_{(A_0,B_0)}\otimes \beta_{(A_1,B_1)}.$$
\end{itemize}

With this information, we can construct a bifunctor,
$\Otimes\colon \vProfGrphf \times \vProfGrphf \to \vProfGrphf $,
which is the starting point for the construction of the symmetric monoidal structure on the bicategory $\vProfGrphf$. The remaining bits of structure look after themselves.

A crucial component for our discussion of the symmetric monoidal structure of the bifunctor,  $\tFQ{\Bc}\colon \tcob{n} \to \vProfGrphf$,  is the discussion in Lemmas \ref{Sym_mon:1} and \ref{def_chi},  and in \S \ref{sec:H-is-sym-monoidal}, which we need to transfer from $\mathbf{2span}(HF)$ to $\tcob{n}$, by using the mapping space construction $\Bc^{(-)}$; see Remark \ref{rem:tcobp} and Subsection \ref{sec:comm_summary} for notation. The notation for the additional bits of structure  that we will give to $\tFQ{\Bc}$ follows the pattern of the notation of Definition \ref{Sym. Mon. Bifunctor}, though we will add a prime to all structure morphisms, to distinguish the notation here from that already used in the context of  $\mathbf{2span}(HF)$.

We first construct a pseudo-natural transformation of bifunctors, fitting into the diagram,
$$\xymatrixcolsep{3cm}
\xymatrix@R=13pt{
( \tcob{n})^2\ar[r]^{\tFQ{\Bc} \times \tFQ{\Bc} }\ar[d]_{ \sqcup } \drtwocell<\omit>{\,\,\boldsymbol{\chi}'}& (\vProfGrphf)^2\ar[d]^\otimes\\
 \tcob{n}\ar[r]_{\tFQ{\Bc}} & \vProfGrphf
.}
$$
Given closed smooth $n$-manifolds, $X$ and $X'$, the cartesian closed structure of $\CGWH$ gives a natural isomorphism of groupoids, $$m'_{(X,X')}\colon  \pi_1(\Bc^{X},\Bc^{X}) \times \pi_1(\Bc^{X'},\Bc^{X'}) \to \pi_1(\Bc^{X \sqcup X'}, \Bc^{X \sqcup X'}).$$
We hence have a profunctor,  using the construction in Example \ref{Ex:maps_to_profunctors},
 $$\boldsymbol{\chi}'_{(X,X')}\colon  \pi_1(\Bc^X,\Bc^X) \times \pi_1(\Bc^{X'},\Bc^{X'}) \bto \pi_1(\Bc^{X \sqcup X'}, \Bc^{X \sqcup X'}),$$ defined as $\boldsymbol{\chi}'_{(X,X')}:=\varphi^{m'_{(X,X')}}$. Furthermore, given cobordisms, $(i,\Sigma,j)\colon X \to Y$ and $(i',\Sigma',j')\colon X' \to Y'$, and hence a cobordism, $\left (i\sqcup i', \Sigma \sqcup \Sigma', j\sqcup j'\right)\colon X \sqcup X'\to Y\sqcup Y'$,
we have a 2-morphism in $\vProfGrphf$,
$$
\xymatrixcolsep{2.5cm}\xymatrix{
\pi_1(\Bc^X)\times \pi_1(\Bc^{X'}) \ar[rr]^{\Hp\big(i^*,\Bc^\Sigma,j^*\big) \otimes \Hp \big({i'}^*,\Bc^{\Sigma'},{j'}^*\big) } \ar[d]_{ \boldsymbol{\chi}'_{(X,X')} } \drrtwocell<\omit>{\hspace*{2cm}\boldsymbol{\chi}'_{\left ( (i,\Sigma,j), (i,\Sigma',j')\right)}}&&
\pi_1(\Bc^Y)\times \pi_1(\Bc^{Y'}) \ar[d]^{\boldsymbol{\chi}'_{(Y,Y')}}\\
\pi_1(\Bc^{X \sqcup X'})
 \ar[rr]_{\Hp  \big( (i\sqcup i')^*, \Bc^{\Sigma \sqcup \Sigma'}, (j\sqcup j')^*\big)  }  &&  \pi_1(\Bc^{Y \sqcup Y'}),}
$$
 where we abbreviated  $\pi_1(B^X,B^X)$, etc, as $\pi_1(B^X)$.
Changing notation, we get
$$\hskip-1.3cm
\xymatrixcolsep{2cm}\xymatrix{
&\tFQ{\Bc}(X)\otimes \tFQ{\Bc}(X')   \ar[rr]^{\tFQ{\Bc} (i,\Sigma,j) \otimes \tFQ{\Bc} (i',\Sigma',j')}\ar[d]_{\boldsymbol{\chi}'_{(X,X')} } \drrtwocell<\omit>{\hspace*{2cm}\boldsymbol{\chi}'_{\left ( (i,\Sigma,j), (i,\Sigma',j')\right)}}&& \tFQ{\Bc}(Y)\otimes \tFQ{\Bc}(Y') \ar[d]^{\boldsymbol{\chi}'_{(Y,Y')} } \\
& \tFQ{\Bc}(X \sqcup X')
 \ar[rr]_{\tFQ{\Bc}\big (  i\sqcup i', \Sigma \sqcup \Sigma', j\sqcup j' \big) }  &&  \tFQ{\Bc}(Y \sqcup Y').}
$$
This natural isomorphism of profunctors is obtained from Lemma \ref{def_chi}.
By applying Lemma \ref{lem:chi-is-natural}, it follows that $\boldsymbol{\chi}'$ is a pseudo-natural transformation, $$\boldsymbol{\chi}'\colon \otimes \circ (\tFQ{\Bc}\times \tFQ{\Bc}) \to \tFQ{\Bc} \circ \sqcup.$$

Let us now  sketch the construction of the rest of the symmetric monoidal structure on $\tFQ{\Bc}$.
We have two  bifunctors $ L,R\colon ( \tcob{n})^3 \to \vProfGrphf$, defined by composition along the boundary left and right paths, in the two diagrams below, as in \cite[page 67]{Gurski_book}, as well as two natural transformations, connecting $L$ and $R$,
\begin{multline*}\allowdisplaybreaks
\xymatrixcolsep{1.5cm}
\xymatrix{ ( \tcob{n})^3 \ddtwocell<\omit>{<-8>\alpha}
\ar[d]_{\id\times \sqcup}\ar[dr]^{\sqcup\times \id} \drrtwocell<\omit>{<-1>\qquad \,\,\,\,\,\boldsymbol{\chi}'\times \tFQ{\Bc}} \ar[rr]^{\tFQ{\Bc} \times \tFQ{\Bc} \times \tFQ{\Bc}} && (\vProfGrphf )^3\ar[d]^{\otimes \times \id}\\
( \tcob{n})^2 \ar[d]_{\sqcup}& ( \tcob{n})^2 \ar[dl]^{\sqcup} \dtwocell<\omit>{ \boldsymbol{\chi}'}
\ar[r]^{\tFQ{\Bc} \times \tFQ{\Bc}} &  (\vProfGrphf )^2\ar[dl]^{\otimes}\\
 \tcob{n} \ar[r]_{\tFQ{\Bc}}& \vProfGrphf
}\\
\xymatrix{& \ar@3[r]_{\boldsymbol{\omega}'} &}  \xymatrixcolsep{1.5cm}
\vcenter{\xymatrix{ ( \tcob{n})^3
\ar[d]_{\id\times \sqcup} \drrtwocell<\omit>{<-1>\,\,\,\, \qquad \tFQ{\Bc}\times \boldsymbol{\chi}' }\ar[rr]^{\tFQ{\Bc} \times \tFQ{\Bc} \times \tFQ{\Bc}} && (\vProfGrphf )^3\ar[d]^{\otimes \times \id}\ar[dl]^{\id\times \otimes}   \ddtwocell<\omit>{<9>\alpha} \\
( \tcob{n})^2 \ar[d]_{\sqcup} \drtwocell<\omit>{\,\,\boldsymbol{\chi}' } \ar[r]^{\tFQ{\Bc} \times \tFQ{\Bc}}   & ( \vProfGrphf)^2   \ar[d]^{\otimes}
 &  (\vProfGrphf )^2\ar[dl]^{\otimes}.\\
 \tcob{n} \ar[r]_{\tFQ{\Bc}}& \vProfGrphf &
}}
\end{multline*}

The constructions, as discussed  in \S \ref{sec:H-is-sym-monoidal}, especially  in Notation \ref{def_omega} and Lemma \ref{Omega_is_natural}, give a modification, $\boldsymbol{\omega}'$, as shown above. Explicitly, given manifolds $X, X'$ and $X''$, and abbreviating, for a topological space $X$, $\pi(X)=\pi_1(X,X)$, we have that  $\boldsymbol{\omega}'_{({{X}},{{X}}',{{X}}'')}$ is a  natural isomorphism of profunctors fitting into the diagram,
$$
\hskip0cm
\xymatrix@R=15pt@C=-14pt{&\big(\pi(\Bc^{{X}})\times \pi(\Bc^{{{X}}'} )\big)\times \pi(\Bc^{{{X}}''})
 \ar[dl]\ar[dr]\\
\pi(\Bc^{{{X}}}) \times \big (\pi(\Bc^{{{X}}'})\times\pi(\Bc^{{{X}}''}) \big )\ar[d]
& \ar@{}[dd]\dtwocell<\omit>{  \boldsymbol{\omega}'_{({{X}},{{X}}',{{X}}'')}} &\pi(\Bc^{{{X}} \sqcup {{X}}' })\times \pi(\Bc^{{{X}}''})\ar[d]  \\
\pi(\Bc^{{X}}) \times \pi(\Bc^{{{X}}' \sqcup  {{X}}''}) \ar[dr]
&&\pi\big ( \Bc^{ ({{X}} \sqcup {{X}}') \sqcup {{X}}''}\big) \ar[dl]\\
&\pi\big (\Bc^{ {{X}} \sqcup ({{X}}' \sqcup {{X}}'')} \big   ).
}
$$
This diagram is the result of applying the pseudo-functor $\varphi^{(-)}\colon \Grp \to \vProfGrp$,  given in Examples \ref{Ex:maps_to_profunctors} and \ref{varphi is pseudo}, in which a functor is converted into a profunctor, to the following commutative diagram of groupoid functors, see the end of \S \ref{sec:pseud-cat-bicat}
$$
\hskip0cm
\xymatrix@R=15pt@C=-6pt{&\big(\pi(\Bc^{{X}})\times \pi(\Bc^{{{X}}'} )\big)\times \pi(\Bc^{{{X}}''})
 \ar[dl]_>>>>{\alpha^{(\Grp,\times)}_{\pi(\Bc^{{X}}),\pi(\Bc^{{{X}}'}),\pi(\Bc^{{{X}}''})} \qquad\,\,\,\,\,\,} \ar[dr]^>>>{\,\,\, m'_{({{X}},{{X}}')} \times \id_{\pi(\Bc^{{{X}}''} )} }\\
\pi(\Bc^{{{X}}}) \times \big (\pi(\Bc^{{{X}}'})\times\pi(\Bc^{{{X}}''}) \big )\ar[d]_{\id_{\pi(\Bc^{{X}})} \times {m'_{({{X}}',{{X}}'')}}}
&  &\pi(\Bc^{{{X}} \sqcup {{X}}' })\times \pi(\Bc^{{{X}}''})\ar[d]^{{m'_{({{X}}\sqcup {{X}}',{{X}}'')}}}  \\
\pi(\Bc^{{X}}) \times \pi(\Bc^{{{X}}' \sqcup  {{X}}''}) \ar[dr]_{{m'_{({{X}},{{X}}'\sqcup  {{X}}'')}}}
&&\pi\big ( \Bc^{ ({{X}} \sqcup {{X}}') \sqcup {{X}}''}\big) \ar[dl]^<<<<<< {{\pi\big(\Bc^{\big(\alpha^{\CGWH,\sqcup }_{{{X}},{{X}}',{{X}}''}\big)^{-1}}\big)}}
\\
&\pi\big (\Bc^{ {{X}} \sqcup ({{X}}' \sqcup {{X}}'')} \big   ).
}
$$
(We applied the associators for $(\Grp,\times)$ and $(\CGWH,\sqcup)$.)
For this reason, the modification, $\boldsymbol{\omega}'$,  satisfies the cocycle equation in \cite[\S 4.3]{Gurski_book}, or  \cite[page 17]{GPS}, if we are given manifolds ${{X}},{{X}}',{{X}}''$ and ${{X}}'''$.

The bifunctor, $\tFQ{\Bc}\colon \tcob{n} \to \vProfGrphf$, is also compatible with the unitor natural transformations in $\tcob{n}$ and in $\vProfGrphf$, as well as with the braiding. This follows from  considerations analogous to those we have just given. 
We  state for the sake of reference:

\begin{Theorem}\label{eQuinn is an eTQFT} The bifunctor, $\tFQ{\Bc}\colon \tcob{n} \to \vProfGrphf,$ of Subsection \ref{sec:once extended}, i.e. what we have called the once-extended Quinn TQFT,  can be upgraded to be a symmetric monoidal bifunctor.
\end{Theorem}
\subsubsection{A symmetric monoidal structure for the $\Bc$-decorated case}
The bicategory $\tcob{n}$ induces an obvious symmetric monoidal structure on the bicategory $\tdcob{n}{\Bc},$ defined in Subsection \ref{sec:fin_ext}. The tensor product of the decorated manifolds $(\Sigma,\fd_\Sigma)$ and $(\Sigma',\gd_{\Sigma'})$ is given by $$(\Sigma,\fd_\Sigma)\otimes (\Sigma',\gd_{\Sigma'})=(\Sigma \sqcup \Sigma', \fd_\Sigma \otimes \gd_{\Sigma'}),$$ where
$$ \fd_\Sigma \otimes \gd_{\Sigma'}:=\big\{ \langle \phi, \phi'\rangle \mid \phi \in \fd_\Sigma \textrm{ and }  \phi' \in \gd_{\Sigma'}\big\}.  $$

The discussion in \S\ref{sec:basicTQFT}  can easily be adapted for the case where we have decorated manifolds, again using the calculations in \S \ref{sec:H-is-sym-monoidal}. We hence have:
\begin{Theorem} \label{Decorated case is eTQFT}The bifunctor, \smash{$\tFQd{\Bc}\colon \tdcob{n}{\Bc} \to \vProfGrpfin$}, defined in  Subsection \ref{sec:fin_ext}, i.e., what we have called the finitary once-extended Quinn TQFT,
can be upgraded to being a  symmetric monoidal bifunctor.
\end{Theorem}

\subsubsection{A symmetric monoidal structure for the Morita valued once-extended Quinn TQFT}\label{sec:Mor_monoidal}
We now  sketch the construction of the symmetric monoidal structure of the Morita valued once-extended Quinn TQFT, 
 defined in  \S \ref{sec:Mor_valued-eTQFT}.
Given that the Morita valued once-extended TQFT  is obtained as a composite of bifunctors,
$$ \tdcob{n}{\Bc} \ra{\tFQd{\Bc}} \vProfGrpfin \ra{\tLin}\Mor ,$$
it will be sufficient to prove that the latter arrow,  from $\vProfGrpfin$ to $\Mor$, can be upgraded to being a  symmetric monoidal bifunctor.

 That $\tLin$ can be given a symmetric monoidal structure is a purely categorical / algebraic exercise, so we will just give  a sketch of that claim. Our main tool is  Mitchell's theorem, here Theorem \ref{Mitchell-Morita} on page \pageref{Mitchell-Morita}, for a linear category $\mathcal{C}$, with finitely many objects, using the approach that we took of constructing a bifunctor $[-]:\mathcal{C}\!-\!Mod\to [\mathcal{C}]\!-\!Mod$ as well as the  two profunctors, $\mathbf{A}_\mathcal{C}$ and $\mathbf{B}_\mathcal{C}$, giving a Morita equivalence; see the discussion in \S \ref{non-functoriality}, starting on page \pageref{non-functoriality}, and especially Remark \ref{Mitchel:map}. In fact, we will show how that theory allows one to define a bifunctor from the bicategory $\vProf$, of $\Vect$-enriched categories and $\Vect$-enriched profunctors between them, to $\Mor$, which is the `reflector' onto the sub-bicategory corresponding to $\Mor$.  We note that no finiteness or other restrictions are needed for this, which is why this construction works on  the whole of $\vProf$, and not just in $\vProfGrpfin$.

We define a bifunctor, $[-]\colon \vProf \to \Mor$, as follows:
 \begin{itemize}[leftmargin=0.6cm]\item  remembering that $\vProf_0$ consists of linear categories, we have that $[-]_0:\vProf_0 \to \Mor_0$ sends $\mathcal{C}$ to the algebra  $[\mathcal{C}]$;
 \item  if $H\colon \mathcal{C}\bto \mathcal{D}$ is a $\Vect$-valued profunctor, then  $[H]\colon [\mathcal{C}]\bto [\mathcal{D}]$ is the composite bimodule, obtained by the following composition of profunctors, as discussed in \S \ref{non-functoriality}, particularly in Remark \ref{Mitchel:map},
 $$\mathbf{B}_\mathcal{C}\bullet H\bullet \mathbf{A}_\mathcal{D}\colon [\mathcal{C}]\bto \mathcal{C}\bto \mathcal{D}\bto [\mathcal{D}].$$
Explicitly,
 $$[H]=\bigoplus_{p\in \mathcal{C}, q\in \mathcal{D}}H(p,q).$$
 The left and right algebra actions of $ [\mathcal{C}]$ and $[\mathcal{D}]$ are as discussed in   \S \ref{non-functoriality}.
 \item  We next assume given $H\colon \mathcal{C}\bto \mathcal{D}$ and $K\mathcal{D}\bto \mathcal{E}$. We do not expect that composition will be preserved by $[-]$, so look at the two ways of producing  things, namely $[H\bullet K]$ and $[H]\bullet [K]$. Firstly $$[H\bullet K]= \mathbf{B}_\mathcal{C}\bullet H \bullet K\bullet\mathbf{A}_\mathcal{E},$$ and then
  $$[H]\bullet [K] =  \mathbf{B}_\mathcal{C}\bullet H \bullet\mathbf{A}_\mathcal{D}\bullet\mathbf{B}_\mathcal{D}\bullet K \bullet\mathbf{A}_\mathcal{E}.$$
 (We will ignore any problems arising from composition being non-associative in $\vProf$, as these can be handled using associators, etc., in a standard way, completely analogous to handling non-associativity of tensor products in $\Mor$, which is a special case of this one.)
 
We recall, from Lemma \ref{A.B and B.A}, that the profunctors $\mathbf{A}_\mathcal{D}\bullet\mathbf{B}_\mathcal{D}$ and $\mathcal{D}(-,-)$, from $\mathcal{D}$ to itself, are naturally isomorphic, so we have an invertible 2-cell $[H]\bullet [K] \Rightarrow [H\bullet K],$ in $\Mor$, as  required, and moreover, this  satisfies the appropriate cocycle identity, given a triple of composable bifunctors.  We note that explicit formulae can be given for this 2-cell  in terms of the direct sum over the objects of $\mathcal{C}$ and $\mathcal{E}$, with tensoring
over the algebra $[\mathcal{D}]$.

 \item  We also need an invertible 2-cell, $id_{[H]}\Rightarrow [id_{H}]$.  The domain of this is the identity morphism on $[H]$, as a bimodule, and the righthand side is the  whiskered composite, $\mathbf{B}_\mathcal{C}\bullet id_{H}\bullet \mathbf{A}_\mathcal{D}$. This is easily checked to be isomorphic to the identity on $[H]$, for instance using that $[H]\cong\oplus_{p\in \mathcal{C}, q\in \mathcal{D}}H(p,q).$
\end{itemize}

 We omit the rest of the verification that this structure gives a bifunctor, $[-]\colon \vProf \to \Mor$, as that verification is quite long and not very insightful.

 We next ask if $[-]$ is monoidal, or more precisely whether it can be given a monoidal structure. We start with $F\colon\mathcal{A}\bto \mathcal{B}$ and $G\colon\mathcal{C}\bto \mathcal{D}$, and also have $F\otimes G\colon\mathcal{A}\times \mathcal{C}\bto \mathcal{B}\times \mathcal{D},$ given by
 $$(F\otimes G)((a,c),(b,c))= F(a,b)\otimes_\kappa G(c,d),$$
 so $$[F\otimes G]= \oplus_{(a,b)}\oplus_{(c,d)}F(a,b)\otimes_\kappa G(c,d).$$
 On the other hand, 
 $$[F]\otimes [G]=\big(\oplus_{(a,b)}F(a,b)\big)\otimes_\kappa \big(\oplus_{(c,d)}G(c,d)\big),$$
 and these two are naturally isomorphic by the standard argument relating tensors and direct sums.

Interaction of $[-]$ with the monoidal units is easy.  In $\vProf$, the monoidal unit is the single object linear category having a 1-dimensional vector space as its endomorphism ring, i.e., it is actually an algebra in its own right, being essentially a copy of $\kappa$ itself, and applying $[-]$ does nothing to it!  We thus have that $[-]$ is a \emph{normalised} monoidal bifunctor.

That $[-]$ respects the \emph{symmetric} monoidal bicategory  structure (up to specified isomorphisms) is then, once again, a result of the natural isomorphisms linking tensor products of direct sums with direct sums of tensor products, when that is suitably interpreted.

 We leave it to the reader  explicitly to write down the full symmetric monoidal bifunctor structure of  $[-]\colon \vProf \to \Mor$, as just outlined,  in the language of Definition \ref{Sym. Mon. Bifunctor}, similarly to that which  we did for the bifunctor $\tFQ{\Bc}\colon \tcob{n} \to \vProfGrphf$ in \S \ref{Quinn_is_sym_mon}.

The implication of the above is that the `$\Mor$-valued' once-extended TQFT, developed in Subsection \ref{sec:Mor_ext}, is, as we claimed, actually a fully-fledged once-extended TQFT,  as it can be given the structure of a symmetric monoidal bifunctor from $\tdcob{n}{\Bc}$ to $\Mor$.

This discussion leads to the following result.
\begin{Theorem} The bifunctor, $\tFQmor{\Bc}\colon \tdcob{n}{\Bc} \to \Mor,$  defined in \S \ref{sec:Mor_valued-eTQFT}, i.e., what we called  the Morita valued once-extended Quinn TQFT,
can be given the structure of a  symmetric monoidal bifunctor, (as explained before this theorem).
\end{Theorem}

\chapter{Calculations for classifying spaces of $\omega$-groupoids}\label{Quinn Calc}

We now have a TQFT and a once-extended TQFT that depend on the choice of a homotopy finite space, $\Bc$.
In this chapter we will show that, if we restrict to those homotopy finite spaces $\Bc$ that are classifying spaces of homotopy finite crossed complexes, then  we have the means for efficiently calculating the values of such TQFTs and once-extended TQFTs.   (We recall that homotopy finite crossed complexes are equivalent to strict $\omega$-groupoids, \cite[\S 13.6]{brown_higgins_sivera}, that are themselves homotopy finite).
Explicit formulae will appear in Section \ref{sec:TQFTS_xcomp}.

As before, let $n$ be a non-negative integer, and  $\Bc$ be a homotopy finite space. We have for  $s$, a complex parameter, Quinn's finite total homotopy TQFT, $$\FQp{\Bc}{s}\colon \cob{n} \to \Vect_\C,$$  (Definition
 \ref{def:quinnTQFT}), as well as the
finitary once-extended Quinn TQFT, $$\tFQd{\Bc}\colon \tdcob{n}{\Bc} \to \vProfGrpfin,$$ (Definition \ref{FinQ}), and the  Morita valued once-extended Quinn TQFT,
$$\tFQmor{\Bc}\colon \tdcob{n}{\Bc} \to \Mor,$$
(Definition \ref{def:finitary_ext_TQFT}). These can all be explicitly computed, combinatorially. This can be achieved in various ways, for example, by passing to the category of simplicial sets, or similar  combinatorial models for homotopy theory.

The calculations of $\FQp{\Bc}{s}$,  $\tFQd{\Bc}$ and   $\tFQmor{\Bc}$, using the category of simplicial sets, require only finite calculations, since, by  Ellis' theorem, \cite{Ellis},  path-connected spaces with a finite number of non-trivial homotopy groups, all of which are finite, can be represented (up to homotopy) by finite simplicial groups. We will examine this in a separate paper, as it requires some development of other techniques.

For  the remainder of this paper, we will outline the interesting special cases of such  explicit combinatorial computations  for the case in which $\Bc$ is the classifying space, $B_\A$, of a homotopy finite crossed complex, $\A$.  We apply the tools of the homotopy theory of crossed complexes, developed  in \cite{brown_higgins_classifying,brown_higgins_sivera,tonks:JPAA:2003}, and show that their use  yields   explicit formulae for $\FQp{B_\A}{s}$ and $\tFQd{B_\A}$, and hence for $\tFQmor{B_\A}$.
The formulae we will obtain for  $\FQp{B_\A}{s}\colon \cob{n} \to \Vect$ extend those of our previous paper, \cite{Martins_Porter}, which only dealt with the case of closed manifolds.

The category of crossed complexes, which strictly includes the category of strict 2-groups, \cite{Baez_Lauda}, and \cite[\S 2.5]{brown_higgins_sivera}, is equivalent to the category of (strict) omega-groupoids. For a precise statement and proof of this see  \cite[\S 13.6]{brown_higgins_sivera}. This latter fact justifies the title of this paper.  Crossed complexes also model strict $\infty$-groupoids via the nerve construction.

It should be noted that homotopy finite crossed complexes do not model all homotopy finite spaces. The specification of the homotopy types thus classified is slightly complicated, and will not be needed here, so will not be recalled in this paper. Such crossed complexes  do, however,  model all 2-types, $X$, that is,  all spaces, $X$, such that $\pi_i(X,x)=0$ if $i\ge 3$, and for all possible choices of base-point.

The explicit formulae we will construct, thus  apply to  the TQFTs $\FQp{\Bc}{s}$, and the extended TQFTs $\tFQd{\Bc}$ and $\tFQmor{\Bc}$, where $\Bc$ is the classifying space of a finite 2-group, as mentioned above, and  as in \cite{Martins_Porter}. In particular, the last part of this paper leads to a construction of TQFTs and extended TQFTs derived from discrete higher gauge theory based on a 2-group; see \cite{Yetter}, and also \cite{loopy,Companion}. In the extended case, the formulae are similar to those derived from the `tube algebras' considered in \cite{Bullivant_Tube} and \cite[Section 3]{Bullivant_Excitations}, in the context of excitations of strict 2-group topological phases.

\

In more detail, in the coming section, Section \ref{Sec:crossed_complexes}, we will review some of the basics of the homotopy theory of crossed complexes, and their classifying spaces, and then prove some refinements of well-known results in the literature, which will lead to the explicit formulae for TQFTs and once-extended TQFTs derived from crossed complexes mentioned above. This latter work will be done in Section \ref{sec:TQFTS_xcomp}.

 We note that subsections \ref{sec:defXcomp} --  \ref{sec:class_space} contain no new results, and essentially follow \cite{brown_higgins_tensor,brown_higgins_classifying,brown_higgins_sivera,tonks:JPAA:2003}. Subsection  \ref{sec:all_crossed_complex_fib}, on fibrations of crossed complexes, revisits definitions and results  from \cite{BrownGolasinski,brown_higgins_classifying}.  A crucial new result, refining the main theorem in \cite{brown_higgins_classifying},  concerns a crossed complex model for the fibre of the restriction map, on function spaces, $(B_\A)^{|S|} \to (B_\A)^{|T|}$, where $\A$ is a crossed complex, and $T$ is a subcomplex of a simplicial set, $S$, ($|-|$ here denoting geometric realisation, of a simplicial set). This result will be  used, later, to write  down, in all detail, TQFTs and extended TQFTs derived from homotopy finite crossed   complexes. Before that, the results in Subsection \ref{sec:homfinXcomp} are essentially in \cite{MartinsCW,Martins_Porter}, and  they allow for a simple calculation of the homotopy content of finite crossed complexes, akin to the well-known formula for the Euler characteristic of a finite CW-complex, as the alternating sum of cardinalities of the sets of $i$-cells.

\section{Crossed complexes: their homotopy theory and  classifying spaces}\label{Sec:crossed_complexes}
\sectionmark{Crossed complexes:  homotopy theory and classifying spaces}

The main sources for this section are \cite{brown_higgins_sivera}, and / or some of the sources already listed, which are summarised therein.
This section will, naturally, consist of lots of definitions, with some commentary.

\subsection{Definition of crossed complexes, and related notions}\label{sec:defXcomp} We first need some useful terminology.
\begin{itemize}[leftmargin=0.6cm]
\item Let $X$ be a set. By a \emph{set over} $X$, we will mean a set, $Y$, together with a surjective map, $\beta\colon Y \to X$. We denote this by $(Y,\beta)$.
\item A groupoid right-action of a groupoid, $\Gamma=(s,t\colon \Gamma_1 \to \Gamma_0)$, on  a set  over $\Gamma_0$, $(Y,\beta)$, is an operation which, given $y \in Y$ and an arrow, $(\beta(y) \ra{\gamma} x) \in \Gamma_1$, associates $y \trl \gamma$, also denoted $y \trl( \beta(y) \ra{\gamma} x)$, in $Y$, with $\beta(y \trl g)=x.$
This is such that if $\beta(y)=a$, we always have:
\begin{align*}
&\big(y \trl (a \ra{\gamma} b)\big) \trl (b \ra{\gamma'} c)=y \trl (a \ra{\gamma\gamma'} c), &\textrm{ and } &&
y \trl ({a \ra{1_a} a})=y.
\end{align*}
\end{itemize}
A groupoid action of $\Gamma$ gives rise to a functor $\Gamma \to \Sets$.

Various equivalent formulations of the definition  of a crossed complex can be found in \cite{baues1,Baues_4D,brown_hha,brown_higgins_sivera,tonks:JPAA:2003}, and many other places in the literature. Note that Baues, in  \cite{baues1}, and \cite{Baues_4D}, preferred to call them \emph{crossed chain complexes}.

For convenience, we recall one form of the definition here, from {\cite[\S 7.1.iii]{brown_higgins_sivera}}.
\begin{Definition}\label{def:crossed_complex}
 A \emph{crossed complex}, $\A=(A_n)_{n \in \mathbb{Z}^+_0}$, is given by:
 \begin{itemize}[leftmargin=1cm]
  \item a set, $A_0$, called the \emph{set of objects of $\A$};
  \item a groupoid, $A_1=(s,t\colon A^1_1 \to A_0)$, with object set $A_0$;
  \item for each integer $n\ge 2$, a totally disconnected groupoid,  with object set $A_0$, $$A_n=(\beta\colon A_n^1 \to A_0);$$
  \item whenever $n\ge 2$, a groupoid map $\d \colon A_n \to A_{n-1}$,  which is  required to restrict to the identity  on  the set of objects $A_0$;
  \item[]\hspace*{-1.2cm} and
  \item a groupoid right-action, $\trl$, of $A_1$ on all the underlying sets, $\beta\colon A_n \to A_0$, over $A_0$, for all $n\ge 2$, which is required to preserve the composition and the identities in each $A_n$.
 \end{itemize}
 Given $x, y \in A_0$, $a=(x \ra{a} x) \in {A_1}(x,x)$, and $g=(x \ra{g} y) \in {A_1}(x,y)$,  we also write $a \trl g:=g^{-1} a g \in A_1(y,y)$.

 This data is to satisfy the following additional conditions:
 \begin{enumerate}[leftmargin=1cm]
  \item for  all $n\ge 3$, given $(x\ra{a} x) \in A_n^1$ with $x \in A_0$, then $\d(\d (x\ra{a} x))=1_x$;
  
  \item  for  $n\ge 2$, and given any $x,y \in A_0$ and $g \in {A_1}(x,y)$,  if $a \in A_n(x,x)$, we have $\d(a \trl g)=\d(a)\trl g$; this is sometimes called the \emph{first Peiffer condition},
  \item for any $x\in A_0$ and $a,b \in A_2(x,x)$, then $a \trl \d(b)=b^{-1} a b$; this is sometimes called the \emph{second Peiffer condition}.
  \item If $n\ge 3$, then given any $x \in A_0$, and $a\in A_2(x,x)$, $b\in A_n(x,x)$, we have $b\trl \d(a)=b$, and so $\d(A_2)\leq A_1$ acts trivially on all $A_n$ for $n \ge 3$.
  \item  If $x \in A_0$ and $n\geq 3$, then each group, $A_n(x,x)$, is abelian.
 \end{enumerate}
 The arrows in $A_n$ will be called \emph{$n$-morphisms} and we may write $A_n(x)$ for $A_n(x,x)$.
 \end{Definition}
The notation, $\A=(A_n)_{n \in \mathbb{Z}^+_0}$, for a crossed complex leaves the boundary maps, $\partial$, and the actions of the groupoid, $A_1$, implicit in the notation.
Another useful way of picturing  a crossed complex is with the commutative diagram:
 \begin{equation}\label{eq:not-Xcomp}
 \vcenter{\xymatrix{  \A= &\cdots\to  A_4^1 \ar[r]^>>>>>{\partial}\ar[drrr]|{\beta} & A_3^1  \ar[drr]|{\beta}\ar[r]^{\partial} & A_2^1 \ar[dr]|{\beta}\ar[r]^{\partial}  & A_1^1\ar@/_0.5pc/[d]_{s} \ar@/^0.5pc/[d]^{t}
 \\
                    &&&& A_0.}}
\end{equation}

If $A_0$ is a singleton, and hence all groupoids $A_n$ have  a single object, as  will often be the case, we will write, again omitting the actions in our notation,
\begin{equation}\label{eq:Ared}
\A=\cdots \ra{\d}  A_{3}   \ra{\d}   A_2 \ra{\d}  A_1.
\end{equation}
Note that, here, we have identified each groupoid, $A^1_i$, with its group of morphisms.

This  latter type of  crossed complex, with a single object, is sometimes referred to as  being \emph{reduced}. The two sources, \cite{baues1,Baues_4D}, restrict attention to such reduced crossed complexes, calling them `crossed chain complexes', but we will need the non-reduced variety as we will be considering the crossed complexes corresponding to function spaces, where the restriction to a reduced case would be very unnatural.

Such a reduced crossed complex is  thus a chain complex of groups, such  that $A_i$ is abelian if $i\ge 3$, together with actions of the group, $A_1$, on all the groups in higher dimensions and, of course, satisfying some other axioms, as above.

 \begin{Definition}
 Let $\A=(A_n)_{n \in \mathbb{Z}^+_0}$ and $\B=(B_n)_{n \in \mathbb{Z}^+_0}$ be crossed complexes. A \emph{crossed complex map}, or \emph{morphism}, $f=(f_n)_{n \in \mathbb{Z}^+_0}\colon \A \to \B$, is given by a set map, $f_0\colon A_0 \to B_0$, and groupoid maps $f_i\colon A_i \to B_i$, where $i\ge 1$, that restrict to $f_0$ on objects. These are required to preserve the actions of $A_1$ and $B_1$, and the boundary maps in $\A$ and $\B$, in the obvious way.
\end{Definition}

Crossed complexes and the maps between them form a category which will be denoted $\Crs$.
As shown in  \cite[\S 7.2]{brown_higgins_sivera},  $\Crs$ is closed under small limits and colimits.

\begin{Definition}
For a positive integer $n$, a crossed complex, $\A=(A_i)_{i \in \mathbb{Z}^+_0}$, is said to be \emph{$n$-truncated} if, for $i>n$, the groupoids, $A_i$, have only identity morphisms.
\end{Definition}
\noindent We thus have that, for $i>n$, each of the groupoids, $A_i$, is a discrete / trivial groupoid on the set of objects, $A_0$.

We have  inclusion functors, $\J_1\colon \Grp \to \Crs$, sending a groupoid  to the obvious 1-truncated complex, and $\J_0\colon \Sets \to \Crs$, sending a set to the obvious $0$-truncated crossed complex, which is thus discrete in all dimensions.
These  inclusion functors have right adjoints, denoted $T_1$ and $T_0$ (respectively). These send $\A=(A_n)_{n \in \mathbb{Z}^+_0}$ to the groupoid, $A_1$,  and to the set,  $A_0$, again respectively.

The functors, $\J_1\colon \Grp \to \Crs$ and $\J_0\colon \Sets \to \Crs$, also  have  very useful left adjoints, $\pi_1\colon \Crs \to \Grp$, the \emph{fundamental groupoid functor}, and $\pi_0\colon \Crs \to \Sets$, the \emph{set of components functor}. In simple terms, these are given as follows:
\begin{Definition}[The functors, $\pi_1(\A)$ and $\pi_0(\A)$]\label{def:pi1andpi0crs} Let $\A=(A_n)_{n \in \mathbb{Z}^+_0}$ be a crossed complex. The \emph{fundamental groupoid of $\A$} is defined as:
$$\pi_1(\A):= A_1/\d(A_2).$$
(We will sometimes denote $\pi_1(\A)$ by the more suggestive $\pi_1(\A,A_0)$.) 
We also put: $$\pi_0(\A):=\pi_0(A_1),$$ the set of connected components of the groupoid, $A_1$, at the base of $\A$.
\end{Definition}
\noindent As the notation indicates, $\pi_1(\A,A_0)$ is a groupoid with one object for each element of $A_0$, and we note that $\pi_0(\A)=\pi_0(\pi_1(\A,A_0))$.

\begin{Definition}\label{def:fibXmap}
 Let $\A=(A_n)_{n \in \mathbb{Z}^+_0}$ and $\B=(B_n)_{n \in \mathbb{Z}^+_0}$ be
crossed complexes. Let $f=(f_n)_{n \in \mathbb{Z}_0^+}\colon \A \to\B$ be a crossed complex map, and let $b \in B_0$.
The \emph{fibre  of $f\colon \A \to \B$, at $b$,} is the sub-crossed complex,
$f^{-1}(b)=(C_n)_{n \in \mathbb{Z}^+_0}$,  of $\A$, with object set 
$C_0=f^{-1}_0(b)$, and such that $C_n$ consists of those elements,  $a \in A_n$, with $f_n(a)=1^{B_n}_b$, the identity of the groupoid, $B_n$, at $b$.
\end{Definition}
 As in the above definition, let $b \in B_0$, then we let $\hat{b}$ denote the sub-crossed complex of $\B$ with
object set $\{b\}$, and only identity arrows.  We clearly have a pullback diagram, where $\inc$ denote the obvious inclusions, of crossed complexes,
$$\xymatrix@R=13pt{
& f^{-1}(b)\ar@{}[dr]|<<<\pullback \ar[d] \ar[r]^{\inc}  & \A \ar[d]^f\\
                & \hat{b} \ar[r]_{\inc} &\B.
} $$

\begin{Definition}[$\PC_x(\A) $ and $\hpiz(\A)$]\label{Def:pizA} Let $\A=(A_n)_{n \in \mathbb{Z}^+_0}$ be a crossed complex.
 Given $x \in A_0$, we define the crossed complex, $\PC_x(\A)=(B_n)_{n \in \mathbb{Z}^+_0}$, as follows.
\begin{itemize}[leftmargin=1cm]\item The set, $B_0$, of objects of $\PC_x(\A)$ consists of all elements in $A_0$ connected to $x$ in the groupoid, $A_1$, so if $a\in B_0$, there is an arrow $a\to x$.

\item For each positive integer $n$, the set of morphisms in $B_n$ consists of the morphisms in $A_n$ connecting elements in $B_0$.
\end{itemize}
We call  the crossed complex $\PC_x(\A)$, the \emph{path-component of $x$ in $\A$.} \label{conn comp crs}

We also  write $\hpiz(\A)=\{\PC_x(\A)\mid x \in A_0\}$, for the collection of these path-components. (Note that different elements $x \in \A_0$ may induce the same $\PC_x(\A)$).

A crossed complex is called \emph{path-connected} if it only has one path-component.
\end{Definition}

Just as with spaces, crossed complexes come with a notion of homotopy groups. They are defined to be the obvious homology groups.

\begin{Definition}[Homotopy groups of crossed complexes]\label{def:pi_n crs}
 Let $\A=(A_n)_{n \in \mathbb{Z}^+_0}$ be a crossed complex. Let $c \in A_0$, and $n\ge 2$. We define $\pi_n(\A,c)$ to be the group,
$$ \ker\big (\d\colon A_n(c,c) \to A_{n-1}(c,c)\big)/ \im\big(\d\colon A_{n+1}(c,c) \to A_{n}(c,c)\big).$$
Also put $\pi_1(\A,c)=A_1(c,c)/\partial\big (A_2(c,c)\big)$.
\end{Definition}

%


\subsection{Fundamental crossed complexes of filtered spaces}\label{sec:Xcomp_Fil}
Many of the prime examples of crossed complexes come from filtered spaces, and, in particular, from CW-complexes considered with their natural skeletal filtration.  For a more complete view of this, see \cite[\S 7.1.i]{brown_higgins_sivera} and the development of related ideas there.

\subsubsection{Filtered spaces and crossed complexes} We first define filtered spaces.
\begin{Definition}
A filtered space, $X_*$, is a CGWH space, $X$, together with an increasing sequence, $X_0 \subseteq X_1 \subseteq X_2\subseteq \dots \subseteq X_n \subseteq \dots \subseteq X$, of subspaces of $X$. 

A filtered map, $f\colon X_*\to Y_*$,  between filtered spaces, is a continuous map, $f\colon X\to Y$, of the ambient spaces, such that $f(X_i) \subseteq Y_i$ for all $i \in \mathbb{Z}^+_0$. We let $\Fil$ denote the category of filtered spaces and filtered maps.

\begin{Notation}\label{not-sqcup-beta}If $X$ is a set and for each $x \in X$ we have a group, $E_x$, we consider the totally disconnected groupoid $\big(\beta \colon \sqcup_{x \in X} E_x  \to X\big)$. Here $\beta\colon \sqcup_{x \in X} E \to X$ is the map that identifies the component of the disjoint union to which an element belongs. The composition in each vertex group, $E_x$, is given by the product in $E_x$.
\end{Notation}

\end{Definition}
We have a crossed complex functor, $\Pi\colon \Fil \to \Crs$, sending a filtered space, $X_*$, to its fundamental crossed complex, $\Pi(X_*)$, defined as follows.

\begin{Definition} The \emph{fundamental crossed complex}, $\Pi(X_*)$, of a filtered space $X_*$, is specified by the following:\begin{itemize}[leftmargin=1cm]
\item the set of objects of $\Pi(X_*)$ is $\Pi(X_*)_0=X_0$;
\item the groupoid, $\Pi(X_*)_1$, is given by the fundamental groupoid, $\pi_1(X_1,X_0)$, of $X_1$, with set of base-points $X_0$;
\item if $n\ge 2$, let $\pi_n(X_n,X_{n-1},x)$ be the usual  relative homotopy group, and then  $$\Pi(X_*)_n=\big (\beta \colon \sqcup_{x \in X_0} \pi_n(X_n,X_{n-1},x) \to X_0\big );$$
\item for each $n\geq 2$ and $x\in X_0$, the boundary map, $\partial\colon \Pi(X_*)_n\to \Pi(X_*)_{n-1},$ is given by the  map appearing at the relevant position in the long homotopy exact sequence of the triple $(X_n,X_{n-1},X_{n-2})$;
 \item[]\hspace*{-1cm}and 
\item the action of $\pi_1(X_1,X_0)$ is the standard one.\end{itemize}\end{Definition} 
We direct the reader to  \cite[\S 7.1.v]{brown_higgins_sivera} for more explanation.

We note that referring to a CW-complex means that the space $X$ comes with a \emph{specified} CW-decomposition, and the cells used,   attaching and characteristic maps, etc.,  are all regarded as part of the structure.
Several filtrations appearing in this paper will be skeletal filtrations of CW-complexes, and these will be denoted $X_{\sk}:=(X^0 \subseteq X^1  \subseteq X^2\subseteq \dots )$, where $X^i$ is the $i$-skeleton of $X$.

As usual, a filtered map, $f\colon X \to Y$, between CW-complexes with the skeletal filtrations, is called \emph{cellular}.

We will give quite a few examples, so as to fix some notation.
\begin{Example}\label{Ex:simple_CW}  

Let $I=[0,1]$, with the standard CW-decomposition with two 0-cells, at $0$ and $1$, and one 1-cell. 

If $n\ge 1$, we let $S^n$ have the CW-decomposition with one $0$-cell, denoted $*$, for concreteness at the south pole, and one $n$-cell.

Let $D^{n+1}$ have the CW-decomposition for which $S^n$ is subcomplex, and we have an additional $(n+1)$-cell attaching along the identity map, $S^n \to S^n$. 

We have:\begin{itemize}[leftmargin=10mm,topsep=0pt] \item $\Pi(I_\sk)\cong\J_1\big(  \pi_1(I,\{0,1\})\big)$. Here, in $\pi_1(I,\{0,1\})$, we have objects $0$ and $1$, and only two non-identity morphisms, ${(0,1)}\colon 0 \to 1$ and ${(1,0)}\colon 1\to 0,$ so this is exactly the `unit interval groupoid,'  that we denoted by $\I$ earlier.
  \item $\Pi(S^1_\sk)\cong\dots \to  0 \to 0 \to 0 \to (\mathbb{Z},+) \to * $.
 \item $\Pi(D^2_\sk)\cong\dots \to  0 \to 0 \to  (\mathbb{Z},+) \ra{\id}  (\mathbb{Z},+) \to * $. Here the action of $ (\mathbb{Z},+)\cong \pi_1(S^1,*)$ on $\pi_2(D^2,S^1,*)\cong  (\mathbb{Z},+)$ is the trivial one.
   \item If $n\ge 2$, then  $\Pi(S^n_\sk)_n= (\mathbb{Z},+)$, $\Pi(D^{n+1}_\sk)_n=\mathbb{Z}$ and $\Pi(D^{n+1}_\sk)_{n+1}\cong  (\mathbb{Z},+) $, and all other groupoids are trivial. Hence:
   \begin{align*} \Pi(S^n_\sk)&\cong\dots  \to 0 \to 0 \,\ra{\,\,\,\,}  (\mathbb{Z},+)\to  0 \to \dots\to 0\to *, \\
\Pi(D^{n+1}_\sk)&\cong\dots \to  0   \to  (\mathbb{Z},+) \ra{\id}  (\mathbb{Z},+) \to  0 \to \dots\to 0\to *.
\end{align*}
\item For a positive integer $n$, let $\Delta(n)$ be the geometric $n$-simplex, the convex hull of $\{e_0,\dots, e_n\}\in \mathbb{R}^{n+1}$, with its obvious CW-decomposition, with $0$-cells at $\{e_0,\dots, e_n\}$. Then for instance:
$$ \Pi\big( \Delta(2)_\sk\big)\cong \dots \to   \bigsqcup_{x \in \{e_0, e_1,e_2 \}} \{0\}  \to  \bigsqcup_{x \in \{e_0, e_1,e_2 \}} (\mathbb{Z},+) \to \pi_1\big(\Delta(2)^1,\{e_0, e_1,e_2 \}\big).$$
Here the fundamental groupoid $\pi_1\big(\Delta(2)^1,\{e_0, e_1,e_2 \}\big)$, of the 1-skeleton of $\Delta(2)$, is the free groupoid on the three 1-cells. The action of $\pi_1$ trivially moves $(\mathbb{Z},+)$ between base-points.
\end{itemize}

\end{Example}

It will be useful to restrict further to what we will call \emph{special CW-complexes}. This ensures neater freeness properties of the fundamental crossed complex.

\begin{Definition}[Special CW-complex]\label{def:specialCW} A \emph{special CW-complex} is a CW-complex,  $X$, for which the attaching maps of all $n$-cells, for $n\ge 2$, are such that the unique $0$-cell of $S^{n-1}$ is sent to a 0-cell of $X^{n-1}$.
\end{Definition}

Let $\CW$ be the category of CW-complexes, (each provided with a specified CW-decomposition), and cellular maps. We let $\sCW$ be the full subcategory, of $\CW$, whose objects are the special CW-complexes.
The fundamental crossed complex functor, $\Pi\colon \Fil \to \Crs$, restricts to functors, $\Pi\colon \CW\to \Crs$, and $\Pi\colon \sCW\to \Crs$.

\subsubsection{Freeness of fundamental crossed complexes of CW-complexes}\label{sec:freeness-Xcomp}
The fundamental crossed complex of a CW-complex  $X$ is `free' on its cells; see \cite[Corollary 7.11]{brown_hha}.   This means that crossed complex maps $f\colon \Pi(X_\sk) \to \A$ can be specified by giving their value on each cell of $X$.  The latter assignment should be compatible with the boundary maps, of $\A$. We need some notation to make this precise.

Let us fix a CW-complex $X$, which, to simplify the exposition, we take to be special.
Given  $n\in \mathbb{Z}^+_0$, we let $C(X,n)$ be the set of $n$-cells of $X$. Given an  $n$-cell, $c\in C(X,n)$, we let $D^n_c=D^n$ and $S^{n-1}_c =S^{n-1}$, (the latter being empty if $n=0$), and let $i_c\colon S^{n-1}_c \to D^n_c$ be the inclusion. Supposing that $n\ge 1$, let $\psi_c\colon S^{n-1} \to X^{n-1}$ be the attaching map of $c$. Let $\phi_c\colon D^n_c \to X^n$ be the characteristic map of $c$. The inclusion of the $(n-1)$-skeleton $X^{n-1}$ into $X^n$, the $n$-skeleton of $X$, will be denoted $\iota_n\colon X^{n-1} \to X^n$, in the two diagrams below.

The following discussion uses \cite[Example  7.3.19 and Corollary 8.3.14]{brown_higgins_sivera}.
Given a positive integer $n$, have a pushout diagram in the category $\CGWH$,
\begin{equation}\label{eq:push-top}
\vcenter{\xymatrix@R=17pt{& {\displaystyle\bigsqcup_{c \in C(X,n)} S^{n-1}_c} \ar[dd]_>>>>>>>{\displaystyle \sum_{c \in C(X,n)} \psi_c} \ar[rrrr]^{\,\,\, \bigsqcup_{c \in C(X,n)} i_c\,\,\,}&&&&  {\displaystyle \bigsqcup_{c \in C(X,n)} D^{n}_c}  \ar[dd]^>>>>>>>>{\displaystyle \sum_{c \in C(X,n)} \phi_c} \\ \\
                    & X^{n-1} \ar[rrrr]_{\iota_n}  &&&& X^n. }}\end{equation}
                    (The vertical arrows arise from universal properties of disjoint unions.)

                    All maps appearing in the diagram above are cellular if $X$ is a special CW-complex, which is the reason for using them. Because the top horizontal arrow is a cofibration, by \cite[Theorem 8.2.5]{brown_higgins_sivera}, the diagram below is a pushout in $\Crs$,
            \begin{equation}\label{eq:push-crs}
\vcenter{\xymatrix@R=17pt{& {\displaystyle\bigsqcup_{c \in C(X,n)} \Pi(S^{n-1}_{c,\sk})} \ar[dd]_>>>>>>>{\displaystyle \sum_{c \in C(X,n)} \Pi(\psi_c)} \ar[rrrr]^{\,\,\,\ \bigsqcup_{c \in C(X,n)} \Pi(i_c)\,\,\,}&&&&  {\displaystyle \bigsqcup_{c \in C(X,n)} \Pi(D^{n}_{c,\sk})}  \ar[dd]^>>>>>>>{\displaystyle \sum_{c \in C(X,n)} \Pi(\phi_c)} \\ \\
& \Pi(X^{n-1}_\sk) \ar[rrrr]_{\Pi(\iota_n)}  &&&& \Pi(X^n_\sk). }}\end{equation}
Here $S^{n-1}_{c,\sk}$ and  $D^{n}_{c,\sk}$  denote the skeletal filtrations of  $S^{n-1}_{c}$ and $D^{n}_{c}$.
Furthermore, we have a natural isomorphism (of functors from $\sCW$ to $\Crs$),
$$\Pi(X_\sk)\cong \colim_{n} \big(\Pi(X^n_\sk),\Pi(\iota_n)\big).$$
This gives the freeness criteria that we mentioned before, see also \cite[page 238]{brown_higgins_sivera}.

Given an $(n+1)$-cell, $c \in C(X,n+1)$,  we have an induced map of pointed spaces, $\psi_c\colon (S^n,*) \to (X^n,\psi_c(*))$, and we let $\iota'(c)\in \pi_n(X,\psi_c(*))$ be the element given by the image of the generating element of $\pi_n(S^n,*)\cong \mathbb{Z}$. This gives an element, $\iota(c)\in \pi_n\big(X,X^{n-1},\psi_c(*)\big)$.

The following follows from the previous discussion, or \cite[page 238]{brown_higgins_sivera}.

\begin{Lemma}\label{lem:maps_on_basis}
Let $\A= \cdots \ra{\d} A_3 \ra{\d}  A_{2}  \ra{\d}   A_1$ be a reduced crossed complex. Crossed complex maps, $f=(f_n)_{n \in \mathbb{Z}^+_0}\colon \Pi(X_\sk) \to \A$, are in one to one correspondence with sequences of maps,  of sets, $\big(f_n'\colon C(X,n) \to A_n)_{n \in \mathbb{Z}^+}$, such that, for each $n$ and $c \in C(X,n)$, we have $f_{n-1}(\iota(c))=\d f'_n(c)$. Here, if we are given $f_i$ and $f'_{i+1}$, then $f_{i+1}$ is determined by the pushout in Equation \eqref{eq:push-crs}.
\end{Lemma}
This result can also be stated in the case when $\A$ is not reduced, but requires some additional conditions that source and target morphisms must match.

\begin{Remark}\label{rem:explicit_freeness} A way to state the freeness of $\Pi(X_\sk)=(\Pi(X_\sk)_n)_{n \in \mathbb{Z}^+_0}$  on the cells of $X$ is as follows.
\begin{itemize}[leftmargin=1cm]
\item
The groupoid, $\Pi(X_\sk)_1=\pi_1(X^1,X^0)$, is the free groupoid on the graph corresponding to the 1-skeleton of $X^1$ of $X$. In other words, $\Pi(X_\sk)_1$ is the free groupoid on the set of $1$-cells of $X$, and their attaching maps in $X^0$. \item The totally disconnected groupoid $\Pi(X_\sk)_2$, with $X_0$ as its set of objects, is the top groupoid of the free crossed  $\pi_1(X^1,X^0)$-module, cf.  \cite[7.3.ii]{brown_higgins_sivera} on the attaching maps for the 2-cells. For an explicit description, see  \cite[\S 3.3]{Companion}.
\item
If $n\ge 3$, then the totally disconnected abelian groupoid, $\Pi(X_\sk)_n$, is a free $\pi_1(X^2,X^0)$-module over the set of $n$-cells of $X$, and the boundary map, $\d\colon \Pi(X_\sk)_n \to \Pi(X_\sk)_{n-1}$, is derived from the attaching maps of the $n$-cells.
\end{itemize}
For more details, see \cite[Definition 7.3.13]{brown_higgins_sivera} and \cite[Chapter III]{Baues_4D}, also reviewed in \cite[\S 2.2.1]{MartinsCW}, and in \cite{Martins-surface} for fundamental crossed modules.
\end{Remark}

\subsubsection{Maps from the fundamental crossed complex of $\Delta(n)$}\label{PiDeltan}
The set of $i$-cells of the geometric $n$-simplex, $\Delta(n)$ in Example \ref{Ex:simple_CW},
 is in bijection with the set of sequences, $(a_0,a_1,\dots, a_i)$, with  $0\leq a_0<a_1<\dots < a_i\leq n$. Let $\A$ be a crossed complex, as in \eqref{eq:not-Xcomp}. By the discussion in \S \ref{sec:freeness-Xcomp}, a crossed  complex  map  $f\colon \Pi\big(\Delta(n)_\sk\big) \to \A$ is uniquely specified by its values on the cells of $\Delta(n)$. Hence $f\colon \Pi\big(\Delta(n)_\sk\big)\to \A$ is given by the following information, and no further compatibility conditions are required (this is as in \cite[p. 99]{brown_higgins_classifying} / \cite[\S 9.9]{brown_higgins_sivera}):
\begin{itemize}[leftmargin=0.6cm]
\item A map, $f_0\colon \{0,\dots,n\}\to  A_0$, so picking out $n+1$  objects of $\A$.
 \item An assignment,  $f_1(a,b)$ of a morphism of $A_1$ to each $1$-cell, $(a,b)$, so with $0\leq a<b \leq n$, that goes between the images of the vertices, $a$ and $b$,  meaning that we have
 ${f_1(a,b)}\colon f_0(a) \to f_0(b)$.
 \item An assignment $f_2(a,b,c)$, of a morphism  in the vertex group, ${A_2}\big(f_0(a),f_0(a)\big)$, to each 2-cell, $(a,b,c)$, such that
$\d(f_2(a,b,c))=f_1(a,b)\, f_1(b,c)\, f_1(a,c)^{-1}$.  The images of $(a,b,c)$ and of its faces thus match together as below:
  $$ \vcenter{\xymatrix{ &b \ar[dr]^{(b,c)}\ar@{}[d]|>>>>>{(a,b,c)}&\\
a\ar[ur]^{(a,b)}\ar[rr]_{(b,c)}&& c}} \qquad  \mapsto \qquad \vcenter{\xymatrix{ &f_0(b) \ar[dr]^{f_2(b,c)}\ar@{}[d]|>>>>>{f_2(a,b,c)}&\\
f_0(a)\ar[ur]^{f_1(a,b)}\ar[rr]_{f_1(b,c)}&& f_0(c).}}$$
 \item An element
  $f_3(a,b,c,d)\in {A_3}(f_0(a),f_0(a))$ for each 3-cell, $(a,b,c,d)$, such that the following holds\footnote{The diagram below shows that we can divide the boundary of $\Ds(3)$ into two parts,
$$\xymatrix{b\ar[r]\ar[dr]^{d_0}_{d_2}&c\ar[d]&&b\ar[r]&c\ar[d] \\
a\ar[r]\ar[u]&d&&a\ar[u]\ar[r]\ar[ur]^{d_3}_{d_1}&d
}$$
This  explains how Equation \eqref{eg:hal-3simplex} arises, as the boundary is the difference between the two parts.}
\begin{equation}\label{eg:hal-3simplex}\d(f_3(a,b,c,d))=\big (f_2(b,c,d)\trl f_1(a,b)^{-1}\big)\,  f_2(a,b,d)\, f_2(a,c,d)^{-1}\,  f_2(a,b,d)^{-1}.\end{equation}

 \item[]\hspace*{-1.1cm} Finally:

 \item for $3<i\leq n$, and each $i$-cell $(a_0,\dots, a_i)$, an element $$ f_i(a_0,\dots, a_i))\in {A_n}\big(f_0(a_0),f_0(a_0)\big),$$  such that, putting $d_j(a_0,\dots, a_i)=(a_0, \dots, \hat{a_j}, \dots a_i)$ we have:
\begin{align*}\d (f_i(a_1,\dots, a_i))=&\\ \big(f_{i-1}(d_i(a_0,&\dots, a_i)) \trl f_1(a_0,a_1)^{-1}\big)  \prod_{{j}=1}^n \big(f_{i-1}(d_{{j}}(a_0,\dots, a_i))\big){}^{(-1)^{{j}}}.\end{align*}
\end{itemize}

\subsubsection{The fundamental crossed complex of a simplicial set}\label{sec:Pi(S)}We freely use the notion of a simplicial set; see e.g,  \cite{Curtis,FP, MaySimplicial} and numerous other places in the literature and on-line. The category of simplicial sets will be denoted by $\Simp$.
The geometric $n$-simplex is, as above, denoted $\Delta(n)$.

Let us explain our convention and notation. Let $\Delta$ be the \emph{simplex category}. The objects of $\Delta$ are non-negative integers, $n$,  or more exactly the finite ordinals,  $[n]= \{0<1<\ldots <n\}$, and the morphisms from  $[m]$ to $[n]$, are the non-decreasing maps, $[m]\to [n]$. A simplicial set, $S$, is then a functor, $S\colon \Delta^{op} \to \Sets$, and, if $n$ is a non-negative integer, the set, $S_n$, of $n$ simplices is the image, $S(n)$, of $[n]$ under $S$. We have face maps $d_i:=d^n_i\colon S_n\to S_{n-1}$, for $0\leq i\leq n$, each of which arises from the unique strictly increasing map $[n-1]\to [n]$ whose image does not contain $i$.
Similarly we have degeneracy maps  $s_j:=s^n_j\colon S_n\to S_{n+1}$,  for $0\leq i\leq n$, each of which arises from the unique surjective non-decreasing map $ [n+1]\to [n]$ that repeats $i$. An $n$-simplex is called \emph{non-degenerate} if it is not in the image of any degeneracy map.

Given a non-negative integer $n$, we  let $\Ds(n)$ be the simplicial $n$-simplex, so $\Ds(n)\colon \Delta^{op} \to \Sets$. This is defined to be the representable functor, $\Ds(n)(m):= \Delta([m],[n])$.
 The set, $\Ds(n)_m$, of $m$-simplices of $\Ds(n)$ is, thus,  the set of non-decreasing maps, $\sigma$, from $[m]$ to $[n]$.  Such an $m$ simplex can be represented by a string, $(a_0,a_1,\dots, a_m)$, where $a_k=\sigma(k)$, and thus we have, $0\leq a_0\leq a_1\leq \dots \leq a_m\leq n$. Each face map $d_i\colon \Ds(n)_m \to \Ds(n)_{m-1}$, where $i=0,\dots, m$, is obtained by omitting the $i^{th}$ entry of $(a_0,a_1,\dots, a_m)$.
An $m$-simplex, $\sigma$, is non-degenerate exactly when the  map, $\sigma$, is injective, so we have $0\leq a_0< a_1< \dots < a_m\leq n$.

We can express any simplicial set as a coend, i.e., as a colimit, of copies of standard simplices. This gives, in its simplest  form,\label{S as a coend}
$$S\cong \int^{n\in \Delta} S_n\times \Ds(n).$$
 This interprets as taking lots of labelled copies of the various standard simplices, and then glueing them along common faces, also taking into account the degeneracies. The geometric realisation of $S$ is then
$$|S|= \int^{n\in \Delta} S_n\times \Delta_n.$$
This $|S|$ is  a special  CW-complex, with one $n$-cell for each non-degenerate $n$-simplex of $S$; see e.g. \cite[Theorem 4.3.5]{FP}.

The geometric realisation of a map between simplicial sets is cellular. Combining geometric realisation with the fundamental crossed complex functor, we, therefore, have a functor, $\Pi\colon \Simp \to \Crs$, which we will refer to as the \emph{fundamental crossed complex functor}, sending a simplicial set $S$ to $\Pi(S):=\Pi(|S|_\sk)$.

We note  the following, which is Proposition 2$\cdot$2 of  \cite{brown_higgins_classifying}.
\begin{Theorem}\label{Th:fundCRS-coend} Let $S$ be a simplicial set. Let $|S|$ be its geometric realisation, with the skeletal filtration. We have a natural isomorphism
 $$\Pi(|S|_\sk)\cong \int^{n\in \Delta}  S_n \times \Pi\big(\Delta(n)_\sk\big).$$
\end{Theorem}
Combining this with the discussion in \S \ref{PiDeltan}, we can see that if $\A$ is a crossed complex then crossed complex maps $f\colon \Pi(S)\to \A$ can be specified combinatorially.

%

\subsection{Homotopy of crossed complexes}\label{sec:xcomphom}A source for much of this review of the homotopy of maps of crossed complexes  is  \cite[\S 9.3]{brown_higgins_sivera}. The particular case of homotopy of crossed modules (of groupoids) is in \cite{Brown_Icen} and, also in  \cite[\S 2.6.1]{Martins_Porter}. The version for reduced crossed complexes is given by Baues, \cite[page 98]{Baues_4D}.

We will give a short description of the notion of homotopy of crossed complex maps, focusing on showing the particular explicit formulae that we will need to write down the TQFTs and once-extended TQFTs derived from finite crossed complexes.

Throughout this subsection, we fix two crossed complexes, $\A$ and $\B$, as below,
$$ \xymatrix@C=22pt{  \A=\cdots  A_4^1 \ar[r]^>>>>{\partial}\ar[drrr]|{\beta} & A_3^1  \ar[drr]|{\beta}\ar[r]^{\partial} & A_2^1 \ar[dr]|{\beta}\ar[r]^{\partial}  & A_1^1,\ar@/_0.5pc/[d]_{s} \ar@/^0.5pc/[d]^{t}
 \\
                    &&& A_0}  \xymatrix@C=22pt{  \B=\cdots  B_4^1 \ar[r]^>>>>{\partial}\ar[drrr]|{\beta} & B_3^1  \ar[drr]|{\beta}\ar[r]^{\partial} & B_2^1 \ar[dr]|{\beta}\ar[r]^{\partial}  & B_1^1\ar@/_0.5pc/[d]_{s} \ar@/^0.5pc/[d]^{t}
 \\
                    &&& B_0.} $$

\subsubsection{Homotopy of crossed complex maps}
There are several equivalent ways of defining the notion of homotopy between morphisms of crossed complexes. Given a crossed complex, $\A$, we can form the tensor product, $\I\otimes \A$.   This gives a model for a `cylinder' on $\A$, so then a homotopy between two maps, $f$ and $f'\colon \A \to \B$, will be a morphism, $h\colon \I\otimes \A\to \B$, satisfying some fairly obvious conditions as in  \cite{brown_higgins_sivera}; see Theorem \ref{homsandotimes}, below.
Alternatively, we can use the internal `hom', $\CRS(-,-)$, in $\Crs$, which  we will meet in \S\ref{sec:internalhom}, and form $\A^\I=\CRS(\I,\A)$. This  leads to a homotopy being seen as a morphism from $\A$ to $\B^\I$.

There is another definition of homotopy, which is the crossed complex analogue of the notion of homotopy of morphisms of chain complexes often given in books on Homological Algebra. This does not need additional constructions to make it work and, in fact, is needed to make sense of the construction, $\CRS(-,-)$, so we start with this. The idea is that we start with both a morphism $f\colon \A \to \B$ and an  `$f$-homotopy', i.e. a homotopy that ends at $f$, and  then obtain the `other end' of the homotopy from that input; see  \cite[\S 7.3]{brown_hha}.

\begin{Definition}\label{def:f-homotopy} 
Consider a crossed complex map, $f=(f_n)_{n \in \mathbb{Z}^+_0}\colon \A \to \B$. An $f$-homotopy, $H=(h_n)_{n \in \mathbb{Z}^+_0}$, or 1-fold $f$-homotopy, is given, in low dimensions, by 
\begin{itemize}[leftmargin=10mm]

 \item a set map, $h_0\colon A_0 \to B_1^1$, such that $t\circ h_0=f_0$,
 \item a set map, $h_1\colon A_1^1 \to B_2^1$, such that $\beta(h_1(g))=f_0(t(g))$.
This map, $h_1$, is furthermore to be a type of derivation, so it is to be such that, if the morphisms 
 $g$ and $g'$, in $A_1$ can be composed, then $$h_1(gg')=\big(h_1(g) \trl f_1(g') \big)\, h_1(g'), $$  We can visualise the conventions for the source and target of $h_0$ and $h_1$ as:
 \begin{equation}\label{conv:h}\turnradius={1pc}
 \vcenter{\xymatrix{ & x \ar[r]^g &y}} \longmapsto \vcenter{\xymatrix@C=30pt{  f_0(x) \ar[r]^{f_1(g)} &f_0(y) \ar@(l,d)[]_<<<<{\,\,\,\,\,\,\,\,\,\,\partial(h_1(g))}\\
 s(h_0(x))\ar[u]^{h_0(x)} &s(h_0(y)).\ar@<-0.3em>[u]_{h_0(y)}}
 }
 \end{equation}

  \item[]\hspace*{-.7cm} In higher dimensions, we have,
 \item if $n\ge 2,$ a groupoid map, $h_n\colon A_n \to B_{n+1}$, which, on objects, restricts to $f_0$, such that,  given $x,y \in A_0$, if $a \in A_n(x)$ and $g\in A_1(x,y)$, then
  $$h_n(a \trl g)=h_n(a) \trl f_1(g).$$
\end{itemize}
We denote the set of $1$-fold $f$-homotopies by $\CRS_1(\A,\B,f)$.
\end{Definition}

In the setting of this  definition, as in \cite[Exercise 7.1.39]{brown_higgins_sivera},  given $$f=(f_n)_{n \in \mathbb{Z}^+_0}\colon \A \to \B\textrm{ and } H=(h_n)_{n \in \mathbb{Z}^+_0}\in \CRS_1(\A,\B,f),$$  it then follows that  we have a crossed complex map, $f'=(f_n')_{n \in \mathbb{Z}^+_0}\colon \A \to \B$, defined by the equation below  (it may be useful to refer to the diagram in \eqref{conv:h}),
\begin{align*}
 f_0'(x)&=s(h_0(x)), \textrm{ if } x \in A_0;\\
 f_1'\big ( x \ra{g} y)&=h_0(x)\,f_1(g) \, \partial\big (h_1(g)\big)  \, h_0(y)^{-1}, \textrm{ if } ( x\ra{g} y ) \in A_1^1; \\
 f_n'(a)&=\big ( f_n(a)\,\,h_{n-1}(\d(a))\,\, \d (h_n(a))\big) \trl h_1(\beta(a))^{-1},  \textrm{ if } n \ge 2 \textrm{ and } a \in A_n^1.
\end{align*}

We denote $s(H,f)=f'$ and $t(H,f)=f$, and frequently write $f'\ra{(H,f)} f$, or $(H,f)\colon f'\to f$. We say that \emph{$(H,f)$ is a crossed complex homotopy from $f'$ to $f$}.

 We put $\CRS_0(\A,\B)=\Crs(\A,\B)$, the set of crossed complex maps from $\A$ to $\B$. As the notation indicates in its use of $s$ and $t$, we have a groupoid $$\CRS_1(\A,\B)=\big (s,t\colon \CRS_1(\A,\B)^1 \to  \CRS_0(\A,\B)\big),$$
  whose objects are the maps, $f\colon \A \to \B$, and the morphisms from $f'$ to $f$ are the $f$-homotopies such that $(H,f)\colon f'\to f$. The composition of $(H',f')\colon f''\to f'$ and $(H,f)\colon f'\to f$, denoted ${(J,f)}\colon f''\to f$, with $J=(j_n)_{ n \in \mathbb{Z}_0^+}$, is such that, if  $H=(h_n)_{n \in \mathbb{Z}^+_0}$ and  $H'=(h_n')_{n \in \mathbb{Z}^+_0}$, then
 \begin{align*}
 j_0(x)&=h'_0(x) \, h_0(x), \textrm{ if } x \in A_0,\\ \intertext{and}
j_n(x \ra{g} y)&=h_n(x \ra{g} y)\,\,\big ( h_n'(x \ra{g} y)\trl h_0(y)\big), \textrm{ if } n \ge 1, \textrm{ and } (x \ra{g} y) \in A^1_n.
\end{align*}

\subsubsection{The internal hom $\CRS(-,-)$}\label{sec:internalhom}
 The groupoid, $\CRS_1(\A,\B)$, can be `extended' to  a crossed complex, $\CRS(\A,\B)$, denoted,
 $$ \xymatrix@C=22pt{  \CRS(\A,\B)=\cdots   \ar[r]^>>>>{\delta}\ar[drrr]|{\beta} & \CRS_3(\A,\B)^1  \ar[drr]|{\beta}\ar[r]^{\delta} & \CRS_2(\A,\B)^1 \ar[dr]|{\beta}\ar[r]^{\delta}  & \CRS_1(\A,\B)^1,\ar@/_0.5pc/[d]_{s} \ar@/^0.5pc/[d]^{t}
 \\
                    &&& \CRS_0(\A,\B),}
$$
by considering $k$-fold homotopies between crossed complex maps for each $k\in\mathbb{Z}^+$.
This construction is explicitly given in  both \cite{brown_higgins_tensor} and  \cite[\S 9.3.i]{brown_higgins_sivera}.
\begin{Definition}
 Consider a crossed complex map, $f=(f_n)_{n \in \mathbb{Z}_0^+}\colon \A \to \B$, and let $k\ge 2$. 
 A $k$-fold $f$-homotopy, $H^k=(h_n^k)_{n \in \mathbb{Z}^+_0}=(h_0^k,h_1^k,h_2^k,\dots)$, is given by:
 \begin{itemize}[leftmargin=1cm]
  \item the choice of an element $h_0^k(x) \in B_k(f_0(x))$ for each $x \in A_0$,
  \item given $(x \ra{g} y)\in A_1^1$, the choice of an element, $h_1^k(x \ra{g} y) \in B_{k+1}(f_0(y)),$  to be such that,
  if $g$ and $g'$ can be composed in $A_1$, then $$h_1^k(gg')=\big(h_1^k(g) \trl f_1(g') \big)\,\, h_1^k(g');$$
  \item given $n\ge 2$, and $x \in A_0$,    a function, $h_n^k\colon A_n^1 \to B_{n+k}^1$, satisfying $$\beta(h_n^k(a))=f_0(\beta(a)), \textrm{ for all } a \in A_n.$$ 
  This mapping, $h_n^k\colon A_n^1 \to A_{n+k}^1$, is to be such that, given any $x \in A_0$, the restriction of $h_n^k$ to $A_n(x)$ is a group homomorphism, $A_n(x) \to B_{n+k}(f_0(x))$, and further, if $x,y \in A_0$, $a \in A_n(x)$ and $(x \ra{g} y) \in A_1^1$, then
  $$ h_n^k\big ( a \trl (x \ra{g} y)\big)= h_n^k\big (a) \trl (f_0(x) \ra{f_1(g)} f_0(y)).$$
 \end{itemize}
We let $\CRS_k(\A,\B,f)$ denote the set of all $k$-fold $f$-homotopies.
 \end{Definition}
 
Let $f\colon \A \to B$ be a crossed complex map. Suppose $k\ge 2$. By using the obvious point-wise product of $k$-fold $f$-homotopies, as in  \cite[Definition 9.3.5]{brown_higgins_sivera}, we have that the set $\CRS_k(\A,\B,f)$ has a group structure, and  that is abelian if $k\ge 3$.

Given $k\ge 2$, we have a totally disconnected groupoid,
$$\CRS_k(\A,\B):=\Big(\beta\colon \bigsqcup_{f\colon \A \to \B }  \CRS_k(\A,\B,f) \to \CRS_0(\A,\B)\Big),$$
with object set, $ \CRS_0(\A,\B)=\Crs(\A,\B)$, the set of crossed complex maps, $f\colon \A \to \B$, and with  the obvious map, $\beta$, that identifies the component of the disjoint union.

\begin{Lemma}Let $f\colon \A \to \B$ be a crossed complex map. Let $H^2=(h_0^2,h_1^2,\dots)$ be a $2$-fold $f$-homotopy, then $\delta(H^2)= (\delta(h_0^2),\delta(h_1^2),\dots)$, defined by:
\begin{itemize}
 \item $\delta(h_0^2)(x):=\d\big ( h_0^1(x)\big)$, for each $x \in A_0$;
 \item $\delta(h_1^2)(x \ra{g} y):=\big (h_0^2(x) \big)^{-1}\trl \big (f_0(x) \ra{f_1(g)} \,\, f_0(y)\big )\, h_0^2(y)\, \d (h_1^2(x\ra{g} y))$, where $(x\ra{g} y)\in A_1^1$,
 \item[]\hspace*{-10mm} and 
 \item given $n\ge 2$ and $a \in A_n^1$,
 $$\delta(h_n^2)(a)=\d \big( h_{n+1}^2 ({a})\big)\,\, h_n^2(\d (a))^{( (-1)^n)},$$
\end{itemize}
is an $f$-homotopy, ${(\delta(H),f)}\colon f\to f$.

We have a groupoid action of $\CRS_1(\A,\B)$ on $\CRS_2(\A,\B)$, by automorphisms, where  given crossed complex maps $f,f'\colon \A \to \B$, an $f'$-homotopy $J=(j_n)_{n \in \mathbb{Z}^+_0}$, connecting $f$ to $f'$ and a 2-fold $f$-homotopy $H^2=(h_0^2,h_1^2,\dots)$, then
$$ \big( f\ra{H^2} f\big )\trl (f \ra{(J,f')} f')= \big(f' \ra{H^2\trl J} f'\big),$$
is the 2-fold $f'$-homotopy such that:
\begin{align*}
(h^2\trl J)_k \big (x \ra{g} y\big )&=h^2_k\big (x \ra{g} y\big) \trl j_0( y).
\end{align*}
\end{Lemma}
\begin{proof}
This is proved by explicit calculations.
 \end{proof}

 Similarly, for  $n\ge 3$, we  have groupoid maps, $\delta \colon \CRS_n(\A,\B) \to \CRS_{n-1}(\A,\B)$, which, again, restrict to the identity on the set of objects, together with actions of the groupoid, $\CRS_1(\A,\B)$, on all of the totally disconnected  groupoids, $\CRS_n(\A,\B)$, for $n \ge 2$. This gives rise to a crossed complex, $\CRS(\A,\B)$, the  \emph{internal hom}  in the category of crossed complexes. For
details, see \cite[\S 9.3.i]{brown_higgins_sivera}.

\begin{Remark}
Let $\A$ and $\B$ be crossed complexes. Following on from the $\CRS(\A,\B)$ construction, we have the  following groupoid, using Definition \ref{def:pi1andpi0crs},
$$\pi_1\big(\CRS(\A,\B)\big)=\big(\CRS(\A,\B),\Crs(\A,\B)\big).$$
The set of objects of  $\pi_1\big(\CRS(\A,\B)\big)$ is, thus, the set, $\CRS_0(\A,\B)=\Crs(\A,\B)$, of crossed complex maps, $f$ from $\A$ to $\B$, and given $f,g\colon \A \to \B$, the set  of arrows from $f$ to $g$ is given by all equivalence classes of homotopies, $[(H,g)]\colon f\to g$, connecting $f$ and $g$, with homotopies considered up to   2-fold homotopy.

This groupoid will play a key role in the description of once-extended TQFTs derived from crossed complexes.
\end{Remark}

\subsubsection{Tensor product and homotopies of crossed complexes}\label{tensorHom}
A crucial property of the crossed complexes, $\CRS(\A,\B)$, where $\A$ and $\B$ are crossed complexes, is that they vary functorially in both positions, so we have a functor, $$\CRS(-,-):\Crs^\op \times \Crs \to \Crs,$$ sending $(\A,\B)$ to $\CRS(\A,\B)$; see  \cite{brown_higgins_tensor}.
This functor,  $\CRS(-,-)$, acts as an `internal hom'. Hence $\CRS(\A,\B)$ behaves like the “object of morphisms” from $\A$ to $\B$, so $\CRS(-,-)$ is analogous to the mapping space functor in $\TOP$ defined on page \pageref{mapping space}.
\begin{Notation}\label{fstar} If $\B$ is a crossed complex,  we have a  functor, $$\CRS(-,\B)\colon \Crs^{\op} \to \Crs.$$ It sends a crossed complex, $\A$, to $\CRS(\A,\B)$, and  a crossed complex map, $f\colon \A' \to \A$, to the crossed complex map, $f^*\colon \CRS(\A,\B) \to \CRS(\A',\B)$, such that:
\begin{enumerate}[leftmargin=1cm]
 \item each crossed complex map, $\phi\colon \A \to \B$, is sent to the composite, $\phi\circ f\colon \A' \to \B;$
 \item given $k\in \mathbb{Z}^+$, a crossed complex map, $\phi\colon \A \to \B$, and a $k$-fold $\phi$-homotopy, $h^k=(h_0^k,h_1^k,h_2^k,\dots)$, then $f^*(h^k):=(h_0^k\circ f,h_1^k\circ f,h_2^k\circ f,\dots)$,
\end{enumerate}
which is a $k$-fold $(\phi\circ f)$-homotopy.
This corresponds to `pre-composition with $f$'.  
\end{Notation}

Given crossed complexes, $\A$ and $\A'$, we can also form their \emph{tensor product}, $\A \otimes \A'$; again, for details, see \cite[\S 9.3.iii]{brown_higgins_sivera} and \cite[Definition 1.4.]{tonks:JPAA:2003}. We have a functor, $\Otimes_\CRS\colon \Crs \times \Crs \to \Crs$, sending $(\A,\A')$ to $\A\otimes \A'$, and an exponential law, $\Crs(\A\otimes \A', \B) \cong \Crs(\A, \CRS(\A',\B))$, that holds naturally in $\A$ and $\B$, showing that the functor `tensor product with $\A'$', i.e.,  $-\otimes \A'$, is left adjoint to the functor, $\CRS(\A',-)$, derived from the internal hom.  This gives $\Crs$ the structure of a monoidal closed category, \cite[Theorem 9.3.17]{brown_higgins_sivera}. In fact, the tensor product is symmetric, so $\Crs$ with the above tensor is a symmetric monoidal closed category, see  \cite[Theorem 9.3.16]{brown_higgins_sivera}.

Given a crossed complex, $\A$, we have  morphisms, which look like the inclusions of the ends of a cylinder, and will here be denoted  $i_0,i_1\colon \A \to \Pi(I_\sk) \otimes \A$; see \cite[page 203]{tonks:JPAA:2003}. (We note that, as this notation, $i_0$, etc.,  is overcharged, occurring in several contexts, often with different meanings, we will sometimes replace $i_0$ by $e_0(\A)$, or $\iota_0^\A$, etc., depending on the other use of symbols in the setting.)  If $k\ge 2$, we also have a canonical inclusion, $i\colon \A \to \Pi(D^k_\sk) \otimes \A $, as $ \Pi(D^k_\sk)$ is a reduced crossed complex.
Moreover, there is a morphism from $\Pi(I_\sk) \otimes \A$ to $A$, which is a partial inverse to the `end inclusion' morphisms. This means that $\Pi(I_\sk) \otimes \A$ behaves exactly like a cylinder on $\A$, and can be used to define a notion of homotopy between morphisms in $\Crs$, which, thankfully, coincides with the one that we introduced earlier, where we used the abbreviation $\I$ for $\Pi(I_\sk)$.  All this is very thoroughly discussed in \cite[\S 9.3.i]{brown_higgins_sivera}, and we note:
\begin{Theorem}\label{homsandotimes}Let $f\colon \A \to \B$ be a crossed complex map. We have
\begin{itemize}[leftmargin=1cm]
\item There is a canonical correspondence between homotopies, ${(f,H)}\colon f'\to f$, and commutative diagrams in $\Crs$ of form,
 $$\xymatrix@R=18pt{  \A \ar[dr]_{f'} \ar[r]^<<<<<{i_0} &\Pi(I_\sk) \otimes \A \ar[d]^{\hat{H}} & \A. \ar[l]_>>>>>>{i_1}  \ar[dl]^{f} \\
 & \B &
 }
 $$
\item If $k\ge 2$, we have  a canonical correspondence between $k$-fold $f$-homotopies, $H^k$, and commutative diagrams in $\Crs$ of form
 $$\xymatrix@R=18pt{  \A \ar[dr]_{f} \ar[r]^<<<<{i} &\Pi(D^k_\sk) \otimes \A \ar[d]^{\hat{H}^k}  \\
 & \B.
 }
 $$
 \end{itemize}
 
\end{Theorem}
%
%

For the geometric interpretation of the tensor product, the following is crucial. See \cite[Theorem 9.8.1]{brown_higgins_sivera}, or  page 92 of \cite{Baues_4D} for the reduced case.
\begin{Theorem}\label{thm:tensor-CW}
 Let $X$ and $Y$ be CW-complexes. Give $X\times Y$ the usual structure of a CW-complex. We have a natural isomorphism of crossed complexes, $$\Pi\big((X\times Y)_\sk\big)\xrightarrow{\cong}\Pi(X_\sk) \otimes \Pi(Y_\sk).$$
\end{Theorem}
\subsubsection{Homotopies and totally free crossed complexes}\label{sec:hom_of_freeXcomp}
To be able to work fairly simply with the above notions of homotopy between crossed complex maps, we will need to be able to construct homotopies in ways analogous to the `induction up the skeleton' methods used in many topological contexts.

The detailed  result in Lemma \ref{lem:derivation_basis},  below, is crucial for what follows and is one such statement. It is given a direct proof in \cite[\S 2.2.6]{MartinsCW}, for the case when the CW-complex, $X$, has a single $0$-cell, and in \cite[\S 2.24.1]{Martins_Porter} for general CW-complexes, but,  there, for the particular case in which $\A$ is 2-reduced.  Related results are also  in  \cite[Corollary 9.6.6]{brown_higgins_sivera}, \cite[Proposition 7.3 I]{Brown_Sivera} and \cite[Chapter III, \S 4]{Baues_4D}.

Fix a reduced crossed complex $\A$ as in \eqref{eq:Ared},  and a special CW-complex  $X$. (These restrictions are there only to simplify the exposition, as for Lemma \ref{lem:maps_on_basis}.) Recall that $C(X,i)$ is the set of $i$-cells of $X$.

\begin{Lemma}\label{lem:derivation_basis} Let $f\colon \Pi(X_\sk)\to \A$ be a crossed complex map.
Let $k\in \mathbb{Z}^+$, then $k$-fold $f$-homotopies are uniquely specified by their value on the elements of $\Pi(X_\sk)$, defined from the cells of $X$.

Explicitly this means that:
\begin{itemize}[leftmargin=1cm]
 \item 
we have a bijection between $1$-fold $f$-homotopies and sequences of maps,
$$(m_i^1\colon  C(X,i) \to A_{i+1})_{i \in \mathbb{Z}^+_0},$$
(of sets) and note that  there are no further compatibility conditions between the maps, $m_i^1\colon C(X,i) \to A_{i+1}$,  and the boundary maps of $\Pi(X_\sk)$ and $\A$.

We will denote this bijection by
\begin{equation}\label{eq:extend}
(m_i^1)_{i \in \mathbb{Z}^+_0} \mapsto \Ext^1_X\big((m_i^1)_{i \in \mathbb{Z}^+_0},f\big)\in \CRS_1(\Pi(X_\sk),\A,f).
\end{equation}
 \item If $k\ge 2$, we have a one-to-one correspondence, between $k$-fold $f$-homotopies and sequences of maps, of sets,
$$(m_i^k\colon  C(X,i) \to A_{i+k})_{i \in \mathbb{Z}^+_0}.$$
We denote this bijection by: 
$$(m_i^k)_{i \in \mathbb{Z}^+_0} \mapsto \Ext^k_X\big((m_i^k)_{i \in \mathbb{Z}^+_0},f\big)\in \CRS_k(\Pi(X_\sk),\A,f).$$
\end{itemize}
\end{Lemma}

\begin{proof}
This follows  from Remark \ref{rem:explicit_freeness}, whose nomenclature we use.

For $k=1$,  given $(m_i^1\colon C(X,i) \to A_{i+1})_{i \in \mathbb{Z}^+_0}$, then $\Ext^1_X\big((m_i^1)_{i \in \mathbb{Z}^+_0},f\big)$
 is the unique 1-fold $f$-homotopy that takes the value, $m_i^1$, on the set of $i$-cells of $X$. The existence and uniqueness of $\Ext^1_X\big((m_i^1)_{i \in \mathbb{Z}^+_0},f\big)$ follows from elementary techniques, since $\Pi(X_\sk)$ is free on the set of cells of $X$, in the sense explained in  Remark \ref{rem:explicit_freeness}. (More details can be found in  \cite[\S 2.2.6]{MartinsCW}, when $X$  is reduced.)

A similar argument  is valid when $k\ge 2$. Given a sequence of maps, $$(m_i^k \colon  C(X,i) \to A_{i+k})_{i \in \mathbb{Z}^+_0},$$ then $\Ext^k_X\big((m_i^k)_{i \in \mathbb{Z}^+_0},f\big)$ is the unique $k$-fold $f$-homotopy that takes the value, $m_i^k$, on the set of $i$-cells of $X$.
\end{proof}
 A perhaps more conceptual proof of Theorem \ref{lem:derivation_basis}  follows by combining Theorem \ref{homsandotimes} and Lemma \ref{lem:maps_on_basis}, using the fact that the crossed complex, $\Pi\big((X\times I)_\sk\big)\cong \Pi(X_\sk)\otimes\Pi(I_\sk)$, is free. The same argument works for $k\ge 2$.

Theorem \ref{lem:derivation_basis}  leads to the following definition. Here $k$ is a positive integer.
\begin{Definition}
 A \emph{$k$-fold homotopy $(X,\A)$-sequence} is a sequence of maps, $$(m_i^k\colon C(X,i) \to A_{i+k} )_{i \in \mathbb{Z}^+_0}.$$
\end{Definition}
\noindent By  Theorem \ref{lem:derivation_basis},  given a crossed complex map, $f\colon \Pi(X_\sk) \to \A$, there is a bijection \eqref{eq:extend} between $k$-fold $f$-homotopies and $k$-fold  homotopy $(X,\A)$-sequences.

We will need a generalisation of Theorem   \ref{lem:derivation_basis} for when $(X,Y)$ is a CW-pair. So we let $Y$ be a subcomplex of $X$,  with  $\iota\colon Y \to X$ denoting the inclusion map.
\begin{Notation}
Let $f\colon \Pi(X_\sk) \to \A$ be a crossed complex map and let   $k\in \mathbb{Z}^+$.
 \begin{itemize}[leftmargin=1cm]
  \item Given a $k$-fold $f$-homotopy, $h^k=(h_j^k)_{j\in \mathbb{Z}},$ we define the $k$-fold homotopy $(Y,\A)$-sequence, denoted $$\Res^k_Y(h^k)=(m^k_j \colon  C(Y,j) \to A_{j+k})_{j \in \mathbb{Z}^+_0},$$
as the restriction of $h^k$ to the elements of $\Pi(X_\sk)$ given by the $j$-cells of $Y$.
  
  \item Given a $k$-fold homotopy $(Y,\A)$-sequence, $n^k=(n_j^k\colon C(Y,j) \to A_{j+k}^1)_{j \in \mathbb{Z}^+_0}$, we define the $k$-fold  homotopy $(X,\A)$-sequence, denoted  $$\Exp^k(n^k,Y,X)=(m_j^k\colon C(X,j) \to A_{j+k}  )_{j \in \mathbb{Z}^+_0},$$
  to be such that, given $j\in \mathbb{Z}^+_0$, then $m_j^k$ coincides with $n_j^k$ over $C(Y,j)\subseteq C(X,j)$, and otherwise $m_j^k$ takes as values the identity element of $A_{j+k}$.
  
 \end{itemize}

\end{Notation}

We use the above notation in the following, where $(X,Y)$ is a CW-pair, with $X$ special,  $\A$ is reduced, and $f\colon \Pi(X_\sk) \to \A$ is a crossed complex map.
\begin{Lemma} \label{lem:fib_direct}
Let $k$ be a positive integer.
\begin{enumerate}[leftmargin=1cm]
 \item Given a $k$-fold $f$-homotopy, $h^k$, we have $\Ext^k_X\big(\Res^k_X(h^k),f\big)=h^k$.
 
 \item \label{Fib-in-practice} 
Let  $\Pi(\iota)\colon \Pi(Y_\sk) \to\Pi(X_\sk)$ be the induced map, which  induces a crossed complex map, going in the other direction, via `restriction', \[\Pi(\iota)^*\colon \CRS\big(\Pi(X_\sk),\A\big) \to \CRS\big(\Pi(Y_\sk),\A\big),\]
(see Notation \ref{fstar}).
 Given a $k$-fold $f\circ \Pi(\iota)$-homotopy, $h^k$, 
 we have that:
 $$\Pi(\iota)^*\Big( 
 \Ext^k\Big (
 \Exp^k \big(
 \Res^k_Y(h^k )
 ,Y,X\big),
 f 
 \Big)
 \Big )=h^k.
 $$
 
\end{enumerate}

\end{Lemma}
\begin{proof}
Follows from the freeness of $\Pi(X_\sk)$ and $\Pi(Y_\sk)$, and  Lemma \ref{lem:derivation_basis}.
\end{proof}

\subsection{The classifying space of  a crossed complex}\label{sec:class_space} Source material for classifying spaces of crossed complexes can be found in \cite{brown_higgins_classifying}, \cite[\S 9.10]{brown_higgins_sivera},  also in \cite{BGPT}  and   \cite{Martins_Porter}, as well as in various other of the sources cited earlier.


\subsubsection{Nerves and classifying spaces of  crossed complexes}\label{sec:simpclass-space}

The classifying space functor, $B\colon \Crs\to \CGWH$, is defined as the composite of the nerve functor, $\N\colon \Crs\to \Simp$, and the geometric realisation functor from $\Simp$ to $\CGWH$.

The main tool to define $\N$ is that we have a functor, $\Pi\circ  \Ds\colon \Delta \to \Crs$, sending $[n]\in \Delta$ to $\Pi(\Ds(n))=\Pi(|\Ds(n)|_\sk)$,
which thus gives a cosimplicial crossed complex.
\begin{Definition}
The nerve of a crossed complex, $\A$, is the simplicial set,  $\N(\A)\colon \Delta \to \Sets$, obtained as the composite below, where $h^\A=\hom_\Crs(-,\A)$,
$$ \Delta \ra{\Pi\circ  \Ds} \Crs \ra{h^\A} \Sets.$$  Hence the set of $n$-simplices is $\N(\A)_n=\Crs(\Pi(\Ds(n)),\A).$ The maps between the different dimensions are induced from  $\Ds\colon \Delta \to \Crs$.

This also defines the  \emph{nerve functor} $\N\colon \Crs \to \Simp$.\end{Definition}
For us, one of the most useful facts about the nerve functor is the following.
\begin{Proposition}[Brown-Higgins]\label{brown_higgins_adjunction} The nerve functor, $\N\colon \Crs \to \Simp$, is right adjoint to the fundamental crossed complex functor, $\Pi\colon \Simp \to \Crs$. \end{Proposition}
\begin{proof} See \cite[Theorem 2.4]{brown_higgins_classifying}, or combine Theorem \ref{Th:fundCRS-coend} (which is also in \textit{loc. cit.}) with the Nerve-Realisation Paradigm in \cite[Proposition 3.2.2]{loregian_2021}. \end{proof}
Given a simplicial set, $S$, and crossed complex, $\A$, we, thus, have a  bijection, $$\phi_{S}^{\A}\colon \Crs\big(\Pi(|S|_\sk),\A\big) \to \Simp\big(S,\N (\A)\big),$$
natural in both $S$ and $\A$. 

A crucial fact that  underpins the use of the nerve is the following; see \cite{brown_higgins_classifying}.
\begin{Proposition}Let $\A$ be a crossed complex, then $\N (\A)$ is Kan.\end{Proposition}

As stated above, the classifying space construction that we will be using is obtained from the nerve by taking geometric realisation.
\begin{Definition}[Classifying space of a crossed complex \cite{brown_higgins_classifying}] The classifying space, $B_\A$, of a crossed complex, $\A$, is defined as the geometric realisation, $|\N(\A)|$.
\end{Definition}
 \begin{Notation}\label{rem:tilde}
Let $\A=(A_n)_{n\in \mathbb{Z}^+_0}$ be a crossed complex.  By construction, $A_0\cong \N(\A)_0$, and so each object $a\in A_0$ of $\A$ gives rise to a $0$-simplex of $\N (\A)$, therefore to a vertex of the CW-complex $|\N(\A)|$, which will be denoted $\tilde{a} \in B_\A$.  \end{Notation}

Given simplicial sets, $K$ and $L$, their simplicial mapping space, i.e., the function complex of \cite[\S 6]{MaySimplicial}, will be  denoted $\SIMP(K,L)$.  We note that $\Simp$ becomes  a cartesian closed category with this function space construction. In particular, we have   $\SIMP(X,Y)_0=\Simp(X,Y)$, the set of simplicial set maps from  $X$ to $Y$. 

As is well known, if $K$ and $L$ are simplicial sets, with $L$ a Kan complex, then we have a weak homotopy equivalence, $|\SIMP(K,L)| \to \TOP(|K|,|L|)=|L|^{|K|}$.
Explicitly, this weak homotopy equivalence sends the equivalence class of
$$(f\colon K \times \underline{\Delta}(n) \to L,s)\in \int^{n \in \Delta} \SIMP(K,L)_n \times |\underline{\Delta}(n)|$$
to the function $|K|\to |L|$, such that $k\mapsto |f|(k,s)$. An explicit proof that this is a weak homotopy equivalence is in \cite[page 131]{Martins_Porter}. Given a simplicial map, $f\colon K \to L$, the weak homotopy equivalence sends the corresponding vertex of $\SIMP(K,L)$ to the geometric realisation, $|f|\colon |K| \to |L|$.

The technical results collected up, for convenience,  in the next theorem are due to Brown--Higgins, \cite{brown_higgins_classifying}, and Tonks, \cite{tonks:1993,tonks:JPAA:2003}. They are discussed in the cubical, as opposed to the simplicial, setting by Brown--Higgins--Sivera in \cite{brown_higgins_sivera}, and, to some extent, in a simplicial setting in \cite{BGPT}.

\begin{Theorem}[Brown--Higgins; Brown--Higgins--Sivera; Tonks]\label{MainCrs}
 As usual, let $\A=(A_n)_{n\in \mathbb{Z}^+_0}$ be a crossed complex, and take $S$ to be a simplicial set, then:
\begin{enumerate}[leftmargin=1cm]
 \item\label{l1} there is  an isomorphism of groupoids, $\pi_1(\A,A_0)\cong \pi_1(B_\A,\widetilde{A_0})$, where $\widetilde{A_0}=\{\tilde{a}\mid a \in A_0\}$, using  Notation \ref{rem:tilde}, which is natural in $\A$, \item[]\hspace*{-5mm}and hence,
 \item\label{l2} there is a  natural bijection, $\pi_0(\A) \cong \pi_0(B_\A)$.
 \item\label{l3} Let $a\in A_0$, and let $n$ be a positive integer, we have a natural isomorphism, $\pi_n(\A,a) \cong \pi_n(B_\A,\tilde{a})$, preserving the  actions of $\pi_1(\A,A_0)$ and $\pi_1(B_\A,\tilde{A_0})$.
   \item\label{l4} There is a weak homotopy equivalence of simplicial sets,
    $$\eta_S^{\A}\colon  \N \big (\CRS(\Pi(|S|_\sk),\A)\big)  \to \SIMP(S, \N (\A)),$$ which, at the level of $0$-simplices, coincides with the bijection, $$\phi_{S}^{\A}\colon \Crs(\Pi(|S|_\sk),\A) \to \Simp\big(S,\N (\A)\big),$$
    given by the adjunction  $\xymatrix{\hskip-0.5mm\Pi\colon \Simp
\ar@/^0.3pc/[r] \ar@{}[r]|\bot \ar@{<-}@/_0.3pc/[r] & \Crs: \hskip-0.5mm{\N},}$  of Proposition \ref{brown_higgins_adjunction}.
    This weak homotopy equivalence is natural in $S$  and also in $\A$. (We note  that $\eta_S^{\A}$ is not simplicially natural in $S$, for which fact see \cite{tonks:JPAA:2003} and \cite{BGPT,BGPT2}.)
  \item\label{l5}\label{simp-maps-to-top-maps} There is a weak homotopy equivalence,
 $$\etab_S^\A\colon  \left |\N \big(\CRS(\Pi(|S|_\sk),\A)\big)\right|  \to \TOP(|S|, B_\A)).$$
 This weak homotopy equivalence is natural in both $S$ and $\A$.
\end{enumerate}
\end{Theorem}
\noindent We will use the results in the previous theorem without giving a proof here, rather we note \eqref{l1} \eqref{l2} and \eqref{l3} form parts of  \cite[Proposition 2.6]{brown_higgins_classifying}; for \eqref{l4}, see \cite[Theorem A]{brown_higgins_classifying} and  \cite[Proposition 3.1.]{BGPT}. For \eqref{l5}, we refer again to  \cite[Theorem A]{brown_higgins_classifying},  \cite[Proposition 3.1.]{BGPT} and \cite[Section 4]{BGPT2}, and then proceed by composing with the canonical weak homotopy equivalence,
  $|\SIMP(S, \N (\A))| \to \TOP(|S|,|\N (\A)|)$.

  We note that item \eqref{simp-maps-to-top-maps}, above, links the classifying space of the crossed complex mapping space (i.e. internal hom) with the topological mapping space, from the realisation of $S$ to the classifying space of $\A$. This is the starting point, in the setting with $\Bc=B_\A$, for computing the Quinn finite total homotopy TQFT, and its extended versions.

\subsubsection{Homotopy classification of maps to classifying spaces}\label{sec:class-maps} Let $S$ be a simplicial set and $\A$ a crossed complex.
Let $f\colon \Pi(|S|_\sk) \to \A$ be a  crossed complex map. The adjunction
 $\xymatrix{\hskip-0.5mm\Pi\colon \Simp
\ar@/^0.3pc/[r] \ar@{}[r]|\bot \ar@{<-}@/_0.3pc/[r] & \Crs: \hskip-0.5mm{\N},}$
gives a simplicial map, $\phi_S^\A(f)\colon S \to \N(\A)$. Its geometric realisation is a continuous map, $$|{\phi_S^\A(f)}|\colon |S| \to B_\A=|\N (\A)|,$$
and  then
$\overline{\eta}_S^\A(\tilde{f})= |{\phi_S^\A(f)}|$, where, following Notation \ref{rem:tilde},
 $\tilde{f}$ is the vertex of the classifying space, $|\N (\CRS(\Pi(|S|_\sk),\A))|$, corresponding to $f$.

The following is essentially as \cite[Theorem A]{brown_higgins_classifying}. We let
$$\Comb(S,\A):=\big\{|\phi_S^\A(f)|  :   f \in \Crs (\Pi(|S|_\sk),\A)\big\} \subseteq \TOP(|S|,B_\A).$$
\begin{Lemma}\label{lem:iso_groupoids}
We have an isomorphism of groupoids, natural in $S$ and $\A$,
 $$\T_{S}^\A\colon \pi_1\big(\CRS(\Pi(|S|_\sk),\A)\big)\to \pi_1\big (\TOP(|S|, B_\A), \Comb(S,\A)\big),
 $$
 and hence we have a bijection between homotopy classes, natural in $S$ and $\A$,
 $$T_{S}^\A\colon \pi_0\big(\CRS(\Pi(|S|_\sk),\A)\big)\to \pi_0\big (\TOP(|S|, B_\A)\big),
 $$
  where on the left we have homotopy classes of crossed complex maps, $\Pi(|S|_\sk) \to \A$, and on the right homotopy classes of continuous maps $|S| \to B_\A$.
\end{Lemma}
\begin{proof}
 This follows from the first point of Theorem \ref{MainCrs}, together with the weak homotopy equivalence,  $\etab_S^\A\colon  \big|\N \big(\CRS(\Pi(|S|_\sk),\A)\big)\big|  \to \TOP(|S|, B_\A))$, since $\etab_S^\A$ is injective on the set $\{\tilde{f} :  f \in\CRS_0(\Pi(|S|_\sk),\A)\}$. (Note that if a weak homotopy equivalence, $g\colon X \to Y$, between spaces, is injective on $X_0\subseteq X$, then $g$ induces an isomorphism of groupoids $\pi_1(X,X_0)\cong \pi_1(Y,g(X_0))$.)
\end{proof}

The previous lemma can also be formulated for CW-complexes. Let $X$ be a CW-complex, and consider the singular set $\Sing(X)$ of $X$. Lemma \ref{lem:iso_groupoids} gives an isomorphism of groupoids,
\begin{multline*}
\T_{\Sing(X)}^\A\colon \pi_1\big(\CRS(\Pi(|\Sing(X)|_\sk),\A\big)\\\to \pi_1\big (\TOP(|\Sing(X)|, B_\A), \Comb(\Sing(X),\A)\big),
\end{multline*}
which is natural with respect to continuous maps $X \to Y$, between CW-complexes. We have a homotopy equivalence,
$P_X\colon |\Sing(X)| \to X.$  Applying the cellular approximation theorem, to the $i$-skeletons of $X$ in increasing dimension, we can find a cellular approximation, $P'_X\colon |\Sing(X)| \to X$, of $P_X$, and a homotopy inverse, $Q_X\colon X \to |\Sing(X)|$, that is cellular. This can be constructed so that if $Y$ is a subcomplex of $X$, the diagrams below commute, where $\iota_Y^X\colon Y \to X$ is the inclusion,
$$\vcenter{\xymatrix@C=45pt{       |\Sing(Y)| \ar[r]^{|\Sing(\iota_{Y}^X)|} &  |\Sing(X)|\\
Y\ar@{<-}[u]^{P'_Y} \ar[r]_{\iota_{Y}^X} &  X  \ar@{<-}[u]_{P'_X} \\
              }}\qquad \textrm{ and } \qquad \vcenter{\xymatrix@C=45pt{ Y\ar[d]_{Q_Y} \ar[r]^{\iota_{Y}^X} &  X  \ar[d]^{Q_X} \\
             |\Sing(Y)| \ar[r]_{|\Sing(\iota_{Y}^X)|} &  |\Sing(X)|.
              } }$$
The cellular maps $Q_X$ and $Q_Y$ are homotopy equivalences, so induce homotopy equivalences of crossed complexes by \cite[Proposition 3.3]{brown_higgins_classifying}/ \cite[Proposition 9.8.3.]{brown_higgins_sivera}.

Now consider the following subset of  $\TOP(X,B_\A)$, $$\overline{\Comb(X,\A)}:= \big \{ |\phi_{S(X)}^\A\big(f\circ \Pi\big(P'_X)\big) |\circ Q_X : f \in \Crs(\Pi(X_{\sk}) ,\A )\big \}.$$
 The previous discussion gives the following, which again is a minor tweaking of \cite[Theorem A]{brown_higgins_classifying}, for $X$ a CW-complex:
\begin{Corollary}\label{cor:CW-class}
There exists an isomorphism of groupoids
  $$\overline{\T}_{X}^\A\colon \pi_1\big(\CRS(\Pi(X_\sk),\A)\big)\to \pi_1\big (\TOP(X, B_\A), \overline{\Comb(X,\A)}\big).$$
  In particular we have a bijection, of homotopy classes of maps,
  $$\overline{T}_{X}^\A\colon \pi_0\big(\CRS(\Pi(X_\sk),\A)\big)\to \pi_0\big (\TOP(X, B_\A)\big),$$
 that is natural with respect to the inclusion of subcomplexes, of $X$.
\end{Corollary}

\subsection{Fibrations of crossed complexes and profunctors}\label{sec:all_crossed_complex_fib}
\subsubsection{Fibrations of crossed complexes}\label{sec:crossed_complex_fib}
We  recall the notion of fibrations of groupoids, which was originally given   in \cite{RB:fibrations:1970}, and  is discussed in  \cite[B.7]{brown_higgins_sivera}. We then turn to fibrations   of crossed complexes, as in  \cite[Definition 12.1.1]{brown_higgins_sivera}.
\begin{Definition}Let $G'=(s,t\colon G_1' \to G_0')$ and  $G=(s,t\colon G_1 \to G_0)$ be groupoids.
A map, $(f_1,f_0)\colon G'\to  G$, of groupoids, is said to be a \emph{fibration of groupoids} if, given any $x \in G_0$, and  $x'$ in $G_0'$ with $f_0(x')=x$, and any arrow in $G$ of  form $(y\ra{g} x) \in G_1$, so ending at $x$, there exists at least one arrow, $$(y'\ra{g'} x') \in G_1',$$ with $$f_1\big (y'\ra{g'} x' \big)=(y \ra{g} x).$$
\end{Definition}
This, then, is a `path lifting' or, more precisely, an `arrow lifting' condition.
\begin{Definition}\label{Def:xcompFib} Let $\A=(A_n)_{n \in \mathbb{Z}^+_0}$ and  $\B=(B_n)_{n \in \mathbb{Z}^+_0}$ be crossed complexes.
A map, $f=(f_n)_{n\in \mathbb{Z}^+_0}\colon \A \to \B$, of crossed complexes, is called a \emph{fibration of crossed complexes} if
\begin{enumerate}[leftmargin=10mm]
\item $f_1\colon A_1 \to B_1$ is a fibration of groupoids,
\item[]\hspace*{-5mm}and
\item given any integer $n\ge 2$, and any $x \in A_0$, the group homomorphism, $$A_n(x,x) \to B_n(f_0(x),f_0(x)),$$ induced by $f_n\colon A_n \to B_n$, is surjective.  
\end{enumerate}
\end{Definition}

An important link with Kan fibrations of simplicial sets is given in  the following result; see \cite[Proposition 6.2]{brown_higgins_classifying}. For a proof in the cubical, as opposed to simplicial set, setting, see \cite[Proposition 12.1.13]{brown_higgins_sivera}. Let $p\colon \A \to \B$ be a crossed complex map.
\begin{Lemma}  The following are equivalent:
\begin{itemize}[leftmargin=1cm]
 \item $p\colon \A \to \B$ is a fibration,
 \item the induced map on nerves, $\N(p)\colon \N (\A) \to \N (\B)$, is a Kan fibration.
\end{itemize}
\end{Lemma}

The next result will play a major role later. It appears to be new, however not unexpected. Recall Notation \ref{fstar} for the induced map in this setting.

Let $\A$ be a crossed complex. Let $X$ be a CW-complex, and $Y$ be a subcomplex of $X$ with $\iota\colon Y \to X$ denoting the inclusion, which induces a crossed complex map, $$\Pi(\iota)\colon \Pi(Y) \to \Pi(X).$$
\begin{Lemma}\label{lem:incs_to_fibs}

The induced crossed complex map between internal homs, $$\Pi(\iota)^*\colon \CRS(\Pi(X),\A) \to\CRS(\Pi(Y),\A),$$ is a fibration of crossed complexes.
\end{Lemma}
\begin{proof}
 When $\A$  has a unique object, which is our main case of interest, a  proof follows directly from the second point of Lemma \ref{lem:fib_direct}. This argument can be easily adapted for the case when $\A$  has more than one object. 
\end{proof}
\begin{Remark}
A model category theoretical proof of Lemma \ref{lem:incs_to_fibs} follows from the fact that the category, $\Crs$, of crossed complexes is a monoidal model category, which was observed by Sauvageot, in  \cite{sauvageot:stabilisation:thesis:2003}, and the fact that the crossed complex map,  $\Pi(Y) \to \Pi(X)$, induced by the inclusion, is a cofibration\footnote{See \cite[section 12.1]{brown_higgins_sivera} and \cite{BrownGolasinski}  for a more detailed discussion of cofibrations and trivial cofibrations of crossed complexes. We just need that inclusions, so cofibrations,  of CW-complexes are sent by $\Pi$ to cofibrations of crossed complexes and similarly for trivial cofibrations.} in that structure; see \cite[Proposition 12.1.4. and Example 7.3.19]{brown_higgins_sivera}.
\end{Remark}
 
By \cite[Proposition 6.2, ii]{brown_higgins_classifying},  given $p\colon \A \to \B$,  a fibration of crossed complexes,  $p$ has the right-lifting property with respect to the map, $\Pi(\{0\}) \to \Pi(I)$, induced by inclusion, where, as usual,  $I=[0,1]$. 

Let $I\times I$ be given the usual product CW-decomposition, and let $U$ be made of the left, right and bottom sides of $I\times I$, with the obvious skeletal filtration. The map,  $\Pi(U) \to \Pi\big( I\times I\big)$, induced by the inclusion, is a trivial cofibration of crossed complexes, and hence has the left-lifting property with respect to all fibrations of crossed complexes; see  \cite[Proposition 2.6]{BrownGolasinski}. From this, we have the following:

\begin{Lemma}[The functor derived from a fibration of crossed complexes]\label{Lem:functCRSfib}
 Let $\A=(A_n)_{n \in \mathbb{Z}^+_0}$ and $\B=(B_n)_{n \in \mathbb{Z}^+_0}$ be crossed complexes, and let $p\colon \A \to \B$ be a fibration  between them. There is a functor, $$\Fc^{p}\colon \pi_1(\B,B_0) \to \Sets,$$ in full $\Fc^{p}=\Fc^{(p\colon \A \to \B)}$, such that $\Fc^{p}$ sends $b \in B_0$ to  $\pi_0(p^{-1}(b))$, where the crossed complex $p^{-1}(b)$ is the fibre of $p\colon \A \to \B$, at $b\in B_0$, as in  Definition \ref{def:fibXmap}.
 
 Given a morphism, $[\gamma]\colon b \to b'$, in $\pi_1(\B,B_0)$, the  map, $$\Fc^p([\gamma])\colon \pi_0(p^{-1}(b)) \to \pi_0(p^{-1}(b')),$$ is defined from the right-lifting property of $p\colon \A \to \B$ with respect to the map $\Pi(\{0\}) \to \Pi(I)$.
\end{Lemma}
This last lemma is a  version of Lemma \ref{main_techfibfunctor} for crossed complexes. Note that a result as strong as  Lemma \ref{main_techfibfunctor} does not hold, since fibrations of crossed complexes, as defined above,  do not necessarily satisfy the full homotopy lifting condition, i.e. they are not necessarily Hurewicz fibrations; see \cite[Proposition 2.2]{BrownGolasinski}.
  
\subsubsection{The profunctor $\Hp^{(X;Y,Z)}_\A$}

 Let $X$ be a special CW-complex,  with $Y$ and $Z$ being two disjoint subcomplexes of $X$. As usual, let $\A$ be a crossed complex.
 
There is a natural isomorphism of crossed complexes, $\Pi(Y_\sk \sqcup Z_\sk) \cong \Pi(Y_\sk) \sqcup \Pi(Z_\sk)$,
 so it follows from
 the closed monoidal structure in $\Crs$, that we have a natural isomorphism of crossed complexes, $$\CRS\big (\Pi(Y_\sk)\sqcup \Pi(Z_\sk),\A\big )\cong\CRS\big (\Pi(Y_\sk),\A\big)\times\CRS\big (\Pi(Z_\sk),\A\big).$$
Note also that given crossed complexes, $\Cc=(C_n)_{n \in \mathbb{Z}^+_0}$ and $\Cc'=(C'_n)_{n \in \mathbb{Z}^+_0}$, there is a natural isomorphism of groupoids,
   $$\pi_1(\Cc \times \Cc',C_0\times C_0') \cong \pi_1(\Cc,C_0) \times \pi_1(\Cc',C_0').$$ 

Applying Lemma \ref{lem:incs_to_fibs}, the map, $\Pi(\iota)\colon \Pi(Y_\sk \sqcup Z_\sk)\to \Pi(X_\sk)$, induced by inclusion, gives a  fibration of crossed complexes,
$$p\colon \CRS(\Pi(X_\sk),\A) \to \CRS(\Pi(Y_\sk \sqcup Z_\sk),\A)\cong \CRS(\Pi(Y_\sk),\A)\times \CRS(\Pi(Z_\sk),\A), $$ where $p:=\Pi(\iota)^*=\CRS(\Pi(\iota),\A)$.
Lemma \ref{Lem:functCRSfib} then gives a functor\footnote{If $\B$ is a crossed complex, we use two notations for the fundamental groupoid of $\B$, namely $\pi_1(\B)$ and $\pi_1(\B,B_0)$. If $W$ is CW-complex, and $\A$ is a crossed complex, we  hence abbreviate $ \pi_1\big(\CRS(\Pi(W_\sk),\A), \Crs(\Pi(W_\sk),\A)\big)$ to $\pi_1\big(\CRS(\Pi(W_\sk),\A)\big).$},
\begin{multline*}
\Fc^{(p)}\colon \pi_1\big(\CRS(\Pi(Y_\sk),\A)\big) \times  \pi_1\big(\CRS(\Pi(Z_\sk),\A) \big) \\ \cong \pi_1 \big ( \CRS(\Pi(Y_\sk),\A)\times \CRS(\Pi(Z_\sk),\A)\big)  \to \Sets.
\end{multline*}

This leads to the following profunctor, whose construction mimics that of the profunctor associated to a fibrant span of CGWH spaces; see Subsection
\ref{proffibspan}.

Let $X$, $Y$, $Z$ and $\A$ be as before.
 \begin{Definition}[The profunctor, $\Hp^{(X;Y,Z)}_\A$]\label{def:crs_profunctor}
 The profunctor, 
 $$\Hp^{(X:Y,Z)}_\A\colon \pi_1\big(\CRS(\Pi(Y_\sk),\A)\big)^\op \times  \pi_1\big(\CRS(\Pi(Z_\sk),\A) \big)\to \Sets,$$
is defined as the composite,
\begin{multline*}
      \pi_1\big(\CRS(\Pi(Y_\sk),\A)\big)^\op \times  \pi_1\big(\CRS(\Pi(Z_\sk),\A) \big)  \ra{((-)^{-1} \times \id)}  \\     \pi_1\big(\CRS(\Pi(Y_\sk),\A)\big) \times  \pi_1\big(\CRS(\Pi(Z_\sk),\A) \big)
      \ra{\Fc^{(p)}}
      \Sets.
\end{multline*}
\end{Definition}

Item \eqref{Fib-in-practice} of Lemma \ref{lem:fib_direct} gives a  way to understand the fibration, $p=\Pi(\iota)^*$, and  hence can be used to  write down, explicitly, the profunctor
 $\Hp^{(X:Y,Z)}_\A$.

 \begin{Example}\label{Prof-for-I} In the previous definition, put $X=[0,1]$, $Y=\{0\}$ and $Z=\{1\}$. Let also $G$ be a groupoid, and take $\A=\J_1(G)$.
Let $*$ be the crossed complex with a unique object, and only identity morphisms, so $\Pi(\{0\})$ and $\Pi(\{1\})$ each are isomorphich to $ *$, and $\CRS \big(*,\J_1(G)\big)\cong \J_1(G)$, canonically.

Let $\I$ be the unit interval groupoid. As observed in Example \ref{Ex:simple_CW}, $\Pi([0,1]_{\sk})\cong \J_1(\I)$. By Theorem \ref{homsandotimes}, morphisms $\J_1(\I) \to \J_1(G)$ are in bijection with arrows $g\colon x \to y$, in $G$. It is easy to see that the arrows in the groupoid $\CRS_1\big (\J_1(\I),\J_1(G)\big)$  are in bijection with commmutative diagrams, in $G$, as below,
\begin{equation}\label{fibration:I}
\vcenter{\xymatrix@C=40pt@R=17pt{& x\ar[d]_{h_L} \ar[r]^g  & y\ar[d]^{h_R}\\ & x' \ar[r]_{h_L^{{-1}}\,\, g\,\, h_R}  & y'.}}
\end{equation}
Here, the two  vertical arrows, $h_L\colon x\to x'$ and $h_R\colon y \to y'$, are taken from the groupoid  $G$.
The groupoid composition in $\CRS_1\big (\J_1(\I),\J_1(G)\big)$ is given by  the obvious vertical composition. For $n\ge 2$, all morphisms in $\CRS_n\big(\J_1(\I),\J_1(G)\big)$ are identity morphisms, since all such $n$-fold homotopies are trivial.

The crossed complex fibration, of crossed complex mapping spaces, given by the inclusion of $\{0,1\}$ in the interval $[0,1]$,
 namely, $$
 p\colon \CRS\big(\Pi([0,1]),\J_1(G)\big) \to \CRS\big(\Pi(\{0, 1\})  ,\J_1(G)\big)\cong \J_1(G) \times \J_1(G),
 $$
is obtained from the groupoid map from  $\CRS_1\big (\J_1(\I),\J_1(G)\big)$ to $G\times G$ that chooses the two  vertical arrows in \eqref{fibration:I}.
Since $\pi_1(\J_1(G))\cong G$,   the profunctor, $$\Hp^{([0,1];\{0\},\{1\})}_{\J_1(G)}\colon \pi_1\big(\CRS(\Pi(\{0\}),\J_1(G))\big)  \bto \pi_1\big(\CRS(\Pi(\{1\}),\J_1(G))\big),$$   is given by the profunctor in Example \ref{identity profunctor},  that is, the identity profunctor on $G$.
 \end{Example}

\subsubsection{Fibrations of crossed complexes and fibrations of mapping spaces}\label{sec:fib-Xcompvs-FibMSpaces}
 
The following theorem, which generalises the last point of Theorem \ref{MainCrs}, will be fundamental in giving explicit calculations of TQFTs from crossed complexes. This result appears to be new, however is not unexpected.

Let $S$ be a simplicial set and $T$ a subcomplex of $S$, with $i_{(T,S)}\colon T\to S$ being the inclusion map. Its geometric realisation, $|i_{(T,S)}|\colon |T|\to |S|$, is  an inclusion of CW-complexes (by, for instance, \cite[Corollary 4.38]{FP}), hence it induces a crossed complex map, $\Pi(i_{(T,S)})\colon \Pi(|T|_\sk)\to \Pi(|S|_\sk)$, in fact a cofibration.  Lemma
\ref{lem:incs_to_fibs}, then, gives    a fibration of crossed complexes  between the appropriate internal homs,
$$\Pi(i_{(T,S)})^*\colon \CRS(\Pi(|S|_\sk),\A) \to \CRS(\Pi(|T|_\sk),\A).  $$

Let $f\colon \Pi(|T|_\sk) \to \A$ be a crossed complex map and, using Definition \ref{def:fibXmap}, consider the crossed complex obtained as the fibre over $f\in \CRS_0(\Pi(|T|_\sk),\A)=\Crs(\Pi(|T|_\sk,\A)$ in this fibration,  $\Pi(i_{(T,S)})^*$. We  denote  this crossed complex by
$$\CRS^{(f)}(\Pi(|S|_\sk),\A):=(\Pi(i_{(T,S)})^*)^{-1}(f).$$
This latter crossed complex fits inside the pull-back diagram below,
where $\hat{f}$ is the crossed complex with object set $\{f\}$, and only identity arrows in all dimensions,
\begin{equation}\label{eq:diag_fib_1}
\vcenter{\xymatrix@R=19pt{
& \CRS^{(f)}(\Pi(|S|_\sk),\A)\ar@{}[dr]|{\pullback} \ar[d] \ar[rr]^{\inc}  && \CRS(\Pi(|S|_\sk),\A) \ar[d]^{\Pi(i_{(T,S)})^*} \\
                & \hat{f} \ar[rr]_{\inc} &&\CRS(\Pi(|T|_\sk),\A)\,.
} }
\end{equation}

Since $|i_{(T,S)}|\colon |T|\to |S|$ is  an inclusion of CW-complexes, hence a cofibration, we have a mapping space fibration of CGWH topological spaces,
$$|i_{(T,S)}|^*\colon \TOP(|S|,B_\A) \to  \TOP(|T|,B_\A).$$
We also have the crossed complex map, $f\colon \Pi(|T|_\sk) \to \A$, and this  gives rise, via the adjunction  $\xymatrix{\hskip-0.5mm\Pi\colon \Simp
\ar@/^0.3pc/[r] \ar@{}[r]|\bot \ar@{<-}@/_0.3pc/[r] & \Crs: \hskip-0.5mm{\N},}$  of Proposition \ref{brown_higgins_adjunction}, followed by geometric realisation,  to a continuous map,
$$|{\phi_T^\A(f)}|\colon |T| \to |B_\A|.$$ 

The fibre of $|{\phi_T^\A(f)}|\colon |T| \to |B_\A|$ under the mapping space fibration will be denoted
$$\TOP^{(|{\phi_T^\A(f)}|)}(|S|,B_\A):=(|i_{(T,S)}|^*)^{-1}(|{\phi_T^\A(f)|}),$$
and we have a pullback diagram in $\CGWH$,
$$\xymatrix{
&  \TOP^{(|{\phi_T^\A(f)}|)}(|S|,B_\A) \ar@{}[dr]|{\pullback} \ar[d] \ar[rr]^{\inc}  && \TOP(|S|,B_\A)
\ar[d]^{|i_{(T,S)}|^*} \\
                & \{|{\phi_T^\A(f)}| \}  \ar[rr]_{\inc} &&\TOP(|T|,B_\A)\,.
}  $$
\begin{Theorem}\label{thm:homotopy_equivalence_restriction}
The Brown--Higgins/Brown--Higgins--Sivera/Tonks weak homotopy equivalence,    $$\etab_S^\A\colon  \left |\N \big( \CRS(\Pi(|S|_\sk),\A)\big) \right |  \to \TOP(|S|, B_{\A}),$$
in Theorem \ref{MainCrs},  restricts to a weak homotopy equivalence,
$$ |\N\big(\CRS^{(f)}(\Pi(|S|_\sk),\A)\big)| \to   \TOP^{(|{\phi_T^\A(f)}|)}(|S|,B_\A) .$$
\end{Theorem}

\begin{proof}
Since the nerve functor $\N\colon \Crs \to \Simp$ is a right adjoint,  it preserves  limits. Applying $\N$ to  diagram in \eqref{eq:diag_fib_1}, we have a pullback diagram in $\Simp$
$$\xymatrix{
&\N \big ( \CRS^{(f)}(\Pi(|S|_\sk),\A) \big)\ar@{}[dr]|\pullback \ar[d] \ar[rr]^{\inc}  &&\N\big( \CRS(\Pi(|S|_\sk),\A)\big) \ar[d]^{\N\big(\Pi(i_{(T,S)})^*\big)} \\
                & \N( \hat{f}) \ar[rr]_{\inc} &&\N\big(\CRS(\Pi(|T|_\sk),\A)\big)\,.
} $$

The geometric realisation functor, $\Simp \to \CGWH$, preserves finite limits,  \cite[Theorem 4.3.16]{FP}, so applying geometric realisations to the previous diagram, yields a  pullback diagram in $\CGWH$,
\begin{equation}\label{eq:keydiagram0}
\vcenter{\xymatrix{
&|\N \big(\CRS^{(f)}(\Pi(|S|_\sk),\A)\big)| \ar@{}[dr]|\pullback\ar[d] \ar[rr]^{\inc}  &&|\N \big(\CRS(\Pi(|S|_\sk),\A)\big)| \ar[d]^{{\big|\N\big(\Pi(i_{(T,S)})^*\big)}\big|} \\
                & |\N (\hat{f})| \ar[rr]_{\inc} &&|\N\big(\CRS(\Pi(|T|_\sk),\A)\big)|
\, .} }\end{equation}

We have another commutative diagram in $\CGWH$, arising from the naturality\footnote{Recall \cite{tonks:JPAA:2003} and \cite[Page 177]{BGPT} that the weak homotopy equivalence is only natural with respect to simplicial maps, but not natural in the enriched sense.}, on  varying the simplicial set, $S$,  see \cite[Proposition 3.1.]{BGPT} and \cite[Theorem A]{brown_higgins_classifying}, of the weak homotopy equivalence, $$\etab_S^\A\colon  |\N \big (\CRS(\Pi(|S|_\sk),\A)\big)|  \to \TOP(|S|, B_\A).$$  This gives a commutative diagram,
\begin{equation}\label{eq:keydiagram}\vcenter{
\xymatrix{
& |\N \big(\CRS(\Pi(|S|_\sk),\A)\big)| \ar[rr]^{\etab_S^\A}
\ar[d]_{{\big|\N\big(\Pi(i_{(T,S)})^*\big)\big|}} && \TOP(|S|,B_\A)\ar[d]^{|i_{(T,S)}|^*} \\
              &|\N\big(\CRS(\Pi(|T|_\sk),\A)\big)| \ar[rr]_{\etab_T^\A}  && \TOP(|T|,B_\A)
\,.} }
\end{equation}
We note that diagrams \eqref{eq:keydiagram0} and \eqref{eq:keydiagram} share one of their vertical arrows.

Next, we note the following.
\begin{itemize}[leftmargin=0.6cm]
 \item Given that $\Pi(i_{(T,S)})^*\colon \CRS(\Pi(|S|_\sk),\A) \to \CRS(\Pi(|T|_\sk),\A)$ is a fibration of crossed complexes (by Lemma \ref{lem:incs_to_fibs}), then its nerve, $$\N\big(\Pi(i_{(T,S)})^*\big)\colon \N\big(\CRS(\Pi(|S|_\sk),\A)\big) \to \N\big(\CRS(\Pi(|T|_\sk),\A)\big),$$ is a fibration of simplicial sets,  and so its geometric realisation is a fibration of CGWH topological spaces;  see \cite[Theorem 4.5.25]{FP}\footnote{Note that geometric realisations of Kan fibrations are only sure to have the homotopy lifting property with respect to homotopies whose domain is a CGWH space, see \cite[Page 185]{FP}.}. The left downwards arrow of diagram \eqref{eq:keydiagram} is, therefore, a fibration in $\CGWH$.
 \item Since the inclusion, $|i_{T,S}|\colon |T| \to |S|$, is a cofibration, in $\CGWH$,  the right downwards arrow of diagram \eqref{eq:keydiagram} is a fibration in $\CGWH$.
 \item The two horizontal maps in diagram \eqref{eq:keydiagram} are weak homotopy equivalences.
 \end{itemize}

Since \eqref{eq:keydiagram} commutes, the map  
$\etab_S^\A\colon  |\N (\CRS(\Pi(|S|_\sk),\A))|  \to \TOP(|S|, B_\A)$  sends fibres to fibres.
As   $f\colon \Pi(|T|_\sk) \to \A$ is a crossed complex map, we have that $\etab_T^\A(\tilde{f})=|{\phi_T^\A(f)}|$; cf. the notation in Remark \ref{rem:tilde}. This follows from the fourth point of Theorem \ref{MainCrs}.
The map, $\etab_S^\A$, restricts to a map on the corresponding fibres, in Equation \eqref{eq:keydiagram}, which we denote  $$g'\colon {\big|\N\big(\Pi(i_{(T,S)})^*\big)\big|}^{-1}(\tilde{f})  \to   \TOP^{(|{\phi_T^\A(f)}|)}(|S|,B_\A) .$$

Making use of the pullback diagram \eqref{eq:keydiagram0},  this map, $g'$, gives rise to another map, in $\CGWH$,
$$g\colon |\N(\CRS^{(f)}(\Pi(|S|_\sk),\A))| \to   \TOP^{(|{\phi_T^\A(f)}|)}(|S|,B_\A) ,$$
arising from the canonical homeomorphism, given by the uniqueness of pullbacks,
$$|\N(\CRS^{(f)}(\Pi(|S|_\sk),\A))|\cong    {\big|\N\big(\Pi(i_{(T,S)})^*\big)\big|}^{-1}(\tilde{f})  .$$
This $g$ is exactly the map we want to prove  is a weak homotopy equivalence.

We will show that $g'$, and hence $g$, is  a weak homotopy equivalence, which yields the statement of the theorem. For this we use the homotopy long exact sequences of the two vertical fibrations of diagram \eqref{eq:keydiagram}. The pair of weak equivalences $\etab_S^\A$, of the total spaces, and $\etab_T^\A$, of the base spaces, together with $g'$, the map on the fibre, maps the first long exact sequence to the latter one. Therefore the five-lemma proves that $g$ induces an isomorphism for all homotopy groups, and hence $g$ is a weak homotopy equivalence.
\end{proof}

\begin{Example}\label{eg:GI} Consider the inclusion, $\iota\colon \{0,1\} \to [0,1]$. Let $G$ be a group, viewed as a crossed complex, via $\J_1(G)$. The ensuing fibration of mapping spaces is denoted $P\colon \TOP\big([0,1],B_G\big) \to \TOP\big(\{0,1\},B_G\big)$.  The corresponding fibration of crossed complexes
$p\colon \CRS\big(\Pi([0,1]_\sk),G\big) \to \CRS\big(\Pi(\{0, 1\})  ,G\big)
$
  was  made explicit in Example \ref{Prof-for-I}.

There is only one crossed complex map $f\colon \Pi(\{0,1\}) \to G$.
The fibre of the projection, $p$, at $f$, is the crossed complex $\CRS^{(f)}(\Pi([0,1]),G)$, whose set of objects is the underlying set of $G$, and with only identity morphisms at all orders. From this fact, we can see that the fibre of the fibration,
$P$, over the map $f'\colon \{0,1\} \to B_{G}$, the map that sends both $0$ and $1$ to the unique vertex of $B_G$,
is homotopic to the classifying space of $\CRS^{(f)}(\Pi([0,1]),G)$. The fibre, $P^{-1}(f')$, is, thus, the disjoint union of contractible spaces, one for each element of $G$, as one should expect, since the classifying space $B_G$ is an aspherical space, and $\pi_1(B_G)\cong G$.
\end{Example}

Note that if $\A$ is a crossed complex, then its classifying space, $B_\A$, is canonically filtered by the classifying spaces of the truncations $T_n(\A)$ that are obtained from $\A$ by making all groupoids, $A_i$ with $i> n$, have only identify morphisms. Letting $(B_\A)_*$ be the corresponding filtered space we have, by \cite[
Theorem 7.15]{brown_hha}, an isomorphism, $\A\to \Pi((B_\A)_*)$, natural in $\A$. Because all $n$-cells of $B_\A$ are contained in $(B_\A)_n$, if $X$ is a CW-complex, a CW approximation of a map $f\colon X \to B_\A$ is, automatically, a filtered map $X_\sk \to (B_\A)_*$.

The following follows as for Corollary \ref{cor:CW-class}.
\begin{Corollary}\label{cor:CW-class-cob}
 Let $(X,Y)$ be a CW-pair and let  $\iota\colon Y \to X$ be the inclusion map. Let $f\colon Y \to B_\A$ be a continuous map that is filtered as a map $Y_{\sk} \to (B_\A)_*$. There is a weak homotopy equivalence,
 $$\TOP^{(f)}(X,B_\A)\cong B_{\CRS^{(\Pi(f))}(\Pi(X_\sk),\A)}, $$
 where again we considered fibres in $\CGWH$ and $\Crs$,  in the pullbacks below:
 $$\vcenter{ \xymatrix@C=21pt{\TOP^{(f)}(X,B_\A) \ar@{}[dr]|<<<<<\pullback\ar[d]\ar[r]^{\inc} & \TOP(X,B_\A) \ar[d]^{i^*} \\
              \{f\} \ar[r]_{\inc} & \TOP(Y,B_\A),
 }}   \vcenter{ \xymatrix@C=21pt{\CRS^{\Pi(f)}(\Pi(X_\sk),\A) \ar@{}[dr]|<<<<<\pullback\ar[d]\ar[r]^{\inc} & \CRS(\Pi(X_\sk),\A) \ar[d]^{\Pi(i)^*} \\
              \{{\Pi(f)}\} \ar[r]_{\inc} & \CRS(\Pi(Y_\sk),\A).
 }} $$
\end{Corollary}


\subsection{Computing the homotopy content of a  finite crossed complex}\label{sec:homfinXcomp}
The results in this subsection are essentially all in \cite{MartinsCW}, or \cite{Martins_Porter}.

If $\A=(A_n)_{n \in \mathbb{Z}^+_0}$  is a finite connected crossed complex, then, in order to determine the homotopy content of the classifying space, $B_\A$, one does not need to compute the homotopy groups of $B_\A$. The computation can be reduced to an alternating product of cardinalities of  sets of certain morphisms in $A_n$. This fact (and its proof) is similar to the fact that the Euler characteristic of a finite CW-complex, $X$, i.e. the alternating sum of the ranks of its homology groups, $H_i(X)$, can be computed as $\sum_{i=1}^\infty (-1)^i n_i$, where $n_i$ is the number of $i$-cells of $X$.
As we will see later, for this reason, the formula for the Quinn TQFT, and its once-extended versions, greatly simplify when $\Bc$ is the classifying space of a finite crossed complex, with a single object.
\subsubsection{Finite and homotopy finite crossed complexes}\label{sec:homotopy finite-CRS}

We will start by defining what it means for a crossed complex to be finite, and more generally homotopy finite, which follows naturally from the framework already introduced.
\begin{Definition}\label{Def:finite_reduced} Let $\A=(A_n)_{n \in \mathbb{Z}^+_0}$ be a crossed complex.
\begin{itemize}[leftmargin=1cm]
 \item We say that $\A$ is \emph{finite} if all the groupoids, $A_n$, are finite, and there exists $m\in \mathbb{N}$ such that the groupoids, $A_n$, have only identity arrows for $n\ge m$.
 \item We say that $\A$ is \emph{homotopy finite} if $\A$ has only a finite number of path components, each of which with a finite number of non-trivial homotopy groups, all of which are finite.
 \end{itemize}
 \end{Definition}
Clearly, if a crossed complex is finite then it is homotopy finite, but not conversely. One can use the results of \cite{Ellis} to show that, if $\A$ is a homotopy finite, then it is weakly homotopy equivalent to a finite one.

The following result will be implicitly used several times.
\begin{Lemma}\label{lem:crs_is_finite} Let $X$ be a special finite CW-complex and let $\A$ be a finite crossed complex. The set, $\Crs\big(\Pi(X_\sk),\A\big)$, of crossed complex maps from $\Pi(X_\sk)$ to $\A$, is finite.
\end{Lemma}
\begin{proof}This follows directly from  Lemma \ref{lem:maps_on_basis}.
\end{proof}

\begin{Definition}
Let $\A$ be a path-connected homotopy finite crossed complex. We define the homotopy content of $\A$ as below, where $c$ is any object of $\A$,
$$\chi^{\pi}(\A):=\prod_{i=1}^{\infty} |\pi_i(\A,c)|^{ (-1)^i }.$$
(Note that this computation does not depend on the object $c$, as $\A$ is path-connected.)
More generally,  if $\B$ is  homotopy finite, but not necessarily path-connected, we define
$$\chi^{\pi}(\B):=\sum_{\A \in \hpiz(\B)}\chi^\pi(\A).$$
Here  $\hpiz(\B)$ is the set of path-components of $\B$; see  Definition \ref{Def:pizA}.
\end{Definition}
Theorem \ref{MainCrs} immediately implies the following lemma.
\begin{Lemma}\label{lem:chiNA}
Let $\A$ be a homotopy finite crossed complex, then its classifying space, $B_\A$, is a homotopy finite space, and $\chi^\pi(\A)=\chi^\pi(B_\A)$.
\end{Lemma}

\begin{Notation}[$\Theta_n^x(\A)$]\label{Big theta}
Let  $\A=(A_n)_{n \in \mathbb{Z}^+_0}$ be a finite crossed complex. Given $x \in A_0$ and $n \in \mathbb{Z}^+$, define
$\Theta_n^x(\A)$, a positive integer, to be the cardinality of the set of morphisms, in the groupoid $A_n$, with source $x$.
\end{Notation}

\begin{Remark}Note that, fixing $n \in \mathbb{Z}^+_0$, then $\Theta_n^x(\A)\in \mathbb{Z}^+$ depends only on the path component in $\A$, or equivalently in $\pi_1(\A)$, to which $x\in A_0$ belongs.
\end{Remark}

The following result appears in \cite[\S 4.2.2]{MartinsCW}.

\begin{Lemma}\label{lem:calc-chi-Euler}
Let  $\A=(A_n)_{n \in \mathbb{Z}^+_0}$ be a finite crossed complex, then
$$ \chi^\pi(\A)=\sum_{x \in A_0}\Big ( \prod_{i=1}^{\infty} \big(\Theta_i^x(\A)\big) ^{ (-1)^i}\Big).
$$
\end{Lemma}
\begin{proof}
 This follows from a  telescopic calculation, similar to the proof of Theorem \ref{main1}. A crucial point in the proof, allowing us to pass from a sum over path components of $\A$ to a sum over objects of $\A$, is that, if  $(s,t\colon G^1 \to G^0)$ is a finite groupoid, then given $x \in G_0$,  the cardinality of the set of morphisms in $G$, with source $x$, is equal to $|G(x,x)|\,|[x]|$. Here $[x]$ is the set of objects of $G$ connected to $x$ by a morphism of $G$.
 Full details are  in \cite[Lemma 4.8]{MartinsCW}.
\end{proof}
\begin{Definition}\label{homogeneous}
 A finite crossed complex, $\A=(A_n)_{n \in \mathbb{Z}^+_0}$, will be called \emph{homogeneous} if, given a non-negative integer $n$, and an object $x$ of $\A$, the value of $\Theta_n^x(\A)$ depends only on $n$. This means that there exists, for each non-negative integer $n$, a positive integer,  $\Theta_n(\A)$, such that, for each $x \in A_0$,  we have $\Theta_n^x(\A) =\Theta_n(\A)$.

 If $\A$ is homogeneous, define
 $$\Theta(\A)=  \prod_{i=1}^{\infty} (\Theta_n(\A))^{ (-1)^i}\in \mathbb{Q}.$$
\end{Definition}
\noindent Note that path-connected finite crossed complexes are automatically homogeneous.
\begin{Corollary}\label{cor-hom}
 Suppose that $\A=(A_n)_{n \in \mathbb{Z}^+_0}$ is homogeneous, and so, in particular, finite, then
 $$\chi^\pi(\A)=\Theta(\A) \, |A_0|. $$
\end{Corollary}

\subsubsection{The homotopy content of $\CRS\big(\Pi(X_\sk),\A\big)$}\label{CRS-euler-formula}
Let the  CW-complex, $X$, be finite and let $Y$ be a subcomplex of $X$.
\begin{Notation}[$\L(n,X)$ and $\L(n,X,Y)$]\label{not:K(X,n)}
 Let  $n$ be a non-negative integer.
 \begin{itemize}[leftmargin=1cm]
\item  We set $\L(n,X)$ to be the number of $n$-cells of $X$.  \item More generally, let $Y$ be a subcomplex of $X$. An $n$-cell, $c$, of $X$ is said to be \emph{internal} to $(X,Y)$, if it is not in $Y$. We let
 $\L(n,X,Y) $
 be the number of $n$-cells of $X$ that are internal to the pair $(X,Y)$.
 \end{itemize}
 
\end{Notation}

By applying the results in \S \ref{sec:hom_of_freeXcomp},  and, in particular,  Lemma \ref{lem:derivation_basis}, together with the notions and notation introduced earlier in this Subsection \ref{sec:homfinXcomp},  we have:
\begin{Lemma}  Let $\A=(A_n)_{n \in \mathbb{Z}^+_0}$ be a finite  crossed-complex with a single object. Let $X$ be a finite special CW-complex. 
The crossed complex $\CRS(\Pi(X_\sk),\A)$ is homogeneous, in the sense of Definition \ref{homogeneous}, and, in particular, finite. Moreover,  we have that, for a positive integer $j$,
$$\Theta_j\big(\CRS(\Pi(X_\sk),\A)\big)=\prod_{i=0}^\infty |A_{i+j}|^{\L(i,X)} .$$
In particular,
$$\chi^\pi\big (\CRS(\Pi(X_\sk),\A)\big)=\big|\Crs(\Pi(X_\sk),\A)\big|\, \prod_{j=1}^\infty \big (\prod_{i=0}^\infty |A_{i+j}|^{\L(i,X)}  \big )^{(-1)^j}. $$\end{Lemma}
This is the special case in which $Y$ is empty, of the following more general result, in which we let $(X,Y)$ be a pair of finite special CW-complexes.
 
 Let $f\colon \Pi(Y_\sk)\to \A$ be a crossed complex map and, as before, let  $\Pi(i)\colon \Pi(Y_\sk) \to \Pi(X_\sk)$ be induced by the inclusion $i\colon Y \to X$. Consider also the induced  fibration  (Lemma \ref{lem:incs_to_fibs}) between the `internal homs',
 $$\Pi(i)^*\colon \CRS\big(\Pi(X_\sk),\A\big) \to \CRS\big(\Pi(Y_\sk),\A\big).$$
 \begin{Lemma} The fibre of $\Pi(i)^*$ at $f$,  Definition \ref{def:fibXmap}, i.e., the crossed complex, $$\CRS^{(f)}\big(\Pi(X_\sk),\A\big):=(\Pi(i)^*)^{-1}(f),$$ is homogeneous, and if $j$ is a positive integer, then
 $$\Theta_j\big(\CRS^{(f)}(\Pi(X_\sk),\A)\big)=\prod_{i=0}^\infty |A_{i+j}|^{\L(i,X,Y)}.$$
 In particular
 \begin{multline*}
 \chi^\pi\big(\CRS^{(f)}(\Pi(X_\sk),\A)\big)=\\ \big|\big\{g\colon \Pi(X_\sk) \to \A\mid g \circ \Pi(i)=f \big \}\big|\,\, \prod_{j=1}^\infty \big (\prod_{i=0}^\infty |A_{i+j}|^{\L(i,X,Y)}  \big)^{(-1)^j}.
 \end{multline*}
 \end{Lemma}
 \begin{proof}This follows from Lemmas  \ref{lem:derivation_basis}, \ref{lem:fib_direct} and \ref{lem:crs_is_finite}. \end{proof}

The particular case when we have a CW-triad, $(X;Y,Z)$, of finite  special CW-complexes will be very useful when we come to write down explicit formulae for TQFTs derived from crossed complexes.  By a CW-triad  $(X;Y,Z)$ here we mean  that  $Y$ and $Z$ are disjoint subcomplexes of a special CW-complex, $X$, which then implies that  $Y\sqcup Z$ is  a subcomplex of $X$. (The example to have in mind is a triangulated cobordism,  $M\colon S \to S'$, where $X=M$, $Y= S$ and $Z = S'$.)
Let $i_{(Y,X)}\colon Y \to X$ and $i_{(Z,X)}\colon Z \to X$ be the inclusion maps. We have a cellular map, $$\soml{ i_{(Y,X)}}{i_{(Z,X)}} \colon Y \sqcup Z  \cong Y \cup Z \to X,$$
which gives the inclusion.
Given crossed complex maps, $f\colon \Pi(Y_\sk) \to \A$ and $f'\colon \Pi(Z_\sk) \to \A$, we can combine them into a crossed complex map, $$\soml{ f}{f'} \colon \Pi\big((Y\cup Z)_\sk\big)\cong \Pi(Y_\sk)\sqcup \Pi(Z_\sk)\to \A.$$

The set of objects of the crossed complex,
$\CRS^{( \soml{ f}{f'})}(\Pi(X_\sk),\A)$, is the set of crossed complex maps, $h\colon \Pi(X_\sk)\to \A$, that make the diagram below commute,
$$\xymatrix@R=12pt
{ && \A\\
& \Pi(Y_\sk) \ar[dr]_{\Pi(i_{(Y,X)})\,\,\,} \ar[ur]^f &
& \Pi(Z_\sk)\ar[ul]_{f'} \ar[dl]^{\,\,\,\Pi(i_{(Z,X)})} \\
&& \Pi(X_\sk).\ar[uu]^h
}$$
In particular, we have  the following.
\begin{Lemma}\label{lem:twosides}
Let $X$ be a finite special CW-complex, with $Y$ and $Z$, two disjoint subcomplexes. Let  $\A=(A_n)_{n \in \mathbb{Z}^+_0}$ be a finite reduced crossed complex. Let $f\colon \Pi(Y_\sk)\to \A$ and $f'\colon \Pi(Z_\sk)\to \A$ be crossed complex maps, then
\begin{multline*}
\chi^\pi\big(\CRS^{( \soml{ f}{f'})}(\Pi(X_\sk),\A)\big )=\\\big |\{h\colon \Pi(X_\sk) \to \A: h \circ \Pi(i_{(Y,X)})=f \textrm{ and }  h \circ \Pi(i_{(Z,X)})=f' \}\big|\\
 \prod_{j=1}^\infty \big (\prod_{i=0}^\infty |A_{i+j}|^{\L(i,X,Y\cup Z)}  \big )^{(-1)^j}.
\end{multline*}
\end{Lemma}
\begin{proof} This follows from the discussion just before the lemma.\end{proof}

Finally, let $X$ be a special CW-complex,  $\A=(A_n)_{n \in \mathbb{Z}^+_0}$ a finite reduced crossed complex, and let $f\colon \Pi(X_\sk) \to \A$ be a crossed complex map. By passing to the path-component, $\PC_f(\CRS(\Pi(X_\sk),\A))$,  of $f$ in the crossed complex $\CRS(\Pi(X_\sk),\A)$, we have, within the same context as before:
\begin{Lemma}\label{path:compcrs}
 Let $f\colon \Pi(X_\sk) \to \A$ be a crossed complex map, then
$$\chi^\pi\big( \PC_f(\CRS(\Pi(X_\sk),\A) \big)=|[f]_{\CRS(\Pi(X_\sk),\A)}| \prod_{j=1}^\infty \big (\prod_{i=0}^\infty |A_{i+j}|^{\L(i,X)}  \big )^{(-1)^j}. $$
Here $[f]_{\CRS(\Pi(X_\sk),\A)}$  denotes the  homotopy class of  $f$ (the set of all crossed complex maps, $\Pi(X_\sk) \to \A$, that are homotopic to $f$).
\end{Lemma}

\subsection{Example computations}\label{sec:computations}
In preparation for the next section, with examples of TQFTs and once-extended TQFTs derived from crossed complexes, we show some key computations of crossed complexes of the form $\CRS(\Pi(X_\sk),\A)$,  as discussed in \S \ref{sec:internalhom},  and, for their  homotopy content, following the discussion in \S \ref{CRS-euler-formula}. Here $X$ is a CW-complex and $\A$ is a crossed complex. We will also show examples, for $Y$  a subcomplex  of $X$, and a crossed complex map $f\colon \Pi(Y_\sk) \to \A$, of the form of the fibre, $\CRS^{(f)}\big(\Pi(X_\sk),\A\big):=(\Pi(\iota)^*)^{-1}(f)$, of the restriction map $\Pi(\iota)^*\colon \CRS\big(\Pi(X_\sk),\A\big)\to \CRS\big(\Pi(Y_\sk),\A\big)$, at $f$.
 Here $\Pi(\iota)\colon \Pi(X_\sk)\to\Pi(Y_\sk)$ arises from the inclusion $\iota\colon Y\to X$.
 We will also show  more examples of the profunctors $\Hp^{(X:Y,Z)}_\A$ of Definition \ref{def:crs_profunctor},  for disjoint subcomplexes, $Y$ and $Z$, of $X$.

We will focus on the easiest  example of crossed complexes that are not 1-truncated, which are those arising from crossed modules (of groups), which, by definition, are reduced 2-truncated crossed complexes.  Some of the computations we show are also in \cite{Porter_HQFT,Porter_Turaev_HQFT}, in the context of homotopy quantum field theories.
\subsubsection{Crossed modules of groups}\label{Crossed modules of groups}

 The  definition of crossed modules of groups is classical, going back at least to Whitehead's original papers on CW-complexes and crossed complexes, \cite{CH1,CH2}. Recent treatments are in \cite{Baez_Lauda,brown_hha,brown_higgins_sivera}. The category of crossed modules is equivalent to the category of 2-groups; see e.g. \cite[\S 2.5]{brown_higgins_sivera} and \cite{Baez_Lauda}.  The homotopy category of crossed modules is equivalent to the homotopy category of 2-types: pointed spaces $(X,*)$ with $\pi_i(X,*)=0$ for $i\ge 3$; see  \cite{Baues_4D,MacLane_Whitehead}.

\begin{Definition} A crossed module, $\Gc=(\d\colon E \to G, \trl)$, of groups is given by:
\begin{itemize}[leftmargin=1cm]
 \item a group homomorphism $\d\colon E \to G$,
 \item[]\hspace{-14mm} together with
 \item a right-action, $\trl$, of $G$ on $E$ by automorphisms. \end{itemize}
This action is such that:
 \begin{enumerate}[leftmargin=1cm]
  \item $\d(a\trl g)=g^{-1}\, \d(a)\, g$, for all $a \in E$, $g \in G$ (called the  first Peiffer relation),
   \item $a\trl \d(e)=e^{-1}\, a\, e$, for all $a,e \in E$ (called the second Peiffer relation).
 \end{enumerate}
\end{Definition}

A crossed module, $\Gc=(\d\colon E \to G, \trl)$, of groups gives rise to a reduced crossed complex, $\J_2(\Gc)$. Explicitly $\J_2(\Gc)$ has the form,
\begin{equation}\label{eq:iota_2}
\J_2(\Gc)=\dots \to \{1\}  \to \{1\} \to E \ra{\d} G\to \{*\}.
\end{equation}
The classifying space, $B_\Gc$, of a crossed module, $\Gc$, is, by definition, the same as the classifying space of the crossed complex $\J_2(\Gc)$.
If the crossed module, $\Gc$, is finite, meaning that both $E$ and $G$ are finite, then $B_\Gc$ will be  homotopy finite.

Our examples of mapping spaces, $\CRS\big(\Pi(X_\sk),\J_2(\Gc)\big)$, will be written down in terms of action groupoids, and semidirect products. Let us explain our conventions.
\begin{Definition}\label{def:act_groupoid}
 Let the group $G$ have a left-action, $\bullet$, on a non-empty set  $X$. The \emph{action groupoid}, $X \sslash G$, or, in full, $X \sslash_\bullet G$,
  has $X$ as  its set of objects. Given $x,y \in X$, the set of morphisms, from  $x$ to $y$, is given by the set of all   $g\in G$, with $g \bullet x = y$. The composition of morphisms in $X\sslash G$ is then as indicated below:  $$\big(x \ra{g} g \bullet x  \ra{h } (hg) \bullet x \big)=\big(x\ra{hg} (hg) \bullet x \big).$$
If $G$ right-acts on $X$, via $\trl$, we will also consider a corresponding action groupoid, $X\sslash_{\trl} G$, where arrows look like $x \ra{g} x\trl g$, with composition: $$\big(x \ra{g}  x\trl g  \ra{h }   x\trl(gh) \big)=\big(x\ra{gh} x\trl (gh) \big).$$
\end{Definition}
By Corollary \ref{cor-hom}, if $X$ is a set and $G$  a group, both finite, then
$\chi^\pi\big (\J_1(X\sslash G)\big)=|X|/|G|,$ where $\J_1\colon \Grp \to \Crs$ is the inclusion, defined in Subsection  \ref{sec:defXcomp}.

If a group $G$ right-acts  on the group $E$ our convention for  $G \ltimes E$ is
 $$(g',e') (g,e):=  \big (g'\, g, e \, ( e' \trl g)\big),\textrm{ hence }  (h,e)^{-1}=(h^{-1},e^{-1}\trl h^{-1}).$$
 This non-standard convention for semidirect products arises from the construction, in Subsection \ref{sec:xcomphom}, of the groupoid, $\CRS_1(\A,\B)$, of crossed complex maps from the crossed complex, $\A$, to the crossed complex, $\B$, and homotopies between them.

We will fix a crossed module $\Gc=(\d\colon E \to G, \trl)$ until the end of this subsection.

\subsubsection{Computation related to $\CRS\big(\Pi(S^1_\sk),\J_2(\Gc)\big)$}\label{sec:G/GG} Consider $S^1$ with a CW-decomposition with a unique $0$-cell, at the south-pole, and one $1$-cell. Let us compute the crossed complex $\CRS\big(\Pi(S^1_\sk),\J_2(\Gc)\big)$.  This construction is also in \cite[Example 9.3.8]{brown_higgins_sivera},  and  \cite{Morton_Picken_2groupactions_and_moduli}, in the language of 2-groupoids.

From Lemma \ref{lem:maps_on_basis},  crossed complex maps from  $\Pi(S^1_\sk)$ to $\J_2(\Gc)$ are in bijection with elements of $G$. This can also be derived from \eqref{eq:iota_2} and that
\begin{equation}\label{eq:pi2S1G}
\Pi(S^1_\sk)\cong \cdots \to \{ 0\} \to \mathbb{Z} \to \{*\}, \quad \J_2(\Gc)=\dots \to \{1\} \to E \ra{\d} G\to \{*\}.
\end{equation}
In order to describe homotopies (of various order) between these two crossed complexes, we can use  Lemma \ref{lem:derivation_basis}, or a direct calculation, to see that, given a map  $f\colon \Pi(S^1_\sk)\to \J_2(\Gc)$, homotopies with target $f$ (i.e. 1-fold $f$-homotopies) are in one-to-one correspondence with elements of $G\times E $, seen as a set. On the other hand, 2-fold $f$ homotopies are in one-to-one correspondence with elements of $E$.

To describe how each such crossed complex $f$-homotopy modifies  $f\colon \Pi(S^1_\sk) \to \J_2(\Gc)$, as explained just after Definition \ref{def:f-homotopy}, we use the following result, which is then also used to determine the rest of  $\CRS\big(\Pi(S^1_\sk),\J_2(\Gc)\big)$. This is motivated by the construction in Subsection \ref{sec:xcomphom}, and the diagram on the left-hand-side of Figure \ref{eq:act-g}, representing a homotopy of maps from  $\Pi(S^1_\sk)$ to $\J_2(\Gc)$, seen as a map $\Pi\big (( S^1\times I)_\sk\big) \to \J_2(\Gc)$; we are using Theorems \ref{homsandotimes} and \ref{thm:tensor-CW} here, and have $g,h \in G$ and $e \in E$.
(The other diagrams in Figure \ref{eq:act-g} consider different CW-decomposition of the annulus $S^1\times I$. They will be addressed latter.)
\begin{figure}
\begin{tikzpicture}[scale=0.5]
\draw[fill=blue!10]  circle (2);
\draw[fill=white] circle (0.8);
\draw[<-,thick] (0,-0.85) -- (0,-2) node[midway,right]{\small{$h$}};
\draw[thick] (0,0) circle (0.8cm);
\draw[<-, thick] (0.1,0.8) -- (0,0.8);
\draw[<-, thick] (0.1,2) -- (0,2);
\draw[thick] (0,0) circle (2cm);
\node at (0,1.4) {\small{$e$}};
\node at (-0.44,0) {\small{$g$}};
\node at (0  ,2.45) {\small{$  h\, g \,  \d( e) \, h^{-1} $}};
\draw[fill=black] (0,-2) circle (0.1cm);
\draw[fill=black] (0,-0.8) circle (0.1cm);
\end{tikzpicture}
\quad
\begin{tikzpicture}[scale=0.5]
\draw[fill=blue!10]  circle (2);
\draw[fill=white] circle (0.8);
\draw[<-,thick] (0,-0.8) -- (0,-2) node[midway,right]{\small{$x$}};
\draw[<-,thick] (0,0.8) -- (0,2) node[midway,right]{\small{$y$}};
\draw[thick] (0,0) circle (0.8cm);
\draw[->,thick] (0.8,0.1) -- (0.8,0);
\draw[->, thick] (2,0.1) -- (2,0);
\draw[<-,thick] (-0.8,0.1) -- (-0.8,0);
\draw[<-, thick] (-2,0.1) -- (-2,0);
\draw[thick] (0,0) circle (2cm);
\node at (-1.4,0) {\small{$a$}};
\node at (1.4,0) {\small{$b$}};
\node at (-0.44,0) {\small{$g$}};
\node at (0.44,0) {\small{$h$}};
\node at (-2.7,1.7) {\small{$  x \, g \,  \d( a) \, y^{-1} $}};
\node at (3.2,1.7) {\small{$  y\, h \,  \d( b) \, x^{-1} $}};
\draw[fill=black] (0,-2) circle (0.1cm);
\draw[fill=black] (0,-0.8) circle (0.1cm);
\draw[fill=black] (0,2) circle (0.1cm);
\draw[fill=black] (0,0.8) circle (0.1cm);
\end{tikzpicture}
\quad
\begin{tikzpicture}[scale=0.5]
\draw[fill=blue!10]  circle (2);
\draw[fill=white] circle (0.8);
\draw[<-,thick] (0,-0.8) -- (0,-2) node[midway,right]{\small{$x$}};
\draw[thick] (0,0) circle (0.8cm);
\draw[->,thick] (0.8,0.1) -- (0.8,0);
\draw[->, thick] (2,0.1) -- (2,0);
\draw[<-,thick] (-0.8,0.1) -- (-0.8,0);
\draw[<-, thick] (-2,0.1) -- (-2,0);
\draw[thick] (0,0) circle (2cm);
\node at (-0.44,0) {\small{$g$}};
\node at (0.44,0) {\small{$h$}};
\node at (-.5,2.4) {\small{\qquad \quad$  x \, g   h \,  \d( a) \, x^{-1} $}};
\draw[fill=black] (0,-2) circle (0.1cm);
\draw[fill=black] (0,-0.8) circle (0.1cm);
\draw[fill=black] (0,0.8) circle (0.1cm);
\node at (0,1.4) {\small{$a$}};
\end{tikzpicture}
\caption{\label{eq:act-g} Three CW-decompositions of the annulus $A=S^1\times I$.}
\end{figure}
\begin{Lemma}\label{lem:exerG//G} The following hold:
\begin{enumerate}[leftmargin=1cm]
 \item\label{it-1H}
We have a left-action, $\bullet$, of $G \ltimes_\trl E$ on the underlying set of $G$,  such that, given $g \in G$ and $(h,e)\in  G \ltimes_\trl E,$ we have $(h,e) \bullet g :=h \, g \, \d( e) \, h^{-1}$.
\item This action restricts to an action of $G\ltimes_\trl E$ on $\d(E)$.
\item \label{it-2H}If $g\in G$, we have a morphism $F_g\colon E \to G \ltimes_\trl E$, with $a \mapsto \big (\partial(a), a^{-1} \triangleleft g \,\, a\big)$, and hence, a right-action of $E$ on the set $G\times E$, with $(p,c)\tl_g e=(p,c) F_g(e)$.
\item \label{it-3H}If $g\in G$ and $a\in E$, then $F_g(a)\bullet g=g$.
\item \label{it-5H} $F_{(h,e)\bullet g}\big( a\triangleleft h \big)=(h,e)^{-1} F_g(a) (h,e)$.
\end{enumerate}
\end{Lemma}

\begin{proof}
These are all standard computations. The hardest of which is \eqref{it-5H}:
\begin{align*}
 F_{(h,e)\bullet g}\big( a\triangleleft h^{-1} \big)&=\big ( h\partial(a) h^{-1}, a^{-1}\triangleleft  (  \, g \, \d( e) \, h^{-1}) \,\, a \triangleleft h^{-1}\big)\\
 &=\big ( h^{-1}\partial(a) h,  e^{-1} \trl h^{-1} \,\, a^{-1}\triangleleft  (  g h^{-1}) \,   e \trl h^{-1} \,  a \trl h^{-1} \big),
\end{align*}
\begin{align*}
 (h,e)  F_g(a) (h,e)^{-1} &= (h,e) (\partial(a), a^{-1} \triangleleft g \,\, a\big)  (h^{-1}, e^{-1}\triangleleft h^{-1}) \\
 &=(h\partial(a) h^{-1},  e^{-1}\triangleleft h^{-1}\,\,  a^{-1} \triangleleft (gh^{-1}) \,\, a\triangleleft h^{-1}\,\, e \triangleleft (\partial(a) h^{-1})\\
  &=(h\partial(a) h^{-1},  e^{-1}\triangleleft h^{-1}\,\,  a^{-1} \triangleleft (gh^{-1}) \,\,  e \triangleleft h^{-1}\,\,  a \triangleleft h^{-1}).
 \end{align*}
 The remaining details are left to the reader.
\end{proof}

 Unpacking the construction in Subsection \ref{sec:xcomphom}, we have:
\begin{Lemma}
  The crossed complex $\CRS\big(\Pi(S^1_\sk),\J_2(\Gc)\big)$ is isomorphic to the following crossed complex, denoted $G\sslash \Gc$, with set of objects $G$,
  $$G\sslash\Gc:=  \cdots\ra{\partial} \bigsqcup_{g \in G}\{1\} \ra{\partial}   \bigsqcup_{g \in G} E \ra{\partial}  G \sslash( G\ltimes E).$$
  Here:
  \begin{enumerate}[leftmargin=1cm]
   \item The groupoid $\bigsqcup_{g \in G} E$ is given by the set map $\beta\colon \bigsqcup_{g \in G}E \to G$, identifying the component of the disjoint union, with composition   via the product in $E$.
  \item The groupoid map $\partial\colon \bigsqcup_{g \in G} E \to   G \sslash( G\ltimes E)$  sends $g\ra{a}g$ to $g\ra{F_g(a)} g$.
   \item $\big( g \ra{a} g\big)\triangleleft\big(  g\ra{(h,e)} (h,e) \bullet  g \big)= (h,e) \bullet  g   \ra{a \triangleleft h^{-1}}  (h,e) \bullet  g $.
 \end{enumerate}
 \end{Lemma}
\noindent The fact that $G\sslash\Gc$ is a crossed complex follows from Lemma \ref{lem:exerG//G}.

\begin{Example}\label{CRSPS^1G} When $E$ is trivial, this gives, for a group $G$,
$$\CRS\big( \Pi(S^1_\sk),\J_1(G)\big)\cong \J_1(G\sslash G), $$
where $G\sslash G$ is the action groupoid of the left-action of $G$ on itself by conjugation, that we met back in Example \ref{ex.Qdoublefingroup}.
\end{Example}

Note that $\pi_0(G \sslash \Gc)$ is the set of orbits, $G/(G\times E)$, of the action $\bullet$. This is in clear bijection with the set of conjugacy classes of the quotient group $G/\partial(E)$.

The fundamental groupoid, $\pi_1(G\sslash\Gc)$, of $G\sslash\Gc$, as in Definition \ref{def:pi1andpi0crs}, has $G$ as its set of objects. Given $g,g'\in G$, morphisms from $g$ to $g'$ are equivalence classes, $[h,e]_g$, of  pairs $(h,e) \in G\ltimes E$, with $(h,e)\bullet g=g'$, where
\begin{equation}\label{eq:cong-pi1}
(h,e) \sim_g (h',e') \textrm{  if there exists }  a \in E, \textrm{ such that } (h',e')= (h,e) \tl_g a.
\end{equation}
Hence $(h,e) \sim_g (h',e')$   if  $h'= h \, \d(a)$ and
$ e' =  ( a^{-1} \trl g ) \, e\, a$, for some $a \in E$.
The composition in  $\pi_1(G\sslash\Gc)$ is inherited from the composition on $G \sslash( G\ltimes E)$.

By Lemma \ref{path:compcrs}, because $S^1$ has unique $0$ and $1$ cells, given $g\in G$, then $$\chi^\pi\big(\PC_g(G\sslash \Gc)\big)= \frac{|\Orb_{G \ltimes E} (g)|}{|G|}.$$ Hence picking representatives $g_i$ of path-components in $G\sslash \Gc$, or by Lemma \ref{sec:homfinXcomp},
$$\chi^\pi(G\sslash \Gc)=\sum_i \frac{|\Orb_{G \ltimes E} (g_i)|}{|G|}  =1.$$

Consider the annulus, $A=S^1\times I$, with a CW-decomposition with two 0-cells, three $1$-cells and one 2-cell, where $S^1$ embeds cellularly as $S^1\times \{0\}$ and $S^1\times \{1\}$. (This is as in the leftmost diagram in Figure \ref{eq:act-g}.) If we combine these inclusions, $\iota_0$ and $\iota_1$, of $S^1$ inside $A$, we get a crossed complex map, $\langle \Pi_1(\iota_0),\Pi_1(\iota_1)\rangle\colon \Pi(S^1_\sk\sqcup S^1_\sk) \to \Pi(A_\sk)$,  which induces a map by composition that we will denote
\begin{multline*}
P\colon \CRS\big (\Pi(A_\sk),\J_2(\Gc)\big)\to \CRS\big(\Pi(S^1_\sk\sqcup S^1_\sk),J_2(\Gc)\big)\cong (G \sslash \Gc)\times (G \sslash \Gc).
\end{multline*}

Given that we only have a 1-cell in $A$ which is neither in $S^1\times \{0\}$ nor in $S^1\times \{1\}$, it follows that the fibre of $P$ at $(g,h) \in G\times G$ is  $P^{-1}(g,h)= \J_1\big(\overline{\Ap}(g,h)\sslash_{\tl_g} E\big )$, where $\overline{\Ap}(g,h):=\{(x,e)\in G\ltimes E: (x,e)\bullet g=h\}$. 
In particular, we have
\begin{equation}\label{eq:fib-A}
 \chi^\pi\big (P^{-1}(g,h)\big) = |\{(x,e)\in G\ltimes E: (x,e)\bullet g=h\}|/|E|.
\end{equation}

The profunctor, arising from Definition \ref{def:crs_profunctor}, namely
$$  \Hp^{(S^1\times I : S^1\times \{0\}, S^1\times \{1\})}_{\J_2(\Gc)}\colon \pi_1\big(\CRS(\pi_1(S^1_\sk),\J_2(\Gc)\big) \bto  \pi_1\big(\CRS(\pi_1(S^1_\sk),\J_2(\Gc)\big), $$
considered as a functor, $\Ap\colon \pi_1(G\sslash \Gc)^\op\times  \pi_1(G\sslash \Gc)\to \Sets$, is such that:
\begin{itemize}[leftmargin=0.6cm]
 \item On objects,
$(g,h)\mapsto \Ap(g,h)$, where $\Ap(g,h):=\overline{\Ap}(g,h)/E$.
\item On morphisms,
\begin{multline*}
\Big((y,b)^{-1}\bullet g \ra{[(y,b)]_{(y,b)^{-1}\bullet g} } g  , h \ra{[(z,c)]_h} (z,c) \bullet h\Big) \\ \stackrel{\Ap}{\longmapsto} \left ( \vcenter{\xymatrix@C=9pt@R=1pt{ \A(g,h) \ar[r]  &\Ap\big( (y,b)^{-1}\bullet g, (z,c)\bullet h\big)  \\
                                            [(x,e)]_g \ar@{|->}[r] & [(y,b) (x,e)(z,c)]_{(y,b)^{-1}\bullet g} }}\right)
.\end{multline*}
\end{itemize}
\noindent In particular $\Ap\colon \pi_1(G\sslash \Gc)\bto \pi_1(G\sslash \Gc)$ is the identity profunctor, in Example \ref{identity profunctor}.

\subsubsection{Computation of $\CRS(\Pi(S^1_{\sk'}),\Gc)$}\label{sec:S2p}

We now briefly discuss the crossed complex $\CRS\big(\Pi(S^1_{\sk'}),\J_2(\Gc)\big)$, where $S^1$ has the CW-decomposition with two $0$-cells, at the south and north poles.  By Lemmas \ref{lem:maps_on_basis} and \ref{lem:derivation_basis}, we have a  bijection between crossed complex maps $\Pi(S^1_{\sk}) \to \J_2(\Gc)$ and pairs $(g,h)\in G\times G$, and furthermore homotopies ending in such crossed complex maps are given by elements $(x,y,a,b)\in G\times G \times E \times E$. Each of the latter, changes a crossed complex map in the way shown in the diagram in the centre of Figure \ref{eq:act-g}.

Consider the group $(G\times G)\ltimes (E\times E)$, with action  $(a,b)\trl (x,y)=(a\trl y, b \trl x)$, for $x,y \in G$ and $a,b \in E$. The group $(G\times G)\ltimes (E\times E)$ acts on $G\times G$ as:
$$(x,y,a,b)\bullet (g,h)=\big(xg\partial(a)y^{-1}, yh\partial(b) x^{-1}\big). $$
Consider, given $(g,h)\in G \times G$, the homomorphism $F_{(g,h)}\colon E \times E \to (G\times G)\ltimes (E\times E)$ with $(a,b)\mapsto \big(\partial(a),\partial(b),(a^{-1} \trl g) \,  b, (b^{-1}\trl h)\, a \big)$, then
 $F_{(g,h)}(a,b)$ stabilises $(g,h)$.

The following lemma follows as for $\CRS\big(\Pi(S^1_\sk),\J_2(\Gc)\big)$.

\begin{Lemma}
 We have that $\CRS\big(\Pi(S^1_{\sk'}),\J_2(\Gc)\big)\cong  (G\times  G)\sslash  \Gc^{(2)}$, where
  $$(G\times  G)\sslash  \Gc^{(2)}:=  \cdots\ra{\partial} \bigsqcup_{(g,h) \in G\times G}\{1\} \ra{\partial}   \bigsqcup_{(g,h) \in G\times G} E\times E \ra{\partial}  G\times G \sslash( G\times G)\ltimes (E\times E).$$
  Here $\partial\big ( (g,h)\ra{(a,b)} (g,h) \big)=(g,h)\ra{F_{(g,h)}(a,b)} (g,h),$ and
 \begin{multline*}
 \big ( (g,h)\ra{(a,b)} (g,h) \big)\trl\big( (g,h)\ra{(x,y,m,n)} (x,y,m,n)\bullet (g,h)\big) \\= \big( (x,y,m,n)\bullet (g,h)\big) \ra{ (a\trl x^{-1}, b \trl y^{-1})} \big ((x,y,m,n)\bullet (g,h)\big).
\end{multline*}
 \end{Lemma}

 There is a crossed complex map $P_1^2\colon (G\times  G)\sslash  \Gc^{(2)}\to G\sslash \Gc$, such that, on objects,  $P_1^2(g,h)=gh$, and on 1- and 2-morphisms, we have:
 \begin{align*}
 P_1^2\big( (g,h)\ra{(x,y,m,n)} (x,y,m,n)\bullet (g,h)\big)& =\big ( (gh)\ra{(x,a\,(b \trl h^{-1})} (x,ab)\bullet (gh)),\\
 P_1^2\big( (g,h)\ra{(m,n)} (g,h) \big)&= \big( gh\ra{m} gh \big ).
 \end{align*}
 In geometric terms, $P_1^2$ is induced by the identity map of $S^1$, which gives a cellular map $S^1_{\sk'}\to S^1_{\sk}$.
 It can be easily seen that  $P_1^2$ induces a bijection between path components of $(G\times  G)\sslash  \Gc^{(2)}$ and $G\sslash \Gc$. Indeed,  both are in bijection with path-components of the function space $\TOP(S^1,B_{\J_2(\Gc)})$, by Corollary \ref{cor:CW-class}.

Consider the annulus $S^1 \times I$, considered with the CW-decomposition, $A_1^2$, shown in the rightmost diagram of Figure  \ref{eq:act-g}. Hence $S^1_\sk$ and $S^1_{\sk'}$ embed cellularly as $S^1\times \{1\}$ and $S^1\times \{0\}$, respectively. The profunctor, as per Definition \ref{def:crs_profunctor}, $$  \Hp^{(A_1^2 : S^1_\sk\times \{0\}, S^1_{\sk'}\times \{1\})}_{\J_2(\Gc))}\colon \pi_1(\CRS(\pi_1(S^1_\sk),\J_2(\Gc)) \bto  \pi_1(\CRS(\pi_1(S^1_{\sk'}),\J_2(\Gc)), $$
considered as a functor, $\Ap_1^2\colon \pi_1(G\sslash \Gc)^\op\times  \pi_1(G\times G)\sslash \Gc^{(2)})\to \Sets$, has the form
$$\Ap_1^2=\Ap\circ ( \id_{\pi_1(G\sslash \Gc)} \times P_1^2 ).$$

\subsubsection{Computations related to $\CRS\big(\Pi(I_\sk),\J_2(\Gc)\big)$}\label{sec:ICRS}
Let $I=[0,1]$, have $0$-cells at $\{0\}$ and $\{1\}$.  By Lemmas \ref{lem:maps_on_basis} and \ref{lem:derivation_basis}, crossed complex maps from $ \Pi(I_{\sk})$ to $\J_2(\Gc)$ are in one-to-one correspondence with elements of $G$. For each such map, $f$, $f$-homotopies correspond to elements of $G\times G\times E$, and each 2-fold $f$-homotopy is then given by an element of $E\times E$.

We will use the following lemma to write down $\CRS\big(\Pi(I_\sk),\J_2(\Gc)\big)$.
\begin{Lemma}
 There exists an action of $G\times (G \ltimes E)$ on $G$ of the form $(x,y,a)\bullet  g = xg \partial(a)y^{-1}$. Furthermore, given $g\in G$, we have a group  map $F'_g\colon E\times E \to G\times (G\ltimes E)$, with  $F'_g(a,b)=\big( \partial(a),\partial(b), a^{-1} \trl g \, b )$, satisfying  $F'_g(a,b)\bullet g=g$.
\end{Lemma}
Given this, the crossed complex $\CRS(\Pi(I_\sk),\Gc)$ is isomorphic to the crossed complex $\Gc^I$ written down below, which has set of objects $G$,
$$\Gc^I:=  \cdots\ra{\partial} \bigsqcup_{g \in G}\{1\} \ra{\partial}   \bigsqcup_{g \in G} E\times E \ra{\partial}  G \sslash\big( G \times (G \ltimes E)\big).$$
We leave it to the reader to write down the action of the groupoid $(\Gc^I)_1$ on $(\Gc^I)_2$.

Noting that $\CRS\big(\Pi(\{0,1\}), \J_2(\Gc)\big)\cong J_2(\Gc)\times J_2(\Gc)$, the map $P_{\partial}\colon \Gc^I \to \J_2(\Gc)\times \J_2(\Gc)$, induced by the inclusion of $\{0,1\}$ into $[0,1]$, is such that all objects $g\in G$ are sent to the unique object, $(*_L,*_R)$ of $\J_2(\Gc)\times\J_2(\Gc)$. Furthermore, on morphisms,
\begin{align*}
P_\partial\big (g \ra{(x,y,a)} (x,y,a)\bullet g \big)& =\big ( (*_L,*_R) \ra{(x,y)}  (*_L,*_R)\big),\\
 P_\partial \big(g \ra{(a,b)} g\big)&=\big ( (*_L,*_R) \ra{(a,b)} (*_L,*_R) \big).
\end{align*}

The fibre $P_\partial^{-1}(\{(*_L,*_R))$ is hence the crossed complex $\J_1(G\sslash_\tl E)$, where the action of $E$ on the underlying set of $G$, is $g\tl e:=g\partial(e)$.
In particular, if $g \in G$,  $$\chi^\pi\left( \PC_g\big(P_\partial^{-1}(*_L,*_R)\big)\right)=\frac{|\Orb_E(g)|}{E}=|\ker(\partial)|.$$

The fundamental groupoid of $\J_2(\Gc)$ is the group $G/\partial(E)$. The profunctor,
$$  \Hp^{([0,1] : \{0\}, \{1\})}_{\J_2(\Gc)}\colon \pi_1\big(\CRS\big(\Pi(\{0\}),\J_2(\Gc)\big)\big) \bto
  \pi_1\big(\CRS\big(\Pi(\{1\}),\J_2(\Gc)\big)\big), $$
 seen as a profunctor, $\Ip\colon G/\partial(E)\bto G/\partial(E)$, is given by the actions of the group $G/\partial(E)$ on its underlying set, by left and right multiplications, so by the identity profunctor in Example \ref{identity profunctor}.

\subsubsection{Computation of $\CRS\big(\Pi(D^2_\sk),\J_2(\Gc)\big)$}\label{sec:eps_eta}
Let $D^2$ have a CW-decomposition where $S^1$, with a unique $0$-cell at the south-pole, is embedded cellularly, with $\mathbf{g}$ being its unique 1-cell,  and with a unique 2-cell, $\mathbf{e}$ that is attached along $\mathbf{g}$. Again by Lemma \ref{lem:maps_on_basis}, crossed complex maps from $\Pi(D^2_\sk)$ to $\J_2(\Gc)$ are in bijection with pairs $(e,g)\in E \times G$, with $\partial(e)=g$, hence with elements of $E$. By Lemma \ref{lem:derivation_basis}, 1-fold homotopies are in one-to-one correspondence with elements of $G\times E$, and 2-fold homotopies with elements of $E$.

To build the crossed complex, $\CRS\big(\Pi(D^2_\sk),\J_2(\Gc)\big)$, first recall \S\ref{sec:G/GG},  mainly the notation in Lemma \ref{lem:exerG//G},  and note:

\begin{Lemma} We have a left-action of $G\ltimes E$ on $E$ where $ (g,e)\bullet a=(a e) \trl g^{-1}$. It satisfies $\partial\big((g,e)\bullet a\big)=(g,e)\bullet \partial(a)$ and $F_{\partial(a)}(n)\bullet a=a$, if $a,n\in E$.
\end{Lemma}

 Unpacking the construction in Subsection \ref{sec:xcomphom}, we have:
\begin{Lemma}
  The crossed complex $\CRS\big(\Pi(D^2_\sk),\J_2(\Gc)\big)$ is isomorphic to the following crossed complex, $E\sslash \Gc$, with set of objects $E$,
  $$E\sslash\Gc:=  \cdots\ra{\partial} \bigsqcup_{e \in E}\{1\} \ra{\partial}   \bigsqcup_{e \in E} E \ra{\partial}  E \sslash( G\ltimes E).$$
  Here:
  \begin{enumerate}[leftmargin=1cm]
  \item The groupoid map $\partial\colon \bigsqcup_{e \in E} E \to   E \sslash( G\ltimes E)$  sends $(a\ra{n}a)$ to $e\ra{F_{\partial(a)}(n)} e$.
   \item $\big( a \ra{n} a\big)\triangleleft\big(  a\ra{(h,e)} (h,e) \bullet  a \big)= (h,e) \bullet  a   \ra{n \triangleleft h^{-1}}  (h,e) \bullet  a$.
 \end{enumerate}
 Moreover, the crossed complex map $\CRS\big(\Pi(D^2_\sk),\J_2(\Gc)\big) \to \CRS\big(\Pi(S^1_\sk),\J_2(\Gc)\big)$ arising from the inclusion $\Pi(S^1_{\sk})\to \Pi(D^2_\sk)$ is given by $P_{S^1}\colon E\sslash \Gc \to G \sslash \Gc$, with $P_{S^1}(e)=\partial(e)$, for any object $e$ of $E\sslash \Gc$, and, on 1-morphisms and 2-morphisms,
 \begin{align*}
  P_{S^1}\big ( a\ra{(h,e)} (h,e) \bullet  a\big)&=\partial(a)\ra{(h,e)} (h,e) \bullet  \partial(a),\\
   P_{S^1}\big( a \ra{n} a\big)&= \partial(a) \ra{n} \partial(a).
 \end{align*}
 \end{Lemma}

The fibre of $P_{S^1}\colon E\sslash \Gc \to G\sslash \Gc$ at $g\in G$ is hence the crossed complex $\J_0\big(\partial^{-1}(g)\big)$, with only identify morphisms. In particular
$\chi^\pi(P_{S^1}^{-1}(g))=|\partial^{-1}(g)|$.

Noting that $\pi_1\big(\CRS\big(\emptyset,\J_2(\Gc)\big)\big)=\{*\}$, let us determine the profunctors,
\begin{align*}
\Hp^{(D^2;  \emptyset, S^1)}_{ \J_2(\Gc)}&\colon \pi_1\big(\CRS\big(\emptyset,\J_2(\Gc)\big)\big) \bto  \pi_1\big(\CRS\big(\Pi(S^1_\sk),\Gc\big)\big),\\
\Hp^{(D^2;   S^1, \emptyset)}_{ \J_2(\Gc)}&\colon  \pi_1\big(\CRS\big(\Pi(S^1_\sk),\Gc\big)\big)\bto \pi_1\big (\CRS\big(\emptyset,\J_2(\Gc)\big)\big),
\end{align*}
considered as a profunctors $\Bp\colon \{*\} \bto \pi_1(G\sslash \Gc)$ and $\Ep\colon \pi_1(G\sslash \Gc) \to \{*\}$, respectively.  These are given by the  functors $\Bp\colon \pi_1(G\sslash \Gc) \to \Sets$ and $\Ep\colon \pi_1(G\sslash \Gc)^\op\to \Sets$  such that, on objects,
$\Bp(g)= \partial^{-1}(g)=\Ep(g)$, and, on morphisms $$\Bp\big (g \ra{[(x,a)]_g} (x,a)\bullet g\big)=\left(\vcenter{\xymatrix@R=0pt@C=9pt{\partial^{-1}(g) \ar[r] & \partial^{-1}\big ( (x,a)\bullet g\big) \\
                      n\ar@{|->}[r] & (x,a)\bullet n         }} \right),$$
$$\Ep\big ( (y,b)^{-1}\bullet h \ra{[y,b]_{(y,b)^{-1}\bullet h}} h\big)=\left(\vcenter{\xymatrix@R=0pt@C=9pt{\partial^{-1}(h) \ar[r] & \partial^{-1}\big ( (y,b)^{-1}\bullet h \big)\\
                      m\ar@{|->}[r] & (y,b)^{-1}\bullet m         }} \right).$$

                      As a consequence, the profunctor $(\Ep\#_ 1 \Bp)\colon \pi_1(G\sslash \Gc) \bto \pi_1( G\sslash \Gc)$ is such that, on objects, $(h,g)\mapsto \partial^{-1}(h)\times \partial^{-1}(g)$.

                      The following lemma will be useful later. The notation $\overline{\Ap}$ is defined in \S \ref{sec:G/GG}.
\begin{Lemma}\label{lem:calc-nat-trans}
Given $(a,b)\in E\times E,$ let $L_{(a,b)}\colon \overline{\Ap}\big(\partial(a),\partial(b)\big) \to E$ be the map, of sets, such that
 $(p,c)\mapsto   (ac)  (b^{-1}\trl p)$, then (recall  $(p,c)\tl_{\partial(a)} e=(p,c) F_{\partial(a)}(e)$),
 \begin{enumerate}[leftmargin=1cm]
  \item $L_{(a,b)}$ takes values in $\ker(\partial)$,
  \item $L_{(a,b)}\big ( (p,c)\tl_{\partial(a)} e\big )=e L(a,b)\big ( (p,c) \big )e^{-1}$,
  \item $L_{\big  (h',f')^{-1}\bullet a,  (h,f)\bullet b\big)}\big ( (h,f) (p,c) (h',f') \big)=\big( L_{(a,b)}(h,f)\big) \trl h'$.
 \end{enumerate}
\end{Lemma}
\begin{proof}
 These are all simple computations. For instance, for the second item:
\begin{align*} L(a,b)\big ( (p,c) F_{\partial(a)}(e)\big )&=L_{a,b}\big(p\partial(e), a^{-1} e^{-1} a  c e\big)\\
                                                          &=   e^{-1} a c e \,\, \big (b^{-1}\trl (p \partial(e))\big)  \\
                                                          &=   e^{-1} a c\, (b^{-1}\trl p ) \, e .
                                                          \end{align*}
For the last item:
\begin{align*}
 L_{\big  (h',f')^{-1}\bullet a,  (h,f)\bullet b\big)}&\big ( (h,f) (p,c) (h',f') \big)\\&=L_{\big(  ( a \,( {f'}^{-1})\trl {h'}^{-1}) \trl h',  (bf)\trl {h}^{-1} \big)} (hph', f'\, c \trl h' \,\, f\trl (ph'))\\
 &= a  \trl h'\,\  c \trl h' \,\, f\trl (ph'))\ (bf)^{-1}\trl (ph')= \big( L_{(a,b)}(h,f)\big) \trl h'.
\end{align*}
The remaining details are left to the reader.
\end{proof}

We have now constructed two different set-valued profunctors $\pi_1(G\sslash \Gc)\bto \pi_1(G\sslash \Gc)$. Namely $\A\colon \pi_1(G\sslash \Gc)\bto \pi_1(G\sslash \Gc)$, written down in \S \ref{sec:G/GG}, which is associated with the two inclusions of $S^1$ into the annulus, $A=S^1\times I$, and  $(\Ep\#_ 1 \Bp)\colon \pi_1(G\sslash \Gc) \bto \pi_1( G\sslash \Gc)$, which is associated with the two inclusions of $S^1$ into $D^2\sqcup D^2$. Consider their linearisation to $\Vect$-profunctors, $\boldsymbol{\A},( \boldsymbol{\Ep}\#_ 0 \boldsymbol{\Bp})\colon \pi_1(G\sslash \Gc)\bto \pi_1(G\sslash \Gc)$, obtained by composition with the free vector space functor.

The previous lemma gives that we have a natural transformations of profunctors
\begin{multline*}{\eta}_{\hat{\mu}}\colon\big( \boldsymbol{\Ep}\#_ 0 \boldsymbol{\Bp}) \colon   \pi_1(G\sslash \Gc)^\op \times  \pi_1(G\sslash \Gc) \to \Vect\big)\\ \implies  \big(  (    \boldsymbol{\A}\colon  \pi_1(G\sslash \Gc)^\op \times  \pi_1(G\sslash \Gc) \to \Vect\big).\end{multline*}
Its matrix elements are, for $g,h$ in the image of $\partial\colon E \to G$, and $a,b\in E$, with $\partial(a)=g$ and $\partial(b)=h$, and $(p,c)$ with $(p,c)\bullet g=h$, as below:
$$\Big\langle g \ra{a\otimes b} h \Big | \eta_{\hat{\mu}} \Big| g \ra{[(p,c)]_g} h \Big\rangle =  \chi^\pi\big (P^{-1}(g,h)\big) \begin{cases}  1 , \textrm{ if } L_{(a,b)}(p,c)=1_E\\
                                                                    0, \textrm{ otherwise}                   \end{cases} $$
That this is well defined and gives a natural transformation of profunctors follows from the previous lemma. The notation  $\chi^\pi\big (P^{-1}(g,h)\big)$ is given in \eqref{eq:fib-A}. The reason for this latter normalising factor is that it agrees with the conventions for the extended Quinn TQFT. This will be explained later in \S \ref{sec:123TQFT_fingroup}

\subsubsection{Computation of $\CRS(\Pi(S^2_\sk),\Gc)$}\label{sec:S2}
Consider $S^2$ with a unique 0-cell, at the south-pole, and a unique 2-cell, $\mathbf{a}$. By Lemma \ref{lem:maps_on_basis}, or by Example \ref{Ex:simple_CW}, crossed complex maps, $f\colon\Pi(S^2_\sk)\to \J_2(\Gc)$, are in one-to-one correpondence with elements of $\ker(\partial)$, and, by Lemma \ref{lem:derivation_basis}, given any such map, $f$,  $f$-homotopies are in bijection with elements of $G$ and 2-fold $f$-homotopies with elements of $E$.

Since $\Gc$ is a crossed module,
the action $\trl$ of $G$ on $E$ restricts to an action of $G$ on $\ker(\partial)$, which descends to an action of the group $G/\partial(E)$ on $\ker(\partial)$. We can form the corresponding action groupoids, $\ker(\partial)\sslash G$ and $\ker(\partial)\sslash \big(G/\partial(E)\big)$.

The following then follows by simple computations:
\begin{Lemma}
  The crossed complex, $\CRS\big(\Pi(S^2_\sk),\J_2(\Gc)\big)$, is isomorphic to the following crossed complex, $\ker(\partial)\| \Gc$, with set of objects $\ker(\partial)$,
  $$\ker(\partial)\|\Gc:=  \cdots\ra{\partial} \bigsqcup_{a \in \ker(\partial)}\{1\} \ra{\partial}   \bigsqcup_{a \in \ker(\partial)} E \ra{\partial}  \ker(\partial) \sslash G.$$
  Here:
  \begin{enumerate}[leftmargin=1cm]
   \item The groupoid map $\partial\colon \bigsqcup_{a \in \ker(\partial)} E \to   \ker(\partial) \sslash  G$  sends $(a\ra{n}a)$ to $a \ra{\partial(n)} a$.
   \item $\big( a \ra{n} a\big)\triangleleft\big(  a\ra{(h,e)} (h,e) \bullet  a \big)= (h,e) \bullet  a   \ra{n \triangleleft h^{-1}}  (h,e) \bullet  a$.
 \end{enumerate}
 Furthermore, $\pi_1(\ker(\partial) \| \Gc)=\ker(\partial)\sslash (G/\partial(E))$.
\end{Lemma}

Note that by Lemma \ref{path:compcrs},  or by Lemma \ref{sec:homfinXcomp}, for $a \in \ker(\partial)$:
$$\chi^\pi\big (\PC_a  (\ker(\partial)\| \Gc) \big)=|\Orb_G(a)|\frac{|E|}{|G|}, \quad \textrm{ and} \quad \chi^\pi(\ker(E)\| \Gc)=\frac{|\ker(\partial) | |E|}{G}.$$

\subsubsection{The profunctor associated to $\Nu$} Consider the 3-manifold $\Nu$ in Figure \ref{fig:Nu}, so $\Nu$ is obtained by removing two disjoint open 3-balls from $D^3$. We have a cobordism $\big(\iota_C, \Nu,\langle \iota_L,\iota_R\rangle\big) \colon S^2  \longrightarrow S^2\sqcup S^2$,  where $\iota_C\colon S^2 \to \partial \Nu$ parametrises the `outside' component, $S^2_C$ of $\partial \Nu$,  and $\iota_L\colon S^2 \to \partial \Nu$ and $\iota_R\colon S^2 \to \partial \Nu$ parametrise the left and right `inside' components, $S^2_L$ and $S^2_R$, of $\partial \Nu$. Give $\Nu$ the CW decomposition in which $\iota_C,\iota_L,\iota_R\colon S^2 \to \Nu$ are all cellular maps, so we have three 0-cells, 2-cells $\mathbf{a}, \mathbf{a'}$ and $\mathbf{b}$, as shown, two 1-cells, $\mathbf{p}$ and $\mathbf{q}$, as shown, and one $3$-cell.
\begin{figure}[ht!]
 \labellist
\pinlabel $S^2$ at 105 113
\pinlabel $\xrightarrow{\iota_C}$ at 193 161
 \pinlabel ${\mathbf{a}}$ at 450 176
\pinlabel $\mathbf{a'}$ at 758 176
\pinlabel $\mathbf{b}$ at 610 205
\pinlabel $\mathbf{p}$ at 488 100
\pinlabel $\mathbf{q}$ at 722  100
\pinlabel ${\Nu}$ at 865 43
\pinlabel $S^2$ at 1210 113
\pinlabel $S^2$ at 1400 113
\pinlabel $\xleftarrow{\langle \iota_L,\iota_R\rangle }$ at 1018 161
\endlabellist
\centering
\includegraphics[scale=0.25]{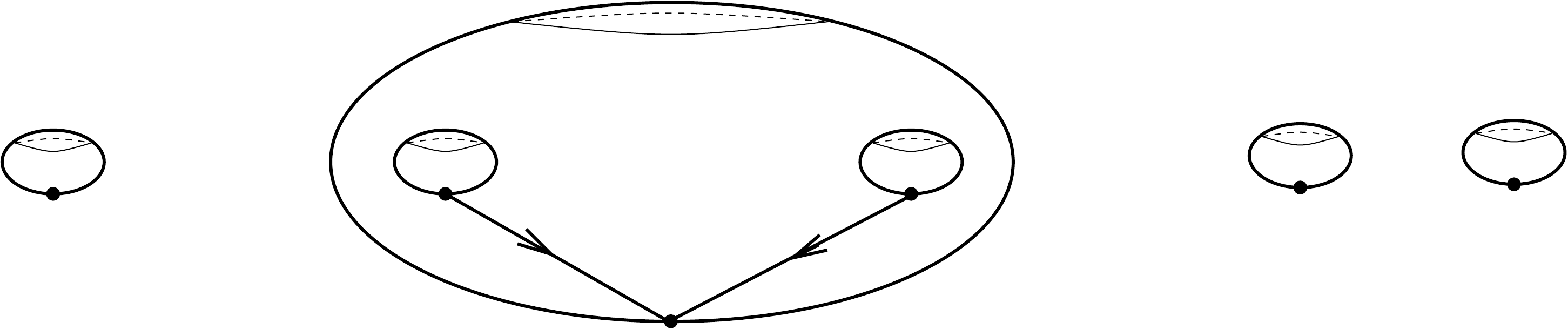}
\caption{The cobordism $(\iota_C, \Nu,\langle \iota_L,\iota_R\rangle) \colon S^2  \longrightarrow S^2\sqcup S^2$.\label{fig:Nu}}
\end{figure}

We wish to determine the set-valued profunctor, following Definition \ref{def:crs_profunctor},
\begin{multline}\label{eq:N}
\Hp^{(\Nu;  S^2_C, S^2_L\cup S^2_R)}_{ \J_2(\Gc)}\colon \pi_1\big(\CRS(\Pi(S^2_\sk),\J_2(\Gc)\big)\\ \bto  \pi_1\big(\CRS(\Pi(S^2_\sk),\J_2(\Gc)\big)\times \pi_1\big(\CRS(\Pi(S^2_\sk),\J_2(\Gc)\big).
\end{multline}
By Lemma \ref{lem:maps_on_basis}, maps $\Pi(\Nu_\sk)\to \J_2(\Gc)$ are in one-to-one correspondence with sequences $(b,a,a',p,q)\in \ker(\partial)^3 \times G^2$ with $b= (a \trl p ) \, ( a' \trl q)$. We have a map,
$$P^\Nu_\partial\colon \CRS\big(\Pi (\Nu_\sk),\J_2(\Gc)\big) \to \CRS\big(\Pi(S^2_\sk), \J_2(\Gc)\big){}^3\cong \big(\ker(\partial)\| \Gc\big){}^3,$$ induced by the inclusion $\langle \iota_C,\iota_L,\iota_R\rangle \colon S^2\sqcup S^2 \sqcup S^2 \to \Nu$. Given  $(b,a,a') \in \ker(\partial)^3$, let us determine the fibre, $({P^\Nu_\partial}){}^{-1}(b,a,a')$.
Define $$\overline{\Np}(b,a,a'):=\{ (p,q)\in G\times G: b= (a \trl p) \, ( a' \trl q)\}.$$ We have a left-action of $E\times E$ on $\overline{\Np}(b,a,a')$ with $(e,e')\bullet (p,q)=(\partial(e) p, \partial(e')q)$. It can easily be shown that:
$$\big(P^\Nu_\partial\big)^{-1}(b,a,a')=\J_1\left(\overline{\Np}(b,a,a')\sslash (E\times E)\right).$$
In particular,
$\chi^\pi\big(\big(P^\Nu_\partial\big)^{-1}(b,a,a')\big)=|\overline{\Np}(b,a,a')|/ |E|^2$.

Finally, put ${\Np}(b,a,a'):=\overline{\Np}(b,a,a')/(E\times E)$.
The profunctor, in \eqref{eq:N}, considered as a functor $$\Np\colon \pi_1( \ker(\partial) \| \Gc)^\op \times \pi_1( \ker(\partial) \| \Gc) \times \pi_1( \ker(\partial) \| \Gc) \to \Sets,$$ is such that, on objects, $(b,a,a')\mapsto \Np(b,a,a')$, where $b,a,a'\in \ker(\partial)$.  On morphisms, given $[y],[x],[x']\in G/\partial(E)$, we have:
\begin{multline*}
\Np\big ( b\trl{y^{-1}}\ra{[y]} b,  a\ra{[x]} a \trl{x},  a'\ra{[x']} a'' \trl{x'}\big)\\=\Bigg(\vcenter{\xymatrix@C=10pt@R=1pt{\Np(b,a,a')\ar[r] &\Np(b \trl{y^{-1}} ,a \trl{x},a'\trl{x'})\\ [(p,q)]\ar@{|->}[r] &(x^{-1}py^{-1},{x'}^{-1}qy^{-1})]}}\Bigg).
\end{multline*}

\subsubsection{Interlude on double groupoids}\label{Interlude-db}
In order to write down more efficiently some of the formulas ahead, let us recall \cite{brown_hha,brown_higgins_sivera,Martins_Picken:2011} that any crossed module of groups, $\Gc$, gives rise to an edge symmetric double groupoid, $\mathcal{D}(\Gc)$, with a single object, $*$.
The squares in $\mathcal{D}(\Gc)$ look like the  diagrams in the figure below,
\begin{equation}\label{eq:squares-G}
\vcenter{\xymatrix@C=13pt@R=13pt{
\ast \ar@{<-}[d]_y \ar@{<-}[r]^x &\ast\ar@{<-}[d]^w \\ \ast \xtwocell[ur]{}\omit{\omit \underline{e}}  \ar@{<-}[r]_z &\ast }}, \textrm{ with } \begin{matrix} &\quad\,\,\, e \in E,  \quad  x,y,z,w \in G,\\ \\ &\d(e)=x^{-1}w^{-1}zy.\end{matrix}
\end{equation}
The  vertical and horizontal compositions in $\mathcal{D}(\Gc)$ are as shown below:
 $$\vcenter{\xymatrix@C=19pt@R=19pt{
\ast \ar@{<-}[d]_y \ar@{<-}[r]^x &\ast\ar@{<-}[d]|{w^{}_{}} \ar@{<-}[r]^{x'} &\ast\ar@{<-}[d]^{w'} \\ \ast \xtwocell[ur]{}\omit{\omit \underline{e}}  \ar@{<-}[r]_z &\ast \xtwocell[ur]{}\omit{\omit \underline{e'}}  \ar@{<-}[r]_{z'} &\ast }}=\vcenter{\xymatrix@C=20pt@R=13pt{
\ast \ar@{<-}[d]_y \ar@{<-}[r]^{x'x} &\ast\ar@{<-}[d]^{w'} \\ \ast \xtwocell[ur]{}\omit{\omit \underline{(e'\trl x)\,\, e  }}  \ar@{<-}[r]_{z'z} &\ast }}, \qquad \vcenter{\xymatrix@C=19pt@R=19pt{
\ast \ar@{<-}[d]_y \ar[r]^x &\ast\ar@{<-}[d]^w \\ \ast \ar@{<-}[d]_{y'}\xtwocell[ur]{}\omit{\omit \underline{e}}  \ar[r]|z &\ast \ar@{<-}[d]^{w'}\\ \ast \xtwocell[ur]{}\omit{\omit \underline{e'}}  \ar[r]_{z'} &\ast  } } = \vcenter{\xymatrix@C=21pt@R=13pt{
\ast \ar@{<-}[d]_{y'y} \ar[r]^x &\ast\ar@{<-}[d]^{w'w} \\ \ast \xtwocell[ur]{}\omit{\omit\underline{ e\,\,( e'\trl y)}}  \ar[r]_{z'} &\ast }}$$
These are associative and  satisfy the interhange law. We abbreviate:
$$\vcenter{\xymatrix@C=13pt@R=13pt{
\ast \ar@{<-}[d]_{1_G} \ar@{<-}[r]^g &\ast\ar@{<-}[d]^{1_G} \\ \ast \xtwocell[ur]{}\omit{\omit \underline{e}}  \ar@{<-}[r]_{g\partial(e)} &\ast }}=\vcenter{\xymatrix@C=18pt@R=13pt{\ast\ar@{<-}@/^1pc/[r]^{g} \xtwocell[r]{}\omit{\omit \underline{e}} \ar@{<-}@/_1pc/[r]_{g\partial(e)} &\ast}},  \vcenter{\xymatrix@C=13pt@R=13pt{
\ast \ar@{<-}[d]_{1_G} \ar@{<-}[r]^g &\ast\ar@{<-}[d]^{1_G} \\ \ast \xtwocell[ur]{}\omit{\omit \underline{1_E}}  \ar@{<-}[r]_{g} &\ast }}=\vcenter{\xymatrix@C=13pt@R=13pt{\ast\ar@{<-}@/^0pc/[r]^{g} &\ast }},
\vcenter{\xymatrix@C=13pt@R=13pt{
\ast \ar@{<-}[d]_{h\partial(e)} \ar@{<-}[r]^{1_G} &\ast\ar@{<-}[d]^{h} \\ \ast \xtwocell[ur]{}\omit{\omit \underline{e}}  \ar@{<-}[r]_{1_G} &\ast }}=
 \vcenter{\xymatrix@R=13pt{\ast
 \ar@{<-}@/_0.7pc/[d]_{h\partial(e)}
 \ar@{<-}@/^0.7pc/[d]^{h}
  \xtwocell[d]{}\omit{\omit \underline{e}}
  \\
  \ast}}
$$

The $\bullet$ action, of $G\ltimes E$ on $G$, in Lemma \ref{lem:exerG//G}, arises from the fact that we have squares as below,
so $G\sslash (G\ltimes E)$  is isomorphic to the  vertical groupoid of $\mathcal{D}(\Gc)$,
$$\vcenter{\xymatrix@C=35pt@R=10pt{
\ast \ar@{<-}[d]_h \ar@{<-}[r]^g &\ast\ar@{<-}[d]^h \\ \ast \xtwocell[ur]{}\omit{\omit \underline{e}}  \ar@{<-}[r]_{(h,e)\bullet g} &\ast }} \qquad \textrm{ and }  \qquad \vcenter{\xymatrix@C=35pt@R=10pt{
\ast \ar@{<-}[d]_h \ar@{<-}[r]^{(h,e)^{-1} \bullet g} &\ast\ar@{<-}[d]^h \\ \ast \xtwocell[ur]{}\omit{\omit \underline{e}}  \ar@{<-}[r]_{ g} &\ast }}$$ The equivalence relation in $G\sslash (G\ltimes E)$ giving $\pi_1(G\sslash \Gc)$, in \eqref{eq:cong-pi1}, is, in this language,
\begin{equation}\label{eq:cong-pi1-DG}
\vcenter{\xymatrix@C=20pt@R=10pt{
\ast \ar@{<-}[d]_h \ar@{<-}[r]^g &\ast\ar@{<-}[d]^h \\ \ast \xtwocell[ur]{}\omit{\omit \underline{e}}  \ar@{<-}[r]_{(h,e)\bullet g} &\ast }} \sim
\vcenter{\xymatrix@C=40pt@R=10pt{
\ast \ar@{<-}[d]^h
\ar@{<-}@/_1.6pc/[d]_{h\partial(a)}
\xtwocell[d]{}\omit{\omit \underline{a} \,\,\,\,\,\,\,\,\,\,\,\,\,}
\ar@{<-}[r]^g &\ast\ar@{<-}[d]_h \ar@{<-}@/^1.6pc/[d]^{h\partial(a)}
\xtwocell[d]{}\omit{\omit \,\,\,\,\,\,\,\,\,\,\,\underline{a^{-1}} }
 \\ \ast \xtwocell[ur]{}\omit{\omit \underline{e}}  \ar@{<-}[r]_{(h,e)\bullet g} &\ast }} = \vcenter{\xymatrix@C=35pt@R=10pt{
\ast \ar@{<-}[d]_{h\partial(a)} \ar@{<-}[r]^g &\ast\ar@{<-}[d]^{h\partial(a)} \\ \ast \xtwocell[ur]{}\omit{\omit \underline{(a^{-1} \trl g)\, e \, a}}  \ar@{<-}[r]_{(h,e)\bullet g} &\ast }}
\end{equation}

Given $(h,g,g')\in G^3$, define $\overline{\Mp}(g,g',h)$ as the set of squares in $\mathcal{D}(\Gc)$ of form:
$$\vcenter{\xymatrix@C=14pt@R=13pt{ \ast \xtwocell[rrrrd]{}\omit{\omit \underline{e}}
   \ar@{<-}[rrrr]^h &&&&
\ast\ar@{<-}[d]^{q} \\
                 \ast  \ar[u]^{p}  \ar@{<-}[r]_{g}  &  \ast  \ar@{<-}[r]_{p^{-1}} &  \ast  \ar@{<-}[r]_{q}  &  \ast  \ar@{<-}[r]_{g'} & \ast,   }}, \quad \textrm{ where} \vcenter{\xymatrix@R=1pt{ e \in E, \quad p,q \in G;\\ \qquad \d(e)=h^{-1}q^{-1}g'qp^{-1} gp.}} $$

We have a left action of $E^\op\times E^\op$ on $\overline{\Mp}(g,g',h')$, as shown below
\begin{align*}
(x,y)\bullet \vcenter{\xymatrix@C=10pt@R=13pt{ \ast \xtwocell[rrrrd]{}\omit{\omit \underline{e}}
   \ar@{<-}[rrrr]|h &&&&
\ast\ar@{<-}[d]_{q} \\
                 \ast  \ar[u]_{p}  \ar@{<-}[r]_{g}  &  \ast  \ar@{<-}[r]_{p^{-1}} &  \ast  \ar@{<-}[r]_{q}  &  \ast  \ar@{<-}[r]_{g'} & \ast   }}   = \hskip-2mm\vcenter{\xymatrix@C=32.5pt{ \ast 
                 \ar@/_1pc/@{{}{ }{}}[rrrr]|{\underline{e}}
   \ar@{<-}[rrrr]|h &&&&\ast \\
                 \ast  \ar@/_1.5pc/[u]_{p} \ar@/^0.0pc/[u]^{p\partial(x)}\xtwocell[u]{}\omit{\omit\,\,\,\,\,\,\,\,\,\,\,\,\, \underline{x}}  \ar@{<-}[r]_{g}  &  \ast \xtwocell[r]{}\omit{\omit \underline{x^{-1} \trl p^{-1}}} \ar@/^0.8pc/@{<-}[r]^{p^{-1}} \ar@/_1.1pc/@{<-}[r]_{\d(x)^{-1} p^{-1}} &  \ast \xtwocell[r]{}\omit{\omit \underline{y}} \ar@/^0.9pc/@{<-}[r]^{q} \ar@/_1.1pc/@{<-}[r]_{q\partial(y)} &  \ast  \ar@{<-}[r]_{g'} & \ast \ar@/^1.6pc/[u]^{q} \ar@/_0.0pc/[u]_{q\d(y)} \xtwocell[u]{}\omit{\omit \underline{y^{-1}} \,\,\,\,\,\,\,\,\,\,\,\,\,\,}   }} \\=
                 \vcenter{\xymatrix@C=39pt{ \ast \xtwocell[rrrrd]{}\omit{\omit \underline{(y^{-1}\trl h)\, e\, x\, \big( (yx^{-1})\trl (p^{-1} g)\big)} }
   \ar@{<-}[rrrr]|h &&&&
\ast\ar@{<-}[d]^{q\partial(y)} \\
                 \ast  \ar[u]^{p\partial(x)}  \ar@{<-}[r]_{g}  &  \ast  \ar@{<-}[r]_{\d(x)^{-1}p^{-1}} &  \ast  \ar@{<-}[r]_{q\d(y)}  &  \ast  \ar@{<-}[r]_{g'} & \ast   }}
                 \end{align*}
That $\bullet$ is a left action, follows from that $\mathcal{D}(\Gc)$ is a double groupoid, and so we have associative horizontal and vertical compositions, satisfying the interchange law.

We have a set-valued profunctor $\overline{\Mp}\colon G\sslash (G\ltimes E) \bto G\sslash (G\ltimes E) \times G\sslash (G\ltimes E)$. It is such that $(h,g,g')\mapsto \overline{\Mp}(h,g,g')$, and given arrows in $G \sslash (G\ltimes E)$, as below
\begin{equation}\label{eq:arrows}
(\alpha,a)^{-1}\bullet h\ra{(\alpha,a)} h, \qquad  g \ra{(\beta,b)} (\beta,b)\bullet g, \qquad g' \ra{(\beta',b')} (\beta',b')\bullet g',
\end{equation}
the associated map $\overline{\Mp}(h,g,g')\to \overline{\Mp}\big ((\alpha,a)^{-1}\bullet h, (\beta,b)\bullet g, (\beta',b')\bullet g')$ is as below,
\begin{multline}\label{eqM}
\vcenter{\xymatrix@C=14pt{ \ast \xtwocell[rrrrd]{}\omit{\omit \underline{e}}
   \ar@{<-}[rrrr]^h &&&&
\ast\ar@{<-}[d]^{q} \\
                 \ast  \ar[u]^{p}  \ar@{<-}[r]_{g}  &  \ast  \ar@{<-}[r]_{p^{-1}} &  \ast  \ar@{<-}[r]_{q}  &  \ast  \ar@{<-}[r]_{g'} & \ast   }}\longmapsto
\vcenter{\xymatrix@C=30pt{ \ast\ar@{<-}[d]_{\alpha} \xtwocell[rrrrd]{}\omit{\omit \underline{a}}
   \ar@{<-}[rrrr]^{(\alpha,a)^{-1}\bullet h} &&&&
\ast\ar@{<-}[d]^{\alpha}\\ \ast \xtwocell[rrrrd]{}\omit{\omit \underline{e}}
   \ar@{<-}[rrrr]|h &&&&
\ast\ar@{<-}[d]^{q} \\
                 \ast  \ar[u]^{p}  \ar@{<-}[r]|{g} \xtwocell[rd]{}\omit{\omit \underline{b}}  &  \ast  \ar@{<-}[r]|{p^{-1}} &  \ast  \ar@{<-}[r]|{q}  &  \ast \xtwocell[rd]{}\omit{\omit \underline{b'}}  \ar@{<-}[r]|{g'} & \ast \\ \ast  \ar[u]^{\beta}  \ar@{<-}[r]_{(\beta,b) \bullet g} &  \ast  \ar[u]_{\beta\,\,\,}  \ar@{<-}[r]_{\alpha p^{-1}\beta^{-1}} &  \ast  \ar@{<-}[r]_{\beta' q \alpha^{-1}}  \ar[u]^{\,\,\,\,\,\,\alpha}  &  \ast  \ar@{<-}[r]_{(\beta',b')\bullet g'}   \ar[u]^{\,\,\,\beta'} & \ast     \ar[u]_{\beta'}}}\\ = \vcenter{\xymatrix@C=30pt{
\ast \ar@{<-}[d]_{\beta p \alpha} \ar@{<-}[rrrr]^{(a,\alpha)^{-1}\bullet h} &&& &\ast\ar@{<-}[d]^{\beta' q \alpha}  \\ \ast \xtwocell[urrrr]{}\omit{\omit \underline{a \,\, ( e \trl \alpha)\,\, \big( b'\trl (qp^{-1}gp\alpha) \big) \,\,\big( b\trl (p\alpha)\big)}}
 \ar@{<-}[r]_{(\beta,b) \bullet g} &  \ast   \ar@{<-}[r]_{\alpha p^{-1}\beta^{-1}} &  \ast  \ar@{<-}[r]_{\beta' q \alpha^{-1}}    &  \ast  \ar@{<-}[r]_{(\beta',b')\bullet g'}    &\ast
 }}
                 \end{multline}

                 Finally we can define a profunctor
$
{\Mp}\colon \pi_1(G\sslash\Gc )\bto \pi_1(G\sslash\Gc)\times \pi_1(G\sslash\Gc)$. On objects, we put $\Mp(h,g,g'):=\overline{\Mp}(h,g,g')/(E^\op\times E^\op)$.
Because $\mathcal{D}(\Gc)$ is a double groupoid,  the map in \eqref{eqM} descends to a map $$\Mp(h,g,g')\to \Mp\big ((\alpha,a)^{-1}\bullet h, (\beta,b)\bullet g, (\beta',b')\bullet g'),$$  depending only on the equivalence class of the arrows in \eqref{eq:arrows}, under the equivalence relation in \eqref{eq:cong-pi1} and \eqref{eq:cong-pi1-DG}. This gives the value of $\Mp$ on morphisms.
\subsubsection{The profunctor associated to $\CRS\big(\Pi(\Mu_\sk),\J_2(\Gc)\big)$}\label{sec:pants}
Let  $\Mu$ be the `pair of pants' manifold, with the CW-decomposition shown in Figure \ref{m}, so with three 0-cells, five 1-cells, $\mathbf{h}$, $\mathbf{g}$, $\mathbf{g'}$, $\mathbf{p}$, and $\mathbf{q}$,  and one 2-cell, $\mathbf{e}$, as shown.
We have inclusions $\iota_L\colon S^1 \to \partial \Mu$, $\iota_R\colon S^1 \to \partial \Mu$, and also $\iota_C\colon S^1 \to \Mu$, parametrising the inside left, $S^1_L$, inside right, $S^1_R$ and outside, $S^1_C$, components of $\partial \Mu$, which are cellular maps. Furthermore, we have 
(1,2)-cobordisms $\big(\langle \iota_L,\iota_R\rangle,\Mu,\iota_C\big)\colon  S^1\sqcup S^1 \to S^1$ and $\big(\iota_C,\Mu,\langle \iota_L,\iota_R\rangle\big)\colon  S^1\to S^1  \sqcup  S^1$.

\begin{figure}[ht!]
 \labellist
\pinlabel ${\mathbf{h}}$ at 75 139
\pinlabel $S^1$ at 108 43
\pinlabel $\xrightarrow{\iota_C}$ at 183 91
 \pinlabel ${\mathbf{g}}$ at 451 98
\pinlabel $\mathbf{g'}$ at 760 98
\pinlabel $\mathbf{h}$ at 610 200
\pinlabel $\mathbf{p}$ at 492 59
\pinlabel $\mathbf{q}$ at 718  59
\pinlabel $\mathbf{e}$ at 587 62
\pinlabel ${\Mu}$ at 929 43
\pinlabel $\mathbf{g}$ at 1141 135
\pinlabel $S^1$ at 1221 43
\pinlabel $S^1$ at 1441 43
\pinlabel $\mathbf{g'}$ at 1350 139
\pinlabel $\xleftarrow{\langle \iota_L,\iota_R\rangle }$ at 1028 91
\endlabellist
\centering
\includegraphics[scale=0.24]{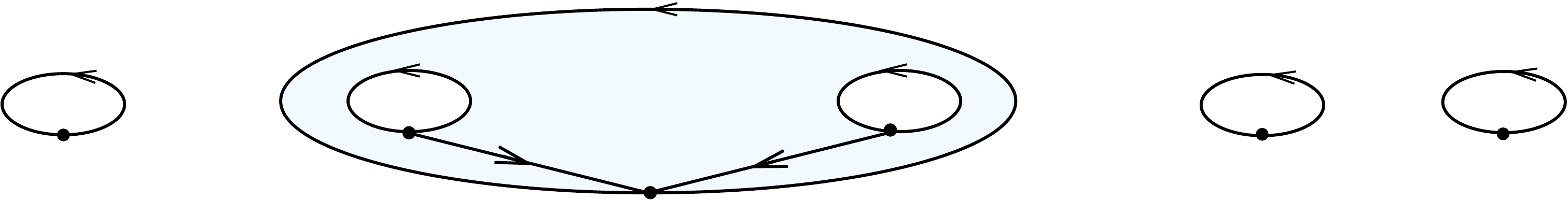}
\caption{{A CW-decomposition of the \emph{pair of pants manifold}, $\Mu$. \label{m}}}
\end{figure}

By Lemma \ref{lem:maps_on_basis}, crossed complex maps, $\Pi(\Mu_\sk)\to \J_2(\Gc)$, are in canonical bijection with 6-tuples, $(h,g,g',p,q,e)\in G^5 \times E$, with $\partial(e)=h^{-1}q^{-1}g'q p^{-1} g p$.
Given a triple $(h,g,g')\in G^3$, the set of objects of $(G \sslash \Gc)^3$, we  identify it with the unique map $\Pi_1(S^1\sqcup S^1 \sqcup S^1)_\sk\to \J_2(\Gc)$, such that $\mathbf{h}\mapsto h$, $\mathbf{g}\mapsto g$ and  $\mathbf{g'}\mapsto g'$.

We will compute the fibre, $\CRS^{(h, g,g')}(\Pi(\Mu_\sk),\iota(\Gc))$, in the pullback below:
\begin{equation}
\vcenter{\xymatrix@R=19pt{
 \CRS^{ (h, g,g')}(\Pi(\Mu_\sk),\iota(\Gc))\ar@{}[dr]|\pullback \ar[d] \ar[rr]^{\inc}  && \CRS(\Pi(\Mu_\sk),\iota(\Gc))  \ar[d]^{ \langle \iota_C,\iota_L, \iota_R\rangle ^*} \\
         \{(h,g,g')\}  \ar[rr]_{\inc} &&\Gc\times \Gc \times \Gc.
} }
\end{equation}
From Lemma \ref{lem:derivation_basis}, we have $\CRS^{(h, g,g')}(\Pi(\Mu_\sk),\iota(\Gc))\cong \overline{\Mp}(h,g,g')\sslash (E^\op \times E^\op)$.
In particular, we can see that:
\begin{align*}
\chi^\pi\big ( \CRS^{ (h, g,g')}(\Pi(\Mu_\sk),&\iota(\Gc))\big)=|\overline{\Mp}(h,g,g)/|E|^2\\&=\big\{ (p,q,e) \in G^2 \times E: \partial(e)=h^{-1}q^{-1}g'q p^{-1} g p\}/|E|^2.
\end{align*}

The discussion just given can be expanded to give that  we furthermore have a natural isomorphism of profunctors, using the notation in Definition \ref{def:crs_profunctor},
\begin{multline}
\Hp^{(\Mu:   S^1_C, S^1_L \sqcup S^1_R)}_{\J_2(\Gc)}\colon \pi_1(\CRS(\Pi(S^1_\sk),\J_2(\Gc))  \bto \pi_1(\CRS(\Pi(S^1_\sk),\J_2(\Gc))^2 \\   \cong {\Mp}\colon \pi_1(G \sslash \Gc )\bto \pi_1 (G \sslash \Gc )^2.
\end{multline}
\subsubsection{The torus}\label{sec:torus}

Consider  the torus $T^2$ with the standard CW-decomposition with a single 0-cell, two 1-cells, $\mathbf{g}$ and $\mathbf{h}$, and one two cell $\mathbf{e}$. Crossed complex maps $\Pi(T^2_\sk) \to \J_2(\Gc)$ are in bijection with elements of the set $$T^2_0(\Gc):=\{ (g,h,e)\in G^2\times E \mid \partial(e)=g^{-1}h^{-1}gh\}.$$ These can be identified with squares of $\mathcal{D}(\Gc)$ as in \eqref{eq:squares-G}, with $x,z = g$ and $y,w=h$. Consider the semidirect product $G\ltimes (E\times E)$, using the diagonal action of $G$. Proceeding as in the previous examples, we have an action of $G\ltimes (E\times E)$ on $T^2_0(\Gc),$ such that:
$$(x,a,b) \bullet
\vcenter{\xymatrix@C=13pt@R=13pt{
\ast \ar@{<-}[d]_h \ar@{<-}[r]^g &\ast\ar@{<-}[d]^h \\ \ast\ar@{}[ur]|{\underline{e}}   \ar@{<-}[r]_g &\ast }} =
\vcenter{\xymatrix@C=49pt@R=13pt{ \ast\ar@{<-}[r]^{1_G} \ar@{<-}[d]_{1_G}\ar@{}[dr]|{\underline{1_E}} & \ar@{<-}[d]_{\,\,\,\,\,\,x^{-1}}\ar@{}[dr]|{\underline{b^{-1}\trl x^{-1}}} \ar@{<-}[r]^{(x,b)\bullet g}\ast &\ast \ar@{<-}[d]^{x^{-1}\,\,\,\,\,\,}  &\ast \ar[l]_{1_G} \ar@{<-}[d]^{1_G}\ar@{}[dl]|{\underline{1_E}} \\
\ast \ar@{<-}[d]_{(x,a)\bullet h}\ar@{<-}[r]|{x^{-1}} & \ast\ar@{<-}[r]|g \ar@{<-}[d]_{\,\,\,\,\,\,h} & \ast \ar@{<-}[d]^{h\,\,\,\,\,\,}  \ar@{<-}[r]|x &\ast\ar@{<-}[d]^{(x,a)\bullet h} \\ \ar@{<-}[r]|{x^{-1}} \ar@{}[ur]|{\underline{a\trl x^{-1}}}
\ast & \ast \ar@{<-}[r]|g \ar@{<-}[d]_{\,\,\,\,\,\,x}\ar@{}[ur]|{\underline{e}} &\ast\ar@{<-}[r]|x \ar@{<-}[d]^{x\,\,\,\,\,\,} \ar@{}[ur]|{\underline{a^{-1}}} &\ast \\ \ast\ar@{<-}[r]_{1_G} \ar[u]^{1_G}\ar@{}[ur]|{\underline{1_E}}  &\ar@{<-}[r]_{(x,b)\bullet g}\ar@{}[ur]|{\underline{b}} \ast & \ast & \ast  \ar@{->}[l]^{1_G} \ar[u]_{1_G}\ar@{}[ul]|{\underline{1_E}}  }}  $$

Similarly to the previous examples, we can prove that:
\begin{Lemma}We have an isomorphism $\CRS\big(\Pi(T^2_\sk), \J_2(\Gc))\cong \T^2(\Gc)$, where
$$\T^2(\Gc):= \cdots \to \bigsqcup_{(g,h,e)\in T^2_0(\Gc)} \{1\}\ra{\partial} \bigsqcup_{(g,h,e)\in T^2_0(\Gc)} E \ra{\partial} T^2_0(\Gc)\sslash (G\ltimes (E\times E)) .$$
Here:
 $$\partial\big( (g,h,e) \ra{c} (g,h,e) \big )=(g,h,e)\ra{\big(\partial(c), F_g(c), F_h(c)\big) } (g,h,e), \textrm{ where } F_g(c)=c^{-1} \trl g \, c,$$
 \begin{multline*}
 \big( (g,h,e) \ra{c} (g,h,e) \big )\trl \big ( (g,h,e)\ra{(x,a,b)} (x,a,b)\bullet (g,h,e) \big)\\= (x,a,b)\bullet (g,h,e)  \ra{c\trl x^{-1}} (x,a,b)\bullet (g,h,e).
\end{multline*}
\end{Lemma}
The set of objects of the fundamental groupoid $\pi_1(T^2(\Gc))$ is  $T^2_0(\Gc)$. Morphisms   have the form ${[(x,a,b)]} \colon (g,h,e)\to(x,a,b)\bullet (g,h,e)$, where given $c \in E$, we have $(x,a,b)\sim  (x,a,b) \big(\partial(c), F_g(c), F_h(c)\big)$,  considering the product in $G\ltimes (E\times E)$.


\section{Once-extended TQFTs derived from  finite crossed complexes}\label{sec:TQFTS_xcomp}
\sectionmark{Once-extended TQFTs from finite crossed complexes}


As before, let $n$ be a non-negative integer and let $\A$ be a homotopy finite (often finite and reduced)
crossed complex, so its classifying space, $B_\A$, is a  homotopy finite space.  As usual, we are working over a subfield, $\kappa$, of $\C$, as we have to be able to invert non-zero integers when working with the homotopy content of spaces.

In this section, we use the techniques of the homotopy theory of crossed complexes, that were recalled and slightly refined in the previous section, to give  formulae for:
\begin{itemize}[leftmargin=0.6cm]
\item Quinn's finite total homotopy TQFT,  $\FQp{B_\A}{s}\colon \cob{n} \to \Vect_\C,$  in Definition  \ref{def:quinnTQFT}, for which we will give explicit formulae in Subsection \ref{Quinn for finite A};
\item[]\hspace{-7mm} and
\item the finitary once-extended Quinn TQFT, $\tFQd{B_\A}\colon \colon \tdcob{n}{B_\A} \to  \vProfGrpfin,$  in Definition \ref{def:finitary_ext_TQFT},
which we will treat in \S \ref{fin_Quinn_xcomp}.
 \item[]\hspace{-7mm}
 This will lead to  formulae for
\item the  Morita-valued once-extended Quinn TQFT, $\tFQmor{B_\A}\colon \tdcob{n}{B_\A} \to \Mor$,
 from Definition \ref{Sec:ext:Mor},
which is discussed in \S\ref{sec:mor-x-complexes}.
\end{itemize}

We will often consider our smooth manifolds, $\Sigma$, to be provided with  what we call \emph{simplicial stratifications}, $\zeta_{\Sigma} \colon |X_\Sigma| \to \Sigma$, where $X_\Sigma$ is a finite simplicial set, with geometric realisation $|X_\Sigma|$, and $\zeta_{\Sigma}$ is a homeomorphism. Most formulae for the Quinn and once-extended Quinn TQFTs will be given in terms of such simplicial stratifications, of closed smooth manifolds, cobordisms and extended cobordism.

Picking $n$-manifolds equipped with simplicial stratifications, leads naturally to another variant, ${\trcob{n}}$, of the bicategory $\tcob{n}$, whose objects are pairs, $(\Sigma, \zeta_{\Sigma} )$, where $ \zeta_{\Sigma} \colon |X_\Sigma| \to \Sigma$ is a simplicial stratification of the  $n$-manifold $\Sigma$, and with the rest of the bicategory structure induced, in the obvious way, from that of $\tcob{n}$. (In particular, cobordisms and extended cobordisms do not come with chosen simplicial stratifications.)

 In \S \ref{eTQFT-colA} and \S\ref{sec:mor-x-complexes}, we show the construction of two symmetric monoidal bifunctors,
 \begin{align}\label{eq:not-tr}
\tFQtr{\A} &\colon {\trcob{n}} \to \vProfGrpfin, \qquad \textrm{ and } \qquad
\tFQmortr{\A} \colon {\trcob{n}} \to \Mor,
\end{align}closely related to $\tFQmor{B_\A}$ and  $\tFQd{B_\A}$,  where $\A$ is a finite crossed complex.
 These latter constructions will allow us to put the non triangulation-invariant bits of the formulae we give for  $\tFQmor{B_\A}$ and  $\tFQd{B_\A}$ on a categorical footing.

The symmetric monoidal bifunctors in Equation \eqref{eq:not-tr} do not attach a value to an $n$-manifold, $\Sigma$, unless it is equipped with  a simplicial stratification, even though the associated groupoids and algebras are unique, up to a canonical invertible profunctor / Morita equivalence. In order to approach the literature on the topic of once-extended TQFTs, such as
\cite{Schommer-Pries,Bartlett_etal,Bartlett_Goosen}, in \S \ref{abs.TQFT}, we will address how to get rid of this latter dependence on the simplicial stratification of the $n$-dimensional manifolds. This step,  however, is non-canonical, and requires the use of the axiom of choice for classes\footnote{{The full force of the choice axiom is not required when the domain bicategory of a once-extended TQFTs is restricted to a `finitary' sub-bicategory of $\tcob{n}$,  for instance those that arise from the finite presentations of the  bicategories $\tcob{0}$ and $\tcob{1}$ developed in \cite{Bartlett_etal,Bartlett_Goosen,Schommer-Pries}.}}. This gives rise to once-extended TQFTs,
\begin{align}\label{eq:not-ch}
\tFQch{\A} &\colon {\tcob{n}} \to \vProfGrpfin, \qquad \textrm{ and } \qquad
 \tFQmorch{\A} \colon {\tcob{n}} \to  \Mor.
\end{align}

The well known $(1,2,3)$-extended TQFT sending $S^1$ to the quantum double of the group algebra of a finite group, \cite{Bartlett_etal,morton:cohomological:2015,Qiu_Wang,MNS}, is an example of this latter construction, and so is the $(0,1,2)$-extended TQFT arising  from the fact that the groupoid algebra of a groupoid $G$ is a `separable symmetric $*$-Frobenius algebra'; see \cite[\S 3.8]{Schommer-Pries} and \cite[Example 5.2]{LaudaPfeiffer}. This paper therefore in particular gives a topological interpretation for these two once-extended TQFTs. However, our construction is considerably more general.

As recalled in the beginning of Chapter \ref{Quinn Calc},  homotopy finite crossed complexes do not model the homotopy types of all homotopy finite spaces $\Bc$. The constructions in this section will not, therefore,  give formulae for all possible Quinn TQFTs, and once-extended Quinn TQFTs. They do, however, provide formulae for those  derived from, for instance,  finite 2-types $\Bc$, most relevant for higher gauge theory; cf., for instance, Baez and Schreiber, \cite{Baez-Schreiber:Streetfest:2007,Baez-Schreiber:2004}, Baez and  Huerta, \cite{Baez_Huerta:2011} or Faria Martins and Picken, \cite{Martins_Picken:2011}, where the links between 2-groups / crossed modules and higher gauge theory are summarised.

The constructions in this section of the paper moreover give a homotopy theoretical explanation for the `tube algebras', considered in \cite{Bullivant-thesis,Bullivant_Tube,Bullivant_Excitations}, in the context of excitations of topological phases. They also prove, as expected from \textit{loc cit.} that those tube algebras can be the starting point for extended TQFTs.

\subsection{Conventions and nomenclature}
In this section, given CGWH spaces, $M$ and $N$,   the space of functions from $M$ to $N$  (with the CGWH topology) will be  denoted  both by  $N^M$ and $\TOP(M,N)$, whichever is more convenient for the formula in question; the overall conventions are otherwise as  in \S\ref{sec:conventions_top}.
If $X$ is a CGWH space, and $x \in X$, then $\PC_x(X)$ denotes the path-component of $x$ in $X$, with the induced CGWH topology, and  $\hpiz(X)$ denotes the set of those $k$-ified path-components of $X$.
 If $f\colon M \to N$ is a map, then $\PC_f(\TOP(M,N))$ will therefore be the space of functions from $M$ to $N$ that are  homotopic to $f$.
\begin{Definition}[Finite Simplicial set]
 A simplicial set, $X$, is called \emph{finite} if it only has a finite number of non-degenerate simplices.
\end{Definition}
If a simplicial set $X$ is finite, then, as recalled earlier,  its geometric realisation, $|X|$, is naturally a finite CW-complex, with one $i$-cell for each non-degenerate $i$-simplex of $X$. Moreover, $|X|$ is a special CW-complex, in the sense of  Definition \ref{def:specialCW}. We also have a relative notion in which $(X,Y)$ is a pair of simplicial sets, meaning that $Y$ is a sub-simplicial set of  $X$,  and then $|Y|$ is naturally a subcomplex of $|X|$.

\begin{Notation}Let $(X,Y)$ be a pair of finite simplicial sets. Following  Notation \ref{not:K(X,n)},  we note that $\L(i,|X|)$ is the number of non-degenerate $i$-simplices of $X$, and $\L(i,|X|,|Y|)$,  the  number of non-degenerate $i$-simplices of $X$ that are not in $Y$.  We will extend the use of the notation, removing the geometric realisations signs for  convenience, so that, from now on:
\begin{itemize}[leftmargin=1cm]
\item $\L(i,X)$ will denote  the number of non-degenerate $i$-simplices of $X$,
\item[]\hspace{-12mm}and
\item $\L(i,X,Y)$ will denote  the  number of those non-degenerate $i$-simplices of $X$, that are not in $Y$.
\end{itemize}

\end{Notation}

Classically,  for instance, see  \cite[pp 107]{Hatcher},  an (abstract) \emph{simplicial complex}, $K$, is defined to be given by a sequence, $K=(K_0,K_1,\dots)$, consisting of  a set, $K_0$, the set of \emph{vertices of $K$}, together with, for each  $i\in \mathbb{Z}^+$, a subset, $K_i$, of the set of subsets of $K_0$ that have cardinality $i$. The elements of $K_i$ are called the \emph{$i$-faces} (or \emph{vertices} if $i=0$) of $K$. By definition, these are to have the property {that if $F$ is  a subset of cardinality $j$ of some $i$-face of $K$,} then $F$ will be  itself a $j$-face of $K$.

\begin{Definition}
A \emph{triangulation} of a  manifold, $M$, is  a homeomorphism, $f\colon |K| \to M$, where $K$ is a simplicial complex and $|K|$ is its geometric realisation.
\end{Definition}

\noindent Note that our definition of a triangulation of $M$ makes no reference to the smooth structure in $M$. This is not required as our construction of TQFTs only makes use of the underlying topological manifold of $M$. Otherwise we would need to consider smooth and regular triangulations of $M$ as, for example, in \cite[Chapter II]{MunkresDiff}.

If a simplicial complex,  $K=(K_0,K_1,\dots)$, is additionally provided with a total order on the set, $K_0$, of its vertices, then $K$ gives rise to a simplicial set, $K'$. This well-known construction appears, for example, in \cite[Examples 1.3]{Curtis} and \cite[\S 1.3.1]{Martins_Porter}. The set of $m$-simplices of $K'$ consists of those sequences $(x_0\leq x_1 \leq \dots \leq x_m)$ of $0$-simplices of $K$ such that $\{x_0,x_1,\dots, x_m\}$ is a simplex of $K$. Also $s_i(x_0\leq x_1 \leq \dots \leq x_m)$ repeats the $i$th entry, whereas $d_i(x_0\leq x_1 \leq \dots \leq x_m)$ removes it.  Given a non-negative integer $i$, we have a bijection between $i$-faces of $K$ and non-degenerate $i$-simplices of $K'$, and  a canonical homeomorphism, $|K|\cong |K'|$.

In this paper, it will be convenient to consider a more general variant of ``triangulations'', defined in the broader context of simplicial sets. We will call these \emph{simplicial stratifications} of the manifold, as in the following definition.

\begin{Definition}\label{def:simp_strat} Consider a closed smooth $n$-manifold $\Sigma$.  A \emph{simplicial stratification} of $\Sigma$ is a homeomorphism, $\zeta_\Sigma\colon  |X_\Sigma| \to \Sigma$, where $X_\Sigma$ is a finite simplicial set. We denote  a simplicial stratification of $\Sigma$ by $(X_\Sigma,\zeta_\Sigma)$.

More generally, consider an $(n+1)$-cobordism, $(i,M,i')\colon \Sigma \to \Sigma'$, between the closed smooth $n$-manifolds, $\Sigma$ and $\Sigma'$.
A \emph{simplicial stratification (of the cobordism)} is given by a triad, $(Y_M; X_\Sigma, X_{\Sigma'})$, of simplicial sets, where  $X_\Sigma$ and $ X_{\Sigma'}$  are subcomplexes  of $Y_M$, and   $X_\Sigma\cap X_{\Sigma'}$ is empty,  together with a map, $(\zeta_\Sigma, \zeta_M, \zeta_{\Sigma'})$, of cospans in $\CGWH$, as  below, where the vertical arrows are homeomorphisms,
\begin{equation}\label{simp-strat}\vcenter{\xymatrix@R=18pt@C=50pt
{ |X_\Sigma|\ar[r]^{|j|}\ar[d]_{\zeta_{\Sigma}} & |Y_M|\ar[d]_{\zeta_M} &|X_{\Sigma'}|\ar[l]_{|j'|}\ar[d]^{\zeta_{\Sigma'}} \\
 \Sigma\ar[r]_i & M &{\Sigma'.}\ar[l]^{i'}
}}\end{equation}
Here,  $j\colon X_\Sigma \to Y_M$ and $j'\colon X'_{\Sigma'} \to Y_M$ denote the obvious simplicial inclusions.
\end{Definition}

Simplicial stratifications of a cobordism, $(i,M,i')\colon \Sigma \to \Sigma'$, arise, for instance, from triangulations of $M$ that extend given triangulations of $\Sigma$ and $\Sigma'$, and which, furthermore, are equipped with a total order on the set of vertices of the triangulation of $M$.
 Since simplicial stratifications can, in general, be chosen to be smaller than triangulations, they have advantages over the triangulated form.

\subsection{TQFTs from homotopy finite and finite  crossed complexes}\label{Quinn for finite A}
In this subsection, we will work over the field $\C$. We also fix $s \in \C$ and  $n\in \mathbb{Z}_{\ge 0}$. The fundamental crossed complex of a simplicial set, $S$, is here denoted $\Pi(S)$; see \S \ref{sec:Pi(S)}.

We now make strong use of the  notation and results from \S \ref{sec:class-maps} and \S \ref{sec:fib-Xcompvs-FibMSpaces}. We have an adjunction
 $\xymatrix{\hskip-0.5mm\Pi\colon \Simp
\ar@/^0.3pc/[r] \ar@{}[r]|\bot \ar@{<-}@/_0.3pc/[r] & \Crs: \hskip-0.5mm{\N},}$ in  Lemma \ref{brown_higgins_adjunction}.
If  $X$ is a simplicial set and $\A$  a  crossed complex, we therefore have a bijection,  natural in $X$ and $\A$,
\[\phi_X^\A\colon \Crs\big(\Pi(X),\A\big)\to \Simp\big(X,\N(\A)\big),\]
and a weak homotopy equivalence,
$
\etab_X^\A\colon \left|\N\big(\CRS(\Pi(X),\A)\big)\right|  \to \TOP(|X|,B_\A)
$.
We hence have a bijection,
$$T_{X}^\A\colon \pi_0\big(\CRS(\Pi(X),\A)\big)\to \pi_0\big (\TOP(|X|, B_\A)\big),$$
   of homotopy classes of maps between crossed complexes and between topological spaces, that is natural with respect to inclusions of subcomplexes of $X$.
\subsubsection{The explicit form of $\FQp{B_\A}{s}\colon \cob{n} \to \Vect_\C$}
We now fix a homotopy finite crossed complex $\A$. By   Lemma \ref{lem:chiNA}, the classifying space, $B_\A$, is homotopy finite. Hence we can consider  Quinn's finite total homotopy TQFT, $\FQp{B_\A}{s}\colon \cob{n} \to \Vect$,  of  \cite[Lecture 4]{Quinn}, as explored here in Subsection \ref{sec:QuinnTQFT}. Our main case of study is  when $\A$ is finite and reduced.  In this case, the formulae for $\FQp{B_\A}{s}$  become particularly simple.

Let $\Sigma$ be a closed smooth $n$-manifold,  with a simplicial stratification, $(X_\Sigma, {\zeta_{\Sigma}})$.
 Given  a crossed complex map, $f\colon \Pi(X_\Sigma) \to \A$, we  define the continuous map, $\underline{f}^{\zeta_\Sigma} \colon \Sigma \to B_\A=|\N(\A)|$, as the composite in the commutative diagram, below,
 \begin{equation}\label{underline-f}
 \vcenter{\xymatrix@R=18pt@C=20pt{ \Sigma\ar@/_0pc/[rd]_{\underline{f}^{\zeta_\Sigma}} \ar[r]^{\zeta_\Sigma^{-1}} &|X_\Sigma|\ar[d]^{|{\phi_X^\A(f)}|} \\ &|\N(\A)|.}}
 \end{equation}

 Given a crossed complex map, $f\colon \Pi(X_\Sigma) \to \A$, we denote its homotopy class  by
 $[f]_{\CRS(\Pi(X_\Sigma),\A)}$.
 Theorem \ref{lem:iso_groupoids} then gives  an isomorphism of vector spaces,
$$\mathbf{T}_{\zeta_{\Sigma}}^\A  \colon \C\big (\pi_0(\CRS(\Pi(X_\Sigma),\A))\big) \to \FQp{\Bc}{s}(\Sigma)= \C \big(\hpi_0( \Bc^\Sigma)\big) $$
where
$\mathbf{T}_{\zeta_{\Sigma}}^\A  \big( [f]_{\CRS(\Pi(X_\Sigma),\A)}\big):=\PC_{\underline{f}^{\zeta_\Sigma}}(\Bc^\Sigma)$.

Consider an
$(n+1)$-cobordism, $(i,M,i')\colon \Sigma \to \Sigma'$,
between  $\Sigma$ and $\Sigma'$, and  a simplicial stratification of  $(i,M,i')$,  as in Equation \eqref{simp-strat}. Consider $f\colon \Pi(X_\Sigma)\to \A$ and $f'\colon \Pi(X_{\Sigma'})\to \A$. Recall Subsection \ref{sec:QuinnTQFT}. Below, $|-|$ is used to denote both the geometric realisation of a simplicial set and the cardinality of a set.

\begin{Theorem}\label{Theo:Quinn_exp_CW}
In the formula, below, for the matrix elements of Quinn's finite total homotopy TQFT,  $\FQp{\Bc}{s}$, 
\begin{multline}\label{eq:matrixXcomp}
 \big \langle  \PC_{\underline{f}^{\zeta_\Sigma}}(\Bc^\Sigma) \mid  \FQp{\Bc}{s}\big ( [(i,M,i')]\big) \mid \PC_{{\underline{f}'{}^{\zeta_{\Sigma'}}}}(\Bc^{\Sigma'})  \big \rangle
\\= \chi^\pi\big( \big\{\underline{f}^{\zeta_\Sigma} |\Bc^M| {\underline{f}'{}^{\zeta_{\Sigma'}}}\big\}\big)\, \big(\chi^\pi\big(\PC_{\underline{f}^{\zeta_\Sigma}}(\Bc^{\Sigma})\big)\big)^s \, \big(\chi^\pi\big(\PC_{{\underline{f}'{}^{\zeta_{\Sigma'}}}}(\Bc^{\Sigma'})\big)\big)^{1-s},
\end{multline}
and  where, as before, we consider the space, with the induced CGWH topology,
$$
\big\{\underline{f}^{\zeta_\Sigma}|\Bc^M|{\underline{f}'{}^{\zeta_{\Sigma'}}}\big\}=
\bigg\{ H\colon M \to \Bc \bigg{|} \hskip-1cm
\vcenter
{
\xymatrix@R=0pt
{ 
 &&\Bc\\
 & \Sigma\ar[ur]^{\underline{f}^{\zeta_\Sigma}} \ar[dr]_{i} && \Sigma' \ar[ul]_{{\underline{f}'{}^{\zeta_{\Sigma'}}}} \ar[dl]^{i'} \\
&& M\ar[uu]^H }
} \textrm{ commutes}\bigg\}{\subseteq} \Bc^M,
 $$
 each factor can be calculated as follows:

 Using the notation of \S \ref{CRS-euler-formula} and \S \ref{sec:fib-Xcompvs-FibMSpaces}, we have, firstly,
  \begin{equation}\label{formula:mels_hfinite}
\chi^\pi \big( \big\{\underline{f}^{\zeta_\Sigma}|\Bc^M|{\underline{f}'{}^{\zeta_{\Sigma'}}}\big\}\big)=\chi^\pi\big(\CRS^{(\soml{ f}{f'})}
(\Pi(Y_M),\A)\big).
 \end{equation}

 If $\A$ is  finite, and reduced (i.e. with a single object), then
 \begin{multline}\label{formula:mels_finite}
\chi^\pi\big(\big\{\underline{f}^{\zeta_\Sigma}|\Bc^M|{\underline{f}'{}^{\zeta_{\Sigma'}}}\big\}\big) =\\
\bigg | \bigg \{ h\colon \Pi(Y_M) \to \A \,\bigg |  \hskip-1cm
\vcenter{ {
\xymatrix@R=0pt{ &&\A\\
& \Pi(X_\Sigma)\ar[ur]^{f} \ar[dr]_{\Pi(j)} && \Pi(X'_{\Sigma'}) \ar[ul]_{f'} \ar[dl]^{\Pi(j')} \\
&& \Pi(Y_M)\ar[uu]^h }}} \textrm{ commutes} \bigg \} \bigg |
 \\ \prod_{k=1}^\infty \big(\prod_{i=0}^\infty |A_{i+k}|^{\L(i,Y_M,X_\Sigma \cup {X'_{\Sigma'}} )}  \big ){}^{(-1)^k}.
      \end{multline} 
        
Continuing with  the evaluation of the other terms in the formula for the matrix elements in \eqref{eq:matrixXcomp}, we have
\begin{equation}
\begin{split}
\chi^\pi\big(\PC_{\underline{f}^{\zeta_\Sigma}}(\Bc^{\Sigma})\big)&=\chi^\pi\big(\PC_f(\CRS(\Pi(X_\Sigma),\A))\big),\\
\chi^\pi\big(\PC_{{\underline{f}'{}^{\zeta_{\Sigma'}}}}(\Bc^{\Sigma'})\big)&=\chi^\pi\big(\PC_{f'}(\CRS(\Pi(X'_{\Sigma'}),\A))\big),
\end{split}
\end{equation}
so (by Lemma \ref{path:compcrs}), if $\A$ is finite and reduced, then
\begin{equation}\label{formula:mels_finite2}
\begin{split}
\chi^\pi\big(\PC_{\underline{f}^{\zeta_\Sigma}}(\Bc^{\Sigma})\big)&=\big|[f]_{\CRS(\Pi(X_\Sigma),\A)}\big| \prod_{k=1}^\infty \big (\prod_{i=0}^\infty |A_{i+k}|^{\L(i,X_\Sigma)}  \big ){}^{(-1)^k},\\
\chi^\pi\big(\PC_{{\underline{f}'{}^{\zeta_{\Sigma'}}}}(\Bc^{\Sigma'})\big)&=\big |[f]_{\CRS(\Pi(X'_{\Sigma'}),\A)}\big| \prod_{k=1}^\infty \big (\prod_{i=0}^\infty |A_{i+k}|^{\L(i,{X'_{\Sigma'}})}  \big ){}^{(-1)^k},
\end{split}
\end{equation}
where, as above, $[f]_{\CRS(\Pi(X_\Sigma),\A)}|$  denotes the homotopy class of the  crossed complex map, $f$, and, analogously, $[f']_{\CRS(\Pi(X'_{\Sigma'}),\A)}|$  is that of $f'$.

\end{Theorem}
\begin{proof}
 This result follows from  the general discussion in \S\ref{sec:simpclass-space}, \S\ref{sec:fib-Xcompvs-FibMSpaces} and \S\ref{CRS-euler-formula}. In particular, the crucial ingredient is the Brown--Higgins--Sivera--Tonks weak homotopy equivalence, $$\etab_S^{\A}\colon  |\N (\CRS(\Pi(S),\A))|  \to \TOP(|S|,B_\A),$$ in item \eqref{simp-maps-to-top-maps} of Theorem \ref{MainCrs}, where  $S$ is a simplicial set and $\A$ is a crossed complex, and its refinement in Theorem \ref{thm:homotopy_equivalence_restriction}, together with the  result that the homotopy groups of a crossed complex coincide with those of its geometric realisation,  for which see again Theorem \ref{MainCrs}.
 
For instance, Equation \eqref{formula:mels_hfinite} follows from Theorem \ref{thm:homotopy_equivalence_restriction}, applied to the simplicial inclusion $\soml{j}{j'}  \colon X_{\Sigma} \sqcup  X'_{\Sigma'} \to Y_M$, and the crossed complex map, $$\soml{ f}{f'} \colon \Pi(X_\Sigma)\sqcup \Pi(X'_{\Sigma'}) \to \A.$$  We then apply Lemma \ref{lem:chiNA}. Finally,  \eqref{formula:mels_finite} follows from Lemma \ref{lem:twosides}, and Equation \eqref{formula:mels_finite2} follows from Lemma \ref{path:compcrs}.
\end{proof}

\begin{Remark}[CW-complexes]\label{rem:CW-TQFT}Even though we stated the previous theorem for simplicial stratifications of manifolds, it also holds for CW-decompositions of manifolds and cobordisms, using Corollaries \ref{cor:CW-class} and
 \ref{cor:CW-class-cob}, by switching from CW-decompositions arising from simplicial stratifications to general CW-decomposition. We will show some examples below.
\end{Remark}

\begin{Remark}[Independence from simplicial stratifications] \label{Independence for triangns}
Note that, by construction, all formulae for Quinn's finite total homotopy TQFT, $\FQp{B_\A}{s}$, in the previous theorem are independent of the chosen simplicial stratifications of the $n$-manifolds, $\Sigma$ and $\Sigma'$, and of the $(n+1)$-cobordism, $(i,M,i')\colon \Sigma \to \Sigma'$. There is no need to  make use of Alexander moves, or, equivalently, of Pachner moves, to prove triangulation-independence\footnote{nor, for this paper, independence from simplicial stratifications.},  as done for instance  in \cite{Barrett_Westbury,Mackaay,Turaev_Viro}. This is because the formulae were directly derived to give quantities that are, by construction, topologically invariant, and related to the homotopy content of function spaces.
\end{Remark}

\begin{Remark}By using the previous theorem together with Lemma \ref{lem:maps_on_basis}, we can see that the calculations of Quinn's finite total homotopy TQFT, $\FQp{B_\A}{s}$, for $\A$ a finite crossed complex, and for  given simplicial stratifications of the manifolds and cobordisms concerned,  could in theory be computed in finite time.
\end{Remark}

We expect that the techniques just shown will also be applicable for computing, explicitly, TQFTs derived from finite crossed complexes, which are, furthermore, equipped with a cohomology class valued in $U(1)$. (The existence of these TQFTs, generalising Dijkgraaf-Witten TQFTs \cite{DW}, was suggested in Remark \ref{coho-twist-Quinn}, and they were treated in \cite{Martins_Porter}, in the particular case of closed manifolds and crossed modules, using similar techniques to those of this paper.) Our approach here can likely also be adapted to give concrete formulae for homotopy quantum field theories derived from (classifying spaces of) crossed complexes, possibly equipped with appropriate cohomology classes. These homotopy quantum field theories are treated in \cite{Sozer-Virelizier:HQFT:ArXiv:2022}, and also in  \cite{Porter_HQFT,Porter_Turaev_HQFT}.

We also expect that similar techniques to those used  in this subsection can be used to give formulae for Quinn's finite total homotopy TQFT, $\FQp{\Bc}{s}$, in the case when $\Bc$ is the classifying space of a finite simplicial group, in which case we would obtain concrete formulae for all types of Quinn's finite total homotopy TQFT. (Since finite simplicial groups model all pointed homotopy finite spaces \cite{Ellis}.)
This would likely yield expressions similar to those in \cite{Porter}.

\subsection{Example: TQFTs from crossed modules of groups} We will now take advantage of the calculations in Subsection \ref{sec:computations}, and show some example computations of TQFTs derived from crossed modules, $\Gc=(\partial\colon E \to G, \trl)$, of groups; see \S \ref{Crossed modules of groups}. When $\Gc$ is finite, then the space $B_{\J_2(\Gc)}$ is homotopy finite. Given  $n$, a non-negative integer,  we thus have an $(n,n+1)$-TQFT, $\FQp{B_{\J_2(\Gc)}}{s}\colon \cob{n} \to \Vect_\C$, as in Definition \ref{def:quinnTQFT}, where  $s \in \C$ is a  parameter, following \cite[Section 4]{Quinn}.

We will only consider the case $s=1$, and we will write, only for this subsection, $$\big(\Z_\Gc\colon \cob{n} \to \Vect_\C\big):=\big( \FQp{B_{\J_2(\Gc)}}{s=1}\colon \cob{n} \to \Vect_\C\big).$$ Given a finite group $G$, we can turn it into a crossed module $\{1\}\to G$, and we put $\Z_G:=\Z_{\{1\} \to G}$. The  TQFT $\Z_G\colon \cob{n} \to \Vect_\C$, is well-known to coincide with Dijkgraaf-Witten TQFT, with a trivial cocycle.

The closed  manifold case of $\Z_\Gc\colon \cob{n} \to \Vect_\C$ was discussed in \cite{Martins_Porter}, in the context of homologically twisted Yetter TQFT, \cite{Yetter}.  A more recent paper, \cite{Sozer-Virelizier:HQFT:ArXiv:2022}, addresses closely related HQFTs derived from finite crossed modules.

 Let us fix a finite crossed module $\Gc=(\partial\colon E \to G, \trl)$.

\subsubsection{$(1,2)$-TQFTs derived from finite crossed modules}\label{sec:123TQFT_finXmod}
For brevity, we will only consider the oriented case, and show computations mainly for CW-complexes. Let $\cob{1}_\oo$ be the category of closed 1-manifolds and oriented cobordisms.  Let $\Z^\oo_\Gc\colon \cob{1}_\oo \to \Vect$ denote the restriction of $\Z_\Gc\colon \cob{1}\to \Vect$ to $\cob{1}_\oo$.

The category, $\cob{1}_\oo$, is generated, as a symmetric monoidal category, see \cite{Kock:TQFT:2003}, by the generator $S^1$, and the cobordisms,
\begin{itemize}
[leftmargin=0.6cm]
\item $\epsilon:=S^1 \ra{(\iota, D^2, 0_{D^2})} \emptyset$, and $\eta:=\emptyset \ra{(0_{D^2}, D^2, \iota)} S^1,
$ where $\iota\colon S^1 \to D^2$ is the obvious inclusion, and $0_{D^2}\colon \emptyset \to D^2$ is the initial map, (cf. Figure \ref{2cobgens}).
\item $\mu:=S^1\sqcup S^1 \ra{(\langle \iota_L,\iota_R\rangle,\Mu,\iota_C)} S^1$ and
$\Delta:=S^1 \ra{(\iota_C,\Mu ,\langle \iota_L,\iota_R\rangle)} S^1 \sqcup S^1,$ where $\Mu$ is as in Figure \ref{m}, and the notation used is defined in \S \ref{sec:pants}.
\end{itemize}

In \S \ref{sec:G/GG},  we considered $S^1$ with a CW-decomposition with a unique $0$-cell, at the south pole. Applying Theorem \ref{Theo:Quinn_exp_CW}, we have, using the notation in \S \ref{sec:G/GG}, that $$\Z^\oo_{\Gc}(S^1)\cong\C\pi_0\big(\CRS\big (\Pi(S^1_\sk),\J_2(\Gc)\big)\big)\cong \C\pi_0\big (G\sslash \Gc ).$$
Hence $\Z^\oo_\Gc(S^1)$ is the free vector space on $G/(G\ltimes E)$, the set of orbits of the $\bullet$ action, of $G\ltimes E$ on the underlying set of $G$. These are in  bijection with conjugacy classes of the quotient group, $G/\partial(E)$.

We can then determine the matrix elements of $\Z^\oo_{\Gc}$ assigned to the generating cobordisms of $\cob{1}_\oo$.  Below, given $g \in G$, we put $[g]={\Orb}_{G\ltimes E}(g)$.
\begin{enumerate}[leftmargin=0.6cm]
\item   $\Z^\oo_{\Gc}(S^1\sqcup S^1)\cong \C\pi_0\big (G\sslash \Gc ) \otimes  \C\pi_0\big (G\sslash \Gc )$,
\item \label{eqeps} $\langle [g] \mid \Z^\oo_{\Gc}(\epsilon) \mid 1\rangle=\chi^\pi\big(P_{S^{1}}^{-1}(g)\big)=|\partial^{-1}(g)|$, using the notation in \S \ref{sec:eps_eta}.\item \label{eqeta}  $\langle 1 \mid \Z^\oo_{\Gc}(\eta ) \mid [g]\rangle=\chi^\pi\big(P_{S^{1}}^{-1}(g)\big)\chi^\pi\big(\PC_g(G\sslash \Gc)\big)= |\partial^{-1}(g)|\frac{|[g]|}{|G|}$.
\item Using the notation and calculations of \S \ref{sec:pants}, we have,
\begin{align*}\langle [g] &\otimes [g']|  \Z_\Gc(\mu)| [h]\rangle= \chi^\pi\big ( \CRS^{ (h, g,g')}(\Pi(\Mu_\sk),\J_2(\Gc)\big)\big)\,\, \chi^\pi\big (\PC_h(G\sslash \Gc)\big) \\ &=\frac{\big\{ (p,q,e) \in G^2 \times E: \partial(e)=h^{-1}q^{-1}g'q p^{-1} g p\}}{|E|^2} \frac{|[h]|}{|G|},
\end{align*}
and
\begin{multline*}\langle [h]\mid  \Z_\Gc(\Delta)|   [g] \otimes [g'|\rangle\\= \chi^\pi\big ( \CRS^{ (h, g,g')}(\Pi(\Mu),\J_2(\Gc))\big) \,\, \chi^\pi\big (\PC_g(G\sslash \Gc)\big)  \chi^\pi\big (\PC_{g'}(G\sslash \Gc)\big) \\ =\frac{\big\{ (p,q,e) \in G^2 \times E: \partial(e)=h^{-1}q^{-1}g'qp^{-1} g p\}}{|E|^2} \frac{|[g]| \, |[g']|}{|G|^2}.
\end{multline*}
\end{enumerate}

We can then determine $\Z^\oo_{\Gc}(S^2)$ either  as, using the notation in \S \ref{sec:S2}, $$\chi^\pi\big(\CRS\big(\Pi(S^2_\sk)\big) ,\Gc\big)=\chi^\pi\big(\ker(E)\| \Gc\big)=\frac{|\ker(\partial) | |E|}{|G|},$$
or as, using the fact that  $\Z^\oo_\Gc\colon \cob{1}_\oo \to \Vect$ is a functor:
\begin{align*}
 \sum_{[g] \in G/  G\ltimes E}\big \langle 1 \mid \Z^\oo_{\Gc}(\eta ) \mid [g]\big\rangle \, \big \langle [g] \mid \Z^\oo_{\Gc}(\epsilon) \mid 1\big\rangle=&  \sum_{[g] \in G/  G\ltimes E}   |\partial^{-1}(g)|\frac{|\Orb_{G \ltimes E} (g)|}{|G|} |\partial^{-1}(g)|\\&=\frac{|\ker(\partial)|^2 |\partial(E)|} {|G|} =\frac{|\ker(\partial) | |E|}{|G|}.
\end{align*}

It is then easy to see that $\Z^\oo_{\Gc}\colon \cob{1}_\oo\to \Vect$ is determined by the restriction of $\Z_{G/\partial(E)}\colon \cob{1}\to \Vect$ to $\cob{1}_\oo$, that is a TQFT defined from the group $G/\partial(E)$, alone; the  additional factors arise from Euler characteristics of cobordisms.

\subsubsection{A quick look at $(2,3)$-TQFTs derived from finite crossed modules} We now address $\Z_\Gc\colon \cob{2}\to \Vect$.
Using the approach in \cite{Martins_Porter}, it follows that $\Z_\Gc\colon \cob{2} \to \Vect$ coincides with Yetter's homotopy 2-type TQFT in \cite{Yetter,Porter}.

Let us show some computations. They will in particular demonstrate that $\Z_\Gc\colon \cob{2}\to \Vect$ is more general than the (2,3)-TQFTs, $\Z_G\colon \cob{2} \to \Vect$, which can be derived from groups alone, unlike for the case $\Z_\Gc\colon \cob{1}_\oo\to \Vect$.

The monoidal category $\cob{2}$ cannot be finitely presented, since there is an infinite number of diffeomorphism classes of closed surfaces. However, given a surface $\Sigma$, we can easily determine $\Z_\Gc(\Sigma)$ by taking a CW-decomposition of $\Sigma$.
It follows, by using  the computations and notation in \S \ref{sec:S2} and  \S \ref{sec:torus}, that:
$$\Z_{\Gc}(S^2)\cong\C \pi_0\big(\ker(\partial) \| \Gc\big), \qquad  \textrm{ and }\qquad \Z_{\Gc}(T^2)\cong\C \pi_0\big(T^2(\Gc)\big),$$
where we have $\pi_0\big(\ker(\partial) \| \Gc\big)=\ker(\partial)/G$ and $\pi_0\big(T^2(\Gc)\big)=T^2_0(\Gc)/\big(G\ltimes (E \times E)\big)$.

The value of $\Z_\Gc\colon \cob{2} \to \Vect$ on the cobordism $\nu:= (\iota_C, \Nu,\langle \iota_L,\iota_R\rangle) \colon S^2  \to S^2\sqcup S^2$, in Figure \ref{fig:Nu}, page \pageref{fig:Nu}, can be calculated by combining the discussion of  \S \ref{sec:S2} with Theorem \ref{Theo:Quinn_exp_CW} / Remark \ref{rem:CW-TQFT}.
This gives that, for $a,a',b \in \ker(\partial)$, and putting   $[a]=\Orb_G(a)$, that \begin{align*}
\big\langle [b] &  |\Z_\Gc(\nu) | [a]\otimes [a']\big\rangle\\&=\chi^\pi\big( (P_\partial^\Nu)^{-1}(b,a,a')\big) \, \chi^\pi(\PC_a(\ker(\partial)\| \Gc))
\, \chi^\pi(\PC_{a'}(\ker(\partial)\| \Gc))\\
                                  &=\big |\{(g,g')\in G\times G: b= a \trl g \,\, a' \trl g'\} \big  |\,  \big  |\Orb_G(a) \big  |\, \big  |\Orb_G(a')\big |.
\end{align*}

For a crossed module, $\Gc$, with trivial $G$, so $\Gc=(E \to \{1\})$, this gives
$
\big\langle [b]   |\Z_\Gc(\nu) | [a]\otimes [a']\big \rangle = \delta(b, a a').$
This  shows that $\Z_{\Gc}$ is strictly more general than the TQFTs that can be obtained from groups, $H$, alone. Indeed, in the latter case $\Z_H(S^2)\cong\C$, and $\Z_H(\nu)\colon \C \to \C \otimes \C \cong \C$ is the identity  map. An easy analysis furthermore shows that $\Z_\Gc$, with $\Gc=(E \to \{1\})$, and $E$ non-trivial, also cannot be obtained as the direct sum of TQFTs, $\Z_H\colon \cob{2} \to \Vect$, derived from groups.

\subsubsection{Why should we  bother with crossed modules?}

The TQFTs, $\Z_\Gc\colon \cob{n} \to \Vect_\C$, obtained from  finite crossed modules, of groups, are strictly more general than the ones that can be obtained from groupoids, equivalently from disjoint unions of finite groups. An easy example showing that this is so arises when $n=4$.

Consider the CW-decomposition of $S^4$ with unique $0$- and  $4$-cells, and no other cells, and the CW-decomposition of $S^2$ with unique $0$- and  $2$-cells. The product CW-decomposition on $S^2\times S^2$, which we will use, then has a unique $0$-cell, no 1-cells, two 2-cells, no $3$-cell, and one $4$-cell. We then have
\begin{align*}
 \Pi(S^4_\sk)&\cong \dots \to \{0\}  \to \{0\} \to \mathbb{Z} \to \{0\} \to\{0\} \to \{1\},\\
 \Pi\big ((S^2\times S^2)_\sk\big)&\cong  \dots \to \{0\}  \to \{0\} \to \mathbb{Z} \to \{0\} \to \mathbb{Z}\oplus \mathbb{Z} \to \{1\}.
\end{align*}

Let $H$ be a finite groupoid. A simple calculation gives that the groupoids
$
\pi_1\big(\CRS\big (\Pi(S^4_\sk),\J_1(H)\big )\big)$ and $\pi_1\big(\CRS\big (\Pi\big((S^2\times S^2)_\sk\big),\J_1(H)\big )\big)$  are both isomorphic to $H$. In  particular, the state spaces, $\Z_H(S^4)$ and $\Z_H(S^2\times S^2)$,
both have dimension  given by the cardinality of the set of components of $H$.

If $\Gc=(\d\colon E \to G, \trl)$ is a  crossed module,  clearly
$\pi_1\big ( \CRS\big (\Pi(S^4_\sk),\J_2(\Gc)\big )\big)\cong G/\d(E), $
a groupoid with a single object. In particular, the dimension of $\Z_\Gc(S^4)$ is always 1, independently of the crossed module $\Gc$.

Using a computation very similar to that of \S \ref{sec:S2}, we have that
$$
\pi_1\big ( \CRS\big (\Pi\big ( (S^2\times S^2)_\sk\big ),\J_2(\Gc)\big )\big) \cong \big (\ker(\partial) \oplus \ker(\partial) \big)\sslash \big (G/\d(E)\big),
$$
where the group $G/\d(E)$ acts on $\ker(\partial)\oplus \ker(\partial)$ via $(a,b)\trl [g]=(a \trl g^{-1}, b \trl g^{-1})$.
In particular,
$\dim\big(\Z_\Gc(S^2 \times S^2)\big)$ is given by the number of orbits of the action of $G$ on $\ker(\partial) \times \ker(\partial)$. This gives, in general, a  value different from $\dim\Z_\Gc(S^4)=1$.

\subsection{The once-extended TQFTs derived from finite  crossed complexes}\label{ETQFT_comp}
 In this section and the next, we will work over the field of rational numbers $\Q$, and also fix a finite (i.e. not just homotopy finite) crossed complex, $\A$. We will freely use the notation and results from Section \ref {sec:def-once-extended}, particularly Subsection \ref{sec:fin_ext}.

Let $\Bc=B_\A$, the classifying space of $\A$, which we recall is a homotopy finite space, by Theorem \ref{MainCrs}. Let $n$ be a non-negative integer. We will give explicit formulae for some instances of the finitary once-extended Quinn TQFT,
$$\tFQd{\Bc}\colon \tdcob{n}{\Bc} \to \vProfGrpfin,$$ and consequently of its Morita version, $$\tFQmor{\Bc}\colon \tdcob{n}{\Bc} \to \Mor.$$
As before, formulae will be given in terms of simplicial stratifications of manifolds, more generally  CW-decompositions, the latter mainly in order to simplify computations.

We will  also treat a few variants of these once-extended TQFTs, as mentioned in the beginning of this section,
which will be done at the same time as we examine the dependence of the formulae that we provide on the choice of the simplicial stratification of an $n$-dimensional manifold that is being used.
\subsubsection{The $B_\A$-decoration of a manifold arising from a simplicial stratification}
Let $\Sigma$ be a closed (and as usual smooth) $n$-manifold. Recall, Definition \ref{def:DecoratedMan}, that, given a HF space, $\Bc$,  a $\Bc$-decoration, $\fd_{\Sigma}$, of  $\Sigma$, is given by a finite subset, $\fd_\Sigma$, of the function space, $\Bc^\Sigma$, of functions from $\Sigma$ to $\Bc$, containing at least one function, $f\colon \Sigma \to \Bc$, from each path-component of $\Bc^\Sigma$. If $\Bc=B_\A$, simplicial stratifications of $\Sigma$ naturally give rise to $\Bc$-decorations of $\Sigma$. We can see this as follows.

 Let  $X_\Sigma$ be a finite simplicial set. By Theorem  \ref{MainCrs}, we have a  weak homotopy equivalence,
  $\etab_{X_\Sigma}^\A\colon |\CRS(\Pi(X_\Sigma),\A)| \to \TOP(|X_\Sigma|,B_\A).$ Given a crossed complex map, $f\colon \Pi(X_\Sigma) \to \A$, we define  $\underline{f}^{\zeta_\Sigma} \colon \Sigma \to B_\A$ by using the commutative diagram in \eqref{underline-f}. There are only finitely many crossed complex maps from  $\Pi(X_\Sigma)$ to $\A$, by Lemma \ref{lem:crs_is_finite}. We thus have a $B_\A$-decoration of $\Sigma$,  defined by
\begin{equation}\label{eq:dec}
  \fd_{\Sigma}(\zeta_{\Sigma},\A):=\big\{ \underline{f}^{\zeta_\Sigma} \mid  f \colon \Pi(X_\Sigma) \to \A) \big \}.
 \end{equation}

\subsubsection{Explicit formulae  for the finitary once-extended Quinn TQFT for $\Bc=B_\A$}\label{fin_Quinn_xcomp}
We now give explicit formulae for the finitary once-extended Quinn TQFT,
 $\tFQd{B_\A}\colon \tdcob{n}{B_\A} \to \vProfGrpfin$, in Definition \ref{def:finitary_ext_TQFT}.
 The formulae will be given, firstly,  in terms of simplicial stratifications of  the manifolds, cobordisms and (reduced; see \S \ref{sec: reduced cobordism}) extended cobordisms. This is analogous, of course, to taking  triangulations, so as to get `lattice models' and `state sum' models, as we mentioned in Remark \ref{Independence for triangns}.
 We will discuss the dependence on these choices in a later subsection.  We will then briefly address the case of CW-decompositions of manifolds.

 In order to write our formulae, we use the form for the matrix elements of  $\tFQd{B_\A}$ in Theorem \ref{thm:red-form-2matrixeles}.

\begin{enumerate}[label=(\roman*), leftmargin=0.6cm]
 \item\label{etqft1} If $\Sigma$ is a closed $n$-manifold, and $\zeta_{\Sigma} \colon |X_\Sigma| \to \Sigma$ is a {simplicial stratification} of $\Sigma$, then we have a canonical isomorphism of groupoids,
 \begin{equation}\label{group:Q}
 \tFQd{B_\A}\big ( \Sigma, \fd_{\Sigma}\big (\zeta_{\Sigma} ,\A)\big)\cong \pi_1\big (\CRS(\Pi(X_\Sigma),\A)\big).
 \end{equation}

\item\label{etqft2} Given an $(n+1)$-cobordism, $(i,M,i')\colon \Sigma \to \Sigma'$, between the closed smooth $n$-manifolds, $\Sigma$ and $\Sigma'$,
consider a simplicial stratification of the cobordism, $(i,M,i')$, derived from a  homeomorphism of cospans in $\CGWH$,
as below
\begin{equation}\label{eq:simp_stract_cob}\vcenter{\xymatrix@R=18pt@C=50pt
{ |X_\Sigma|\ar[r]^{|j|}\ar[d]_{\zeta_{\Sigma}} & |Y_M|\ar[d]_{g_M} &|X'_{\Sigma'}|\ar[l]_{|j'|}\ar[d]^{\zeta'_{\Sigma'}} \\
 \Sigma\ar[r]_i & M &{\Sigma'.}\ar[l]^{i'}
}}\end{equation}
We have $B_\A$-decorations,
$\fd_{\Sigma}(\zeta_{\Sigma},\A),$  of $\Sigma$, and
$\fd_{\Sigma'}(\zeta'_{\Sigma'}, \A), $ of $\Sigma'$,
giving  the associated 1-morphism in the bicategory $\tdcob{n}{B_\A}$,
 $$ \big(\Sigma, \fd_{\Sigma}(\zeta_{\Sigma} ,\A)\big) \ra{(i,M,i')} \big(\Sigma', \fd_{\Sigma'}(\zeta'_{\Sigma'},\A)\big).$$
Using Definitions \ref{def:finitary_ext_TQFT} and  \ref{def:crs_profunctor}, we have a natural isomorphism of profunctors,
\begin{multline}
\label{profunctor:Q}
\tFQd{B_\A}\big(  \big(\Sigma, \fd_{\Sigma}(\zeta_{\Sigma} ,\A)\big) \ra{(i,M,i')}  \big(\Sigma', \fd_{\Sigma'}(\zeta'_{\Sigma'} ,\A)\big) \big)  \\ \cong  \Lin\circ  \Hp^{\big (|Y_M|_\sk;|X_\Sigma|_\sk,|X'_{\Sigma'}|_\sk\big)}_\A.
\end{multline}(Recall also that $\Lin\colon \Sets \to \Vect$ is the free vector space functor.)
\item\label{etqft3} Finally, at the level of 2-morphisms, consider an $(n + 2)$-extended cobordism, $\K$, and also its reduction, $\hat{\K}$, as defined in \S\ref{sec: reduced cobordism}, both written in \eqref{diag_for_K} below,
\begin{equation}\label{diag_for_K}
\K=\vcenter{\xymatrix@C=2.6pc@R=20pt{
  \Sigma  \ar[r]^{i}\ar[d]_{\iota_0}  &  {M}\ar[d]^{i_N}   & \Sigma'\ar[l]_{i'} \ar[d]^{\iota_0'}\\
  \Sigma\times I\ar[r]^{i_E}   &  K & \Sigma'\times I, \ar[l]_{i_W} \\
  \Sigma \ar[u]^{\iota_1} \ar[r]_{j}  &  {M}'\ar[u]_{i_S}    & \Sigma'\ar[l]^{j'} \ar[u]_{\iota_1'}
 }} \,\, \,\, \hat{\K} =\vcenter{\xymatrix@C=2.8pc@R=2pc{
   &  M\ar[d]^{\hat{i_N}}   & \\
  \Sigma\ar[r]^{\hat{i_E}}\ar[ur]^{{i}}  \ar[rd]_{{j}}  &  \hat{K} & \Sigma'\ar[l]_{\hat{i_W}} \ar[ul]_{{i'}}\ar[ld]^{{j'}} \\
    &  M'\ar[u]_{\hat{i_{S}}}     &
 }}.
  \end{equation}

Consider also a diagram,  of finite simplicial sets, $\hat{W}_\K$, in \eqref{simp_stract_for_K}, below,
together with a homeomorphism of diagrams,  $g\colon |\hat{W}_\K| \to \hat{\K}$, in $\CGWH$, where $|\hat{W}_\K|$ is obtained by applying geometric realisation to all components of $\hat{W}_\K$, and, in order to simplify notation, all `components' of $g$  will be  denoted  $g$,
\begin{equation}\label{simp_stract_for_K}
\hat{W}_\K=\vcenter{\xymatrix@C=3pc@R=20pt{
     &  Y_{M}\ar[d]^{\overline{i_N}}   &   \\
  X_{\Sigma} \ar[ru]^{\overline{i}}  \ar[rd]_{\overline{j}}  \ar[r]^{\overline{i_E}}   &  Z_{\hat{K}} & X'_{\Sigma'}, \ar[lu]_{\overline{i'}}  \ar[l]_{\overline{i_W}} \ar[ld]^{\overline{j'}}\\
   &  Y'_{M'}\ar[u]_{\overline{i_S}}    &
 }} \quad\textrm{ and } g\colon |\hat{W}_\K| \to \hat{\K}.
   \end{equation}
  Note that $g\colon |\hat{W}_\K| \to \hat{\K}$ gives simplicial stratifications, $\zeta_\Sigma$ and $\zeta'_{\Sigma'}$, for $\Sigma$ and $\Sigma'$, extending to simplicial stratifications  of the  $(n+1)$-cobordisms, $(i,M,i')\colon \Sigma \to \Sigma'$ and $(j,M',j')\colon \Sigma \to \Sigma'$, and for the reduced extended cobordism $\hat{K}$.

 Also, consider the following pushout of simplicial sets, $$\fr(\hat{W}_\K):=Y_M\sqcup_{(X_\Sigma \sqcup X'_{\Sigma'})} Y'_{M'},$$ and let  $f_{\hat{W}_\K}\colon  \fr(\hat{W}_\K)\to    Z_{\hat{K}}$ be defined from the universal property of pushouts. Since we have cofibrations,
  $$\langle \overline{i},   \overline{i'} \rangle \colon X_\Sigma \sqcup X'_{\Sigma'}\to Y_M, \quad  \textrm{ and }
 \quad \langle  \overline{j},  \overline{j'} \rangle \colon X_\Sigma \sqcup X'_{\Sigma'}\to Y'_{M'},$$
 the higher homotopy van Kampen Theorem \cite[8.2.i]{brown_higgins_sivera}, or the discussion in \S \ref{sec:freeness-Xcomp}, gives that  $$\Pi\big(\fr(\hat{W}_\K)\big)\cong \Pi(Y_M)\sqcup_{(\Pi(X_\Sigma) \sqcup \Pi(X'_{\Sigma'}))} \Pi(Y'_{M'}).$$

Consider a  commutative diagram of crossed complex maps, as  below,
\begin{equation}\label{eq:2cobCRS}
\vcenter{\xymatrix@C=2pc@R=20pt{
     &  \Pi(Y_{M})\ar[d]^{H}   &   \\
  \Pi(X_{\Sigma}) \ar[ru]^{\Pi(\overline{i})}  \ar[rd]_{\Pi(\overline{j})}  \ar[r]^{f}   &  \A & \Pi(X'_{\Sigma'}) \ar[lu]_{\Pi(\overline{i'})}  \ar[l]_{f'} \ar[ld]^{\Pi(\overline{j'})}\\
   &  \Pi(Y'_{M'})\ar[u]_{H'}    &
 }}.
 \end{equation}
  Let  $[H,H']\colon \Pi\big(\fr(\hat{W}_K)\big)\to \A$ be defined by the universal property of pushouts.
 As in \S\ref{sec:fib-Xcompvs-FibMSpaces},  let $\CRS^{( [H,H'])}(\Pi (Z_{\hat{K}}),\A \big )$ denote the fibre of the restriction map $(f_{\hat{W}_\K})^*\colon \CRS\big (\Pi (Z_{\hat{K}}),\A \big ) \to \CRS\big (\Pi (\fr(\hat{W}_\K)),\A \big )$, at $[H,H']$.
 The fibre of the restriction map $\CRS(\Pi(Y'_{M'}),\A)\to \CRS\big(\Pi(X_\Sigma) \sqcup \Pi(X'_{\Sigma'}),\A\big)$ at the crossed complex map $\langle f,f' \rangle\colon \Pi(X_\Sigma) \sqcup \Pi(X'_{\Sigma'}) \to \A$  is denoted $\CRS^{(\langle f,f' \rangle)}(\Pi(Y'_{M'}),\A)$.

Applying the adjunction, $\xymatrix{\hskip-0.5mm\Pi\colon \Simp
\ar@/^0.3pc/[r] \ar@{}[r]|\bot \ar@{<-}@/_0.3pc/[r] & \Crs: \hskip-0.5mm{\N},}$ in  Lemma \ref{brown_higgins_adjunction} to the diagram in \eqref{eq:2cobCRS} gives rise to the commutative diagram of simplicial sets,

\begin{equation}\label{eq:2cobSIMP}
 \quad \vcenter{\xymatrix@C=4.9pc@R=24pt{
     &  Y_{M} \ar[d]|{\phi^\A_{Y_M}(H)}   &   \\
X_{\Sigma} \ar[ru]^{\overline{i}}  \ar[rd]_{\overline{j}}  \ar[r]|{{\phi^\A_{X_\Sigma}(f)}}   &  \N(\A) & X'_{\Sigma'}. \ar[lu]_{\overline{i'}}  \ar[l]|{\phi^\A_{X'_{\Sigma'}}(f')} \ar[ld]^{\overline{j'}}\\
   &  Y'_{M'}\ar[u]|{\phi^\A_{Y'_{M'}}(H')}    &
 }}
 \end{equation}

Continuing from \eqref{underline-f}, we apply geometric realisation to \eqref{eq:2cobSIMP}, and compose with $g^{-1}\colon \hat{\K} \to  |\hat{W}_\K|$, to get  the following diagram of spaces and continuous maps,
$$\vcenter{\xymatrix@C=4.8pc@R=2pc{
   &  M\ar[d]|{\underline{H}}   & \\
  \Sigma\ar[r]|{\underline{f}^{\zeta_\Sigma}}\ar[ur]^{{i}}  \ar[rd]_{{j}}  &  B_\A & \Sigma',\ar[l]|{{\underline{f}'{}^{\zeta_{\Sigma'}}}}  \ar[ul]_{{i'}}\ar[ld]^{{j'}} \\
    &  M'\ar[u]|{\underline{H}'}     &
 }} $$
where (we recall that all components of  $g\colon |\hat{W}_\K| \to \hat{\K}$ are denoted $g$),
 \begin{align*}
 \underline{H}&=\big |\phi_{Y_{M}}^\A(H)\big|\circ g^{-1}\colon M \to B_\A,
 &\underline{H}'&=\big|\phi_{Y'_{M'}}^\A(H')\big|\circ g^{-1}\colon M' \to B_\A,\\
 \underline{f}^{\zeta_\Sigma}&=\big |\phi_{X_\Sigma}^\A(f)\big|\circ g^{-1}\colon \Sigma \to B_\A,
 &\underline{f}'{}^{\zeta_{\Sigma'}}&=\big |\phi_{X'_{\Sigma'}}^\A(f')\big|\circ g^{-1}\colon \Sigma' \to B_\A.
 \end{align*}

This leads to the following description of the corresponding matrix entries of $\tFQd{B_\A}\colon \tdcob{n}{B_\A} \to \vProfGrpfin,$ where, to  simplify expressions slightly, we have  written $\Bc$ for $B_\A$.
\begin{multline}\label{formula_matrix_2els}
 \big \langle \PC_{\underline{H}}\big ({\{\underline{f}^{\zeta_\Sigma}|\Bc^{{M}}|\underline{f}'{}^{\zeta_{\Sigma'}}\}}\big ) \mid \big(\tFQ{\Bc}^2([\K])\big)_{(\underline{f}^{\zeta_\Sigma},\underline{f}'{}^{\zeta_{\Sigma'}})}\mid  \PC_{\underline{H'}}\big ({\{\underline{f}^{\zeta_\Sigma}|\Bc^{M'}|\underline{f}'{}^{\zeta_{\Sigma'}}\}} \big )\big \rangle\\
= \chi^{\pi}\big( \CRS^{( [H,H']) }\big(\Pi( Z_{\hat{K}}),\A\big)\big)\, \chi^\pi\big(\PC_{H'}\big(\CRS^{({\soml{ f}{f'}})} \big(\Pi(Y'_{M'}),\A\big)\big) \big).
 \end{multline}
\end{enumerate}

 \begin{Theorem}\label{Thm:fin_Quinn_xcomp}
 Let $\A$ be a finite crossed complex, and $n$   a non-negative integer. The structures specified in \ref{etqft1},  \ref{etqft2}, and  \ref{etqft3}, above, give the finitary once-extended Quinn TQFT, in Definition \ref{FinQ}, namely
  $\tFQd{B_\A}\colon \tdcob{n}{B_\A} \to \vProfGrpfin$,
 if we restrict to the objects of  $\tdcob{n}{B_\A}$  of the form 
   $\big(\Sigma, \fd_{\Sigma}(\zeta_{\Sigma} ,\A)\big)$, where  $\zeta_{\Sigma} \colon |X_\Sigma| \to \Sigma$ is a simplicial stratification of a closed smooth $n$-manifold $\Sigma$.
 \end{Theorem}
 
 \begin{proof}
  For the most part, the proof is essentially as in Theorem \ref{Theo:Quinn_exp_CW}. For instance,  Equation \eqref{group:Q} follows from Lemma \ref{lem:iso_groupoids}, and Equation \eqref{formula_matrix_2els} follows from Lemma \ref{thm:homotopy_equivalence_restriction}.  Equation \eqref{profunctor:Q} follows from the fact that we have a weak homotopy equivalence,
  $$\etab_{Y_M}^\A\colon  \big |\N (\CRS(\Pi(Y_M),\A))\big|  \to \TOP(|Y_M|, B_\A),$$
by (Brown--Higgins--Sivera--Tonks) Theorem \ref{MainCrs}.
 \end{proof}

 If we assume, furthermore, that $\A$ is reduced, then the crossed complexes appearing in \eqref{formula_matrix_2els} are homogeneous. This follows  from the discussion in  Subsection \ref{sec:homfinXcomp}. In particular (as in Theorem \ref{Theo:Quinn_exp_CW}), we can obtain very simple formulae for their homotopy content, similar to  \eqref{formula:mels_finite} and  \eqref{formula:mels_finite2}, using  Lemmas \ref{lem:twosides}.

 We will show the latter for the case of CW-decompositions of manifolds. Again, consider an $(n + 2)$-extended cobordism, $\K$, and also its reduction, $\hat{\K}$,  in \eqref{diag_for_K}. Suppose that we are given simplicial stratifications $\zeta_{\Sigma}\colon |X_\Sigma| \to \Sigma$, and ${\zeta_{\Sigma'}}\colon |X'_{\Sigma'}| \to \Sigma$, of $\Sigma$ and $\Sigma'$.
   Choose  CW-decompositions of $M$ and $M'$ extending the CW-decompositions of $\Sigma$ and $\Sigma'$ induced by their simplicial stratifications, and also consider a CW-decomposition of $\hat{K}$, extending that of $M$ and that of $M'$.

   Recall the notation in Corollary \ref{cor:CW-class}. Choose crossed complex morphisms $H\colon \Pi(M_\sk)\to \A$ and $H'\colon \Pi(M'_\sk)\to \A$, extending, respectively, $f\colon \Pi(\Sigma_\sk)\to \A$ and   $f'\colon \Pi(\Sigma'_\sk)\to \A$. Put $\Bc=B_\A$, and consider the following continuous maps,
   \begin{align*}&\hat{\underline{H}}=\overline{\T}_{M_\sk}^\A(H)\colon M \to \Bc, & \textrm{ and } && \hat{\underline{H}}{}'=\overline{\T}_{M'_\sk}^\A(H')\colon M' \to \Bc.\end{align*}
  We can suppose that both $\underline{\hat{H}}$ and $\hat{\underline{H}}{}'$  extend $\langle \underline{f}^{\zeta_\Sigma}, {\underline{f}'{}^{\zeta_{\Sigma'}}}\rangle \colon \Sigma \sqcup \Sigma' \to \Bc$.
 \begin{Corollary}\label{cor:matrix-els-CW}If $\A$ is finite and reduced, and $\Bc=B_\A$, then the respective matrix elements of $\tFQd{\Bc}\colon \tdcob{n}{B_\A} \to \vProfGrpfin$ can be computed as:
 \begin{multline*}
 \Big \langle \PC_{\underline{\hat{H}}}\big ({\{\underline{f}^{\zeta_\Sigma}|\Bc^{{M}}|\underline{f}'{}^{\zeta_{\Sigma'}} \}}\big ) \mid \big(\tFQ{\Bc}^2([\K])\big)_{(\underline{f}^{\zeta_\Sigma},\underline{f}'{}^{\zeta_{\Sigma'}} )}\mid  \PC_{\hat{\underline{H}}{}'}\big ({\{\underline{f}^{\zeta_\Sigma}|\Bc^{M'}|\underline{f}'{}^{\zeta_{\Sigma'}}\}} \big )\Big \rangle\\
=\left | \left \{T\colon \Pi(\hat{K}_\sk) \to \A \, \big | \,
T_{|\Pi(\Sigma_\sk)}=f,\,\, T_{|\Pi(\Sigma'_\sk)}=f', \,\,  T_{|\Pi(M_\sk)}=H,  \,\,  T_{|\Pi(M'_\sk)}=H'\right\}\right |
\\
\left | \left \{T'\colon \Pi({M'}_\sk) \to \A \, \big | \, T'_{|\Pi(\Sigma)}=f, \quad T'_{|\Pi(\Sigma')}=f' \right\}\right |
\\\, \prod_{n=1}^\infty \big (\prod_{m=0}^\infty |A_{m+n}|^{\L\big(m, \hat{K}, \fr(\hat{W}_\K)\big) }  \big)^{(-1)^n}
\prod_{n=1}^\infty \big (\prod_{m=0}^\infty |A_{m+n}|^{\L\big(m,M',\Sigma \sqcup \Sigma'\big)}  \big)^{(-1)^n}.
 \end{multline*}
 \end{Corollary}

\subsubsection{Dependence of  the formulae on  simplicial stratifications}\label{rem:2-triangulation_independence} We now address the dependence of the formulae in Theorem \ref{Thm:fin_Quinn_xcomp} on the choice of simplicial stratifications. We freely use \S\ref{Sec-Decorations:ind}, where  the dependence of $ \tFQd{\Bc}\colon \tdcob{n}{\Bc} \to \vProfGrpfin$  on decorations of $n$-manifolds was discussed.

Let $\Sigma$ be a closed (and, as usual, smooth) $n$-manifold. If we choose different simplicial stratifications,  $\zeta_{\Sigma} \colon |X_\Sigma| \to \Sigma$ and $\zeta'_\Sigma \colon |X'_\Sigma| \to \Sigma$, of $\Sigma$, then the corresponding $B_\A$-decorations of $\Sigma$,
$\fd_{\Sigma}(\zeta_{\Sigma},\A)$ and  $\fd_{\Sigma}(\zeta'_\Sigma,\A)$, as defined in Equation \eqref{eq:dec},
will, in general, be different. Nevertheless, we have an invertible profunctor,
$$\Psi\big( \fd_{\Sigma}(\zeta_{\Sigma},\A), \fd_{\Sigma}(\zeta'_\Sigma,\A) \big)\colon \tFQd{B_\A}\big (\Sigma, \fd_{\Sigma}( \zeta_{\Sigma},\A)\big)\bto
\tFQd{B_\A}\big(\Sigma, \fd_{\Sigma}\big (\zeta'_\Sigma,\A)\big).$$
This profunctor is natural with respect to the profunctors associated to cobordisms.

On the other hand, the formula, in Equation \eqref{profunctor:Q}, for the profunctors associated to an $(n+1)$-cobordism, $(i,M,i')\colon \Sigma \to \Sigma'$, does not depend on the simplicial stratification, $g_M\colon |Y_M| \to M$, of $M$,  extending that of $\Sigma$ and  $\Sigma'$, as shown in Equation \eqref{eq:simp_stract_cob}. Note that the simplicial stratifications of $\Sigma$ and $\Sigma'$  were part of the given data and so are themselves fixed.

Likewise,  in Equation \eqref{formula_matrix_2els}, the formula for matrix elements associated to the natural transformation of profunctors provided by an $(n+2)$-extended cobordism, $\K$ in \eqref{diag_for_K}, required the choice of a diagram $\hat{W}_\K$ of simplicial sets, and a homeomorphism of diagrams, $g\colon |\hat{W}_\K| \to \hat{\K}$, as shown in \eqref{simp_stract_for_K}.
The value in \eqref{formula_matrix_2els}, however, depends neither on the simplicial stratifications of $M$ and $M'$, extending those of $\Sigma$ and $\Sigma'$,  nor on the simplicial stratification of the reduced $(n+2)$-cobordism $\hat{K}$, extending those of $M$ and $M'$ that  $g\colon |\hat{W}_\K| \to \hat{\K}$ gives.

 \subsubsection{The bifunctor, $\tFQtr{\A} \colon {\trcob{n}} \to  \vProfGrpfin $ }\label{eTQFT-colA}
Some of discussion concerning  the (in)dependence  of the formulae for the finitary once-extended Quinn TQFT, with respect to the simplicial stratifications, can be repackaged inside a new  version of the finitary once-extended Quinn TQFT, that we now address.

We first define  a variant, ${\trcob{n}}$, of the bicategory $\tcob{n}$.
\begin{itemize}[leftmargin=0.6cm] \item The objects of  ${\trcob{n}}$ are pairs, $(\Sigma, \zeta_{\Sigma} )$, where $ \zeta_{\Sigma} \colon |X_\Sigma| \to \Sigma$ is a simplicial stratification of the (closed and smooth) $n$-manifold $\Sigma$;
\item  1-morphisms,
$(\Sigma, \zeta_{\Sigma}) \to
  (\Sigma', \zeta'_{\Sigma'}), $
  are $(n+1)$-cobordisms, $(i,M,j)\colon \Sigma \to \Sigma'$,
   \item  2-morphisms,  in $\trcob{n}$, denoted,
$$\big ( (i,M,j)\colon (\Sigma, \zeta_{\Sigma} )  \to  (\Sigma', \zeta'_{\Sigma'} ) \big) \\ \Longrightarrow  \big ( (i',M',j')\colon (\Sigma, \zeta_\Sigma )  \to  (\Sigma', \zeta'_{\Sigma'}) \big),
 $$
 are given by equivalence classes of extended cobordisms, as for $\tcob{n}$,
$$\K\colon \big ((i,M,j)\colon \Sigma  \to \Sigma' \big)  \Longrightarrow \big ((i',M',j')\colon \Sigma  \to \Sigma' \big) .
  $$  
  \end{itemize}
 The rest of the bicategory structure for $\trcob{n}$ is induced from that of the  bicategory $\tcob{n}$, in the obvious way, as in the construction of the bicategory $\tdcob{n}{\Bc}$, in Definition \ref{FinQ}.
  
  Given a finite crossed complex, $\A$, we therefore have a bifunctor,
  $$\mathcal{V}^\A\colon \trcob{n} \to \tdcob{n}{B_\A},$$ which, on objects, is such that, using the notation  $\fd(\zeta_{\Sigma},\A )$ of Equation \eqref{eq:dec},
  $$ \big (\Sigma, \zeta_{\Sigma}\big)\stackrel{\mathcal{V}^\A}{\longmapsto} \big (\Sigma, \fd(\zeta_{\Sigma},\A ) \big),$$  on 1-morphisms,   $$\mathcal{V}^\A\big ( (i,M,j)\colon (\Sigma, \zeta_{\Sigma} )  \to  (\Sigma', \zeta'_{\Sigma'})\big)= (i,M,j)\colon \big(\Sigma, \fd(\zeta_{\Sigma},\A )  \big)  \to  \big(\Sigma', \fd(\zeta'_{\Sigma'},\A )  \big) ,$$
  and analogously for 2-morphisms.

  The symmetric monoidal structure of $\tdcob{n}{B_\A}$, which is naturally derived from that of $\tcob{n}$, was briefly explained at the end of Subsection \ref{Quinn_is_sym_mon}. In particular, the tensor product of two $B_\A$-decorated $n$-manifolds is $$(\Sigma,\fd_\Sigma)\otimes (\Sigma',\gd_{\Sigma'})=(\Sigma \sqcup \Sigma', \fd_\Sigma \otimes \gd_{\Sigma'}),$$ where
$$ \fd_\Sigma \otimes \gd_{\Sigma'}:=\big\{ \langle \phi, \phi'\rangle \mid \phi \in \fd_\Sigma \textrm{ and }  \phi' \in \gd_{\Sigma'}\big\}.  $$
(Here, given $\phi\colon \Sigma \to B_\A$ and $\phi'\colon \Sigma'\to B_\A$, $\langle \phi, \phi'\rangle \colon \Sigma \sqcup \Sigma' \to B_\A$ is defined from the universal property of disjoint unions.)

We can define a symmetric monoidal structure in \smash{$\trcob{n}$}, where the tensor product of two closed, smooth, $n$-manifolds, $\Sigma$ and $\Sigma'$, provided with simplicial stratifications, $\zeta_\Sigma\colon |X_\Sigma| \to \Sigma$ and $\zeta'_{\Sigma'}\colon |X'_{\Sigma'}| \to \Sigma'$, is given by
$$(\Sigma, \zeta_\Sigma) \otimes (\Sigma', \zeta'_{\Sigma'}):= \big(\Sigma \sqcup \Sigma',  (\zeta_\Sigma\sqcup' \zeta'_{\Sigma'})\colon |X_\Sigma \sqcup X'_{\Sigma'}| \to \Sigma \sqcup \Sigma' \big), $$
where, explicitly, the homeomorphism $\zeta_\Sigma\sqcup' \zeta'_{\Sigma'}$ is defined as the composite
$$ |X_\Sigma \sqcup X'_{\Sigma'}| \ra{\cong} |X_\Sigma| \sqcup |X'_{\Sigma'}| \ra{\zeta_\Sigma\sqcup \zeta'_{\Sigma'}} \Sigma\sqcup \Sigma' .$$
From the fact that $\Pi(X_\Sigma \sqcup X'_{\Sigma'})\cong \Pi(X_\Sigma) \sqcup \Pi(X'_{\Sigma'})$, it can moreover be proved that $\mathcal{V}^\A$ is compatible with the symmetric monoidal structures of $\trcob{n}$ and $\tdcob{n}{B_\A}$.

This discussion leads to the following:
\begin{Theorem}\label{Thm:e-triang_independence} Let $\A$ be a finite crossed complex. There is a (symmetric monoidal) bifunctor, denoted
$$\tFQtr{\A} \colon {\trcob{n}} \to \vProfGrpfin,$$
which is defined as the following composite of  bifunctors,
$$\trcob{n} \ra{\,\, \mathcal{V}^\A \,\, } \tdcob{n}{B_\A} \ra{ \,\, \tFQd{B_\A} \,\,} \vProfGrpfin.$$
 \end{Theorem}
 Note that $\tFQtr{\A}$ is now decorated with a crossed complex $\A$, rather than with its classifying space $B_\A$. This is because the step
$\mathcal{V}^\A$ depends on the crossed complex $\A$, and not only on its classifying space.
\subsubsection{Morita-valued  once-extended TQFTs from finite crossed complexes}\label{sec:mor-x-complexes} As above, $\A$ denotes a fixed finite crossed complex.
 Explicit formulae for the Morita-valued version of the once-extended Quinn TQFT, in \S\ref{sec:Mor_valued-eTQFT},
   $\tFQmor{B_\A}\colon \tdcob{n}{B_\A} \to \Mor$,
can be derived from Theorem \ref{Thm:fin_Quinn_xcomp}, by applying the general constructions from Subsection \ref{sec:Mor_ext}. Let us give some brief details.

Passing from groupoids, $\Gamma$, to their groupoid algebras, $\tLin(\Gamma)$,  as in \S \ref{Groupoid algebra}, and with   $\Sigma$ a closed smooth $n$-manifold, $\zeta_{\Sigma} \colon |X_\Sigma| \to \Sigma$ being a simplicial stratification of $\Sigma$,  we have a canonical isomorphism of finite dimensional algebras,
 \begin{equation}
 \tFQmor{B_\A}\big ( \Sigma, \fd_{\Sigma}(\zeta_{\Sigma},\A)\big)\cong \tLin\big (\pi_1\big(\CRS(\Pi(X_\Sigma),\A)\big )\big).
 \end{equation}

These finite dimensional algebras  associated to a closed $n$-manifold, $\Sigma$, with a simplicial stratification, depend, explicitly,  on the chosen simplicial stratification of $\Sigma$. This dependence is, however, in  a quite `mild' way,  exactly as for the case of $\tFQd{B_\A}$ outlined in  \S \ref{rem:2-triangulation_independence} and Remark \ref{Decorations:indMod}, each using both the ideas and approach, given in \cite[10.3]{Bullivant-thesis}.
If we choose two simplicial stratifications,  $\zeta_{\Sigma} \colon |X_\Sigma| \to \Sigma$ and  $\zeta'_\Sigma \colon |X'_\Sigma| \to \Sigma$, of $\Sigma$, then there exists a canonically defined invertible bimodule,
 $$ \overline{\Psi}\big( \fd_{\Sigma}(\zeta_{\Sigma},\A), \fd_{\Sigma}(\zeta'_\Sigma,\A) \big)\colon \tFQmor{B_\A}\big(\Sigma, \fd_{\Sigma}(\zeta_{\Sigma},\A)\big) \bto   \tFQmor{B_\A} \big  (\Sigma, \fd_{\Sigma}(\zeta'_\Sigma,\A)\big),$$
connecting the  algebras thus obtained.
By construction, these bimodules compose well if we make further changes to the simplicial stratification and  are natural with respect to the bimodules associated to cobordisms, $(i,M,j)\colon \Sigma \to \Sigma'$, where both $\Sigma$ and $\Sigma'$ have a simplicial stratification, and hence a given $B_\A$-decoration.

As before, we have,
\begin{Theorem}\label{barMorTQFT} We have a, symmetric monoidal, bifunctor,
$$\tFQmortr{\A} \colon {\trcob{n}} \to \Mor,$$
 obtained by the following composite of bifunctors,
$$\trcob{n} \ra{\,\,\mathcal{V}^\A\,\,} \tdcob{n}{B_\A} \ra{\, \tFQd{B_\A}\,} \vProfGrpfin \ra{\,\tLin\,} \Mor.$$
 \end{Theorem}

\subsubsection{Absolute once-extended TQFTs derived from finite crossed complexes}\label{abs.TQFT}
Let  $\A$  be a finite crossed complex.
 As was noted earlier, the once-extended TQFTs,
\begin{align*}
&\tFQtr{\A} \colon {\trcob{n}} \to \vProfGrpfin &&& \textrm{ and } &&&&
\tFQmortr{\A} \colon {\trcob{n}} \to  \Mor,
\end{align*}
do not give a value to a closed smooth $n$-manifold $\Sigma$, unless $\Sigma$ is given a simplicial stratification.
In order to construct bifunctors whose domain is $\tcob{n}$, and whose target, unlike that of the once-extended Quinn TQFT, of Definition \ref{def:once_ext_Quinn}, $\tFQ{B_\A}\colon \tcob{n} \to\vProfGrphf,$ only outputs finite groupoids and finite dimensional algebras, we must specify a symmetric monoidal bifunctor, from   $\tcob{n}$ to $\trcob{n}.$

If $n=0$, this is  easy to do, as $0$-dimensional manifolds have only one simplicial stratification,
so we have a symmetric monoidal bifunctor  $\tcob{0}\to \trcob{0}$.

For $n\ge 1$, in order to construct a symmetric monoidal bifunctor  $\tcob{n} \to \trcob{n}$, we pick a simplicial stratification of each connected compact smooth manifold $\Sigma$, and then, if $\Sigma'$ is a not-necessarily connected manifold, the decomposition  of $\Sigma'$ into path-components provides a simplicial stratification of $\Sigma'$. (As we mentioned at the beginning of Section \ref{sec:TQFTS_xcomp}, this step is non-canonical, as we are using the choice axiom for classes.) This gives a symmetric monoidal bifunctor,  $\tcob{n} \to \trcob{n}$.

 \begin{Theorem}\label{thm.abs} For a non-negative integer $n$, we have once-extended TQFTs,
 \begin{align*}
 &\tFQch{\A} \colon {\tcob{n}} \to \vProfGrpfin &&\textrm{ and }
 &&&\tFQmorch{\A} \colon {\tcob{n}} \to  \Mor.
\end{align*}
 They can be `normalised' so that, if $\{\Sigma_k\}_{k \in \mathcal{K}}$ is any chosen set of path-connected closed smooth $d$-manifolds, and we have selected simplicial stratifications of each manifold $\Sigma_k$, namely $\zeta_{\Sigma_k}\colon |X_{\Sigma_k}| \to \Sigma_k$, then, for each $k$,
\begin{align*}
  \tFQch{\A}(\Sigma_k)& \cong \pi_1\big(\CRS(\Pi(X_{\Sigma_k}),\A)\big) \,\, \textrm{ and } \,\,
  \tFQmorch{\A}(\Sigma_k  ) \cong \tLin\big(\pi_1\big(\CRS(\Pi(X_{\Sigma_k}),\A)\big)\big).
\end{align*}
 \end{Theorem}
\begin{proof}
We compose the chosen symmetric monoidal bifunctor $\tcob{n} \to \trcob{n}$   with either
$\tFQtr{\A} \colon {\trcob{n}} \to \vProfGrpfin$  or
$ \tFQmortr{\A} \colon {\trcob{n}} \to  \Mor$. The remaining details are left to the reader.
 \end{proof}

\begin{Remark}\label{rem:Quinnfrom2Quinn} Under the conditions of the previous theorem, we note that we will always have that the state space of Quinn's finite total homotopy TQFT, $\FQp{B_\A}{s}\colon \tcob{n} \to \Vect_\C$, on $\Sigma$, is canonically isomorphic to the free vector space on the set of components of the groupoid $\pi_1\big(\CRS(\Pi(X_\Sigma) \big)$. In other words,
$$  \FQp{B_\A}{s}(\Sigma)=\C\big (\pi_0\big(\CRS(\Pi(X_\Sigma),\A)\big)\big);$$
see Theorem \ref{Theo:Quinn_exp_CW}. This makes it again clear in what sense the once-extended Quinn TQFT is a categorification of Quinn's finite total homotopy TQFT.
\end{Remark}

\subsection{Some explicit calculations for the once-extended  TQFTs derived from finite groupoids and 2-groups} \label{sec:exp_formulae_Quinn_groups}
For the rest of this paper we work over $\Q$. We give examples of the profunctors and algebras that the once-extended TQFTs,
\begin{align*}
\tFQtr{\A} \colon {\trcob{n}} &\to \vProfGrpfin, \quad \textrm{and} \quad
\tFQmortr{\A} \colon {\trcob{n}} \to  \Mor,
\end{align*}
assign to some $n$-dimensional manifolds.   These will be for low dimensions, $n=0,1,2$, and when $\A$ is the crossed complex given by a  finite group,  a finite groupoid, or  a crossed module of finite groups (equivalently a finite 2-group). The algebras we assign to loops and surfaces are particular cases of `tube algebras' considered in \cite{Bullivant_Tube}, \cite[Chapters 10 and 13]{Bullivant-thesis} and \cite[Section 3]{Bullivant_Excitations}, in the context of models for excitations of topological phases, derived from discrete higher gauge theory.
\subsubsection{The simplest example: the $(0,1,2)$-extended TQFT derived from a finite groupoid}\label{sec:012TQFT_fingroup}
Recall from Subsection \ref{sec:defXcomp} that we can think of a groupoid, $G$,  as a 1-truncated crossed complex, leading to a functor,  $\J_1\colon \Grp \to \Crs$.
Each $0$-dimensional manifold is trivially diffeomorphic to the disjoint union of copies of the singleton manifold, $\{*\}$.  On unpacking the construction in Subsection \ref{sec:xcomphom},  we have an isomorphism of groupoids $\pi_1\big(\CRS\big(\Pi(\{*\}),\J_1(G)\big)\big)\cong G$.

We can then consider the identification, $\trcob{1}\cong \tcob{1}$,  mentioned in  \S \ref{abs.TQFT}, and compose it with either of the bifunctors below,
 \begin{align*}
&\tFQtr{\J_1{(G)}} \colon \trcob{0}\to \vProfGrpfin &&&\textrm{ or }
&&&&\tFQmortr{\J_1{(G)}} \colon \trcob{0}\to  \Mor.
\end{align*}
Applying Theorem \ref{thm.abs}, we have the following result, essentially in \cite{LaudaPfeiffer,LaudaPfeiffer2,Teleman}.
\begin{Theorem}\label{thm:012-Grp}
For $G$ a finite groupoid,
there are once-extended TQFTs,
\begin{align*}
&\tFQch{\J_1(G)}\colon \tcob{0}
\to  \vProfGrpfin &&&\textrm{ and }
&&&& \tFQmorch{\J_1(G)}&\colon \tcob{0}  \to \Mor,
\end{align*}
such that
\begin{align*}&\tFQch{\J_1(G)}\big(\{*\}\big) \cong G &&\textrm{ and } 
&&&\tFQmorch{\J_1(G)}\big (\{*\}\big)\cong \tLin(G).
\end{align*}
Here $\tLin(G)$ is the groupoid algebra of $G$.
\end{Theorem}

\noindent The remaining parts of the specification of these $(0,1,2)$-extended TQFTs can be obtained from Theorem \ref{Thm:fin_Quinn_xcomp}. We will explain this for crossed modules in \S\ref{sec:012-xmod}.
\subsubsection{Frobenius algebras}\label{sec:frob} 
 That the once-extended TQFT, $\tFQmorch{\J_1(G)} \colon \tcob{0} \to \Mor $, in  Theorem \ref{thm:012-Grp}, exists is well known, see e.g. \cite[3.9 Remark]{Teleman}.  This follows from \cite[Theorems 3.52 / 3.5.4, in \S 3.8]{Schommer-Pries}, since  such $(0,1,2)$-extended  TQFTs are given by \emph{separable symmetric stellar Frobenius algebras}, and groupoid algebras of finite groupoids can be given such a structure; see \cite[Examples 5.1 and 5.2.]{LaudaPfeiffer}.
We give some details, following the conventions of \cite[
\S3.8]{Schommer-Pries}, and in the particular case when the stellar structure is given by a $*$-structure, as in \cite[Example 3.79]{Schommer-Pries}.

A \emph{symmetric $*$-Frobenius algebra}, $\A=(A,\cdot,1_A,\lambda, e,\dagger)$, by convention here over $\Q$, is given by an associative $*$-algebra,  $(A,\cdot,1_A,\dagger)$, with $1$, together with:
 \begin{itemize}[leftmargin=0.6cm]
  \item A $\Q$-linear map, $\lambda \colon A \to \Q$, satisfying that $\lambda ( a\cdot b)=\lambda(b\cdot a)$, for all $a,b\in A$, and also that $\lambda(a^\dagger)=\lambda(a)$, for each $a \in A$.
  \item An element, $e=\sum_{i} x_i \otimes y_i\in A\otimes A$, satisfying that given any $w \in A$, we have
  $$\sum_{i}(w  \cdot x_i) \otimes y_i= \sum_{i} x_i  \otimes ( y_i \cdot w).$$
 \end{itemize}
Moreover, $\lambda$ and $e$ should satisfy the following compatibility condition,
 $$\sum_{i}\lambda(x_i) \otimes y_i=1_A= \sum_{i} x_i  \otimes \lambda( y_i).$$

By direct application of the axioms, the following bilinear operation, $\circ$, on $A$,
$$a\circ b=\sum_i a \cdot x_i\cdot  b \cdot y_i,$$
gives an associative algebra, $(A,\circ)$, and the subspace, $[A,A]=\{a\cdot b-b\cdot a: a,b \in A\}$, is a bilateral ideal of $(A,\circ)$. Furthermore, we have an algebra homomorphism, $F_{\A}\colon (A,\circ) \to (A,\cdot)$, with $F_\A(a):=\sum_i x_i \cdot a \cdot y_i$, which factors through the centre $Z(A,\cdot)$, of $(A,\cdot)$, and which vanishes on $[A,A]$. In particular, $F_\A$ descends to an algebra homomorphism $F_\A'\colon (A,\circ)/[A,A]\to Z(A,\cdot)$.

Recall that an algebra with $1$, $(A,\cdot,1_A)$, is called \emph{separable} if there exists
 an element  $\overline{e}=\sum_{i} x'_i \otimes y'_i\in A\otimes A$,  satisfying that
  $\sum_i x_i' \cdot y_i' =1_A$, and also that
  $$\sum_{i}(w \cdot x'_i) \otimes y'_i= \sum_{i} x'_i  \otimes ( y'_i \cdot w),
  \textrm{ for any $w\in A$.}$$

 From the discussion in \cite[\S 3.8.5.]{Schommer-Pries}, it follows that if
 the underlying algebra with $1$, $(A,\cdot,1_A)$, of the symmetric Frobenius algebra, $\A$, is separable, then $ F_\A'\colon (A,\circ)/[A,A]\to Z(A,\cdot)$ is an isomorphism of associative  algebras.

After these preliminaries, we  now sketch how a separable symmetric $*$-Frobenius algebra $\A$ gives a (0,1,2)-extended TQFT, $\Z_\A\colon \tcob{0}\to \Mor$. Our conventions are non-standard, so as to have a clearer match with $\tFQmorch{\J_1(G)}\colon \tcob{0}  \to \Mor$.
\begin{enumerate}[leftmargin=0.6cm]
 \item If $\{*\}$ is a singleton manifold, then $\Z_\A(\{*\})=(A,\cdot,1_A)$. We shall associate $(A,\cdot,1_A)^{\otimes n}$ to all $n$-fold disjoint unions of $\{*\}$.
\item The bimodules given by the cobordisms, $\id_{\{*\}}\colon \{*\}\to \{*\}$, $\theta_{\{*\}}\colon \{*\}\to \{*\}$ and $( \id_{\{*\} }\sqcup \id_{\{*\}})\colon \{*\} \sqcup \{*\} \to  \{*\} \sqcup \{*\}$, depicted below,
$$\begin{tikzpicture}[xscale=0.6,yscale=0.5]
\draw plot [smooth] coordinates {(0,0) (3.5,0)};
\node at (-0.6,0) {$\{*\}$};
\node at (4.1,0) {$\{*\},$};
\node at (1.8,0.5) {$\id_{\{*\}}$};
\draw plot [smooth] coordinates {(7,0) (9.8,-1) (7.5,-1) (8.4,-0.65)};
\draw plot [smooth] coordinates {(10.5,0) (9.0,-0.5)};
\node at (6.5,0) {$\{*\}$};
\node at (11,0) {$\{*\},$};
\node at (8.8,0.3) {$\theta_{\{*\}}$};
\draw plot [smooth] coordinates {(14,0) (17.5,0)};
\draw plot [smooth] coordinates {(14,-1) (17.5,-1)};
\node at (13.4,0) {$\{*\}$};
\node at (18.1,0) {$\{*\}$};
\node at (13.4,-1) {$\{*\}$};
\node at (18.1,-1) {$\{*\}$};
\node at (15.8,0.5) {$\id_{\{*\}}\sqcup \id_{\{*\}}$};
\node at (18.6,-1) {$,$};
\end{tikzpicture}
$$
are given, respectively, by $A$, $A$ and $A\otimes A$, with actions, $a_1 \tr b \tl a_2=a_1 \cdot b \cdot a_2$, $ a_1 \tr b \tl a_2=a_2^\dagger\cdot b \cdot a_1^\dagger$, and  $(a'_1\otimes  a_1) \tr (b' \otimes b)\tl (a'_2\otimes  a_2)  =(a'_1  \cdot b' \cdot  a'_2) \otimes (a_2^\dagger \cdot  b  \cdot a_1^\dagger) $. In order that components are more clearly identified, we denote these actions graphically, for instance as below,
$$\vcenter{\xymatrix@C30pt{a_1 \ar@{-}[r]^{b} & a_2}}=\vcenter{\xymatrix@C30pt{1_A \ar@{-}[r]^{a_1\cdot b \cdot a_2 } &1_A}}, \quad \textrm{ and } \quad \vcenter{\xymatrix@=40pt@R=2pt{a'_1 \ar@{-}[r]^{b'} & {a_2'} \\ a_1 \ar@{-}[r]_{b} & {a_2} }} = \vcenter{\xymatrix@=40pt@R=2pt{1_A \ar@{-}[r]^{ a'_1 \cdot b'\cdot  a_2' } & 1_A \\ 1_A\ar@{-}[r]_{a_2^\dagger \cdot b\cdot  a_1^\dagger } & 1_A. }}$$

 \item We have cobordisms,  $$\emptyset \xrightarrow{{\cup}} \{*\} \sqcup \{*\}, \textrm{ graphically }\quad \vcenter{\xymatrix@R=1pt{\ast\ar@{-}@/_1pc/[d] \\\ast}}, \quad \textrm{ and } \quad   \{*\} \sqcup \{*\}  \xrightarrow{{\cap}} \emptyset, \textrm{ graphically } \vcenter{\xymatrix@R=1pt{\ast\ar@{-}@/^1pc/[d] \\\ast}} \quad .$$
The associated bimodules,
$$\Z_\A(\cup)\colon \Q \bto (A,\cdot,1_A)\otimes (A,\cdot,1_A), \quad \textrm{  and } \quad \Z_\A(\cap)\colon (A,\cdot,1_A)\otimes (A,\cdot,1_A) \to \Q$$ are both given by $A$, with actions,
\begin{align*}
b \tl (a_1\otimes a_2)&= a_2^* \cdot  b \cdot  a_1, \qquad \textrm{ graphically  }
 &\vcenter{\xymatrix@R=1pt@C=30pt{a_1 \ar@{-}@/_2pc/[d]_b \\ a_2}}&= \vcenter{\xymatrix@R=1pt@C=30pt{1_A \ar@{-}@/_2pc/[d]_{a_2^\dagger \cdot b \cdot a_1} \\ 1_A}},\\
 (a_1\otimes a_2) \tr b &=a_1 \cdot  b \cdot  a_2^*, \qquad \textrm{ graphically } &\vcenter{\xymatrix@R=1pt@C=30pt{a_1 \ar@{-}@/^2pc/[d]^b \\ a_2}} &= \vcenter{\xymatrix@R=1pt@C=30pt{1_A \ar@{-}@/^2pc/[d]^{a_1 \cdot b \cdot a_2^\dagger} \\ 1_A}} .\end{align*}

 \item  $\Z_A(S^1)=A/[A,A]$, considered as a vector space, and $\Z_\A(\varepsilon)\colon \Z_\A(S^1)\to \Q$ is induced by $\lambda\colon A \to \Q$. (The cobordism $\epsilon\colon S^1 \to \emptyset$, is defined in  \S \ref{sec:123TQFT_finXmod}.)

 \item The saddle,
 \begin{multline*}S\colon \big( (\id_{\{*\} }\sqcup \id_{\{*\}}) \colon \{*\} \sqcup \{*\} \to   \{*\} \sqcup \{*\} \big) \\ \implies \big( ( \cap \#_0 \cup) \colon \{*\} \sqcup \{*\} \to   \{*\} \sqcup \{*\} \big),
 \end{multline*}
 is sent to the following map of $\big((A,\cdot,1_A)\otimes (A,\cdot,1_A), (A,\cdot,1_A)\otimes (A,\cdot,1_A)\big)$-bimodules, each with underlying vector space $A \otimes A$: $b'\otimes b \mapsto \sum_i b'\cdot x_i \otimes b \cdot y_i$. We can visualise this map with the following diagram:
 $$\vcenter{\xymatrix@R2pt{ 1_A \ar@{-}[r]^{b'} & 1_A \\ 1_A \ar@{-}[r]^{b} & 1_A }}\mapsto \sum_{i} \vcenter{\xymatrix@R=2pt@C=20pt{  1_A \ar@{-}@/^2pc/[d]^{b' \cdot x_i} &&& 1_A \ar@{-}@/_2pc/[d]_{b \cdot  y_i} \\ 1_A  &&& 1_A . }}$$
\end{enumerate}

  Groupoid algebras, $\tLin(G)$, of finite groupoids are separable algebras, e.g. via  $$ \overline{e}=   \sum_{g\colon x \to  y} \frac{1}{N_x} (x\ra{g} y) \otimes (y\ra{g^{-1}} x),$$
where the sum is extended to all morphisms $g\colon x \to  y$ in $G$, and given an object $x$, in $G$, $N_x$ is the number of morphisms in $G$ with source $x$.

 Given a finite groupoid, $G$, the data that makes $\tLin(G)$ a symmetric Frobenius $*$-algebra, \cite{LaudaPfeiffer}, is as shown below.
 \begin{itemize}[leftmargin=1cm]
\item $\lambda \colon \tLin(G)  \to \Q$ is defined by $ \lambda(x\ra{g} y)=\begin{cases}1, \textrm{ if } (x\ra{g} y)= (x \ra{\id_x} x)\\ 0, \textrm{ otherwise},  \end{cases}$
\item[]\hspace{-10mm}whilst, as above,
\item $e=  \displaystyle \sum_{g\colon x \to  y} (x\ra{g} y) \otimes (y\ra{g^{-1}} x) \in \tLin(G) \otimes \tLin(G)$,
\item $(x\ra{g} y)^\dagger=(y\ra{g^{-1}} x)$.
 \end{itemize}

Using the explicit formulae in Theorem \ref{Thm:fin_Quinn_xcomp}, the $(0,1,2)$-extended TQFT constructed, as in \cite[\S 3.8.5 and 3.8.6]{Schommer-Pries}, from the separable symmetric Frobenius $*$-algebra  $\big( \tLin(G),\lambda, e,\dagger\big)$ coincides with $\tFQmorch{\J_1{(G)}}\colon \tcob{0} \to \Mor$. We will show this in the coming section, in the  more general case of crossed modules.
\subsubsection{(0,1,2)-extended TQFTs from finite crossed modules}\label{sec:012-xmod}
Let $\Gc=(\d\colon E \to G, \trl)$ be a crossed module  (of groups). Recall that $\d(E)$ is a normal subgroup of $G$. We continue denoting the elements of $G/\partial(E)$ by $[g]$, where $g \in G$.

The following follows immediately. (We use the notation of Subsection \ref{sec:xcomphom}.)

\begin{Lemma}  We have isomorphisms of groupoids, $\CRS_1\big(\Pi(\{*\}),\J_2(\Gc)\big)\cong G$, and  $\pi_1\big (\CRS \big(\Pi(\{*\}),\J_2(\Gc) \big)\big) \cong G/\partial(E)$.
\end{Lemma}
As a consequence, in this context, we have the following result.
\begin{Theorem}Suppose that $\Gc$ is finite. The once-extended TQFTs,
 \begin{align*}
&\tFQch{\J_2(\Gc)} \colon\tcob{0}
\to  \vProfGrpfin &&& \textrm{ and } &&&&
\tFQmorch{\J_2(\Gc)}&\colon \tcob{0} \to \Mor ,
\end{align*}
are such that
\begin{align*}
&\tFQch{\J_2(\Gc)}(\{*\}) \cong   G/\d(E) &&& \textrm{ and } &&&& \tFQmorch{\J_2(\Gc)}(\{*\}) \cong \mathbb{Q}\big(G/\d(E)\big).
\end{align*}
\end{Theorem}

The remaining parts of the specification of $\tFQch{\J_2(\Gc)}$  and $\tFQmorch{\J_2(\Gc)}$ can be obtained from Theorem \ref{Thm:fin_Quinn_xcomp} / Corollary \ref{cor:matrix-els-CW}.
Let us give some details.

Define $\Z^0_\Gc\colon \tcob{0}
\to  \vProfGrpfin$ as below: $$\big(\Z^0_\Gc\colon \tcob{0}
\to  \vProfGrpfin\big)=\big(\tFQch{\J_2(\Gc)} \colon\tcob{0}
\to  \vProfGrpfin\big).$$
In order to determine the value of  $\Z_\Gc^0$ on  $\id_{\{*\}}\colon  \{*\}\to \{*\}$, $\cup\colon\emptyset\to \{*\}\sqcup \{*\},$ $\cap\colon \{*\}\sqcup \{*\} \to\emptyset$, $\id_{\{*\} \sqcup \{*\}}\colon \{*\} \sqcup \{*\}\to  \{*\} \sqcup \{*\}$, and  on the saddle, $S$, we consider their CW-decompositions indicated in the figure below, in the same order,
$$\vcenter{\xymatrix@R=1pt{ \ast\ar[r] &\ast}},  \quad  \vcenter{\xymatrix@R=1pt{ \ast\ar@/^1pc/[d]  \\\ast}}\quad \quad \quad \vcenter{\xymatrix@R=1pt{ \ast\ar@{<-}@/_1pc/[d] \\\ast}},  \qquad\vcenter{ \xymatrix@R=1pt{\ast\ar[r] & \ast\\ \ast & \ast\ar[l]}}, \quad \quad \vcenter{\xymatrix@C=40pt{  \ast \ar@/^1pc/[d]_{\boldsymbol{g}}\ar[r]^{\boldsymbol{p}}  &\ast \ar@{<-}@/_1pc/[d]^{\boldsymbol{h}} \\ \ast \ar@{<-}[r]_{\boldsymbol{q}}\ar@{{}{}{}}[ur]|{\boldsymbol{e}}  &\ast.  }}$$
In particular, our CW-decomposition of the saddle, $S$, has four 0-cells, four 1-cells, $\boldsymbol{g}$, $\boldsymbol{h}$, $\boldsymbol{p}$, and $\boldsymbol{q}$, and one 2-cell, $\boldsymbol{e}$, attaching along $\boldsymbol{ph^{-1}q g^{-1}}$.

Note that if $H$ and $H'$ are groups, then $\Vect$-profunctors, $H\bto H'$, are nothing but vector spaces with a left-$H$-representation and a right-$H'$-representation that are compatible.  Using the calculations at the end of \S\ref{sec:ICRS}, it follows that $\id_{*}\colon \{*\} \to \{*\}$ is sent to $\Q\big(G/\partial(E)\big)$, with the left and right actions of $G/\partial(E)$, by left and right multiplications. For the same reason, the cobordisms, $$ (\cap \#_0 \cup) \colon (\{*\}\sqcup \{*\}) \to (\{*\}\sqcup \{*\})\textrm{ and } \id_{\{*\} \sqcup \{*\}}\colon (\{*\}\sqcup \{*\}) \to (\{*\}\sqcup \{*\}),$$ are both sent to $\Q ( G/\partial E)\otimes \Q(G/\partial E)$, with actions, respectively,
\begin{align*}
 ([a]\otimes [b])\tr ([g]\otimes [h]) \tl ([c]\otimes [d])&=([a][ g] [b]^{-1})\otimes ([d]^{-1} [b] [c]),\\
 ([a]\otimes [b])\tr ([g]\otimes [h]) \tl ([c]\otimes [d])&=([a][ g] [c])\otimes ([c]^{-1} [b] [b]^{-1}).\end{align*}

 To compute the matrix element of $\Z_\Gc^0(S)$, note that crossed complex maps $\Pi(S_\sk)\to \J_2(\Gc)$ are in one-to-one correspondence with elements of $\{(g,h,p,q,e)\in G^4\times E: \partial(e)=ph^{-1}q g^{-1}\}$. We then have, by Corollary \ref{cor:matrix-els-CW}, that
\begin{align*}
\big\langle [g]\otimes [h ]& | \Z_\Gc^0(S) | [p]\otimes [q] \big \rangle\\&=\big|\{(g,h,p,q,e)\in G^4\times E: \partial(e)=ph^{-1}q g^{-1}\}\big| \frac{|\Orb_{E^{\op}}(p)| \, |\Orb_{E^{\op}}(q) |}{|E^2|}\\
&=\delta([q] [g], [h][p]\}/|\ker(\partial)|.
\end{align*}

The value of $\Z_\Gc^0\big ( \epsilon\colon S^1 \to \emptyset\big)$ was  determined in \S\ref{sec:123TQFT_finXmod}. In the current notation,
$$
 \langle [g] \mid \Z_{\Gc}(\epsilon) \mid 1\rangle=\delta([g],1_{G/\partial(E)})\,  |\ker(\partial)|.$$

From these calculations, it can then be seen that  $ \tFQmorch{\J_2(\Gc)}\colon \tcob{0}\to \Mor$ is obtained from the following  symmetric Frobenius algebra structure on the group algebra $\Q\big (G/\partial(E)\big)$ of $G/\partial(E)$,  with the star structure $[g]^\dagger=[g^{-1}]$:
\begin{itemize}
\item $\lambda \colon \tLin(G)  \to \Q$ is defined by $\lambda([g])=|\ker(\partial)|\, \delta\big([g],1_{G/\partial(E)}\big)$,
\item $e=  \frac{1}{|\ker(\partial)|}\displaystyle \sum_{[g] \in G/\partial(E)} [g] \otimes  [g]^{-1} \in \Q\big(G/\partial(E)\big) \otimes \Q\big(G/\partial(E)\big).$
 \end{itemize}

In particular, the $(0,1,2)$-extended TQFTs $\tFQmorch{\J_2(\Gc)}$ and  $\tFQmorch{\J_1(G/\partial(E))}$ are equivalent. Hence,
$(0,1,2)$-extended TQFTs derived from finite crossed modules are not more general than the  $(0,1,2)$-extended TQFTs derived from finite groups. It is an open problem whether (0,1,2)-extended TQFTs derived from  general homotopy finite spaces can similarly  be reduced to the finite group case.

\subsubsection{$(1,2,3)$-extended TQFT derived from finite crossed modules}\label{sec:123TQFT_fingroup}
 We now address $(1,2,3)$-extended TQFTs derived from a finite crossed module of groups, $\Gc=(\d\colon E \to G, \trl)$.
For simplicity, we only look at the oriented case.

Consider $S^1$ with a simplicial stratification, $\zeta_{S^1}\colon |X_{S^1}|\to S^1$, where the simplicial set, $X_{S^1}$, has a single $0$-simplex and a single non-degenerate 1-simplex, and hence such that the associated CW-decomposition of $S^1$ has unique 0- and 1-cells.
This gives  simplicial stratifications for arbitrary disjoint unions of $S^1$, by using the obvious disjoint unions of this  simplicial stratification.

In \S \ref{sec:G/GG}, we defined the crossed complex $G\sslash \Gc$, proved it to be isomorphic to $\CRS\big(\Pi(S^1_\sk),\J_2(\Gc)\big)$, and we also computed the groupoid $\pi_1(G\sslash \Gc)$.
 In particular, $\tFQtr{\J_2(\Gc)}\colon  \trcob{1}
\to  \vProfGrpfin$ and $\tFQmortr{\J_2(\Gc)}\colon \trcob{1} \to \Mor$ satisfy that,
\begin{align*}
&\tFQtr{\J_2(\Gc)}(S^1,\zeta_{S^1})\cong \pi_1(G\sslash \Gc ) &&& \textrm{ and } &&&&\tFQmortr{\J_2(\Gc)}(S^1,\zeta_{S^1})\cong \tLin\big (\pi_1 (G\sslash \Gc)\big).
\end{align*}

As a consequence, applying Theorem \ref{thm.abs}, we have the following.
\begin{Theorem}\label{Thm:TQFT_Qdouble}
We have once-extended TQFTs,
\begin{align}
&\tFQch{\J_2(\Gc)}\colon \tcob{1}
\to  \vProfGrpfin, &&\textrm{ and } &
&\tFQmorch{\J_2(\Gc)}\colon \tcob{1} \to \Mor.\label{eq:123calG}
\intertext{Their values in $S^1$ are, respectively,}
& \tFQch{\J_2(\Gc)}(S^1) \cong \pi_1(G\sslash \Gc ), &&\textrm{ and } &&\tFQmorch{\J_2(\Gc)}(S^1)\cong \tLin\big (\pi_1(G\sslash \Gc )\big).
 \end{align}
 \end{Theorem}

 An important special case is when $E$ is trivial. Combining Theorem \ref{Thm:TQFT_Qdouble} with Example \ref{CRSPS^1G}, it follows that, given a finite group $G$, we  have once-extended TQFTs,
\begin{align}
&\tFQch{\J_1{(G)}}\colon \tcob{1} \to  \vProfGrpfin, &&\textrm{ and } &&\tFQmorch{\J_1(G)} \colon \tcob{1} \to \Mor,\intertext{such that}
 &\tFQch{\J_1(G)}(S^1) \cong G\sslash G,
 &&\textrm{ and } &&\tFQmorch{\J_1(G)}(S^1)\cong \tLin(G\sslash G).
\end{align}

As we  recalled in Example \ref{ex.Qdoublefingroup}, the algebra $\tLin(G\sslash G)$ coincides with the quantum double of the group algebra of $G$; some extra discussion on this  is found in \cite{Willerton}, and also in \cite{loopy}. In particular, the argument leading to Theorem \ref{Thm:TQFT_Qdouble} gives another proof of (and provides a homotopy theoretical underpinning for) the fact that, if $G$ is a finite group, then there exists a Morita-valued (1,2,3)-extended TQFT sending $S^1$ to the quantum double of the group algebra of $G$, see \cite{Bartlett_etal,morton:cohomological:2015,MNS}.

The rest of the structures of the once-extended TQFTs in Theorem \ref{Thm:TQFT_Qdouble} can be obtained from the discussion in \S \ref{fin_Quinn_xcomp}. For proof-of-principle, let us compute,
$$\tFQch{\J_2(\Gc)}\colon \tcob{1}
\to  \vProfGrpfin,$$
on
some of the generators of $\tcob{1}$ in \cite{Bartlett_etal,Bartlett_Goosen}. \begin{figure}[ht!]
 \labellist
 \pinlabel $\Du=$ at -10 206
 \pinlabel $\Mu=$ at 595 206
 \pinlabel $\epsilon=$ at 1183 127
 \pinlabel $\beta=$ at 1783 127
 \pinlabel $\Au=$ at 2300 127
 \endlabellist
\centering
\includegraphics[scale=0.12]{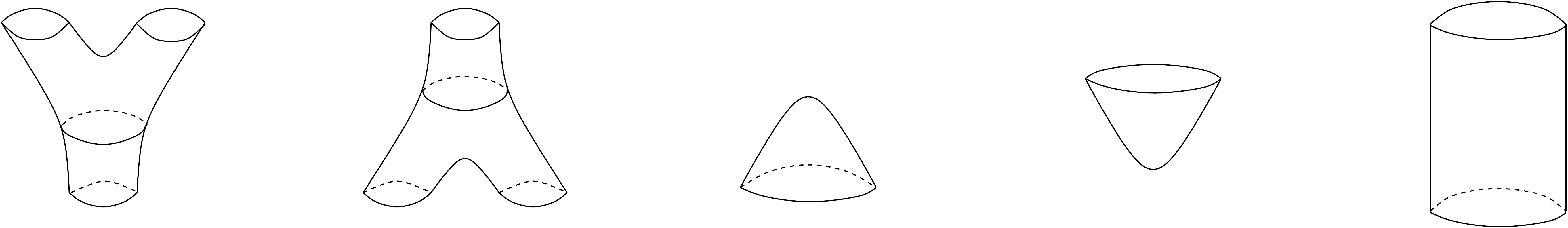}
\caption{Some 1-morphism of the bicategory $\tcob{1}$. \label{2cobgens}}
\end{figure}

On objects we already know that $S^1\mapsto \pi_1(G \sslash \Gc)$.
On 1-morphisms, for the cobordisms shown in Figure \ref{2cobgens} we have (recall that given a profunctor $F\colon G^\op\times H\to \Sets$, its linearisation, to a $\Vect$-profunctor, is denoted $\boldsymbol{F}\colon G^\op\times H\to \Vect$):
\begin{itemize}[leftmargin=0.6cm]
 \item $\Mu\colon  S^1\sqcup S^1 \to S^1$ is sent to $\boldsymbol{\Mp}\colon\pi_1(G \sslash \Gc)\times \pi_1(G \sslash \Gc) \bto \pi_1(G \sslash \Gc)$; see \S \ref{Interlude-db}.
 \item $\Du\colon  S^1\sqcup S^1 \to S^1$ is sent to $\boldsymbol{\Mp}^\dagger\colon  \pi_1(G \sslash \Gc) \bto\pi_1(G \sslash \Gc)\times \pi_1(G \sslash \Gc)$, where given a profunctor $P\colon \Gamma \to \Gamma'$, $P^\dagger\colon \Gamma' \bto \Gamma$ is $P^\dagger:=P\circ ( (-)^{-1},(-)^{-1})$.
 \item $\epsilon\colon S^1 \to \emptyset$ is sent to
 $\boldsymbol{\Ep}\colon \pi_1(G\sslash \Gc) \bto \{*\}$; see \S \ref{sec:eps_eta},
 \item $\beta\colon \emptyset \to S^1$ is sent to  $\boldsymbol{\Bp}\colon \{*\} \bto \pi_1( G \sslash \Gc)$, also in \S \ref{sec:eps_eta}.
 \item $\id_{S^1}\colon S^1 \to S^1$ is sent to $\boldsymbol{\Ap}\colon \pi_1(G \sslash \Gc)\bto \pi_1(G \sslash \Gc)$, defined in  \S\ref{sec:G/GG}.
 \end{itemize}

 As an example of the calculations associated to 2-morphisms, we show how to compute the natural transformation of profunctors given by the extended cobordism, $\mu\colon ( \epsilon \#_1 \beta\colon S^1 \to S^1  )\Rightarrow (\id_{S^1}\colon S^1 \to S^1)$,  on the left-hand-side of Figure \ref{2cob}.
\begin{figure}[ht!]
 \labellist
  \pinlabel $\mu=$ at -55 350
    \pinlabel $\hat{\mu}=$ at 1270 350
 \pinlabel ${\mathbf{a}}$ at 1429 484
  \pinlabel ${\mathbf{b}}$ at 1430 232
  \pinlabel ${\mathbf{g}}$ at 1673 662
  \pinlabel ${\mathbf{h}}$ at 1575 75
\pinlabel ${\mathbf{p}}$ at 1935 430
  \pinlabel ${\mathbf{c}}$ at 1850 336
\endlabellist
\centering
\includegraphics[scale=0.15]{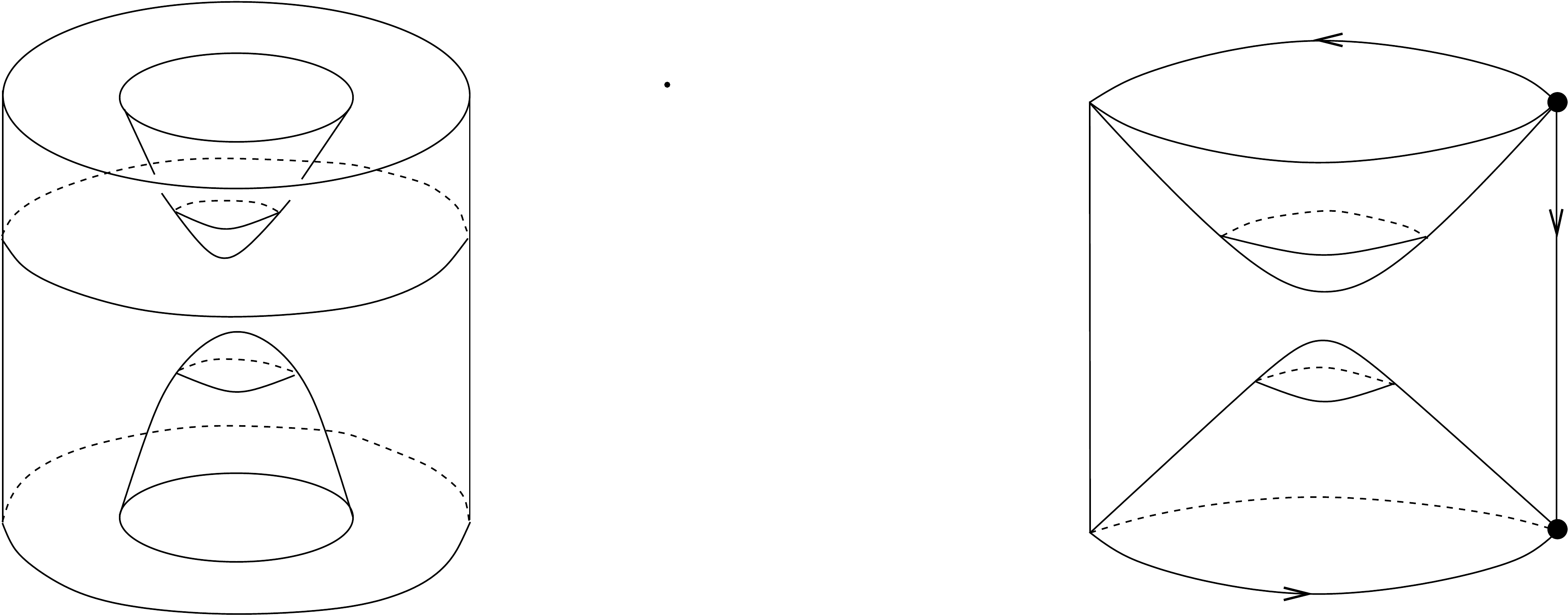}
\caption{The extended cobordism $\mu\colon ( \epsilon \#_1 \beta\colon S^1 \to S^1  )\Rightarrow (\id_{S^1}\colon S^1 \to S^1)$, and a CW-decomposition of its reduction, $\hat{\mu}$, so with two 0-cells, three 1-cells, $\mathbf{h}$, $\mathbf{g}$ and $\mathbf{p}$, three 2-cells, $\mathbf{a}$, $\mathbf{b}$ and $\mathbf{c}$, and one 3-cell. \label{2cob}}
\end{figure}

\noindent This gives a natural transformation
\begin{multline*}
\tFQch{\J_2(\Gc)}(\mu)\colon \big (\boldsymbol{\Ep\#_1\Bp}\colon  \pi_1(G\sslash \Gc)^\op \times  \pi_1(G\sslash \Gc)^\op \to \Vect\big )  \\ \Longrightarrow  \big( \boldsymbol{\Ap}\colon  \pi_1(G\sslash \Gc)^\op \times  \pi_1(G\sslash \Gc)^\op \to \Vect\big).
\end{multline*}

Morphisms $\Pi(\hat{\mu}_\sk) \to \J_2(\Gc)$ are in bijection with sequences $(g,h,p,a,b,c)\in G^3 \times E^3$, with $\partial(a)=g$, $\partial(b)=h$,  and $ac=b\trl p$.  Using Corollary \ref{cor:matrix-els-CW}, we have:
\begin{multline*}
 \Big\langle g \ra{a\otimes b} h \Big | \tFQch{\J_2(\Gc)}(\mu)  \Big| g \ra{[(p,c)]_g} h \Big\rangle \\=                                                    \frac{|\{(x,e)\in G\ltimes E: (x,e)\bullet g=h\}|}{|E|} \begin{cases}  1 , \textrm{ if } \partial(a)=g, \, \partial(b)=h, \,  ac=b\trl p,  \\
                                                                  0, \textrm{ otherwise}.                   \end{cases}
\end{multline*} This is exactly the natural transformation of profunctors, $\eta_{\hat{\mu}}$, discussed in  \S \ref{sec:eps_eta}


\begin{Remark}
We could consider CW-decompositions of $S^1$ with more than one 0-cell. The algebras and groupoids thereby obtained would then, in general, be different, however a natural Morita equivalence connects them all; see \S \ref{sec:mor-x-complexes}. For instance, the profunctor associated to moving from a cell decomposition of $S^1$ with one 0-cell to one with two 0-cells is the profunctor $\Ap_1^2\colon \pi_1(G\sslash \Gc)\bto \pi_1(G\times G)\sslash \Gc^{(2)})$ in \S \ref{sec:S2p}.
The case of $S^1$ decomposed by using multiple $0$-cells is reminiscent of the calculations in \cite[II-d]{lan-wen:topological:2014} and in \cite{Bullivant_Tube,Bullivant_Excitations}. Credit is due here to discussions with Alex Bullivant, including \cite[Theorem 10.3.2]{Bullivant-thesis}.
\end{Remark}

\subsubsection{$(2,3,4)$-extended TQFTs derived from finite crossed modules}\label{sec:234QFT_fingroup}
We will briefly discuss $(2,3,4)$-extended TQFTs derived from finite crossed modules, $\Gc$.

There is an infinite number of diffeomorphism classes of surfaces, thus an infinite amount of data\footnote{This issue will likely disappear if one further categorification level is introduced.} is \textit{a priori} required to write down the once-extended TQFTs,
\begin{align*}
&\tFQtr{\J_2{(\Gc)}} \colon\trcob{2}
\to  \vProfGrpfin &&& \textrm{ and }\quad  
&&&&\tFQmortr{\J_2(\Gc)}\colon \trcob{2}\to \Mor.
\end{align*}

In this paper, we will focus on the groupoids and  algebras assigned to $S^2$ and $T^2=S^1\times S^1$. (Recall they are only defined up to canonical profunctor / Morita equivalence.) We could consider all other surfaces (orientable and non-orientable), e.g. by choosing their usual CW-decompositions, with unique $0$ and $2$-cells.

Consider a simplicial stratification, $\zeta_{S^2}\colon |X_{S^2}|\to S^2$, of $S^2$, where  $X_{S^2}$ has a single $0$-simplex and a single non-degenerate 2-simplex. Let $S^2_\sk$ be the induced CW-decomposition of $S^2$, which has unique  $0$- and $1$-cells. In \S \ref{sec:S2}, we defined the crossed complex $\ker(\partial)\| \Gc$,  which is isomorphic to $\CRS(\Pi(S^2_\sk),\J_2(\Gc))$. Since $\pi_1(\ker(\partial) \| \Gc)=\ker(\partial)\sslash \big(G/\partial(E)\big)$, we have:
\begin{align*}
\tFQtr{\J_2(\Gc)}(S^2,\zeta_{S^2})&\cong \ker(\partial)\sslash (G/\partial(E)).
\end{align*}
This implies that if $G$ is a finite group, with group-algebra $\Q(G)$,
$$
\tFQtr{\J_1{(G)}}(S^2,\zeta_{S^2})\cong G, \ \quad  \text{ and }\ \quad \tFQmortr{\J_1{(G)}}(S^2,\zeta_{S^2})\cong \mathbb{Q}(G).
$$

We now determine the groupoid associated to the 2-torus, $T^2=S^1 \times S^1$.
There is a simplicial stratification of the 2-disk, $D^2$, with two non-degenerate 2-simplices, meeting along a diagonal edge. Identifying boundary edges in the usual way, this gives a simplicial stratification of the 2-torus, $T^2$,  here denoted $\zeta_{T^2}\colon |X_{T^2}|\to T^2$. It induces a CW-decomposition of the torus, with one 0-cell, three 1-cells and two 2-cells. We have in \S \ref{sec:torus} computed
$\CRS(\Pi(T^2_\sk),\J_2(\Gc))$ for the standard CW-decomposition of the torus, which however does not arise from a simplicial stratification. Let $T^2_{\sk'}$ be the torus, with the CW-decomposition induced by $\zeta_{T^2}$.
 Using analogous computations as in \S \ref{sec:torus}, we can see that $\pi_1(\CRS(\Pi(T^2_{\sk'}),\J_2(\Gc)))$ is isomorphic to the groupoid,  $\hat{T}^2(\Gc)$, below.
\begin{itemize}[leftmargin=0.6cm]
 \item The  objects of $\hat{T}^2(\Gc)$ are diagrams of the form below,
 (these can be interpreted as fake-flat discrete 2-gauge $\Gc$-configurations in $T^2_\sk$, see \cite[\S 3.5]{Companion}),
 $$\vcenter{\xymatrix@R=7pt@C=7pt{&\ast\ar[rr]^g\ar@{<-}[dd]_h  &&\ast \ar@{<-}[dd]^h \\ &&& \\
            &\ast\ar[rr]_g \ar[rruu]|{ {}_{\quad}^{\quad}v{}^{\quad}_{\quad} } \ar@{{}{ }{}}@/_1pc/[rruu]|{e'} \ar@{{}{ }{}}@/^1pc/[rruu]|{e}  &&\ast }}\quad \quad \vcenter{\xymatrix@R=1pt{&g,h,v \in G,   &e,e'\in E, \\
&\d(e)=g^{-1} h^{-1} v,  & \d(e')=v^{-1} gh.}}$$
\item The 1-morphisms of $\hat{T}^2(\Gc)$ are equivalence classes of arrows of  the form  below
where $g,h,v \in G$,  $e,e'\in E$, $x\in G$, and $a,b,c \in E$.
(These can be seen as gauge transformations between fake-flat 2-gauge $\Gc$-configurations, see \cite[\S 4.3.1]{Companion}.)
$$\vcenter{
\xymatrix@R=7pt@C=7pt{
&\ast\ar[rr]^g\ar@{<-}[dd]_h  &&\ast \ar@{<-}[dd]^h \\ &&& \\
            &\ast\ar[rr]_g \ar[rruu]|{{}_{\quad}^{\quad}v{}_{\quad}^{\quad}} \ar@{{}{ }{}}@/_1pc/[rruu]|{e'} \ar@{{}{ }{}}@/^1pc/[rruu]|{e}  &&\ast
       }
        }
        \ra{(x,a,b,c)} \hskip-0.8cm
            \vcenter{\xymatrix@C=5pc@R=1pc{&\ast \ar[rr]^{xg\d(a) x^{-1}}\ar@{<-}[dd]_{xh\d(b)x^{-1}}  &&\ast \ar@{<-}[dd]^{x\,h\,\d(b)\, x^{-1}} \\ &&& \\
            &\ast\ar[rr]_{x\, g\,\d(a)\, x^{-1}} \ar[rruu]|{x\, v \, \d(c)\,x^{-1}} \ar@{{}{ }{}}@/_1pc/[rruu]|>>>>>>>>>>>>>>{\qquad \qquad \qquad\big(c^{-1} \, e'\, (a\trl h) \,b \big)\trl x^{-1} } \ar@{{}{ }{}}@/^1pc/[rruu]|>>>>>>>>>>>>{\big (a^{-1} \ ( b^{-1} \trl g) \, e\, c\big) \trl x^{-1}\qquad \qquad }  &&\ast }}.$$
            \item Two arrows, given by $(x,a,b,c)$ and $(x',a',b',c')$ in $\hat{T}^2(\Gc)$, where $x,x' \in G$ and $a,a',b,b',c,c' \in E$, with the same source and target, as in the example above, are said to be equivalent if there exists a $ p \in E$ such that
$$(x',a',b',c')=\big(x \, \d(p), (p^{-1} \trl g) \, a \, p,  ( p^{-1} \trl h) \, b \, p,  (p^{-1} \trl v) \, c \, p\big).   $$
(In terms of discrete higher gauge theory, we are here identifying two gauge transformations when they differ by a 2-gauge transformation; see \cite{loopy}.)
\item[]\hspace{-12mm}Finally,
\item  the composition in the groupoid, $\mathcal{T}^2(\Gc)$, is induced by the semi-direct product, $G \ltimes_\trl ( E\times E \times E)$, with the product action of $G$, $(a,b,c)\trl g=(a\trl g, b\trl g, c \trl g)$.
\end{itemize}

\noindent In the particular case,  when $E$ is trivial, we have, for  $G$ a group,
$$\pi_1\big(\CRS( \Pi(T^2_{\sk'}), \J_1(G))\big)\cong \{(a,b)\in G^2: [a,b]=1 _G \}\sslash G,$$
where $g\trl (a,b)=(gag^{-1}, gbg^{-1})$.

We thus have, applying Theorem \ref{thm.abs}.
 \begin{Theorem} We have $(2,3,4)$-extended TQFTs,
\begin{align*}&\tFQch{{\J_2(\Gc)}}\colon \tcob{2}
\to  \vProfGrpfin &&& \textrm{ and } &&&& \tFQmorch{{\J_2(\Gc)}}\colon \tcob{2} \to \Mor.
\end{align*}
These can be normalised such that, for the 2-sphere $S^2$,
\begin{align*}
\tFQch{{\J_2(\Gc)}}(S^2)&\cong \ker(\d)\sslash (G/\d(E)), \quad
 &&\tFQmorch{{\J_2(\Gc)}}(S^2)\cong \tLin\big ( \ker(\d)\sslash \big (G/\d(E)\big ) \big ),
\end{align*}
and, on the 2-torus, $T^2$,
\begin{align*}
&\tFQch{{\J_2(\Gc)}} (T^2)\cong \hat{T}^2(\Gc), && \textrm{ and } &&
\tFQmorch{{\J_2(\Gc)}} (T^2) \cong \tLin\big (\hat{T}^2(\Gc)\big).
\end{align*}
\end{Theorem}
The remaining parts of the specification of the (2,3,4)-extended TQFTs, $\tFQch{\J_2(\Gc)}$  and $\tFQmorch{\J_2(\Gc)}$, can be obtained from Theorem \ref{Thm:fin_Quinn_xcomp} / Corollary \ref{cor:matrix-els-CW}.


\backmatter

\backmatter
\bibliographystyle{plainurl}

\phantomsection\bibliography{Quinn.bib}

\end{document}